\PassOptionsToPackage{table}{xcolor}
\documentclass[12pt,a4paper]{report}
\let\openright=\clearpage

\usepackage[a-2u]{pdfx}

\usepackage[utf8]{inputenc}

\usepackage{lmodern}

\usepackage{amsmath}        
\usepackage{amsfonts}       
\usepackage{amsthm}         
\usepackage{bbding}         
\usepackage{bm}             
\usepackage{graphicx}       
\usepackage{fancyvrb}       
\usepackage{natbib}         
\usepackage[nottoc]{tocbibind} 

\usepackage{dcolumn}        
\usepackage{booktabs}       
\usepackage{paralist}       
\usepackage{bbm}          
\usepackage{accents}      
\usepackage{float}        

\usepackage{tikz-cd}      
\newcommand{\nospaceperiod}{\makebox[0pt][l]{\,.}} 

\usepackage{amssymb}				

\usepackage{pinlabel}				

\usepackage{stackrel}             

\usepackage{nameref}				

\usepackage{todonotes}				

\usepackage{subcaption}				

\usepackage{mathtools}				

\usepackage{extarrows}				

\DeclareMathOperator{\sign}{sign}
\DeclareMathOperator{\degr}{deg}

\usepackage{hyperref}
\hypersetup{unicode}
\hypersetup{breaklinks=true}

\makeatletter
\def\@makechapterhead#1{
  {\parindent \z@ \raggedright \normalfont
   \Huge\bfseries \thechapter. #1
   \par\nobreak
   \vskip 20\p@
}}
\def\@makeschapterhead#1{
  {\parindent \z@ \raggedright \normalfont
   \Huge\bfseries #1
   \par\nobreak
   \vskip 20\p@
}}
\makeatother

\theoremstyle{plain}
\newtheorem{thm}{Theorem}
\newtheorem{cor}[thm]{Corollary}
\newtheorem{lemma}[thm]{Lemma}

\newtheorem{claim}[thm]{Claim}
\newtheorem{defn}[thm]{Definition}
\newtheorem{defn_thm}[thm]{Definition/Theorem}
\newtheorem{defn_lemma}[thm]{Definition/Lemma}
\newtheorem{conj}[thm]{Conjecture}

\newtheorem{Darboux}[thm]{Darboux's Theorem}
\newtheorem{Pfaff}[thm]{Pfaff's Theorem}
\newtheorem{Morse}[thm]{Morse Index Theorem}

\newtheorem{Thom}[thm]{Thom's Transversality Theorem}
\newtheorem{Sternberg}[thm]{Sternberg's Linearization Theorem}

\newtheorem{unstable}[thm]{Unstable manifold Theorem}

\newtheorem{Gronwall}[thm]{Gronwall's inequality}
\newtheorem{ExchangeA}[thm]{$(\ell_u+1)$-Exchange Lemma}
\newtheorem{ExchangeB}[thm]{$(\ell_u+\sigma+1)$-Exchange Lemma}

\newtheorem{stability}[thm]{Stability Lemma}
\newtheorem{Ehresmann}[thm]{Ehresmann's Theorem}
\newtheorem{ThomIsotop}[thm]{Thom's first isotopy Lemma}
\newtheorem{KupkaSmale}[thm]{Kupka-Smale Theorem}
\newtheorem{MorseHom}[thm]{Morse Homology Theorem}

\newtheoremstyle{named}{}{}{\itshape}{}{\bfseries}{.}{.5em}{\thmnote{#3's }#1}
\theoremstyle{named}

\theoremstyle{remark}
\newtheorem{rem}[thm]{Remark}
\newtheorem{example}[thm]{Example}

\newtheorem{rem_not}[thm]{Remark/Notation}
\newtheorem{n_example}[thm]{Non-Example}
\newtheorem{notat}[thm]{Notation}
\newtheorem{assump}[thm]{Assumption}
\newtheorem{setup}[thm]{Set-up}
\newtheorem{conv}[thm]{Convention}

\newenvironment{myproof}[2]{
  \par\medskip\noindent
  \textit{Proof of  {#1} {#2}}.}{\vspace{-\baselineskip}
\newline\rightline{$\qedsymbol$}
}

\newenvironment{sproof}{%
  \proof}{\endproof}
  
\newenvironment{pproof}{%
  \proof}{\endproof}
  
\newenvironment{hproof}{%
  \proof}{\endproof}
  
\newenvironment{iproof}{%
  \proof}{\endproof}

\DefineVerbatimEnvironment{code}{Verbatim}{fontsize=\small, frame=single}

\newcommand{\R}{\mathbb{R}}
\newcommand{\N}{\mathbb{N}}

\DeclareMathOperator{\codim} {codim}

\DeclareMathOperator{\End} {End}

\newcommand{\rot}{\hbox{rot}}

\newcommand{\ev}{{\rm ev}}

\newcommand{\fs}{\cdot_{f\hbox{-}s\,}}
\newcommand{\ffs}{{f\hbox{-}s\,}}

\newcommand\at[2]{\left.#1\right\vert_{#2}}

\newcommand{\evK}{\ev_{\bm{\mathcal{R}}}}

\newcommand{\ioo}{i^{\text{\textit{out-out}}}}

\begin{document}
\pagenumbering{roman}
     \setcounter{page}{1}
     \thispagestyle{empty}
     \begin{center}
	\centering
	\hspace{0pt}
    \vfill
    {\LARGE\bfseries A cord algebra for tori in three-space\par}
	\vspace{2cm}
	{\Large\bfseries Dissertation\par}
	\vspace{1.5cm}
	{\large zur Erlangung des akademischen Grades\par Dr. rer. nat.\par}
	\vspace{1.5cm}
	{\large eingereicht an der\par Mathematisch-Naturwissenschaftlich-Technischen Fakultät\par der Universität Augsburg\par}
	\vspace{1.5cm}
	{\large von\par}
	{\large\bfseries Marián Poppr\par}
	\vspace{2cm}
	{\large Augsburg, July 2025\par}
	\vspace{2cm}
\begin{figure}[!htbp]
\centering
\includegraphics[scale=0.45]{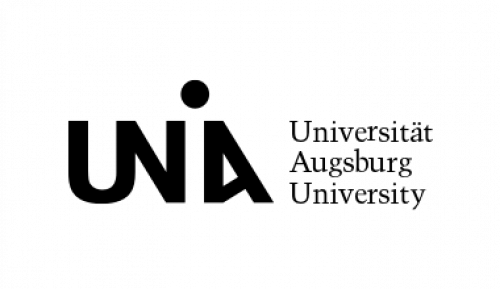}
\vspace{0.3cm}
\end{figure}
\end{center}
{
    \centering
    \hspace{0pt}
    \vfill
    \begin{tabular}{ r l }
        Betreuer:                   & Prof. Dr. Kai Cieliebak, Universität Augsburg \\
        Gutachter:                  & Prof. Dr. Tobias Ekholm, Uppsala University \\
    \\ \\
    \end{tabular}
    \newline
    \begin{tabular}{ r l }
        Tag der mündlichen Prüfung: & 14.11.2025
    \end{tabular}
}
\chapter*{Abstract}
Given a thin torus $T_{K, \varepsilon}$ around a knot $K\subset \R^3$, we construct Morse models of cord algebra $Cord(T_{K, \varepsilon})$ with $\mathbb{Z}$ and loop space coefficients. Using the Multiple time scale dynamics we identify $Cord(T_{K, \varepsilon}; \mathbb{Z})$ with $Cord(K; \mathbb{Z})$. In combination with the work of \cite{Cieliebak2016KnotCH, petrak2019definition} this indirectly relates $Cord(T_{K, \varepsilon})$ to $0$-th degree Legendrian contact homology $LCH_0(\mathcal{L}^\ast_+ T_{K, \varepsilon})$ of one component of unit conormal bundle over $T_{K, \varepsilon}$. Our definition of $Cord(T_{K, \varepsilon})$ is motivated by $J$-holomorphic curves with boundary on the Lagrangian submanifold $L^\ast_+ T_{K, \varepsilon}\cup\R^3$ with arboreal singularity along the torus $T_{K, \varepsilon}$.
\chapter*{Acknowledgements}

First and foremost, I would like to thank my supervisor, Kai Cieliebak, for many inspiring discussions and for encouraging me throughout all these years. Even though I kept calling our meetings ``consultations'', which sounded sort of official, it was especially a place where I could share any of my thoughts and always obtain very human advice.

I am grateful to my wife Zuzka for being with me on this adventure abroad. Without her unceasing support and love, it would not be possible. Also, I thank my parents, Andrea and Roman, for always being here for me. I cannot forget to express my gratitude to my grandparents Jarmila and Emil, brothers Vojta and Vratislav, and to Milada; this is also for you, grandma!

Next, I am thankful to Urs Frauenfelder for always having time for my questions, sometimes even missing lunch, as the canteen closed. I thank Yuchen Wang for fruitful discussions about dynamical systems. I also want to thank all other members of our symplectic group in Augsburg for all the fun and valuable math inputs. Namely (among others), I thank Shuaipeng Liu, Frederic Wagner, Emilia Konrad, Kevin Ruck, Zhen Gao, Miguel Pereira, Hanna Häußler, Filip Bro\'{c}i\'{c}, Evgeny Volkov, Airi Takeuchi, Milan Zerbin, Jennifer Gruber, and Lei Zhao.

I thank Shah Faisal, Francisco Javier Martínez Aguinaga, and Filip Strako\v{s} for friendly talks about math. I also thank Ji\v{r}\'{i} Zeman and Pavel H\'{a}jek for introducing me to life in Augsburg from a ``Czech perspective''. I would like to thank Yukihiro Okamoto for the discussions about various aspects of knot contact homology, it was enjoyable to talk about the topics so closely related to my thesis. I thank Pavel Nov\'{a}k for helping me with the English grammar. I thank my office mate, David Buchberger, for allowing me to expand to other tables and for explaining to me the basics of functional analysis.

During my studies, I also greatly benefited from the discussions with other people: I thank Roman Golovko, Christian Kuehn, Alberto Abbondandolo, Frank Sottile, and Joa Webber. I also thank Vivek Shende for suggesting to study $J$-holomorphic curves with the switches along the arboreal singularity $T_{K}$ of the Lagrangian submanifold $L^\ast_+ T_{K}\cup \R^3$. I am grateful to Tobias Ekholm, who agreed to be a reviewer of my thesis.

Next, I thank my floorball teammates from Augsburg for keeping me (somehow) fit. And finally, I thank my Czech friends for helping me not to forget the life in my home country.

\tableofcontents
\newpage

\pagenumbering{arabic}

\chapter{Introduction}
\label{ch:intro}

\section{Algebraic invariants from conormal lifts}
With any submanifold $N\subset \R^n$ we can associate its unit conormal bundle $\mathcal{L}^\ast N$ inside the unit co-sphere bundle $S^\ast \R^n$. $S^\ast\R^n$ together with the restricted Liouville form $pdq$ becomes a contact manifold and $\mathcal{L}^\ast N$ its Legendrian submanifold. In addition, an ambient isotopy of $N$ induces a Legendrian isotopy of its (unit) conormal lift $\mathcal{L}^\ast N$. Hence, for the study of $N$ we can use the toolkit of $J$-holomorphic curves type invariants associated to $\mathcal{L}^\ast N$ such as Legendrian contact homology (LCH).

Legendrian contact homology was proposed in \cite{SFT00} as the simplest invariant of the general framework of (relative) symplectic field theory. It arises as the homology of a certain differential graded algebra generated by words of Reeb chords, and with the differential counting $J$-holomorphic curves with punctures asymptotic to Reeb chords. Later, it was proven in \cite{EES05, Ekholm2005LegendrianCH} that Legendrian contact homology is well-defined for Legendrian submanifolds for $1$-jet spaces, and hence also of $J^1(S^{n-1})\cong S^\ast\R^n$.

In the last decades, $LCH(\mathcal{L}^\ast N)$ was extensively studied in the case when $N$ is a knot $K\subset\R^3$. Namely, L. Ng gave in a series of papers \cite{ng2005knot, Ng_2005, ng2007framedknotcontacthomology} a combinatorial description of $LCH(\mathcal{L}^\ast K)$ and among others introduced the cord algebra $Cord^{top}(K)$. The cord algebra $Cord^{top}(K)$ is a computable invariant of $K$ which is isomorphic to $LCH_0(\mathcal{L}^\ast K)$ and is generated by paths with endpoints on $K$ modulo certain skein relations. In particular, for framed oriented knots, the cord algebra detects the unknot. 

Later, it was verified in \cite{ekholm2013knot} that the combinatorial description agrees with the count of $J$-holomorphic curves from $LCH(\mathcal{L}^\ast K)$. The argument involved a representation of the knot $K$ in a more suitable form, as a braid along the unknot. This allowed to identify $LCH(\mathcal{L}^\ast K)$ with a count of certain gradient flow trees in the front projection $J^0(S^2)$.	

Next, in \cite{Cieliebak2016KnotCH} a refined version of $LCH_0(\mathcal{L}^\ast K)\cong Cord^{top}(K)$ was shown. Very roughly said, this more direct approach works as follows. We consider certain $J$-holomorphic curves with one puncture in a Reeb chord and the boundary mapped to two cleanly intersecting Lagrangian submanifolds - the conormal bundle $L^\ast K$ and $\R^3$. Then $Cord^{top}(K)$ can be translated from the topological boundary of these $J$-holomorphic curves.  

The approach of \cite{Cieliebak2016KnotCH} outlines a potential framework for the identification of $LCH(\mathcal{L}^\ast K)$ with its topological counterpart for more general manifolds than only $K$. As a consequence, in \cite{Ekholm_2017} it was shown that if $p$ is a point in $\R^3\setminus K$, then the enhanced Legendrian contact homology $LCH_0(\mathcal{L}^\ast K\cup \mathcal{L}^\ast p)$ is a complete knot invariant. We stress that this result uses the ``fully noncommutative'' version of $LCH$ with loop space coefficients.

In \cite{petrak2019definition}, a Morse-theoretic model of $Cord^M(K)$ with loop space coefficients was constructed. This model counts certain Morse flow trees given by the gradient flow of the standard energy functional on the space $K\times K$. Also, $Cord^M(K)$ is isomorphic to $Cord^{top}(K)$ (even though the whole proof did not appear in the literature, to experts the isomorphism is known).

Recently, in \cite{okamoto2022toward}, Y. Okamoto, using de Rham chains, proposed a topological model of $LCH(\mathcal{L}^\ast N)$ with $\R$-coefficients. In this case, $LCH(\mathcal{L}^\ast N)$ is of arbitrary degree and $N\subset \R^n$ is a closed oriented submanifold of arbitrary codimension. Then the same author studied in \cite{okamoto2024legendrian} the singular and Morse models of $LCH(\mathcal{L}^\ast N)$ with $\mathbb{Z}_2$ coefficients, provided that $N$ has a high codimension.

We would like to mention some other works on the knot case. Namely, $LCH(\mathcal{L}^\ast K)$ can be related to the framework of Ooguri-Vafa large N duality, see for example \cite{ekholm2020physics, aganagic2014topological, ekholm2018knot}. A different technique is microlocal sheaf theory, which was applied in \cite{shende2019conormal} to the proof that $\mathcal{L}^\ast K$ is a complete knot invariant. A relation between microlocal simple sheaves and augmentations of the DGA underlying $LCH(\mathcal{L}^\ast K)$ was described in \cite{gao2018simple}.

To the end, we mention also the work of J. Asplund \cite{asplund2019fiber} which studied $\mathcal{L}^\ast N$ using partially wrapped Floer homology in the spirit of wrapped Fukaya categories of Liouville sectors from \cite{Ganatra2018SectorialDF}.

\section{Main results}
To any knot $K\subset \R^3$ we can associate an $\varepsilon$-radius torus $T_{K, \varepsilon}$ which arises as the boundary of the $\varepsilon$-tubular neighborhood of $K$. Since $T_{K, \varepsilon}$ is a codimension $1$ submanifold, its (unital) conormal lift has two connected components $\mathcal{L}^\ast_\pm T_{K, \varepsilon}$. Let us take only one component; $\mathcal{L}^\ast_+ T_{K, \varepsilon}$. It is not hard to see that there is a canonical Legendrian isotopy $\mathcal{L}^\ast K\cong\mathcal{L}^\ast_+ T_{K, \varepsilon}$ (Lemma \ref{lemma_isotop_conormal}). Hence, in particular, $LCH_0(\mathcal{L}^\ast K)\cong LCH_0(\mathcal{L}^\ast_+ T_{K, \varepsilon})$.
As we stated above, by \cite{Cieliebak2016KnotCH, petrak2019definition} there is an isomorphism $LCH_0(\mathcal{L}^\ast K)\cong Cord^{M}(K)$. Thus, it is natural to expect that there exists a Morse-theoretic $Cord^{M}(T_{K, \varepsilon})$ which matches its symplectic counterpart $LCH_0(\mathcal{L}^\ast_+ T_{K, \varepsilon})$. In short, in this thesis, we are going to study the following diagram
\[
\begin{tikzcd}
LCH_0(\mathcal{L}^\ast K) \arrow[r, "\cong"] \arrow[d, "\phi_K"]           & LCH_0(\mathcal{L}^\ast_+ T_{K, \varepsilon})    \arrow[d, "\phi_{T_{K, \varepsilon}}"]     \\
Cord^M(K)\arrow[r, "\Theta^M"] & Cord^M(T_{K, \varepsilon})
\end{tikzcd}
\]

The core of the thesis is definition of the Morse model $Cord^M(T_{K, \varepsilon})$ with $\mathbb{Z}$ and \textit{loop space coefficients} and its relation to $Cord^M(K)$. In addition, this indirectly relates $Cord^M(T_{K, \varepsilon})$ to $LCH_0(\mathcal{L}^\ast_+ T_{K, \varepsilon})$. The map $\phi_{T_{K, \varepsilon}}$ will be only outlined in Chapter \ref{ch:rel_symplectic}. However, leaving first a word about $\phi_{T_{K, \varepsilon}}$ will enlight the definition of $Cord^M(T_{K, \varepsilon})$.

By $L^\ast_\pm T_K$ let us denote two half's of the conormal bundle $L^\ast T_K$ such that $\mathcal{L}^\ast_\pm T_K\subset L^\ast_\pm T_K$. Next, similarly to $\phi_K$, the map $\phi_{T_{K, \varepsilon}}$ will also use punctured $J$-holomorphic curves with one puncture asymptotic to a Reeb chord and from whose topological boundary one can expect to reconstruct $Cord^M(T_{K, \varepsilon})$. However, now the boundary of the $J$-holomorphic curves will be mapped to the Lagrangian submanifold $L^\ast_+ T_{K, \varepsilon}\cup \R^3$ with a simple arboreal singularity along $T_{K, \varepsilon}$ in the sense of \cite{nadler2017arboreal}. This reflects the choice of the single component of $\mathcal{L}^\ast T_{K, \varepsilon}$. As a consequence, the boundary of $J$-holomorphic curves $u$ will satisfy certain jet conditions along the zero section $\R^3$. Let us denote by $c_u$ the oriented path given by one connected component of $\partial u\cap \R^3$. In particular, after a choice of appropriate conventions, we conclude that the path $c_u$ will be outward-pointing from $T_{K, \varepsilon}$ at both of its endpoints.

Let us return, for now, back to $Cord^M(T_{K, \varepsilon})$. We aim to define $Cord^M(T_{K, \varepsilon})$ by the dynamics of a gradient-like flow of the standard energy functional $E_\varepsilon$ on $T_{K, \varepsilon}\times T_{K, \varepsilon}$. However, motivated by $J$-holomorphic curves we do not consider as an ambient space for the Morse dynamics of $Cord^M(T_{K, \varepsilon})$ the full product $T_{K, \varepsilon}\times T_{K, \varepsilon}$,  but rather only the subspace 
$$M_{K, \varepsilon}\subset T_{K, \varepsilon}\times T_{K, \varepsilon},$$
which represents the spaces of outward-pointing chords on $T_{K, \varepsilon}$, i.e., the space of oriented line segments on $T_{K, \varepsilon}$ that are pointing out from the torus at their endpoints.

\begin{thm}\label{thm_intro_MK}Let $K$ be generic and $\varepsilon>0$ be small. Then, outside of a small neighborhood of the diagonal, the space $M_{K, \varepsilon}$ is a $4$-dimensional manifold with corners.
\end{thm}

In addition, we will be able to describe the surprisingly rich topology of $M_{K, \varepsilon}$. In more detail, outside of certain special and diagonal points, $M_{K, \varepsilon}$ will look just like a fiber bundle over $K\times K$ with fibers diffeomorphic to the square. If we also include the special points, we obtain a broken fibration. We will compute $H_1^{sing}$ of the broken fibration and conclude that the space $M_{K, \varepsilon}$ itself remembers nontrivial information about the underlying knot $K$.

The topology of $M_{K, \varepsilon}$ becomes more dramatic near the diagonal $\Delta_{full}$, where $M_{K, \varepsilon}$ will become singular. As we will see later, this will significantly influence $Cord^{M}(T_{K, \varepsilon})$ with loop space coefficients. We will show that if $\varepsilon>0$ is small, then near the diagonal the space $M_{K, \varepsilon}$ looks approximately as a stratified fiber bundle over a half-torus with cuspical fibers.

In order to define $Cord^{M}(T_{K, \varepsilon}; \mathbb{Z})$ we construct a chain complex in degrees $1, 0$ which is generated by critical points of $E_\varepsilon\vert_{M_{K, \varepsilon}\setminus\Delta_{full}}$ of Morse indices $2$ and $1$, respectively. Then the differential will be counting certain Morse flow trees $\clubsuit_{X_{E_\varepsilon}}^{out-out}\subset M_{K, \varepsilon}\setminus\Delta_{full}$ given by the flow of a gradient-like vector field $X_{E_\varepsilon}$. Finally, $Cord^{M}(T_{K, \varepsilon}; \mathbb{Z})$ will be $0$-th homology of this chain complex. We remark that this version is not using the loop space coefficients, since we are not dealing with the diagonal. 

Similarly, A. Petrak's \cite{petrak2019definition} cord algebra $Cord^{M}(K; \mathbb{Z})$ was defined by chains generated by critical points of the energy functional $E_0$ on $K\times K$. The critical points of indices $1$ and $0$ were considered, and the differential counted Morse flow trees $\clubsuit_{X_{E_0}}$. Since there is a bijection between the critical points entering the cord algebras for $K$ and $T_{K, \varepsilon}$, one can guess that the algebras might be isomorphic. In fact, we will prove the following result.

\begin{thm}\label{thm_intro_isom} For generic $K$, generic gradient-like vector fields $X_{E_0}, X_{E_\varepsilon}$ and $\varepsilon>0$ sufficiently small the cord algebras $Cord^{M}(K; \mathbb{Z})$ and $Cord^{M}(T_{K, \varepsilon}; \mathbb{Z})$ are isomorphic on the chain level. 
\end{thm} 

\begin{iproof} Theorem \ref{thm_intro_isom} is an adiabatic limit problem with $\varepsilon\rightarrow 0$. Instead of working with collapsing tori $T_{K, \varepsilon}\times T_{K, \varepsilon}\subset \R^3\times \R^3$ we work rather with configuration spaces given by parametrizations of $K$ and $T_{K, \varepsilon}$, namely with $(\R/T\mathbb{Z})^2$ and $(\R/T\mathbb{Z}\times S^1)^2$  (note that $M_{K, \varepsilon}\subset(\R/T\mathbb{Z}\times S^1)^2$ also makes sense). Then the main tool for the description of the dynamics of $X_{E_\varepsilon}$ on $(\R/T\mathbb{Z}\times S^1)^2$ will be multiple time scale dynamics, see for example \cite{Kuehn2015MultipleTS}. In particular,
\begin{equation}\label{eqn_intro_fs}
\dot{u}_\varepsilon=X_{E_\varepsilon}
\end{equation}
is a fast-slow flow system with two possible time scales $t, \tau$, related by $\tau=t/\varepsilon$. Sending $\varepsilon$ to $0$ for different times scales, $t, \tau,$ we will obtain the slow and the fast subsystem. Ideally, for $\varepsilon>0$ small, one will be able to recover the dynamics $(\ref{eqn_intro_fs})$ from the interplay between the fast and slow subsystems, provided that $(\ref{eqn_intro_fs})$ is not too singular. Since now we are not dealing with the diagonal and we work with critical points of low degrees, we will be able to avoid the ``too singular'' behavior of $(\ref{eqn_intro_fs})$. Additionally, the dynamics of the slow subsystem will almost coincide with the dynamics of $X_{E_0}$ on $(\R/T\mathbb{Z})^2.$

In order to keep track of the dynamics of the trees $\clubsuit_{X_{E_\varepsilon}}^{out-out}$ (and $\clubsuit_{X_{E_0}}$) it will be not enough to use the standard techniques of Fenichel theory \cite{Fenichel1979GeometricSP} and Exchange Lemmata \cite{Brunovsk1999CrInclinationTF, Schecter2008ExchangeL2}. Indeed, the trees have bifurcations, which depend also on the exterior geometry of $T_{K, \varepsilon}$ and $K$ inside $\R^3$. So we will also need to mix the Fenichel theory and Exchange Lemmata with the perturbations of these bifurcations.

In the end, we remark that we also prove a weaker version of Theorem \ref{thm_intro_isom} using the Conley index, \cite{conley1978isolated, Salamon_1985conley}. In more detail, we show a $\hbox{mod }2$ correspondence between the heteroclinic orbits entering $Cord^{M}(K; \mathbb{Z})$ and $Cord^{M}(T_{K, \varepsilon}; \mathbb{Z})$, where we will additionally be able to capture their homotopy classes (as paths with fixed endpoints). An advantage of the Conley index argument is that it reflects well the compactness argument that we also use in a more involved fashion use also in the multiple time scale dynamics for the trees.
\end{iproof}

Now, we turn to the cord algebra with loop space coefficients $\mathbb{Z}[\lambda^{\pm1}, \mu^{\pm 1}]$ ($\cong\mathbb{Z}\pi_{1}(T_K)\cong H_\bullet(\Omega_\ast T_{K})$), where $\lambda$ corresponds to a longitude and $\mu$ to the meridian. Contrary to the knot case, the loop space coefficients in $Cord^M(T_{K, \varepsilon}; \mathbb{Z}[\lambda^{\pm1}, \mu^{\pm1}])$ appear much more naturally than in the knot case, because now the ends of the chords are actually living on $T_{K, \varepsilon}$. This becomes even more evident when we approach the diagonal. In the knot case, the chords at the diagonal are replaced by the slightly mysterious algebraic expression $1-\mu$. This cannot be seen from the Morse model and has to be justified by the behavior of $J$-holomorphic curves on $\mathcal{L}^\ast K$. On the other hand, in the $T_{K, \varepsilon}$ case, the algebraic expression $1-\mu$ will be dictated by the cuspidal behavior of $M_{K, \varepsilon}$ near the diagonal.

In order to define $Cord^{M}(T_{K, \varepsilon}; \mathbb{Z}[\lambda^{\pm1}, \mu^{\pm1}])$ we again have to first construct a chain complex $C_1^M(T_{K, \varepsilon})\xrightarrow{D_\varepsilon}C_0^M(T_{K, \varepsilon})$. This will consist of the following:
\begin{itemize}
\item  $C_1^M(T_{K, \varepsilon})$ is precisely the same as in the $\mathbb{Z}$-coefficient case - the $\mathbb{Z}$ vector space generated by $Crit_2(E_\varepsilon\vert_{M_{K, \varepsilon}\setminus \Delta_{full}})$.
\item However, $C_0^M(T_{K, \varepsilon})$ is a bit more involved. By $m_{\Delta_\varepsilon}$ we denote the unique global minimum of an auxiliary Morse function $h_{\Delta_\varepsilon}$ on the diagonal $\Delta_{full}$. Then $C_0^M(T_{K, \varepsilon})$ is the free unital noncommutative $\mathbb{Z}$-algebra generated by $Crit_1(E_\varepsilon\vert_{M_{K, \varepsilon}\setminus \Delta_\varepsilon})\cup\lbrace m_{\Delta_\varepsilon}\rbrace$ and extra generators $\lambda, \mu$ such that $\lambda$ and $\mu$ have inverses and commute together. Note that the Morse index of $m_{\Delta_\varepsilon}$ is $0$ and not $1$!
\item Finally, since $\Delta_{full}$ is a Bott-nondegenerate critical manifold of $E_\varepsilon$, the differential $D_\varepsilon$ will count cascade Morse flow trees $\clubsuit^{c, out-out}_{X_{E_\varepsilon}}$ together with the sign count of $\lambda, \mu$ as the chord-endpoints from the flow trees cross the meridian or the longitude.
\end{itemize}  

\begin{defn}$$Cord^M(T_{K, \varepsilon}; \mathbb{Z}[\lambda^{\pm 1}, \mu^{\pm 1}])=C_0^M(T_{K, \varepsilon}; \mathbb{Z}[\lambda^{\pm 1}, \mu^{\pm 1}])/\mathcal{I}_{T_{K, \varepsilon}},$$
where $\mathcal{I}_{T_{K, \varepsilon}}$ is the two-sided ideal of $C_0^M(T_{K, \varepsilon}; \mathbb{Z}[\lambda^{\pm 1}, \mu^{\pm 1}])$ generated by the image of the linear map $D_{\varepsilon}$ and $m_{\Delta_\varepsilon}-1$.
\end{defn}

Now, we conceptually explain how $Cord^M(T_{K, \varepsilon}; \mathbb{Z}[\lambda^{\pm 1}, \mu^{\pm 1}])$ counts the Morse flow trees near the diagonal. We remark that the following discussion will be only partially proven, but the author believes that the claims are supported by enough evidence. For more details, we refer to Chapter \ref{ch:cord}.

For simplicity, let us take $K$ to be the standard unknot $K_U$, which lies in the $xy$-plane. Let $p_\varepsilon\in Crit_2(E_\varepsilon\vert_{M_{K, \varepsilon}\setminus\Delta_{full}})$. Then the crucial claim is the following.

\begin{claim}\label{intro_eat}
$$\dim\big(W^u_{-\nabla E_\varepsilon}(p_\varepsilon)\cap W^s_{-\nabla E_\varepsilon}(\Delta_{full})\cap M_{K, \varepsilon}\big)=1$$
(instead of the expected $2$-dimensional intersection).
\end{claim}

\begin{figure}[H]
\labellist
\pinlabel $c_{p_\varepsilon}$ at 330 120
\pinlabel $\times$ at 178 255
\pinlabel $c_{triv}$ at 190 276
\endlabellist
\centering
\includegraphics[scale=0.4]{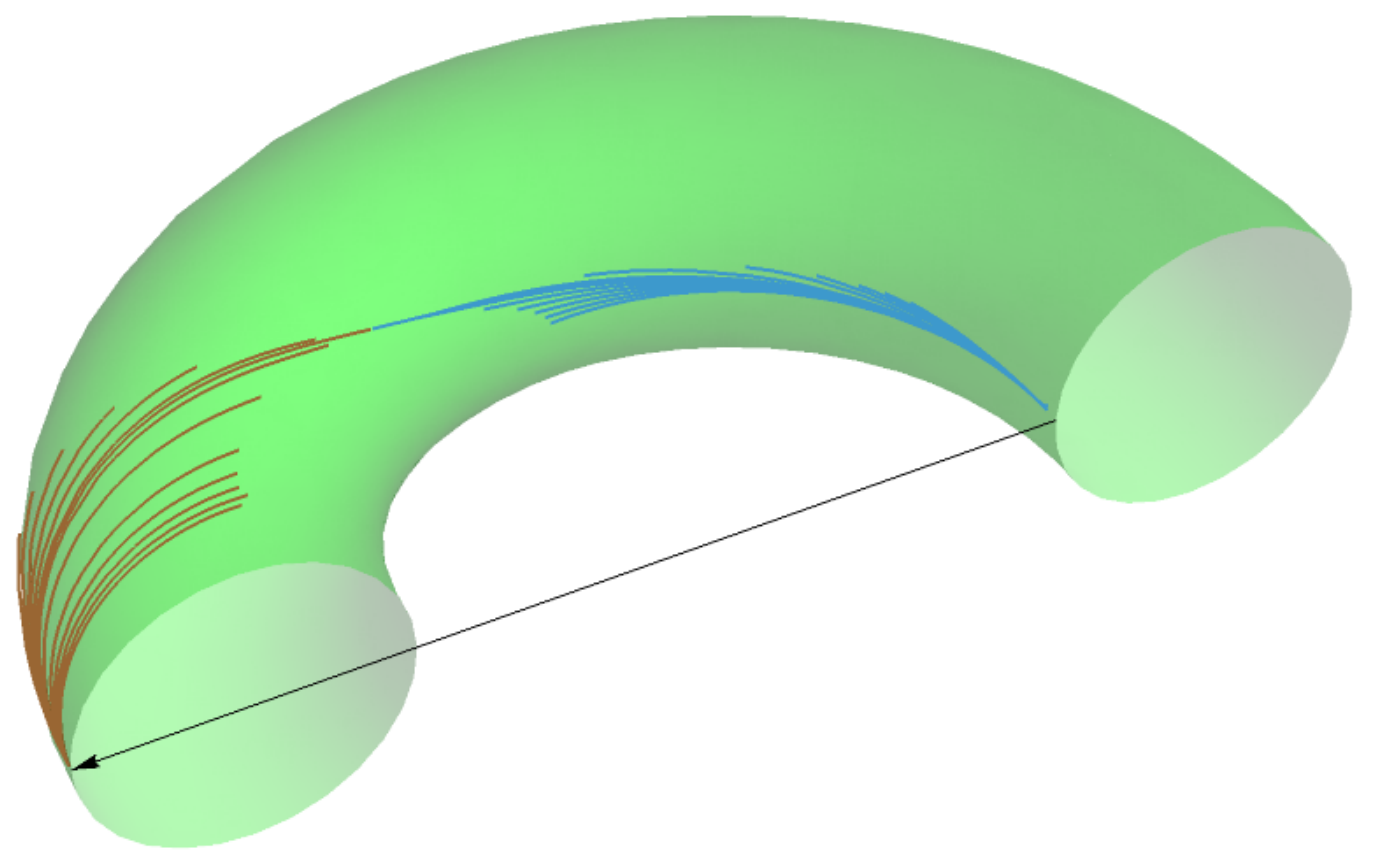}
\vspace{0.1cm}
\caption{A visualization of $W^u_{-\nabla E_\varepsilon}(p_\varepsilon)\cap M_{K_U, \varepsilon}$ on $T_{K_U, \varepsilon}$. The \textcolor{red}{red} curves depict the endpoints of chords emanating from the chord $c_{\varepsilon}$ under the negative gradient flow constrained to $M_{K_U, \varepsilon}$. Similarly, the \textcolor{blue}{blue} curves depict the starting points of chords emanating from the chord $c_{\varepsilon}$ under the negative gradient flow constrained to $M_{K_U, \varepsilon}$. 
Note also a single pair of red and blue curves that are connected by a trivial chord $c_{triv}$; they represent a trajectory $u_{\varepsilon}$ from $p_\varepsilon$ to $\Delta_\varepsilon$ that stays the whole time in $M_{K_U, \varepsilon}$.}
\label{figure_unstable_diagonal_intro}
\end{figure}

In other words, the cuspical singularity of $M_{K, \varepsilon}$ is decreasing the expected dimension of the intersection of stable and unstable manifolds by $1$. This also depends on two other properties of $-\nabla E_\varepsilon$ near the diagonal. First, the eigenvalues of $D\nabla E_\varepsilon(\Delta_{full})$ are so fortunate that by Sternberg's Linearization Theorem the flow $-\nabla E_\varepsilon$ is approaching $\Delta_{full}$ from all directions ``uniformly''. The second property is that $-\nabla E_\varepsilon$ is strictly outward-pointing from $\partial M_{K, \varepsilon}$ near $\Delta_{full}$. This property will be only computed on one example, but is expected to be true. See also the \textit{Mathematica} model in Figure \ref{figure_unstable_diagonal_intro}.

\begin{figure}[!htbp]
\labellist
\pinlabel $c_{p_\varepsilon}$ at 330 120
\pinlabel $u_{\varepsilon}^{trace}$ at 270 300
\pinlabel $\widetilde{u}_{\varepsilon}^{trace}$ at 360 200
\endlabellist
\centering
\includegraphics[scale=0.4]{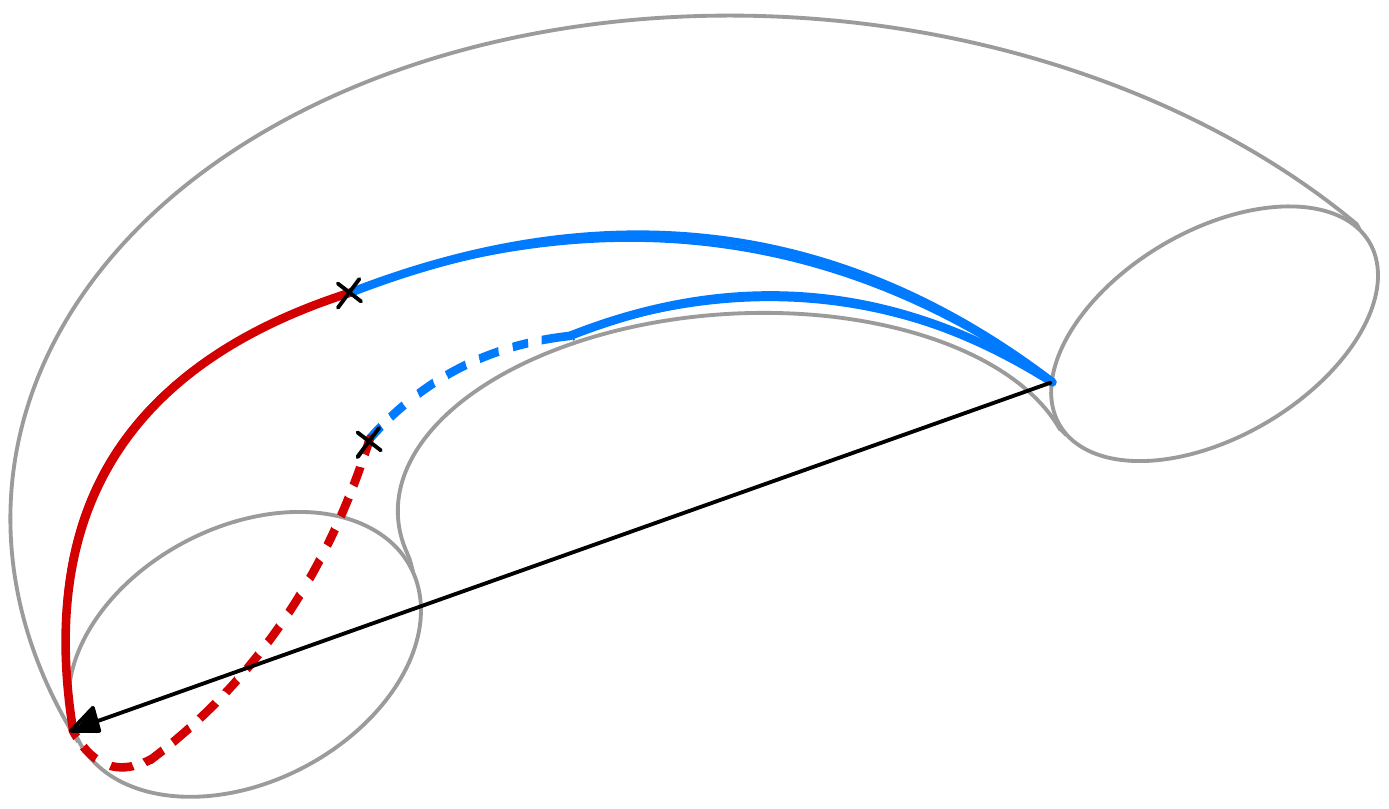}
\vspace{0.3cm}
\caption{Traces of $u_{\varepsilon}$ and $\widetilde{u}_{\varepsilon}$.}
\label{figure_two_trajectories_intro}
\end{figure}

In Figure \ref{figure_unstable_diagonal_intro}, we see the trace of a single trajectory $u_\varepsilon$ emanating from $p_\varepsilon$ that reaches $\Delta_{full}$. However, since the torus is symmetric along the $xy$-plane, there is also another trajectory $\widetilde{u}_\varepsilon$ with the trace symmetric to the trace of $u_\varepsilon$, see Figure \ref{figure_two_trajectories_intro}. Since the contribution of these two paths to the loop space coefficients differs by the meridian, we obtain the desired algebraic expression $1-\mu$.

The last topic of the thesis is about the map $\phi_{T_{K, \varepsilon}}$ from $LCH_0(\mathcal{L}^\ast_+ T_{K, \varepsilon})$ to $Cord^M(T_{K, \varepsilon}, \mathbb{Z}[\lambda^{\pm1}, \mu^{\pm1}])$ which will be only outlined. In more detail, we analyze the moduli space of $J$-holomorphic curves with boundary on the Lagrangian submanifold $L^\ast_+ T_{K, \varepsilon}\cup \R^3$. We stress the geometric interpretation of the compactness phenomena of $J$-holomorphic curves arising from the arboreal singularity. From that we propose a model of $0$-th string homology $H^{string}_0(T_{K, \varepsilon})$ which is expected to be a middle step between $LCH_0(\mathcal{L}^\ast_+ T_{K, \varepsilon})$ and $Cord^M(T_{K, \varepsilon}, \mathbb{Z}[\lambda^{\pm1}, \mu^{\pm1}])$. In short, our string homology $H^{string}_0(T_{K, \varepsilon})$ will be defined by a chain complex consisting of relative singular chains of broken strings.

In the end, we remark that we included in the Appendix also the computation of the interplay between the Morse and Maslov indices of the corresponding Reeb and binormal chords. Such a formula appeared independently in \cite{okamoto2024legendrian}. However, we believe that our approach focuses on a slightly different aspect of the problem and it is meaningful to include it.

\chapter{Knots and tori}
\label{ch:file6}
In the thesis, we would like to associate an algebraic structure with a thin torus $T_K\subset\R^3$ which arises as the boundary of a tubular neighborhood of the knot $K$. So the aim of this short chapter will be to inspect the map $K\mapsto T_K$ between the spaces of embedded knots and tori. For example, a priori it is not clear how many knots are ``bounded'' by a single torus or even $T_K$. For this, we use classical results from low-dimensional topology.

\begin{defn}Let $Q$ be an boundaryless submanifold of $N$. Let $U\subset \mathcal{N}(Q, N)$ be an open neighborhood of the zero section $0_{\mathcal{N}(Q, N)}$ of the normal bundle $\mathcal{N}(Q, N)$. Then a \textbf{tubular neighborhood $\nu Q$ of $Q$ in $N$} is a mapping $f:\overline{U}\rightarrow N$ satisfying the following:
\begin{itemize}
\item[$(i.)$] $f$ is an embedding.
\item[$(ii.)$] $f|_Q=Id_Q$ after the identification of $Q$ with the zero section $0_{\mathcal{N}(Q, N)}.$
\item[$(iii.)$] $f(U)$ is an open neighborhood of $Q$ in $N$.
\end{itemize}
\end{defn}

\begin{rem}\label{rem_tub} Let us consider a Riemannian metric $g$ on $N$, this induce a $\mathcal{N}_g(Q, N)$. Then a tubular neighborhood can be constructed as an exponential map that is precomposed with the radial shrinking along each fiber to a sufficiently small neighborhood of $0_{\mathcal{N}(Q, N)}$, for details see \cite[Thm 7.1.5]{mukherjee_2015}. Also, all tubular neighborhoods are (ambient) isotopic relative to the zero section \cite[Thm 5.3]{hirsch_1997}.
\end{rem} 

\begin{rem_not} Now we restrict ourselves to $N:=S^3$. By \cite{moise_1952} $S^3$ as a topological $3$-manifold admits unique smooth and piecewise-linear structures.

Now, a \textbf{\textit{knot}} $K$ is a smooth embedding $S^1\hookrightarrow S^3$.  Two knots $K_1$ and $K_2$ are \textbf{\textit{equivalent}} if there is an orientation--preserving homeomorphism $h:S^3\rightarrow S^3$ such that $h\circ K_1=K_2$.

Analogously, a \textbf{\textit{torus}} $T$ is a smooth embedding $S^1\times S^1\hookrightarrow S^3$. Two tori $T_1$ and $T_2$ are \textbf{\textit{equivalent}} if there is an orientation--preserving homeomorphism $h:S^3\rightarrow S^3$ such that $h\circ T_1=T_2$.
\end{rem_not}

\begin{rem}\label{rem_tub_neighb_invariant} By Remark \ref{rem_tub}, knots are completely determined by their tubular neighborhoods.

Next, let $K_1$ and $K_2$ be equivalent knots. Then, by restriction, there is an orientation-preserving homeomorphism between the \textbf{\textit{knot complements}} $\overline{S^3\setminus\nu K_1}$ and $\overline{S^3\setminus\nu K_2}$. In particular, it is an orientation-preserving homeomorphism between the boundaries of these complements. Thus, the boundary of a tubular neighborhood of the knot is a knot invariant.	
\end{rem}

\begin{defn} For any knot $K$ we define two simple curves $\mathfrak{m}, \mathfrak{L}\subset\partial(\nu K)$. The \textbf{meridian} $\bm{m}$ satisfies $\mathfrak{m}\not\sim 0$ in $\pi_1(\partial\nu K)$ and $\mathfrak{m}\sim 0$ in $\pi_{1}(\nu K)$. And the \textbf{longitude} $\mathfrak{L}$ represents a generator of $\pi_{1}(\nu K)\cong \mathbb{Z}$.
\end{defn}

The following lemma is an adapted version of the Solid torus theorem from \cite{rolfsen2003knots}.

\begin{lemma}\label{lemma_knot_complement} Let $T$ be a torus in $S^3$. Let $A$ be one connected component of $S^3\setminus T$. Then the closure $\overline{A}$ is a tubular neighborhood of some knot if and only if the natural inclusion $\pi_{1}(T)\rightarrow\pi_{1}(\overline{A})$ has a non-trivial kernel. Moreover, if $\overline{A}$ is a tubular neighborhood of some knot $K$, then $K$ is unique.
\end{lemma}
\begin{proof} Let the inclusion $\pi_{1}(T)\rightarrow\pi_{1}(\overline{A})$ has a non-trivial kernel. Then by Loop theorem \cite{papa_dehn_lemma} there is a simple non-contractible closed curve in $T$ that bounds a disk $D$ in $\overline{A}$. Next, note that $\overline{A}$ without a bicollared neighborhood $N$ of $D$ bounds (piecewise-linear) $S^2$ given by the union of two discs and the annulus $T\setminus N$. Hence, by the polyhedral Schoenflies theorem \cite{rolfsen2003knots}, $\overline{A}\setminus N$ is $D^3$. Finally, canonical gluing of $N$ to $D^3$ is homeomorphic to $S^1\times D^2$, and hence it is a tubular neighborhood of some knot. The opposite implication is trivial.

It remains to show the uniqueness of $K$. We need to inspect in how many ways we can embed $S^1\times D^2$ into $S^3$ such that the boundary is glued to the boundary of $\overline{S^3\setminus A}$ via an orientation--preserving homeomorphism. This is precisely the \textit{knot complement problem}, see \cite{Tietze1908berDT}, which was solved by Gordon and Luecke in \cite{Gordon1989KnotsAD}. They showed that there is only a trivial gluing (mapping the meridian to the meridian) or $\overline{A}$ is a tubular neighborhood of the unknot (in this case, we can map the meridian to the meridian plus some number of twists by the longitude). The lemma follows.
\end{proof}

\begin{lemma}\label{lem_torus_bounds}Any torus $T$ bounds exactly one tubular neighborhood of some knot, unless $T$ bounds a tubular neighborhood of an unknot, in which case it bounds two tubular neighborhoods of two unknots.
\end{lemma}
\begin{proof}

Let $A_1, A_2$ be connected components of $S^3\setminus T$. Also, $\pi_{1}(T)\cong \mathbb{Z}\oplus\mathbb{Z}$ and $\pi_{1}(S^3)\cong 0$. Hence by Van Kampen's theorem at least one of inclusions $\pi_1(T)\rightarrow\pi_1(\overline{A}_1)$ and $\pi_1(T)\rightarrow\pi_1(\overline{A}_2)$ has a nontrivial kernel. Then by Lemma \ref{lemma_knot_complement} the torus $T$ bounds at least one tubular neighborhood of some knot.

Let $T$ bounds $\nu K$. Then $\pi_{1}(T)\cong\pi_{1}(\partial\nu K)$ is generated by the meridian $\mathfrak{m}$ and the longitude $\mathfrak{L}$ of $K$. Also, by \cite{rolfsen2003knots} we know two observations: $\mathfrak{m}$ is not contractible in $\overline{S^3\setminus \nu K}$; and $\mathfrak{L}$ is contractible in $\overline{S^3\setminus \nu K}$ if and only if $K$ is an unknot. Then Lemma \ref{lemma_knot_complement} finishes the proof.
\end{proof}
\begin{cor}In $S^3$ there is a bijection between equivalence classes of knots and tori.
\end{cor}

\begin{rem} Let us fix an orientation of $S^3$. We would like to discuss the choices of \textbf{orientations} of knots and tori. 

Let $T\subset S^3$ which is not the boundary of a tubular neighborhood of an unknot. Then by Lemma \ref{lem_torus_bounds} there is a unique $K$ and a $\nu K$ such that $T=\partial(\nu K)$. The open subset $Int(\nu K)$ inherits the orientation of $S^3$. In particular, $\nu K$ is oriented, and hence $\partial(\nu K)$ is oriented too. This induces an orientation of $T$. However, such an orientation of $T$ does \textit{not} determine an orientation of $K$.

We note that an orientation of $\mathfrak{m}$ determines an orientation of $K$ (even if $K$ is an unknot). Indeed, by the definition of $\mathfrak{m}$ we can take a disk $D\subset\nu K$ such that $\partial D=\mathfrak{m}$ and $K\pitchfork D$ at the unique point $x$. Then
$$T_x S^3=T_x K\oplus T_x D.$$
Since $S^3$ is oriented and $D$ inherits the orientation from $\mathfrak{m}$, we obtain the orientation of $K$.
\end{rem}

\begin{defn} A \textbf{framed knot $K$} is a knot together with a nowhere-vanishing smooth section of the normal bundle $\mathcal{N}(K, S^3)$ called a \textbf{framing}. Moreover, we enhance $K$ with a base point $\ast$ on $K$.
\end{defn}

\begin{rem} Note that the choice of framing defines a (non-oriented) longitude $\mathfrak{L}\subset\partial(\nu K)$. An orientation of $K$ induces canonically an orientation of $\mathfrak{L}$. There is also a canonically chosen base point $\ast$ on $\mathfrak{L}$ corresponding to the base point on $K\subset \mathcal{N}(K, S^3)$.
\end{rem}

\begin{rem}Let $p\in S^3$. Then the stereographic projection gives us a diffeomorphism between $S^3\setminus\lbrace p\rbrace$ and $\R^3$.
\end{rem}

\noindent{\textit{From now on, we will consider knots and tori only in $\R^3$.}}

\begin{defn} Let $K$ be a knot in $(\R^3, g_{euc})$ and $\varepsilon>0$ small. Then we take a tubular neighborhood $\nu K$ given by an embedding $\overline{U_\varepsilon}\hookrightarrow\R^3$; $(q, v_q)\mapsto q+v_q$, where $U_\varepsilon:=\lbrace (q, v_q)\in \mathcal{N}_{g_{Euc}}(K, \R^3)\,|\,\,\Vert\,v_q\,\Vert_{g_{euc}}^2< \varepsilon\rbrace$.

Then the boundary $\partial (\nu K)$ will be denoted as $T_{K}.$
\end{defn}

\begin{rem}\label{rem_god_rad} There is a constant $\varepsilon_{good}>0$ such that $T_{K}$ exists for all $\varepsilon\in(0, \varepsilon_{good}]$. Moreover, they are all in the same (ambient) isotopy class, see Remark \ref{rem_tub}. Later, we would like often to stress out the size of the radius $\varepsilon$, and then we will write $T_{K, \varepsilon}$.
\end{rem}

\begin{cor} Let us consider a torus $T_K$.  Then the orientation classes of $K$ are in a bijection with the orientation classes of $\mathfrak{m}\subset T_K$. Moreover, the homotopy classes of framings are in bijection with homotopy classes of $\mathfrak{L}\subset T_K$.
\end{cor}

\chapter{Space of outward-pointing chords $M_{K, \varepsilon}$}

Given a thin torus $T_{K, \varepsilon}\subset (\R^3, g_{Euc})$, we consider the space of chords - oriented straight lines with endpoints on $T_{K, \varepsilon}$. In fact, we would like to consider a smaller space - $M_{K, \varepsilon}$ which consists of those chords that, at their endpoints, point out from the tubular neighborhood $\nu K$.

We emphasize that our primary motivation for studying this space is given by certain $J$-holomorphic curves. Such $J$-holomorphic curves have a boundary on the Lagrangian manifold $\R^3\cup L^\ast_+ T_{K, \varepsilon}$ with the arboreal singularity along $T_{K, \varepsilon}$ and at $+\infty$ are asymptotic to Reeb chords of $\R_+\times\mathcal{L}^\ast_+ T_{K, \varepsilon}$. It turns out that these $J$-holomorphic curves, when restricted to the zero section $\R^3$, behave as collections of paths with endpoints on $T_{K, \varepsilon}$, where they have \textit{outward-pointing conditions}. In other words, one can think nonformally about $J$-holomorphic curves as a potential chain map from Legendrian contact homology of $\mathcal{L}^\ast_+ T_{K, \varepsilon}$ to some Morse-theoretic algebra of chords on $T_{K, \varepsilon}$ (called cord algebra). And the $J$-holomorphic curves indicate that $M_{K, \varepsilon}$ is a good candidate for the ambient space of the cord algebra. For more details, we refer to Chapter \ref{ch:rel_symplectic}. On the other hand, the space $M_{K, \varepsilon}$ itself is interesting due to its nontrivial topology. For example, from its first homology we will be able to recollect some nontrivial information about the underlying knot.

The main goal of this chapter is to inspect the manifold structure of $M_{K, \varepsilon}$ with special care for the singular behavior near the constant chords, which will later play an important role in the cord algebra. We also study $E_\varepsilon$; the restriction of the standard energy function to $T_{K, \varepsilon}\times T_{K, \varepsilon}$.
\section{A bit of differential topology}
In this section, we present some rather standard properties of manifolds with corners that will be used throughout the subsequent sections.
\begin{defn} A \textbf{manifold with corners $M$} of dimension $n$ is a Hausdorff space such that:
\begin{itemize}
\item[$(i.)$] For every $x\in M$ there is a chart $(U_x, \varphi_x),$ i. e. an open neighbourhood $U_x\subset M$ and $\varphi_x:U_x\rightarrow \R^j\times[0, \infty)^{n-j}$ a homeomorphism onto an open subset of $\R^j\times[0, \infty)^{n-j}$ for some $j\in\mathbb{N}_0$.
\item[$(ii.)$] For any two charts $(U_x, \varphi_x)$ and $(U_y, \varphi_y)$ the composition $\varphi_x\circ\varphi_y^{-1}:\varphi_y(U_x\cap U_y)\rightarrow\varphi_x(U_x\cap U_y)$ extends to a smooth map between open sets of $\R^n$.
\end{itemize}
\end{defn}

\begin{defn_lemma}\cite{Nielsen1981TransversalityAT}Let $M$ be a manifold with corners and $x\in M$. We say that $j$ is the \textbf{index $\iota(M, x)$ of $x$ in $M$} if there is a chart $(U_x, \varphi_x)$ on $M$ such that $\varphi:U_x\rightarrow \R^j\times[0, \infty)^{n-j}$ and $\varphi_x(x)=0$. $(U_x, \varphi_x)$ will be called an \textbf{index chart for $x$}. The index is independent of the choice of an index chart.

$M$ is a disjoint union of (boundaryless) manifolds $\lbrace M_j\rbrace_{j\in\lbrace 0,\dots, n\rbrace}$ such that
$$M_j=\lbrace x\in M\,|\,\iota(M, x)=j\rbrace.$$
A connected component of $M_j$ is called an \textbf{$j$-stratum}, $M_n$ is the interior $Int(M)$ of $M$, and the set $M\setminus M_n$ is the boundary of $M$. The closure of each $M_j$ is called a \textbf{$j$-face}. We assume that each face is also a manifold with corners. 
\end{defn_lemma}


\begin{defn}Let $M$ be a manifold with corners, $N$ be a boundaryless manifold and $Q\subset N$ is a submanifold with corners. We say that a smooth map $f:M\rightarrow N$ is \textbf{stratum transverse to $Q$} if for any two strata $\mathcal{M}$ and $\mathcal{Q}$ the restriction $f|_{\mathcal{M}}:\mathcal{M}\rightarrow N$ is transverse to $\mathcal{Q}$. In more detail, for any $x\in \mathcal{M}$ holds that $f(x)\notin\mathcal{Q}$ or 
$$df(x)(T_x\mathcal{M})+T_{f(x)}\mathcal{Q}=T_{f(x)}N.$$
\end{defn}

\begin{thm}\cite{Nielsen1981TransversalityAT}\label{thm_invers_im} Let $M$ be a manifold with corners, $N$ be a boundaryless manifold and $Q\subset N$ be a submanifold with corners. Let $f:M\rightarrow N$ be a stratum transverse map to $Q$. Then either $f^{-1}(Q)=\emptyset$, or
\begin{itemize}
\item[$(i.)$]$f^{-1}(Q)$ is a submanifold with corners of $M$, and
\item[$(ii.)$]$f^{-1}(Q)$ is stratified by connected components of $(f|_{\mathcal{M}})^{-1}(\mathcal{Q})$ for all strata $\mathcal{M}$ and $\mathcal{Q}$ (of $M$ and $Q$, respectively), and
\item[$(iii.)$]$\dim M-\dim f^{-1}(Q)=\dim N-\dim Q$, and
\item[$(iv.)$]$\iota(M, x)-\iota(f^{-1}(Q), x)=\dim N-\iota(Q, f(x))$ for all $x\in f^{-1}(Q).$
\end{itemize}
\end{thm}

\begin{rem} The following theorem is a generalization of \cite[thm 1.4]{hirsch_1997} from manifolds with boundary to manifolds with corners.
\end{rem}

\begin{thm}\cite{hirsch_1997}\label{thm_open_embedd} Let $k\geq1$. Let $M, N$ be manifolds with corners. Then the set of $C^k$ embeddings from $M$ to $N$ is open in $C^k(M, N)$ (with strong $C^k$ topology).
\end{thm}

\begin{thm}\cite{Trotman1978StabilityOT, goresky2012stratified}\label{thm_open_tran}Let $k\geq1$. Let $M$ be a manifold with corners, $N$ be a boundaryless manifold, and $Q\subset N$ be a submanifold with corners which is \emph{closed} in $N$. Then the set
$$\lbrace f\in C^{k}(M, N)\,|\,f\hbox{ is stratum transverse to }Q\rbrace$$
is open (in the strong $C^{k}$ topology) in $C^{k}(M, N)$.
\end{thm}

\begin{n_example}\cite{guillemin2010differential} Typically, we would like to apply Theorem \ref{thm_open_tran} in the following scenario. We have a transverse map, and we would like to extend this map to a smooth homotopy of transverse maps. However, we will see that this might fail if the domain of the map is not compact. 

For this, we choose an auxiliary smooth function $\rho:\R\rightarrow\R$ such that 
$$\rho(t)=
\begin{cases}
1, &\text{ for }|t|<1,\\
0, &\text{ for }|t|>2.\\
\end{cases}
$$
Now we define a family of functions $\lbrace f_t:\R\rightarrow\R\rbrace_{t\in\R}$ by $$f_t(x)=x\rho(tx).$$
Let us put $Q=\lbrace 0\rbrace$. Then $f_0\pitchfork Q$, but for any $|t|>0$ it holds that $f_t\not\pitchfork Q$.
\end{n_example}

\begin{rem}The following theorem is a generalization of Ehresmann's Theorem \cite{2020IntroductionTL} for manifolds with boundary to manifolds with corners. See also Thom's first isotopy Lemma \ref{thm_thom_isot}, which can be seen as Ehresmann's Theorem for Whitney stratified spaces.
\end{rem}

\begin{Ehresmann}\cite{2020IntroductionTL}\label{thm_ehrsm} Let $M$ be a manifold with corners, $N$ be a boundaryless manifold. Let $\pi:M\rightarrow N$ be a proper surjection that is a submersion on each stratum of $M$. Then $\pi$ defines a locally trivial \textbf{fiber bundle}, i.e., for every $x\in N$ there is an open neighbourhood $U_x$ in $N$ and a (strata preserving) diffeomorphism $\psi: U_x\times\pi^{-1}(x)\rightarrow\pi^{-1}(U_x)$ such that $\pi\circ\psi(u, v)=u$ for every $(u, v)\in U_x\times\pi^{-1}(x)$.
\end{Ehresmann}

\begin{rem}Similarly to vector bundles, as discussed in Appendix \ref{sect_vect}, also any locally trivial fiber bundle $M\rightarrow N$ over a contractible space $N$ is trivial.

In short. We pick a smooth homotopy $F:N\times [0, 1]\rightarrow N$ from the identity to a constant map. Then the pull-back bundle $F^\ast M\rightarrow N\times [0, 1]$ possesses a complete Ehresmann's connection \cite{Hoyo2015CompleteCO} (\cite{Hoyo2015CompleteCO} also generalizes to the case for manifolds with corners). Then the time $1$ flow of the horizontal lift of $(0, \partial_t)$ gives a global trivialization of the fiber bundle $M$.
\end{rem}

\begin{stability}\label{lem_stability}Let $M$ be a compact manifold with corners, $N$ be a boundaryless manifold and $Q\subset N$ be a submanifold with corners which is closed in $N$. Assume that there are two maps $f, g:M\rightarrow N$ such that $f$ is stratum transverse to $Q$ and $f(M)\cap Q\neq\emptyset$. If $g$ is sufficiently close to $f$ in $C^\infty$ topology, then there is a smooth homotopy $F:M\times[0, 1]\rightarrow N$ such that
\begin{itemize}
\item[$(i.)$]$F|_{M\times\lbrace 0\rbrace}=f$ and $F|_{M\times\lbrace 1\rbrace}=g$.
\item[$(i.)$]For any $t\in[0, 1]$ it holds that $F|_{M\times\lbrace t\rbrace}$ is stratum transverse to $Q$.
\end{itemize}
Moreover, $f^{-1}(Q)$ and $g^{-1}(Q)$ are isotopic submanifolds with corners of $M$.
\end{stability}
\begin{proof}
Let $U$ be a small open neighbourhood of $f$ (in strong $C^\infty$ topology) in $C^\infty(M, N)$, which consists of stratum transverse maps to $Q$. Such a neighbourhood exists by Theorem \ref{thm_open_tran}.

Next, by the Strong Whitney embedding theorem \cite{hirsch_1997}, there is an embedding $\varphi:N\rightarrow\R^k$ for some $k$. Let us consider the normal bundle $\mathcal{N}(N, \R^k)$ with the projection $\pi_{\mathcal{N}}$. Put $I_\varepsilon:=(-\varepsilon, 1+\varepsilon)$ for some $\varepsilon>0$. Now we define a smooth homotopy $\widetilde{F}:M\times I_\varepsilon\rightarrow N$ by 
$$\widetilde{F}|_{M\times\lbrace t\rbrace}=(\varphi^{-1}\circ\pi_{\mathcal{N}})\big(t\varphi(g)+(1-t)\varphi(f)\big).$$
Since $M$ is compact, we see that if $\varepsilon>0$ is sufficiently small and $g$ is sufficiently close to $f$, then each $\widetilde{F}|_{M\times\lbrace t\rbrace}$ lies also in $U$. In particular, $\widetilde{F}$ is also stratum transverse to $Q$.  

Now, by Theorem \ref{thm_invers_im} $(i.)$, $\widetilde{F}^{-1}(Q)$ is a manifold with corners. We would like to show that $\pi_t:\widetilde{F}^{-1}(Q)\rightarrow I_\varepsilon$ defines a locally trivial fiber bundle. For this, we would like to apply Ehresmann's Theorem \ref{thm_ehrsm}. First, $\pi_t$ is smooth. Since $f$ has a non-empty intersection with some stratum of $Q$, the Simplified Newton's method \cite{webber99} on $F$ gives us the surjectivity of $\pi_t$ (provided that $\varepsilon>0$ is sufficiently small and $g$ is sufficiently close to $f$). And since $M$ is compact and $Q$ closed, the fibers are compact and $\pi_t$ is proper. 

It remains to show that $\pi_t$ is a submersion on each stratum of $\widetilde{F}^{-1}(Q)$. By Theorem \ref{thm_invers_im} $(ii.)$ $\widetilde{F}^{-1}(Q)$ is stratified by connected components of $\widetilde{F}|_{\mathcal{M}\times I_\varepsilon}^{-1}(\mathcal{Q})$ for all strata $\mathcal{M}\subset M$ and $\mathcal{Q}\subset Q$.

Let $x=(m, t_0)\in\widetilde{F}|_{\mathcal{M}\times I_\varepsilon}^{-1}(\mathcal{Q})$. We would like to show that there is a tangent vector $V=(v, \partial_t)\in T_x(\widetilde{F}|_{\mathcal{M}\times I_\varepsilon}^{-1}(\mathcal{Q}))$. It holds that
$$T_x(\widetilde{F}|_{\mathcal{M}\times I_\varepsilon}^{-1}(\mathcal{Q}))=\big\lbrace (w_1, w_2)\in T_x(\mathcal{M}\times I_\varepsilon)\,|\,d\widetilde{F}|_{\mathcal{M}\times I_\varepsilon}(x)(w_1, w_2)\in T_{\widetilde{F}(x)}\mathcal{Q}\big\rbrace.$$
Moreover, by the construction of $\widetilde{F}$ it holds that
\begin{equation}\label{eqn_split_lin}
d\widetilde{F}|_{\mathcal{M}\times I_\varepsilon}(x)\big(T_m\mathcal{M}, 0\big)+T_{\widetilde{F}(x)}\mathcal{Q}=T_{\widetilde{F}(x)}N.
\end{equation}
Let $\lbrace\partial_{z_i}, \partial_{z_j^\bot}\rbrace$ be a basis of $T_{\widetilde{F}(x)}N$ such that $\lbrace \partial_{z_i}\rbrace$ is a basis of $T_{\widetilde{F}(x)}\mathcal{Q}.$ Then 
$$d\widetilde{F}|_{\mathcal{M}\times I_\varepsilon}(x)(0, \partial_t)=a^i\partial_{z_i}+b^j\partial_{z_j^\bot},$$
for some $a^i, b^i\in \R$. By relation $(\ref{eqn_split_lin})$ we see that there exists $v\in T_m\mathcal{M}$ such that 
$$d\widetilde{F}|_{\mathcal{M}\times I_\varepsilon}(x)(v, 0)=-b^j\partial_{z_j^\bot}.$$
This gives us the desired vector $V$, and consequently $\pi_t$ is also a submersion on each stratum of $\widetilde{F}^{-1}(Q)$.

Then Ehresmann's Theorem \ref{thm_ehrsm} gives us a locally trivial fiber bundle. The bundle has a contractible base, so in particular it is a trivial bundle. Since each $\widetilde{F}|_{M\times\lbrace t\rbrace}^{-1}(Q)$ is a submanifold of $M$, $f^{-1}(Q)$ and $g^{-1}(Q)$ are isotopic submanifolds of $M$.
\end{proof}

\begin{defn}\label{def_outward_point} Let $M$ be an $n$-dimensional manifold with corners, $x\in M$ with $\iota(M, x):=j<n$ and $v\in T_x M$. Let $(U_x, \varphi_x: U_x\rightarrow\R^j\times [0, \infty)^{n-j})$ be an index chart for $x$. We consider coordinates $(q_1,\dots, q_n)$ on $\R^n\times\R^{n-j}$. Then we say that the vector $v$ is \textbf{outward-pointing at $x$} if
\begin{equation}\label{eqn_out_cond}
dq_i\big[(\varphi_x)_\ast v\big](0)\leq0
\end{equation}
for every $i\in\lbrace j+1,\dots, n\rbrace$. If the signs in $(\ref{eqn_out_cond})$ are all $\geq$, we say that the vector $v$ is \textbf{inward-pointing at $x$}. If the signs are all sharp, we say that $v$ is \textbf{strictly outward/inward-pointing at $x$}. 

The definitions are independent on the choice of the index chart. Moreover, if $v=0$, then $v$ is both outward-pointing and inward-pointing.
\end{defn}

\section{Definition of $M_{K, \varepsilon}$}
We are going to introduce $M_{K, \varepsilon}$; the space of outward-pointing chords,  together with some associated terms.

\begin{Thom}[\cite{mukherjee_2015}]\label{thm_thom}
Let $r\in\mathbb{N}_0$. Suppose $M, N, Q$ are boundaryless manifolds, where $Q$ is a submanifold of $J^r(M, N)$. Then the set $\lbrace f\in C^\infty(M, N)\,|\,j^r(f)\pitchfork Q\rbrace$ is a dense subset of $C^\infty(M, N)$ (with strong topology) and open if $Q$ is closed.
\end{Thom}

\begin{rem}\label{rem_stratified}A closed subset $\Sigma$ of a manifold $N$ is called \textbf{stratified}, if it is a disjoint locally finite union of connected submanifolds of $M$ called \textbf{strata}. The dimension of $\Sigma$ is the maximal dimension of its strata.

The stratified subset $\Sigma$ is called \textbf{Whitney}, if the following holds. Let $\mathcal{Q}, \mathcal{M}$ are two strata with $\dim(\mathcal{Q})<\dim(\mathcal{M})$. Suppose that there are two sequences of points $x_i\in \mathcal{M}$ and $y_i\in \mathcal{Q}$ that converge to some $y\in\mathcal{Q}$. Let the secant lines $\overline{x_i, y_i}$ converge to some limiting line $\ell$ and the planes $T_{x_i}\mathcal{M}$ converge to some limiting plane $\tau$. Then the \textbf{Whitney conditions} holds, that is $T_y\mathcal{Q}\subset\tau$ and $\ell\subset\tau$. The limiting plane $\tau$ is called a \textbf{generalized tangent space at $y$}.

For example, every real analytic set can be (Whitney) stratified \cite{goresky2012stratified}. 

A map $f:M\rightarrow N$ is transverse to a (Whitney) stratified subset $\Sigma$ if it is transverse to each stratum. In addition, Theorem \ref{thm_thom} holds also in the case $Q=\Sigma$, see \cite[Thm 2.3.2]{Eliashberg2002IntroductionTT} and \cite{Trotman1978StabilityOT}. And if $f\pitchfork\Sigma$, then $f^{-1}(\Sigma)$ is a canonically (Whitney) stratified subspace of $M$ of the same codimension as $\Sigma$, see \cite[Thm 1.4]{Gibson1976TopologicalSO}. Also, the transverse intersection of two (Whitney) stratified subspaces is a (Whitney) stratified subspace \cite{Gibson1976TopologicalSO}.
\end{rem}

\begin{lemma}\label{lemma_curv}A generic knot in $\R^3$ has nowhere-vanishing curvature.
\end{lemma}
\begin{proof}Let $\gamma:S^1\rightarrow\R^3$ be a knot. If $s\in S^1$, then the curvature $\kappa$ at the point $s$ is defined as $\kappa(s):=||\dot{\gamma}(s)\times\ddot{\gamma}(s)||/||\dot{\gamma}(s)||^3$.

Let us consider coordinates $(s, a, b, c)$ on $J^2(S^1, \R^3)\cong S^1\times\R^3\times\R^3\times\R^3$. If $s\in S^1$, then $j^2\gamma(s)=(s, \gamma(s), \dot{\gamma}(s), \ddot{\gamma}(s))\in J^2(S^1, \R^3)$.

Let us put $\Sigma:=\Sigma_1\cup\Sigma_2$, where
$$\Sigma_1=\big\lbrace(b, c)\in\R^3\setminus\lbrace 0\rbrace\times\R^3\,|\,b\times c= 0\big\rbrace\hbox{ and }\Sigma_2=\lbrace 0\rbrace\times\R^3.$$
We are going to show that $\Sigma$ is a $\codim 2$ stratified subset of $\R^6$.

From $b\times c=0$ we see that $b$ and $c$ are linearly dependent. So $\Sigma_1$ has a structure of real line bundle over $\R^3\setminus\lbrace 0\rbrace$, which is a $4$-dimensional submanifold of $\R^6$. Moreover, we see that $\Sigma$ is closed. Alternatively, we can see $\Sigma$ as an algebraic set $b\times c=0$. The map $b\times c$ has rank $2$ on $\Sigma\setminus\lbrace(0, 0)\rbrace$, and thus $(0, 0)$ is the only singular point of $\Sigma$. Then the submanifold structure of $\Sigma\setminus\lbrace(0, 0)\rbrace$ will follow from \cite[Prop 3.3.10-11]{Bochnak1992RealAG}.

So $\widetilde{\Sigma}:=S^1\times\R^3\times\Sigma$ is a codimension $2$ stratified subset of $J^2(S^1, \R^3)$. By Thom's Transversality Theorem \ref{thm_thom} we have that $j^2\gamma\pitchfork \widetilde{\Sigma}$ for a generic knot. Hence, from the dimension argument, the lemma follows.
\end{proof}

\begin{lemma}\label{lemma_tangent_gen}For a generic knot $K$ in $\mathbb{R}^3$ it holds
\begin{itemize}
\item[$(i.)$]Let $(p, q)\in K\times K$ such that $p\neq q$. The set of all pairs (points) $(p, q)$ such that the tangent line to the knot from $p$ meets the knot again at $q$ is finite. Any such pair $(p, q)$ is called \textbf{$p$-special}. 
\item[$(ii.)$]Any $p$-special $(p, q)$ pair has the property that tangent line at $q$ does not lie in the osculating plane of $p$. Recall that the osculating plane of $p$ is the affine plane spanned by the tangent and normal vectors to $K$ at $p$.
\item[$(iii.)$]The subset of those pairs, such that $p$ and $q$ have the same tangent line, is empty.
\end{itemize}
The analogous statement holds also for \textbf{$q$-special pairs (points) $(p, q)$.}
\end{lemma}

\begin{proof}For $(i.)$ and $(iii.)$, see the proof of Lemma $7.10 (b)$ in \cite{Cieliebak2016KnotCH}. The part $(ii.)$ is a straightforward application of Thom's Transversality Theorem \ref{thm_thom_isot}.
\end{proof}

\begin{rem} From now on, we will consider only knots in $\R^3$ that satisfy Lemmata \ref{lemma_tangent_gen} and \ref{lemma_curv}.
\end{rem}

\begin{rem_not}\label{rem_tor_param}From now on, until otherwise said, $K$ will be a knot in $\mathbb{R}^3$ with an arc length parametrization $\gamma:\mathbb{R}/T\mathbb{Z}\rightarrow\mathbb{R}^3$. 

Now, we introduce a function
\begin{align*}
v:\mathbb{R}/T\mathbb{Z}\times S^1&\rightarrow\mathbb{R}^3\\
(s, \theta)&\mapsto \cos(\theta) n(s)+\sin(\theta) b(s),
\end{align*}
where $n(s)=\ddot{\gamma}(s)/||\ddot{\gamma}(s)||$ is the normal vector at $\gamma(s)$ and $b(s)=\dot{\gamma}(s)\times n(s)$ is the binormal vector at $\gamma(s)$. Later, we will also use the notion of the torsion $\tau(s):=\langle\dot{\gamma}(s)\times\ddot{\gamma}(s), \dddot{\gamma}(s)\rangle/||\dot{\gamma}(s)\times\ddot{\gamma}(s)||^2$.

Let $\varepsilon\in(0, \varepsilon_{good}]$. Recall that $\varepsilon_{good}$ was defined in Remark \ref{rem_god_rad} as a radius of a tubular neighbourhood of $K$. Then the embedded torus $T_{K, \varepsilon}$ can be naturally parametrized as
\begin{align*}
\Gamma_\varepsilon:\mathbb{R}/T\mathbb{Z}\times S^1&\rightarrow\mathbb{R}^3\\
(s, \theta)&\mapsto \gamma(s)+\varepsilon v(s, \theta).
\end{align*} 
Note that the map $\Gamma_{\varepsilon=0}$ is also well-defined and its image coincides with $K$.

Let us recall the Frenét equations
$$\ddot{\gamma}=\kappa n,\quad \dot{n}=-\kappa \dot{\gamma}+\tau b,\quad\dot{b}=-\tau n,$$
which we use for the auxiliary computations 
\begin{align}
\begin{split}\label{eqn_torus_aux_deriv}
    \partial_{s}v(s, \theta)&=-\cos(\theta)\kappa(s)\dot{\gamma}(s)+\cos(\theta)\tau(s)b(s)-\sin(\theta)\tau(s)n(s),\\
    \partial_{\theta}v(s, \theta)&=-\sin(\theta)n(s)+\cos(\theta)b(s),\\
        \partial_{s}\Gamma_\varepsilon(s, \theta)&=\big(1-\varepsilon\cos(\theta)\kappa(s)\big)\dot{\gamma}(s)+\varepsilon\cos(\theta)\tau(s)b(s)-\varepsilon\sin(\theta)\tau(s)n(s),\\
    \partial_{\theta}\Gamma_\varepsilon(s, \theta)&=-\varepsilon\sin(\theta)n(s)+\varepsilon\cos(\theta)b(s).
\end{split}
\end{align}
\end{rem_not}

\begin{rem}\label{rem_eps_bound}
Since $T_{K, \varepsilon}$ is an embedded surface, the rank of $d\Gamma_\varepsilon$ is $2$ at each point $(s, \theta)$. Hence, by relations (\ref{eqn_torus_aux_deriv}), we obtain that $1-\varepsilon\cos(\theta)\kappa(s)>0$. This give us an upper bound on $\varepsilon_{good}:$
$$\varepsilon_{good}<\min\left\lbrace\left\lvert\frac{1}{\kappa(s)} \right\rvert\,|\,s\in\R/T\mathbb{Z}\right\rbrace.$$
\end{rem}

\begin{rem}\label{rem_normal_vect} We are going to show that $v(s, \theta)$ is outward-pointing normal vector to $T_{K, \varepsilon}$ at $\Gamma(s, \theta)$, as expected. For this, we compute
\begin{align*}
\partial_{s}\Gamma_\varepsilon(s, \theta)\times\partial_{\theta}\Gamma_\varepsilon(s, \theta)=&\varepsilon\big(1-\varepsilon\cos(\theta)\kappa(s)\big)\dot{\gamma}(s)\times\partial_{\theta}v(s, \theta)+\varepsilon^2\tau(s)\partial_{\theta}v(s, \theta)\times\partial_{\theta}v(s, \theta)\\
=&\varepsilon\big(1-\varepsilon\cos(\theta)\kappa(s)\big)\dot{\gamma}(s)\times\partial_{\theta}v(s, \theta)\\
=&-\varepsilon\big(1-\varepsilon\cos(\theta)\kappa(s)\big)v(s, \theta).
\end{align*}
Note that the vector product is nonzero by the definition of $\varepsilon$ from Remark \ref{rem_tor_param} and the bound on $\varepsilon_{good}$ from Remark \ref{rem_eps_bound}. So $v(s, \theta)$ is normal to $T_{K, \varepsilon}$. Outward-pointing property is immediate from the definition of $\Gamma_\varepsilon$.
\end{rem}

\begin{defn}\label{defn_funct_out}
For $i=1, 2$, we define two families of functions $\lbrace F^{[\varepsilon]}_i\rbrace_{\varepsilon\in[0, \varepsilon_{good}]}$ which are given by
\begin{align*}
    F_i^{[\varepsilon]}:(\mathbb{R}/T\mathbb{Z}\times S^1)^2&\rightarrow\mathbb{R}\\
    (s_1, \theta_1, s_2, \theta_2)&\mapsto\langle\Gamma_\varepsilon(s_2, \theta_2)-\Gamma_\varepsilon(s_1, \theta_1), v(s_i, \theta_i)\rangle.
\end{align*}
Here $\langle\cdot,\cdot\rangle$ is the standard Euclidean metric on $\R^3$. See also Figure \ref{figure_torus_chords}.
\begin{figure}[!htbp]
\labellist
\pinlabel ${\Gamma_\varepsilon(s_2, \theta_2)-\Gamma_\varepsilon(s_1, \theta_1)}$ at 250 225
\pinlabel ${v(s_1, \theta_1)}$ at 336 183
\pinlabel ${v(s_2, \theta_2)}$ at 395 310
\endlabellist
\centering
\includegraphics[scale=0.95]{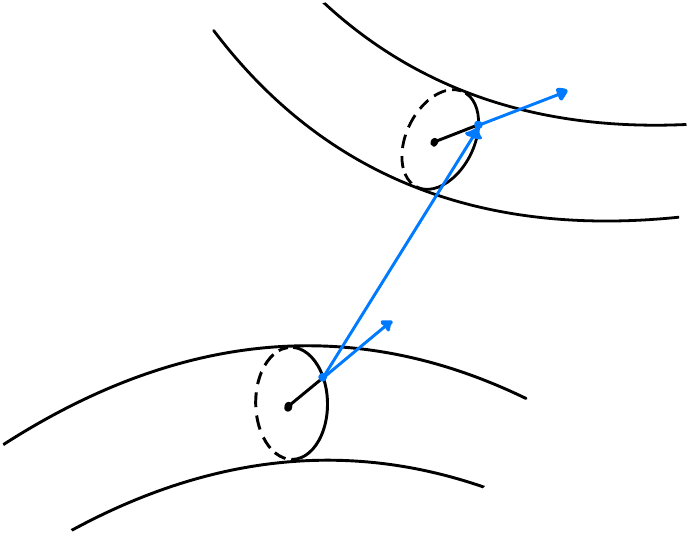}
\vspace{0.3cm}
\caption{The vectors $\Gamma_\varepsilon(s_2, \theta_2)-\Gamma_\varepsilon(s_1, \theta_1), v(s_1, \theta_1)$ and $v(s_2, \theta_2)$ appearing in the functions $F_1^{[\varepsilon]}$ and $F_2^{[\varepsilon]}$ in the case $\varepsilon>0$.}
\label{figure_torus_chords}
\end{figure}
\end{defn}

\begin{defn} \label{defn_out_out}Let $\varepsilon\in[0, \varepsilon_{good}]$. The set $$M_{K, \varepsilon}:=\overline{\left\{
(s_1, \theta_1, s_2, \theta_2)\in(\mathbb{R}/T\mathbb{Z}\times S^1)^2 \Bigm| \begin{array}{l}
F_1(s_1, \theta_1, s_2, \theta_2)\geq 0,\\
F_2(s_1, \theta_1, s_2, \theta_2)\geq 0,\\
\hspace{1.2cm} s_1\neq s_2
\end{array}
\right\}}.$$ is called the \textbf{set of outward-pointing chords}.
 
Moreover, if $\varepsilon>0$, then the set $M_{K, \varepsilon}|_{s_1=s_2}$ is called \textbf{$\varepsilon$-diagonal} and denoted by $\Delta_\varepsilon$. The set $\lbrace(\mathbb{R}/T\mathbb{Z}\times S^1)^2\,|\,s_1=s_2\wedge \theta_1=\theta_2\rbrace$ is called the \textbf{full diagonal} and denoted by $\Delta_{full}$.
\end{defn}

\begin{rem_not}\label{notation_on_circle} Let $r\in\R_+$. An identification of $S^1$ with $\R/r\mathbb{Z}$ gives $S^1$ a structure of an abelian group. Let $a, b, c\in S^1$. 

Under the expression $a\pm b$ we will understand the unique equivalence class $[a\pm b]$.

The linear order on $\R$ induce on $S^1$ a ternary relation - \textbf{cyclic order} $[a, b, c]$, which roughly means ``after $a$, one reaches $b$ before $c$'', for details see \cite{nlab:cyclic_order}. This allows us to define intervals on $S^1$. For example, the open interval $(a, c)$ is defined as the set of all $b$ such that $[a, b, c]$.

Next, $\widetilde{\langle\cdot,\cdot\rangle}$ will denote the flat metric on $S^1$, which is defined as the pull-back of the Euclidean metric on $\R$. This will induce the distance function $\widetilde{d}(\cdot,\cdot)$ on $S^1$.
\end{rem_not}

\begin{defn} Let $T_{K, \epsilon}$ be a torus in $\R^3$, then for each $\overline{s}\in\R/T\mathbb{Z}$ we define $\mathcal{C}_{\overline{s}, \varepsilon}$ as follows. It will be the infinite cylinder in $\R^3$ that is spanned by the lines with directional vector $\dot{\gamma}(\overline{s})$, which are intersecting the circle $\Gamma_\varepsilon|_{s=\overline{s}}$.
\end{defn}

\begin{defn_lemma}\label{lemma_aux_diag} Let $K$ be a knot in $\R^3$. Then for each $\overline{s}\in\R/T\mathbb{Z}$ there is a $\delta_{\overline{s}}\in(0, T/2)$ such that for each $\delta_0\in(0, \delta_{\overline{s}}]$ there is $\varepsilon_{0, \overline{s}}\in(0, \delta_{\overline{s}})$ with the following property. If $\varepsilon\in(0, \varepsilon_{0, \overline{s}}]$, then
\begin{itemize}
\item[$(i.)$]There exist the unique $a\in(\overline{s}-\delta_0, \overline{s})$ and $b\in(\overline{s}, \overline{s}+\delta_0)$ such that the circles $\Gamma_\varepsilon|_{s=a}, \Gamma_\varepsilon|_{s=b}$ are tangent to the infinite cylinder $\mathcal{C}_{\overline{s}, \varepsilon}$ at one point.

Moreover, it holds that $a\in(\overline{s}-\delta_0, \overline{s}-\varepsilon)$ and $b\in(\overline{s}+\varepsilon, \overline{s}+\delta_0)$, and for any $s\in[\overline{s}-\delta_0, a)\cup(b, \overline{s}+\delta_k]$ the circle $\Gamma|_{s}$ has an empty intersection with $\mathcal{C}_{\overline{s}, \varepsilon}$.

See Figure \ref{figure_numeric_diagonal} (top).
\item[$(ii.)$]There exist the unique $c\in(\overline{s}-\delta_0, \overline{s})$ and $d\in(\overline{s}, \overline{s}+\delta_0)$ such that the infinite cylinders $\mathcal{C}_{c, \varepsilon}, \mathcal{C}_{d, \varepsilon}$ are tangent to the circle $\Gamma_\varepsilon|_{s=\overline{s}}$ at one point.

Moreover, it holds that $c\in(\overline{s}-\delta_0, \overline{s}-\varepsilon)$ and $d\in(\overline{s}+\varepsilon, \overline{s}+\delta_0)$, and for any $s\in[\overline{s}-\delta_0, a)\cup(b, \overline{s}+\delta_k]$ the circle $\Gamma|_{s=\overline{s}}$ has an empty intersection with $\mathcal{C}_{s, \varepsilon}$.

See Figure \ref{figure_numeric_diagonal} (bottom).
\end{itemize}

Now we define the following terms.

Let $\delta\in(0, \inf\lbrace\delta_{\overline{s}}\,|\,\overline{s}\in \R/T\mathbb{Z}\rbrace)$. Let us consider a set 
$$\lbrace(s_1, s_2)\in\R/T\mathbb{Z}\times\R/T\mathbb{Z}\,|\,\widetilde{d}(s_1, s_2)\leq\delta\rbrace$$ 
or equivalently
$$F_{dist}^{-1}([-\delta, \delta]),$$
where the smooth function $F_{dist}:\R/T\mathbb{Z}\times\R/T\mathbb{Z}\rightarrow \R/T\mathbb{Z}$ is given by $(s_1, s_2)\mapsto s_2-s_1$. 
Such a set will be called a \textbf{weak diagonal set}, and elements of the weak diagonal set are called \textbf{weak diagonal pairs (points)}.

Put $\varepsilon_{diag}=\inf\lbrace \varepsilon_{0, \overline{s}}\,|\,\overline{s}\in\R/T\mathbb{Z}\wedge \delta_0=\delta]\rbrace$. Let $\varepsilon\in(0, \varepsilon_{diag})$. The subset of the weak diagonal set
$$\lbrace(s_1, s_2)\in\R/T\mathbb{Z}\times\R/T\mathbb{Z}\,|\,\widetilde{d}(s_1, s_2)<\varepsilon\rbrace$$
is called an \textbf{almost diagonal set} and elements of the almost diagonal set are called \textbf{almost diagonal pairs}. Finally, pairs $(s_1, s_2)$ with $s_1=s_2$ are called \textbf{$0$-diagonal} and form the set denoted by $\Delta_0$.
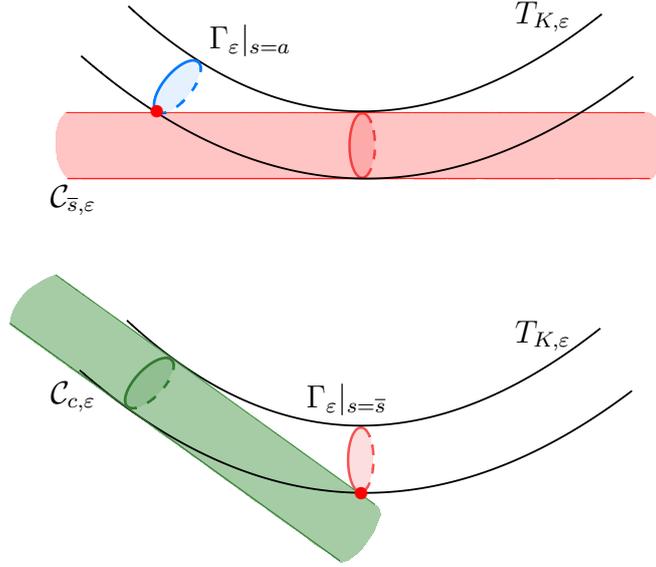
\begin{figure}[!htbp]
\labellist
\pinlabel $\textcolor{red}{\bullet}$ at 219 531
\pinlabel $\textcolor{red}{\bullet}$ at 92 770
\pinlabel $\mathcal{C}_{\overline{s}, \varepsilon}$ at 40 715
\pinlabel $\Gamma_\varepsilon|_{s=a}$ at 150 815
\pinlabel $T_{K, \varepsilon}$ at 332 830
\pinlabel $\Gamma_\varepsilon|_{s=\overline{s}}$ at 210 590
\pinlabel $\mathcal{C}_{c, \varepsilon}$ at 40 592
\pinlabel $T_{K, \varepsilon}$ at 332 630
\endlabellist
\centering
\includegraphics[scale=0.6]{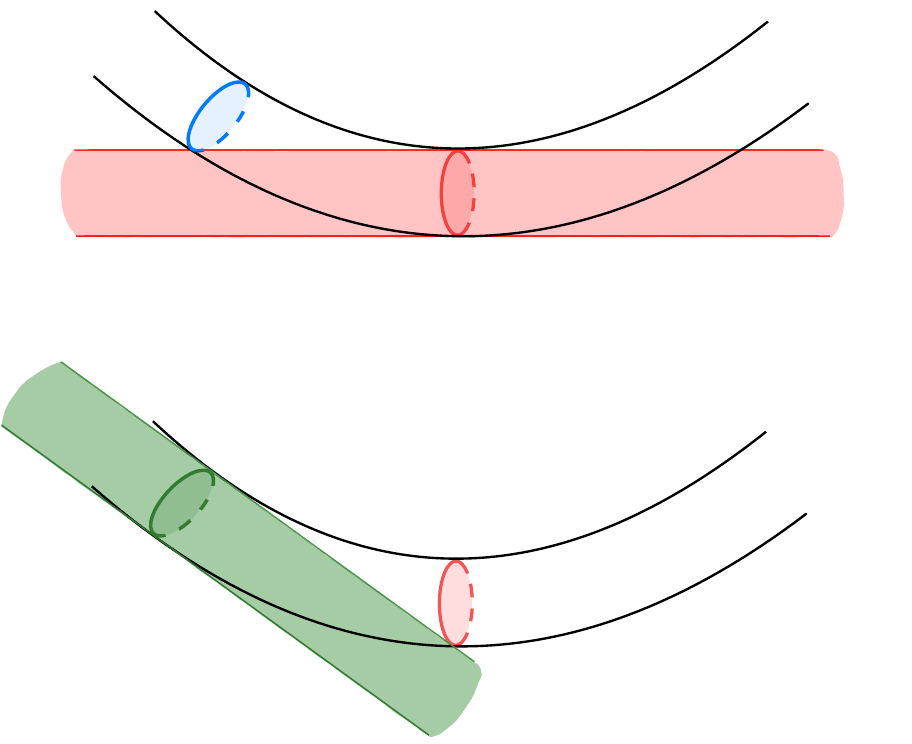}
\vspace{0.3cm}
\caption{\textit{Top:} the intersection (red point) of the cylinder $\mathcal{C}_{\overline{s}, \varepsilon}$ with the circle $\Gamma_\varepsilon|_{s=a}$. \textit{Bottom:} the intersection (red point) of the cylinder $\mathcal{C}_{c, \varepsilon}$ with the circle $\Gamma_\varepsilon|_{s=\overline{s}}$.}
\label{figure_numeric_diagonal}
\end{figure}
\end{defn_lemma}

\begin{proof} Even though the lemma is ``evident'' from Figure \ref{figure_numeric_diagonal}, the actual proof is a bit technical, since one has to consider all possible local perturbations of the underlying knot $K$. Thus, we leave the straightforward, but technical, proof to Appendix \ref{apx:proof_lemma}.
\end{proof}

\begin{defn}\label{lemma_special_delta} Let $K$ be a knot in $\R^3$ and $\overline{s}_2$-special pair $(\overline{s}_1, \overline{s}_2)\in\R/T\mathbb{Z}\times\R/T\mathbb{Z}$. Let also $\delta_{(\overline{s}_1, \overline{s}_2)}>0$. Then the set
$$W_{(\overline{s}_1, \overline{s}_2)}=\lbrace (s_1, s_2)\in[\overline{s}_1-\delta_{(\overline{s}_1, \overline{s}_2)}, \overline{s}_1+\delta_{(\overline{s}_1, \overline{s}_2)}]\times[\overline{s}_2-\delta_{(\overline{s}_1, \overline{s}_2)}, \overline{s}_2+\delta_{(\overline{s}_1, \overline{s}_2)}]\rbrace$$
is called the \textbf{weakly $\overline{s}_2$-special set} and its elements as the \textbf{weakly $\overline{s}_2$-special pairs (points)}.

The analogous statement also holds for $\overline{s}_1$-special pairs.
\end{defn}

\begin{rem}\label{rem_delta_standard} Let $K$ be a generic knot. There is 
$$\delta_0\in\left(0, \hbox{min}\big\lbrace \delta_{diag}, \delta_{(\overline{s}_1, \overline{s}_2)}\,\vert\, (\overline{s}_1, \overline{s}_2)\hbox{ special pair}\big\rbrace\right)$$
such that if $\delta_K\in(0, \delta_0)$, then weak diagonal and weak special sets are disjoint. Let us fix in this Chapter some $\delta_K$. The complement of the corresponding weak special and weak diagonal sets is called the \textbf{standard set $S_K$}. 
\end{rem}

\section{Topology of $M_{K, \varepsilon}$}
In this section, we would like to inspect the manifold structure of $M_{K, \varepsilon}\setminus\Delta_\varepsilon$ provided that $K$ is generic and the radius $\varepsilon>0$ is small. We aim to write $M_{K, \varepsilon}$ as a broken fibration over $K\times K$ and describe the fibers.

\begin{defn}\label{defn_fake}Let $M$ be a $k$-dimensional boundaryless manifold and $N$ be a $n$-product of copies of $\R$ or $\R/T\mathbb{Z}$. Let $F=(F_1,\dots,F_n):M\rightarrow N$ be a continuous map and let $Q=I_1\times\dots\times I_n\subset N$ be a product of intervals. Also let us consider a subset $\lbrace i_1,\dots, i_j\rbrace\subset\lbrace 1, \dots, n\rbrace$. Then we put
$$\mathcal{Q}_{i_1,\dots, i_j}=\big\lbrace x\in Q\,|\,\pi_{\ell}(x)\in\partial I_j\hbox{ iff }\ell\in\lbrace i_1,\dots, i_j\rbrace\big\rbrace,$$
where $\pi_\ell$ is the canonical projection onto $\ell$-th component of $N$. In particular, $\mathcal{Q}$ is the interior $Int(Q)$.

The set
$$\partial^{fake}F^{-1}(Q)=\bigcup_{i_1,\dots, i_j\in\lbrace k-1, \dots, n\rbrace}F^{-1}(\mathcal{Q}_{i_1,\dots, i_j})$$ 
is called the \textbf{fake boundary of $F^{-1}(Q)$}. Here recall that $\dim(M)=k$. Next, the set
$$\partial^{real}F^{-1}(Q)=F^{-1}(Q)\setminus \big( F^{-1}(\mathcal{Q})\cup\partial^{fake}F^{-1}(Q)\big)$$
is called the \textbf{real boundary of $F^{-1}(Q)$}.
\end{defn}

\begin{rem} In our applications, the manifold $M$ from Definition \ref{defn_fake} will be a submanifold of $(\R/T\mathbb{Z}\times S^1)^2$ or $(\R/T\mathbb{Z}\times S^1)\times\R/T\mathbb{Z}$. In particular, $k$ will be $4$ or $3$.

We will aim to split the space $M_{K, \varepsilon}$ into several pieces, which will be easier to describe. So the notion of the real and the fake boundary will help us to distinguish which strata were made artificially by the splittings.
\end{rem}

\begin{thm}\label{thm_mfld_coners}
Let us restrict ourselves outside of the weakly diagonal set (i.e., to the set $J:=\lbrace (s_1, \theta_1, s_2, \theta_2)\in(\R/T\mathbb{Z}\times S^1)^2\,|\,\widetilde{d}(s_1, s_2)\geq\delta_K\rbrace$) and put $\widehat{M}_{K, \varepsilon}:=M_{K, \varepsilon}\vert_{J}$. Then there is a $\varepsilon_0\in(0, \varepsilon_{good}]$ such that for each $\varepsilon\in[0, \varepsilon_0]$ it holds that $\widehat{M}_{K, \varepsilon}$ is a compact $4$-dimensional submanifold with corners. Moreover, all of these $\widehat{M}_{K, \varepsilon}$ are isotopic.

The stratification of each $\widehat{M}_{K, \varepsilon}$ is induced by the stratum transverse intersections of $F^{[\varepsilon]}=(F^{[\varepsilon]}_1, F^{[\varepsilon]}_2, F_{dist})$ and $Q=[0, \infty)\times[0, \infty)\times[\delta_K, -\delta_K]$.
\end{thm}

\begin{proof}
If we verify that $F^{[0]}$ is stratum transverse to $Q$, then the theorem follows from Stability Lemma \ref{lem_stability}. 

But first, we would like to do a few auxiliary computations and describe the matrix $\left(dF_1^{[\varepsilon]}, dF_2^{[\varepsilon]}\right)^T$ (at some general point of $(\R/T\mathbb{Z}\times S^1)^2$). We will use the computation also for the study of the diagonal $\Delta_\varepsilon$ (Definition \ref{defn_out_out}), and hence we do not put $\varepsilon:=0$ immediately.

We compute few partial derivatives for $F^{[\varepsilon]}_1$
\begin{align*}
 \partial_{s_1}F^{[\varepsilon]}_1&=\langle\Gamma_\varepsilon(s_2, \theta_2)-\gamma(s_1), \partial_{s_1}v(s_1, \theta_1)\rangle,\\
 \partial_{\theta_1}F^{[\varepsilon]}_1&=\langle\Gamma_\varepsilon(s_2, \theta_2)-\gamma(s_1), \partial_{\theta_1}v(s_1, \theta_1)\rangle\\
 \partial_{s_2}F^{[\varepsilon]}_1&=\langle\partial_{s_2}\Gamma_\varepsilon(s_2, \theta_2), v(s_1, \theta_1)\rangle,\\
 \partial_{\theta_2}F^{[\varepsilon]}_1&=\langle\partial_{\theta_2}\Gamma_\varepsilon(s_2, \theta_2), v(s_1, \theta_1)\rangle.
\end{align*}
One can show that
\begin{align*}
 \partial_{s_1}F^{[\varepsilon]}_1-\tau(s_1)\partial_{\theta_1}F^{[\varepsilon]}_1&=-\cos(\theta_1)\kappa(s_1)\big\langle\gamma(s_2)-\gamma(s_1)+\varepsilon\cos(\theta_2)n(s_2)+\varepsilon\sin(\theta_2)b(s_2), \dot{\gamma}(s_1)\big\rangle,\\
 \partial_{s_2}F^{[\varepsilon]}_1-\tau(s_2)\partial_{\theta_2}F^{[\varepsilon]}_1&=(1-\varepsilon\cos(\theta_2)\kappa(s_2))\big\langle\cos(\theta_1) n(s_1)+\sin(\theta_1)b(s_1), \dot{\gamma}(s_2)\big\rangle.
\end{align*}
Here we recall that $1-\varepsilon\cos(\theta_2)\kappa(s_2)\neq 0$ due to our choice of $\varepsilon$.

The case of $F^{[\varepsilon]}_2$ is completely analogous; we only point out that
\begin{align*}
 \partial_{s_1}F^{[\varepsilon]}_2&=-\langle\partial_{s_1}\Gamma_\varepsilon(s_1, \theta_1), v(s_2, \theta_2)\rangle,\\
 \partial_{\theta_1}F^{[\varepsilon]}_2&=-\langle\partial_{\theta_1}\Gamma_\varepsilon(s_1, \theta_1), v(s_2, \theta_2)\rangle,\\
 \partial_{s_2}F^{[\varepsilon]}_2&=-\langle\Gamma_\varepsilon(s_1, \theta_1)-\gamma(s_2), \partial_{s_2}v(s_2, \theta_2)\rangle,\\
 \partial_{\theta_2}F^{[\varepsilon]}_2&=-\langle\Gamma_\varepsilon(s_1, \theta_1)-\gamma(s_2), \partial_{\theta_2}v(s_2, \theta_2)\rangle.
\end{align*}

Let us introduce the notation
$$v_i:=v(s_i, \theta_i)\hbox{ and }v_i^\bot:=\partial_{\theta_i}v_i,$$
$$d_i:=1-\varepsilon\cos(\theta_i)\kappa(s_i),$$
where $i=1, 2$. Last, we will denote $$P:=\gamma(s_2)-\gamma(s_1).$$

Hence the matrix $\left(dF_1^{[\varepsilon]}, dF_2^{[\varepsilon]}\right)^T$ is of the form
\begin{equation}
\hspace*{-3cm}\begin{pmatrix}\label{eqn_matrix_epsioln_nonzero}
-\cos(\theta_1)\kappa(s_1)\langle P+\varepsilon v_2, \dot{\gamma}(s_1)\rangle+\tau(s_1)\langle P+\varepsilon v_2, v_1^{\bot}\rangle & -d_1\kappa(s_1))\langle v_2, \dot{\gamma}(s_1)\rangle-\tau(s_1)\langle v_2, \varepsilon v_1^{\bot}\rangle \\
\langle P+\varepsilon v_2, v_1^{\bot}\rangle & -\langle v_2, \varepsilon v_1^{\bot}\rangle \\
d_2\langle v_1, \dot{\gamma}(s_2)\rangle+\tau(s_2)\langle v_1, \varepsilon v_2^{\bot}\rangle & \cos(\theta_2)\kappa(s_2)\langle -P+\varepsilon v_1, \dot{\gamma}(s_2)\rangle-\tau(s_2)\langle -P+\varepsilon v_1, v_2^{\bot}\rangle \\
\langle v_1, \varepsilon v_2^{\bot}\rangle & -\langle -P+\varepsilon v_1, v_2^{\bot}\rangle
\end{pmatrix}
\end{equation}

And thus $\left(dF_1^{[0]}, dF_2^{[0]}\right)^T$ is equal to the matrix
\begin{align}
\begin{split}\label{eqn_matrix_epsioln_zero}
A&=\\
&\begin{pmatrix}
-\cos(\theta_1)\kappa(s_1)\langle P, \dot{\gamma}(s_1)\rangle+\tau(s_1)\langle P, v_1^{\bot}\rangle & -\langle v_2, \dot{\gamma}(s_1)\rangle \\
\langle P, v_1^{\bot}\rangle & 0 \\
\langle v_1, \dot{\gamma}(s_2)\rangle & -\cos(\theta_2)\kappa(s_2)\langle P, \dot{\gamma}(s_2)\rangle+\tau(s_2)\langle P, v_2^{\bot}\rangle \\
0 & \langle P, v_2^{\bot}\rangle
\end{pmatrix}
\end{split}
\end{align}

We will denote the elements of the matrix $\left(dF_1^{[0]}, dF_2^{[0]}\right)^T$ as $\lbrace a_{i, j}\rbrace_{i, j}$.

Let us show that $F^{[0]}\pitchfork\mathcal{Q}_{1, 2}$. So we need to show that $A$ has rank $2$ along the solution set of 
\begin{align}
\langle P, v_1\rangle&=0, \label{eqn_morse1_red}\\
\langle P, v_2\rangle&=0, \label{eqn_morse2_red}\\
\widetilde{d}(s_1, s_2)&>\delta_K. \label{eqn_morse3_red}
\end{align}

Arguing by contradiction, we assume that $\hbox{rank}(A)<2$.

First, we inspect the cases where the second column is $0$. If $\cos(\theta_2)\neq 0$, then by $a_{3, 2}$ and $a_{4, 2}$ the vector $P$ is parallel to $v_2$. Which is a contradiction by equation (\ref{eqn_morse2_red}) and the fact that $s_1\neq s_2$.

If $\cos(\theta_2)=0$, then $v_2=\pm b(s_2)$. Also, by $a_{4, 2}$ and equation (\ref{eqn_morse2_red}), the vector $P$ is parallel to $\dot{\gamma}(s_2)$. So by our generic assumptions on $K$, we have that $P$ is not parallel to $\dot{\gamma}(s_1)$. Also $s_1\neq s_2$. Thus $(s_1, s_2)$ is a $s_2$-special pair. However, in that case we have assumed that $\dot{\gamma}(s_1)$ is not lying in the osculating plane of $\gamma(s_2)$. But that is a contradiction with $a_{1, 2}$.

Hence, there exists $t\in\R$ such that $t$-times the second column of $A$ is the first column. If $t\neq0$, then from the second and fourth row of $A$ and equations (\ref{eqn_morse1_red}) and (\ref{eqn_morse2_red}) we obtain that $P$ is parallel to both $\dot{\gamma}(s_1)$ and $\dot{\gamma}(s_2)$. But that is not possible for a generic $K$.

If $t=0$, then the first column is $0$, which is analogous to the case when the second column is $0$. Hence we obtained contradiction and thus $$F^{[0]}\pitchfork\mathcal{Q}_{1, 2}.$$

Moreover, observe that when we assumed that one column of $A$ is zero, then we actually used for the contradiction always only one of the equations (\ref{eqn_morse1_red}) and (\ref{eqn_morse2_red}). Hence, in fact, we also showed that
$$F^{[0]}\pitchfork\mathcal{Q}_{1}\hbox{ and }F^{[0]}\pitchfork\mathcal{Q}_{2}.$$

Also, trivially
$$F^{[0]}\pitchfork\mathcal{Q}.$$

Now, it remains to inspect the strata that contains the boundary of $[\delta_K, -\delta_K]$. Note that at each point of $(\R/T\mathbb{Z}\times S^1)^2$
$$dF_{dist}=(-1, 0, 1, 0).$$

Let us show that $F^{[0]}\pitchfork\mathcal{Q}_{1, 2, 3}.$ Let $x\in (F^{[0]})^{-1}(\mathcal{Q}_{1, 2, 3})$. Note that by the construction of weakly diagonal sets, it holds that $\pi_{s_1, s_2}(x)$ is not a special pair. In particular, $a_{2, 1}$ and $a_{4, 2}$ are non-zero. Hence from the description of $dF_{dist}(x)$ it follows that $\hbox{rank}(dF^{[0]}(x))=3$, and thus
$$F^{[0]}\pitchfork\mathcal{Q}_{1, 2, 3}.$$

The cases $F^{[0]}\pitchfork\mathcal{Q}_{1, 3}$ and $F^{[0]}\pitchfork\mathcal{Q}_{2, 3}$ follow from the same argument as $F^{[0]}\pitchfork\mathcal{Q}_{1, 2, 3}$ did. 

Finally, the case $F^{[0]}\pitchfork\mathcal{Q}_{3}$ is immediate.
\end{proof}

\begin{notat}\label{notat_torus_first} We recall that in the proof of Theorem \ref{thm_mfld_coners} we introduced the notation
$$v_i:=v(s_i, \theta_i)\hbox{ and }v_i^\bot:=\partial_{\theta_i}v_i,$$
$$d_i:=1-\varepsilon\cos(\theta_i)\kappa(s_i),$$
where $i=1, 2$. Last, we denoted $$P:=\gamma(s_2)-\gamma(s_1).$$
\end{notat}

\begin{rem} The proofs of the following lemmata - Lemma \ref{lem_mfld_special} and Lemma \ref{lem_mfld_standard} are analogous to the proof of the Theorem \ref{thm_mfld_coners}.
\end{rem}

\begin{lemma}\label{lem_mfld_special}
Let us restrict ourselves over some weakly special set $W_{(s_1, s_2)}$ (i.e. to the set $J:=\lbrace (s_1, \theta_1, s_2, \theta_2)\in(\R/T\mathbb{Z}\times S^1)^2\,|\,(s_1, s_2)\in W_{(s_1, s_2)}\rbrace$) and put $\widehat{M}_{K, \varepsilon}=M_{K, \varepsilon}\vert_J$. Then there is a $\varepsilon_0\in(0, \varepsilon_{good}]$ such that for each $\varepsilon\in[0, \varepsilon_0]$ it holds that $\widehat{M}_{K, \varepsilon}$ is a $4$-dimensional compact submanifold with corners. Moreover, all of these $\widehat{M}_{K, \varepsilon}$ are isotopic.

The stratification of each $\widehat{M}_{K, \varepsilon}$ is induced by the stratum transverse intersections of $F^{[\varepsilon]}=(F^{[\varepsilon]}_1, F^{[\varepsilon]}_2, Id, Id)$ and 
$$Q=[0, \infty)\times[0, \infty)\times W_{(s_1, s_2)}.$$
\end{lemma}

\begin{lemma}\label{lem_mfld_standard}
Let us restrict ourselves over $\overline{S_K}$ (i.e. to the set $J:=\lbrace (s_1, \theta_1, s_2, \theta_2)\in(\R/T\mathbb{Z}\times S^1)^2\,|\,(s_1, s_2)\in \overline{S_K}\rbrace$, where $S_K$ was introduced in Remark \ref{rem_delta_standard}). We put and put $\widehat{M}_{K, \varepsilon}=M_{K, \varepsilon}\vert_J$. Then there is a $\varepsilon_0\in(0, \varepsilon_{good}]$ such that for each $\varepsilon\in[0, \varepsilon_0]$ it holds that $\widehat{M}_{K, \varepsilon}$ is a $4$-dimensional compact submanifold with corners. Moreover, all of these $\widehat{M}_{K, \varepsilon}$ are isotopic.

The stratification of each $\widehat{M}_{K, \varepsilon}$ is induced by the stratum transverse intersections of $F^{[\varepsilon]}=(F^{[\varepsilon]}_1, F^{[\varepsilon]}_2, F_{dist}, Id,\dots,Id)$ and 
$$Q=[0, \infty)\times[0, \infty)\times[\delta_K, -\delta_K]\times \prod_{(s_1, s_2) \text{ special pair}}\big((\R/T\mathbb{Z}\times S^1)^2\setminus Int(W_{(s_1, s_2)})\big).$$
\end{lemma}

\begin{rem} In order to understand $M_{K, \varepsilon}$ outside of a weakly diagonal set, it will be enough by Lemmata \ref{lem_mfld_special} and \ref{lem_mfld_standard} to study only $M_{K, 0}$.
\end{rem}

\begin{lemma}\label{lemma_standard_fibre} Let $\widehat{M}_{K, 0}:=M_{K, 0}\vert_{\pi_{s_1, s_2}^{-1}(S_K)}$. Then the canonical projection $\pi_{s_1,s_2}:\widehat{M}_{K,0}\rightarrow S_{K}$ defines a locally trivial fiber bundle.
\end{lemma}

\begin{proof} Let $\widehat{F}^{[0]}=(\widehat{F}_1^{[0]}, \widehat{F}^{[0]}_2)$ denotes the restriction of $F^{[0]}=(F_1^{[0]}, F^{[0]}_2)$ to $\lbrace (s_1, \theta_1, s_2, \theta_2)\in(\R/T\mathbb{Z}\times S^1)^2\,|\,(s_1, s_2)\in \overline{S_K}\rbrace$. By Lemma \ref{lem_mfld_standard} $\widehat{M}_{K, 0}$ is a $4$-manifold with corners and the stratification is induced by the stratum intersections of $\widehat{F}^{[0]}=(\widehat{F}_1^{[0]}, \widehat{F}^{[0]}_2)$ and $Q=[0, \infty)\times[0, \infty)$. Our aim is to apply Ehresmann's Theorem \ref{thm_ehrsm} to the projection $\pi_{s_1, s_2}:\widehat{M}_{K, 0}\rightarrow S_{K}$.

$\pi_{s_1, s_2}$ is smooth and surjective. Also, since the fibers are compact, $\pi_{s_1, s_2}$ is proper. We need to show that $\pi_{s_1, s_2}$ is a submersion on each stratum. 

On $(\widehat{F}^{[0]})^{-1}(\mathcal{Q})$, the projection $\pi_{s_1, s_2}$ is a submersion, since there is locally no condition on $s_1$ and $s_2$.

Let $x=(s_1, \theta_1, s_2, \theta_2)\in(\widehat{F}^{[0]})^{-1}(\mathcal{Q}_1)$. Then 
$$T_x\big((\widehat{F}^{[0]})^{-1}(\mathcal{Q}_1)\big)=\big\lbrace v\in T_x(\R/T\mathbb{Z}\times S^1)^2\,|\, d\widehat{F}^{[0]}_1(x)v=0\big\rbrace.$$
We would like to show that there are $w_1, w_2\in  T_x\big((\widehat{F}^{[0]})^{-1}(\mathcal{Q}_1)\big)$ such that $d\pi_{s_1, s_2}(x)w_1=\partial_{s_1}$ and $d\pi_{s_1, s_2}(x)w_2=\partial_{s_2}$. Hence, it will be sufficient to verify that
$$\hbox{rank}(d\widehat{F}^{[0]}_1|_{\pi_{s_1, s_2}^{-1}(s_1, s_2)}( \theta_1, \theta_2))=1.$$
Since we are outside of (weakly) special and (weakly) diagonal pairs, $\langle P, v_1^\bot\rangle\neq 0$. Hence, the rank condition immediately follows from matrix (\ref{eqn_matrix_epsioln_zero}).

We treat the remaining two cases similarly.
\end{proof}

\begin{rem}We are going to replace Lemma \ref{lemma_standard_fibre} with a stronger result - Lemma \ref{lemma_standard_square}. However, the proof of Lemma \ref{lemma_standard_fibre} is still conceptually useful, since we would like to describe globally $M_{K, 0}$ as a broken fibration - Definition \ref{defn_broken}.
\end{rem}

\begin{lemma}\label{lemma_standard_square} Over $\overline{S_{K}}$ the manifold $M_{K, 0}$ is diffeomorphic to $\overline{S_{K}}\times[-\pi/2, \pi/2]\times[-\pi/2, \pi/2]$.
\end{lemma}

\begin{proof}
Similarly to Lemma \ref{lemma_standard_fibre}, let $\widehat{M}_{K, 0}$ denotes the restriction of $M_{K, 0}$ over $\overline{S_K}$, which is by Lemma \ref{lem_mfld_standard} a $4$-manifold with corners.

Observe that over $\overline{S_{K}}$ the projections of $P$ into normal planes at $\gamma(s_1)$ and $\gamma(s_2)$ are nowhere-vanishing vectors. Hence, for $i=1, 2$, we can make the following changes of smooth orthonormal frames $\lbrace n(s_i), b(s_i)\rbrace$ to some $\lbrace n^{\ast}_i(s_1, s_2), b^{\ast}_i(s_1, s_2)\rbrace$ which are given by
\begin{align*}
n(s_i)&\longmapsto n^{\ast}_i(s_1, s_2):=D_i^{-1/2}\big(\langle P, n(s_i)\rangle n(s_i)+\langle P, b(s_i)\rangle b(s_i)\big),\\
b(s_i)&\longmapsto b^{\ast}_i(s_1, s_2):=D_i^{-1/2}\big(-\langle P, b(s_i)\rangle n(s_i)+\langle P, n(s_i)\rangle b(s_i)\big),
\end{align*}
where $D_i=\langle P, n(s_i)\rangle^2+\langle P, b(s_i)\rangle^2.$ The frames $\lbrace n^{\ast}_i(s_1, s_2), b^{\ast}_i(s_1, s_2)\rbrace$ are also smooth.

Let $\theta^\ast_i\in S^1$ be angular coordinates defined by (positively oriented) orthonormal frames $\lbrace n^{\ast}_i(s_1, s_2), b^{\ast}_i(s_1, s_2)\rbrace$. Hence for each $(s_1, s_2)\in \overline{S_K}$ the change of frames gives us an automorphims on $S^1$
$$\varphi_{i; (s_1, s_2)}:\theta_i\longmapsto\theta^\ast_i(s_1, s_2).$$

Note that $\langle P, \pm b^{\ast}_i(s_1, s_2)\rangle=0$. In other words, for $\theta_i=\varphi_{i; (s_1, s_2)}^{-1}(\pm \pi/2)$ it holds that $F^{[0]}_i=0$. We remark that each of $F^{[0]}_i$ depends, besides $s_1, s_2,$ only on one of angular coordinates and that is $\theta_i$. Moreover, for $\theta_i\in\varphi_{i; (s_1, s_2)}^{-1}((-\pi/2, \pi/2))$ holds that $F^{[0]}_i>0$. Altogether, there is a (strata preserving) diffeomorphism $\widehat{M}_{K, 0}\rightarrow\overline{S_{K}}\times[-\pi/2, \pi/2]\times[-\pi/2, \pi/2]$ which is given by
$$(s_1, \theta_1, s_2, \theta_2)\mapsto(s_1, s_2, \varphi_{1; (s_1, s_2)}(\theta_1), \varphi_{2; (s_1, s_2)}(\theta_2)).$$
\end{proof}

\begin{lemma}\label{lemma_line_product} Let $(\overline{s}_1, \overline{s}_2)$ be a $\overline{s}_2$-special pair. Then over $W_{(\overline{s}_1, \overline{s}_2)}$ the $4$-manifold $M_{K, 0}$ is diffeomorphic to $N_{K, 0}\times [-\pi/2, \pi/2]$, where $N_{K, 0}$ is a compact $3$-manifold with corners that depends only on $(s_1, \theta_2, s_2)$. In more detail, $N_{K, 0}$ is stratified by the stratum transverse intersection of $F^{[0]}=(F^{[0]}_2, Id, Id)$ and $Q=[0, \infty)\times W_{(\overline{s}_1, \overline{s}_2)}$ (here $F^{[0]}$ is the function in variables $(s_1, \theta_2, s_2)$).

The analogous statement holds for $\overline{s}_1$-special pairs.
\end{lemma}

\begin{proof}
Let $(\overline{s}_1, \overline{s}_2)$ be a $\overline{s}_2$-special pair. Let $\widehat{M}_{K, 0}$ denotes the restriction of $M_{K, 0}$ over $W_{(\overline{s}_1, \overline{s}_2)}$, which is a $4$-manifold with corners by Lemma \ref{lem_mfld_special}. Similarly to Lemma \ref{lemma_standard_square} we can write $N_{K, 0}\times[-\pi/2, \pi/2]$ as the image of the embedding of $\widehat{M}_{K, 0}$ which is given as
$$(s_1, \theta_1, s_2, \theta_2)\mapsto(s_1, s_2, \theta_2, \varphi_{1; (s_1, s_2)}(\theta_1)).$$
\end{proof}

\begin{defn}\label{defn_broken}Let $M$ be a manifold with corners and $f:M\rightarrow \R$ be a smooth function.

A point $p\in M$ is called a \textbf{critical point of $f$} if $df(p)|_{T_p \mathcal{Q}}=0$, where $\mathcal{Q}$ is the stratum containing $p$. By $Crit(f)$, we will denote the \textbf{set of critical points of $f$}.

The function $f$ is called \textbf{Morse} if:
\begin{itemize}
\item[$(i.)$] $f$ is proper.
\item[$(ii.)$] For each stratum $\mathcal{Q}$, the critical points of $f|_\mathcal{Q}$ are nondegenerate.
\item[$(iii.)$] If $p\in \mathcal{Q}$ is a critical point, then $df(p)V\neq 0$, where $V$ is any vector space that satisfies the following property. There is a stratum $\mathcal{M}$, higher then the stratum $\mathcal{Q}$, and a sequence of points $p_j\in\mathcal{M}$ converging to $p$ such that $\lim_{p_j\rightarrow p}T_{p_j}\mathcal{M}=V$.
\end{itemize}

Let $I\subset \R$ be an interval. In addition, if the Morse function $f$ is surjective to $I$, then the restriction $f:f^{-1}(I)\rightarrow I$ is called a \textbf{broken fibration}.
\end{defn}

\begin{rem}See also \cite{goresky2012stratified} for the notion of Morse functions for (Whitney) stratified spaces. There is also a terminology of broken Lefschetz fibrations, for example, in \cite{lekili2012heegaardfloerhomologybroken}.
\end{rem}

\begin{lemma}\label{lemma_crit} Let us consider any $\overline{s}_2$-special pair $(\overline{s}_1, \overline{s}_2)$. And put $\widehat{N}_{K, 0}:=N_{K, 0}\vert_{\pi_{s_1}^{-1}(\overline{s}_1-\delta_K, \overline{s}_1+\delta_K)}$. Then the base projection $\pi_{s_1}:\widehat{N}_{K, 0}\rightarrow (\overline{s}_1-\delta_K, \overline{s}_1+\delta_K)$ is a broken fibration with exactly two critical points. In more detail, the critical points are $(\overline{s}_1, \overline{s}_2, \pm\pi/2)$, they lie in the real boundary of $N_{K, 0}$ and are of the Morse index $1$.

The analogous statement also holds for $s_1$-special pairs.
\end{lemma}

\begin{rem} The geometric picture of Lemma \ref{lemma_crit} will be described in Remark \ref{rem_pict_special} and Figure \ref{fibre_phase}.
\end{rem}

\begin{proof} 
Let us describe $\widehat{N}_{K, 0}$ as $(F^{[0]})^{-1}(Q)$, where $F^{[0]}=(F^{[0]}_2, Id, Id)$ and $Q=[0,\infty)\times (\overline{s}_1-\delta_K, \overline{s}_1+\delta_K)\times [\overline{s}_2-\delta_K, \overline{s}_2+\delta_K]$. Recall that here we see $F^{[0]}$ only as a function of variables $(s_1, s_2, \theta_2)$. Also $F^{[0]}$ and $Q$ are intersecting stratum transversely, see Theorem \ref{thm_mfld_coners} or Lemma \ref{lem_mfld_special}.

Let us find the critical points. The strategy will be similar to the proof of Lemma \ref{lemma_standard_fibre}, only now we look for the points where $\pi_{s_1}$ fails to be a fibration.

So we immediately obtain that there are no critical points on $(F^{[0]})^{-1}(\mathcal{Q})$ or $(F^{[0]})^{-1}(\mathcal{Q}_3)$. 

Let us consider $(F^{[0]})^{-1}(\mathcal{Q}_1)$. We pick $p\in(F^{[0]})^{-1}(\mathcal{Q}_1)$. Then the matrix $(dF^{[0]}_2(p))^T$ is given by
\begin{align*}
\begin{pmatrix}
 -\langle v_2, \dot{\gamma}(s_1)\rangle \\
 -\cos(\theta_2)\kappa(s_2)\langle P, \dot{\gamma}(s_2)\rangle+\tau(s_2)\langle P, v_2^{\bot}\rangle \\
\langle P, v_2^{\bot}\rangle
\end{pmatrix}.
\end{align*}
Here we used the computation (\ref{eqn_matrix_epsioln_zero}). Hence, $p$ is a critical point iff it satisfies
\begin{align}
-\cos(\theta_2)\kappa(s_2)\langle P, \dot{\gamma}(s_2)\rangle+\tau(s_2)\langle P, v_2^{\bot}\rangle&=0, \label{eqn_spec_1}\\
\langle P, v_2^{\bot}\rangle&=0, \label{eqn_spec_2}\\
\langle P, v_2\rangle&=0, \label{eqn_spec_3}\\
\pi_{s_1, s_2}(p)&\in Int(W_{(\overline{s}_1, \overline{s}_2)}). \label{eqn_spec_4}
\end{align}
By (\ref{eqn_spec_4}) we are outside of the weak diagonal set, and in particular $P\neq 0$. Since $\lbrace\dot{\gamma}(s_2), v_2, v_2^\bot\rbrace$ is an orthonormal basis, by (\ref{eqn_spec_2}) and (\ref{eqn_spec_3}) $P$ is parallel to $\dot{\gamma}(s_1)$. So $\pi_{s_1, s_2}(p)=(\overline{s}_1, \overline{s}_2)$, and in particular $\langle P, \dot{\gamma}(s_2)\rangle\neq 0$. Hence, by (\ref{eqn_spec_1}) $\cos(\theta_2)=0$, and thus $\theta_2=\pm\pi/2$ and $v_2=\pm b_2(s_2)$.

Let $p=(\overline{s}_1, \overline{s}_2, \pi/2)$, the other case will be treated analogously. Now, we are going to compute the Morse index of $\pi_{s_1}$ at $p$. Let us consider the implicit equation
$$F^{[0]}_2(s_1 ,s_2, \theta_2)=0.$$

Since $F^{[0]}$ and $Q$ intersect stratum transversely and $p$ is a critical point of $\pi_{s_1}$, it holds that $\partial_{s_1}F^{[0]}_2(p)\neq 0$.

So we can, by the Implicit function theorem, describe the Hessian matrix $\widetilde{H}[\pi_{s_1}](p)$. Here the tilde describes that the Hessian matrix is with respect to the flat metric $\widetilde{\langle\cdot,\cdot\rangle}$ on $\R/T\mathbb{Z}\times S^1$. Thus, we obtain
\begin{equation*}
\widetilde{H}[\pi_{s_1}](p)=\frac{-1}{\partial_{s_1}F^{[0]}_2(p)}\begin{pmatrix}
\partial_{s_2, s_2}F^{[0]}_2(p) & \partial_{s_2, \theta_2}F^{[0]}_2(p)\\
\partial_{\theta_2, s_2}F^{[0]}_2(p) & \partial_{\theta_2, \theta_2}F^{[0]}_2(p)\\
\end{pmatrix}.
\end{equation*}
In order to compute the sign of the $\det(\widetilde{H}[\pi_{s_1}](p))$, we can ignore the fraction before $\widetilde{H}[\pi_{s_1}](p)$.

Let us compute 
\begin{align*}
\partial^2_{s_2}F_2(p)&=\at{\partial_{s_2}\langle P, -\cos(\theta_2)\kappa(s_2)\dot{\gamma}(s_2)+\tau(s_2)v_2^\bot\rangle}{(s_1, s_2, \theta_2)=p} \\
&=\at{\partial_{s_2}\langle P, \tau(s_2)n(s_2)\rangle}{(s_1, s_2, \theta_2)=p} \\
&=\langle\dot{\gamma}(\overline{s}_2), \tau(\overline{s}_2)n(\overline{s}_2)\rangle+\left\langle P, \at{\partial_{s_2}\big(\tau(s_2)n(s_2)\big)}{(s_1, s_2, \theta_2)=p}\right\rangle \\
&=\langle P, \dot{\tau}(\overline{s}_2)n(\overline{s}_2)-\tau(\overline{s}_2)\kappa(\overline{s}_2)\dot{\gamma}(\overline{s}_2)+\tau(\overline{s}_2)^2\kappa(\overline{s}_2)b(\overline{s}_2)\rangle,\\
&=-\tau(\overline{s}_2)\kappa(\overline{s}_2)\langle P, \dot{\gamma}(\overline{s}_2)\rangle,\\
\\
\partial_{\theta_2}\partial_{s_2}F_2(p)&=\at{\partial_{\theta_2}\langle P, -\cos(\theta_2)\kappa(s_2)\dot{\gamma}(s_2)+\tau(s_2)v_2^\bot\rangle}{(s_1, s_2, \theta_2)=p} \\
&=\langle P, \sin(\theta_2)\kappa(s_2)\dot{\gamma}(s_2)-\tau(s_2)v_2\rangle\vert_{(s_1, s_2, \theta_2)=p} \\
&=\langle P, \kappa(\overline{s}_2)\dot{\gamma}(\overline{s}_2)-\tau(\overline{s}_2)b(\overline{s}_2)\rangle, \\
&=\kappa(\overline{s}_2)\langle P, \dot{\gamma}(\overline{s}_2)\rangle,\\
\\
\partial^2_{\theta_2}F_2(p)&=\at{\partial_{\theta_2}\langle P, v_2^\bot\rangle}{(s_1, s_2, \theta_2)=p} \\
&=-\langle P, b(\overline{s}_2)\rangle,\\
&=0.
\end{align*}
Hence, 
$$\det(\widetilde{H}[\pi_{s_1}](p))=-\kappa(\overline{s}_2)^2\langle P, \dot{\gamma}(\overline{s}_2)\rangle^2<0$$
due to our generic assumptions on $K$. Thus the Morse index of $\pi_{s_1}$ at $p$ is equal to $1$.

Since the fibers are compact, $\pi_{s_1}$ is proper. $\pi_{s_1}$ is also surjective.

It remains to verify that the condition $(iii.)$ from Definition \ref{defn_broken}. But $\widehat{N}_{K, 0}$ is a $\codim 0$ submanifold of $(\R/T\mathbb{Z}\times S^1)\times \R/T\mathbb{Z}$. Hence from the existence of a collar neighbourhood near any critical point $p\in (F^{[0]})^{-1}(\mathcal{Q}_1)$ we obtain that for any sequence of points $p_j\in(F^{[0]})^{-1}(\mathcal{Q})$ converging to $p$ it holds
$$\lim_{p_j\rightarrow p}T_{p_j}\left((F^{[0]})^{-1}(\mathcal{Q}_1)\right)=T_p\left((\R/T\mathbb{Z}\times S^1)\times \R/T\mathbb{Z}\right)\supset\partial_{s_1}\vert_p.$$

The lemma follows.  
\end{proof}

\begin{rem} Now, we would like to depict $N_{K, 0}$. However, the author finds it easier to visualize, instead of $N_{K, 0}$, rather its certain perturbation: $N_{K, \varepsilon}$.
\end{rem}

\begin{defn} Let $\varepsilon\in[0, \varepsilon_{good}]$ and $(\overline{s}_1, \overline{s}_2)$ be a $\overline{s}_2$-special pair. Then we put $N_{K, \varepsilon}:=(\widetilde{F}^{[\varepsilon]})^{-1}(Q)$, where $\widetilde{F}^{[\varepsilon]}=(F^{[0]}_2+\varepsilon, Id, Id)$ and $Q=[0, \infty)\times W_{(\overline{s}_1, \overline{s}_2)}$ (note that $F^{[0]}_2+\varepsilon=\langle P+\varepsilon v_2, v_2\rangle$).

An analogous statement also holds for $\overline{s}_1$-special pairs.
\end{defn}

\begin{rem} Let $\varepsilon>0$ be small and $(\overline{s}_1, \overline{s}_2)$ be a $\overline{s}_2$-special pair. Then by Stability Lemma \ref{lem_stability} and Lemma \ref{lemma_line_product} $N_{K, 0}$ is isotopic to $N_{K, \varepsilon}$. By Stability Lemma \ref{lem_stability} one can show that $\pi_{s_1}:N_{K, \varepsilon}\rightarrow(\overline{s}_1-\delta_K, \overline{s}_1+\delta_K)$ is also broken fibration with precisely two critical points. The critical points will be the (transversely cut out) solutions of the system
\begin{equation}\label{eqn_crit_system}
\lbrace\theta_2=\pm\pi/2\wedge\langle P, v_2\rangle=\mp\varepsilon\wedge\langle P, v_2^\bot\rangle=0\rbrace.
\end{equation}
Geometrically, the critical points will correspond to the intersections of the piece of the knot $K\vert_{(\overline{s}_1-\delta_K, \overline{s}_1+\delta_K)}$ with the lines
$$\lbrace\gamma(s_2)\pm\varepsilon b(s_2)+r\dot{\gamma}(s_2)\,|\,r\in\R\rbrace,$$
where $s_2\in(\overline{s}_2-\delta_K, \overline{s}_2+\delta_K)$.

Also, by the Implicit function theorem applied to (\ref{eqn_crit_system}), one can show $\pi_{s_1}$ will now have no distinct critical values.
\end{rem}

\begin{rem}\label{rem_pict_special} Let $\varepsilon>0$ be small and $(\overline{s}_1, \overline{s}_2)$ be a $\overline{s}_2$-special pair. We would like to present several visualizations of $N_{K, \varepsilon}$.

Let us shrink the radius $\varepsilon$ of $T_{K ,\varepsilon}$ to $0$ in a neighborhood of $\gamma(\overline{s}_1)$. Hence on $(\overline{s}_1-\delta_K, \overline{s}_1+\delta_K)$ the torus $T_{K, \varepsilon}$ will become $T_{K, 0}=K$, see Figure \ref{special_pair}. So $N_{K, \varepsilon}$ represents the set of chords from $T_{K, 0}\vert_{(\overline{s}_1-\delta_K, \overline{s}_1+\delta_K)}$ to $T_{K, \varepsilon}\vert_{(\overline{s}_2-\delta_K, \overline{s}_2+\delta_K)}$ that are outward-pointing in their endpoints.

Next, the first indication that the function $\pi_{s_1}$ consists of saddle-type critical points can be seen in Figure \ref{fig_fibers}.

The standard $\R^3$ is not suited for the visualization of the whole $N_{K, \varepsilon}$. For this, will rather use $\R/T\mathbb{Z}\times \R/T\mathbb{Z}\times S^1$, see Figure \ref{fibre_phase}.

It is worth recalling that all boundary components shown in Figure \ref{fibre_phase} are not the real boundary components of $N_{K, \varepsilon}$. In other words, in Figure \ref{fibre_phase}, we also draw fake boundary components, which come from the fact that we restricted the $s$-coordinates to $W_{(\overline{s}_1, \overline{s}_2)}$. Fake boundary components (of $N_{K, \varepsilon}$, or more precisely of $N_{K, 0}$) will tell us, how $N_{K, 0}\times[0, 1]$ is glued to the rest of $M_{K, 0}$. It turns out that $N_{K, \varepsilon}$ is homeomorphic to the solid torus $S^1\times D^2$. Moreover, this homeomorphism maps the fake boundary to a thin strip on $S^1\times \partial D^2$ that is homotopic to the concatenation of meridian and longitude (given by Seifert framing). See Figure \ref{fibre_homeo}.
\begin{figure}[!htbp]
\centering
\includegraphics[scale=0.7]{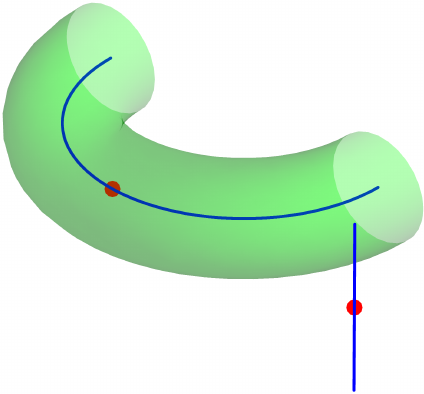}
\vspace{0.3cm}
\caption{A neighborhood of $\overline{s}_2$-strongly special pair in $\R^3$. The red points represent the $\overline{s}_2$-strongly special pair $(\overline{s}_2, \overline{s}_1).$ The blue curves describe the knot $K$, where the vertical part is viewed as $T_{K, 0}$. Around $\gamma(\overline{s}_2)$ the torus $T_{K, \varepsilon}$ is depicted in green.}
\label{special_pair}
\end{figure}

\begin{figure}
\centering
\begin{subfigure}[b]{0.55\textwidth}
\includegraphics[width=0.7\linewidth]{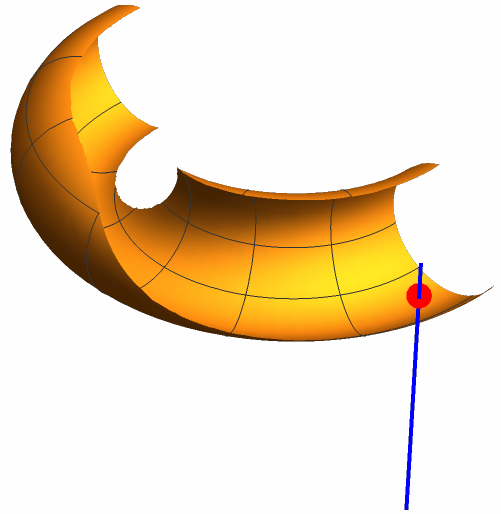}
\end{subfigure}
\begin{subfigure}[b]{0.55\textwidth}
\includegraphics[width=0.7\linewidth]{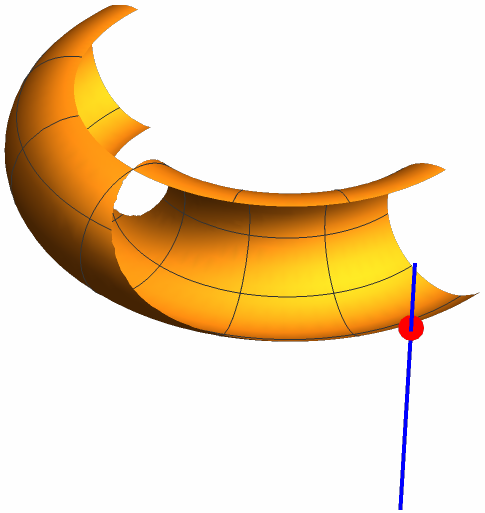}
\end{subfigure}
\begin{subfigure}[b]{0.55\textwidth}
\includegraphics[width=0.7\linewidth]{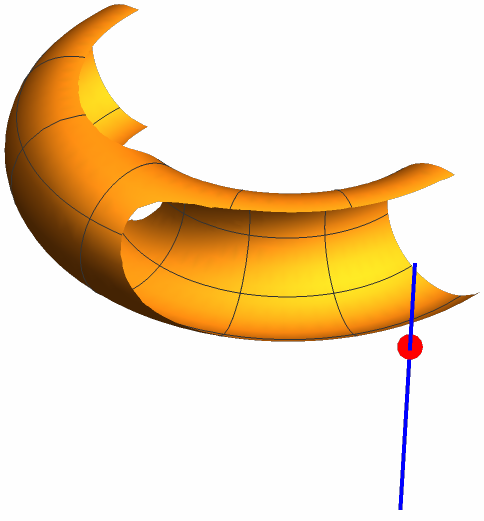}
\end{subfigure}
\caption[Fibers]{A visualization in $\R^3$ of fibers $\pi^{-1}_{s_1}(s_1)\subset N_{K, \varepsilon}$ as $s_1$ pass a singular fiber (the middle figure). The blue lines represent the torus $T_{K, 0}$ of zero radius, which is parametrized by $s_1$. Next, the orange sets describe the sets of all ends of outward-pointing chords that start in the red points. Orange sets are parametrized by $(s_2, \theta_2)$-coordinates.}
\label{fig_fibers}
\end{figure}

\begin{figure}[!htbp]
\labellist
\pinlabel $s_1$ at 20 25
\pinlabel $s_2$ at 192 87
\pinlabel $\theta_2$ at 115 5
\pinlabel $0$ at 65 -10
\pinlabel $2\pi$ at 175 23
\endlabellist
\centering
\includegraphics[scale=0.89]{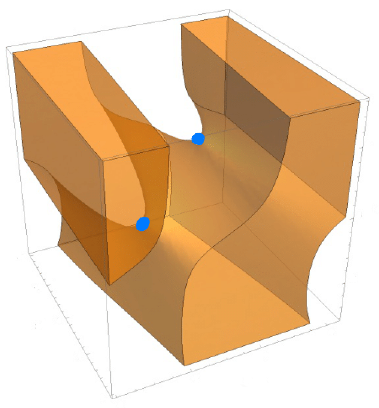}
\vspace{0.3cm}
\caption{A visualization of $N_{K, \varepsilon}$ in coordinates $(s_1, s_2, \theta_2)$. The blue dots represent two critical points of $\pi_{s_1}$ (Lemma \ref{lemma_crit}).}
\label{fibre_phase}
\end{figure}

\begin{figure}[!htbp]
\labellist
\pinlabel $\cong$ at 220 87
\endlabellist
\centering
\includegraphics[scale=0.89]{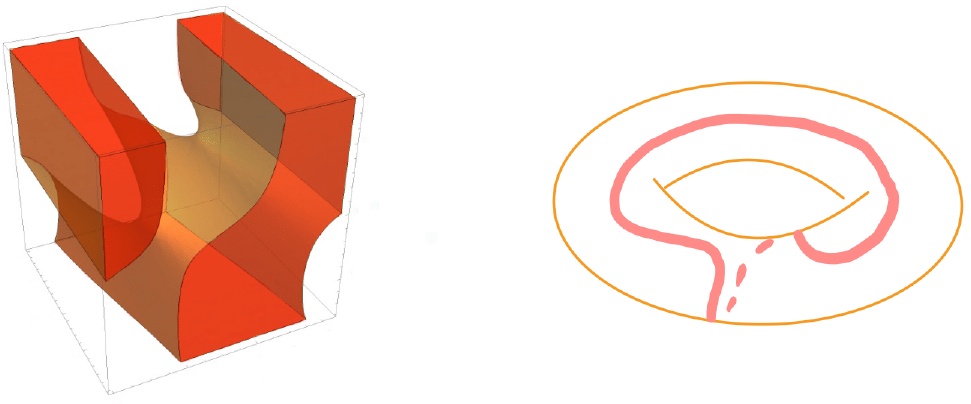}
\vspace{0.3cm}
\caption{The homeomorphism between $N_{K, \varepsilon}$ and $S^1\times D^2$. On the left, the red surface describes the fake boundary, and on the right, we see the image of the fake boundary. The image is homotopic to the concatenation of the meridian and the longitude.}
\label{fibre_homeo}
\end{figure}
\end{rem}

\newpage

The following result is probably originally due to T. Tsuboi, \cite{LoktevaKnots}.

\begin{thm}\label{thm_tsub}If a knot has an empty set of special pairs, then the knot is an unknot.
\end{thm}

\begin{cor}\label{cor_tsub} Let $K$ be generic and $\varepsilon\geq0$ be small. If $H_1^{sing}(M_{K, \varepsilon}|_{\widetilde{d}(s_1, s_2)>\delta_{K}}; \mathbb{Z})\cong\mathbb{Z}$, then $K$ is an unknot.
\end{cor}

\begin{proof}
Let us compute
\begin{align*}
H_1^{sing}(M_{K, \varepsilon}|_{\widetilde{d}(s_1, s_2)>\delta_{K}}; \mathbb{Z})&\stackrel{(1)}{\cong}H_1^{sing}(M_{K, 0}|_{\widetilde{d}(s_1, s_2)>\delta_{K}}; \mathbb{Z})\\
&\stackrel{(2)}{\cong} H_1^{sing}\big((\R/T\mathbb{Z}\times\R/T\mathbb{Z})\setminus(\Delta_0\cup\bigcup_{\substack{x\text{ is a special}\\ \text{point}}}W_x); \mathbb{Z}\big)\\
&\cong \mathbb{Z}^{1+\sharp\text{ special points}},
\end{align*}
where the step $(1)$ follows from Theorem \ref{thm_mfld_coners} and in the step $(2)$ we applied the following two deformations retracts.

First, we deformation retract each $N_{K, 0}\times\lbrace t\rbrace \in N_{K, 0}\times[-1, 1]$ to the fake boundary of the whole product. To obtain such a deformation retract, we identify each $N_{K, 0}\times\lbrace t\rbrace\cong N_{K, 0}$ with $S^1\times D^2$ by the (orientation preserving) homeomorphism from Remark \ref{rem_pict_special}. Then the deformation retraction of the solid torus is straightforward.

Second, by Lemma \ref{lemma_standard_square}, the remaining space is diffeomorphic to the product of $\overline{S_K}$ and the square. So the radial shrinking of the square will give us the second deformation retraction. Last, we can potentially slightly extend the space closer to $\Delta_0$.

Hence, we showed that $H_1^{sing}$ counts the special points. Thus, we can apply Theorem \ref{thm_tsub}.
\end{proof}

\newpage
\begin{conj}\label{conj_diag} Let $\varepsilon\in(0, \varepsilon_{diag})$. Let $\widehat{M}_{K, \varepsilon}$ denotes the restriction of $M_{K, \varepsilon}$ to the set $\lbrace(s_1, \theta_1, s_2, \theta_2)\in(\R/T\mathbb{Z}\times S^1)^2\,|\,\varepsilon\leq\widetilde{d}(s_1, s_2)\leq\delta_K\rbrace$. Then the following holds
\begin{itemize}
\item[$(i.)$]$\widehat{M}_{K, \varepsilon}$ is a $4$-manifold with corners and the projection $\pi_{s_1}:\widehat{M}_{K, \varepsilon}\rightarrow\R/T\mathbb{Z}$ induce a locally trivial fibration. 

In more detail, $\widehat{M}_{K, \varepsilon}=(F^{[\varepsilon]})^{-1}(Q)$, where $F^{[\varepsilon]}=(F^{[\varepsilon]}_1, F^{[\varepsilon]}_2, F^{dist})$ and $Q=[0, \infty)\times[0, \infty)\times\lbrace[-\delta_K, -\varepsilon]\cup[\varepsilon, \delta_K]\rbrace$.
\item[$(ii.)$]Let $\overline{s}_1, \overline{s}_2\in\R/T\mathbb{Z}$ be such that $\overline{s}_1=\overline{s}_2$, i.e. $(\overline{s}_1, \overline{s}_2)$ is a diagonal pair. Then $\pi_{s_1}^{-1}(\overline{s}_1)$ is a $3$-manifold with corners and the projection $\pi_{s_2}:\pi_{s_1}^{-1}(\overline{s}_1)\rightarrow (\overline{s}_2-\delta_K, \overline{s}_2-\varepsilon)\cup(\overline{s}_2+\varepsilon, \overline{s}_2+\delta_K)$ induce a broken fibration.

In more detail, $\pi_{s_1}^{-1}(\overline{s}_1)=(\widetilde{F}^{[\varepsilon]})^{-1}(\widetilde{Q})$, where $$\widetilde{F}^{[\varepsilon]}=(F^{[\varepsilon]}_1\vert_{\pi_{s_1}^{-1}(\overline{s}_1)}, F^{[\varepsilon]}_2\vert_{\pi_{s_1}^{-1}(\overline{s}_1)}, Id)$$ and $$\widetilde{Q}=[0, \infty)\times[0, \infty)\times\lbrace[\overline{s}_2-\delta_K, \overline{s}_2-\varepsilon]\cup[\overline{s}_2+\varepsilon, \overline{s}_2+\delta_K]\rbrace.$$ And the only critical points of $\pi_{s_2}$ are two critical points in $(\widetilde{F}^{[\varepsilon]})^{-1}(\widetilde{\mathcal{Q}}_1)$, each of them has Morse index equal to $1$.
\end{itemize}
\end{conj}

\begin{rem} The geometric picture of Conjecture \ref{conj_diag} will be described in Figure \ref{figure_diagonal_example}.
\end{rem}

\begin{pproof}
Recall from the proof of Theorem \ref{thm_mfld_coners} that the matrix $\left(dF_1^{[\varepsilon]}, dF_2^{[\varepsilon]}\right)^T$ (at some general point of $(\R/T\mathbb{Z}\times S^1)^2$) is given by
\begin{equation*}
\hspace*{-3cm}\begin{pmatrix}
-\cos(\theta_1)\kappa(s_1)\langle P+\varepsilon v_2, \dot{\gamma}(s_1)\rangle+\tau(s_1)\langle P+\varepsilon v_2, v_1^{\bot}\rangle & -d_1\kappa(s_1))\langle v_2, \dot{\gamma}(s_1)\rangle-\tau(s_1)\langle v_2, \varepsilon v_1^{\bot}\rangle \\
\langle P+\varepsilon v_2, v_1^{\bot}\rangle & -\langle v_2, \varepsilon v_1^{\bot}\rangle \\
d_2\langle v_1, \dot{\gamma}(s_2)\rangle+\tau(s_2)\langle v_1, \varepsilon v_2^{\bot}\rangle & \cos(\theta_2)\kappa(s_2)\langle -P+\varepsilon v_1, \dot{\gamma}(s_2)\rangle-\tau(s_2)\langle -P+\varepsilon v_1, v_2^{\bot}\rangle \\
\langle v_1, \varepsilon v_2^{\bot}\rangle & -\langle -P+\varepsilon v_1, v_2^{\bot}\rangle
\end{pmatrix}.
\end{equation*}
We will denote the elements of the matrix by $\lbrace a_{i, j}\rbrace_{i, j}.$
\newline
\newline
\noindent{\textit{Ad $(i.)$:}}

First, we would like to verify the stratum transversality of $F^{[\varepsilon]}$ and $Q$ and then apply Theorem \ref{thm_invers_im}.

Let us show that $F^{[\varepsilon]}\pitchfork\mathcal{Q}_1$. We need to show that 
$$\hbox{rank}(dF_1^{[\varepsilon]}(p))=1$$
for any $p\in (F^{[\varepsilon]})^{-1}(\mathcal{Q}_1)$. Assume for contradiction that this is not true. Then by $a_{3, 1}$ and $a_{4, 1}$ we obtain that $v_1=\pm v_2$. However, we also know that $F^{[\varepsilon]}_1(p)=0$ and $F^{[\varepsilon]}_1(p)>0$, which is a contradiction.

Similarly, we can show that $F^{[\varepsilon]}\pitchfork\mathcal{Q}_2$. To show that $F^{[\varepsilon]}\pitchfork \mathcal{Q}_{3}$ is immediate, since at each point it holds that the matrix $dF^{dist}$ is given by $(1, 0, 1, 0)$.

Let us show that $F^{[\varepsilon]}\pitchfork \mathcal{Q}_{1, 2}$. It will be enough to show that for any $p\in(F^{[\varepsilon]})^{-1}(\mathcal{Q}_{1, 2})$ it holds that at least one of $a_{2, 1}$ and $a_{4, 1}$ is non-zero. Assume for contradiction, that $a_{2, 1}=a_{4, 1}=0$. We also know that $F^{[\varepsilon]}_1(p)=0$. Hence, geometrically the cylinder $\mathcal{C}_{\pi_{s_1}(p), \varepsilon}$ is tangent to the circle $\Gamma_\varepsilon\vert_{\pi_{s_2}(p)}$. Then by Lemma \ref{lemma_aux_diag} $(i.)$ we obtain that $\varepsilon<\widetilde{d}(\pi_{s_1}(p), \pi_{s_2}(p))<\delta_K$. Which is a contradiction, since $\widetilde{d}(\pi_{s_1}(p), \pi_{s_2}(p))$ is equal to $\varepsilon$ or $\delta_K$.

Similarly, we can show that $F^{[\varepsilon]}\pitchfork\mathcal{Q}_{2, 3}$. To show that $F^{[\varepsilon]}\pitchfork \mathcal{Q}$ is immediate.

It remains to show that $F^{[\varepsilon]}\pitchfork\mathcal{Q}_{1, 2}$ and $F^{[\varepsilon]}\pitchfork\mathcal{Q}_{1, 2, 3}$. Here we \textit{conjecture} that 
\begin{equation}\label{eqn_assump_rank}
\hbox{rank}\left(\begin{pmatrix}
a_{2, 1} & a_{2, 2} \\
a_{4, 1} & a_{4, 2}
\end{pmatrix}\right)=2
\end{equation}
at any point $p\in (F^{[\varepsilon]})^{-1}(\mathcal{Q}_{1, 2}\cup\mathcal{Q}_{1, 2, 3})$. See also the discussions in Remark \ref{rem_enum_alg} and Example \ref{example_diag}. 

Hence, with the assumption $(\ref{eqn_assump_rank})$, we verified that $F^{[\varepsilon]}$ is stratum transverse to $Q$. In fact, with the assumption $(\ref{eqn_assump_rank})$, we also showed that $\pi_{s_1}:\widehat{M}_{K, \varepsilon}\rightarrow \R/T\mathbb{Z}$ is a submersion on each stratum. Moreover, $\pi_{s_1}$ is proper and surjective. So we can apply Ehresmann's Theorem \ref{thm_ehrsm}, which finishes the part $(i.)$.
\newline
\newline
\noindent{\textit{Ad $(ii.)$:}} 

Let us assume $(\ref{eqn_assump_rank})$. 

Then, in fact, we showed already in part $(i.)$ that also $\pi_{s_1}^{-1}(\overline{s}_1)$ is a $3$-manifold with corners.

We are going to inspect the critical points of $\pi_{s_2}:\pi_{s_1}^{-1}(\overline{s}_1)\rightarrow (\overline{s}_2-\delta_K, \overline{s}_2-\varepsilon)\cup(\overline{s}_2+\varepsilon, \overline{s}_2+\delta_K)$.

Let $p\in (\widetilde{F}^{[\varepsilon]})^{-1}(\mathcal{Q}_2)$. Similar to the proof of Lemma \ref{lemma_crit}, $p$ is a critical point of $\pi_{s_2}$ iff
$$a_{2, 2}=0\hbox{ and }a_{4, 2}=0.$$
By the definition of $p$ it holds that $F^{[\varepsilon]}_2\vert_{\pi_{s_1}^{-1}(\overline{s}_1)}(p)=0$. Hence geometrically,  the cylinder $\mathcal{C}_{\pi_{s_2}(p), \varepsilon}$ is tangent to the circle $\Gamma_\varepsilon\vert_{\overline{s}_1}$. Then by Lemma \ref{lemma_aux_diag} $(ii.)$ (see Figure \ref{figure_numeric_diagonal} (bottom)) we obtain that $F^{[\varepsilon]}_1(p)<0,$ which is a contradiction, since $F^{[\varepsilon]}_1(p)<0.$ So there is no critical point on $(\widetilde{F}^{[\varepsilon]})^{-1}(\mathcal{Q}_2)$.

Let $p\in (\widetilde{F}^{[\varepsilon]})^{-1}(\mathcal{Q}_1)$. Then $p$ is a critical point iff
$$a_{2, 1}=0\hbox{ and }a_{4, 1}=0.$$
By the definition of $p$ it holds that $F^{[\varepsilon]}_1\vert_{\pi_{s_1}^{-1}(\overline{s}_1)}(p)=0$. Hence geometrically,  the cylinder $\mathcal{C}_{\overline{s}_1, \varepsilon}$ is tangent to the circle $\Gamma_\varepsilon\vert_{\pi_{s_2}(p)}$. Then by Lemma \ref{lemma_aux_diag} $(i.)$ there are precisely two such critical points.

Let $p=(\theta_1, \overline{s}_2, \theta_2)$ be one of these critical points of $\pi_{s_2}$. We would like to compute its Morse index. Let us consider an implicit equation
$$F^{[\varepsilon]}_1\vert_{\pi_{s_1}^{-1}(\overline{s}_1)}(\theta_1, s_2, \theta_2)=0.$$
Since $\widetilde{F}^{[\varepsilon]}\pitchfork\mathcal{Q}_1$, it holds that $\partial_{s_2}F^{[\varepsilon]}_1\vert_{\pi_{s_1}^{-1}(\overline{s}_1)}(p)\neq0$.

So we can, similarly to the proof of Lemma \ref{lemma_crit}, apply the Implicit function theorem and describe the Hessian matrix $\widetilde{H}[\pi_{s_2}](p)$ as
\begin{equation*}
\frac{-1}{\partial_{s_2}F^{[\varepsilon]}_2\vert_{\pi_{s_1}^{-1}(\overline{s}_1)}(p)}\begin{pmatrix}
\partial_{\theta_1, \theta_1}F^{[\varepsilon]}_2\vert_{\pi_{s_1}^{-1}(\overline{s}_1)}(p) & \partial_{\theta_1, \theta_2}F^{[\varepsilon]}_2\vert_{\pi_{s_1}^{-1}(\overline{s}_1)}(p)\\
\partial_{\theta_2, \theta_1}F^{[\varepsilon]}_2\vert_{\pi_{s_1}^{-1}(\overline{s}_1)}(p) & \partial_{\theta_2, \theta_2}F^{[\varepsilon]}_2\vert_{\pi_{s_1}^{-1}(\overline{s}_1)}(p)\\
\end{pmatrix}.
\end{equation*}
In order to compute the sign of the $\det(\widetilde{H}[\pi_{s_2}](p))$, we can ignore the fraction before $\widetilde{H}[\pi_{s_2}](p)$. Hence, we are left with the matrix
\begin{equation*}
\begin{pmatrix}
-\langle P+\varepsilon v_2, v_1\rangle & \varepsilon\langle v_1^{\bot}, v_2^{\bot}\rangle\\
\varepsilon\langle v_2^{\bot}, v_1^{\bot}\rangle & -\varepsilon\langle v_1, v_2\rangle\\
\end{pmatrix}.
\end{equation*}

Next, since $F^{[\varepsilon]}_1\vert_{\pi_{s_1}^{-1}(\overline{s}_1)}(p)=0$, we obtain that $-\langle P+\varepsilon v_2, v_1\rangle=-\varepsilon\langle v_1, v_1\rangle=-\varepsilon<0$. Finally, from the geometric picture of Lemma \ref{lemma_aux_diag} $(ii.)$, we see that $-\varepsilon\langle v_1, v_2\rangle>0$.

Hence all together, $\det(\widetilde{H}[\pi_{s_2}](p))<0$. And thus, the Morse index of $\pi_{s_2}$ at $p$ is equal to $1$.

Now, it is straightforward that there are no other critical points of $\pi_{s_2}$ on $(\widetilde{F}^{[\varepsilon]})^{-1}(\mathcal{Q}\cup\mathcal{Q}_{1, 2})$. Consequently, $\pi_{s_2}:\pi_{s_1}^{-1}(\overline{s}_1)\rightarrow(\overline{s}_2-\delta_K, \overline{s}_2-\varepsilon)\cup(\overline{s}_2+\varepsilon, \overline{s}_2+\delta_K)$ induce a broken fibration.
\end{pproof}

\begin{rem}\label{rem_enum_alg} The reason why we did not prove fully the Conjecture \ref{conj_diag} is that it is quite difficult to solve
$$F^{[\varepsilon]}_1=0\wedge F^{[\varepsilon]}_2=0$$
on the set $\lbrace(s_1, \theta_1, s_2, \theta_2)\in(\R/T\mathbb{Z}\times S^1)^2\,|\,\varepsilon\leq\widetilde{d}(s_1, s_2)\leq\delta_K\rbrace$. Unfortunately, now we cannot shrink $\varepsilon\rightarrow 0$ and apply the same trick as in Lemma \ref{lemma_standard_square}.

Let us aim for a weaker result - to show that close to the diagonal $\Delta_{\varepsilon}$ the system $F^{[\varepsilon]}_1=0\wedge F^{[\varepsilon]}_2=0$ has exactly $4$ real solutions. However, after tangent half-angle substitutions, we will need to find the intersections of two biquadratics. Such a problem is expected to have eight complex solutions and does not have a solution in radicals. Moreover, one can show that the geometric description of $F^{[\varepsilon]}_1=0\wedge F^{[\varepsilon]}_2=0$ leads to two biquadratics. In other words, independently of the choice of parametrization, we have to study intersections of two biquadratics. Still, some techniques from Enumerative real algebraic geometry might lead to a proof of the existence of four distinct real solutions. See also Example \ref{example_diag}.
\end{rem}

\newpage
\section{Closer look on the $\varepsilon$-diagonal}
\label{sec:closer_diag}
In this section, we inspect the singular behavior of $M_{K, \varepsilon}$ near the diagonal $\Delta_\varepsilon$. In other words, how, or even whether, the outward-pointing chords can approach the constant chords. This turns out to be a local question about the geometry of $T_{K, \varepsilon}$ near points with non-positive Gaussian curvature. After slightly laborious computation, we conclude that for $\varepsilon>0$ small, $M_{K, \varepsilon}$ looks near $\Delta_\varepsilon$ roughly as a stratified fiber bundle over a half-torus with cuspidal fibers.  

\begin{rem} We assume that the reader is familiar with basic notions from the geometry of surfaces, such as the first and the second fundamental forms, principal directions, and Gaussian curvature. For more on this theory, see, for example \cite{Carmo1976DifferentialGO}.
\end{rem}

\begin{conv}\label{conv_princip} Let $M\subset\R^3$ be a closed oriented surface. Then there is the unique compact $3$-submanifold $\widetilde{M}\subset\R^3$ such that $\partial\widetilde{M}=M$. This gives us a canonical outward-pointing normal vector field $N^{out}$ on $M$. Let $p\in M$ be such that the principal curvatures $\kappa_M(p)\geq\kappa_m(p)$ are not equal. Then we choose the (normalized) principal directions $w_M(p), w_m(p)$ such that
$$w_M(p)\times w_m(p)=N^{out}_p.$$ 
\end{conv}

\begin{lemma}\label{lemma_curvature} The Gaussian curvature $\kappa_G$ of $T_{K, \varepsilon}$ at $p=(s, \theta)$ is given by

$$\kappa_G(p)=\frac{-\kappa(s)\cos(\theta)}{\varepsilon d},$$
where $d=(1-\varepsilon\cos(\theta)\kappa(s)).$

Also, the principal curvatures $\kappa_M, \kappa_m$ and the (normalized) principal directions $w_M, w_m$ are given by
$$\kappa_M(p)=\frac{\kappa(s)\cos(\theta)}{d}\hbox{ and }w_M(p)=-d^{-1}(\partial_s-\tau(s)\partial_\theta),$$
and
$$\kappa_m(p)=-\varepsilon^{-1}\hbox{ and }w_m(p)=\varepsilon^{-1}\partial_\theta.$$
\end{lemma}

\begin{proof}
The strategy of the computation is the following. First, we need to know the first and the second fundamental forms with respect to the basis $\lbrace\partial_s, \partial_\theta\rbrace$. We denote the first and second fundamental forms as $G$ and $H$, respectively. Then, we obtain the Gaussian curvature and the mean curvature $\kappa_{mean}$ at $p\in T_K$ as
$$\kappa_G=\frac{\det(H)}{\det(G)}\hbox{ and }\kappa_{mean}=\frac{h_{11}g_{22}-2 h_{12}g_{12}+h_{22}g_{11}}{2\det(G)}.$$
Then, the principal curvatures are determined as
$$\kappa_1=\kappa_{mean}+\sqrt{\kappa_{mean}^2-\kappa_G}\hbox{ and }\kappa_2=\kappa_{mean}-\sqrt{\kappa_{mean}^2-\kappa_G}.$$
From which we finally compute the principal directions as the eigenvectors of $G^{-1}H$.

Let us compute $G$. We recall relations $(\ref{eqn_torus_aux_deriv})$, i.e.
\begin{align*}
\partial_s\Gamma_\varepsilon&=d\dot{\gamma}(s)+\varepsilon\tau(s)v^\bot,\\
\partial_\theta\Gamma_\varepsilon&=\varepsilon v^\bot.
\end{align*}
Thus, we obtain that
\begin{equation}\label{eqn_metric_torus}
G=\begin{pmatrix}
d^2+\varepsilon^2 \tau^2(s) & \varepsilon^2 \tau(s)\\
\varepsilon^2 \tau(s) & \varepsilon^2
\end{pmatrix}.
\end{equation}
Now, we continue with the computation of $H$. It holds that
\begin{align*}
\partial_{s}v^\bot&=\partial_s\big(-\sin(\theta)n(s)+\cos(\theta)b(s)\big)\\
&=-\sin(\theta)\big(-\kappa(s)\dot{\gamma}(s)+\tau(s)b(s)\big)+\cos(\theta)\big(-\tau(s)n(s)\big)\\
&=\sin(\theta)\kappa(s)\dot{\gamma}(s)-\tau(s) v
\end{align*}
Hence, we compute
\begin{align}
\begin{split}\label{eqn_second_deriv}
\partial_s^2\Gamma_\varepsilon &=\big(-\varepsilon\cos(\theta)\dot{\kappa}(s)+\varepsilon\tau(s)\sin(\theta)\kappa(s)\big)\dot{\gamma}(s)+d\cdot\kappa(s) n(s)+\varepsilon\dot{\tau}(s)v^\bot-\varepsilon\tau^2(s)v,\\
\partial_s\partial_\theta\Gamma_\varepsilon &=\varepsilon\sin(\theta)\kappa(s)\dot{\gamma}(s)-\varepsilon\tau(s)v,\\
\partial_\theta^2\Gamma_\varepsilon &=-\varepsilon v.
\end{split}
\end{align}
By Remark \ref{rem_normal_vect} we know that $v$ are the unital normal vectors to $T_{K, \varepsilon}$ and thus the second fundamental form is given as the dot product of the second derivatives from (\ref{eqn_second_deriv}) with $v$. Therefore, $H$ is given as
$$\begin{pmatrix}
\kappa(s)\cos(\theta)d-\varepsilon\tau^2(s) & -\varepsilon\tau(s)\\
-\varepsilon\tau(s) & -\varepsilon
\end{pmatrix}.$$
Then $\kappa_{mean}(p)=\frac{1-2\varepsilon\kappa(s)\cos(\theta)}{-2\varepsilon d}$ and the rest of the computation follows.
\end{proof}

\begin{rem}$\kappa_G$ is negative for $\theta\in(-\pi/2, \pi/2)$.
\end{rem}

\begin{rem}$w_M=-\dot{\gamma}(s)$ and $w_m=v^\bot$. And in particular, it holds that $w_M\times w_m=v$.
\end{rem}

\begin{notat}\label{not_deriv} Let $f:\R^2_{x, y}\rightarrow\R$ be a smooth function and $f(x, y)=0$ an implicit equation which is locally solved by the function $y(x)$. Then, we will use the following notation for the partial and implicit derivatives of $f$ at the point $(x, y(x))$. That is
$$\partial_x f(x, y(x)):=\partial_x f(x, y)\vert_{(x, y)=(x, y(x))}$$
and
$$\dot{f}(x, y(x)):=\frac{d}{dx}f(x, y(x)).$$
\end{notat}

\begin{defn} Let $w_M, w_m$ be two orthogonal vectors in $\R^2$. Let $Q$ be a quadrant that is spanned by axes $\langle w_M\rangle, \langle w_m\rangle$ and let $c:[0, T]\rightarrow Q$ be a smooth curve such that $c(0)=(0, 0)$ and $c((0, T])\subset Int(Q).$ $c$ is called \textbf{$w_M$-attracted} if there is $\delta; 0<\delta<T,$ such that $c((0, \delta])$ is strictly concave graph over the axis $\langle w_M\rangle$ (here the positive direction of the functions over $\langle w_M\rangle$ is given by the half-axis $\langle w_m\rangle\cap Q$).

Analogously for $w_m$-attracted curves.
\end{defn}

\begin{defn_lemma}\label{lemma_convexity}Let $M$ be a closed oriented surface and let $p\in M$ such that $\kappa_G(p)<0$. Let $w_M, w_m$ be (normalized) principal directions $w_M, w_m$ of principal curvatures $\kappa_M>\kappa_m$ (recall Convention \ref{conv_princip}). 
We also consider the two planes spanned by the normal vector to $M$ at $p$ together with $w_M$ or $w_m$. Then, these planes split a neighbourhood of $p$ in $M$ into four sectors $\lbrace Q_j\rbrace_{j\in\lbrace1,\dots, 4\rbrace}$, see Figure \ref{figure_diagonal_convexity}. Let also $c_j$ be the unique curves in $Q_j$ given by the intersection of the affine plane $T_pM$ and $M$ and defined in a small neighbourhood of $p$.

Now, we consider $Q_j$ for any $j$. If the curve $c_j$ is $w_M$-attracted in $T_pM$ for some $\delta>0$, then for $t\in[0, \delta]$ the chords $\overline{p, c_j(t)}$ are outward-pointing in the sense of Definition \ref{defn_out_out} and there is a smooth curve $\widehat{c}_j:[0, \delta]\rightarrow Q_j$ such that
\begin{itemize}
\item[$(i.)$] $\widehat{c}_j(0)=p$,
\item[$(ii.)$] $c_j((0, \delta])\cap\widehat{c}_j((0, \delta])=\emptyset$,
\item[$(iii.)$] the chords $\overline{p, \widehat{c}_j(t)}$ are outward-pointing.
\end{itemize}
See also Figure \ref{figure_diagonal_convexity}.

The analogous statement holds in the $w_m$-attracted case, only then the chords need to be inward-pointing on both ends. We will call such a quadrant $Q_j$ \textbf{outward} or \textbf{inward pointing}.

Let $\widetilde{Q_j}$ be the quadrants of $T_pM$ that naturally correspond to the quadrants $Q_j$; see also Figure \ref{figure_diagonal_convexity}. Then $\rho^j_a, \rho^j_b$ be the following functions on $\widetilde{Q_j}$. For $j=1, 3$, $\rho^j_a$ is $+1$ and else $-1$. For $j=1, 2$, $\rho^j_b$ is $+1$ and else $-1$.

If at the point $p$ holds that
\begin{equation}\label{eqn_attracted}
\rho^j_b\left(2\rho^j_a\nabla_{w_M}(\kappa_M)+2\nabla_{w_m}(\kappa_m)\left\lvert\frac{\kappa_M}{\kappa_m} \right\rvert^{3/2}+6\nabla_{w_m}(\kappa_M)\left\lvert\frac{\kappa_M}{\kappa_m} \right\rvert^{1/2}-6\rho^j_a\nabla_{w_M}(\kappa_m)\frac{\kappa_M}{\kappa_m}\right)<0,
\end{equation}
then $Q_j$ is outward-pointing. If the inequality is opposite, $Q_j$ is inward-pointing.
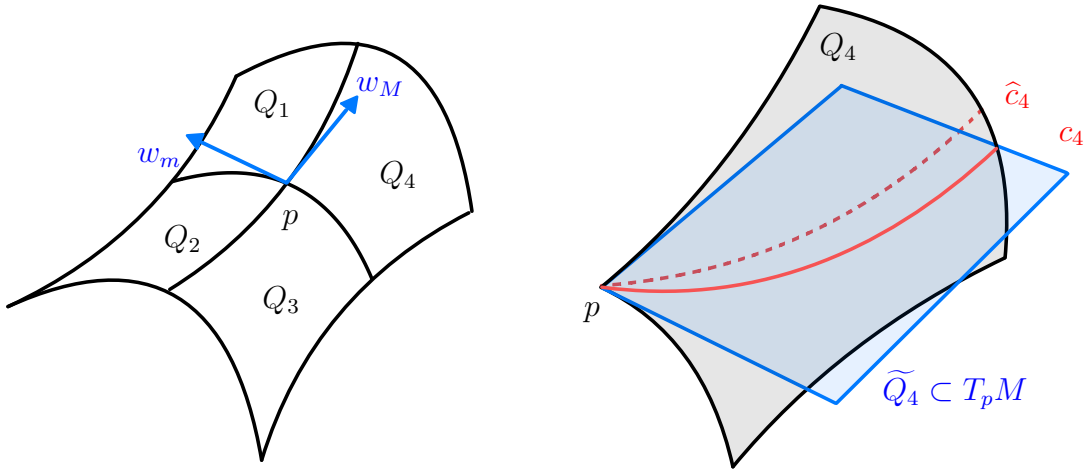
\begin{figure}[!htbp]
\labellist
\pinlabel $Q_1$ at 150 710
\pinlabel $Q_2$ at 100 634
\pinlabel $Q_3$ at 155 600
\pinlabel $Q_4$ at 220 670
\pinlabel $\textcolor{blue}{w_M}$ at 210 720
\pinlabel $\textcolor{blue}{w_m}$ at 87 680
\pinlabel $p$ at 160 645
\pinlabel $p$ at 328 593
\pinlabel $Q_4$ at 465 742
\pinlabel $\textcolor{blue}{\widetilde{Q_4}\subset T_pM}$ at 530 550
\pinlabel $\textcolor{red}{c_4}$ at 595 692
\pinlabel $\textcolor{red}{\widehat{c}_4}$ at 565 715
\endlabellist
\centering
\includegraphics[scale=0.68]{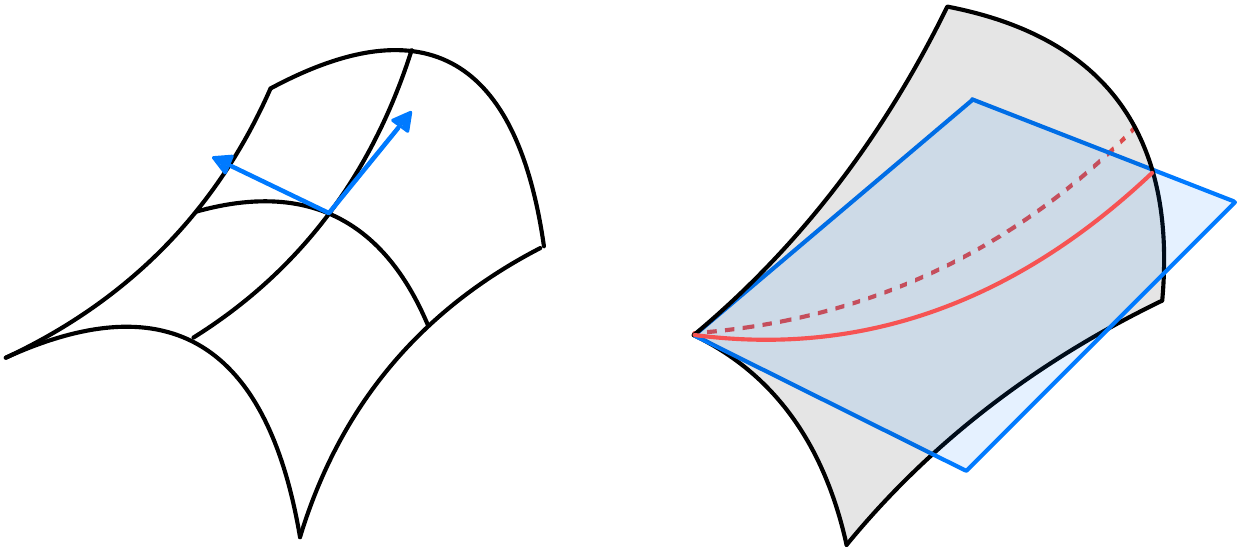}
\vspace{0.3cm}
\caption{On the left: the splitting of the neighbourhood of $p$ in $M$ onto quadrants $Q_j$ using principal directions $w_M,w_m$. By our convention, $w_M\times w_m$ is equal to the outward-pointing normal vector. On the right: $w_M$-attracted curve $c_4$ together with possible $\widehat{c}_4$. Observe also the quadrant $\widetilde{Q_4}$ corresponding to $Q_4$.}
\label{figure_diagonal_convexity}
\end{figure}
\end{defn_lemma}

\begin{n_example} Let $M$ be a hyperbolic paraboloid. Then still $\kappa_G(p)<0$. However, since $M$ is a ruled surface, the curves $c_j$ are neither $w_M$- nor $w_m$-attracted. In particular, there are no curves $\widehat{c}_j$ that satisfy the conditions $(i.)-(iii.)$ from Lemma \ref{lemma_convexity}.
\end{n_example}

\begin{proof}
Let $c_j$ is $w_M$-attracted, then for each $t\in(0, \delta]$ chords $\overline{p, c_j(t)}$ are outward-pointing. Indeed, $\overline{p, c_j(t)}$ lie in $T_pM$, so they are tangent to $M$ at $p$. At $c_j(t)$ they are also outward-pointing, since $c_j$ is $w_M$-attracted and $\kappa_M>\kappa_m$, see Figure \ref{figure_diagonal_chord}. Now, let us we make a small $C^0$ perturbation $\widetilde{c}_j$ of $c_j$ inside $\widetilde{Q}_j$ such that $c_j(0)=\widetilde{c}_j$. If the perturbation is small enough then the lift ($=:\widehat{c}_j$) of $\widetilde{c}_j$ to $Q_j$ has a property that the chords $\overline{p, \widehat{c}_j(t)}$ are outward-pointing in $\widehat{c}_j(t)$. If we perturb $c_j$ such that $\widetilde{c}_j\subset \widetilde{Q}_j$ lies in the convex hull given by $c_j$ and $\langle w_M\rangle\cap\widetilde{Q}_j$, the the lift $\widehat{c}_j$ satisfies the property $(iii.)$. In Lemma \ref{cor_convexity_tangent} we will construct a special example of smooth $\widehat{c_j}$ such that the chords $\overline{p, \widehat{c}_j(t)}$ are tangent to $M$ at $\widehat{c}_j(t)$.

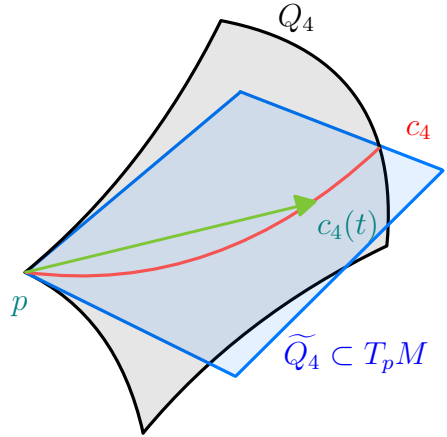
\begin{figure}[!htbp]
\labellist
\pinlabel $\textcolor{teal}{p}$ at 0 660
\pinlabel $Q_4$ at 172 840
\pinlabel $\textcolor{blue}{\widetilde{Q_4}\subset T_pM}$ at 209 630
\pinlabel $\textcolor{red}{c_4}$ at 248 770
\pinlabel $\textcolor{teal}{c_4(t)}$ at 204 710
\endlabellist
\centering
\includegraphics[scale=0.60]{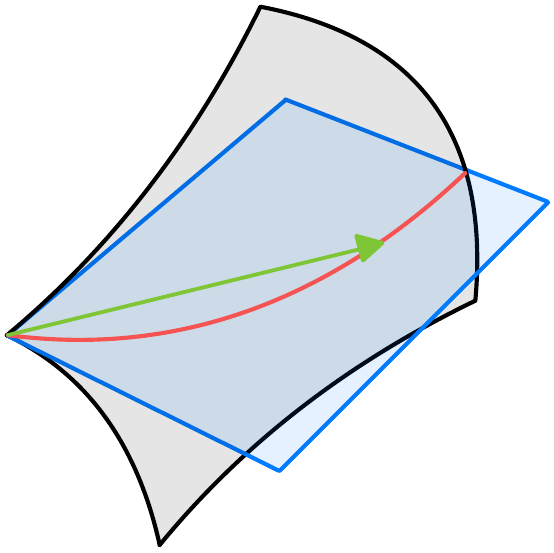}
\vspace{0.3cm}
\caption{In green a chord $\overline{p, c_4(t)}$ that is outward-point.}
\label{figure_diagonal_chord}
\end{figure}

Now, the inequality (\ref{eqn_attracted}). After rigid transformations of $\R^3$ we can assume that $M$ is locally a graph of a function $f:\R_{x, y}^2\rightarrow\R$ such that
\begin{itemize}
\item[$(i.)$] $(q, f(q))=(0, 0, 0)=p$ for $q=(0, 0)$.
\item[$(ii.)$] $\partial_xf(q)=\partial_yf(q)=0.$
\item[$(iii.)$] $\partial^2_xf(q)>0$, $\partial^2_yf(q)<0$ and $\partial_x\partial_yf(q)=0$.
\end{itemize}

We are going to study the set $f(x, y)=0$. By the Morse Lemma \cite{hirsch_1997} there is a coordinate change $\varphi:\R^2_{x, y}\rightarrow \R^2_{\widetilde{x}, \widetilde{y}}$, which is defined in a small neighbourhood of $q$ and satisfy that
$$d\varphi(q)=\mathbbm{1}\hbox{ and }f(\widetilde{x}, \widetilde{y})=[\partial^2_xf(q)]\cdot\widetilde{x}^2+[\partial^2_yf(q)]\cdot\widetilde{y}^2.$$
Since $q$ is a nondegenerate critical point, $f(x, y)=0$ is diffeomorphic to a cone, and each of the two smooth components can be written as a graph of $y(x)$. In particular, $(x, y(x))$ will lie inside $Q_1$ and $Q_3$ or $Q_2$ and $Q_4$. Note that we can see immediately the slope $\dot{y}(0)$ from the Morse Lemma, but for further calculations, it will be illustrative to make the following computation
\begin{align*}
\dot{y}(0)&\stackrel{(1)}{=}\lim_{x\rightarrow 0}-\frac{\partial_x f(x, y(x))}{\partial_y f(x, y(x))}\\
&\stackrel{(2)}{=}\lim_{x\rightarrow 0}-\frac{(\dot{\partial_x f})(x, y(x))}{(\dot{\partial_y f})(x, y(x))}\\
&\stackrel{(3)}{=}-\frac{\partial^2_xf(q)+\partial_x\partial_yf(q)\cdot\dot{y}(0)}{\partial_x\partial_yf(q)+\partial^2_yf(q)\cdot\dot{y}(0)}\\
&=-\frac{\partial^2_xf(q)}{\partial^2_yf(q)\cdot\dot{y}(0)},
\end{align*}
where we used in step $(1)$ implicit differentiation and in step $(2)$ the L$^\prime$ Hospital rule. We applied in step $(3)$ implicit differentiation and in step the assumptions on $f$. We also recall that we used Notation \ref{not_deriv} for various differentiating. Hence, as expected, we obtain that
\begin{equation}\label{eqn_aux_attract}
\lim_{x\rightarrow 0}\frac{\partial_x f(x, y(x))}{\partial_y f(x, y(x))}=-\rho_a\left\lvert\frac{\partial^2_xf(q)}{\partial^2_yf(q)} \right\rvert^{1/2},
\end{equation}
where $\rho_a=\hbox{sgn}(\dot{y}(0))$.

In the following, we write for simplicity $f(x, y(x))$ as $f$. Also, we would like to denote partial derivatives as subscripts, for example, $\partial_x\partial_yf=f_{xy}$.

We are interested in the sign of $\ddot{y}(0)$, so we can compute
\begin{align*}
\hspace*{-1cm}\ddot{y}(0)&=\lim_{x\rightarrow 0}\ddot{y}(x)\\
&\stackrel{(1)}{=}\lim_{x\rightarrow 0}\frac{-\dot{(f_{x})}f_y+f_x\dot{(f_{y})}}{f_y^2}\\
&\stackrel{(2)}{=}\lim_{x\rightarrow 0}\frac{-f_{xx}f_y^2+2f_xf_yf_{xy}-f_{yy}f_x^2}{f_y^3}\\
&\stackrel{(3)}{=}\lim_{x\rightarrow 0}\frac{\splitfrac{f_x\left[2f_{xxx}f_{yy}+f_{yyy}\left[3f_{xx}\frac{f_x}{f_y}+f_{yy}\left(\frac{f_x}{f_y}\right)^3\right]+3f_{xxy}\left[f_{xx}\left(\frac{f_x}{f_y}\right)^{-1}-f_{yy}\frac{f_x}{f_y}\right]-6f_{xyy}f_{xx}\right]}{+O\big((f_x+f_y)(f_x+f_y+f_{xy})\big)}}{3f_y\big(f_{yy}^2\left(\frac{f_x}{f_y}\right)^2-f_{yy}f_{xx}\big)+O\big((f_x+f_y)(f_x+f_y+f_{xy})\big)}\\
&\stackrel{(4)}{=}\lim_{x\rightarrow 0}\frac{\frac{f_x}{f_y}\left[2f_{xxx}f_{yy}+f_{yyy}\left[3f_{xx}\frac{f_x}{f_y}+f_{yy}\left(\frac{f_x}{f_y}\right)^3\right]+3f_{xxy}\left[f_{xx}\left(\frac{f_x}{f_y}\right)^{-1}-f_{yy}\frac{f_x}{f_y}\right]-6f_{xyy}f_{xx}\right]}{3\big(f_{yy}^2\left(\frac{f_x}{f_y}\right)^2-f_{yy}f_{xx}\big)}\\
&\stackrel{(5)}{=}\lim_{x\rightarrow 0}\frac{\frac{f_x}{f_y}\left[2f_{xxx}+f_{yyy}\left[3\frac{f_{xx}}{f_{yy}}\frac{f_x}{f_y}+\left(\frac{f_x}{f_y}\right)^3\right]+3f_{xxy}\left[\frac{f_{xx}}{f_{yy}}\left(\frac{f_x}{f_y}\right)^{-1}-\frac{f_x}{f_y}\right]-6f_{xyy}\frac{f_{xx}}{f_{yy}}\right]}{3\left(f_{yy}\left(\frac{f_x}{f_y}\right)^2-f_{xx}\right)}.
\end{align*}
In equality $(1)$ we used implicit differentiation (and the chain rule), and in $(2)$ we applied implicit differentiation. For $(3)$, we applied the L$^\prime$ Hospital rule twice. In $(4)$ we got rid of $O$-terms; for their closer description, see below. This was possible since $f_x, f_y, f_{xy}\rightarrow 0$ as $x\rightarrow 0$, $\lim_{x\rightarrow0}f_x/f_y\in\R\setminus\lbrace 0\rbrace$, and $f^2_{yy}f_x^2/f_y^2-f_{yy}f_{xx}>0$. And in $(5)$ we cancelled $f_{yy}$.

For completeness, we write the description of the numerator and the denominator in step $(3)$. So the numerator is of the form
\begin{align*}
&f_x\left[2f_{xxx}f_{yy}+f_{yyy}\left[3f_{xx}\frac{f_x}{f_y}+f_{yy}\left(\frac{f_x}{f_y}\right)^3\right]+3f_{xxy}\left[f_{xx}\left(\frac{f_x}{f_y}\right)^{-1}-f_{yy}\frac{f_x}{f_y}\right]-6f_{xyy}f_{xx}\right]\\
&-\frac{f_{yyyy}f_x^4}{f_y^2}-\frac{f_{yyy}f_x^3f_{xy}}{f_y^2}+\frac{6f_{xyy}f_x^2f_{xy}}{f_y}+\frac{4 f_{xyyy}^3f_x^3}{f_y}-6f_x^2f_{xxyy}\\
&-f_yf_{xy}f_{xxx}+4f_yf_x f_{xxxy}-f_y^2f_{xxxx}
\end{align*}
and the denominator is equal to
\begin{align*}
&3\left(\frac{f_{yy}^2f_x^2}{f_y}-f_yf_{yy}f_{xx}\right)\\
&+3\big(f_x^2f_{yyy}-2f_xf_{xy}f_{yy}+2f_y f_{xy}^2-2f_xf_yf_{xyy}+f_y^2f_{xxy}\big).
\end{align*}

Hence, by relation (\ref{eqn_aux_attract}) and assumptions on $f$ at $q$, we obtain
\begin{align*}
\hspace*{-2cm}\hbox{sgn}(\ddot{y}(0))&=\hbox{sgn}\,\lim_{x\rightarrow 0}-\frac{f_x}{f_y}\left[2f_{xxx}+f_{yyy}\left[3\frac{f_{xx}}{f_{yy}}\frac{f_x}{f_y}+\left(\frac{f_x}{f_y}\right)^3\right]+3f_{xxy}\left[\frac{f_{xx}}{f_{yy}}\left(\frac{f_x}{f_y}\right)^{-1}-\frac{f_x}{f_y}\right]-6f_{xyy}\frac{f_{xx}}{f_{yy}}\right]\\
&=\hbox{sgn}\,\lim_{x\rightarrow 0}\rho_a\left\lvert\frac{f_{xx}}{f_{yy}} \right\rvert^{\frac{1}{2}}\Bigg[2f_{xxx}+f_{yyy}\left[-3\rho_a\frac{f_{xx}}{f_{yy}}\left\lvert\frac{f_{xx}}{f_{yy}} \right\rvert^{\frac{1}{2}}-\left\lvert\frac{f_{xx}}{f_{yy}} \right\rvert^{\frac{3}{2}}\right]\\
&\quad\quad\quad\quad\,\,\,\,+3f_{xxy}\left[-\rho_a\frac{f_{xx}}{f_{yy}}\left\lvert\frac{f_{xx}}{f_{yy}} \right\rvert^{-\frac{1}{2}}+\rho_a\left\lvert\frac{f_{xx}}{f_{yy}} \right\rvert^{\frac{1}{2}}\right]-6f_{xyy}\frac{f_{xx}}{f_{yy}}\Bigg]\\
&=\hbox{sgn}\,\lim_{x\rightarrow 0}\rho_a\left[2f_{xxx}+\rho_a 2f_{yyy}\left\lvert\frac{f_{xx}}{f_{yy}} \right\rvert^{\frac{3}{2}}+6\rho_af_{xxy}\left\lvert\frac{f_{xx}}{f_{yy}} \right\rvert^{\frac{1}{2}}-6f_{xyy}\frac{f_{xx}}{f_{yy}}\right]\\
&=\hbox{sgn}\,\lim_{x\rightarrow 0}\left[2\rho_a f_{xxx}+2f_{yyy}\left\lvert\frac{f_{xx}}{f_{yy}} \right\rvert^{\frac{3}{2}}+6f_{xxy}\left\lvert\frac{f_{xx}}{f_{yy}} \right\rvert^{\frac{1}{2}}-6\rho_a f_{xyy}\frac{f_{xx}}{f_{yy}}\right].
\end{align*}

Now, it remains to relate $f_{xxx}, f_{xxy}, f_{yyx}, f_{yyy}$ at $q$ with the directional derivatives $\nabla_{w_{M}}\kappa_{M}$, etc. It is not hard to see that at $q$ it holds $f_{yy}=\kappa_M$, $\partial_x=w_M$, $f_{yy}=\kappa_m$, $\partial_y=w_m$. However, the quantity $\nabla_{w_m}\kappa_M(q)$ depends pointwise only on $w_M$; for $\kappa_M$, we need to know the values in a neighborhood $U_q$.

Hence, let us first find the descriptions of $H$ and $G^{-1}$ around $q$ at graphical coordinates $(x, y)$.

Since
\begin{equation}\label{eqn_metr_graph}
G=\begin{pmatrix}
f_x^2+1 & f_x f_y\\
f_x f_y & f_y^2+1
\end{pmatrix},
\end{equation}
we have that
\begin{equation}\label{eqn_inv_metr_graph}
G^{-1}=\frac{1}{f_x^2+f_y^1+1}\begin{pmatrix}
f_y^2+1 & -f_x f_y\\
-f_x f_y & f_x^2+1
\end{pmatrix}.
\end{equation}

Also at any $r=(x, y, f(x, y))$ we choose a normal vector $N_r$ to the graph of $f$ as
$$N_r=\frac{(-f_x, -f_y, 1)^T}{(f_x^2+f_y^2+1)^{1/2}}.$$
Then
$$H=\frac{1}{(f_x^2+f_y^2+1)^{1/2}}\begin{pmatrix}
f_{xx} & f_{xy}\\
f_{xy} & f_{yy}
\end{pmatrix}.$$

Next,
$$\kappa_G=\frac{f_{xx}f_{yy}-f_{xy}^2}{(f_x^2+f_y^2+1)^{3}}\hbox{ and }\kappa_{mean}=\frac{f_{xx}(f_y^2+1)-2f_{xy}f_x f_y+f_{yy}(f_{xx}^2+1)}{2(f_x^2+f_y^2+1)^{3}},$$
and hence
$$\kappa_G(q)=f_{xx}f_{yy}\hbox{ and }\kappa_{mean}(q)=\frac{f_{xx}+f_{yy}}{2}.$$
We can also compute that
\begin{align*}
\partial_x \kappa_G(q)&=(f_{xxx}f_{yy}+f_{xx}f_{yyx}-2 f_{xy}f_{xxy})(f_x^2+f_y^2+1)^{-3}\\
&\quad-3(f_{xx}f_{yy}-f_{xy}^2)(f_x^2+f_y^2+1)^{-4}(2 f_x f_{xx} +2 f_y f_{xy})\\
&=f_{xxx}f_{yy}+f_{xx}f_{yyx}
\end{align*}
and similarly
$$\partial_x \kappa_{mean}(q)=\frac{f_{xxx}+f_{yyx}}{2}.$$

Now we compute
\begin{align*}
\nabla_{w_M}\kappa_M(q)&=\partial_x\kappa_M(q)\\
&=\partial_x\big[\kappa_{mean}+(\kappa_{mean}^2-\kappa_G)^{1/2}\big]\\
&=\partial_x(\kappa_{mean})+\frac{2 \kappa_{mean}\partial_x(\kappa_{mean})-\partial_x(\kappa_G)}{2(\kappa_{mean}^2-\kappa_G)^{1/2}}\\
&=\frac{f_{xxx}+f_{yyx}}{2}+\frac{\frac{(f_{xx}+f_{yy})(f_{xxx}+f_{yyx})}{2}-(f_{xxx}f_{yy}+f_{xx}f_{yyx})}{((f_{xx}+f_{yy})^2-4f_{xx}f_{yy})^{1/2}}\\
&\stackrel{(1)}{=}\frac{f_{xxx}+f_{yyx}}{2}+\frac{(f_{xx}-f_{yy})(f_{xxx}-f_{yyx})}{2(f_{xx}-f_{yy})}\\
&=f_{xxx},
\end{align*}
where in $(1)$ we used the assumption that $f_{xx}>0>f_{yy}$.

We obtain similarly that $\nabla_{w_m}\kappa_M(q)=f_{xxy}, \nabla_{w_M}\kappa_m(q)=f_{yyx}, \nabla_{w_m}\kappa_m(q)=f_{yyy}$. This concludes the lemma.
\end{proof}

\begin{lemma}\label{cor_convexity_tangent} Assume that the quadrant $Q_j\subset M$ from Lemma \ref{lemma_convexity} is outward-pointing. Then the curve $\widehat{c}_j(t)$ can be chosen such that the chords $\overline{p, \widehat{c}_j(t)}$ are tangent to $M$ at $\widehat{c}_j(t)$. Moreover, the curves $c_j$ and $\widehat{c}_j$ have a common tangent line at $p$.

The analogous statement holds for the inward-pointing case.
\end{lemma}

\begin{proof}
We will prove the outward-pointing case. As in the proof of Lemma \ref{lemma_convexity}, let us consider $M$ locally as a graph of $f(x, y)$.

Now we define two functions $F_{1, p}, F_{2, p}:M\rightarrow\R$ which are similar to the functions from Definition \ref{defn_funct_out}. More precisely, at $r=(x, y, f(x, y))$ we put
\begin{equation}\label{eqn_out_funct_simple}
F_{1, p}(r):=\langle r-p, N_p\rangle\hbox{ and }F_{2, p}(r):=\langle r-p, N_r\rangle.
\end{equation}
Recall that $p=(0, 0, 0)$ and $N_r=(f_x^2+f_y^2+1)^{-1/2}(-f_x, -f_y, 1)^T$. Note that now, unlike in Definition \ref{defn_funct_out}, $F_{1, p}$ and $F_{2, p}$ have fixed variables in the ``start-point.''

Hence, in graphical coordinates
$$F_{1, p}(x, y)=f\hbox{ and }F_{2, p}(x, y):=\frac{-x f_x-y f_y+f}{(f_x^2+f_y^2+1)^{1/2}}.$$

As in the proof of Lemma \ref{lemma_convexity}, we can write around $q=(0, 0)$ $y$ as the unique function of $x$ constrained to $F_{1, p}=0\wedge \pi_{x, y}(Q_j)$. Such a function gives us $\lbrace q^1_t\rbrace$ and in particular
$$\dot{y}(0)=-\rho_a^j\left\lvert\frac{\partial^2_xf(q)}{\partial^2_yf(q)} \right\rvert^{1/2}.$$

Now, in $\pi_{x, y}(Q_j)$, we would like to write similarly $y$ as the locally unique function of $x$ constrained to $F_{2, p}=0$. Also along the set $\lbrace y=0\rbrace\cap\pi_{x, y}(Q_j)$ it holds that $F_{1, p}\leq 0$. Then, by the Intermediate value theorem and $w_M$-attraction of $c_j$, $y(x)$ will contribute by a family of chords that are not only outward pointing in $(x, y(x))$, but also in $p$. For this, it will be sufficient to show that at $p$ the functions $F_{1, p}$ and $F_{2, p}$ have the same value, the same Jacobi matrix, and, up to a multiple of $-1$, the same Hessian matrix. Note that this will, in particular, also imply that the curves $c_j$ and $\widehat{c}_j$ have a common tangent line at $p$.

First, by the choice of $f$, it holds clearly that $F_{1, p}(q)=0=F_{2, p}(q)$. Then, let us compute
\begin{align*}
\partial_x F_{2, p}(x, y)&=-F_{2, p}(x, y)\frac{f_x f_{xx}+f_y f_{xy}}{f_x^2+f_y^2+1}-\frac{xf_{xx}+y f_{xy}}{(f_x^2+f_y^2+1)^{1/2}},\\
\partial_y F_{2, p}(x, y)&=-F_{2, p}(x, y)\frac{f_y f_{yy}+f_x f_{xy}}{f_x^2+f_y^2+1}-\frac{yf_{yy}+x f_{xy}}{(f_x^2+f_y^2+1)^{1/2}}.
\end{align*}
Hence both Jacobi matrices vanish. Let us continue with
\begin{align*}
\partial^2_x F_{2, p}(x, y)&=-\partial_x F_{2, p}(x, y)g_1(x, y)- F_{2, p}(x, y) g_2(x, y)\\
&\quad\quad-\frac{\textcolor{blue}{f_{xx}}+xf_{xxx}-yf_{xxy}}{(f_x^2+f_y^2+1)^{1/2}}-\frac{(f_x f_{xx}+f_{y}f_{yy})^2}{(f_x^2+f_y^2+1)^{3/2}},\\
\partial_x\partial_y F_{2, p}(x, y)&=-\partial_x F_{2, p}(x, y)g_3(x, y)- F_{2, p}(x, y) g_4(x, y)\\
&\quad\quad-\frac{y f_{xyy}++f_{xy}+xf_{xxy}}{(f_x^2+f_y^2+1)^{1/2}}+\frac{(yf_{yy}+xf_{xy})(f_xf_{xx}+f_yf_{xy})}{(f_x^2+f_y^2+1)^{3/2}},\\
\partial^2_y F_{2, p}(x, y)&=-\partial_y F_{2, p}(x, y)g_5(x, y)- F_{2, p}(x, y) g_6(x, y)\\
&\quad\quad-\frac{\textcolor{blue}{f_{yy}}+yf_{yyy}-xf_{xyy}}{(f_x^2+f_y^2+1)^{1/2}}-\frac{(f_x f_{xx}+f_{y}f_{yy})^2}{(f_x^2+f_y^2+1)^{3/2}},
\end{align*}
where $\lbrace g_i\rbrace_{i=1,\dots, 6}$ are some functions that are bounded in a small neighborhood of $q$. Note that the \textcolor{blue}{blue} terms are the only nonzero contributions to the Hessian matrix of $F_{2, p}$ at $q$. In particular, $H[F_{1, p}](q)=-H[F_{2, p}](q)$, and the lemma follows.
\end{proof}

\begin{lemma}\label{lemma_aux_differential}Let $M$ and $p$ be as in Lemma \ref{lemma_convexity}. For $i=1, 2$ let $F_{i, p}$ be the functions from (\ref{eqn_out_funct_simple}). Then there is a small collar neighborhood $U^{\circ}_p$ of $p$ in $M$ such that
$$\hbox{rank}(dF_{i, p}(r))=1$$
for any $r\in U^{\circ}_p$ such that $F_{i, p}(r)=0.$
\end{lemma}
\begin{proof}
We consider the graphical coordinates $(x, y, f(x, y))$ as in Lemma \ref{lemma_convexity}. Then for each quadrant $Q_j$ we can find the locally unique solution $(x, y(x))$ of $F_{i, p}(x, y(x))=0$. We know that $\partial_x F_{i, p}(q)=0$, where $q=(0, 0)$. Then one can show that the implicit derivative $\dot{(\partial_x F_{i, p})}(q)\neq 0$, since $\kappa_G(p)\neq 0$. So on a small neighborhood of $q$ it holds that $\partial_x F_{i, p}(q)$ is strictly increasing or decreasing along $(x, y(x))$.
\end{proof}

\begin{ThomIsotop}\cite{goresky2012stratified, Mather2012}\label{thm_thom_isot} Let $Z$ be a Whitney stratified subset of a smooth manifold $M$ and $\pi:M\rightarrow N$ be a smooth map onto a smooth manifold $N$. Assume that $\pi\vert_Z$ is a proper surjection and $\pi$ is a submersion on each stratum of $Z$. Then $\pi\vert_Z$ defines a locally trivial \textbf{stratified fiber bundle}, i.e., for every $x\in N$ there is an open neighborhood $U_x$ in $N$ and a strata preserving homeomorphism $\psi: U_x\times\pi\vert_Z^{-1}(x)\rightarrow\pi\vert_Z^{-1}(U_x)$ which is smooth on each stratum and $\pi\vert_Z\circ\psi(u, v)=u$ for every $(u, v)\in U_x\times\pi\vert_Z^{-1}(x)$.
\end{ThomIsotop}

\begin{lemma}\label{conj_local_fibr_diag} Let $(M, g)$ be an oriented surface in $(\R^3, g_{Euc})$ with induced metric and $p\in M$ with $\kappa_G(p)<0$. Let us consider a neighborhood $U_p$ of $p$ in $M$ and the diagonal $\Delta_{U_p}\subset U_p\times U_p$. Next, $\overline{\nu_\delta}(\Delta_{U_p})\subset M\times M$ be a $\delta$-radius closed tubular neighborhood induced by the normal bundle $\mathcal{N}_{g\times g}(\Delta_{U_p}, U_p\times U_p)$. By $M^{out-out}$ we denote the set of outward-pointing chords $\overline{\nu_\delta}(\Delta_{U_p})\subset M\times M$. Assume that each quadrant $Q_j$ of $p$ is either outward- or inward-pointing. If $U_p$ and $\delta>0$ are sufficiently small, then $M^{out-out}$ is $\codim 0$ Whitney stratified subset of $\overline{\nu}_\delta(\Delta_{U_p})$ and 
$$M^{out-out}\xrightarrow{\pi_{\mathcal{N}}}\Delta_{U_p}$$
is a locally trivial stratified fiber bundle. Here the projection $\pi_{\mathcal{N}}$ is determined by the bundle projection in identification of $\overline{\nu_\delta}(\Delta_{U_p})$ with a disk subbundle of $\mathcal{N}_{g\times g}(\Delta_{U_p}, U_p\times U_p)$.

In particular, any generalized tangent space of $\pi^{-1}_\mathcal{N}((p, p))$ at $(p, p)$ lies in the span of $w_M(p)\times -w_M(p)$ and $w_m(p)\times-w_m(p)$ (for the definition of the generalized tangent space see Remark \ref{rem_stratified}).
\end{lemma}

\begin{proof}
For $i=1, 2$ let us consider functions $F_i(\cdot,\cdot):=F_{i,\cdot}(\cdot):\overline{\nu_\delta}(\Delta_{U_p})\setminus\Delta_{U_p}\rightarrow\R$ and put $F=(F_1, F_2)$. If $U_p$ together with $\varepsilon>0$ are small, then by Lemma \ref{lemma_aux_differential} $F$ is stratum transverse to $[0, \infty)^2$. Hence, if there is at least one outward-pointing quadrant, then $F^{-1}([0, \infty)^2)$ is a $4$-manifold with corners. Then $M^{out-out}=F^{-1}([0, \infty)^2)\cup \Delta_{U_p}$, so we can consider on $M^{out-out}$ a stratification given by the canonical stratifications of the $4$-manifold with corners $F^{-1}([0, \infty)^2)$ and $\Delta_{U_p}$. We would like to show that this stratification is Whitney.

Let us canonically extend $F_1$ and $F_2$ to $\Delta_{U_p}$. Then we claim that $\Delta_{U_p}$ is a Bott nondegenerate manifold (Definition \ref{defn_morse_bott}) for the functions $F_1$ and $F_2$. Moreover, we claim that in both cases the metric Hessian has eigenvectors $w_M \times -w_M$ and $w_m \times -w_m$ with the eigenvalues equal to $\kappa_M$ and $\kappa_m$, respectively. For this, one has to consider for each $(\widetilde{p}, \widetilde{p})\in\Delta_{U_p}$ the graphical coordinates $(x_1, y_1, f(x_1, y_1), x_2, y_2, f(x_2, y_2))$ centered at $(\widetilde{p}, \widetilde{p})$, see coordinates in Lemma \ref{lemma_convexity}. Then the rest follows from a straightforward computation as in Lemma \ref{cor_convexity_tangent}; we will leave the computation to the reader.

Now, similarly to the proof of Lemma \ref{lemma_convexity}, we can apply the Morse-Bott Lemma \cite{banyaga2004lectures, hirsch_1997} to each of $F_1$ and $F_2$ and conclude the Whitney conditions along the stratum $\Delta_{U_p}$.

Then the lemma follows from Thom's first isotopy Lemma \ref{thm_thom_isot}. In particular, to show that $\pi_\mathcal{N}$ is a submersion on each strata, one has to use Lemma \ref{lemma_aux_differential}.
\end{proof}

\begin{lemma}\label{lemma_appl_convex}For each $(s, \theta)\in \R/T\mathbb{Z}\times(-\pi/2, \pi/2)$ we consider $\lbrace Q_j\rbrace_{j\in\lbrace1,\dots, 4\rbrace}$ - quadrants on a $O(\varepsilon)$-small neighbourhood of $\Gamma(s, \theta)=p\in T_{K, \varepsilon}$, which were defined in Lemma \ref{lemma_convexity}.

There is a positive constant $C_s$ such that for each sufficiently small $\varepsilon>0$ the following holds
\begin{itemize}
\item If $\theta\in[\varepsilon^{1/2} C_s,\pi/2-\varepsilon^{1/2} C_s]$, then $Q_1, Q_2$ are outward-pointing and $Q_3, Q_4$ are inward-pointing.
\item If $\theta\in[-\pi/2+\varepsilon^{1/2} C_s, -\varepsilon^{1/2} C_s]$, then $Q_1, Q_2$ are inward-pointing and $Q_3, Q_4$ are outward-pointing.
\end{itemize}
See also Figure \ref{figure_diagonal_convexity2}.
\begin{figure}[!htbp]
\labellist
\pinlabel $w_M$ at 250 630
\pinlabel $w_m$ at 280 748
\pinlabel $T_{K, \varepsilon}$ at 100 850
\pinlabel $p$ at 200 695
\endlabellist
\centering
\includegraphics[scale=0.68]{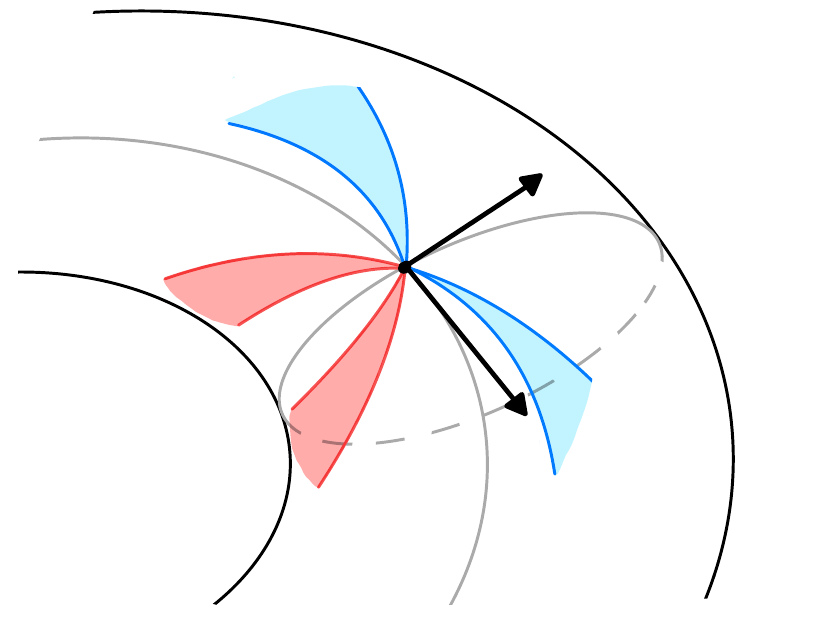}
\vspace{0.3cm}
\caption{Local description of chords emanating from $p=\Gamma_\varepsilon(s, \theta)$ for $\theta\in[\varepsilon^{1/2} C_s, \pi/2- \varepsilon^{1/2} C_s]$. Blue regions describe the possible endpoints for \textit{``out-out''} chords. Red regions describe the possible endpoints for \textit{``in-in''} chords. Grey lines are integral curves of principal directions from $p$.}
\label{figure_diagonal_convexity2}
\end{figure}
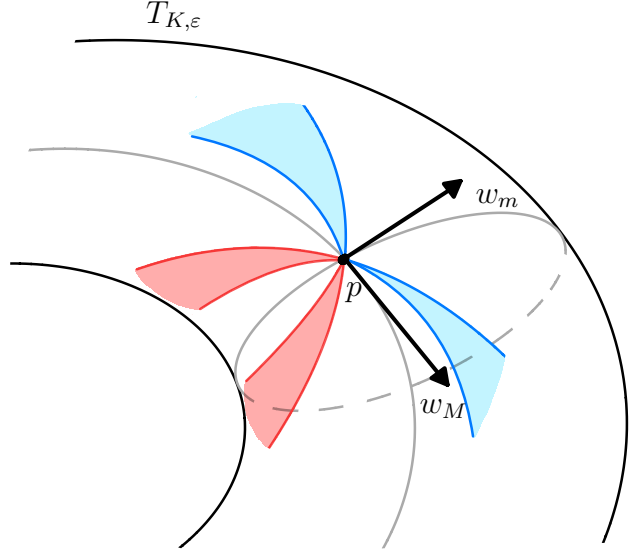
\end{lemma}

\begin{proof}
The lemma follows from Lemmata \ref{lemma_curvature} and \ref{lemma_convexity}. Observe that at $(s, \theta)$ it holds:
\begin{align*}
\nabla_{w_M}\kappa_M&=-\frac{1}{d}(\partial_s-\tau(s)\partial_\theta)\left(\frac{\kappa(s)\cos(\theta)}{d}\right)\\
&=-\frac{\dot{\kappa}(s)\cos(\theta)+\tau(s)\kappa(s)\sin(\theta)}{d^3},\\
\nabla_{w_m}\kappa_M&=\frac{1}{\varepsilon}\partial_\theta\left(\frac{\kappa(s)\cos(\theta)}{d}\right)\\
&=-\frac{\kappa(s)\sin(\theta)}{\varepsilon d^2},\\
\nabla_{w_M}\kappa_m&=0,\\
\nabla_{w_m}\kappa_m
&=0.
\end{align*}
Hence, on the $Q_j$ quadrant inequality (\ref{eqn_attracted}) becomes equivalent to
$$\rho_b^j\big(-2\rho_a^j d^{1/2}(\dot{\kappa}(s)\cos(\theta)+\tau(s)\kappa(s)\sin(\theta))-(1/\varepsilon^{1/2})6\kappa^{3/2}(s)\sin(\theta)|\cos(\theta)|^{1/2}\big)< 0.$$
Which can be rewritten as
\begin{equation}\label{eqn_bound_diag}
\varepsilon^{1/2}A(s, \theta)-\rho_b^j\sin(\theta)|\cos(\theta)|^{1/2}< 0
\end{equation}
for some bounded function $A(s, \theta)$. And the lemma follows.
\end{proof}

\begin{rem}\label{rem_descr_diag} The diagonal $\Delta_\varepsilon$ is diffeomorphic to $S^1\times [0, 1]$, where $S^1\times\lbrace0, 1\rbrace$ correspond to points of $T_{K, \varepsilon}$ with $\kappa_G=0$. In particular, $\Delta_\varepsilon\subsetneq \Delta_{full}$.
\end{rem}

\begin{cor}\label{lem_cusp_fibre} Put $C_K:=2\sup \lbrace C_{s_1}\,|\, s_1\in\R/T\mathbb{Z}\rbrace$. In a small closed tubular neighborhood of 
\begin{equation}\label{def_diag_stand}
\Delta_\varepsilon^{cusp}:=\Delta_{\varepsilon}\setminus\big\lbrace (s_1, \theta_1,s_2, \theta_2)\,|\, \theta_1\in (k\pi/2-\varepsilon^{1/2} C_{K}, k\pi/2+\varepsilon^{1/2} C_{K}), k\in\mathbb{Z}\big\rbrace
\end{equation}
it holds that $M_{K, \varepsilon}$ has a structure of a locally trivial stratified fibration. The fibers are identified with the union of two cuspidal regions. More precisely, the map identifying each fiber $\pi^{-1}_\mathcal{N}(p)$ with 
$$\big\lbrace (x, y)\in\R^2\,|\,y^2\leq|x^3|\wedge x^2+y^2\leq 1\big\rbrace$$ is a stratum preserving homeomorphism, which is a diffeomorphism outside of $p$.
\end{cor}

\begin{proof} The corollary follows from the combination of Lemmata \ref{lemma_appl_convex}, \ref{conj_local_fibr_diag} and \ref{cor_convexity_tangent}.
\end{proof}

\begin{example}\label{example_diag} In a simplified case, we are going to illustrate the behavior of $M_{K, \varepsilon}$ near the diagonal. We will consider $K$ with $\tau=\dot{\kappa}=0$ (note that this will, in particular, also simplify inequality $(\ref{eqn_bound_diag})$  to $\rho_b^j\sin(\theta)|\cos(\theta)|^{1/2}>0$).

Let $T_{K, \varepsilon}$ be a standard torus parameterized by $$\Gamma_\varepsilon(s, \theta)=\big(\cos(s)(3+\cos(\theta)), \sin(s)(3+\cos(\theta)), \sin(\theta)\big).$$ 
Let $\widehat{M}_{K, \varepsilon}$ denotes the restriction of $M_{K, \varepsilon}$ to the set $\lbrace(s_1, \theta_1, s_2, \theta_2)\in(\R/T\mathbb{Z}\times S^1)^2\,|\,\widetilde{d}(s_1, s_2)\leq\delta_K\rbrace$ and $\pi_{s_1}:\widehat{M}_{K, \varepsilon}\rightarrow\R/T\mathbb{Z}$ be the canonical projection. We would like to inspect the fiber $\pi_{s_1}^{-1}(3\pi/2)$.

To enlight Conjecture \ref{conj_diag} and Remark \ref{rem_enum_alg} we can explicitly compute the solutions of $F^{[\varepsilon]}_1=0\wedge F^{[\varepsilon]}_2=0$ in $(\theta_1, s_2, \theta_2)$-coordinates. If we ignore the trivial solutions coming from $\Delta_\varepsilon$, we obtain four (piecewise smooth) curves $\lbrace c_j\rbrace_{j\in\lbrace 1,\dots, 4\rbrace}:$
\begin{align*}
c_{1, 2}(t)&:=\left(2\pi-t, \pi\pm\arcsin\left(\frac{3-\cos(t)+2\sec(t)}{3+\cos(t)}\right), t\right),\\
c_{3, 4}(t)&:=(\pi\pm\pi/2, t, \pi\pm\pi/2).
\end{align*}
The curves $c_1$ and $c_2$ intersect (touch) at the point $(\pi, 3\pi/2, \pi)$, and their union is homeomorphic to a cone and consists of two smooth curves. See Figure \ref{figure_diagonal_example}.
\begin{figure}[!htbp]
\labellist
\pinlabel $\theta_1$ at 27 180
\pinlabel $s_2$ at 162 35
\pinlabel $0$ at 308 -9
\pinlabel $\theta_2$ at 358 87
\pinlabel $2\pi$ at 398 155
\pinlabel $\textcolor{blue}{\bullet}$ at 90 182
\pinlabel $\textcolor{blue}{\bullet}$ at 251 142
\endlabellist
\centering
\includegraphics[scale=0.77]{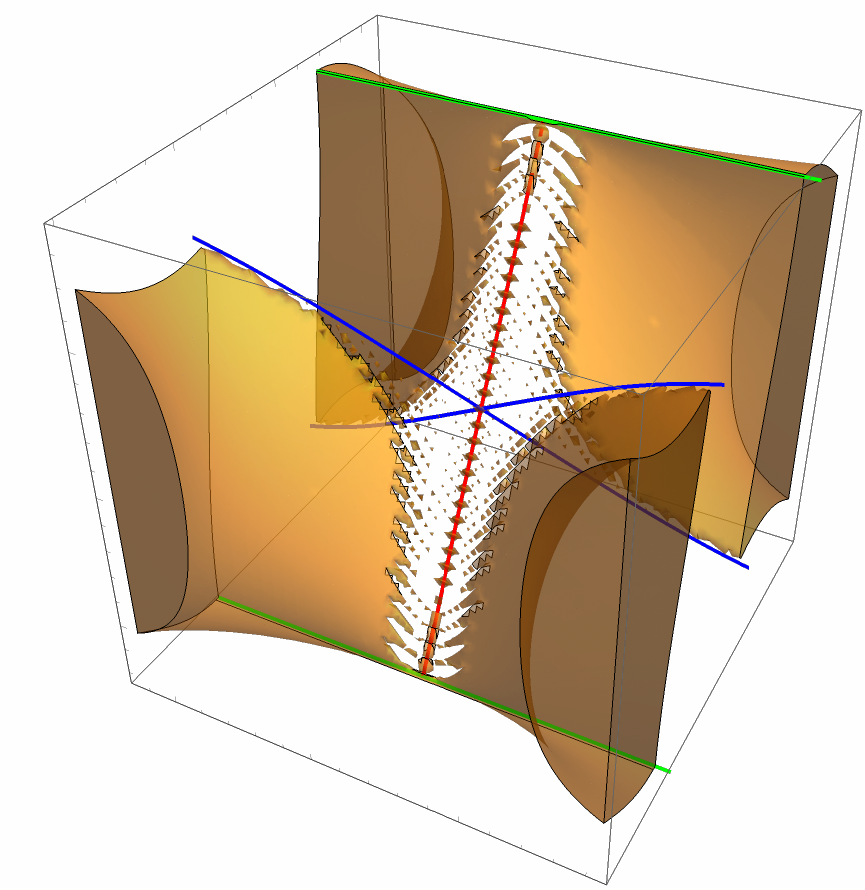}
\vspace{0.3cm}
\caption{A visualisation of fiber $\pi_{s_1}^{-1}(3\pi/2)$ in $(\theta_1, s_2, \theta_2)$-coordinates. The red line corresponds to points on $\Delta_\varepsilon$, the blue curves denote $c_1$ and $c_2$, and the green lines represent $c_3$ and $c_4$. Two blue dots represent critical points of $\pi_{s_2}$ (Conjecture \ref{conj_diag} $(ii.)$).
Note that the singular behavior of $\Delta_\varepsilon$ causes lower quality of the visualization in a neighborhood of $\Delta_\varepsilon$. But still, we can see the cuspidal fibration from Conjecture \ref{lem_cusp_fibre}.}
\label{figure_diagonal_example}
\end{figure}
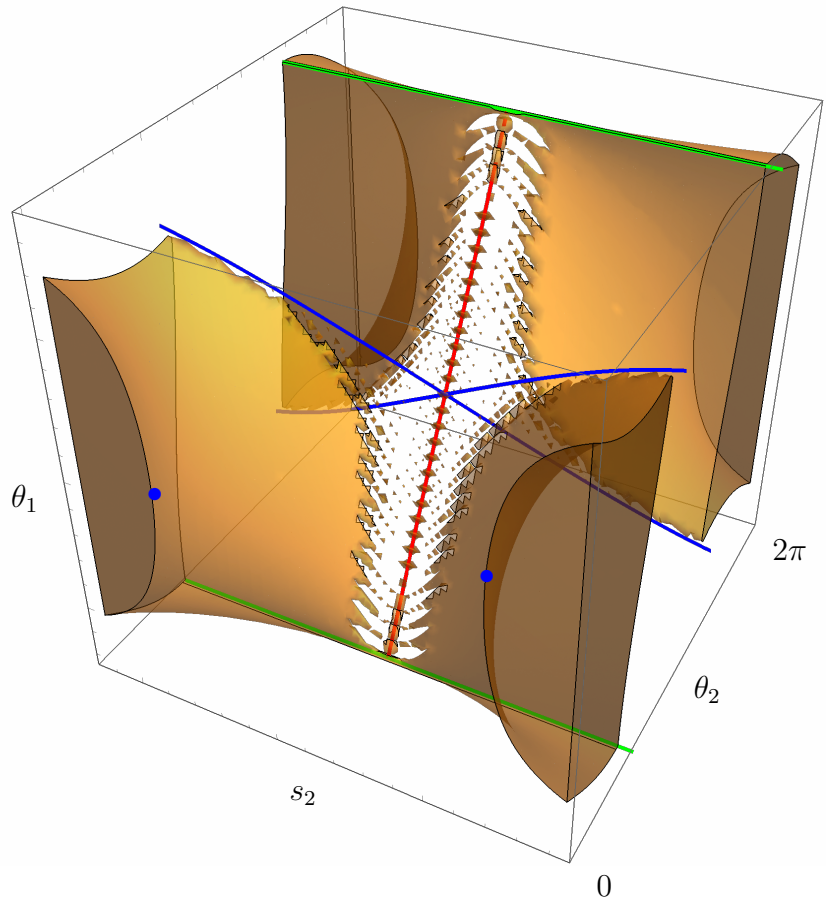
\end{example}

\newpage
\section{Energy functions}
\label{s:energy}
We are going to study $E_0, E_\varepsilon$; the restrictions of the standard energy function $E:\R^3\times\R^3\rightarrow\R$ to $K\times K$ and $T_{K, \varepsilon}\times T_{K, \varepsilon}$  $(M_{K, \varepsilon})$, respectively. In particular, we show that for generic $K$ and $\varepsilon>0$ the functions $E_0, E_\varepsilon$ are Morse-Bott, and we also relate their critical points. In the end of the section we inspect how the cuspical geometry of $M_{K, \varepsilon}$ restricts the number of potential $-\nabla E_\varepsilon$-trajectories in $M_{K, \varepsilon}$ that can reach the minimum $\Delta_\varepsilon$ before flowing out of $M_{K, \varepsilon}$. This will also depend on the non-resonance of eigenvalues of $-D\nabla E_{\varepsilon}(\Delta_\varepsilon)$ and the behavior of $-\nabla E_\varepsilon$ along $\partial M_{K, \varepsilon}$. 

\begin{defn}\label{defn_morse_bott}Let $(M, g)$ be a closed oriented Riemanninan manifold. A smooth function $f:M\rightarrow\R$ is called \textbf{Morse-Bott}, if the set of critical points $Crit(f)$ is a disjoint union of connected submanifolds $\lbrace Q_j\rbrace_j$ of $M$ with the following property: the Hessian $H^g[f]$, when restricted fiberwise to the normal bundle of $Q_j$ in $M$ $(:=\mathcal{N}^g(Q_j, M))$, defines a nondegenerate bilinear form.
The $Q_j$ will be called \textbf{Bott-nondegenerate (critical) submanifold}.
 
In addition, if $\mathcal{N}_{\pm}^g(Q_j, M)$ denotes the splinting of $\mathcal{N}^g(Q_j, M)$ onto positive and negative eigenspaces of $H^g[f]$, then we define the \textbf{Morse index of $Q_j$} as $$Ind(Q_j)=\hbox{rank}(\mathcal{N}_{-}^g(Q_j, M)).$$

$Crit_k(f)$ will denote the subset of $Crit(f)$ that consists of Bott-nondegenerate critical submanifolds of the Morse index $k$.
\end{defn}

\begin{defn}Let $f$ be a Morse-Bott function on the Riemannian manifold $(M, g)$. Then the (negative) \textbf{gradient flow} is the flow of $-\nabla^g f$ and will be denoted by $\phi^t_f$ or $\cdot_f$. 
\end{defn}

\begin{defn} Let $f:M\rightarrow\R$ be a Morse-Bott function and $Q\subset Crit(f)$. Then we define the \textbf{unstable manifold} $W_f(Q_j)$ as
$$W_f^{u}(Q):=\lbrace x\in M\,|\,\lim_{t\rightarrow-\infty}\phi^{t}_f(x)\in Q\rbrace.$$
and the \textbf{stable manifold} $W_f(Q)$ as
$$W_f^s(Q):=\lbrace x\in M\,|\,\lim_{t\rightarrow+\infty}\phi^{t}_f(x)\in Q\rbrace.$$
\end{defn}

\begin{rem} \cite[Prop 3.2]{Austin1995MorseBottTA} \label{rem_stable_dim}$W_f^{u}(Q)$ is the image of an injective immersion of $\mathcal{N}_{-}^g(Q, M)$ into $M$ and the map $x\mapsto\lim_{t\rightarrow-\infty}\phi^t_f(x)$ gives to $W_f^{u}(Q)$, in a neighborhood of $Q$, a structure of locally trivial fiber bundle over $Q$. The analogous statement also holds for $W_f^{s}(Q)$.
\end{rem}

\begin{defn}On $K$ we define
\begin{align*}
E_0:\hspace{1cm}R/T\mathbb{Z}\times R/T\mathbb{Z}\hspace{0.85cm}&\longrightarrow\hspace{2cm}\mathbb{R}\\
(s_1, s_2)&\longmapsto\frac{1}{2}||\gamma(s_2)-\gamma(s_1)||^2,
\end{align*}
where $\gamma:\R/T\mathbb{Z}\rightarrow\R^3$ is an arc length parametrization of $K$ and $\langle\cdot, \cdot\rangle$ is the induced metric on $K$ from $(\R^3, g_{Euc})$.
\end{defn}

\begin{rem}\label{rem_grad_0}We can compute
\begin{align*}
\partial_{s_1}E_0&=-\langle P, \dot{\gamma}(s_1)\rangle,\\
\partial_{s_2}E_0&=\langle P, \dot{\gamma}(s_2)\rangle.
\end{align*}
Hence, the critical points of $E_0$ are the solutions of 
$$\langle P, \dot{\gamma}(s_i)\rangle=0, \hbox{ for }i=1,2,$$
which are the binormal chords on $K$ and points from $\Delta_0$.
Next, we can compute
\begin{align*}
\partial^2_{s_1}E_0&=1-\kappa(s_1)\langle P, n(s_1)\rangle,\\
\partial_{s_1}\partial_{s_2}E_0&=-\langle \dot{\gamma}(s_1), \dot{\gamma}(s_2)\rangle,\\
\partial^2_{s_2}E_0&=1+\kappa(s_2)\langle P, n(s_2)\rangle.
\end{align*}
We also point out that, by the construction, $\gamma$ is an isometry between the flat metric on $\R/T\mathbb{Z}$ and $\langle\cdot, \cdot\rangle$ on $K$.
\end{rem}

\begin{lemma}\label{lemma_generic_morse_knot} For a generic $K$ in $\R^3$ the function $E_0$ is Morse outside of the Bott nondegenerate manifold $\Delta_0$.
\end{lemma}

\begin{proof}
For any $K$ the Bott-nondegenaracy of $\Delta_0$ is immediate from Remark \ref{rem_grad_0}, see also \cite[Lemma 7.10]{Cieliebak2016KnotCH}.

Let 
\begin{align*}
\mathfrak{s}:K\times K&\rightarrow T^\ast(K\times K)\\
p&\mapsto dE_0(p)
\end{align*}
be a section with zeros as critical points of $E_0$. Recall from Remark \ref{rem_split_tangent} that along the zero section, $p\in\mathfrak{s}^{-1}(0)$, we have the canonical splitting of $T_{(p, 0)}T^\ast(K\times K)$. In particular, we consider the vertical linearized map
\begin{align*}
d^v\mathfrak{s}(p):T_p(K\times K)&\rightarrow T^\ast_p(K\times K)\\
v&\mapsto H[E_0](v,\cdot)\quad\big(=d^2E_0(v,\cdot)\big).
\end{align*}
By the Rank-nullity theorem, we obtain that $d^v\mathfrak{s}(p)$ is a Fredholm map with the index zero. In particular, $d^v\mathfrak{s}(p)\pitchfork 0$ iff $d^v\mathfrak{s}(p)$ is injective. And by the construction of $d^v\mathfrak{s}(p)$, this is equivalent to $p$ being a nondegenerate critical point of $E_0$.

Now, let us consider a section
\begin{align*}
S:C^k(S^1, \R^3)\times\big((S^1\times S^1)\setminus \Delta_0\big)&\rightarrow\R^2\\
(\gamma, s_1, s_2)&\mapsto dE_0(s_1, s_2),
\end{align*}
where $C^k(S^1, \R^3)$ is a Banach space with $C^k$ topology and $k$ is big enough. We note that $\gamma$ is here just an arbitrary element of $C^k(S^1, \R^3)$ and that $E_0$ is here understood as $(\gamma\times\gamma)^\ast E$, where $E$ is the standard energy function on $(\R^3\times\R^3, g_{Euc})$. 

We are going to show the surjectivity of the linearised map $dS(\gamma, s_1, s_2):C^k(S^1, \R^3)\times\R^2\rightarrow\R^2$ at the zero section. 

We compute that
\begin{align*}
dS(\gamma, s_1, s_2)\cdot(\hat{\gamma}, 0, 0)=\begin{pmatrix}
  -\langle\hat{\gamma}(s_2)-\hat{\gamma}(s_1), \dot{\gamma}(s_1)\rangle-\langle\gamma(s_2)-\gamma(s_1), \dot{\hat{\gamma}}(s_1)\rangle \\
    \langle\hat{\gamma}(s_2)-\hat{\gamma}(s_1), \dot{\gamma}(s_2)\rangle+\langle\gamma(s_2)-\gamma(s_1), \dot{\hat{\gamma}}(s_2)\rangle
\end{pmatrix}.
\end{align*}
We choose $\hat{\gamma}_{i}, i=1, 2,$ as follows. We put $\hat{\gamma}_1(s_2)-\hat{\gamma}_1(s_1)=\dot{\hat{\gamma}}_1(s_1)=\gamma(s_2)-\gamma(s_1)$ and $\dot{\hat{\gamma}}_1(s_2)=\dot{\gamma}(s_2)$. And analogously, $\hat{\gamma}_2(s_2)-\hat{\gamma}_2(s_1)=\dot{\hat{\gamma}}_2(s_2)=\gamma(s_2)-\gamma(s_1)$ and $\dot{\hat{\gamma}}_2(s_1)=\dot{\gamma}(s_1)$.

Hence, we obtain that
\begin{align*}
dS(\gamma, s_1, s_2)\cdot(\hat{\gamma}_1, 0, 0)=\begin{pmatrix}
  -||\gamma(s_2)-\gamma(s_1)||^2 \\
    0
\end{pmatrix}
\end{align*}
and
\begin{align*}
dS(\gamma, s_1, s_2)\cdot(\hat{\gamma}_2, 0, 0)=\begin{pmatrix}
  0 \\
     ||\gamma(s_2)-\gamma(s_1)||^2
\end{pmatrix}.
\end{align*}
Since $||\gamma(s_2)-\gamma(s_1)||^2\neq 0$, we obtain the surjectivity of $dS$ along the zero section. Now, by the Implicit function theorem, see for example \cite{cieliebak_2018}, $S^{-1}(0)$ is a Banach manifold. Let us consider the projection $\pi_1:S^{-1}(0)\rightarrow C^k(S^1, \R^3)$, which is, for $k$ big enough, smooth. By Sard-Smale theorem, \cite{cieliebak_2018}, the set of regular values $(=C_{reg})$ of $\pi_1$ is residual in $C^k(S^1, \R^3)$. By \cite[Cor 1.6, Thm 2.13]{hirsch_1997}, $C^k$-embeddings $(=C_{embedd})$ are dense and open in $C^k(S^1, \R^3)$. Hence, for each $\gamma$ from the generic subset $C_{reg}\cap C_{embedd}$ of $C^\infty(S^1, \R^3)$ the corresponding section $\mathfrak{s}\pitchfork 0$. And, in particular, $E_0$ is Morse-Bott.
\end{proof}

\noindent{\textit{From now on, we will work with knots that are generic in the spirit of Lemma \ref{lemma_generic_morse_knot}}}

\begin{defn}\label{def_energ_2} Let $K$ be a knot in $\R^3$ and $\varepsilon\in (0, \varepsilon_{good}]$. Then we define
\begin{align*}
E_\varepsilon:\hspace{1cm}(\R/T\mathbb{Z}\times S^1)^2\hspace{0.85cm}&\longrightarrow\hspace{2cm}\mathbb{R}\\
(s_1, \theta_1, s_2, \theta_2)&\longmapsto\frac{1}{2}||\Gamma_\varepsilon(s_2, \theta_2)-\Gamma_\varepsilon(s_1, \theta_1)||^2.
\end{align*}
Also, by $E^{out-out}_\varepsilon$ we will denote the restriction of $E_\varepsilon$ to $M_{K, \varepsilon}$.
\end{defn}

\begin{rem}\label{rem_grad_comp}We would like to describe $-\nabla E_\varepsilon$.

We recall our notation that $\langle\cdot,\cdot\rangle$ is the standard metric in $\R^3$ or $\R^3\times\R^3$ and $\widetilde{\langle\cdot,\cdot\rangle}$ is the standard metric in the coordinate charts  $(s_1, s_2)$  or $(s_1, \theta_1, s_2, \theta_2)$, see Notation \ref{notation_on_circle}. Hence, the gradients with respect to those metrics will be denoted as $\nabla$ and $\widetilde{\nabla}$, respectively, and so on.

Let us consider the map $$\bm{\Gamma}_\varepsilon:(s_1, \theta_1, s_2, \theta_2)\mapsto\left(\Gamma_\varepsilon(s_1, \theta_1),\Gamma_\varepsilon(s_2, \theta_2)\right)\subset\R^3\times\R^3$$
which is not, unlike in the knot case, an isometry.

First, we compute
\begin{align*}
\partial_{s_1}E_\varepsilon&=-\langle \Gamma_\varepsilon(s_2, \theta_2)-\gamma(s_1), \partial_{s_1}\Gamma_\varepsilon(s_1, \theta_1) \rangle,\\
\partial_{\theta_1}E_\varepsilon&=-\langle \Gamma_\varepsilon(s_2, \theta_2)-\gamma(s_1), \partial_{\theta_1}\Gamma_\varepsilon(s_1, \theta_1) \rangle,\\
\partial_{s_2}E_\varepsilon&=\langle \gamma(s_2)-\Gamma_\varepsilon(s_1, \theta_1), \partial_{s_2}\Gamma_\varepsilon(s_2, \theta_2) \rangle,\\
\partial_{\theta_2}E_\varepsilon&=\langle \gamma(s_2)-\Gamma_\varepsilon(s_1, \theta_1), \partial_{\theta_2}\Gamma_\varepsilon(s_2, \theta_2) \rangle.
\end{align*}
where we used, as in Theorem \ref{thm_mfld_coners}, the fact that $$\langle \partial_{s_i}\Gamma_\varepsilon(s_i, \theta_i), v_i\rangle=\langle \partial_{\theta_i}\Gamma_\varepsilon(s_i, \theta_i), v_i\rangle=0.$$ And hence
\begin{equation}\label{eqn_energ_deriv}
\begin{split}
\partial_{s_1}E_\varepsilon&=-\langle P+\varepsilon v_2, d_1\dot{\gamma}(s_1)+\varepsilon\tau(s_1)v_1^{\bot}\rangle,\\
\partial_{\theta_1}E_\varepsilon&=-\langle P+\varepsilon v_2,\varepsilon v_1^{\bot} \rangle,\\
\partial_{s_2}E_\varepsilon&=\langle P-\varepsilon v_1, d_2)\dot{\gamma}(s_2)+\varepsilon\tau(s_2)v_2^{\bot} \rangle,\\
\partial_{\theta_2}E_\varepsilon&=\langle P-\varepsilon v_1,\varepsilon v_2^{\bot}\rangle,
\end{split}
\end{equation}
where $d_i=(1-\varepsilon \cos(\theta_i)\kappa(s_i))$.

Next, let us consider the pull-back metric $g=\bm{\Gamma}^{\ast}_\varepsilon\langle\cdot,\cdot\rangle$. $g$ is a product metric on the product of tori $(\mathbb{R}/T\mathbb{Z}\times S^1)\times(\mathbb{R}/T\mathbb{Z}\times S^1)$ and recall that by (\ref{eqn_metric_torus}) it holds
\begin{align}
\begin{split}\label{eqn_first_fund}
g(\partial_{s_i}, \partial_{s_i})&=\varepsilon^2 \tau^2(s_i)+(1-\varepsilon \cos(\theta_i)\kappa(s_i))^2,\\
g(\partial_{s_i}, \partial_{\theta_i})&=\varepsilon^2\tau(s_i),\\
g(\partial_{\theta_i}, \partial_{\theta_i})&=\varepsilon^2. 
\end{split}
\end{align}
Hence the inverse metric $g^{-1}$ is given by the matrix 
\begin{equation}\label{eqn_inv_prod_metric}
G^{-1}=\begin{pmatrix}
    1/d_1^2 & -\tau(s_1)/d_1^2 & 0 & 0 \\
    -\tau(s_1)/d_1^2 & \frac{d_1^2+\varepsilon^2\tau^2(s_1)}{\varepsilon^2 d_1^2} & 0 & 0\\
        0      & 0 & 1/d_2^2 & -\tau(s_2)/d_2^2 \\
            0 & 0 &  -\tau(s_2)/d_2^2 & \frac{d_2^2+\varepsilon^2\tau^2(s_2)}{\varepsilon^2 d_2^2}\\
\end{pmatrix}.
\end{equation}

Next, by the definition
\begin{align*}
\nabla E_\varepsilon&=\big(\partial_{s_1}E_\varepsilon(g^{-1})_{1, 1}+\partial_{\theta_1}E_\varepsilon(g^{-1})_{1, 2}\big)\partial_{s_1}
+\big(\partial_{s_1}E_\varepsilon(g^{-1})_{2, 1}+\partial_{\theta_1}E(g^{-1})_{2, 2}\big)\partial_{\theta_1}\\
&\quad+\big(\partial_{s_2}E_\varepsilon(g^{-1})_{3, 3}+\partial_{\theta_2}E_\varepsilon(g^{-1})_{3, 4}\big)\partial_{s_2}
+\big(\partial_{s_2}E_\varepsilon(g^{-1})_{4, 3}+\partial_{\theta_2}E_\varepsilon(g^{-1})_{4, 4}\big)\partial_{\theta_3}.
\end{align*}
Thus
\begin{equation}\label{eqn_grad_E}
\nabla E_\varepsilon=\begin{pmatrix}
      -(1/d_1)\langle P+\varepsilon v_2, \dot{\gamma}(s_1)\rangle\\
\tau(s_1)(1/d_1)\langle P+\varepsilon v_2, \dot{\gamma}(s_1)\rangle-(1/\varepsilon)\langle P+\varepsilon v_2, v_1^\bot\rangle\\
(1/d_2)\langle P-\varepsilon v_1, \dot{\gamma}(s_2)\rangle\\
-\tau(s_2)(1/d_2)\langle P-\varepsilon v_1, \dot{\gamma}(s_2)\rangle+(1/\varepsilon)\langle P-\varepsilon v_1, v_2^\bot\rangle    
\end{pmatrix}.
\end{equation}
\end{rem}

\begin{lemma}\label{lemma_morse_corresp} Let us restrict $E_0$ and $E_\varepsilon$ outside of the diagonals $\Delta_0$ and $\Delta_{full}$, respectively. We denote these restrictions by $\widehat{E}_0$ and $\widehat{E}_\varepsilon$. Then it holds that if $\widehat{E}_0$ is Morse and $\varepsilon>0$ is sufficiently small, then $\widehat{E}_\varepsilon$ is also Morse.

In addition, there is a bijection between $Crit_{k}(\widehat{E}_0)$ and $Crit_{k+1}(\widehat{E}_\varepsilon^{out-out})$, for any $k\in \mathbb{N}_0$. More precisely, the sets $Crit(\widehat{E}_0)$ and $Crit(\widehat{E}_\varepsilon)\cap Int(M_{K, \varepsilon})$ are in a bijection and $$Ind_{E_0}(p_0)=Ind_{E_\varepsilon}(p_\varepsilon)-1,$$
where $p_0$ and $p_\varepsilon$ are the corresponding critical points from $Crit(\widehat{E}_0)$ and $Crit(\widehat{E}_\varepsilon)\cap Int(M_{K, \varepsilon})$, respectively. 
\end{lemma}

\begin{proof}
From Remark \ref{rem_grad_0} we know that the critical points of $\widehat{E}_0$ are the binormal chords on $K$, i.e. the solutions of
$$\langle P, \dot{\gamma}(s_i)\rangle=0.$$

Since our ambient space is $\R^3$ with the flat metric, we can relate the energy functionals from Definitions \ref{def_energ_1} and \ref{def_energ_2}, $\widehat{\bm{E}}_{\varepsilon}$ and $\widehat{E}_{\varepsilon}$, respectively. In particular, by Remark \ref{rem_crit_path}, the critical points $p_\varepsilon$ of $\widehat{E}_{\varepsilon}$ represent the binormal chords of $T_{K, \varepsilon}$, that is
$$P+\varepsilon v_2-\varepsilon v_1\,||\, v_i.$$
Alternatively, we can see this directly from relations (\ref{eqn_energ_deriv}). Altogether, at $p_\varepsilon$ it holds that $v_1=\pm v_2$ and each $p_\varepsilon$ is uniquely determined as the solution of 
\begin{equation}\label{eqn_of_crit_pt}
\langle P, \dot{\gamma}(s_i)\rangle=0\wedge\langle P, v_i^\bot\rangle=0.
\end{equation}
So, for each $p_\varepsilon\in Crit(\widehat{E}_\varepsilon)$ there is a $p_0\in Crit(\widehat{E}_0)$ such that $\pi_{s_1, s_2}(p_\varepsilon)=p_0$.
Now, we would like to compare the signs of the Hessian matrices of these $p_0$ and $p_\varepsilon$. 

But first, let us make the following auxiliary computations
\begin{align*}
\partial_{s_1}\langle P+\varepsilon v_2, \dot{\gamma}(s_1)\rangle&=-1+\kappa(s_1)\langle P+\varepsilon v_2, n(s_1)\rangle,\\
\partial_{\theta_1}\langle P+\varepsilon v_2, \dot{\gamma}(s_1)\rangle&=0,\\
\partial_{s_2}\langle P+\varepsilon v_2, \dot{\gamma}(s_1)\rangle&=\langle d_2\dot{\gamma}(s_2)+\varepsilon\tau(s_2)v_2^\bot, \dot{\gamma}(s_1)\rangle,\\
\partial_{\theta_2}\langle P+\varepsilon v_2, \dot{\gamma}(s_1)\rangle&=\varepsilon\langle v_2^\bot, \dot{\gamma}(s_1)\rangle
\end{align*}
and 
\begin{align*}
\partial_{s_1}\langle P+\varepsilon v_2, v_1^\bot\rangle&=-\tau(s_1)\langle P+\varepsilon v_2, v_1\rangle+\sin(\theta_1)\kappa(s_1)\langle P+\varepsilon v_2, \dot{\gamma}(s_1)\rangle,\\
\partial_{\theta_1}\langle P+\varepsilon v_2, v_1^\bot\rangle&=-\langle P+\varepsilon v_2, v_1\rangle,\\
\partial_{s_2}\langle P+\varepsilon v_2, v_1^\bot\rangle&=\langle d_2\dot{\gamma}(s_2)+\varepsilon\tau(s_2)v_2^\bot, v_1^\bot\rangle,\\
\partial_{\theta_2}\langle P+\varepsilon v_2, v_1^\bot\rangle&=\varepsilon\langle v_2^\bot, v_1^\bot\rangle.
\end{align*}

Now at any critical point $p_\varepsilon$ of $E_\varepsilon$ we have that
\begin{align}\label{eqn_sec_deriv_energ}
\begin{split}
\partial^2_{s_1}E_\varepsilon&=-d_1\big(-1+\kappa(s_1)\langle P+\varepsilon v_2, n(s_1)\rangle\big)+\varepsilon\tau^2(s_1)\langle P+\varepsilon v_2, v_1\rangle,\\
\partial_{\theta_1}\partial_{s_1}E_\varepsilon&=\varepsilon\tau(s_1)\langle P+\varepsilon v_2, v_1\rangle,\\
\partial_{s_2}\partial_{s_1}E_\varepsilon&=-d_1d_2\langle\dot{\gamma}(s_2), \dot{\gamma}(s_1)\rangle-\varepsilon\tau(s_1)d_2\langle \dot{\gamma}(s_2), v_1^\bot\rangle\\
&\quad-\varepsilon\tau(s_2)d_1\langle v_2^\bot, \dot{\gamma}(s_1)\rangle-\varepsilon^2\tau^2(s_2)\langle v_1^\bot, v_2^\bot\rangle,\\
\partial_{\theta_2}\partial_{s_1}E_\varepsilon&=-\varepsilon d_1 \langle v_2^\bot, \dot{\gamma}(s_1)\rangle-\varepsilon^2\tau(s_1)\langle v_2^\bot, v_1^\bot\rangle,\\
\\
\partial^2_{\theta_1}E_\varepsilon&=\varepsilon\langle P+\varepsilon v_2, v_1\rangle,\\
\partial_{s_2}\partial_{\theta_1}E_\varepsilon&=-\varepsilon d_2\langle \dot{\gamma}(s_2), v_1^\bot\rangle-\varepsilon^2\tau(s_1)\langle v_2^\bot, v_1^\bot\rangle,\\
\partial_{\theta_2}\partial_{\theta_1}E_\varepsilon&=-\varepsilon^2\langle v_2^\bot, v_1^\bot\rangle,\\
\\
\partial^2_{s_2}E_\varepsilon&=d_2\big(1+\kappa(s_2)\langle P-\varepsilon v_1, n(s_2)\rangle\big)-\varepsilon\tau^2(s_2)\langle P-\varepsilon v_1, v_2\rangle,\\
\partial_{\theta_2}\partial_{s_2}E_\varepsilon&=-\varepsilon\tau(s_2)\langle P-\varepsilon v_1, v_2\rangle,\\
\\
\partial^2_{\theta_2}E_\varepsilon&=-\varepsilon\langle P-\varepsilon v_1, v_2\rangle.\\
\end{split}
\end{align}

Hence, by Remark \ref{rem_grad_0} $\partial_{s_i}\partial_{s_j}E_0(p_0)=\partial_{s_i}\partial_{s_j}E_\varepsilon(p_\varepsilon)+O(\varepsilon)$. And also all other second derivatives of $E_\varepsilon$ are of type $O(\varepsilon)$.

This motivates us to the following comparison of the Hessian matrices of $\widehat{E}_0$ and $\widehat{E}_\varepsilon$. 

Because $K$ is generic, by Lemma \ref{lemma_generic_morse_knot} $\widehat{E}_0$ is Morse. Since nondegenerate critical points are isolated, we can assume that critical points of $\widehat{E}_0$ lie outside of the weakly special and weakly diagonal points. Let us pick any $p_\varepsilon\in Crit(\widehat{E}_\varepsilon)$ and put $p_0:=\pi_{s_1, s_2}(p_\varepsilon)$.

Recall that at the critical points, the signs of the Hessian matrix do not depend on the pull-back metric and can be described in coordinates just as second derivatives.

Now we write the Hessian matrix $\widetilde{H}[E_0](p_0)$ as $4\times 4$ matrix, where the even rows and columns are filled with zeros and we implicitly considered the standard metric $\widetilde{\langle\cdot,\cdot\rangle}$ on coordinates $(s_1, \theta_1, s_2, \theta_2)$. Hence $\widetilde{H}[E_\varepsilon](p_\varepsilon)$ can be seen as $\varepsilon$-perturbation of elements of $\widetilde{H}[E_0](p_0)$. By Kato's selection theorem, this gives us a $\varepsilon$-perturbation of the eigenvalues of matrices.

By our assumptions on $p_0$, $\widetilde{H}[E_0](p_0)$ has two eigenvalues $\lambda_{0; 1}, \lambda_{0; 2}$ such that $$\delta\leq\max\lbrace|\lambda_{0; 1}|, |\lambda_{0; 2}|\rbrace$$ for some $\delta\geq0$. Then $\widetilde{H}[E_\varepsilon](p_\varepsilon)$ has $4$ eigenvalues  $\lambda_{\varepsilon; 1}, \lambda_{\varepsilon; 2}, \lambda_{\varepsilon; 3}, \lambda_{\varepsilon; 4}$ such that
\begin{align*}
\lambda_{\varepsilon; 1}&=\lambda_{0; 1}+O(\varepsilon),\\
\lambda_{\varepsilon; 2}&=\lambda_{0; 2}+O(\varepsilon),\\
\lambda_{\varepsilon; 3}&=O(\varepsilon),\\
\lambda_{\varepsilon; 4}&=O(\varepsilon).
\end{align*}
Hence, if $0<\varepsilon\ll\delta$, the eigenvalues $\lambda_{\varepsilon; 1}, \lambda_{\varepsilon; 2}$ have the same signs as $\lambda_{\varepsilon; 1}, \lambda_{\varepsilon; 2}$, respectively.

Now we aim to inspect the sign of $\det\big(\widetilde{H}[E_\varepsilon](p_\varepsilon)\big)$. We can compute
\begin{align*}
\det\big(\widetilde{H}[E_\varepsilon](p_\varepsilon)\big)&=\big(\lambda_{0; 1}\lambda_{0; 2}\big)(-1)\varepsilon^2\langle P+\varepsilon v_2, v_1\rangle\langle P-\varepsilon v_1, v_2\rangle+O(\varepsilon^3)\\
&=\big(\lambda_{0; 1}\lambda_{0; 2}\big)(-1)\varepsilon^2\langle P, v_1\rangle\langle P, v_2\rangle+O(\varepsilon^3).
\end{align*}

Since $p_\varepsilon$ is a binormal chord, from (\ref{eqn_of_crit_pt}) we have that $P\,||\,v_i$. Moreover, since $p_\varepsilon\notin \Delta_{full}$, we conclude that for $\varepsilon>0$ sufficiently small it holds
$\det\big(\widetilde{H}[E_\varepsilon](p_\varepsilon)\big)\neq 0$.
In particular, we showed that if $\varepsilon>0$ is small and $\widehat{E}_0$ is Morse, then $\widehat{E}_\varepsilon$ is Morse too.

At the critical point $p_\varepsilon\notin \Delta_{full}$, we can deduce from $P$ being parallel to $v_i$ even more. It holds that
\begin{equation}\label{ref_eqn_short_lift}
v_i=\pm D_i^{-1/2}\langle P, n(s_i)\rangle n(s_i)+\langle P, b(s_i)\rangle b(s_i),
\end{equation}
where $D_i=\langle P, n(s_i)\rangle^2+\langle P, b(s_i)\rangle^2.$ For the geometric picture, see the proof of Lemma \ref{lemma_standard_square}. Note that since $p_\varepsilon$ does not project onto any diagonal or special point, $D_i>0$. So for each $p_0\in Crit(\widehat{E}_0)$ there are four possible $p_\varepsilon$ with $\pi_{s_1, s_2}(p_\varepsilon)=p_0$. Then it holds
\begin{align}
\begin{split}\label{comput_binom}
\hspace{-2cm}\langle P, v_i\rangle&=\pm D_i^{-1/2}\big\langle \langle P, \dot{\gamma}(s_i)\rangle \dot{\gamma}(s_i)+\langle P, n(s_i)\rangle n(s_i)+\langle P, b(s_i)\rangle b(s_i), \langle P, n(s_i)\rangle n(s_i)+\langle P, b(s_i)\rangle b(s_i)\big\rangle\\
&=\pm D_i^{-1/2}||\langle P, n(s_i)\rangle n(s_i)+\langle P, b(s_i)\rangle b(s_i)\rangle||^2.
\end{split}
\end{align}
Hence for only one of the four critical points $p_\varepsilon$, that project to $p_0$, it holds that 
\begin{equation}\label{eqn_sign_conv}
\langle P, v_1\rangle>0\wedge \langle P, v_2\rangle>0.
\end{equation}
Let $p_\varepsilon$ be the critical point satisfying (\ref{eqn_sign_conv}). Since $p_0$ is outside of weakly special and weakly diagonal points, for $\varepsilon>0$ sufficiently small it holds that $p_\varepsilon\in Int(M_{K, \varepsilon})$. Moreover, for $\varepsilon>0$ sufficiently small the sign of 
$\det\big(\widetilde{H}[E_\varepsilon](p_\varepsilon)\big)$ is \textit{opposite} to the sign of the product $\lambda_{0; 1}\lambda_{0; 2}$. Hence we obtain that $Ind_{E_0}(p_0)=Ind_{E_\varepsilon}(p_\varepsilon)-1$.
\end{proof}

\begin{rem} Later, we will see Lemma \ref{lemma_morse_corresp} also as a consequence of Multiple-time scale dynamics.
\end{rem}

\begin{lemma}\label{lem_ener_bott}$\Delta_{full}$ is a Bott nondegenerate critical manifold of $E_\varepsilon$. Moreover, if $p\in\Delta_{full}$, then the nonzero eigenvalues of $H[E_\varepsilon](p)$ are \textcolor{blue}{both equal to $2$} (we took $H[E_\varepsilon](p)$ with respect to the induced metric $g=\bm{\Gamma}^{\ast}_\varepsilon\langle\cdot,\cdot\rangle$).
\end{lemma}

\begin{proof}
Let $p\in\Delta_{full}$. By the proof of Lemma \ref{lemma_morse_corresp} we quickly obtain that $\widetilde{H}[E_\varepsilon](p)$ is given by
$$\begin{pmatrix}
    d^2+\varepsilon^2\tau^2(s)      & \varepsilon^2\tau(s) & -d^2-\varepsilon^2\tau^2(s) & -\varepsilon^2\tau(s) \\
    \varepsilon^2\tau(s) & \varepsilon^2 &  -\varepsilon^2\tau(s) & -\varepsilon^2\\
        -d^2-\varepsilon^2\tau^2(s)      & -\varepsilon^2\tau(s) & d^2+\varepsilon^2\tau^2(s) & \varepsilon^2\tau(s) \\
            -\varepsilon^2\tau(s) & -\varepsilon^2 &  \varepsilon^2\tau(s) & \varepsilon^2\\
\end{pmatrix}.$$
Hence, by Remark \ref{rem_grad_comp} we can compute that
\begin{align*}
H[E_\varepsilon](p)&=G^{-1}(p)\cdot \widetilde{H}[E_\varepsilon](p)\\
&=\varepsilon^{-2}d^{-2}\begin{pmatrix}
    \varepsilon^2      & -\varepsilon^2\tau(s) & 0 & 0 \\
    -\varepsilon^2\tau(s) & d^2+\varepsilon^2\tau(s)^2 &  0 & 0\\
        0      & 0 & \varepsilon^2      & -\varepsilon^2\tau(s) \\
            0 & 0 &  -\varepsilon^2\tau(s) & d^2+\varepsilon^2\tau(s)^2\\
\end{pmatrix}\\
&\hspace{1.6cm}\cdot
\begin{pmatrix}
    d^2+\varepsilon^2\tau^2(s)      & \varepsilon^2\tau(s) & -d^2-\varepsilon^2\tau^2(s) & -\varepsilon^2\tau(s) \\
    \varepsilon^2\tau(s) & \varepsilon^2 &  -\varepsilon^2\tau(s) & -\varepsilon^2\\
        -d^2-\varepsilon^2\tau^2(s)      & -\varepsilon^2\tau(s) & d^2+\varepsilon^2\tau^2(s) & \varepsilon^2\tau(s) \\
            -\varepsilon^2\tau(s) & -\varepsilon^2 &  \varepsilon^2\tau(s) & \varepsilon^2\\
\end{pmatrix}\\
&=\begin{pmatrix}
    1&0&-1&0\\
    0&1&0&-1\\
    -1&0&1&0\\
    0&-1&0&1
\end{pmatrix}.
\end{align*}
Observe that the matrix $H[E_\varepsilon](p)$ has a $2$-dimensional kernel spanned by $\partial_{s_1}+\partial_{s_2}$ and $\partial_{\theta_1}+\partial_{\theta_2}$. Moreover $H[E_\varepsilon](p)$ is positive definite in directions $\partial_{s_1}-\partial_{s_2}$ and $\partial_{\theta_1}-\partial_{\theta_2}$ with the same positive eigenvalues equal to $2$. Note also that the $0$-eigenspace is orthogonal to the $2$-eigenspace with respect to the metric $g$. This finishes the Bott nondegeneracy of $\Delta_{full}$.
\end{proof}

\begin{lemma}\label{lem_grad_standard} Let $x_\varepsilon\in M_{K, \varepsilon}$ such that $y_\varepsilon:=\pi_{s_1, s_2}(x_\varepsilon)\in \overline{S_K}$ (see Remark \ref{rem_delta_standard} for the notion of the standard set $S_K$). Provided that $\varepsilon>0$ is sufficiently small, it holds that
\begin{itemize}
\item[$(i.)$]If $F^{[\varepsilon]}_1(x_\varepsilon)=0$ and $F^{[\varepsilon]}_2(x_\varepsilon)>0,$ then $-\nabla E_\varepsilon(x_\varepsilon)$ is strictly inward-pointing into $M_{K, \varepsilon}$.
\item[$(ii.)$]If $F^{[\varepsilon]}_1(x_\varepsilon)>0$ and $F^{[\varepsilon]}_2(x_\varepsilon)=0$, then $-\nabla E_\varepsilon(x_\varepsilon)$ is strictly outward-pointing from $M_{K, \varepsilon}$.
\item[$(iii.)$] If $F^{[\varepsilon]}_1(x_\varepsilon)=F^{[\varepsilon]}_1(x_\varepsilon)=0$, then there is $\delta>0$ such that $$\lbrace x_\varepsilon\cdot_{E_\varepsilon}[-\delta, \delta]\rbrace\cap M_{K, \varepsilon}=x_\varepsilon$$ and $-\nabla E_\varepsilon(x_\varepsilon)\neq 0$.
\end{itemize}
See also Figure \ref{figure_square}.
\begin{figure}[!htbp]
\centering
\includegraphics[scale=0.68]{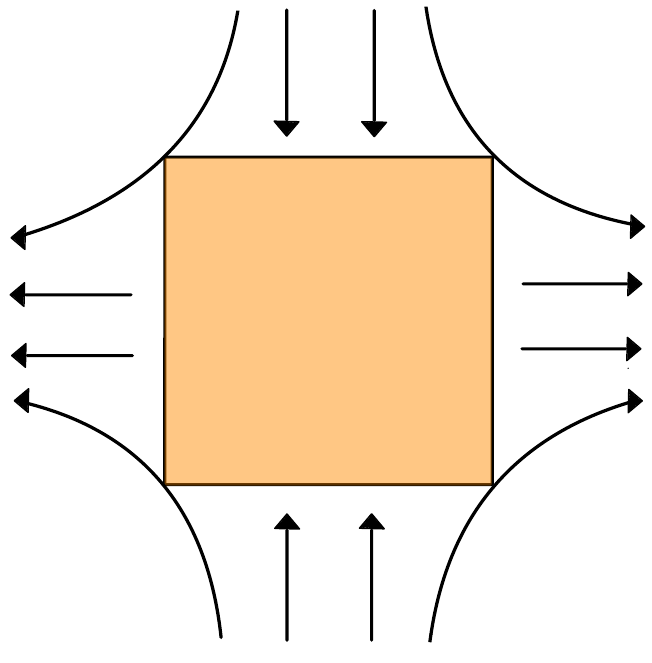}
\vspace{0.3cm}
\caption{The square diffeomorphic to $M_{K, \varepsilon}|_{\pi_{s_1, s_2}^{-1}(y_\varepsilon)}$, where the horizontal sides correspond to $F^{[\varepsilon]}_1=0$ and the vertical sides correspond to $F^{[\varepsilon]}_2=0$. Moreover, the arrows describe the behavior of the gradient flow of $-\nabla E_\varepsilon|_{\pi_{s_1, s_2}^{-1}(y_\varepsilon)}$ on the boundary $\partial M_{K, \varepsilon}|_{\pi_{s_1, s_2}^{-1}(y_\varepsilon)}$.}
\label{figure_square}
\end{figure}
\end{lemma}

\begin{proof}
We are going to prove only the case $(i.)$, since the other cases are analogous.

Let us consider a smooth family of functions 
\begin{equation}\label{enq_funct_aux_grad}
\big\lbrace f_\varepsilon:\lbrace x\in (\R/T\mathbb{Z}\times S^1)^2\,|\,\pi_{s_1, s_2}(x)\in \overline{S_K}\rbrace\rightarrow\R\big\rbrace_{\varepsilon\in(0, \varepsilon_{good}]}
\end{equation}
such that 
$$f_\varepsilon(x)=\varepsilon \widetilde{\langle \widetilde{\nabla} F^{[\varepsilon]}_1(y), -\nabla E_\varepsilon(x)\rangle}.$$
Note that by relation $(\ref{eqn_grad_E})$ the family of $f_\varepsilon$ extends smoothly to $0$.

Let $\lbrace\widehat{M}_{K, \varepsilon}\rbrace_{\varepsilon\in[0, \varepsilon_0]}$ be the compact manifolds of outward-pointing chords over $\overline{S_K}$ from Lemma \ref{lem_mfld_standard}. By Lemma \ref{lem_mfld_standard} all of these $\widehat{M}_{K, \varepsilon}$ are isotopic submanifolds of $(\R/T\mathbb{Z})^2$.

Let $\varepsilon\in(0, \varepsilon_0]$ and $x_\varepsilon\in\widehat{M}_{K, \varepsilon}$ be such that $F^{[\varepsilon]}_1(x_\varepsilon)=0$ and $F^{[\varepsilon]}_2(x_\varepsilon)>0$. By $x_0\in\widehat{M}_{K, 0} $ we denote the isotopic point to $x_\varepsilon$. In particular, $F^{[0]}_1(x_0)=0$ and $F^{[0]}_2(x_0)>0$.

Note that
$$f_0(x_0)=\langle P, v_1^\bot\rangle^2>0,$$
since we consider points over $\overline{S_K}$, and moreover $F^{[0]}_1(x_0)=\langle P, v\rangle=0$. Hence if $\varepsilon>0$ is sufficiently small, then also
$$f_\varepsilon(x_0)>0$$
by the smoothness of the family $(\ref{enq_funct_aux_grad}).$ Moreover, $x_\varepsilon$ and $x_0$ are isotopic. So if $\varepsilon>0$ is sufficiently small, then also
$$f_\varepsilon(x_\varepsilon)>0.$$
Finally, since $\varepsilon>0$, it holds that 
\begin{equation*}
\widetilde{\langle \widetilde{\nabla} F^{[\varepsilon]}_1(x_\varepsilon), -\nabla E_\varepsilon(x_\varepsilon)\rangle}>0
\end{equation*}
which finishes the case $(i.)$.
\end{proof}

\begin{conj}\label{conj_grad_diag}Let $x_\varepsilon \in M_{K, \varepsilon}\setminus \Delta_\varepsilon$ such that $\pi_{s_1, s_2}(x_\varepsilon)$ is an almost diagonal pair for $\varepsilon>0$ small (recall Definition \ref{lemma_aux_diag} for almost diagonal pairs). If $$F^{[\varepsilon]}_1(x_\varepsilon)=0\vee F^{[\varepsilon]}_1(x_\varepsilon)=0,$$
then $-\nabla E_\varepsilon(x_\varepsilon)$ is strictly outward-pointing from $M_{K, \varepsilon}$. I.e. $-\nabla E_\varepsilon(x_\varepsilon)\neq 0$ and there is $\delta>0$ such that $x_\varepsilon\cdot_{E_\varepsilon}(0, \delta]\cap M_{K, \varepsilon}=\emptyset$ and $x_\varepsilon\cdot_{E_\varepsilon}[-\delta, 0)\subset Int(M_{K, \varepsilon})$. See also Example \ref{exm_diag_grad}.
\end{conj}

\begin{example}\label{exm_diag_grad} We are going to describe Conjecture \ref{conj_grad_diag} for the outward-pointing chords on a closed oriented surface $M\subset \R^3$ in a neighborhood of a point $p\in M$ with the negative Gaussian curvature. 

Assume that $M$ is locally given as a graph $\lbrace (x, y, x^2-y^2-x^2 y)\rbrace$, where $p:=(0, 0, 0)$. So $(x, y)$ gives us the graphical coordinates of $M$. We can assume that the vector $(0, 0, 1)^T$ at $p$ is an outward-pointing normal vector, see Convention \ref{conv_princip}.

Let us split a neighborhood $U_p$ of $p$ into quadrants $Q_i$ as in Lemma \ref{lemma_convexity}. Then the quadrants $Q_1, Q_2$ are outward-pointing and $Q_3, Q_4$ are inward-pointing. Hence by Lemma \ref{conj_local_fibr_diag} the space of outward-pointing chords $M^{out-out}\subset U_p\times U_p$ has near $(p, p)$ a structure of a stratified fiber bundle $M^{out-out}\xrightarrow{\pi_\mathcal{N}}U_p\times U_p$. Here, the outward-pointing conditions are given by functions $F_1, F_2$ of coordinates $(x_1, y_1, x_2, y_2)$. 

By $E$ we will understand the restriction of the standard energy function from $(\R^3\times\R^3, g_{Euc})$ to $U_p\times U_p$.

Let us consider the set
$$L_p=\lbrace (x_1, y_1, x_2, y_2)\in M^{out-out}\,|\,(x_1, y_1)=(0, 0)\rbrace,$$
i.e., the set of outward-pointing chords starting at $p$. Recall that $L_p$ is not the fiber $\pi^{-1}_\mathcal{N}(p, p)$, but $M^{out-out}$ can be covered by the sets $\lbrace L_q\,|\,q\in\widetilde{U}_p\rbrace$ for some $\widetilde{U}_p\supset U_p$.

We would like to inspect the behavior of $-\nabla E$ and $-\nabla E|_{L_p}$ on $L_p$. Let us parametrize the set $\lbrace F_1=0\rbrace\cap L_p$ by
$$y_2(x_2)=\frac{1}{2}(-x_2^2\pm x_2\sqrt{4+x_2^2}),$$
and the set $\lbrace F_2=0\rbrace\cap L_p$ by
$$x_2(y_2)=\frac{\mp y_2}{\sqrt{1+2 y_2}}.$$

To understand $-\nabla E$ along the boundary of $L_p$ in some neighborhood of $p$, we use the same condition as in Lemma \ref{lem_grad_standard}. Using \textit{Mathematica} we compute that
\begin{align*}
\widetilde{\langle \widetilde{\nabla}F_1, -\nabla E\rangle}|_{Q_1}=-x_2^3+O(x_2^4)<0,\\
\widetilde{\langle \widetilde{\nabla}F_1, -\nabla E\rangle}|_{Q_2}=x_2^3+O(x_2^4)<0
\end{align*}
along $\lbrace F_1=0\rbrace\setminus\lbrace(0, 0, 0, 0)\rbrace$. And also
$$\widetilde{\langle \widetilde{\nabla}F_2, -\nabla E\rangle}|_{Q_1\cup Q_2}=-y_2^3+O(y_2^4)<0$$
along $\lbrace F_2=0\rbrace\setminus\lbrace(0, 0, 0, 0)\rbrace$.

Hence, if $\widetilde{V}_p^{\circ}$ is some small collar neighborhood of $p$, then $-\nabla E$ is strictly outward-pointing from $M^{out-out}$ along the set $L_p\cap \widetilde{V}_p^{\circ}\cap\lbrace F_1=F_2=0\rbrace.$

It is interesting to consider also the semi-fixed case, i.e., $-\nabla E|_{L_p}$. We then obtain that 
\begin{align*}
\widetilde{\langle \widetilde{\nabla}F_1, -\nabla E|_{L_p}\rangle}|_{Q_1}=x_2^3+O(x_2^4)>0,\\
\widetilde{\langle \widetilde{\nabla}F_1, -\nabla E|_{L_p}\rangle}|_{Q_2}=-x_2^3+O(x_2^4)>0
\end{align*}
along $\lbrace F_1=0\rbrace\setminus\lbrace(0, 0, 0, 0)\rbrace$. And also
$$\widetilde{\langle \widetilde{\nabla}F_2, -\nabla E|_{L_p}\rangle}|_{Q_1\cup Q_2}=-2y_2^3+O(y_2^4)<0$$
along $\lbrace F_2=0\rbrace\setminus\lbrace(0, 0, 0, 0)\rbrace$. So now, the negative gradient is strictly outward-pointing only from the boundary corresponding to $F_2=0$. To the last, note that the dot product condition for $F_i$ gave twice as big value and with the opposite sign when restricted to $(x_i, y_i)$, then when restricted to $(x_j, y_j)$ for $j\neq i$.
\end{example}

\begin{rem}We expect that over any weakly $\overline{s}_2$-special set the following holds. There is a point $x_\varepsilon\in M_{K, \varepsilon}$ such that $F^{[\varepsilon]}_1(x_\varepsilon)>0\vee F^{[\varepsilon]}_2(x_\varepsilon)=0$ and $-\nabla E_\varepsilon(x_\varepsilon)$ is now inward-pointing into $M_{K, \varepsilon}$ In other words, the relative direction of $-\nabla E_\varepsilon(x_\varepsilon)$ is opposite to the direction in Lemma \ref{lem_grad_standard} $(i.)$. For geometric intuition, the reader might find it convenient to see Figure \ref{fig_fibers}. We expect the analogous phenomena also for weakly $\overline{s}_1$-special sets.
\end{rem}

\begin{defn} Let $X$ be a $C^1$ vector field on the Riemannian manifold $(M, g)$ that vanishes at a point $p$. Let $x\in M\setminus p$. We say that $\ell_p\in T_p M$ is an \textbf{asymptotic limit vector of $x$} iff there is a partial flow line $u:[0, \infty)\rightarrow M$ of $X$ such that
\begin{itemize}
\item[$(i.)$] $u(0)=x$.
\item[$(ii.)$] $\lim_{t\rightarrow\infty}u(t)=p.$
\item[$(iii.)$] $\lim_{t\rightarrow\infty}\frac{\dot{u}(t)}{||\dot{u}(t)||}=\ell_p.$
\end{itemize}
\end{defn}

\begin{rem}\cite[Lemma B.5]{Schwarz1993MorseHT} Let $f$ be a Morse function on $(M, g)$. Then for the gradient flow $\phi^t_f$ it holds that every point $x\in M\setminus Crit(f)$ has an asymptotic limit vector.
\end{rem}

\begin{Sternberg}\cite{Nelson_dynamics}\label{thm_sternberg}
Let $n\geq2$ and $X$ be a smooth vector field on $\R^n$ that vanishes at $0$. Let $\chi_1,\dots,\chi_n$ be (complex) eigenvalues of $DX(0)$ that satisfy the following non-resonant assumption:
\begin{itemize}
\item For $i=1,\dots, n$ and all $m_j$ non-negative integers such that $2\leq\sum_{j=1}^n m_j$ holds
$$\chi_i\neq\sum_{j=1}^n m_j\chi_j.$$
\end{itemize}
Then there is a local smooth diffeomorphism around $0$ conjugating $X$ with its linear part on $T_0\R^n$.
\end{Sternberg}

\begin{lemma}\label{lemma_assympt_diag} Let $p\in \Delta_{full}$ and $\ell_p$ be a unit vector of $T_p(\R/T\mathbb{Z}\times S^1)^2.$ Let $\nu_\delta(\Delta_{full})$ be a $\delta$-radius tubular neighborhood of $\Delta_{full}$ for some $\delta>0$. If $\delta>0$ is small, then for the flow $\phi^t_{E_\varepsilon}$ on $(\R/T\mathbb{Z}\times S^1)^2$ it holds that the set
\begin{equation*} 
\lbrace x_\varepsilon\in \nu_\delta(\Delta_{full})\,|\,\ell_p\hbox{ is the asymtotic limit vector of }x_\varepsilon\rbrace
\end{equation*}
is flow invariant and diffeomorphic to an open interval.
\end{lemma}

\begin{proof}
By Lemma \ref{lem_ener_bott} $\Delta_{full}$ is a Bott nondegenerate critical manifold for $E_\varepsilon$. Let $\widehat{W}^s_{E_\varepsilon}(\Delta_{full})$ be the restriction of $W^s_{E_\varepsilon}(\Delta_{full})$ to $\nu_\delta(\Delta_{full})$. By Remark \ref{rem_stable_dim} the flow endpoint map $\pi_{E_\varepsilon}:\widehat{W}^s_{E_\varepsilon}(\Delta_{full})\rightarrow \Delta_{full}$ induce a locally trivial fiber bundle.

Let $p\in\Delta_{full}$. We consider the fiber $\mathcal{F}=\pi^{-1}_{E_\varepsilon}(p)$ and the flow of $-\nabla E_\varepsilon\vert_{\mathcal{F}}$. In particular, $\mathcal{F}$ is flow invariant.  By Lemma \ref{lem_ener_bott} the spectrum of $-D\nabla E_\varepsilon\vert_{\mathcal{F}}(p)$ consists of two eigenvalues equal to $2$. Hence if $\delta>0$ is small, then by Sternberg's Linearization Theorem \ref{thm_sternberg} $-\nabla E_\varepsilon\vert_{\mathcal{F}}$ is conjugated by a smooth diffeomorphism with $-D\nabla E_\varepsilon\vert_{\mathcal{F}}$ (restricted to some small neighborhood $U_0$ of $0$ in $T_p\mathcal{F}$).

Now, by the spectrum of $-D\nabla E_\varepsilon\vert_{\mathcal{F}}(0)$ the following holds for the flow of $-D\nabla E_\varepsilon\vert_{\mathcal{F}}$. If $\ell_0\in T_0(T_p\mathcal{F})$ is a unit vector, then the set of $x\in U_0$, with the asymptotic limit vector $\ell_0$, is diffeomorphic to an open interval. 

Since the flows are conjugated by a smooth diffeomorphism, the lemma follows. 
\end{proof}

\begin{conj}\label{conj_eating_cusp} For every $\delta>0$ small there is a subset $\Delta_{\varepsilon, \delta}\subset\Delta_{\varepsilon}^{cusp}$ which is $O(\delta)$-close and diffeomorphic to $\Delta_{\varepsilon}^{cusp}$ and the following holds.

If $p\in\Delta_{\varepsilon, \delta}$, then the set
$$
A
_p=\left\{ x_\varepsilon\in\nu_\delta(\Delta_{\varepsilon, \delta})\setminus\Delta_{\varepsilon, \delta} \ \middle\vert \begin{array}{l}
    x_\varepsilon\cdot_{E_\varepsilon}[0, \infty)\subset M_{K, \varepsilon}, \\
    p\hbox{ is the omega limit of }x_\varepsilon
  \end{array}\right\}
$$
is flow invariant and diffeomorphic to a disjoint union of two open intervals. See also Figure \ref{figure_eat_trajectory}.
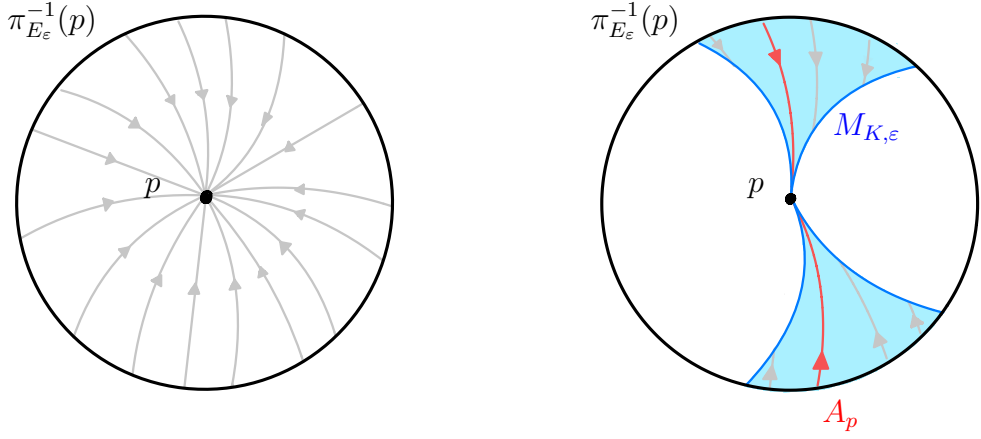
\begin{figure}[!htbp]
\labellist
\pinlabel $p$ at 80 744
\pinlabel $\pi_{E_\varepsilon}^{-1}(p)$ at 30 828
\pinlabel $p$ at 380 744
\pinlabel $\pi_{E_\varepsilon}^{-1}(p)$ at 320 828
\pinlabel $\textcolor{blue}{M_{K, \varepsilon}}$ at 435 774
\pinlabel $\textcolor{red}{A_p}$ at 422 630
\endlabellist
\centering
\includegraphics[scale=0.75]{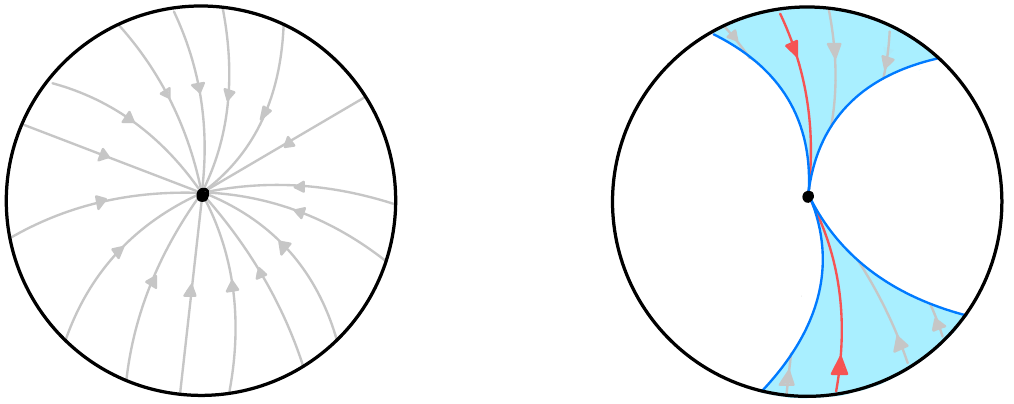}
\vspace{0.3cm}
\caption{\textit{On the left:} the fiber $\pi_{E_\varepsilon}^{-1}(p)$ of the fibration $\widehat{W}^s_{E_\varepsilon}(\Delta_{full})\xrightarrow{\pi_{E_{\varepsilon}}}\Delta_{full}$ inside $\nu_\delta(\Delta_{full})$. The gradient flow trajectories of $-\nabla E_\varepsilon$ are in grey. \textit{On the right:} the restriction of the gradient flow trajectories from $\pi_{E_\varepsilon}^{-1}(p)$ to $\pi_{E_\varepsilon}^{-1}(p)\cap M_{K, \varepsilon}$. The set $A_p$ is depicted in red.}
\label{figure_eat_trajectory}
\end{figure}
\end{conj}

\begin{pproof}Inside the tubular neighborhood $\nu_\delta(\Delta_{full})$ we have two (stratified) fiber bundles $\widehat{W}^s_{E_\varepsilon}(\Delta_{full})\xrightarrow{\pi_{E_{\varepsilon}}}\Delta_{full}$ and $\widehat{M}_{K, \varepsilon}\xrightarrow{\pi_\mathcal{N}}\Delta_\varepsilon^{cusp}
$ given by Remark \ref{rem_stable_dim}, Lemma \ref{conj_local_fibr_diag} and Corrollary \ref{lem_cusp_fibre}. Recall that $\pi_{E_\varepsilon}$ is the flow end-point map and $\pi_\mathcal{N}$ is induced from the normal bundle $\mathcal{N}_{\bm{\Gamma}^{\ast}\langle\cdot, \cdot\rangle}(\Delta^{cusp}_\varepsilon, (\R/T\mathbb{Z}\times S^1)^2)$.

Now, let us assume that $p\in \Delta_\varepsilon^{cusp}$ is such that 
\begin{equation}\label{eqn_fit_fibre}
\pi_{\mathcal{N}}\circ\pi_{E_\varepsilon}^{-1}(p)\subset\Delta_\varepsilon^{cusp}.
\end{equation}
So, if $\delta>0$ is small, then by Conjecture \ref{conj_grad_diag} and Corollary \ref{lem_cusp_fibre} the set $A_p$ consists of points with one of two possible asymptotic limit vectors at $p$. Hence, by Lemma \ref{lemma_assympt_diag}, it holds that $A_p$ is diffeomorphic to a disjoint union of two open intervals.

It remains to show that the relation (\ref{eqn_fit_fibre}) can be satisfied by a sufficiently large subset $\Delta_{\varepsilon, \delta}\subset\Delta_\varepsilon^{cusp}$. By Lemma \ref{conj_local_fibr_diag} we see that any generalized tangent space of $\pi_{\mathcal{N}}^{-1}(p)$ at $p$ is contained at $T_p(\pi_{E_\varepsilon}^{-1}(p))$. Finally, as we shrink $\delta$-radius of $\nu_\delta(\Delta_{full})$, $\Delta_{\varepsilon, \delta}$ can be arbitrarily large as a subset of $\Delta_\varepsilon^{cusp}$. 
\end{pproof}

\chapter{Adiabatic limit with Conley index}
\label{ch:adiab_conley}
From Section \ref{s:energy} we know the correspondence between the critical points $Crit(E_0)$ and $Crit(E_\varepsilon\vert_{M_{K, \varepsilon}})$ and their degree shift. Hence, it is natural to ask about the correspondence between the trajectories on $K\times K$ and $M_{K, \varepsilon}$, at least for $\varepsilon>0$ small. However, as we shrink $T_{K, \varepsilon}$ to $K$, the gradient $\nabla E_\varepsilon$ explodes. Hence, due to the singular behavior, we would like to apply an Adiabatic limit argument as in \cite{frauenfelder2022lagrange, Frauenfelder2022LagrangeMA}, see also \cite{webber99}. In order to avoid a lot of analysis, we would like to perform instead of two Implicit function theorems rather a certain Conley index argument. As a consequence, such a correspondence will not identify the particular trajectories, but rather only certain equivalence classes of trajectories. However, we will be able to close each equivalence class of trajectories to an arbitrarily thin strip (isolating neighborhood), and in particular know reasonably well the ``shape'' of the trajectories. In addition, these thin strips and a certain argument by contradiction will be used to trap the trajectories and obtain the adiabatic compactness. 

We remark that in our argument, we want the trajectories on $K\times K$ to generically avoid special and diagonal points. For this, we restrict the adiabatic limit only for trajectories between critical points in low degrees that moreover do not belong to diagonals.

To the end, we note that in the literature, some Conley index arguments were used for the singular dynamical systems, see for example \cite{Gedeon1999TheCI}. However, to the author\textsc{\char13}s knowledge, only in the task of the existence of orbits. This is probably due to the fact that our dynamics is given by gradient(-like) vector fields instead of some general dynamical system.

\section{Gradient-like vector fields}
We are going to introduce gradient-like vector fields and discuss some of their perturbation techniques.

\begin{defn}\label{defn_grad_like}Let $f$ be a Morse-Bott function on closed oriented Riemannian manifold $(M, g)$. A (smooth) vector field $X_f$ is called \textbf{gradient-like} adapted to $f$ if it holds that
\begin{itemize}
\item[$(i.)$]$df(x)X_f|_x<0$ for every $x\notin Crit(f)$.
\item[$(ii.)$] $-\nabla^g f= X_f$ in some open neighborhood of $Crit(f)$.
\end{itemize}
The flow of $X_f$ will be denoted as $\phi^t_{X_f}$ or $\cdot_{X_f}$. The stable/unstable manifolds of $Q\subset Crit(f)$ will be denoted as $W^{s/u}_{X_f}(Q)$. Moreover, 
$$W_{X_f}(Q_i, Q_j)=W^u_{X_f}(Q_i)\cap W^s_{X_f}(Q_j),$$
where $Q_i, Q_j\subset Crit(f)$.

In addition, we put 
$$\mathcal{M}_{X_f}(Q_i, Q_j)=W_{X_f}(Q_i, Q_j)/\R,$$
where the quotient is given by the free action of $\R$, which shifts the time on the flow lines.
\end{defn}

\begin{rem} Due to the condition $(i.)$ in Definition \ref{defn_grad_like}, we can see any gradient-like vector field as a negative gradient vector field with respect to some different metric. In more detail, let $X_f$ be a gradient-like vector field adapted to a Morse-Bott function $f$ on $(M, g)$. Then there is a Riemannian metric $\widehat{g}$ on $M$ such that $X_f=-\nabla^{\widehat{g}}f$.

To see this, one has to find a suitable metric locally in nice charts and then use the partition unity argument. Around $Crit(f)$, we already have a good metric. Outside of $Crit(f)$ we can always take local coordinates $(x_0,\dots, x_n)$ such that $f$ is a constant in coordinates $(x_1,\dots, x_n)$. Then, we define point-wise a matrix $A$ by columns which consist of bases vectors $\partial_{x_1},\dots,\partial_{x_n}$ and $-X_{f}/\sqrt{|df(X)|}$. Finally, $(AA^T)^{-1}$ will define point-wise the desired local metric.
\end{rem}

\begin{defn}\cite{hurtubise2013approaches}Let $f:M\rightarrow\R$ be a Morse-Bott function and $X_f$ be a gradient-like vector field. Then the vector field $X_f$ is called \textbf{Morse-Bott-Smale} if for any two $Q_{i,j}\in Crit(f)$ it holds that
\begin{equation}\label{eqn_Morse_Bott_Smale}
W^u_{X_f}(q)\pitchfork W^s_{X_f}(Q_j),
\end{equation}
for all $q\in Q_i$. Specially, if $f$ is Morse and $X_f$ satisfy $(\ref{eqn_Morse_Bott_Smale})$, then $X_f$ is called \textbf{Morse-Smale}.
\end{defn}

\begin{rem}\label{rem_Morse_smale_same_ind}Let $X_f$ be a Morse-Smale vector field adapted to a Morse function $f$ on $(M, g)$ and $p, q\in Crit(f)$. Then by Remark \ref{rem_stable_dim} and Theorem \ref{thm_invers_im} $W_{X_f}(p, q)$ is an (embedded) submanifold of $M$ and $\dim(W_{X_f}(p, q))=Ind_f(p)-Ind_f(q)$.

Note that if moreover $W_{X_f}(p, q)\neq \emptyset$, then $W_{X_f}(p, q)$ consists of at least one gradient flow line, and thus $\dim(W_f(p, q))>0$ . Hence, there are no gradient flow lines between the critical points of the same Morse index.
\end{rem}

\begin{rem} Let $f$ be a Morse function on the $n$-dimensional manifold $M$ and $q_0, q_1, q_{n-1}, q_n$ be its critical points with Morse indices equal to the subscripts. Assume that $X_f$ is a Morse-Smale vector field adapted to $f$. Then the manifolds $W_{X_f}(q_1, q_0), W_{X_f}(q_{n-1}, q_n)$ both contain up to $2$ (unparametrized) gradient-like flow lines.
\end{rem}

\begin{rem} If two stable and unstable manifolds intersect transversely at some point $x\in M$, then they intersect transversely along the whole flow line containing $x$. 
\end{rem}

\begin{lemma}\cite{Smale1961OnGD}\label{lemma_kupka_aux}Let $f$ be a Morse-Bott function on closed oriented Riemannian manifold $(M, g)$. Let $X_f$ be a gradient-like vector field adapted to $f$. Let $Q\in Crit(f)$ such that $f(Q)\neq\max(f)$ and $U_Q$ be an open neighborhood of $Q$ such that $U_Q\cap Crit(f)=Q$. If $\delta>0$ is small, then $W^s_{X_f}(Q)\cap U_Q\cap f^{-1}\big((f(Q)+\delta, f(Q)+2\delta)\big)\neq\emptyset$ and there is a vector field $Y$ such that
\begin{itemize}
\item[$(i.)$] $Y$ is compactly supported in $U_Q\cap f^{-1}\big((f(Q)+\delta, f(Q)+2\delta)\big)$.
\item[$(ii.)$] $Y$ is $C^1$ $O(\delta)$-close to the zero section $0_M$ of $TM$.
\item[$(iii.)$] For any $Q_j\in Crit(f)$ it holds that
$$W^u_{X_f+Y}(Q_j)\pitchfork W^s_{X_f+Y}(Q).$$
\end{itemize} 
\end{lemma}

\begin{KupkaSmale}\cite{Smale1961OnGD}\label{thm_kupka} Let $X_f$ be a gradient-like vector field adapted to the Morse-Bott function $f$ on the closed oriented Riemannian manifold $(M, g)$. Let $y_m$ be the minimum value of $f$. Assume that $Crit(f)\setminus f^{-1}(y_m)$ consists of $0$-dimensional critical manifolds. Then there exists a Morse-Bott-Smale vector field $Y_f$ that is $C^1$ close to $X_f$.
\end{KupkaSmale}

\begin{proof}
Recall that (nontrivial) gradient-like flow lines flow ``downhill'' with respect to the level sets of the Morse function. Hence, we can partially order $Crit(f)$ by the heights of critical values (from the highest to the lowest). We aim to inductively apply Lemma \ref{lemma_kupka_aux}. For this, we will reorder the critical manifolds with the same critical values by different choices of small $\delta>0$, so that the perturbations from Lemma \ref{lemma_kupka_aux} will be made in different heights. The points of $Crit(f)$ with maximal height can be reordered arbitrarily, since for them we do not need to call Lemma \ref{lemma_kupka_aux}. So the set $Crit(f)$ is now ordered and ready for inductive application of Lemma \ref{lemma_kupka_aux}.

To the end, we point out that since only the higher-dimensional critical manifolds lie in $f^{-1}(y_m)$, the Morse-Bott-Smale transversality at the unstable manifolds of these critical manifolds is automatic.
\end{proof}

\begin{rem}Presented Lemma \ref{lemma_kupka_aux} and Kupka-Smale Theorem \ref{thm_kupka} slightly extend the work \cite{Smale1961OnGD} which was done for Morse-Smale flows adapted to Morse functions with distinct critical values. For another proof, see also \cite{Palis1982GeometricTO, Peixoto1967OnAA}.
\end{rem}

\begin{lemma}\cite{petrak2019definition}\label{lem_perturb_petrak}Let $X_{E_0}$ be a gradient-like vector field. Next, let $p_0\in Crit_1(E_0)$ and $x\in(\R/T\mathbb{Z})^2$ be a special point. Let $U_x$ be an open neighborhood of $x$. Then there is a vector field $Y$ such that
\begin{itemize}
\item[$(i.)$]$Y$ is compactly supported on $U_x.$ 
\item[$(ii.)$]$Y$ is $C^1$ close to the zero section $0_{(\R/T\mathbb{Z})^2}$ of $T(\R/T\mathbb{Z})^2$.
\item[$(iii.)$]$W^u_{X_{E_0}+Y}(p_0)\pitchfork x$, i.e. there is no intersection of $W^u_{X_{E_0}+Y}(p_0)$ and $x$.
\item[$(iv.)$]$W^u_{X_{E_0}}(p_0)=W^u_{X_{E_0}+Y}(p_0)$ on $(\R/T\mathbb{Z})^2\setminus U_x.$
\end{itemize}
\end{lemma}
\begin{proof}
Note that in Lemma \ref{lemma_kupka_aux}, we made a local perturbation of the flow that globally deforms the stable and unstable manifolds. However, here we demand the property $(iv.)$, so now the perturbation technique will be a bit different than the perturbation in Lemma \ref{lemma_kupka_aux}.

First, we make a local compact perturbation of the $1$-dimensional manifold $W^u_{X_{E_0}}(p_0)$ such that the transversality with $x$ is achieved. This can be done by Relative Thom's transversality Theorem, see \cite[thm A.4]{petrak2019definition}. Then we find a new vector field $Y$ such that the perturbed $1$-manifold is an integral of $Y$, see \cite[lem 2.19]{petrak2019definition}. If the perturbation was small enough, then $Y$ is also a gradient-like vector field adapted to $E_0$. 
\end{proof}

In this chapter, we will consider the following \textbf{generic approximations} of $-\nabla E_0$ and $-\nabla E_\varepsilon$.

\begin{cor}\label{cor_perturb_grad_like}There is a $C^1$ small vector field $Y$ on $(\R/T\mathbb{Z})^2$ such that the vector field 
$$X_{E_0}:=-\nabla E_0+Y$$
is Morse-Bott-Smale and for every $p_0\in Crit(E_0)$ it holds that $W^u_{X_{E_0}}(p_0)$ is transverse to all special points.

Let $\varepsilon\in(0, \varepsilon_{good}]$ and $Y_0$ be a canonical lift of $Y$ to $(\R/T\mathbb{Z})\times S^1)^2$, i.e. $(\pi_{s_1, s_2})_\ast Y_0=Y$ and $(\pi_{\theta_1, \theta_2})_\ast Y_0=0$. Then there is a $C^1$ $O(\varepsilon^2)$-small vector field $Z_\varepsilon$ on $(\R/T\mathbb{Z}\times S^1)^2$ such that the vector field
$$X_{E_\varepsilon}:=-\nabla E_\varepsilon+Y_0+Z_\varepsilon$$
is Morse-Bott-Smale and $X_{E_\varepsilon}$ and $-\nabla E_\varepsilon+Y_0$ agree outside of a $O(\varepsilon^2)$-small neighborhood of $Crit(E_\varepsilon)$.
\end{cor}

\begin{proof}
Recall that by Lemma \ref{lemma_tangent_gen} there is only a finite number of special points. Moreover, these points are clearly not critical points of $E_0$ (or inside $\Delta_0)$. In order to obtain $X_{E_0}$, we would like to inductively apply to $-\nabla E_0$ the finite number of perturbations from the Kupka-Smale theorem \ref{thm_kupka} and Lemma \ref{lem_perturb_petrak}.

We consider the ordered set of $Crit(E_0)$ from the Kupka-Smale Theorem \ref{thm_kupka} and augment the set with the special points. The special points will be ordered using their values in $E_0$. The resulting set $S$ is now partially ordered. In particular, if for a special point $x$ and a critical point $p_0\in Crit(E_0)$ holds that $E_0(x)=E(p_0)\neq\hbox{max}(E_0)$, then we chose a convention $x<p_0$. This is due to the fact that the perturbation for $p_0$ will be made at the level $(f(p_0)+\delta, f(p_0)+2\delta)$ and not in $(f(p_0)-\delta, f(p_0)+\delta)$ for some $\delta>0$. If $E_0(x)=E(p_0)=\hbox{max}(E_0)$, then we reorder by $x<p_0$. The special points with the same values of $E_0$ can be reordered among themselves arbitrarily. So the set $S$ is now ordered.

Let us consider a step of the induction, where we want to perform a small local perturbation around a special point $x$. By taking a sufficiently small neighborhood $U_x$, we can achieve that this perturbation deforms (locally) only the unique $1$-dimensional unstable manifold that intersects the special point, and that this perturbation is not made in the same levels as the perturbations for the critical points. In particular, the perturbed unstable manifold will have no intersection with any special or critical point on the same $E_0$-level as $x$. Since the gradient-like flows flow downhill, the induction will produce the desired $X_{E_0}$.

Now, the perturbation of $-\nabla E_\varepsilon+Y_0$ is an immediate consequence of the Kupka-Smale Theorem \ref{thm_kupka}.
\end{proof}

\begin{defn}Let $X_{E_0}$ be a gradient-like vector field and $p_0, q_0\in Crit(E_0)$. A gradient-like flow line $u$ of the flow $\phi^t_{X_{E_0}}$ is called a \textbf{$0$-solution from $p_0$ to $q_0$} if $u\subset W_{X_{E_0}}(p_0, q_0)$.

Let $\varepsilon\in(0, \varepsilon_{good}]$, $X_{E_\varepsilon}$ be a gradient-like vector field and $p_\varepsilon, q_\varepsilon\in Crit(E_\varepsilon)$. A gradient-like flow line $u_\varepsilon$ of the flow $\phi^t_{X_{E_\varepsilon}}$ is called an \textbf{$\varepsilon$-solution from $p_\varepsilon$ to $q_\varepsilon$} if $u_\varepsilon\subset W_{X_{E_\varepsilon}}(p_\varepsilon, q_\varepsilon)$.
\end{defn}

\section{Uniform bounds}
In this section, we provide uniform bounds for the $\varepsilon$-trajectories that will be needed for the adiabatic compactness.

\begin{thm}\label{thm_unif_bound}
There is a $\varepsilon_0\in(0, \varepsilon_{good}]$ and positive constants $\lbrace c_j\rbrace_{j=1,\dots, 8}$ such that for every $\varepsilon\in(0, \varepsilon_0]$ the following holds.

Let $p_\varepsilon\in Cr_2(E_\varepsilon)$ and $q_\varepsilon\in Cr_1(E_\varepsilon)$ such that $p_\varepsilon, q_\varepsilon\in M_{K, \varepsilon}\setminus \Delta_\varepsilon$. Let $X_{E_\varepsilon}$ be a generic approximation of $-\nabla E_\varepsilon$ (see Corollary \ref{cor_perturb_grad_like}). Let $u_\varepsilon(t)$ be an $\varepsilon$-solution from $p_\varepsilon$ to $q_\varepsilon$ which is described in coordinates as $$u_\varepsilon(t)=(s_1(t), \theta_1(t), s_2(t), \theta_2(t)).$$
Then there are bounds
$$|\dot{s}_1(t)|\leq c_1, |\dot{s}_2(t)|\leq c_2, |\ddot{s}_1(t)|\leq c_3, |\ddot{s}_2(t)|\leq c_4$$
for $t\in \R$.

Moreover, let us assume that $u_\varepsilon$ avoids weakly special and weakly diagonal points, i.e., $\pi_{s_1, s_2}\circ u_\varepsilon\subset S_K$. Then $u_\varepsilon\subset M_{K, \varepsilon}$. And, in addition, there are bounds
$$|\dot{\theta}_1(t)|\leq c_1, |\dot{\theta}_2(t)|\leq c_2, |\ddot{\theta}_1(t)|\leq c_3, |\ddot{\theta}_2(t)|\leq c_4$$
for $t\in \R$.
\end{thm}

\begin{rem}We are going to show Theorem \ref{thm_unif_bound} only for $-\nabla E_\varepsilon$. Then, provided that $\varepsilon>0$ is sufficiently small, Theorem \ref{thm_unif_bound} will also hold for any generic approximation of $-\nabla E_\varepsilon$ which is given by Corollary \ref{cor_perturb_grad_like}. 
\end{rem}

\begin{rem} The proof of Theorem \ref{thm_unif_bound} will be the aim of the rest of this section. But first we would like to enlight the main idea of the proof. We would like to construct certain more and more delicate ``Lyapunov-like traps'' for $u_\varepsilon$. Unlike for the standard level sets of a Lyapunov function, in our case, a general flow line could potentially escape from the trap (and even flow back). However, due to the global geometry of $T_{K, \varepsilon}\times T_{K, \varepsilon}$, we will see that $u_\varepsilon$ will not be able to escape. 
\end{rem}

\begin{rem} Let $\varepsilon\in(0, \varepsilon_{good}]$. With the help of (\ref{eqn_grad_E}), we are going to describe in coordinates the first two derivatives of the components of $u_\varepsilon$. 

The first derivatives are immediately given by
\begin{align}
\dot{s}_1&=-d_1^{-1}\langle P+\varepsilon v_2, \dot{\gamma}(s_1)\rangle,\label{eqn_unif_1}\\
\dot{s}_2&=d_2^{-1}\langle P-\varepsilon v_1, \dot{\gamma}(s_2)\rangle,\label{eqn_unif_2}\\
\varepsilon\dot{\theta_1}&=\varepsilon\tau(s_1)d_1^{-1}\langle P+\varepsilon v_2, \dot{\gamma}(s_1)\rangle-\langle P+\varepsilon v_2, v_1^{\bot}\rangle,\label{eqn_unif_3}\\
\varepsilon\dot{\theta_2}&=-\varepsilon\tau(s_2)d_2^{-1}\langle P-\varepsilon v_1, \dot{\gamma}(s_2)\rangle+\langle P-\varepsilon v_1, v_2^{\bot}\rangle,\label{eqn_unif_4}
\end{align}
since $u_\varepsilon$ is an integral of $-\nabla E_\varepsilon$. Recall that $d_i=1-\varepsilon\cos(\theta_i)\kappa(s_i).$

The second derivatives will be described with the help of certain (bounded) functions $\lbrace f_{k, \varepsilon}:(\R/T\mathbb{Z}\times S^1)^2\rightarrow\R\rbrace_{k=1,\dots, 4}$. We point out that since we are in coordinates, we are computing the derivatives of $\dot{u}_\varepsilon$ with respect to the flat connection (not $\nabla^g$).
\begin{align}
\begin{split}\label{eqn_unif_5}
\ddot{s}_1&=-\partial_{\theta_1}(d_1^{-1})\langle P+\varepsilon v_2, \dot{\gamma}(s_1)\rangle\dot{\theta}_1
-\varepsilon d_1^{-1}\langle v_2^{\bot}, \dot{\gamma}(s_1)\rangle\dot{\theta}_2\\
&\quad\,-\partial_{s_1}(d_1^{-1})\langle P+\varepsilon v_2, \dot{\gamma}(s_1)\rangle\dot{s}_1
-d_1^{-1}\big(1+\kappa(s_1)\langle P+\varepsilon v_2, n(s_1)\rangle\big)\dot{s}_1\\
&\quad\,-d_1^{-1}\langle d_2\dot{\gamma}(s_2)+\varepsilon\tau(s_2)v_2^{\bot}, \dot{\gamma}(s_1)\rangle\dot{s}_2\\
&=-\frac{\sin(\theta_1)\kappa(s_1)}{d_1^2}\langle P+\varepsilon v_2, \dot{\gamma}(s_1)\rangle\big(\varepsilon\tau(s_1)d_1^{-1}\langle P+\varepsilon v_2, \dot{\gamma}(s_1)\rangle-\langle P+\varepsilon v_2, v_1^{\bot}\rangle\big)\\
&\quad\,-\varepsilon d_1^{-1}\langle v_2^{\bot}, \dot{\gamma}(s_1)\rangle\big(-\varepsilon\tau(s_2)d_2^{-1}\langle P-\varepsilon v_1, \dot{\gamma}(s_2)\rangle+\langle P-\varepsilon v_1, v_2^{\bot}\rangle\big)\\
&\quad\,+\frac{\varepsilon\cos(\theta_1)\dot{\kappa}(s_1)}{d_1^3}\langle P+\varepsilon v_2, \dot{\gamma}(s_1)\rangle\langle P+\varepsilon v_2, \dot{\gamma}(s_1)\rangle\\
&\quad\,-d_1^{-1}(1+\kappa(s_1)\langle P+\varepsilon v_2, n(s_1)\rangle)\dot{s}_1\\
&\quad\,-d_1^{-1}\langle d_2\dot{\gamma}(s_2)+\varepsilon\tau(s_2)v_2^{\bot}, \dot{\gamma}(s_1)\rangle\dot{s}_2\\
&=:f_{1, \varepsilon}(s_1, \theta_1, s_2, \theta_2),
\end{split}
\end{align} 
analogously also 
\begin{equation}\label{eqn_unif_6}
\ddot{s}_2=:f_{2, \varepsilon}(s_1, \theta_1, s_2, \theta_2).
\end{equation}
Moreover,
\begin{align}
\begin{split}
\varepsilon\ddot{\theta}_1&=\frac{d}{dt}\big[-\varepsilon\tau(s_1)\dot{s}_1-\langle P+\varepsilon v_2, v_1^\bot\rangle\big]\\
&=-\varepsilon\dot{\tau}(s_1)(\dot{s}_1)^2-\varepsilon\tau(s_1)\ddot{s}_1-\varepsilon\frac{d}{dt}\langle v_2, v_1^\bot\rangle-\frac{d}{dt}\langle P, v_1^\bot\rangle\\
&=\varepsilon\Big[-\dot{\tau}(s_1)(\dot{s}_1)^2-\tau(s_1)\ddot{s}_1-\frac{d}{dt}\langle v_2, v_1^\bot\rangle\Big]+\widetilde{\langle\widetilde{\nabla}\langle P, v_1^\bot\rangle, \nabla E_\varepsilon\rangle}\\
&=:f_{3, \varepsilon}(s_1, \theta_1, s_2, \theta_2)+\widetilde{\langle\widetilde{\nabla}\langle P, v_1^\bot\rangle, \nabla E_\varepsilon\rangle}\\
&=f_{3, \varepsilon}(s_1, \theta_1, s_2, \theta_2)\\
&\quad\,+\varepsilon\Big[-\partial_{s_1}\langle P, v_1^\bot\rangle\langle v_2, \dot{\gamma}_1\rangle d_1^{-1}+\partial_{s_2}\langle P, v_1^\bot\rangle\langle v_1, \dot{\gamma}_2\rangle d_2^{-1}\Big]\\
&\quad\,+\frac{d_2\sin(\theta_1)\kappa(s_1)\langle P, \dot{\gamma}(s_1)\rangle^2-d_1d_2\langle P, v_1\rangle\langle v_2, v_1^\bot\rangle-d_1\langle\dot{\gamma}(s_2), v_1^\bot\rangle\langle  \dot{\gamma}(s_2), v_1\rangle}{d_1d_2}\\
&\quad\,-\frac{\langle P, v_1\rangle\langle P, v_1^\bot\rangle}{\varepsilon}
\end{split}
\end{align}
And also analogously, we have that
\begin{align}\label{eqn_unif_8}
\varepsilon\ddot{\theta}_2=:f_{4, \varepsilon}(s_1, \theta_1, s_2, \theta_2)-\widetilde{\langle\widetilde{\nabla}\langle P, v_2^\bot\rangle, \nabla E_\varepsilon\rangle},
\end{align}
where $-\widetilde{\langle\widetilde{\nabla}\langle P, v_2^\bot\rangle, \nabla E_\varepsilon\rangle}$ also splits in the terms $O(\varepsilon)+O(1)+O(1/\varepsilon)$. Specially, $O(1/\varepsilon)$ terms consists of $$\frac{\langle P, v_2\rangle\langle P, v_2^\bot\rangle}{\varepsilon}.$$
\end{rem}

\begin{lemma}\label{lem_bound_first} There is $\varepsilon_0$ such that $c_1$ and $c_2$ exist.
\end{lemma}
\begin{proof}
This follows from equations (\ref{eqn_unif_1}) and (\ref{eqn_unif_2}). For $i=1, 2$, let us consider a family of functions $\lbrace h_{i, \varepsilon}:=(\pi_{s_i})_\ast(-\nabla E_\varepsilon):(\R/T\mathbb{Z}\times S^1)^2\rightarrow\R\rbrace_{\varepsilon\in(0, \varepsilon_{good}]}$. 

Note that each $h_{i, \varepsilon}$ is a smooth function on a compact set, and hence bounded. Also $\lim_{\varepsilon\rightarrow 0}d_i^{-1}=1$. In particular, the smooth family $\lbrace h_{i, \varepsilon}\rbrace_{\varepsilon\in(0, \varepsilon_{good}]}$ smoothly extends to $\varepsilon=0$. By the compactness of the domain $(\R/T\mathbb{Z}\times S^1)^2$ and smoothness of the family, it will be enough to verify the uniform bounds for $|h_{i, 0}|$. Then there will be $\varepsilon_0>0$ small such that the uniform bounds hold also for each $|h_{i, \varepsilon}|$ with $\varepsilon\in[0, \varepsilon_0]$. And in particular, we will obtain the uniform bounds for $|\dot{s}_i|$. Since the uniform bound for $|h_{i, 0}|$ is straightforward, the lemma follows.
\end{proof}

\begin{lemma} There is $\varepsilon_0$ such that $c_3$ and $c_4$ exist.
\end{lemma}

\begin{proof}
Let us see the relations (\ref{eqn_unif_5}) and (\ref{eqn_unif_6}). Our aim is to bound the functions $\lbrace f_{i, \varepsilon}\rbrace_{i=1, 2}$ analogously to the proof of Lemma \ref{lem_bound_first}. Only potential difficulty might be that $f_{i, \varepsilon}(u_\varepsilon(t))$ contains a term $\dot{\theta}_i(t)$ which might blow up as $\varepsilon\rightarrow 0$. However, this term is paired with the scalar $\varepsilon$, so that the (potential) blowing up is canceled.
\end{proof}

\begin{lemma}\label{lemma_aux_strip}For $i=1, 2,$ we define a function $\overline{G}_i$ on the set $\lbrace x=(s_1, \theta_1, s_2, \theta_2)\in(\R/T\mathbb{Z}\times S^1)^2\,\vert\, \pi_{s_1, s_2}(x)\in\overline{S_K}\rbrace$ by
$$\overline{G}_i:x\mapsto\langle P, v_i^{\bot}\rangle.$$ 
Then there is a $\varepsilon_0>0$ and $\widehat{\delta}>0$ such that for every $\varepsilon\in[0, \varepsilon_0)$ it holds that
$$\vert\langle P, v_1\rangle\vert>\widehat{\delta}\wedge \vert\langle P, v_2\vert\rangle>\widehat{\delta}$$
at each $x\in \lbrace |\overline{G}_1(x)|\leq\varepsilon^{1/2}\wedge|\overline{G}_2(x)|\leq\varepsilon^{1/2}\rbrace.$
\end{lemma}

\begin{proof}
Since we are outside of the special and diagonal points, we can find a bound $\widehat{\delta}>0$ for the compact set $\lbrace \overline{G}_1=\overline{G}_2=0\rbrace$. Hence, we can also find a bound $\widehat{\delta}>0$ for a small open neighborhood $U$ of the set $\lbrace \overline{G}_1=\overline{G}_2=0\rbrace$. Finally, we can find $\varepsilon_0>0$ small such that $U\supset\lbrace |\overline{G}_1|\leq\varepsilon_0^{1/2}\wedge|\overline{G}_2|\leq\varepsilon_0^{1/2}\rbrace$, see the diffeomorphism from Lemma \ref{lemma_standard_square}. So the lemma follows.
\end{proof}

\begin{rem}\label{rem_sheets_of_G} By the diffeomorphism from Lemma \ref{lemma_standard_square} the set $\lbrace \overline{G}_1=\overline{G}_2=0\rbrace$ is described as a $4$-sheet graph over $\overline{S_K}$ with the sheets given by $$(\pm n^\ast_1(s_1, s_2), \pm n^\ast_2(s_1, s_2))\in S^1\times S^1.$$
It is natural to inspect what happens to the sheets $\lbrace \overline{G}_1= \overline{G}_2=0\rbrace$, if we extend the domain of $\overline{G}_1, \overline{G}_2$ also over weakly special and weakly diagonal points. We will be denote such a extension by $\overline{\overline{G}}_1, \overline{\overline{G}}_2$.  Over a special point, two sheets of $\lbrace \overline{\overline{G}}_1=\overline{\overline{G}}_1=0\rbrace$ get connected by a circle $S^1$ and over diagonal points they get all connected by the torus, i.e., the diagonal $\Delta_{full}$. We will study these sets later in Section \ref{sec_appl_chord}.
\end{rem}

\begin{rem_not}Let us introduce the following non-negative constants
$$\widehat{P}:=\hbox{sup}_{(\R/T\mathbb{Z})^2}\lbrace \vert P\vert\rbrace+1, \widehat{\kappa}:=\hbox{sup}_{\R/T\mathbb{Z}}\lbrace\vert \kappa(s)\vert\rbrace, \widehat{\tau}:=\hbox{sup}_{\R/T\mathbb{Z}}\lbrace\vert \tau(s)\vert\rbrace.$$
The constants are real numbers by the compactness of $\R/T\mathbb{Z}$. Note also that $d_i^{-1}<2$ provided that $\varepsilon<\min\lbrace\varepsilon_{good}, 1\rbrace.$
\end{rem_not}

\begin{defn}\label{defn_g_strip}Let $\varepsilon>0$. Then we define a \textbf{$G_{K, \varepsilon}$-strip} as the set $$\big\lbrace x\in M_{K, \varepsilon}\,\vert\,\vert G_1(x)\vert\leq\varepsilon c_G\wedge\vert G_2(x)\vert\leq\varepsilon c_G\big\rbrace,$$
where
$$c_G=\frac{2\widehat{\kappa}\widehat{P}^2+4\widehat{\tau}\widehat{P}^2+3\widehat{P}+1}{\widehat{\delta}}$$
and $G_i$ are the restrictions of $\overline{G}_i$ to $\pi_{s_1, s_2}^{-1}(S_K)$.
\end{defn}

\begin{thm}\label{thm_g}There is a $\varepsilon_0>0$ such that for each $\varepsilon\in(0, \varepsilon_0)$ the $G_{K, \varepsilon}$-strip has the following properties: 
\begin{itemize}
\item[$(i.)$] The $G_{K, \varepsilon}$-strip is a submanifold with corners of $Int(M_{K, \varepsilon})$.
\item[$(ii.)$] $\pi_{s_1, s_2}$ induces on $G_{K, \varepsilon}$-strip a structure of a locally trivial fiber bundle with fibers diffeomorphic to the square.
\item[$(iii.)$] $-\nabla E_\varepsilon$ is strictly outward pointing on $\lbrace |G_1|<\varepsilon c_G\wedge |G_2|=\varepsilon c_G\rbrace$, strictly inward pointing on $\lbrace |G_2|<\varepsilon c_G\wedge |G_1|=\varepsilon c_G\rbrace$. If $x\in\lbrace |G_2|<\varepsilon c_G\wedge |G_1|=\varepsilon c_G\rbrace$, then there is $\delta>0$ such that $\lbrace x_\varepsilon\cdot_{E_\varepsilon}[-\delta, \delta]\rbrace\cap G_{K, \varepsilon}\hbox{-strip}=x_\varepsilon$ and $-\nabla E_\varepsilon(x_\varepsilon)\neq 0$. (The analogous statement holds also for a generic approximation $X_{E_\varepsilon}$ of $-\nabla E_\varepsilon$.)
\item[$(iv.)$] As $\varepsilon\rightarrow 0$, the width of $G_{K, \varepsilon}$-strip (in $\theta_1, \theta_2$ directions and with respect to $\widetilde{\langle\cdot,\cdot\rangle}$ metric) converge to $0$.
\item[$(v.)$] The $G_{K, \varepsilon}$-strip contains in its interior all critical points (except $\Delta_\varepsilon$) of $E_\varepsilon$ on $M_{K, \varepsilon}$.
\end{itemize}
\end{thm}

\begin{proof}
Ad $(i.)$: Let us show that the map $(G_1, G_2):\lbrace x\in(\R/T\mathbb{Z}\times S^2)^2\,\vert\,\pi_{s_1, s_2}\in S_K\rbrace\rightarrow\R^2$ is stratum transverse to $[-\varepsilon c_G, \varepsilon c_G]^2$ for $\varepsilon\in(0, \varepsilon_0)$, where $\varepsilon_0>0$ is small. Hence, let us describe the matrix $B:=(dG_1, dG_2)^T$ at any point $x\in (G_1, G_2)^{-1}([-\varepsilon c_G, \varepsilon c_G]^2)$:
\begin{equation*}
B=\begin{pmatrix}
\sin(\theta_1)\kappa(s_1)\langle P, \dot{\gamma}(s_1)\rangle-\tau(s_1)\langle P, v_1\rangle& -\langle v_1^{\bot}, \dot{\gamma}_2(s_2)\rangle   \\
-\langle P, v_1\rangle & 0  \\
\langle v_1^{\bot}, \dot{\gamma}_2(s_2)\rangle & \sin(\theta_2)\kappa(s_2)\langle P, \dot{\gamma}(s_2)\rangle-\tau(s_2)\langle P, v_2\rangle \\
0 & -\langle P, v_2\rangle 
\end{pmatrix}
\end{equation*}
Let $\varepsilon_0>0$ be small. Then, by Lemma \ref{lemma_aux_strip} we can bound $|\langle P, v_i\rangle|$ from bellow by $\overline{\delta}$. Thus we have that $\hbox{rank}(B)=2$. So we obtain the stratum transversality and by Theorem \ref{thm_invers_im} $(G_1, G_2)^{-1}([-\varepsilon c_G, \varepsilon c_G]^2)$ is a submanifold of $(\R/T\mathbb{Z}\times S^2)^2$. In addition, if $\varepsilon_0>0$ is small enough, we obtain by Lemma \ref{lemma_aux_strip} that the $G_{K, \varepsilon}$-strip is in fact a submanifold of $Int(M_{K, \varepsilon})$.

Ad $(ii.)$: Recall from above that at $x\in G_{K, \varepsilon}$-strip $\partial_{\theta_i}G_i=\pm \langle P, v_i\rangle\neq 0$. So the fiber bundle structure immediately follows from Ehresmann's Theorem \ref{thm_ehrsm}. The structure of each fiber follows from the diffeomorphism from Lemma \ref{lemma_standard_square}. 

Ad $(iii.)$: Let us compute the quantity $\widetilde{\langle \widetilde{\nabla}\overline{G}_i, -\nabla E_\varepsilon\rangle}$  along $\lbrace M_{K, \varepsilon}\cap|\overline{G}_i|=\varepsilon c_G\rbrace$. Observe that there is a term of the form $\pm\frac{1}{\varepsilon}\langle P, v_i\rangle\langle P, v_i^\bot\rangle.$ Now, the constant $c_G$ was chosen such that this term is always leading, and hence its sign determines the behavior of the flow of $-\nabla E_\varepsilon$ along the boundary of the $G_{K, \varepsilon}$-strip. And the part $(iii.)$  for $-\nabla E_\varepsilon$ follows. For the case of $X_{E_\varepsilon}$ observe that $\widetilde{\langle \widetilde{\nabla}\overline{G}_i, -\nabla E_\varepsilon\rangle}$ is $O(1)$ along the compact set $\lbrace M_{K, \varepsilon}\cap|\overline{G}_i|=\varepsilon c_G\rbrace$.

Ad $(iv.)$: It is clear that if $\varepsilon_a, \varepsilon_b\in[0, \varepsilon_0)$ and $\varepsilon_a> \varepsilon_b$, then $\lbrace |G_i|\leq \varepsilon_b\rbrace\subset\lbrace |G_i|\leq \varepsilon_a\rbrace.$

Ad $(v.)$: This immediately follows from (\ref{eqn_of_crit_pt}), provided that $\delta_K$ is sufficiently small, i.e., weakly special and weakly diagonal sets are chosen sufficiently small.
\end{proof}

\begin{lemma}\label{lemma_solution_strip} Let $u_\varepsilon$ be a $\varepsilon$-solution from $p_\varepsilon$ to $q_\varepsilon$ ($p_\varepsilon, q_\varepsilon\in M_{K, \varepsilon}$) such that $\pi_{s_1, s_2}\circ u_\varepsilon\subset S_K$. Then $u_{\varepsilon}\subset G_{K, \varepsilon}\hbox{-strip}$.
\end{lemma}

\begin{proof}
Until now, we were studying the space of ``\textit{out-out}'' chords - $M_{K, \varepsilon}$. This space turned out to be a submanifold of $(\R/T\mathbb{Z}\times S^1)^2$, which has a nice topology outside of weakly special and weakly diagonal pairs, see Theorem \ref{thm_mfld_coners} and Lemma \ref{lemma_standard_square}. Also, the behavior of $-\nabla E_\varepsilon$ was described in Lemma \ref{lem_grad_standard}. 

However, one can consider $M_{K, \varepsilon}$ as one part of the splitting of the whole space $(\R/T\mathbb{Z}\times S^1)^2$. We can naturally split $(\R/T\mathbb{Z}\times S^1)^2$ (or more precisely $(\R/T\mathbb{Z}\times S^1)^2\setminus\Delta_{full}$) in four sets of chords that are classified by the conditions on their endpoints $\lbrace F_1^{[\varepsilon]}\lessgtr 0, F_2^{[\varepsilon]}\lessgtr 0\rbrace$, i.e., the sets of ``\textit{out-out, in-out, out-in, in-in}'' chords. It is not hard to construct for them analogous statements as Theorem \ref{thm_mfld_coners}, Lemma \ref{lemma_standard_square}, and Lemma \ref{lem_grad_standard}. That is, all of them have a manifold structure, outside of weakly special and weakly diagonal pairs, they are diffeomorphic to the product of $S_K$ and the square. In addition, the behavior of $-\nabla E_\varepsilon$ on the lower strata of these manifolds is induced from the knowledge of the behavior of $-\nabla E_\varepsilon$ on the lower strata of $M_{K, \varepsilon}$.

Next, we construct a trap for the global dynamics of $u_\varepsilon$. Here, recall that $u_\varepsilon$ can be not only a flow line corresponding to $-\nabla E_\varepsilon$, but also to a generic approximation $X_{E_\varepsilon}$. Assume that $u_\varepsilon$ flows out at some point from \textit{out-out} manifold ($M_{K, \varepsilon}$), then at that point $u_\varepsilon$ will enter the \textit{out-in} manifold. But from the \textit{out-in} manifold, nothing can escape outside of weakly special and weakly diagonal pairs. So also $u_\varepsilon$ can not escape from the \textit{out-in} manifold, and hence could not flow to $q_\varepsilon$. Contradiction. Thus $u_\varepsilon\subset M_{K, \varepsilon}$. See Figure \ref{figure_laypunov_trap} on the left.

\begin{figure}[!htbp]
\centering
\includegraphics[scale=0.68]{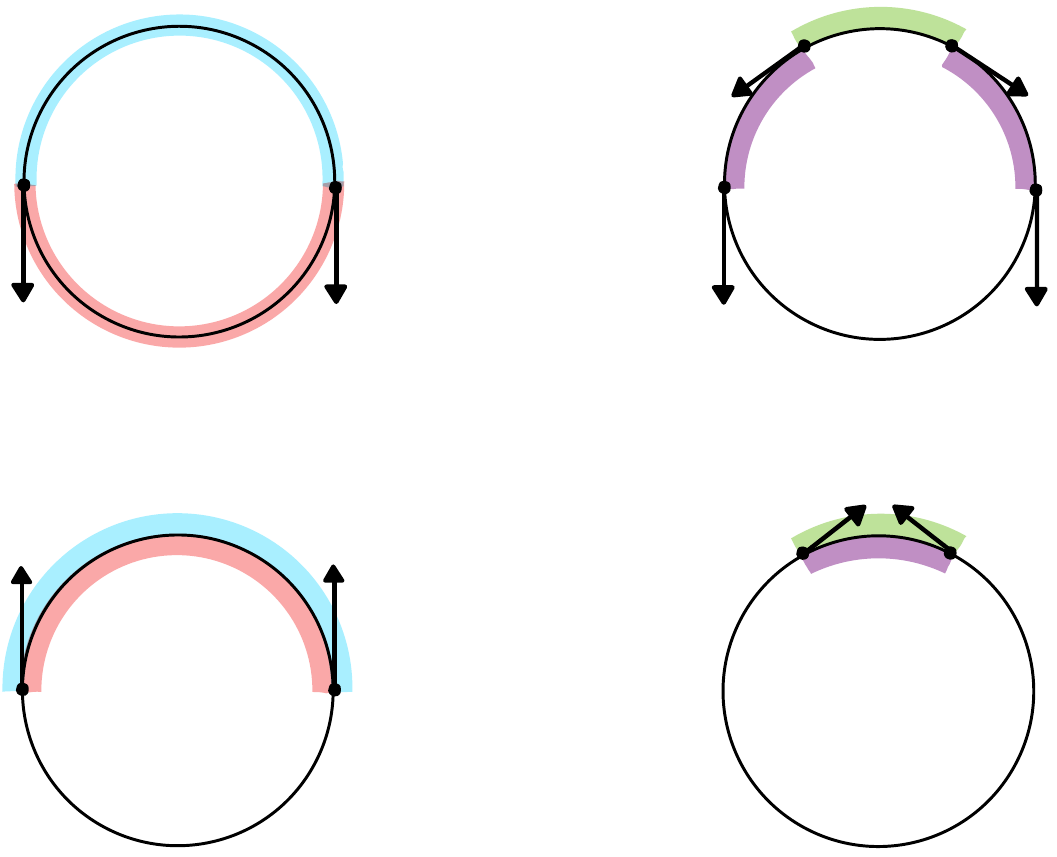}
\vspace{0.3cm}
\caption{Cartoon pictures of the discussed spaces restricted to fibers $\pi^{-1}_{s_1, s_2}((s_1, s_2))$. On the left: In blue is the \textit{out-out} manifold ($M_{K, \varepsilon}$) and in red is the \textit{out-in} manifold. Black arrows describe $-\nabla E_\varepsilon$ at the boundaries. On the right: In green is the set $\lbrace |G_2|\leq\varepsilon c_G\cap|G_1|\leq\varepsilon c_G\cap M_{K, \varepsilon}\rbrace$ ($G_{K, \varepsilon}$-strip). In deep purple is the set $\lbrace |G_2|\geq\varepsilon c_G\cap|G_1|\leq\varepsilon c_G\cap M_{K, \varepsilon}\rbrace$}
\label{figure_laypunov_trap}
\end{figure}
 
In fact, we can improve the localization of $u_\varepsilon$. From Theorem \ref{thm_g} we know that $u_\varepsilon$ starts and ends at the $G_{K, \varepsilon}$-strip, and we also know the description of $-\nabla E_\varepsilon$ along the boundary of the $G_{K, \varepsilon}$-strip. 

Now, we will use a similar argument as above. Let us assume that $u_\varepsilon$ escapes at some point from the $G_{K, \varepsilon}$-strip. Then $u_\varepsilon$ will flow through the set (manifold) $\lbrace |G_2|\geq\varepsilon c_G\cap|G_1|\leq\varepsilon c_G\cap M_{K, \varepsilon}\rbrace$. But from this manifold $u_\varepsilon$ can flow only into \textit{out-in} manifold, and hence will never get to $q_\varepsilon$, contradiction. See also Figure \ref{figure_laypunov_trap} on the right.
\end{proof}

\begin{lemma} There is $\varepsilon_0$ such that $c_5$ and $c_6$ exist.
\end{lemma}

\begin{proof}
Observe that the only potential problematic terms in equations (\ref{eqn_unif_3}) and (\ref{eqn_unif_4}) are of the form $\pm\langle P, v_i^\bot\rangle$. However, by the assumptions on $u_\varepsilon$ and Lemma \ref{lemma_solution_strip} we know that $u_\varepsilon$ cannot leave the $G_{K, \varepsilon}$-strip, which gives us the desired bounds on those problematic terms. More precisely,
$$|\dot{\theta}_i(t)|\leq 2\widehat{\tau}\widehat{P}+1+c_G,$$
provided that $\varepsilon<\min\lbrace\varepsilon_{good}, 1\rbrace.$
\end{proof}

\begin{lemma}\label{lem_aux_h_strip}Let $\varepsilon>0$. For $i=1, 2$ we define a function $H_i$ on the set $\lbrace x=(s_1, \theta_1, s_2, \theta_2)\in(\R/T\mathbb{Z}\times S^1)^2\,\vert\, \pi_{s_1, s_2}(x)\in {S_K}\rbrace$ by
$$H_i:x\mapsto\widetilde{\langle \widetilde{\nabla}\langle P, v_i^{\bot}\rangle, -\nabla E_\varepsilon(x)\rangle}.$$ 
Then there is a $\varepsilon_0>0$ and $\widehat{c}_{H}>0$ such that for every $\varepsilon\in(0, \varepsilon_0)$ and every $x\in Int(G_{K, \varepsilon}\hbox{-strip})$ holds that
$$\vert\widetilde{\langle\widetilde{\nabla}H_i(x), -\nabla E_\varepsilon(x)\rangle}\vert\leq\widehat{c}_{H}+\left\vert\frac{\langle P, v_i\rangle}{\varepsilon}H_i(x)\right\vert.$$
\end{lemma}

\begin{proof}
See (\ref{eqn_unif_6}) for the splitting of $H_i(x)$ in the terms $O(\varepsilon)+O(1)+O(1/\varepsilon)$. Then the lemma follows from the compactness of $(\R/T\mathbb{Z}\times S^1)$ and the fact that $x\in G_{K, \varepsilon}\hbox{-strip}$ (so the terms $\langle P, v^\bot_i\rangle/\varepsilon$ will not blow up).
\end{proof}

\begin{defn}\label{defn_h_strip} Let $\varepsilon>0$ be small. Then we define a \textbf{$H_{K, \varepsilon}$-strip} as the set of $x\in Int(G_{K, \varepsilon}$-$strip)$ such that it is satisfied that
$$|H_i(x)|\leq\varepsilon c_H,$$
where $i=1, 2,$ and $c_H=(\widehat{c}_{H}+1)/\widehat{\delta}$ (recall Lemma \ref{lemma_aux_strip} for the definition of $\widehat{\delta}$).
\end{defn}

\begin{thm}\label{thm_h}There is a $\varepsilon_0>0$ such that for each $\varepsilon\in(0, \varepsilon_0)$ the $H_{K, \varepsilon}$-strip has the following properties: 
\begin{itemize}
\item[$(i.)$] The $H_{K, \varepsilon}$-strip is a submanifold with corners of $Int(G_{K, \varepsilon}$-$strip)$.
\item[$(ii.)$] $-\nabla E_\varepsilon$ is strictly outward pointing on $\lbrace |H_1|<\varepsilon c_H\wedge |H_2|=\varepsilon c_H\rbrace$, strictly inward pointing on $\lbrace |H_2|<\varepsilon c_H\wedge |H_1|=\varepsilon c_H\rbrace$ and never tangent to the corners. If $x\in\lbrace |H_2|<\varepsilon c_H\wedge |H_1|=\varepsilon c_H\rbrace$, then there is $\delta>0$ such that $\lbrace x_\varepsilon\cdot_{E_\varepsilon}[-\delta, \delta]\rbrace\cap H_{K, \varepsilon}\hbox{-strip}=x_\varepsilon$ and $-\nabla E_\varepsilon(x_\varepsilon)\neq 0$. (The analogous statement holds also for a generic approximation $X_{E_\varepsilon}$ of $-\nabla E_\varepsilon$.)
\item[$(iii.)$] The $H_{K, \varepsilon}$-strip contains in its interior all critical points (except $\Delta_\varepsilon$) of $E_\varepsilon$ on $M_{K, \varepsilon}$.
\end{itemize}
\end{thm}

\begin{proof}
Ad $(i.)$: Let $x\in Int(G_{K, \varepsilon}\hbox{-strip})$. By (\ref{eqn_unif_6}), we see that the following terms are the only terms of $\partial_{\theta_i}H_i(x)$ that can potentially blow up:
$$\frac{\langle P, v^\bot_i\rangle^2}{\varepsilon}\hbox{ and }\frac{\langle P, v_i\rangle^2}{\varepsilon}.$$
However, we have that 
$$\frac{\langle P, v_i^\bot\rangle^2}{\varepsilon}\leq \varepsilon c_G^2$$
and
$$\frac{\langle P, v_i\rangle^2}{\varepsilon}\geq\frac{\widehat{\delta}^2}{\varepsilon}.$$
So the case $(i.)$ follows from Theorem \ref{thm_mfld_coners}.

Ad $(ii.)$: This case is obtained from Lemma \ref{lem_aux_h_strip} and the definition of the boundary of the $H_{K, \varepsilon}$-strip.

Ad $(iii.)$: This case is analogous to the Theorem \ref{eqn_unif_6} $(v.)$.
\end{proof}

\begin{lemma}\label{lemma_solution_strip_H}  Let $u_\varepsilon$ be a $\varepsilon$-solution from $p_\varepsilon$ to $q_\varepsilon$ ($p_\varepsilon, q_\varepsilon\in M_{K, \varepsilon}$) such that $\pi_{s_1, s_2}\circ u_\varepsilon\subset S_K$. Then $u_{\varepsilon}\subset H_{K, \varepsilon}\hbox{-strip}$.
\end{lemma}

\begin{proof}
The lemma is based on a similar trick as in Lemma \ref{lemma_solution_strip}. Let us assume that $u_\varepsilon$ leaves at some point $x$ the $H_{K, \varepsilon}$-strip. This can be done in two ways. 

First, if $x\in\partial(H_{K, \varepsilon}$-strip). Then $u_\varepsilon$ will flow through the set (manifold)  $\lbrace |H_2|\geq\varepsilon c_H\cap|G_1|\leq\varepsilon c_H\cap G_{K, \varepsilon}\hbox{-strip}\rbrace$. From $\lbrace |H_2|\geq\varepsilon c_H\cap|G_1|\leq\varepsilon c_H\cap G_{K, \varepsilon}\hbox{-strip}\rbrace$ $u_\varepsilon$ cannot flow directly back to the $H_{K, \varepsilon}$-strip. So the $u_\varepsilon$ can only flow out from the $G_{K, \varepsilon}$-strip. But this is not possible by Lemma \ref{lemma_solution_strip}.

Second, let $x$ be on the boundary of the $G_{K, \varepsilon}$-strip. But again, this could not happen by Lemma \ref{lemma_solution_strip}. Contradiction.
\end{proof}

\begin{lemma} There is $\varepsilon_0$ such that $c_6$ and $c_8$ exist.
\end{lemma}

\begin{proof}
Let us consider the equations (\ref{eqn_unif_6}) and (\ref{eqn_unif_8}). The functions $\lbrace f_{i, \varepsilon}\rbrace_{i=3, 4}$ are straightforward to bound. The bounds for the remaining terms follow from Definition \ref{defn_h_strip} and Lemma \ref{lemma_solution_strip_H}.
\end{proof}

\section{Conley index}
We are going to recollect basic notions of Conley index theory and give an explicit construction of a regular index pair that will be used later.

\begin{defn} Let $\varphi=\lbrace\varphi^t\rbrace_{t\in\R}$ be a flow on locally compact metric space $M$. For a compact subset $N\subset M$ we put
$$I(N, \varphi):=\lbrace x\in N\,|\,\varphi^t(x)\in N\hbox{ for all }t\in\R\rbrace.$$ 
$N$ is called an \textbf{isolating neighborhood} if $I(N, \varphi)\subset Int(N)$.

A (compact) subset $\mathcal{S}\subset M$ is called \textbf{isolated invariant set} if $\mathcal{S}=I(N, \varphi)$ for some isolating neighborhood $N$.
\end{defn}

\begin{defn}\label{defn_index_pair} An \textbf{index pair} $(N, L)$ for an isolated invariant set $\mathcal{S}$ is a pair of compact sets $L\subset N\subset M$ such that
\begin{itemize}
\item[$(i.)$] $\mathcal{S}=I(\overline{N\setminus L}, \varphi)\subset Int(N\setminus L)$.
\item[$(ii.)$] $L$ is an \textbf{exit set} for $N$. That is; for all $x\in N$, if $\varphi^t(x)\notin N$ for some $t>0$, then there is $\tau\in[0, t]$ with $\varphi^\tau(x)\in L$.
\item[$(iii.)$] L is \textbf{positively invariant} in $N$. That is; if $x\in L$ and $t>0$ such that $\varphi([0, t],x)\subset N$, then $\varphi([0, t],x)\subset L$.
\end{itemize}
The \textbf{Conley index} associated to an isolated invariant set $\mathcal{S}$ is defined to be the pointed homotopy type $(N/L,[L]).$
\end{defn}

\begin{thm}\cite{conley1978isolated, Salamon_1985conley} Let $\mathcal{S}$ be an isolated invariant set, then:
\begin{itemize}
\item[$(i.)$] $\mathcal{S}$ admits an index pair.
\item[$(ii.)$] The Conley index is invariant under the choice of the index pair.
\item[$(iii.)$] The Conley index is invariant under continuation: Let $\lbrace\varphi_\lambda\rbrace_{\lambda\in[0, 1]}$ be a homotopy of flows on $M$ such that $N$ is an isolating neighborhood for each $\varphi_\lambda,\lambda\in[0, 1]$. Then the Conley indices of $I(N, \varphi^0)$ and $I(N, \varphi^1)$ are the same.  
\end{itemize} 
\end{thm}

\begin{defn} An index pair $(N, L)$ is called \textbf{regular} if the inclusion $L\hookrightarrow N$ is a cofibration, i.e. whenever we are given a topological space $M$, continuous maps $f:N\rightarrow M$ and $H:L\times[0, 1]\rightarrow M$ such that $f(x)=H(x, 0)$ for all $x\in L$, then there exists a continuous map $F:N\times [0, 1]\rightarrow M$ extending $H$ such that $F(x, 0)=f(x)$ for all $x\in X$.
\end{defn}

\begin{rem} For a regular index pair $(N, L)$ it holds that $H_\ast^{sing}(N, L; \mathbb{Z}_2)\cong H_\ast^{sing}(N/L, [L]; \mathbb{Z}_2)$. Also, the canonical inclusion of a subcomplex of a CW-complex is a cofibration.
\end{rem}

\begin{lemma}\cite{Salamon_1985conley}\label{lem_cond_reg}An index pair $(N, L)$ is regular if it holds that
$$\varphi([0, \varepsilon], x)\not\subset \overline{N\setminus L}$$
for every $x\in L$ and every $\varepsilon>0$. 
\end{lemma}

\begin{lemma} \cite{Salamon_1985conley} Let $(N, L)$ be an index pair for an isolated invariant set $\mathcal{S}$, then there is a continuous Lyapunov function $h:N\rightarrow[0, 1]$ such that
\begin{itemize}
\item[$(i.)$]$h(x)=1$, iff $\varphi([0, \infty), x)\subset N$ and the omega limit of $x$ is in $\mathcal{S}$.
\item[$(ii.)$]$h(x)=0$, iff $x\in L$.
\item[$(iii.)$] If $t>0, 0<h(x)<1, \varphi([0, t], x)\subset N$, then $h(\varphi^t(x))<h(x)$.
\end{itemize}
Moreover, we obtain the following. Let $\varepsilon\in(0, 1)$. If we replace $L$ by $L_\varepsilon=\lbrace x\in N\,\vert\, h(x)\leq\varepsilon\rbrace$, then $(N, L_\varepsilon)$ is a regular index pair for $\mathcal{S}$.
\end{lemma}

\begin{example}\label{example_smooth_index_pair}
Let $X_{E_0}$ be a gradient-like vector field adapted to $E_0$. Let $p_0\in Crit_1(E_0)\setminus \Delta_0$, $q_0\in Crit_0(E_0)\setminus \Delta_{0}$ and $u$ be a $0$-solution from $p_0$ to $q_0$. We would like to construct a smooth regular index pair $(N, L)$ for the invariant set $\lbrace p_0, q_0, u\rbrace.$ 

Since $Ind_{E_0}(q_0)=0$, $q_0$ is a local minimum. Let $\delta>0$ and $U_{q_0}$ be a open neighborhood of $q_0$. We put $N_{q_0}:=U_{q_0}\cap E^{-1}_0[E(q_0), E(q_0)+\delta]$ which is an embedded $2$-disc, provided that $\delta>0$ and $U_{q_0}$ are sufficiently small. Then $(N_{q_0}, \emptyset)$ is a (regular) index pair for $q_0$.

From the index reasons, $p_0$ has $1$-dimensional stable and unstable manifolds; $W^s_{X_{E_0}}(p_0), W^u_{X_{E_0}}(p_0).$ Let $\delta>0$ be small. Then $W^s_{X_{E_0}}(p_0)\pitchfork\partial E^{-1}_0([E_0(p_0)-\delta, E_0(p_0)+\delta])$ at two points $A, B$. Next, there are two small closed neighborhoods $D_A, D_B$ of points $A, B$ in $\partial E^{-1}_0([E_0(p_0)-\delta, E_0(p_0)+\delta])$, respectively. Now we put
$$N_{p_0}:=\big(\varphi([0, \infty),D_A\cup D_B)\cup W_{E_0}^u(p_0)\big)\cap E^{-1}_0([E_0(p_0)-\delta, E_0(p_0)+\delta])$$
and 
$$L_{p_0}:=N_{p_0}\cap E^{-1}_0(E_0(p_0)-\delta).$$
See also Figure \ref{figure_Conely}. By the construction, if $\delta>0$, $D_A, D_B$ are small enough, then $N_{p_0}$ is a (smooth) $2$-manifold with corners and $(N_{p_0}, L_{p_0})$ is an index pair of $p_0$. Also note that $L_{p_0}$ has two connected components. However, the index pair is not necessarily regular. We also remark that by Hartman-Grobman Theorem \cite{Palis1993HyperbolicityAS} there is a homeomorphism conjugating $X_{E_0}$ with $D X_{E_0}(p_0)$ around $p_0$. Which describes the dynamics on $N_{p_0}$. We will leave $(N_{p_0}, L_{p_0})$ for now and construct an index pair for $\lbrace p_0, q_0, u\rbrace$. 

\begin{figure}[!htbp]
\labellist
\pinlabel $p_0$ at 120 770
\pinlabel $A$ at 102 846
\pinlabel $B$ at 102 667
\endlabellist
\centering
\includegraphics[scale=0.75]{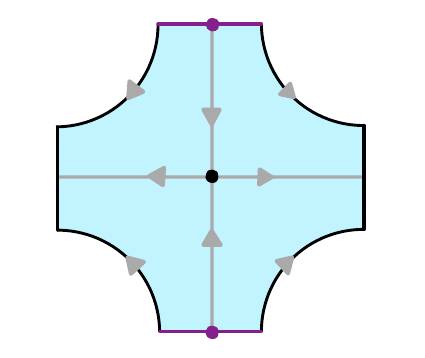}
\vspace{0.3cm}
\caption{The neighbourhood $U_{p_0}$ of the critical point $p_0$. In purple, there are neighbourhoods $D_A, D_B$ of points $A, B$, respectively.}
\label{figure_Conely}
\end{figure}

Let $C$ be the unique point in $u\cap \partial N_{q_0}$. Let $A_{1, q_0}\in\partial N_{q_0}\setminus C$. If $A_{1, q_0}$ is sufficiently close to $C$, then the set 
$$\phi^t_{X_{E_0}}\big((-\infty, 0], A_{1, q_0}\big)\cap (D_A\cup D_B)$$
consists of the unique point $(=: A_{1, p_0})$ and we can assume that $A_{1, p_0}\in D_A$. This follows from the continuous dependence of the vector fields on the initial conditions (Remark \ref{rem_continuity_ode}) and the Hartman-Grobman Theorem. Similarly, we can take another $A_{2, q_0}\in N_{q_0}$ close to $C$ such that the $A_{2, p_0}$ is on the same connected component of $L_{p_0}$ as $A_{1, p_0}$. We can choose $A_{2, q_0}$ such that the arc $\overline{A, A_{2, p_0}}$ inside $D_A$ is shorter then $\overline{A, A_{1, p_0}}$.

This gives us a rectangular tube $\overline{A_{1, q_0}, A_{1, p_0}, A_{2, p_0}, A_{2, q_0}}$ which is from $D_A$ to $\partial N_{q_0}$ smoothly foliated by the flow lines. Here we recall that because $A_{1, p_0}, A_{2, p_0}$ are close to $C$, the tube contains no critical point of $E_0$.

Then we can consider a smooth section $\ell_A$ of this foliation from $A_{2, p_0}$ to $A_{1, q_0}$ that is transverse to the leaves. In another words, $\ell_A\pitchfork X_{E_0}$. See Figure \ref{figure_Conley_2} left.

Now we construct analogously a tube from $D_B$ to $\partial N_{q_0}$ and a section $\ell_B$. Then we obtain a hull around $\lbrace p_0, q_0, u\rbrace$ connecting $N_{p_0}$ and $N_{q_0}$ with $\ell_A$ and $\ell_B$, which give us almost $N$. We only need to deform the remaining two components of the hull, which are not transverse to $X_{E_0}$. This can be done by the same trick as above, and hence we obtain a regular index pair $(N, L)$. See Figure \ref{figure_Conley_2}. Observe that by the construction, $N$ is a smooth manifold with corners and $L$ is contractible. 

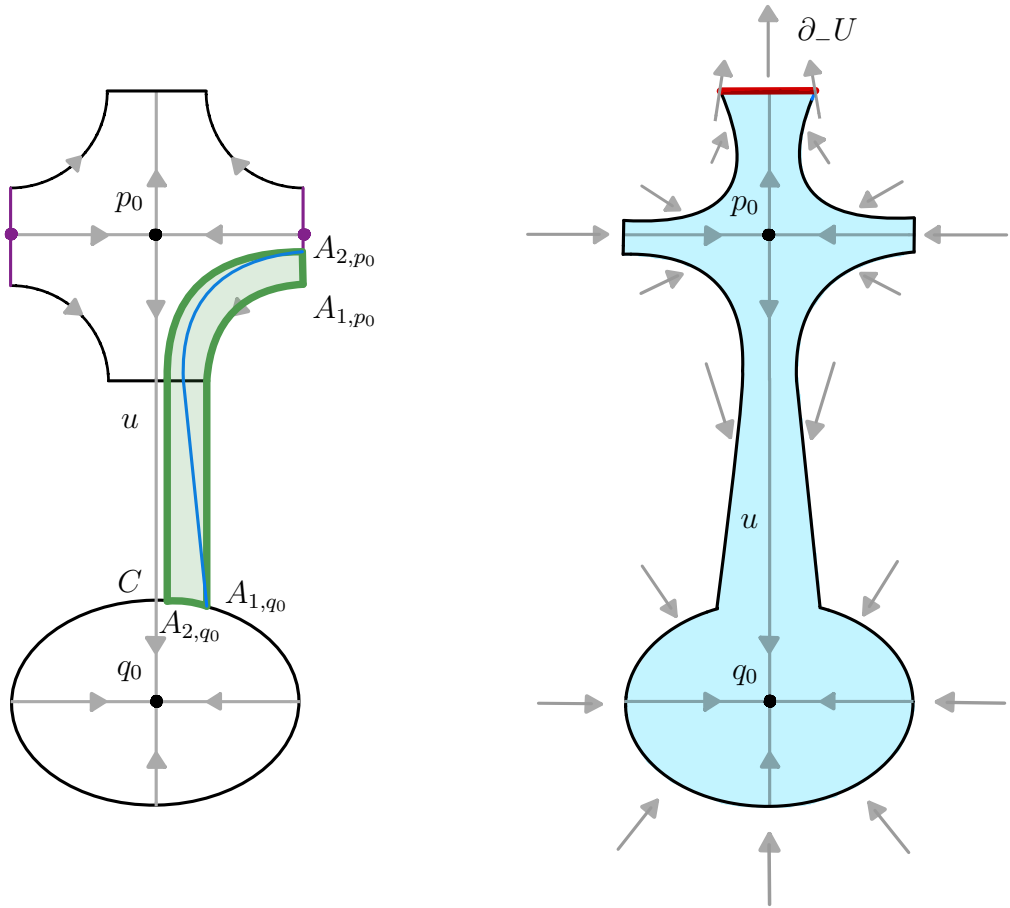
\begin{figure}[!htbp]

\labellist
\pinlabel $p_0$ at 120 730
\pinlabel $q_0$ at 120 495
\pinlabel $A_{1, p_0}$ at 227 675
\pinlabel $A_{2, p_0}$ at 227 705
\pinlabel $A_{1, q_0}$ at 183 532
\pinlabel $A_{2, q_0}$ at 150 517
\pinlabel $p_0$ at 426 728
\pinlabel $q_0$ at 426 493
\pinlabel ${\partial_- U}$ at 467 815
\pinlabel $u$ at 428 570
\pinlabel $u$ at 120 620
\pinlabel $C$ at 120 540
\endlabellist
\centering
\includegraphics[scale=0.75]{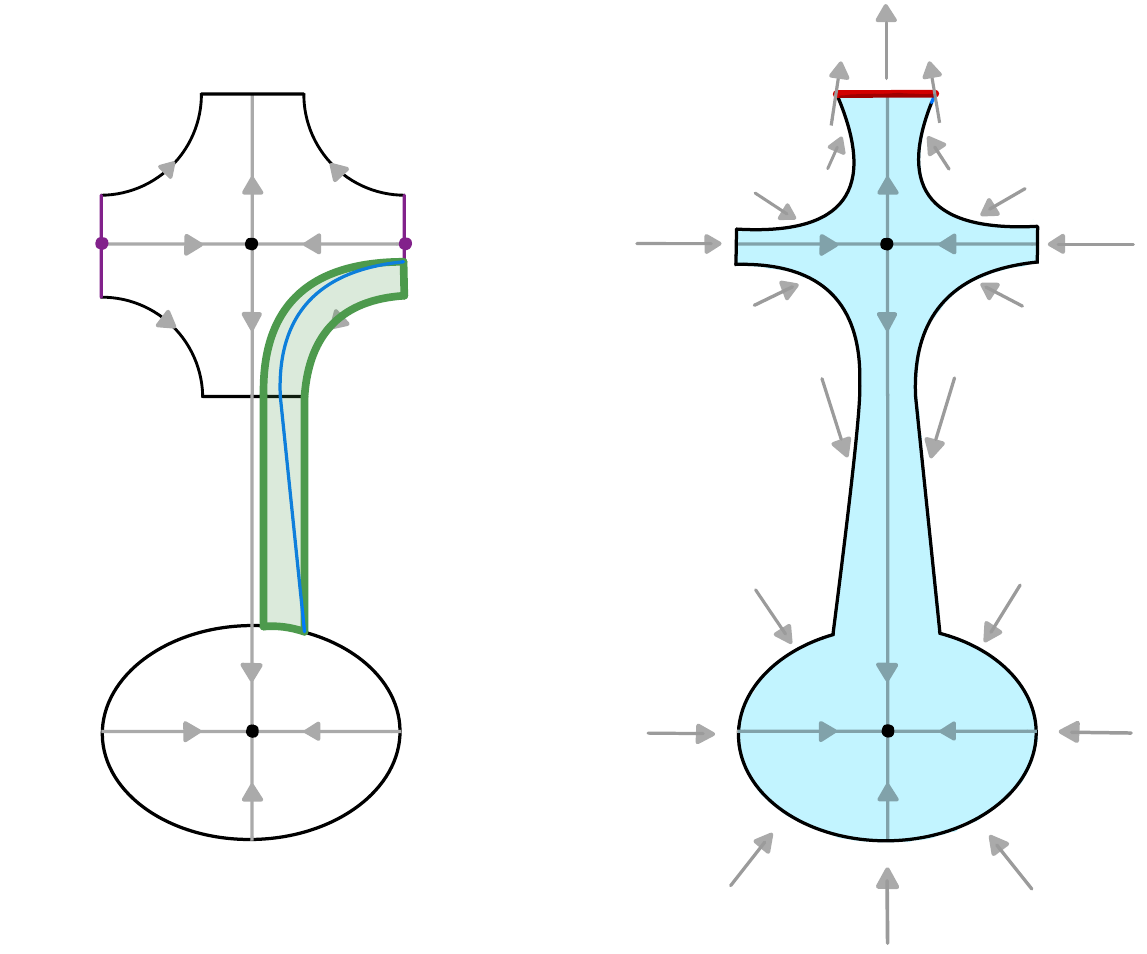}
\vspace{0.3cm}
\caption{\textit{On the left:} In \textcolor{teal}{green} the rectangular tube $\overline{A_{1, q_0}, A_{1, p_0}, A_{2, p_0}, A_{2, q_0}}$ connecting $N_{p_0}$ and $N_{q_0}$. In \textcolor{blue}{blue} the section $\ell_A$. \textit{On the right:} The regular index pair $(N, L)$ for $\lbrace p_0, q_0, u\rbrace$. Here in \textcolor{purple}{red} is the exit set $L$. }
\label{figure_Conley_2}
\end{figure}
\end{example}

\newpage
\section{Morse homology}
Using the Conley index theory, we define a Morse homology for manifolds with corners and introduce the Morse Homology Theorem \ref{thm_Morse_homology}.

\begin{assump}In this section, we will assume the following. $(M, g)$ will be a compact Riemannian $n$-manifold with corners. Next, $f:M\rightarrow\R$ will be a Morse function such that $Crit(f)\cap\partial M=\emptyset$.
\end{assump}

\begin{defn}A gradient-like vector field $X_f$ (adapted to $f$) is said to be \textbf{regular on $\partial M$} if $X_f$ is not tangent to any $(n-1)$-face of $M$. Then $\partial_- M$ will denote the closure of the set of all $x\in\partial M$ such that $X_f$ is strictly outward-pointing at $x$.
\end{defn}

\begin{rem} Let $X_f$ be a gradient-like vector field which is regular on $\partial M$. Then $\partial_- M$ is a manifold with corners with stratification canonically induced from $M$.
\end{rem}

\begin{defn} Let $X_f$ be a gradient-like vector field on $M$, then we put
$$\mathcal{S}_{X_f}:=\bigcup_{p, q\in Crit(f)} W_{X_f}(p, q).$$
(In particular, $Crit(f)\subset \mathcal{S}_{X_f}$.)
\end{defn}

\begin{lemma} Let $X_f$ be a gradient-like vector field which is regular on $\partial M$. Then $(M, \partial_- M)$ is a regular index pair for the isolated invariant set $\mathcal{S}_{X_f}$. And in particular, $\mathcal{S}_{X_f}$ is compact. 
\end{lemma}
\begin{proof}
By the assumption of $X_f$ on $\partial M$ it holds that $\mathcal{S}_{X_f}=I(M, \phi_{X_f})$. Since $M$ is compact, $\mathcal{S}_{X_f}$ is compact too. $\partial_- M$ is a positively invariant exit set of $M$, so $(M, \partial_-M)$ is an index pair for $\mathcal{S}_{X_f}$.  Then by Lemma \ref{lem_cond_reg}, the index pair $(M, \partial_- M)$ is regular.
\end{proof}

\begin{lemma}\label{lemma_isolat}Let $X_f$ be a gradient-like vector field which is regular on $\partial M$. If a gradient-like vector field $Y_f$ is sufficiently $C^0$-close to $X_f$, then $Y_f$ is also regular on $\partial M$ with the same $\partial_- M$. And in particular, $(M, \partial_-M)$ is an regular index pair for $\mathcal{S}_{Y_f}$.
\end{lemma}
\begin{proof}
The lemma immediately follows from the compactness of $M$ and the fact that the assumption for $X_f$ on $\partial M$ is open.
\end{proof}

\begin{defn}Let $X_f$ be a Morse-Smale vector field which is regular on $\partial_-M$. Then we define a \textbf{Morse chain complex} as the following graded abelian groups $C_\ast(X_f)$ together with a differential $\partial^m_\ast:C_\ast(X_f)\rightarrow C_{\ast-1}(X_f)$.

For $i\in\mathbb{N}_0$, we generate $C_i(X_f)$ by the critical points as
$$C_i(Y_f)=\mathbb{Z}\langle Crit_i(f)\rangle\otimes \mathbb{Z}_2.$$
Next for $p\in Crit_i(f)$ we put
\begin{equation}\label{rem_eqn_morse}
\partial_i^m(p):=\sum_{q\in Crit_{i-1}(f)} \#_2\mathcal{M}_{X_f}(p, q)\,q,
\end{equation}
where $\#_2\mathcal{M}_{X_f}(p, q)$ means modulo $2$ count of elements of $\mathcal{M}_{X_f}(p, q)$. Then we extend $\partial^m$ by $\mathbb{Z}_2$-linearity on the whole $C_i(X_f)$.

Then the \textbf{Morse homology} is defined by
$$HM_i(X_f; \mathbb{Z}_2)=\frac{\hbox{Ker} \,\,\partial^m_i}{\hbox{Im} \,\,\partial^m_{i+1}}.$$
\end{defn}

\begin{rem}\label{rem_cmp_lines} Let $X_f$ be a Morse-Smale vector field, $p\in Crit_\ast(f)$ and $q\in Crit_{\ast-1}(f)$. By Remark \ref{rem_Morse_smale_same_ind}, there are no flow lines between the critical points of the same index, so the set
\begin{equation}\label{eqn_inv_isolat_for_pair}
S_{X_f}(p, q):=W_{Y_f}(p, q)\cup\lbrace p, q\rbrace
\end{equation}
 is closed. Since $M$ is compact, $W_{X_f}(p, q)\cup\lbrace p, q\rbrace$ is compact too. It follows that $\mathcal{M}_{Y_f}(p, q)$ is a finite set, hence the sum in $(\ref{rem_eqn_morse})$ is finite too.
\end{rem}

\begin{lemma}\cite{Salamon1990MorseTT,Rot_2014}Let $X_f$ be a Morse-Smale vector field which is regular on $\partial M$. Then
\begin{equation}\label{eqn_square_to_zero}
\partial^m_\ast\circ\partial^m_{\ast-1}=0.
\end{equation}
\end{lemma}

\begin{sproof} We present two different ways how to show $(\ref{eqn_square_to_zero})$. Both of them strongly depend on the isolating property of $\mathcal{S}_{X_f}$.

\textit{First approach:}

Let $p\in Crit_i(f)$ and $q\in Crit_{i-2}(f)$. Let us consider any sequence of gradient-like trajectories of $W_{X_f}(p, q)$. Then the sequence converges subsequentially in Floer-Gromov $C^1_{loc}$ sense \cite{frauenfelder2020moduli} to a broken gradient-like flow line from $p$ to $q$ that passes through a single element of $Crit_{i-1}(f)$. This induces the compactification $\overline{\mathcal{M}_{X_f}(p, q)}$ of $\mathcal{M}_{X_f}(p, q)$ in the Floer-Gromov topology, see \cite{frauenfelder2020moduli}.

Hence, $\overline{\mathcal{M}_{X_f}(p, q)}$ has a structure of a compact $1$-dimensional manifold with boundary. In particular, $\overline{\mathcal{M}_{X_f}(p, q)}$ is diffeomorphic to a union of circles and closed intervals, which has an even number of boundary points.

Since
$$\partial^m_{i-1}\circ\partial^m_i(p)=\sum_{\widehat{q}\in Crit_{i-2}(f)}\#_2\partial\overline{\mathcal{M}_{X_f}(p, \widehat{q})}\,\widehat{q},$$
the relation $(\ref{eqn_square_to_zero})$ follows.

\textit{Second approach:}

Let $r\in Crit_j(f)$. Then by Hartman-Grobman Theorem \cite{Palis1993HyperbolicityAS} we can quickly find a regular index pair $(N_r, L_r)$ for $r$, see also Example \ref{example_smooth_index_pair}. In particular, $H_j^{sing}(N_r, L_r; \mathbb{Z}_2)\cong \mathbb{Z}_2$ which gives an identification
$$C_j(X_f)\cong\bigoplus_{r\in Crit_j(f)}H_j^{sing}(N_r, L_r; \mathbb{Z}_2).$$

Let $p\in Crit_i(f)$ and $q\in Crit_{i-1}(f)$. Next, let $\mathcal{S}_{X_f}(p, q)$ be the isolated invariant set from $(\ref{eqn_inv_isolat_for_pair})$ and $(N_2, N_0)$ be some regular index pair for $\mathcal{S}_{X_f}(p, q)$. Then we define $N_1:=N_1\cup(N_2\cap M_a)$, where $M_a=\lbrace x\in M\,\vert\,f(x)\leq a\rbrace$ for some $a\in(f(q), f(p))$. Now we introduce an isomorphism
$$\Delta_i(p, q):H_{i}^{sing}(N_p, L_p; \mathbb{Z}_2)\rightarrow H_{i-1}^{sing}(N_q, L_q; \mathbb{Z}_2)$$
as the composition
$$H_{i}^{sing}(N_p, L_p; \mathbb{Z}_2)\rightarrow H_{i}^{sing}(N_2, N_1; \mathbb{Z}_2)\xrightarrow{\partial} H_{i-1}^{sing}(N_1, N_0; \mathbb{Z}_2)\rightarrow H_{i-1}^{sing}(N_q, L_q; \mathbb{Z}_2),$$
where the first and the third map are isomorphisms induced by the homotopy independence of the Conley index. The map $\partial$ is given by the LES for the triple $(N_0, N_1, N_2).$

In fact, using a certain deformation argument for $X_f$ outside of $\mathcal{S}_{X_f}(p, q)$, one can choose a nice $N_2, N_0$ and show that $\Delta_i(p, q)$ actually counts the gradient-like trajectories from $p$ to $q$.

Now, we would like to inspect the global properties of the map $\Delta_\ast$.

Let us define
$$\mathcal{S}_{X_f}^{i, j}=\bigcup \lbrace W_{X_f}(p, q)\,\vert\,i\leq Ind_f(p)\leq Ind_f(q)\leq j\rbrace$$
which are compact isolated invariant subsets of $\mathcal{S}_{X_f}$. Then by \cite{conley1978isolated} there is a filtration
\begin{equation}\label{eqn_filter_conley_index}
\partial_-M=M_{-1} \subseteq M_0\subseteq\dots \subseteq M_{n-1}\subseteq M_n=M
\end{equation}
such that each $(M_i, M_{j-1})$ is a regular index pair for $\mathcal{S}_{X_f}^{i, j}$.

Hence, we obtain a commutative diagram
\[
\begin{tikzcd}
C_i(f) \arrow[d, "\cong"] \arrow[r, "\partial^m_i"] & C_{i-1}(f) \arrow[d, "\cong"] \arrow[r, "\partial^m_{i-1}"] & C_{i-2}(f) \arrow[d, "\cong"] \\
H_{i}^{sing}(N_i, N_{i-1}; \mathbb{Z}_2) \arrow[r, "\Delta_i"]     & H_{i-1}^{sing}(N_{i-1}, N_{i-2}; \mathbb{Z}_2) \arrow[r, "\Delta_{i-1}"]             & H_{i-2}^{sing}(N_{i-2}, N_{i-3}; \mathbb{Z}_2)\nospaceperiod              
\end{tikzcd}.
\]
and in particular also relation $(\ref{eqn_square_to_zero}).$
\end{sproof}

\begin{MorseHom}\cite{Salamon1990MorseTT,Rot_2014}\label{thm_Morse_homology}Let $X_f$ be a Morse-Smale vector field which is regular on $\partial M$. Then
\begin{equation*}
HM_\ast(X_f; \mathbb{Z}_2)\cong H_\ast^{sing}(M, \partial_-M; \mathbb{Z}_2).
\end{equation*}
\end{MorseHom}

\begin{sproof}The theorem follows from the standard argument which combines LESs of the triples from filtration $(\ref{eqn_filter_conley_index})$.
\end{sproof}

\section{Map}
In this section, we finally combine the uniform bounds together with Conley index theory and show the correspondence between the $0$-solutions and the $\varepsilon$-solutions.

\begin{defn}Let $p_\varepsilon, q_\varepsilon\in Crit(E_\varepsilon)$ and $X_{E_\varepsilon}$ be a gradient-like vector field. Then 
$$\mathcal{M}_{X_{E_\varepsilon}}^{out-out}(p_\varepsilon, q_\varepsilon)$$
will denote the restriction of the space of (unparametrized) trajectories from $\mathcal{M}_{X_{E_\varepsilon}}(p_\varepsilon, q_\varepsilon)$ to the trajectories that lie completely in $M_{K, \varepsilon}$.
\end{defn}

\begin{thm}\label{thm_almost_bijection} Let $p_0\in Cr_1(E_0)\setminus\Delta_0, q_0\in Cr_0(E_0)\setminus\Delta_0$ and $p_\varepsilon, q_\varepsilon$ be their unique corresponding critical points of $E_\varepsilon$ in $M_{K, \varepsilon}$. Let $X_{E_0}$ and $X_{E_\varepsilon}$ be generic approximations of $-\nabla E_0$ and $-\nabla E_\varepsilon$, respectively (see Corollary \ref{cor_perturb_grad_like}). Then for $\varepsilon>0$ sufficiently small, there is a multivalued function
$$\Phi^\varepsilon_{p_0, q_0}:\mathcal{M}_{X_{E_0}}(p_0, q_0)\rightarrow \mathcal{M}_{X_{E_\varepsilon}}(p_\varepsilon, q_\varepsilon).$$
with the following properties:
\begin{itemize}
\item[$(i.)$] For each $u\in \mathcal{M}_{X_{E_0}}(p_0, q_0)$ it holds that $\#_2\lbrace\Phi^\varepsilon_{p_0, q_0}(u)\rbrace=1$.
\item[$(ii.)$] All elements in $\Phi^\varepsilon_{p_0, q_0}(u)$ are homotopic through curves in $M_{K, \varepsilon}$ with fixed endpoints $p_\varepsilon, q_\varepsilon$.
\item[$(iii.)$] For each two different $u, \widehat{u}\in \mathcal{M}_{X_{E_0}}(p_0, q_0)$ it holds that $\Phi^\varepsilon_{p_0, q_0}(u)\cap\Phi^\varepsilon_{p_0, q_0}(\widehat{u})=\emptyset$ (as elements in the set $\mathcal{M}_{X_{E_\varepsilon}}(p_\varepsilon, q_\varepsilon)$) .
\item[$(iv.)$] $\Phi^\varepsilon_{p_0, q_0}$ is surjective.
\end{itemize}
Moreover, 
\begin{equation}\label{eqn_no_new_traj}
\mathcal{M}_{X_{E_\varepsilon}}(p_\varepsilon, q_\varepsilon)=\mathcal{M}_{X_{E_\varepsilon}}^{out-out}(p_\varepsilon, q_\varepsilon).
\end{equation}
\end{thm}
\begin{proof}
Follows from Lemmata \ref{thm_ambient_ift} and \ref{thm_base_ift}.
\end{proof}

\begin{rem} The count over $\mathbb{Z}_2$ in Theorem \ref{thm_almost_bijection} $(i.)$ is due to the fact that we considered for simplicity the Morse homology over $\mathbb{Z}_2$. However, one can expect that the count will also work over $\mathbb{Z}$, see \cite{Salamon1990MorseTT,Rot_2014} for the Morse homology over $\mathbb{Z}$. Later, using the techniques of Multiple-time scale dynamics in Chapter \ref{ch:file9}, we will show that $\Phi^\varepsilon_{p_0, q_0}$ is in fact a bijection.
\end{rem}

\begin{lemma}\label{thm_ambient_ift} If $\varepsilon>0$ is sufficiently small, then there is a multivalued function $\Phi^\varepsilon_{p_0, q_0}: \mathcal{M}_{X_{E_0}}(p_0, q_0)\rightarrow \mathcal{M}_{X_{E_\varepsilon}}(p_\varepsilon, q_\varepsilon)$ which is$\mod 2$ injective, i.e. $\Phi^\varepsilon_{p_0, q_0}$ satisfies properties $(i.)-(iii.)$ from Theorem \ref{thm_almost_bijection}.
\end{lemma}

\begin{proof}
Let $u\in \mathcal{M}_{X_{E_0}}(p_0, q_0)$. The strategy is the following. First, we will find a compact submanifold with corners $N\subset (\R/T\mathbb{Z})^2$ such that $(N, \partial_-N)$ is a regular index pair for $\lbrace p_0, q_0, u\rbrace.$ Then, we construct a lift of $N$ to $N_\varepsilon\subset M_{K, \varepsilon}$, where in particular $p_\varepsilon, q_\varepsilon\in N_\varepsilon$. Next, due to the genericity assumptions of $X_{E_0}$ we will be able to apply Morse homology Theorem \ref{thm_Morse_homology} to $X_{E_0}\vert_{N_\varepsilon}$. Hence, we obtain that
$$HM_1(X_{E_\varepsilon}\vert_{N_\varepsilon}; \mathbb{Z}_2)\cong H_1^{sing}(N_\varepsilon, \partial_-N_\varepsilon; \mathbb{Z}_2)\cong0,$$
where the last equality follows from the fact that the constructed $4$-manifold $N_\varepsilon$ is contractible and $\partial_-N_\varepsilon$ is contractible too. Then, we focus on how the relative Morse homology was computed. Since $p_\varepsilon, q_\varepsilon$ will be the only critical points of $E_\varepsilon$ inside $N_\varepsilon$, it means that by the homology reasons 
$$\partial^m(p_\varepsilon)=1\cdot q_\varepsilon \hbox{ over }\mathbb{Z}_2.$$
In other words, for $u$ we find $1 \mod 2$ of elements in $\mathcal{M}_{X_{E_\varepsilon}}(p_\varepsilon, q_\varepsilon)$. This gives us the assignment $u\mapsto \Phi^\varepsilon_{p_0, q_0}(u)$ and also the property $(i.)$. Then $(ii.)$ is clear from the contractibility of $N_\varepsilon$. Since $(N, \partial_- N)$ is a regular index pair for $\lbrace p_0, q_0, u\rbrace$, there could not be any other element of $\mathcal{M}_{X_{E_0}}(p_0, q_0)$ inside $N$ other than $u$ (recall that in $\mathcal{M}_{X_{E_0}}(p_0, q_0)$ we modded out the reparametrizations). Then, from the canonical choice of the lift $N_\varepsilon$, the case $(iii.)$ will follow. 

So it remains to do few things: to construct $N$, to show that $N_\varepsilon$ and $X_{E_\varepsilon}\vert_{N_\varepsilon}$ satisfy the assumptions of Morse homology Theorem \ref{thm_Morse_homology}, to verify the contractibility of $N_\varepsilon, \partial_-N_\varepsilon$ and the case $(iii.)$.

\textit{Construction of $N$:}

For the construction of the regular index pair $(N, \partial_- N)$, we consider precisely the index pair $(N, L)$ from Example \ref{example_smooth_index_pair}. We remark that due to the generic assumptions on $X_{E_0}$, the submanifold $N$ lies inside the standard set $S_K$.

\textit{Assumptions of Morse homology Theorem \ref{thm_Morse_homology}:}

We define $N_\varepsilon$ as the subset of the $G_{K, \varepsilon}$-strip (Definition \ref{defn_g_strip}) that projects under $\pi_{s_1, s_2}$ to $N$. Since $N\subset S_K$ and $N$ is contractible, by Theorem \ref{thm_g} we obtain that $N_\varepsilon\xrightarrow{\pi_{s_1, s_2}} N$ is a trivial fiber bundle with fibers diffeomorphic to the square $[-1, 1]^2$. In particular, $N_\varepsilon$ is a manifold with corners.

Let us inspect $X_{E_{\varepsilon}}$ along the $3$-faces of $N_\varepsilon$. They are of two kinds. Elements of the first group are given by the equations $G_i=\varepsilon c_G$; here the behavior of $X_{E_\varepsilon}$ is already known from Theorem \ref{thm_g}. The second group consists of $3$-faces arising from the boundary sides of $N$. So let us focus on this group.

By the definition, the compact boundary of $N$ can be finitely covered by the regular $0$-sets of some smooth locally supported functions $f_j(s_1, s_2)$. Then the intersections of these level sets give us corners of $N$. In particular, by the construction of $N$, it holds that 
$$\widetilde{\langle\widetilde{\nabla}f_j, X_{E_0} \rangle}\neq 0$$
along the level sets of $f_j$. Here, recall that $\widetilde{\langle\cdot,\cdot\rangle}$ is our notation for the standard metrics on coordinate charts $(s_1, s_2)$ or $(s_1, \theta_1, s_2, \theta_2)$.

Moreover, the functions $f_j$ give us canonically the covering of the remaining boundary sides of $N_\varepsilon$ by the regular $0$-sets of the functions $f_j\circ \pi_{s_1, s_2}(s_1, \theta_1, s_2, \theta_2)$.

We need to verify that the relations
\begin{equation}\label{eqn_sign_lift_0}
\hbox{sgn}\widetilde{\langle\widetilde{\nabla}(f_j\circ\pi_{s_1, s_2}), X_{E_\varepsilon} \rangle}\neq 0
\end{equation}
and
\begin{equation}\label{eqn_sign_lift}
\hbox{sgn}\widetilde{\langle\widetilde{\nabla}f_j, X_{E_0} \rangle}=\hbox{sgn}\widetilde{\langle\widetilde{\nabla}(f_j\circ\pi_{s_1, s_2}), X_{E_\varepsilon} \rangle}
\end{equation}
hold along the corresponding points of the level sets of $f_j$ and $f_j\circ\pi_{s_1, s_2}$. First, observe that clearly $\widetilde{\nabla}f_j$ and $\widetilde{\nabla}(f_j\circ\pi_{s_1, s_2})$ depend only on variables $s_1, s_2$, and $\widetilde{\nabla}f_j=(\pi_{s_1, s_2})_\ast\widetilde{\nabla}(f_j\circ\pi_{s_1, s_2})$. Also, by Remarks \ref{rem_grad_0} and \ref{rem_grad_comp} and Corrollary \ref{cor_perturb_grad_like} the $s_i$-components of $\nabla E_0$ and $\nabla E_\varepsilon$ differ only by $O(\varepsilon)$-terms.

Since $\partial N$ is compact, we can bound each $\vert\langle\widetilde{\nabla}f_j, X_{E_0} \rangle\vert$ from bellow by some constant $c_0>0$. Moreover, since $(\R/T\mathbb{Z}\times S^1)^2$ is also compact, it holds that
$$\Vert X_{E_0}-(\pi_{s_1, s_2})_\ast X_{E_\varepsilon}\Vert_{\infty}\ll c_0,$$
provided that $\varepsilon>0$ is sufficiently small. And the relations (\ref{eqn_sign_lift_0}) and (\ref{eqn_sign_lift}) hold. Hence, the assumptions of Morse homology Theorem \ref{thm_Morse_homology} are satisfied.

\textit{Contractibility:}

By the above, it holds that
$$\partial_-N_\varepsilon\cong\partial_-(N\times[-1, 1]^2)=(\partial_-N\times[-1, 1]^2)\cup(N\times[-1, 1]\times\lbrace\pm1\rbrace).$$

In particular, $\partial_-{N_\varepsilon}$ is contractible, which is also the case for $N_\varepsilon$.

\textit{The case $(iii.)$:}
Now, this is immediate from the construction of $N$ and $N_\varepsilon$.
\end{proof}

\begin{lemma}\label{thm_base_ift} If $\varepsilon>0$ is sufficiently small, $\Phi^\varepsilon_{p_0, q_0}$ satisfies property $(iv.)$ from Theorem \ref{thm_almost_bijection}. In addition, the equality of moduli spaces $(\ref{eqn_no_new_traj})$ holds.
\end{lemma}

\begin{proof}
By contradiction. Let us assume that there is a sequence $\lbrace u_{\varepsilon_n}\rbrace$ of $\varepsilon_n$-solutions with $\varepsilon_n\rightarrow 0$ that are not images of the multivalued maps $\Phi^{\varepsilon_n}_{p_0, q_0}$. Which is by the proof Lemma \ref{thm_ambient_ift} equivalent to the assumption that $u_{\varepsilon_n}$ are not in any of the sets $N_{\varepsilon_n}$.

On the other hand, recall the uniform bounds from Theorem \ref{thm_unif_bound} on $u_{\varepsilon_n}$. Hence, by \cite[Thm 4.7]{frauenfelder2020moduli}, the $\pi_{s_1, s_2}(u_{\varepsilon_n})$ converge in Floer-Gromov $C^1_{loc}$ sense to some broken flow line from $p_0$ to $q_0$. But since $X_{E_0}$ is Morse-Bott-Smale, by the index reason, the limit is in fact a $0$-solution. We will denote it by $u$.  Hence, for $n\gg 0$ it holds that $\pi_{s_1, s_2}(u_{\varepsilon_n})\subset N$. Recall also that $u$ is the only $0$-solution inside $N$.

Now, since $\pi_{s_1, s_2}(u_{\varepsilon_n})\subset N$, in particular $u_{\varepsilon_n}$ avoids special and diagonal points. Hence, we can apply Lemma \ref{lemma_solution_strip} and thus, each $u_{\varepsilon_n}\subset G_{K, \varepsilon_n}\hbox{-strip}$. But, by the construction, $N_{\varepsilon_n}=\pi_{s_1, s_2}^{-1}(N)\cap G_{K, \varepsilon_n}\hbox{-strip}$. Hence, for $n\gg 0$, $u_{\varepsilon_n}\subset N_{\varepsilon_n}$. Contradiction.

Since $G_{K, \varepsilon_n}\hbox{-strip}\subset M_{K, \varepsilon}$, we in fact also showed $(\ref{eqn_no_new_traj})$.
\end{proof}

\begin{rem} Note that in the proof of Lemma \ref{thm_base_ift} for the Floer-Gromov convergence of $\pi_{s_1, s_2}(u_{\varepsilon_n})$ we needed in fact only the uniform bounds for derivatives of $s_i$-coordinates.

On the other hand, the use of the uniform bounds for derivatives of $\theta_i$-coordinates inside the strips tells us interesting information about the geometry of the trajectories $u_{\varepsilon_n}$ in $M_{K, \varepsilon}$. They Floer-Gromov $C^1_{loc}$ converge to the algebraic ODE, which turns out to be the slow trajectory between $\mathfrak{p}_0$ and $\mathfrak{q}_0$, see Lemma \ref{lemma_identif}. This can be seen as an indication for the use of the Multiple-time scale dynamics to study our problem, see Chapter \ref{ch:file9}.
\end{rem}

\chapter{Morse flow trees}
\label{ch:file8}
In this chapter, we introduce certain Morse flow trees that will be counted later in cord algebras for knots and tori. In \cite{Ng_2005}, the notion of Morse flow trees for knots was outlined. Later, in \cite{petrak2019definition}, the Morse flow trees were used in more detail for the Morse model of Cord algebra for knots with loop space coefficients. See also \cite{okamoto2024legendrian} for Morse chains with a coproduct that were motivated by the language of \cite{abouzaid2013symplectic}.

For further purposes, it will be useful for us to have various models of Morse flow trees and describe the translations between them. In more detail, motivated by \cite{mescher2016perturbedgradientflowtrees}, we define trees as preimages of certain evaluation maps. This will allow us later to control the trees under the deformations in the adiabatic limit problem. Next, in the configuration space $(\R/T\mathbb{Z})^2$, we also use an iterative construction of trees under the forward-flow, see \cite{petrak2019definition}. This will be suitable for achieving the transversality conditions. However, for tori we use in the configuration space $(\R/T\mathbb{Z}\times S^1)^2$ an iterative construction of trees under the backward-flow. This is slightly technical, but it will allow us to track the trees better and keep the dimension of tracked trajectories constant. 

\section{Trees for knots}

\begin{rem_not}\label{rem_graphs}
We would like to introduce few notions from graph theory. For more about graphs, see, for example \cite{Matouek2000KapitolyZD}.

A \textbf{rooted tree} $\mathcal{T}$ is a connected acyclic graph, where one vertex has been designated the \textbf{root}. $V(\mathcal{T}), E(\mathcal{T})$ will denote the sets of all vertices/edges of $\mathcal{T}$, respectively.

\textbf{Inner vertices} are the vertices of $\mathcal{T}$ that have the valence higher than $1$. We will assume that the \textit{valence of interior vertices is equal to $3$ and the valence of the root is equal to $1$}. The set of interior vertices will be denoted by $V_{int}(\mathcal{T})$. A \textbf{leaf} is any vertex of $\mathcal{T}$ that is not an inner vertex or the root. The set of leaves will be denoted by $V_{leaf}(\mathcal{T})$. Let $m=\vert V_{leaf}(\mathcal{T})\vert$.

The paths from the root induce the canonical orientation of $\mathcal{T}$.

We consider the tree $\mathcal{T}$ with an \textbf{ordering} that is an assignment of each leaf to the unique number from $\lbrace 1,\dots,m\rbrace$. Let $\mathcal{T}_1$ and $\mathcal{T}_2$ be two rooted trees with orderings. We say that $\mathcal{T}_1$ and $\mathcal{T}_2$ are \textbf{isomorphic} if there is a graph isomorphism from $\mathcal{T}_1$ to $\mathcal{T}_2$ that preserves leaves together with their ordering. An equivalence class of isomorphic rooted trees is called an \textbf{ordered tree} and will be still denoted by $\mathcal{T}$. The set of all ordered trees (with $m$ leaves) will be denoted by $\clubsuit$ ($\clubsuit_m$).

An edge $\emph{e}$ is called \textbf{interior} if the endpoints of $\emph{e}$ are interior vertices. The set of all interior edges will be denoted by $E_{int}(\mathcal{T})$. Using the orientation on $\mathcal{T}$, we will say that one of the endpoints of any edge $\emph{e}$ is the \textbf{incoming vertex} $v^{in}_\emph{e}$ and the other is the \textbf{outcoming vertex} $v^{out}_\emph{e}$.

For any vertex $v\in V_{int}(\mathcal{T})$ the orientation on $\mathcal{T}$ induces the unique \textbf{incoming edge} $e_v$ such that $v$ belongs to $\emph{e}_v$. Moreover, using the ordering on $\mathcal{T}$, we will distinguish the two remaining edges with endpoint at $v$ - the \textbf{lower edge} $\emph{e}^L_v$ and the \textbf{upper edge} $\emph{e}^U_v$. The convention is that the unique path from the root to the leaf, which goes through $\emph{e}^U_v$, will end in the leaf with a higher number than the path through $\emph{e}^L_v$. By $\emph{e}_r$ we will denote the unique edge that contains the root.
\end{rem_not}

\begin{rem}\label{rem_flag}Let $\mathcal{T}\in\clubsuit$. A \textbf{flag} in $\mathcal{T}$ is a pair $(v, \emph{e})$, where $v\in V_{int}(\mathcal{T})$ and $\emph{e}\in E(\mathcal{T})$. We can naturally collect flags in two ways: by the common edges and by the common interior vertices. This will motivate us for the reordering functions in Definitions \ref{defn_reord} and \ref{ev_eps_tree}. See also \cite{cieliebak2023chern}.
\end{rem}
 
\begin{rem_not}We will be using three notions for \textbf{chords on $K$}. If $\gamma:\R/T\mathbb{Z}\rightarrow \R^3$ is an arc reparametrization of $K$, then any chord can be seen as a pair $(s_1, s_2)\in(\R/T\mathbb{Z})^2$ or a vector $P=\gamma(s_2)-\gamma(s_1)$. Moreover, the chord can also be seen as a smooth map $c_{s_1, s_2}:[0, 1]\rightarrow\R^3$ defined by 
$$c_{s_1, s_2}(\ell):=(1-\ell)\gamma(s_1)+ \ell\gamma(s_2).$$
\end{rem_not}

\begin{defn}\label{defn_ribbon_tree_K}\cite{Ng_2005} Let $p_0\in Crit_1(E_0)\setminus \Delta_0$ and $q_0^1, \dots, q_0^m\in Crit_0(E_0)\setminus \Delta_0$. Let $\mathcal{T}\in\clubsuit_m.$ Let $X_{E_0}$ be a gradient-like vector field adapted to $E_0$. A \textbf{Morse flow tree from $p_0$ to $q_0^1, \dots, q_0^m$ modeled on $\mathcal{T}$ (in the configuration space $(\R/T\mathbb{Z})^2$)} is an ordered tree $\mathcal{T}$ together with the following assignment:
\begin{itemize}
\item[$(i.)$] If $m=1$ then the unique edge corresponds to a standard flow line of $\phi_{X_{E_0}}^t$ with $t\in\R.$

If $m>1$, then the interior edges correspond to \textcolor{blue}{positive finite length} partial flow lines of $\phi_{X_{E_0}}^t$. Then the other edges correspond to partial flow lines of $\phi_{X_{E_0}}^t$ parametrized by $t\in[0, \infty)$ or $t\in(-\infty, 0]$. 
\item[$(ii.)$] The root corresponds to $p_0$.
\item[$(iii.)$] $i$-th leaf is identified with $q_0^i$.
\item[$(iv.)$] The tree is allowed to bifurcate at the points where the flow lines consist of nontrivial chords intersecting $K$ by their interior.

Let $(s_1, s_2)$ be such a chord and $s_3$ corresponds to the intersection of $c_{s_1, s_2}\vert_{(0, 1)}$ with $K$. Then we jump to the chords $(s_1, s_3), (s_3, s_2)$ and continue to flow with $\phi_{X_{E_0}}^t$. Now, the triple $\lbrace(s_1, s_2),(s_1, s_3), (s_3, s_2)\rbrace$ is assigned to the corresponding interior vertex $v$ of the tree $\mathcal{T}$. The partial flow-line starting from $(s_1, s_3)$ will correspond to the edge $\emph{e}^L_v$ and the partial flow-line starting from $(s_3, s_2)$ will correspond to the edge $\emph{e}^U_v$.  See also Figure \ref{figure_bifurcation_oftree}.
\end{itemize}

The \textbf{set of all Morse flow trees from $p_0$ to $q_0^1, \dots, q_0^m$ modeled on $\mathcal{T}$} is denoted by $$\clubsuit_{X_{E_0}}(\mathcal{T}; p_0; q_0^1, \dots, q_0^m).$$
The set of all Morse flow trees for all choices of $m, \mathcal{T}, p_0\in Crit_1(E_0)\setminus\Delta_0$ and $q_0^1, \dots, q_0^m\in  Crit_0(E_0)\setminus\Delta_0$ will be denoted by $\clubsuit_{X_{E_0}}$.

If $m=1$, we mod out the natural $\R$-action, then 
$$\clubsuit_{X_{E_0}}(\mathcal{T}; p_0; q_0^1)=\mathcal{M}_{X_{E_0}}(p_0; q_0^1).$$

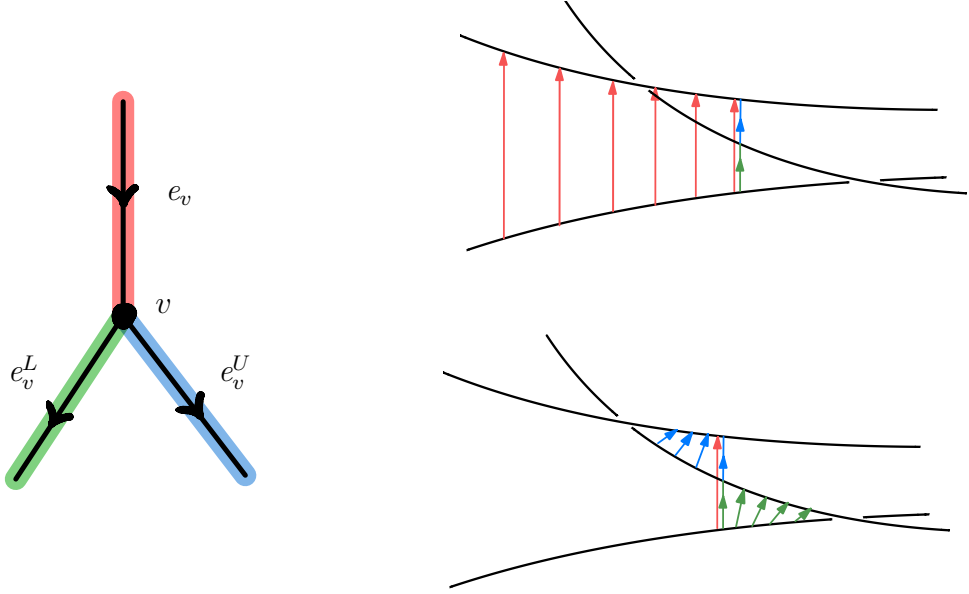
\begin{figure}[!htbp]
\labellist
\pinlabel $v$ at 148 700
\pinlabel $\emph{e}_v^L$ at 85 670
\pinlabel $\emph{e}_v^U$ at 180 670
\pinlabel $\emph{e}_v$ at 155 750
\endlabellist
\centering
\includegraphics[scale=0.84]{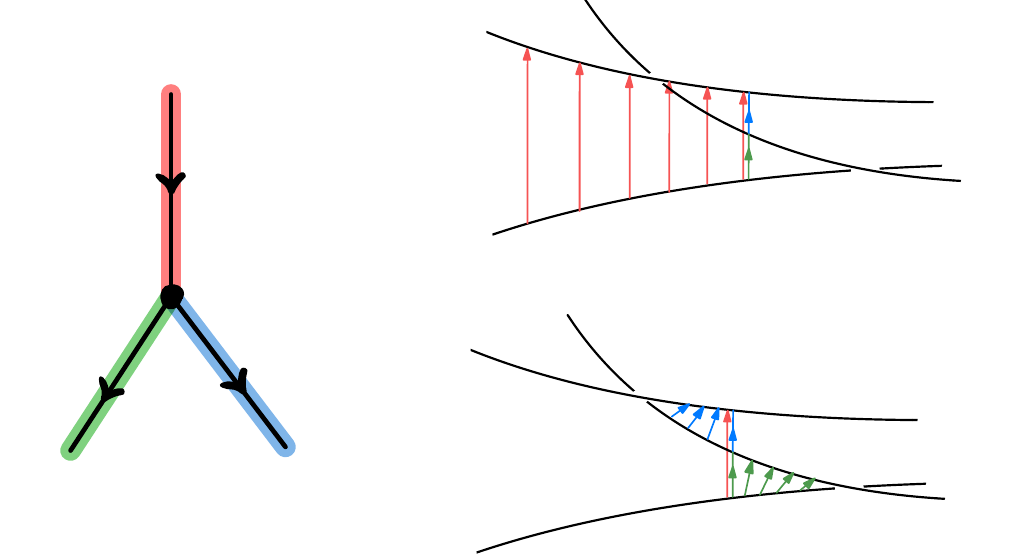}
\vspace{0.3cm}
\caption{On the left: a bifurcation of a Morse flow tree in the vertex $v$. $\emph{e}_v$ is the incoming edge and $\emph{e}_v^L, \emph{e}_v^U$ are the lower and upper edges, respectively. Recall that the orientation of the edges agrees with the flow of $X_{E_0}$. On the right: the realization of the bifurcation of the tree on the knot $K$ in $\R^3$.}
\label{figure_bifurcation_oftree}
\end{figure}
\end{defn}

\begin{rem}\label{rem_brokenK} If in Definition \ref{defn_ribbon_tree_K} we demand that at least one edge corresponds to a constant or a broken (partial) flow line of $\phi_{X_{E_0}}^t$, then the resulted tree will be called \textbf{broken Morse flow tree}.

We stress that even for the broken Morse flow trees we do not allow that $q_0^i\in\Delta_0$ for some $i\in\lbrace 1,\dots, m\rbrace.$
\end{rem}

\begin{lemma}\label{lemma_dim_tree_knot} Let $p_0\in Crit_1(E_0)\setminus \Delta_0$ and $q_0^1, \dots, q_0^m\in Crit_0(E_0)\setminus \Delta_0$.  Let $\mathcal{T}\in\clubsuit_m.$ If $K$ is a generic knot and $X_{E_0}$ is a generic perturbation of $-\nabla E_0$, then $\clubsuit_{X_{E_0}}(\mathcal{T}; p_0; q_0^1, \dots, q_0^m)$ is a zero dimensional compact manifold. In addition, there are no broken Morse flow trees, and for $m\gg 0$ the set $\clubsuit_{X_{E_0}}(\mathcal{T}; p_0; q_0^1, \dots, q_0^m)$ is empty.
\end{lemma}

\begin{rem}  Lemma \ref{lemma_dim_tree_knot} was done  in \cite[lem 2.23]{petrak2019definition} with a small mistake.  First, we are going to show Lemma \ref{lemma_dim_tree_knot} heuristically and later present (hopefully) a correct argument.
\end{rem}

\begin{hproof}
Let $m=1$. Since $\Delta_0$ is the global minimum of $E_0$, we can restrict the flow $\phi_{X_{E_0}}^t$ to $(\R/T\mathbb{Z})^2$ without a small open neighborhood of $\Delta_0$ such that $X_{E_0}$ is outward-pointing from the boundary. Then the case $m=1$ is precisely Remark \ref{rem_cmp_lines}. For $m>1$ we obtain $\clubsuit_{X_{E_0}}(\mathcal{T}; p_0; q_0^1, \dots, q_0^m)$ as an intersection of certain transversality conditions; let us compute
\begin{align*}
\dim \clubsuit_{X_{E_0}}&(\mathcal{T}; p_0; q_0^1, \dots, q_0^m)\\
\stackrel{(1)}{=}&Ind_{E_0}(p_0)+\sum_{i=1}^m(2-Ind_{E_0}(q_0^i))+\sum_{\emph{e}\in E_{int}(\mathcal{T})}(2+1)-\sum_{v\in V_{int}(\mathcal{T})} (2\cdot 2\, \textcolor{blue}{+\,1})\\
=&Ind_{E_0}(p_0)+\sum_{i=1}^m(2-Ind_{E_0}(q_0^i))+(m-2)(2+1)- (m-1)(2\cdot 2+1)\\
=&Ind_{E_0}(p_0)-\sum_{i=1}^m Ind_{E_0}(q_0^i)-1\\
=&0,
\end{align*}
where the terms in $(1)$ are obtained as follows. The first two summands are coming from the dimensions of the unstable and stable manifolds, respectively. The third summand is the contribution of interior edges. They contribute by the spaces of positive finite length trajectories. Now only the terms corresponding to constraints remain. The first group of these terms is realized by the transversality conditions between the stable and unstable manifolds, and the space of finite-length trajectories. The terms in the last group are determined by the fact that the condition for a chord to intersect the knot in its interior is of \textcolor{blue}{$\codim 1$}.

The ambient space is compact and by the dimension reasons there will be generically no nontrivial subsequential convergence phenomena (see \cite[thm 4.22]{mescher2016perturbedgradientflowtrees}), and hence $\clubsuit_{X_{E_0}}(\mathcal{T}; p_0,; q_0^1, \dots, q_0^m)$ is subsequentially compact. 

In more detail, the phenomenon of broken trees will be of too high codimension, so it will be avoidable for trees that satisfy certain generic conditions. To achieve the transversality, it will be more suitable for us to impose the perturbations of $-\nabla_{E_0}$ geometrically in the configuration space, which we will do later.

To see that there are no trees with too many bifurcations, one has to estimate how much any branch is shorter after a bifurcation.
\end{hproof}

\begin{defn}\label{defn_evK} Let us introduce the evaluation map 
$$\evK:[0, 1]\times (\R/T\mathbb{Z})^2\times (\R/T\mathbb{Z})^2\times (\R/T\mathbb{Z})^2\rightarrow (\R/T\mathbb{Z})^3\times \R^3$$ 
as
$$\evK:\begin{pmatrix}
\ell\\
s_1\\
s_2\\
y_1\\
y_3\\
z_2\\
z_3
\end{pmatrix}
\longmapsto
\begin{pmatrix}
s_1-y_1\\
s_2-z_2\\
y_3-z_3\\
\gamma(y_3)-(1-\ell)\gamma(s_1)-\ell\gamma(s_2)
\end{pmatrix},$$
where the coordinates $(s_1, s_2, y_1, y_3, z_2, z_3)$ are induced by the single arc length parametrization $\gamma:\R/T\mathbb{Z}\rightarrow K$. By Remark \ref{notation_on_circle}, $\evK$ is a well-defined smooth map. See Figure \ref{figure_R_coincidence}.
 
Next, $\widehat{\evK}$ will denote the restriction of $\evK$ to $(0, 1)\times((\R/T\mathbb{Z})^2\setminus\Delta_0)^3$. Then the subset $\mathcal{R}\subset (\R/T\mathbb{Z})^2$ is defined as
$$\mathcal{R}:=\pi_{s_1, s_2}\big(\widehat{\evK}^{-1}(0)\big)$$
and called the \textbf{space of chords meeting $K$ in their interior}. We also define two other sets
$$\mathcal{R}^L:=\pi_{y_1, y_3}\big(\widehat{\evK}^{-1}(0)\big)\hbox{ and }\mathcal{R}^U:=\pi_{z_2, z_3}\big(\widehat{\evK}^{-1}(0)\big).$$
\begin{figure}[!htbp]
\labellist
\pinlabel $c_{s_1, s_2}$ at 57 715
\pinlabel $c_{y_1, y_3}$ at 127 670
\pinlabel $c_{z_3, z_2}$ at 150 750
\endlabellist
\centering
\includegraphics[scale=0.84]{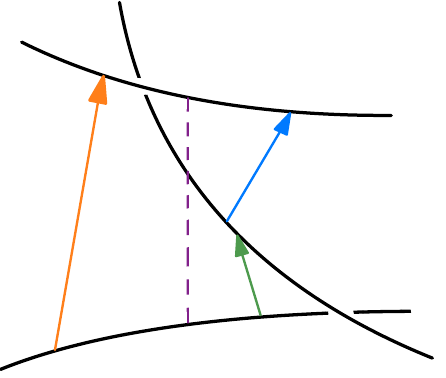}
\vspace{0.3cm}
\caption{Geometrically, we can see a zero of $\evK$ (or $\widehat{\evK}$) as a coincidence of $c_{s_1, s_2}$ with the union of $c_{y_1, y_3}$ and $c_{z_3, z_2}$. The coincidence is depicted as the dashed chord.}
\label{figure_R_coincidence}
\end{figure}
\end{defn}

\begin{rem} Due to the choice of coordinates, we can canonically see all $\mathcal{R}, \mathcal{R}^L, \mathcal{R}^U$ as subsets of the single $(\R/T\mathbb{Z})^2$. In particular, it is reasonable to admit that we chose more coordinates in Definition \ref{defn_evK} than we needed. However, the presented definition of $\evK$ will be better for various notions of trees $\clubsuit_{X_{E_0}}$.
\end{rem}

\begin{rem}\label{rem_d_evK}For further discussions, it will be useful to have described the differential 
$d\evK$ at a point $x=(\ell, s_1, s_2, y_1, y_3, z_2, z_3)$
$$d\evK(x)=
\begin{pmatrix}
0 & 1 & 0 & -1 & 0 & 0 & 0\\
0 & 0 & 1 & 0 & 0 & -1 & 0\\
0 & 0 & 0 & 0 & 1 & 0 & -1\\
\ell\gamma(s_1)-\ell\gamma(s_2) & (\ell-1)\dot{\gamma}(s_1) & -\ell\dot{\gamma}(s_2) & 0 & \dot{\gamma}(y_3) & 0 & 0\\
\end{pmatrix}.$$
Note that the last row of the matrix is in fact a $3\times 7$ matrix.
\end{rem}

The following lemma is a small extension of \cite[lem 7.10]{Cieliebak2016KnotCH}

\begin{lemma}\label{space_R}For a generic $K$ the set $\widehat{\evK}^{-1}(0)$ is a $1$-dimensional submanifold and the maps $\pi_{s_1, s_2}, \pi_{y_1, y_3}, \pi_{z_2, z_3}$ are immersions, when restricted to $\widehat{\evK}^{-1}(0)$.

Moreover, the closure $\overline{\mathcal{R}}$ is a compact immersed $1$-dimensional submanifold of $(\R/T\mathbb{Z})^2\setminus \Delta_0$ with finitely many transverse self-intersections and boundary such that:
\begin{itemize}
\item[$(i.)$]The self-intersections consist of pairs $(s_1, s_2)\in \mathcal{R}$ such that $c_{s_1, s_2}|_{(0, 1)}$ intersects $K$ twice. We will denote the set of these pairs by $\lozenge \mathcal{R}$.
\item[$(ii.)$] The boundary of $\overline{\mathcal{R}}$ are the special pairs (at special pairs $\overline{\mathcal{R}}$ is not self-intersecting).
\end{itemize} 
\end{lemma}

\begin{lemma}\label{cor_maps_R} Let $K$ be generic. Then the closures $\overline{\mathcal{R}^L}, \overline{\mathcal{R}^U}$ are compact immersed $1$-dimensional submanifolds of $(\R/T\mathbb{Z})^2$ with boundaries given by transverse intersections with $\Delta_0$ together with $y_1$- and $z_2$-special points, respectively. Along $\Delta_0$, each of $\overline{\mathcal{R}^L}, \overline{\mathcal{R}^U}$ has no self-intersections.

If we canonically represent $\overline{\mathcal{R}^L}, \overline{\mathcal{R}^U}$ as submanifolds of the single $(\R/T\mathbb{Z})^2$, then they are symmetric along the axis $\Delta_0$. See Figure \ref{figure_R_curves}.

Moreover, there are well-defined maps
\begin{align*}
\begin{split}
\varphi^L: \mathcal{R}\setminus\lozenge\mathcal{R} &\longrightarrow \mathcal{R}^L \hspace{3.9cm} \varphi^U: \mathcal{R}\setminus\lozenge\mathcal{R} \longrightarrow \mathcal{R}^U\\
\pi_{s_1, s_2}(x) \,&\longmapsto\pi_{y_1, y_3}(x)\,\,\quad\quad\quad\raisebox{16pt}{ and }  \quad\quad \quad\,  \pi_{s_1, s_2}(x) \,\longmapsto\pi_{z_2, z_3}(x),
\end{split}
\end{align*}
where $x\in\evK^{-1}(0)\setminus\pi^{-1}_{s_1, s_2}(\lozenge\mathcal{R})$.
\begin{figure}[!htbp]
\vspace{0.2cm}
\labellist
\pinlabel $0$ at -10 -10
\pinlabel $T$ at 220 -10
\pinlabel $T$ at -10 212
\pinlabel $(\R/T\mathbb{Z})^2$ at 220 235
\endlabellist
\centering
\includegraphics[scale=0.84]{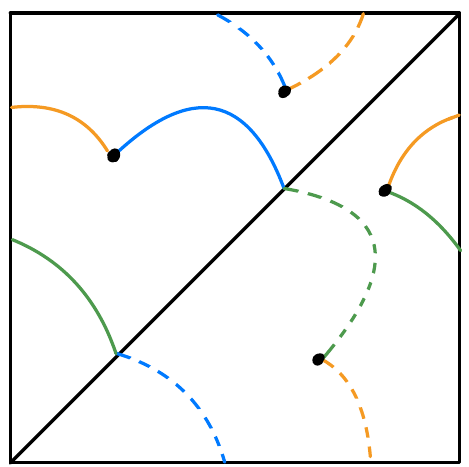}
\vspace{0.4cm}
\caption{A visualization of $\textcolor{orange}{\mathcal{R}}, \textcolor{teal}{\mathcal{R}^L}, \textcolor{blue}{\mathcal{R}^U}$ in the single configuration space $(\R/T\mathbb{Z})^2.$ Black dots represents special pairs. The groups of full and dashed curves correspond to two different connected components of $\evK^{-1}(0)$.}
\label{figure_R_curves}
\end{figure}
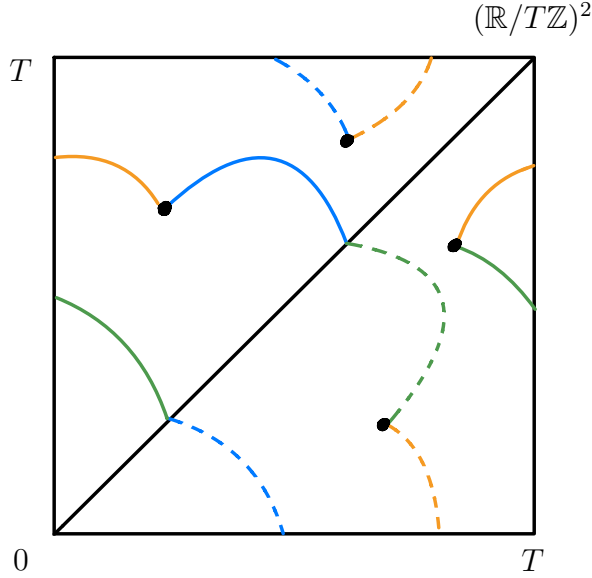
\end{lemma}

\begin{proof}
The lemma is a consequence of Lemma \ref{space_R} (i.e. \cite[lem 7.10]{Cieliebak2016KnotCH}). We remark that in the proof of Lemma \ref{space_R}, a small neighborhood of $\partial\overline{\mathcal{R}}$ in $\overline{\mathcal{R}}$ was explicitly constructed as a \textit{single} regular curve emanating from each special pair. This will in fact imply the transvesrality of $\overline{\mathcal{R}^L}, \overline{\mathcal{R}^U}$ with respect to $\Delta_0$.
\end{proof}

\begin{lemma}\label{cor_transverse_hitR} Let $K$ be generic. If $x\in \widehat{\evK}^{-1}(0)$, then $\hbox{rank}(d\evK(x))=6$.

Moreover, let $v\in T_{\pi_{s_1, s_1}(x)}(\R/T\mathbb{Z})^2$. Then $v\pitchfork\mathcal{R}$, iff the vectors 
$$\lbrace (\evK)_\ast v, \dot{\gamma}(y_3), \ell\gamma(s_2)-\ell\gamma(s_1)\rbrace$$ are linearly independent.
\end{lemma}

\begin{proof}
The lemma can be seen as a corollary of the proof of Lemma \ref{space_R} (i.e. \cite[lem 7.10]{Cieliebak2016KnotCH}), where one had used twice Thom's Transversality Theorem \ref{thm_thom}.

In more detail. The condition ``four vectors from $\R^3$ lie in a plane'' is of codim $2$. Hence, from the matrix $d\evK(x)$ (Remark \ref{rem_d_evK}), we see that the condition ``$\hbox{rank}(d\evK(x))<6$, where $x\in\widehat{\evK}^{-1}(0)$'' is of codim $3+3+2=8$. Thus, it is avoidable for generic knots.

In the second statement, one obtains analogously that for a generic knot the vectors $\dot{\gamma}(y_3), \ell\gamma(s_2)-\ell\gamma(s_1)$ are linearly independent. Indeed, the condition ``two vectors from $\R^3$ are linearly dependent'' is of codim $2$. Finally, since $T_x(\evK^{-1}(0))=\ker(d\evK(x))$, $(\evK)_\ast v$ could not lie in the $span$ of $\dot{\gamma}(y_3), \ell\gamma(s_2)-\ell\gamma(s_1)$. So the second statement follows.
\end{proof}

\begin{lemma}\label{lem_R_crit} For a generic knot $K$ it holds that 
$$Crit(E_0)\cap \overline{\mathcal{R}}=Crit(E_0)\cap \mathcal{R}^L=Crit(E_0)\cap \mathcal{R}^U=\emptyset.$$
\end{lemma}

\begin{proof}
By Lemma \ref{lemma_generic_morse_knot} we know that for a generic $K$ the set $Crit(E_0)\setminus\Delta_0$ consists of a finite number of nondegenerate critical points. In particular, $Crit(E_0)\setminus\Delta_0$ is stable under perturbations of $K$ (in the sense of Stability Lemma \ref{lem_stability} for the perturbations of the map $(\gamma\times\gamma)^\ast E$). Hence, for a generic knot, we can assume that any two binormal chords are not parallel. In particular, if $x\in\widehat{\evK}^{-1}(0)$, then at most one point of $\pi_{s_1, s_2}(x), \pi_{y_1, y_3}(x), \pi_{z_2, z_3}(x)$ lies in $Crit(E_0)$. Then the lemma follows from perturbations of $K$ as in \cite[lem 2.17]{petrak2019definition}.
\end{proof}

\begin{lemma}\label{lem_R_perturb}Let $K$ be a generic knot and $X_{E_0}$ be a gradient-like vector field. Then there is a gradient-like vector field $Y_{E_0}$ such that $Y_{E_0}$ is tangent to $\overline{\mathcal{R}}$ at only a finite number of the points and $X_{E_0}$ and $Y_{E_0}$ are $C^1$ close.
\end{lemma}

\begin{proof}
By Lemma \ref{space_R} we know that $\overline{\mathcal{R}}$ is a compact immersed $1$-manifold with boundary and a finite number of transverse self-intersections of the covering degree $2$. Hence there is a finite covering $\lbrace R_i\rbrace_{i=1, \dots,k}$ of $\overline{\mathcal{R}}$ by compact connected $1$-manifolds $R_i$ embedded in $(\R/T\mathbb{Z})^2$. Moreover, there are open subsets $U_i\subset (\R/T\mathbb{Z})^2$ such that $R_i\subset U_i$ and $\partial\overline{U_i}\pitchfork \overline{\mathcal{R}}$. 

Let us put $Y_0:=X_{E_0}$. By an induction in $i$ we would like to construct perturbed vector fields $Y_i$ such that $Y_i\pitchfork \overline{\mathcal{R}}$ along $R_i$ and $Y_i=Y_{i-1}$ on $(\R/T\mathbb{Z})^2\setminus \overline{U_i}$ (see below). Then we put $Y_{E_0}:=Y_k$. If for each $i$ the perturbation is sufficiently small in $C^\infty$ topology (relative to the perturbations at the steps $1,\dots, i-1$), then by Stability Lemma \ref{lem_stability} the vector field $Y_{E_0}$ has the desired tangency properties and is $C^1$ close to $X_{E_0}$. Moreover, by Lemma \ref{lem_R_crit} $Y_{E_0}$ is gradient-like.

Let us discuss the $i$-th perturbation. For simplicity of the notation, we assume that $R_i$ does not contain any self-intersection of $\overline{\mathcal{R}}$. Let $R^{ext}_i$ be a small open extension of $R_i$ inside $\overline{\mathcal{R}}$ (if $R_i\cap \partial \overline{\mathcal{R}}\neq\emptyset$, then we consider, instead of $\overline{\mathcal{R}}$, any small auxiliary extension of $\overline{\mathcal{R}}$ beyond the boundary $\partial\overline{\mathcal{R}}$). Now,  we can assume that a small tubular neighborhood $\nu(R^{ext}_i)$ induces a chart $\varphi:\nu(R^{ext}_i)\rightarrow\R_{x, y}$ such that 
$$\overline{U_i}=[-2, 2]\times[-1, 1],\quad R^{ext}_i\cap\overline{U_i}=[-2, 2]\times\lbrace0\rbrace,\quad R_i=[-1, 1]\times\lbrace 0\rbrace.$$
Let $\widetilde{Y}_{i-1}:=\varphi_\ast Y_{i-1}\varphi^{-1}$. Then by Relative Thom's transversality Theorem \cite[thm A.4]{petrak2019definition} there is a $C^1$ small function $\rho:[-2, 2]\times\lbrace0\rbrace\rightarrow \R$ such that $dy(\widetilde{Y}_{i-1})+\rho\pitchfork 0$ along $[-2, 2]\times\lbrace0\rbrace$ and $\rho(\lbrace\pm 2\rbrace\times\lbrace0\rbrace)=0.$ Now we can extend $\rho$ to the whole $\overline{U}_i$ as $\sigma(y)\rho(x)$, where $\sigma(y)$ is a small bump function which is supported on $[-1, 1]$. Then $\sigma\rho\partial_y$ determines the desired perturbation of $Y_{i-1}$.
\end{proof}

\begin{myproof}{Lemma}{\ref{lemma_dim_tree_knot}} Lemma \ref{space_R} and Corrollary \ref{cor_maps_R} give us a recipe, how to use the \textit{forward flow} $\phi^t_{X_{E_0}}$ and iteratively construct in the configuration space $(\R/T\mathbb{Z})^2$ all trees of $\clubsuit_{X_{E_0}}$.

That is the following. We consider $p_0\in Crit_1(E_0)\setminus\Delta_0$ and $x\in W^u_{X_{E_0}}(p_0)$ that is close to $p_0$. Then we continue flowing with $\phi^t_{X_{E_0}}(x)$ until we hit $\overline{\mathcal{R}}$, ideally in $\mathcal{R}\setminus\lozenge\mathcal{R}$. Then we can bifurcate and jump with the maps $\varphi^{L/U}$. After that, we continue by flowing and repeat. 

Hence, what we need to achieve is that each potential intersection of the (partial) flow trajectory with $\overline{\mathcal{R}}$ is transverse and avoid $\partial \overline{\mathcal{R}}, \lozenge\mathcal{R}$. Moreover, we need to avoid $Crit(E_0)$ besides the root and the leaves. About $\Delta_0$ and $\partial\overline{\mathcal{R}}$, we do not need to worry. Indeed, we are not counting trees with leaves in $\Delta_0$, and $\Delta_0$ is the global minimum of $E_0$. So, we can take a small regular index pair $(N, \emptyset)$ for $\Delta_0$ and study the flow only on $(\R/T\mathbb{Z})^2\setminus Int(N).$ In particular, any bifurcation will shrink the chords at least by some positive constant. Since the gradient-like vector field flows downhill, the number of potential bifurcations is bounded from above.

By Lemma \ref{lem_R_perturb}, we can assume that there is only a \textit{finite} number of problematic points along the whole $\overline{\mathcal{R}}$. Now, a similar \textit{finite} inductive argument as in Corollary \ref{cor_perturb_grad_like} will result in a desired perturbed gradient-like vector field. We remark that the local perturbations around the intersections with $\overline{\mathcal{R}}$ are described in \cite[lem 2.20]{petrak2019definition}, see also Lemma \ref{lem_perturb_petrak}.
\end{myproof}

\begin{rem}\label{rem_tree_standard_set} Since by Lemma \ref{space_R} $\partial\overline{\mathcal{R}}$ corresponds to the special points, we achieved in Lemma \ref{lemma_dim_tree_knot} that the trees $\clubsuit_{X_{E_0}}$ (when realized in $(\R/T\mathbb{Z})^2$) all lie in a standard set $S_K$. Here we implicitly assumed that the constant $\delta_K>0$ for $S_K$ is sufficiently small, see Remark \ref{rem_delta_standard}.
\end{rem}

\begin{setup} Let us fix $m\in \mathbb{N}, \mathcal{T}\in\clubsuit_m, p_0\in Crit_1(E_0)\setminus\Delta_0$ and $q_0^1, \dots, q_0^m\in Crit_0(E_0)\setminus \Delta_0$ and a gradient-like vector field $X_{E_{0}}$.
\end{setup}

The following exposition is motivated by \cite{mescher2016perturbedgradientflowtrees}, see also \cite{abouzaid2010topologicalmodelfukayacategories, mazuir2022higheralgebraainftyomega}.

\begin{rem}Recall that $m=\vert V_{leaf}(\mathcal{T})\vert=\vert V_{int}(\mathcal{T})\vert+1=\vert E_{int}(\mathcal{T})\vert+2$.
\end{rem}

\begin{defn}We define a map
\begin{align*}
i_{time}:(\R/T\mathbb{Z})^2\times \R_+&\longrightarrow (\R/T\mathbb{Z})^2\times (\R/T\mathbb{Z})^2\\
(s_1, s_2, t)\hspace{0.8cm}&\longmapsto \hspace{0.1cm}(s_1, s_2, \phi^t_{X_{E_{0}}}(s_1, s_2)).
\end{align*}
(Here $\R_+=(0, \infty)$.) The domain $(\R/T\mathbb{Z})^2\times \R_+$ will be called as the \textbf{space of positive finite length trajectories}.
\end{defn}

\begin{rem} Since $\frac{d}{dt}\vert_{t=t_0}\phi^t_{X_{E_{0}}}(s_1, s_2)=X_{E_{0}}(\phi^{t_0}_{X_{E_{0}}}(s_1, s_2))$, the map $i_{time}$ has full rank, when restricted outside of $Crit(E_0)\times\R_+$. In fact, because $X_{E_0}$ is a gradient-like vector field, the map $i_{time}$ is an embedding when restricted outside of $Crit(E_0)\times\R_+$.
\end{rem}

\begin{defn}\label{defn_incl_can}The map
\begin{align*}
i_{can}:\quad&\prod_{v\in V_{int}(\mathcal{T})}[0, 1]\times W^{u}_{X_{E_0}}(p_0)\times\prod_{\emph{e}\in E_{int}(\mathcal{T})}\big[(\R/T\mathbb{Z})^2\times \R_+\big]\times\prod^{m}_{i=1}W^s_{X_{E_0}}(q_0^i) \longrightarrow\\
 &\prod_{v\in V_{int}(\mathcal{T})}[0, 1]\times (\R/T\mathbb{Z})^2\times\prod_{\emph{e}\in E_{int}(\mathcal{T})}\big[(\R/T\mathbb{Z})^2\times (\R/T\mathbb{Z})^2\big]\times\prod^{m}_{i=1}(\R/T\mathbb{Z})^2 
\end{align*}
is induced by $i_{time}$ and canonical inclusions on the remaining terms. (We assume that $E_{int}(\mathcal{T})$ and $V_{int}(\mathcal{T})$ are arbitrarily ordered.)
\end{defn}

\begin{defn}\label{defn_reord}The \textbf{reordering map} $\rho_{\mathcal{T}}:$
\begin{align*}
\hspace{-1cm}\prod_{v\in V_{int}(\mathcal{T})}[0, 1]\times(\R/T\mathbb{Z})^2\times\prod_{\emph{e}\in E_{int}(\mathcal{T})}(\R/T\mathbb{Z})^4\times\prod^{m}_{i=1}(\R/T\mathbb{Z})^2 &\longrightarrow\prod_{v\in V_{int}(\mathcal{T})}\big[[0, 1]\times(\R/T\mathbb{Z})^6\big],\\
((\ell_v)_{v\in V_{int}(\mathcal{T})}, a^0, (a_{\emph{e}}^{in}, a_{\emph{e}}^{out})_{\emph{e}\in E_{int}(\mathcal{T})}, a^1,\dots, a^m)\hspace{0.8cm}&\longrightarrow \hspace{0.5cm}(\ell_v, b_v, b_v^L, b_v^U)_{v\in V_{int}(\mathcal{T})}
\end{align*}
will be defined by permuting the coordinates as follows. Here letters ``$a$'' or ``$b$'' represent, in fact, pairs of coordinates.

To each pair of $a$-coordinates or $b$-coordinates, we would like to assign a flag. If two pairs of $a$-coordinates and $b$-coordinates will have assigned the same flag, then they will be canonically identified under $\rho_{\mathcal{T}}$. For this, we will use the notions from Remarks \ref{rem_graphs} and \ref{rem_flag}.

\begin{itemize}
\item To each of $a^0, a^1,\dots, a^m$ we assign the unique edge $\emph{e}$ such that the root or the $i$-th leaf belongs to $\emph{e}$ $(i\in\lbrace 1,\dots, m\rbrace)$. The assigned vertex will be the unique interior vertex that belongs to $\emph{e}$.
\item To each $a^{in/out}_\emph{e}$ we assign the edge $\emph{e}$ and the vertex $v^{in/out}_\emph{e}$.
\item To each of $b_v, b_v^L, b_v^U$ we assign the vertex $v$ and the edges $\emph{e}_v, \emph{e}_v^L, \emph{e}_v^U$, respectively.
\end{itemize}
\end{defn}

\begin{defn} Let $X_{E_0}$ be a gradient-like vector field. Then we define the evaluation map $\emph{Ev}_{\mathcal{T}; p_0; q_0^1, \dots, q_0^m}:$
$$\hspace{-0.9cm}\prod_{v\in V_{int}(\mathcal{T})}[0, 1]\times W^{u}_{X_{E_0}}(p_0)\times\prod_{\emph{e}\in E_{int}(\mathcal{T})}\big[(\R/T\mathbb{Z})^2\times \R_+\big]\times\prod^{m}_{i=1}W^s_{X_{E_0}}(q_0^i)\longrightarrow\prod_{v\in V_{int}(\mathcal{T})}\big[(\R/T\mathbb{Z})^3\times\R^3]$$
as
$$\emph{Ev}_{\mathcal{T}; p_0; q_0^1, \dots, q_0^m}:=\big(\prod_{v\in V_{int}}\evK\big)\circ\rho_{\mathcal{T}}\circ i_{can}.$$
\end{defn}

\begin{cor}\label{cor_bij_treesT}For $K$ generic, there is a canonical bijection 
$$\clubsuit_{X_{E_0}}(\mathcal{T}; p_0; q_0^1, \dots, q_0^m)\xlongleftrightarrow{1-1}(\emph{Ev}_{\mathcal{T}; p_0; q_0^1, \dots, q_0^m})^{-1}(0).$$

In particular, if $X_{E_0}$ is a generic gradient-like vector field from Lemma \ref{lemma_dim_tree_knot} (which was constructed by certain perturbations in the configuration space $(\R/T\mathbb{Z})^2$), then
\begin{equation}\label{eqn_trans_Big_evK}
\emph{Ev}_{\mathcal{T}; p_0; q_0^1, \dots, q_0^m}\pitchfork 0.
\end{equation} 
\end{cor}

\begin{proof}
The bijection follows from the construction of \textit{$\emph{Ev}_{\mathcal{T}; p_0; q_0^1, \dots, q_0^m}$}. Here we remark that \textit{$\emph{Ev}_{\mathcal{T}; p_0; q_0^1, \dots, q_0^m}$} does not count trees with edges of constant paths. Indeed, for interior edges it is clear, and for the remaining edges it follows from $\mathcal{R}\cap Crit(E_0)=\emptyset$. Also for every $v\in V_{int}(\mathcal{T})$ it holds that \textit{$(\emph{Ev}_{\mathcal{T}; p_0; q_0^1, \dots, q_0^m})^{-1}(0)\cap\pi_{\ell_v}^{-1}(\lbrace0, 1\rbrace)=\emptyset$}. This follows from the fact that $p_0, q_0^1,\dots, q_0^m\in Crit(E_0)\setminus\Delta_0$, see also Remark \ref{rem_tree_standard_set}. 

The transversality (\ref{eqn_trans_Big_evK}) then follows from Lemma \ref{cor_transverse_hitR}, see also Remark \ref{rem_d_evK} for the description of $d\evK$.
\end{proof}

We finish the section with a small discussion about the evaluation map $\evK$.

\begin{rem}\label{rem_restr_R}It is easy to see that $\pi_{s_1, s_2}(\evK^{-1}(0))=(\R/T\mathbb{Z})^2,$ which is too big space. We have also introduced the restriction $\widehat{\evK}$, where $\pi_{s_1, s_2}(\widehat{\evK}^{-1}(0))=\mathcal{R}$. By Lemma \ref{space_R} $\mathcal{R}$ has a nice geometric picture, and moreover allows us to describe all bifurcations of $\clubsuit_{X_{E_0}}$. However, for further work, it will be a bit inconvenient that the domain of $\widehat{\evK}$, and in particular also $\mathcal{R}$, are both non-compact. Luckily, we have seen in Remark \ref{rem_tree_standard_set} that in fact all bifurcations of $\widehat{\evK}$ have to appear outside of a small neighborhood of $\partial\overline{\mathcal{R}}.$ So it will be enough to consider only certain restrictions of $\mathcal{R}$. This motivates the following definition.
\end{rem}

\begin{defn}$\widehat{\widehat{\evK}}$ will denote the restriction of $\evK$ to $[0, 1]\times((\R/T\mathbb{Z})^2\setminus\nu_{\delta_{\mathcal{R}}}(\Delta_0))^3$, where $\nu_{\delta_{\mathcal{R}}}(\Delta_0)$ is a $\delta_{\mathcal{R}}$-radius (open) tubular neighborhood of $\Delta_0$ for some $\delta_{\mathcal{R}}>0$ small.

Then the canonical projections $\pi_{s_1, s_2}, \pi_{y_1, y_3}, \pi_{z_2, z_3}$ on $\widehat{\widehat{\evK}}^{-1}(0)$ will define the sets $\widehat{\widehat{\mathcal{R}}}, \widehat{\widehat{\mathcal{R}^L}}, \widehat{\widehat{\mathcal{R}^U}},$ respectively.
\end{defn}

\begin{rem}\label{rem_shorter_R_K} If $\delta_{\mathcal{R}}>0$ is small and $K$ generic, then by Lemma \ref{space_R} and Lemma \ref{cor_maps_R} the following holds. $\widehat{\widehat{\mathcal{R}}}$ is a compact immersed $1$-dimensional submanifold of $(\R/T\mathbb{Z})^2$ with boundary. Moreover, outside of small neighborhoods of special points $\mathcal{R}$ and $\widehat{\widehat{\mathcal{R}}}$ coincide and $\widetilde{d}_{Haus}(\mathcal{R}, \widehat{\widehat{\mathcal{R}}})\rightarrow 0$ as $\delta_{\mathcal{R}}\rightarrow 0$.
\end{rem}

\section{Trees for tori}
We would like to define Morse flow trees on $T_{K, \varepsilon}\times T_{K, \varepsilon}$ (on $M_{K, \varepsilon}$) in a similar flavor as was done in the knot case.

\begin{rem_not}Similarly to the knot case, we will be using three notions for \textbf{chords on $T_{K, \varepsilon}$}. If $\Gamma_\varepsilon:(\R/T\mathbb{Z}\times S^1)\rightarrow \R^3$ is the parametrization of $T_{K, \varepsilon}$ from Remark \ref{rem_tor_param}, then any chord can be seen as a quadruple $(s_1, \theta_1, s_2, \theta_2)\in(\R/T\mathbb{Z}\times S^1)^2$ or a as vector $\Gamma_\varepsilon(s_2, \theta_2)-\Gamma_\varepsilon(s_1, \theta_1)$. Moreover, the chord can also be seen as a smooth map $c_{s_1, \theta_1, s_2, \theta_2}^\varepsilon:[0, 1]\rightarrow\R^3$ defined by 
$$c_{s_1, \theta_1, s_2, \theta_2}^\varepsilon(\ell):=(1-\ell)\gamma(s_1)+ \ell\gamma(s_2)+\varepsilon\big[(1-\ell)v(s_1, \theta_1)+ \ell v(s_2, \theta_2)\big].$$
Here recall that in Remark \ref{rem_tor_param} we introduced $v(s_i,  \theta_i)=\cos(\theta_i)n(s_i)+\sin(\theta_i)b(s_i)$.
\end{rem_not}

\begin{defn}\label{defn_ribbon_tree_T} Let $\varepsilon\in(0, \varepsilon_{good}]$. Let $p_\varepsilon\in Crit_2(E_\varepsilon\vert_{M_{K, \varepsilon}\setminus\Delta_\varepsilon})$ and $q_\varepsilon^1, \dots, q_\varepsilon^m\in Crit_1(E_\varepsilon\vert_{M_{K, \varepsilon}\setminus\Delta_\varepsilon})$ such that $p_\varepsilon, q_\varepsilon^1, \dots, q_\varepsilon^m\in M_{K, \varepsilon}$. Let $\mathcal{T}\in\clubsuit_m.$ Let $X_{E_\varepsilon}$ be a gradient-like vector field adapted to $E_\varepsilon$. A \textbf{Morse flow tree from $p_\varepsilon$ to $q_\varepsilon^1, \dots, q_\varepsilon^m$ modeled on $\mathcal{T}$ (in the configuration space $(\R/T\mathbb{Z}\times S^1)^2$)} is an ordered tree $\mathcal{T}$ together with the following assignment:
\begin{itemize}
\item[$(i.)$] If $m=1$ then the unique edge corresponds to a standard flow line $\phi_{X_{E_0}}^t$ with $t\in\R.$

If $m>1$, then the interior edges correspond to \textcolor{blue}{positive finite length} partial flow lines of $\phi_{X_{E_\varepsilon}}^t$. The other edges correspond to partial flow lines of $\phi_{X_{E_\varepsilon}}^t$ parametrized by $t\in[0, \infty)$ or $t\in(-\infty, 0]$.
\item[$(ii.)$] The root corresponds to $p_\varepsilon$.
\item[$(iii.)$] $i$-th leaf is identified with $q_\varepsilon^i$.
\item[$(iv.)$] The tree is allowed to bifurcate at the points where the flow lines consist of nontrivial chords intersecting $T_{K, \varepsilon}$ by their interior.

Let $(s_1, \theta_1, s_2, \theta_2)$ be such a chord and $(s_3, \theta_2)$ corresponds to the intersection of $c_{s_1, \theta_1, s_2, \theta_2}^\varepsilon\vert_{(0, 1)}$ with $T_{K, \varepsilon}$. Then, we jump to the chords $(s_1, \theta_1, s_3, \theta_3)$ and $(s_3, \theta_3, s_2, \theta_2)$ and continue to flow with $\phi_{X_{E_\varepsilon}}^t$. Now, the triple $$\lbrace(s_1, \theta_1, s_2, \theta_2),(s_1, \theta_1, s_3, \theta_3), (s_3, \theta_3, s_2, \theta_2)\rbrace$$ is assigned to the corresponding interior vertex $v$ of the tree $\mathcal{T}$. The partial flow-line starting from $(s_1, \theta_1, s_3, \theta_3)$ corresponds to the edge $\emph{e}^L_v$ and the partial flow-line starting from $(s_3, \theta_3, s_2, \theta_2)$ corresponds to the edge $\emph{e}^U_v$.
\end{itemize}

The \textbf{set of all Morse flow trees from $p_\varepsilon$ to $q_\varepsilon^1, \dots, q_\varepsilon^m$ modeled on $\mathcal{T}$} is denoted by $$\clubsuit_{X_{E_\varepsilon}}(\mathcal{T}; p_\varepsilon; q_\varepsilon^1, \dots, q_\varepsilon^m).$$
The set of all Morse flow trees for all choices of $m, \mathcal{T}, p_\varepsilon\in Crit_2(E_\varepsilon\vert_{M_{K, \varepsilon}\setminus\Delta_\varepsilon})$ and $q_\varepsilon^1, \dots, q_\varepsilon^m\in  Crit_1(E_\varepsilon\vert_{M_{K, \varepsilon}\setminus\Delta_\varepsilon})$ will be denoted by $\clubsuit_{X_{E_\varepsilon}}$.

If $m=1$, we mod out the natural $\R$-action, then 
$$\clubsuit_{X_{E_\varepsilon}}(\mathcal{T}; p_\varepsilon; q_\varepsilon^1)=\mathcal{M}_{X_{E_\varepsilon}}(p_\varepsilon; q_\varepsilon^1).$$
\end{defn}

\begin{rem}\label{rem_brokenT} If in Definition \ref{defn_ribbon_tree_T} we demand that at least one edge corresponds to a constant or broken (partial) flow line of $\phi_{X_{E_\varepsilon}}$, then the resulted tree will be called \textbf{broken Morse flow tree}.

Again, as in Remark \ref{rem_brokenK}, we stress out that even for the broken Morse flow trees we are not allowing that $q_\varepsilon^i\in\Delta_{full}$ for some $i\in\lbrace 1,\dots, m\rbrace.$
\end{rem}

\begin{defn} We put
$$\clubsuit_{X_{E_\varepsilon}}^{out-out}:=\clubsuit_{X_{E_\varepsilon}}\cap M_{K, \varepsilon}.$$
I. e. $\clubsuit_{X_{E_\varepsilon}}^{out-out}$ consists of Morse flow trees of $\clubsuit_{X_{E_\varepsilon}}$ that lie completely in $M_{K, \varepsilon}$.
\end{defn}

\begin{lemma}\label{lemma_dim_tree_tree}Let $\varepsilon\in(0, \varepsilon_{good}]$. Let $p_\varepsilon\in Crit_2(E_\varepsilon\vert_{M_{K, \varepsilon}\setminus\Delta_\varepsilon})$ and $q_\varepsilon^1, \dots, q_\varepsilon^m\in Crit_1(E_\varepsilon\vert_{M_{K, \varepsilon}\setminus\Delta_\varepsilon})$ such that $p_\varepsilon, q_\varepsilon^1, \dots, q_\varepsilon^m\in M_{K, \varepsilon}$. If $K$ is a generic knot, $\varepsilon>0$ is small, $X_{E_\varepsilon}$ is a generic perturbation of $-\nabla E_\varepsilon$, then $\clubsuit_{X_{E_\varepsilon}}(\mathcal{T}; p_\varepsilon; q_\varepsilon^1, \dots, q_\varepsilon^m)$ is a zero dimensional compact manifold. In addition, there are no broken Morse flow trees, $\clubsuit_{X_{E_\varepsilon}}=\clubsuit_{X_{E_\varepsilon}}^{out-out}$ and for $m\gg 0$ the set $\clubsuit_{X_{E_\varepsilon}}(\mathcal{T}; p_\varepsilon; q_\varepsilon^1, \dots, q_\varepsilon^m)$ is empty.
\end{lemma}

\begin{rem}  Now, we are going to show a heuristic proof of the dimension count for $\clubsuit_{X_{E_\varepsilon}}(\mathcal{T}; p_\varepsilon; q_\varepsilon^1, \dots, q_\varepsilon^m)$. The precise proof of Lemma \ref{lemma_dim_tree_knot} will be the aim of Chapter \ref{ch:file9}, where we will compare $\clubsuit_{X_{E_\varepsilon}}$ with $\clubsuit_{X_{E_0}}$ using the Multiple time scale dynamics together with the compactness argument motivated by Chapter \ref{ch:adiab_conley}.
\end{rem}

\begin{hproof}
For $m=1$, the argument is the same as in Lemma \ref{lemma_dim_tree_knot}. For $m>1$ we obtain $\clubsuit_{X_{E_\varepsilon}}(\mathcal{T}; p_\varepsilon; q_\varepsilon^1, \dots, q_\varepsilon^m)$ as an intersection of certain transversality conditions; let us compute
\begin{align*}
\dim \clubsuit_{X_{E_\varepsilon}}&(\mathcal{T}; p_\varepsilon; q_\varepsilon^1, \dots, q_\varepsilon^m)\\
\stackrel{(1)}{=}&Ind_{E_\varepsilon}(p_\varepsilon)+\sum_{i=1}^m(4-Ind_{E_\varepsilon}(q_\varepsilon^i))+\sum_{\emph{e}\in E_{int}(\mathcal{T})}(4+1)-\sum_{v\in V_{int}(\mathcal{T})} (2\cdot 4\, \textcolor{blue}{+\,0})\\
=&Ind_{E_\varepsilon}(p_\varepsilon)+\sum_{i=1}^m(4-Ind_{E_\varepsilon}(q_\varepsilon^i))+(m-2)(4+1)- (m-1)(2\cdot 4)\\
=&Ind_{E_\varepsilon}(p_\varepsilon)-\sum_{i=1}^m Ind_{E_\varepsilon}(q_\varepsilon^i)+m-2\\
=&0,
\end{align*}
The step $(1.)$ is similar to that in Lemma \ref{lemma_dim_tree_knot} except that now, the constraint for a chord to intersect $T_{K, \varepsilon}$ with its interior is of \textcolor{blue}{$\codim 0$}.
\end{hproof}

\begin{defn}\label{defn_evE} Let $\varepsilon\in[0, \varepsilon_{good}]$ (i.e. $\varepsilon$ can be even $0$). Let us introduce the evaluation map 
$$\ev_\varepsilon: [0, 1]\times (\R/T\mathbb{Z}\times S^1)^2\times (\R/T\mathbb{Z}\times S^1)^2\times (\R/T\mathbb{Z}\times S^1)^2\rightarrow (S^1)^3\times(\R/T\mathbb{Z})^3\times \R^3$$ 
as
$$\begin{pmatrix}
\ell\\
s_1\\
\theta_1\\
s_2\\
\theta_2\\
y_1\\
\alpha_1\\
y_3\\
\alpha_3\\
z_2\\
\beta_2\\
z_3\\
\beta_3\\
\end{pmatrix}
\longmapsto
\begin{pmatrix}
\theta_1-\alpha_1\\
\theta_2-\beta_2\\
\alpha_3-\beta_3\\
s_1-y_1\\
s_2-z_2\\
y_3-z_3\\
\gamma(y_3)-(1-\ell)\gamma(s_1)-\ell\gamma(s_2)+\varepsilon\big[v(y_3, \alpha_3)-(1-\ell)v(s_1, \theta_1)-\ell v(s_2, \theta_2)\big]
\end{pmatrix}$$
where the coordinates $((s_1, \theta_1, s_2, \theta_2), (y_1, \alpha_1, y_3, \alpha_3), (z_2, \beta_2, z_3, \beta_3))$ are induced by a single parametrization $\Gamma_\varepsilon:\R/T\mathbb{Z}\times S^1\rightarrow T_{K, \varepsilon}$. Also, by Remark \ref{notation_on_circle} $\ev_\varepsilon$ is a well-defined smooth map. 

Next, $\widehat{\widehat{\ev_\varepsilon}}$ will denote the restriction of $\ev_\varepsilon$ to 
$$[0, 1]\times\big((\R/T\mathbb{Z}\times S^1)^2\setminus(\nu_{\delta_{\mathcal{R}}}(\Delta_{0})\times (S^1)^2)\big)^3$$ for some $\delta_{\mathcal{R}}>0$ small. Then the subset $\widehat{\widehat{\mathcal{R}_\varepsilon}}\subset (\R/T\mathbb{Z}\times S^1)^2$ is defined as
$$\widehat{\widehat{\mathcal{R}_\varepsilon}}=\pi_{s_1, \theta_1, s_2, \theta_2}\big(\widehat{\widehat{\ev_\varepsilon}}^{-1}(0)\big).$$
We also define two sets
$$\widehat{\widehat{\mathcal{R}_\varepsilon^L}}:=\pi_{y_1, y_3}\big(\widehat{\widehat{\ev_\varepsilon}}^{-1}(0)\big)\hbox{ and }\widehat{\widehat{\mathcal{R}_\varepsilon^U}}:=\pi_{z_2, z_3}\big(\widehat{\widehat{\ev_\varepsilon}}^{-1}(0)\big).$$
See also Figure \ref{figure_R_eps_coincidence}.
\begin{figure}[!htbp]
\labellist
\pinlabel $c_{s_1, s_2}$ at 57 715
\pinlabel $c_{y_1, y_3}$ at 127 670
\pinlabel $c_{z_3, z_2}$ at 150 750
\endlabellist
\centering
\includegraphics[scale=0.84]{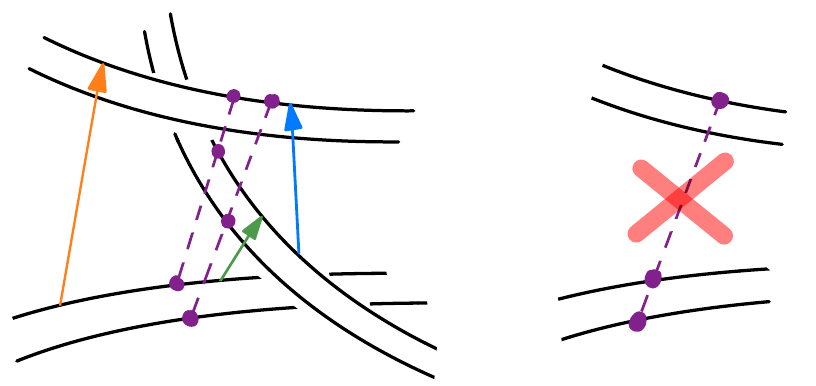}
\vspace{0.3cm}
\caption{\textit{On the left:} for $\varepsilon>0$ we can see geometrically a zero of $\widehat{\widehat{\ev_\varepsilon}}$ as a coincidence of \textcolor{orange}{$c^\varepsilon_{s_1, \theta_1, s_2, \theta_2}$} with the union of \textcolor{teal}{$c^\varepsilon_{y_1, \alpha_1, y_3, \alpha_3}$} and \textcolor{blue}{$c^\varepsilon_{z_3, \beta_3, z_2, \beta_2}$}. Two possible coincidences are depicted as the dashed chords with three marked points. Note that we are not imposing any outward-pointing condition.
\textit{On the right:} an example of a not allowed coincidence, since we are not counting intersections that are close in knot coordinates.}
\label{figure_R_eps_coincidence}
\end{figure}
\end{defn}

\begin{rem}\label{rem_simpl_ev0} Note that the family $\lbrace \ev_{\varepsilon}\rbrace_{\varepsilon\in[0, \varepsilon_{good}]}$ is smooth. Moreover
\begin{equation}\label{eqn_ev0}
\ev_0=\begin{pmatrix}
\theta_1-\alpha_1\\
\theta_2-\beta_2\\
\alpha_3-\beta_3\\
\evK\\
\end{pmatrix}
\end{equation}
and
$$\mathcal{R}_0\cong\mathcal{R}\times (S^1)^2.$$
\end{rem}

\begin{lemma}\label{lem_preim_treeT}For a generic knot and $\delta_{\mathcal{R}}>0$ sufficiently small, there is a $\varepsilon_0>0$ such that $\lbrace\widehat{\widehat{\ev_\varepsilon}}^{-1}(0)\rbrace_{\varepsilon\in[0, \varepsilon_0]}$ is a smooth family of isotopic $4$-dimensional compact submanifolds with boundary. 
\end{lemma}

\begin{proof}
Our aim is to apply the Stability Lemma \ref{lem_stability}. Hence, we need to verify that for a generic knot the map $\widehat{\widehat{\ev_0}}$ is stratum transverse to $0$. 

If we restrict $\widehat{\widehat{\ev_0}}$ to $Int(Dom(\widehat{\widehat{\ev_0}}))$, then the transversality is straightforward by relation $(\ref{eqn_ev0})$ and Lemma \ref{cor_transverse_hitR}. By Lemmata \ref{space_R} and \ref{cor_maps_R}, provided that $\delta>0$ is sufficiently small, it holds
\begin{itemize}
\item $\overline{\mathcal{R}^L}, \overline{\mathcal{R}^U}\pitchfork \overline{\nu_{\delta_{\mathcal{R}}}(\Delta_0)}.$
\item $\lbrace x\in\widehat{\widehat{\ev_\varepsilon}}^{-1}(0)\,\vert\,\pi_{y_1, y_3}(x)\in\overline{\nu_{\delta_{\mathcal{R}}}(\Delta_0)}\wedge\pi_{z_2, z_3}(x)\in\overline{\nu_{\delta_{\mathcal{R}}}(\Delta_0)}\rbrace=\emptyset$.
\item $\overline{\mathcal{R}}\cap \overline{\nu_{\delta_{\mathcal{R}}}(\Delta_0)}=\emptyset.$
\end{itemize}
In particular, $\widehat{\widehat{\ev_0}}^{-1}(0)\cap\pi_\ell^{-1}(\lbrace0, 1\rbrace)=\emptyset$. Which implies the transversality for the remaining cases.
\end{proof}

\begin{lemma}\label{lem_no_crit_bifurT}For a generic knot and $\delta_{\mathcal{R}}, \varepsilon>0$ sufficiently small it holds that
$$Crit(E_\varepsilon)\cap\widehat{\widehat{\mathcal{R}_\varepsilon}}=Crit(E_\varepsilon)\cap\widehat{\widehat{\mathcal{R}_\varepsilon^L}}=Crit(E_\varepsilon)\cap\widehat{\widehat{\mathcal{R}_\varepsilon^U}}=\emptyset.$$
\end{lemma}

\begin{proof}
By Lemma \ref{lem_R_crit} we know that for generic $K$ it holds that $Crit(E_0)\cap\widehat{\widehat{\mathcal{R}}}=\emptyset$. Hence by Lemma \ref{lemma_morse_corresp} (the relations $(\ref{eqn_of_crit_pt})$) we obtain that $Crit(E_\varepsilon)\cap\widehat{\widehat{\mathcal{R}_0}}=\emptyset$. Since $E_\varepsilon$ is Morse-Bott, the critical submanifolds are isolated. So the isotopy from Lemma \ref{lem_preim_treeT} finishes the proof, provided that $\varepsilon>0$ is sufficiently small. Other cases are analogous.
\end{proof}

\begin{lemma}\label{lemma_collect_intersect_treeT} For a generic knot and $\delta_{\mathcal{R}}, \varepsilon>0$ sufficiently small it holds that each bifurcation of $\clubsuit_{X_{E_\varepsilon}}$ appears at $\widehat{\widehat{\mathcal{R}_\varepsilon}}$.
\end{lemma}

\begin{proof}
Note that there is a $\delta_{K}>0$ small such that for any $\varepsilon>0$ sufficiently small
$$\big(\lbrace(s_1, \theta_1, s_2, \theta_2)\in(\R/T\mathbb{Z})^2\,\vert\,\widetilde{d}(s_1, s_2)\leq\delta_{K}\rbrace, \emptyset\big)$$
is a regular index pair of $\Delta_{full}$ under the flow $\phi_{X_{E_\varepsilon}}$. Also, similarly to the knot case, bifurcations are shortening the chords, so the forward-flow $\phi_{X_{E_\varepsilon}}$ flows along $\clubsuit_{X_{E_\varepsilon}}$ downhill. Hence, since elements of $\clubsuit_{X_{E_\varepsilon}}$ have no leaves in $\Delta_{full}$, for any $\varepsilon>0$ sufficiently small $\pi_{s_1, s_2}(\clubsuit_{X_{E_\varepsilon}})$  avoids a $\delta_K$-neighborhood of any special point. Now, we take $\delta_{\mathcal{R}}>0$ such that $\delta_{\mathcal{R}}\ll\delta_K$ which determines  $\widehat{\widehat{\mathcal{R}_\varepsilon}}$ with the desired property.
\end{proof}

\begin{defn}\label{ev_eps_tree} Let $\varepsilon\in(0, \varepsilon_{good}]$. Let $p_\varepsilon\in Crit_2(E_\varepsilon\vert_{M_{K, \varepsilon}\setminus\Delta_\varepsilon})$ and $q_\varepsilon^1, \dots, q_\varepsilon^m\in Crit_1(E_\varepsilon\vert_{M_{K, \varepsilon}\setminus\Delta_\varepsilon})$ such that $p_\varepsilon, q_\varepsilon^1, \dots, q_\varepsilon^m\in M_{K, \varepsilon}$. Let $\mathcal{T}\in\clubsuit_m.$ Let $X_{E_\varepsilon}$ be a gradient-like vector field adapted to $E_\varepsilon$.  

Then we define the evaluation map $\emph{Ev}_{\mathcal{T}; p_\varepsilon; q_\varepsilon^1, \dots, q_\varepsilon^m}^\varepsilon:$
\begin{align*}
\prod_{v\in V_{int}(\mathcal{T})}[0, 1]\times W^{u}_{X_{E_\varepsilon}}(p_\varepsilon)\times\prod_{\emph{e}\in E_{int}(\mathcal{T})}\big[(\R/T\mathbb{Z}\times S^1)^2\times \R_+\big]\times\prod^{m}_{i=1}W^s_{X_{E_\varepsilon}}(q_\varepsilon^i)\longrightarrow\\\prod_{v\in V_{int}(\mathcal{T})}\big[(S^1)^3\times(\R/T\mathbb{Z})^3\times\R^3]
\end{align*}
as
$$\emph{Ev}_{\mathcal{T}; p_\varepsilon; q_\varepsilon^1, \dots, q_\varepsilon^m}^\varepsilon=\big(\prod_{v\in V_{int}(\mathcal{T})}\ev_\varepsilon\big)\circ\rho_{\mathcal{T}}\circ i_{can}$$
where the definitions of the reordering function $\rho_\mathcal{T}$ and $i_{can}$ are analogous to Definitions \ref{defn_incl_can} and \ref{defn_reord}.
\end{defn}

\begin{cor}\label{cor_coresp_tree_tor}For $K$ generic, $\delta_{\mathcal{R}}, \varepsilon>0$ sufficiently small there is a canonical bijection 
$$\clubsuit_{X_{E_\varepsilon}}(\mathcal{T}; p_\varepsilon; q_\varepsilon^1, \dots, q_\varepsilon^m)\xlongleftrightarrow{1-1}(\emph{Ev}_{\mathcal{T}; p_\varepsilon; q_\varepsilon^1, \dots, q_\varepsilon^m})^{-1}(0).$$
\end{cor}

\begin{proof} Similarly as in Corollary \ref{cor_bij_treesT}, the bijection follows from the construction of \textit{$\emph{Ev}_{\mathcal{T}; p_\varepsilon; q_\varepsilon^1, \dots, q_\varepsilon^m}$}. Here we remark that \textit{$\emph{Ev}_{\mathcal{T}; p_\varepsilon; q_\varepsilon^1, \dots, q_\varepsilon^m}$} also does not count trees with edges of constant paths. Indeed, for interior edges it is clear, and for the remaining edges it follows from Lemmata \ref{lem_no_crit_bifurT} and \ref{lemma_collect_intersect_treeT}. Also by Lemma \ref{lemma_collect_intersect_treeT}, for every $v\in V_{int}(\mathcal{T})$ it holds that \textit{$(\emph{Ev}_{\mathcal{T}; p_\varepsilon; q_\varepsilon^1, \dots, q_\varepsilon^m})^{-1}(0)\cap\pi_{\ell_v}^{-1}(\lbrace0, 1\rbrace)=\emptyset$}.
\end{proof}

\begin{rem}\label{rem_backward_tree} Recall that in the knot case we also used an iterative construction of $\clubsuit_{X_{E_0}}$ under the \textit{forward-flow} $\phi_{X_{E_0}}^t$. I. e., roughly speaking, we track $W^u_{X_{E_0}}(p_0)$ for some $p_0\in Crit_1(E_0)$. When we hit $\mathcal{R}$, we bifurcate and continue flowing along the new branches and repeat.

However, for the torus case, this method is not very convenient. If we start naively to flow along some trajectory from $p_\varepsilon\in Crit_2(E_\varepsilon)$, then we realize that suddenly after the first bifurcation we have to follow two $1$-dimensional families of trajectories and so on. Hence, it might look that the dimensions of the tracked trajectories are exploding after each bifurcation. Fortunately, this all will get corrected once we reach the leaves. Contrary to the knot case, now the constraints to reach the stable manifolds at leaves are of $\codim 1$, which will reduce the dimension of the tracked trajectories. 

In order to avoid these temporary explosions of the dimension of tracked trajectories, we will be considering in the torus case an iterative construction under the \textit{backward-flow} $\phi_{X_{E_\varepsilon}}^{-t}$. I. e., roughly speaking, we start to track two $3$-dimensional stable manifolds of two leaves. When they hit $\widehat{\widehat{\mathcal{R}_\varepsilon^L}}$ and $\widehat{\widehat{\mathcal{R}_\varepsilon^U}}$, then they induce a subset $R_\varepsilon\subset\widehat{\widehat{\mathcal{R}_\varepsilon}}$. Our wish will be that $R_\varepsilon$ is a $2$-dimensional manifold, which is transverse to the flow $\phi_{X_{E_\varepsilon}}^{-t}.$ Then we will follow $\phi_{X_{E_\varepsilon}}^{(-\infty, 0]}(R_\varepsilon)$, which will be only $3$-dimensional manifold. Altogether, we started tracking two $3$-dimensional stable manifolds and continue to track one $3$-dimensional flow invariant manifold. This induces the recursion under the backward-flow. We will use this method in Chapter \ref{ch:file9}.
\end{rem}

We finish the chapter with one auxiliary lemma.

\begin{lemma}\label{lem_aux_intersect_tor_chords}Let $(\ell, (s_1, \theta_1, s_2, \theta_2), (y_1, \alpha_1, y_3, \alpha_3), (z_2, \beta_2, z_3, \beta_3))\in \widehat{\widehat{\ev_\varepsilon}}^{-1}(0)$, then the vector $\gamma(s_2)-\gamma(s_1)$ can be represented in two ways as
$$\frac{\Vert\gamma(s_2)-\gamma(s_1)\Vert}{\Vert\gamma(y_3)-\gamma(y_1)\Vert}\left[\big(\gamma(y_3)-\gamma(y_1)\big)+O(\varepsilon)\right]$$
or
$$\frac{\Vert\gamma(s_2)-\gamma(s_1)\Vert}{\Vert\gamma(z_2)-\gamma(z_3)\Vert}\left[\big(\gamma(z_2)-\gamma(z_3)\big)+O(\varepsilon)\right].$$
\end{lemma}
\begin{proof}
The lemma immediately follows from the geometric picture of $\widehat{\widehat{\ev_\varepsilon}}^{-1}(0)$.
\end{proof}
\chapter{Multiple time scale dynamics}
\label{ch:file9}
In this chapter, we will see the dynamics of $-\nabla E_\varepsilon$ as a fast-slow system on the configuration space $(\R/T\mathbb{Z}\times S^1)^2$ which in different time scales approaches the fast or slow subsystem as $\varepsilon\rightarrow 0$. If we ignore the singularities at special and diagonal points, the slow subsystem will look like $-\nabla E_0$ dynamics on $4$ sheets of $(\R/T\mathbb{Z})^2$ embedded into $(\R/T\mathbb{Z}\times S^1)^2$. In addition, the fast system looks like the flow in $\theta$-directions between the $4$ sheets, which are now normally hyperbolic critical manifolds (again, the situation is a bit more complicated once also the special and diagonal points get involved). For $\varepsilon>0$ small, we will be able to recover the dynamics of the fast-slow system from the combination of the fast and slow subsystems. The Fenichel theory and Exchange Lemmata will be the techniques that we will be using for this task. 

In addition, one of the sheets, where the slow dynamics is contained, lies in $M_{K, \varepsilon}.$ Hence, we will be able to construct the bijective correspondence between the $-\nabla E_0$-heteroclinic orbits on $(\R/T\mathbb{Z})^2$ and the $-\nabla E_\varepsilon$-heteroclinic orbits $M_{K, \varepsilon}$, provided that the rest points have low Morse indices and do not lie on diagonals. In fact, we will be able to show the bijective correspondence also for the Morse flow trees. For the surjectivity, we will adapt the compactness from Chapter \ref{ch:adiab_conley}, where we did the adiabatic limit with the Conley index. Later, in Chapter \ref{ch:cord}, the correspondence of trees will allow us to relate various cord algebras.

We remark that in Morse theory, the Multiple-time scale dynamics naturally appeared in the following places and was used for the bijective correspondence of (broken) heteroclinic orbits. It was used in \cite{schecter2014morse} for Lagrangian multipliers and in \cite{banyaga2013cascades} for the cascades between Bott-critical manifolds.

In Sections \ref{s:track}-\ref{s:exch} we recollect the notions of Fenichel theory and Exchange Lemmata. Then, finally, in Section \ref{s:appl_trees} we apply the theory to our setup.

\section{Basic tracking of trajectories}
\label{s:track}

\begin{Gronwall}\cite{hirsch2004differential}\label{thm_gronwall} Let $u\in C^0([0, T], \R_{\geq0})$. Assume that $C, L$ are non-negative constants such that
$$u(t)\leq C+\int_0^t L u(s) ds$$
for all $t\in[0, T]$. Then, for each such $t$ holds that
$$u(t)\leq C e^{L t}.$$
\end{Gronwall}

\begin{rem}\cite{hirsch2004differential}\label{rem_continuity_ode} A useful consequence of Gronwall's inequality \ref{thm_gronwall} is that the flow of $C^1$ vector field depends continuously on initial conditions.
\end{rem}

\begin{lemma}\cite[lem 3.2.7]{Wiggins1994NormallyHI}\label{lemma_general_track}Let $f$ and $f^{pert}$ be two $C^1$ vector fields on $\R^n$. Let us consider an open bounded set $U\subset\R^n$ and $\varepsilon>0$. Assume that $f$ and $f^{pert}$ are $C^0$ $\varepsilon$-close on $\overline{U}$, i. e.
\begin{equation}\label{eqn_assump_track}
\sup_{x\in\overline{U}}\Vert f(x)-f^{pert}(x)\Vert \leq \varepsilon.
\end{equation}
We fix $T>0$ and denote by $\phi^t$ and $\phi^t_{pert}$ the flows generated by $f$ and $f^{pert}$, respectively. Then for each $x\in \overline{U}$ and $|t|\leq T$ holds
$$\Vert \phi^t(x)-\phi^t_{pert}(x)\Vert =O(\varepsilon T e^{LT}),$$
provided that $\phi^t(x)$ and $\phi^t_{pert}(x)$ remain in the compact set $\overline{U}$. Here $L$ is the Lipschitz constant of $f$.
\end{lemma}

\begin{proof}
We compute
\begin{align*}
\Vert\phi^t(x)-\phi^t_{pert}(x)\Vert&=\left\Vert\phi^0(x)-\phi^0_{pert}(x)+\int_0^t f(\phi^s(x))-f^{pert}(\phi^s_{pert}(x))ds \right\Vert\\
&\leq\int_0^t\Vert f(\phi^s(x))-f^{pert}(\phi^s_{pert}(x))\Vert ds\\
&\leq\int_0^t\Vert f(\phi^s(x))-f(\phi^s_{pert}(x))+f(\phi^s_{pert}(x))-f^{pert}(\phi^s_ {pert}(x))\Vert ds\\
&\stackrel{(1)}{\leq} \varepsilon T+\int_0^t\Vert f(\phi^s_{pert}(x))-f^{pert}(\phi^s_{pert}(x))\Vert ds\\
&\stackrel{(2)}{\leq} \varepsilon T+\int_0^t L\Vert \phi^s_{pert}(x)-\phi^s_{pert}(x)\Vert ds\\
&\stackrel{(3)}{\leq} \varepsilon T e^{LT},
\end{align*}
where in the step $(1.)$ we used the assumption (\ref{eqn_assump_track}), in $(2.)$ we used the Lipschitz continuity of $f$, and finally in $(3.)$ we applied Gronwall's inequality \ref{thm_gronwall}.
\end{proof}

\begin{rem}\label{rem_smooth_data} Alternatively, one can use an Implicit function theorem argument and show a smooth dependence of $\phi^{t}(x)$ on the initial data. Here, one has to consider the vector fields from a bounded open subset of the space of bounded vector fields with the uniform $C^k$-norm. See \cite[appx A]{Eldering2013NormallyHI}.
\end{rem}

\section{Fenichel theory}

\textit{From now on in this chapter, we will consider that $k$ is a natural number bigger than $3$.}

\begin{defn} Let
\begin{equation}\label{eqn_gen_ODE}
\dot{z}=f(z)
\end{equation}
be a general autonomous ODE on $\R^n$ given by a $C^k$ function $F:\R^n\rightarrow\R^n$. Then $\phi^t$ will denote the flow induced by (\ref{eqn_gen_ODE}). Let $M$ be a connected compact $C^k$-manifold in $\R^n$ with corners ($M$ is potentially of $\codim>0$). $M$ is called:
\begin{itemize}
\item \textbf{Overflowing} if $\phi^t(p)\in M$ for every $p\in M$ and $t\leq0$, and in addition the vector field of (\ref{eqn_gen_ODE}) is strictly outward-pointing from $\partial M$ (recall the terminology for vectors on the boundary in Definition \ref{def_outward_point}).
\item \textbf{Inflowing} if $\phi^t(p)\in M$ for every $p\in M$ and $t\geq0$, and in addition the vector field of (\ref{eqn_gen_ODE}) is strictly inward-pointing from $\partial M$.
\item \textbf{Invariant} if $\phi^t(p)\in M$ for every $p\in M$ and $t\in \R$.
\item \textbf{Locally invariant} if for every $p\in M^{\circ}$ there exists some open interval $I_p$ such that $0\in I_p$ and $\phi^t(p)\in M$ for $t\in I_p$.
\end{itemize}
\end{defn}

\begin{rem}\cite{Fenichel71} We would like to inspect the stability of the \textit{overflowing} manifold $M$ under a perturbation of $f$. For this, we will first introduce the rate of growth of normal vectors under the linearised dynamics in backward time. And then we will impose a function measuring the flattening of orbits relative to $M$ by relating the normal and tangential directions in backward time.

More precisely, let $\mathcal{N}:=N(M, \R^n)$ be the normal bundle of $M$. Let $p\in M$. Then for any $v_0\in T_p M$ and $w_0\in \mathcal{N}_p$ we put
$$v_{-t}:=d\phi^{-t}(p)v_0\hbox{ and }w_{-t}:=(\pi_M\circ d\phi^{-t})w_0.$$

Now, we consider any Riemannian metric on $\R^n$ and define the following \textbf{generalized Lyapunov-type numbers} as
$$\nu(p)=\inf\left\lbrace a\,\Big|\,\lim_{t\rightarrow \infty}\frac{\Vert w_0\Vert }{\Vert w_{-t}\Vert }/a^t=0,\hbox{ for each }w_0\in \mathcal{N}_p\right\rbrace$$
and, if $\nu(p)<1$,
$$\sigma(p)=\inf\left\lbrace b\,\Big|\,\lim_{t\rightarrow \infty}\big(\Vert w_0\Vert ^b/\Vert v_0\Vert \big)/\big(\Vert w_{-t}\Vert ^b/\Vert v_{-t}\Vert \big)=0,\hbox{ for each }v_0\in T_p M, w_0\in \mathcal{N}_p\right\rbrace.$$
\end{rem}

\begin{rem}\label{rem_layp1}\cite{Fenichel71, Wiggins1994NormallyHI} $\nu$ and $\sigma$ are independent of choices of the Riemannian metric on $\R^n$ and $\mathcal{N}$. Moreover, $\nu$ and $\sigma$ are constant along orbits.

Also, for the computation (especially at equilibrium points), there are useful relations
\begin{align*}
\nu(p)&=\overline{\lim_{t\rightarrow\infty}}\Vert \pi_\mathcal{N}\circ d\phi^{t}(\phi^{-t}(p))|_{\mathcal{N}}\Vert^{1/t},\\
\sigma(p)&=\overline{\lim_{t\rightarrow\infty}}\frac{\log\Vert d\phi^{-t}(p)|_{TM}\Vert }{-\log\Vert \pi_\mathcal{N}\circ d\phi^{t}(\phi^{-t}(p))|_{\mathcal{N}}\Vert }.
\end{align*}
\end{rem}

\begin{thm}\label{thm_overlfow}\cite{Fenichel71, Wiggins1994NormallyHI, Fenichel1979GeometricSP} Let $M$ be an overflowing manifold for the system (\ref{eqn_gen_ODE}). Assume that $\nu(p)<1$ and $\sigma(p)<1/k$ for each $p\in M$.

Let the $C^k$-vector field $f^{\varepsilon}$ be $C^1$ $O(\varepsilon)$-close to $f$, where $0<\varepsilon\ll 1$. Then there is an overflowing $C^k$-manifold $M_\varepsilon$ for $f^\varepsilon$ which is diffeomorphic to $M$. In addition, $M$ and $M^\varepsilon$ are $O(\varepsilon)$-close in Hausdorff distance.
\end{thm}

\begin{sproof}The proof is based on the original idea of Hadamard from $1901$ of the graph transform, see \cite{Hasselblatt2017ErgodicTA}. First, we consider the disc bundle $\mathcal{D}_\delta$ over $M$ consisting of normal vectors of norm $\leq\delta$. This gives us local coordinates $\R^{\dim M}\times\R^{n-\dim M}$ that identify locally $M$ with $\R^{\dim M}\times \lbrace 0\rbrace$. For simplicity, we assume the coordinates are global. In particular, any section $v\in \mathcal{D}_\delta$ can be seen as a graph ($=:\text{graph}(v)$) over $\R^{\dim M}$.

Let $\mathcal{S}_\delta$ be the space of sections in $\mathcal{D}_\delta$ with Lipschitz constant $\delta$. Then we define a graph transform $G$ on $\mathcal{S}_\delta$ by
$$\text{graph}(v)\mapsto \phi^T_\varepsilon(\text{graph}(v)),$$
for some fixed $T>1$. If we ignore what overflows from the bundle $\mathcal{D}_\delta$ and $\delta, \varepsilon$ are sufficiently small, then $G$ is well-defined. Next, the generalized Lyapunov-type numbers $\nu, \sigma$ imply that $\phi^T$ is contracting $\delta$ of $D_\delta$. However, if $\delta$ and $\varepsilon$ are sufficiently small, then $G$ is also a contraction. Hence, by the Banach contraction mapping theorem, $G$ has a unique fixed point $v$. In particular, $\hbox{graph}(v)\subseteq\phi^T_\varepsilon(\hbox{graph}(v))$. Then an inductive argument, starting from small $t>0$, shows that $\hbox{graph}(v)\subseteq\phi^t_\varepsilon(\hbox{graph}(v))$ for all $t\geq0$.

In the end, we remarked that we worked in a slightly simplified setting. In order to show the regularity of $M^\varepsilon$, we should not have started with $D_\delta\subset \mathcal{N}(M, \R^n)$. We should have rather considered $D_\delta$ inside a so-called transversal bundle. That is a certain $C^k$ $\dim M$-bundle in $T\R^n|_{M}$ that is close to $\mathcal{N}(M, \R^n)$ (which is only $C^{k-1}$). For more details on the regularity of $M^\varepsilon$, see \cite{Wiggins1994NormallyHI}.
\end{sproof}

\begin{rem}\label{rem_ext_overflow} Reversing time brings a straightforward extension of Theorem \ref{thm_overlfow} to inflowing manifolds and, in particular, to boundary-less invariant manifolds.

In the case of invariant manifolds with boundaries, the situation is a bit more delicate. By the definition, the vector field $f$ is also defined on a small neighborhood $U$ of $M$. Now, we deform the vector field $f$ by adding certain small bump functions supported in $U\setminus M$. The purpose of this modification is that we will obtain a compact set $\widetilde{M}\supset M$, where $\widetilde{M}$ is inflowing or overflowing. The downside is that for the perturbation $f^\varepsilon$, the manifold $M^\varepsilon$ will be only locally invariant due to the potential leaving of the flow through the boundary. Luckily, in our applications, we will be able to recover a lot of the structure of the dynamics restricted to $M^\varepsilon$. For more details on this construction, see \cite{Fenichel1979GeometricSP, Wiggins1988GlobalBA}.
\end{rem}

\begin{rem}\label{rem_layp2}\cite{Fenichel71, Wiggins1994NormallyHI} Continuing the discussion from Remark \ref{rem_layp1}, we introduce new \textbf{generalized Lyapunov-type numbers}. Assume we have a $(C^k)$ splitting
$$T\R^n|_M=\mathcal{N}^u\oplus TM\oplus\mathcal{N}^s,$$
where we assume that $\mathcal{N}^u\oplus TM$ and $\mathcal{N}^s\oplus TM$ are invariant under $D\phi^t$ from the system (\ref{eqn_gen_ODE}) for $t\leq0$. Let us pick a Riemannian metric on $\R^n$.

Now, we put
\begin{align*}
\nu^u(p)&=\overline{\lim_{t\rightarrow\infty}}\Vert \pi_{\mathcal{N}^u}\circ d\phi^{t}(\phi^{-t}(p))|_{\mathcal{N}^u}\Vert ^{1/t},\\
\nu^s(p)&=\overline{\lim_{t\rightarrow\infty}}\Vert \pi_{\mathcal{N}^s}\circ d\phi^{-t}(\phi^{t}(p))|_{\mathcal{N}^s}\Vert^{1/t},\\
\sigma^u(p)&=\overline{\lim_{t\rightarrow\infty}}\frac{\log\Vert d\phi^{-t}(p)|_{TM}\Vert }{-\log\Vert \pi_{\mathcal{N}^u}\circ d\phi^{t}(\phi^{-t}(p))|_{\mathcal{N}^u}\Vert },\\
\sigma^s(p)&=\overline{\lim_{t\rightarrow\infty}}\frac{\log\Vert d\phi^{t}(p)|_{TM}\Vert }{-\log\Vert \pi_{\mathcal{N}^s}\circ d\phi^{-t}(\phi^{t}(p))|_{\mathcal{N}^s}\Vert }.
\end{align*}

These numbers are also independent of the Riemannian metric and constant along orbits.
\end{rem}

\begin{unstable}\label{thm_unstable}\cite{Fenichel71, Wiggins1994NormallyHI} Let $M$ be an overflowing $C^k$-manifold for the system (\ref{eqn_gen_ODE}) such that $\hbox{sup}_{p\in M}\hbox{max}\lbrace \nu^u(p), \nu^s(p), \sigma^u(p)\rbrace<1.$ Then there is an overflowing $C^k$-manifold $W^u$ tangent to $\mathcal{N}^u$ along $M$.

$W^u$ persists under perturbation by $C^k$-vector field $f^\varepsilon$ $C^1$ $O(\varepsilon)$-close to $f$. In more detail, there is an overflowing $C^k$-manifold $W^u_\varepsilon$ diffeomorphic to $W^u$ and $O(\varepsilon)$-close in Hausdorff distance. 
\end{unstable}

\begin{sproof}For the existence, we study the graph transform on sections of $\mathcal{N}^u\oplus\mathcal{N}^s\rightarrow\mathcal{N}^u$ given by the linearised flow. Hence, the Banach contraction mapping theorem gives us the existence of $W^u$ in a small neighborhood $U$ of $M$.

For the persistence, we do not need to compute generalized Lyapunov-type numbers for $W^u$. But we would rather use the fact that in $U$ the manifolds $W^u$ and $\mathcal{N}^u$ are close (in the coordinates induced by the splitting of $T\R^n|_M$). Hence, if $U$ is sufficiently small, then we obtain the appropriate bounds for the Banach contraction mapping theorem and thus $W^u_\varepsilon$. 
\end{sproof}

\begin{defn} Let $M$ be a compact invariant manifold without boundary. The splitting $T\R^n|_M=\mathcal{N}^u\oplus TM\oplus\mathcal{N}^s$ is called \textbf{normally hyperbolic}, if $$\hbox{sup}_{p\in M}\hbox{max}\lbrace \nu^u(p), \nu^s(p), \sigma^u(p), \sigma^s(p)\rbrace<1.$$ Then $M$ is called a \textbf{normally hyperbolic manifold}. We will denote by $\ell_u$ the dimension of fibers of $\mathcal{N}^u$ and by $\ell_s$ the dimension of fibers of $\mathcal{N}^s$.
\end{defn}

\begin{rem}\cite{Fenichel71} The Unstable Manifold Theorem \ref{thm_unstable} extends also to normally hyperbolic invariant manifolds. In this case, we obtain the existence of overflowing $W^u$ and inflowing $W^s$. Moreover, not only $W^u$ and $W^s$ persist under perturbation of $f$, but also $M$ persists. In this case $W^u_\varepsilon\cap W^s_\varepsilon=M_\varepsilon$, $W^u_\varepsilon$ is overflowing, $W^s_\varepsilon$ is inflowing and $M_\varepsilon$ is invariant. We will call $W^{s}_\varepsilon$ and $W^{u}_\varepsilon$ the \textbf{local stable and unstable manifolds} of $M_\varepsilon$, and denote them rather as $W^{s, loc}_\varepsilon(M_\varepsilon)$ and $W^{u, loc}_\varepsilon(M_\varepsilon)$, respectively. The unions of the images of $W^{s/u, loc}_\varepsilon(M_\varepsilon)$ under the flow for all times give us \textbf{global stable and unstable manifolds} of $M_\varepsilon$.
\end{rem}




\begin{example}\label{example_hyperbolicity} Let $F$ be a Morse function on a closed Riemannian manifold $(M, g).$ Let $X_F$ be a gradient like-vector field adapted to $F$ and $p_0\in Crit(F)$. Then, with respect to $X_F$, the point $p_0$ is a normally hyperbolic invariant manifold with the unique local unstable manifold $W^{u, loc}_{X_F}(p_0)$ (the uniqueness is with respect to the fixed radius of the corresponding disk in $\mathcal{N}^u$). 

Hence, if $Y^\varepsilon$ is a $C^k$ $O(\varepsilon)$-perturbation of $X_F$, then the unique perturbations of $p_\varepsilon$ and $W^{u, loc}_{Y^\varepsilon}(p_\varepsilon)$ are $C^k$ $O(\varepsilon)$-close to $p_0$ and $W^{u, loc}_{X_F}(p_0)$, respectively.

Let us moreover assume that the perturbation $Y^\varepsilon$ of $X_F$ is $C^k$ smooth in $\varepsilon\in[-\varepsilon_0, \varepsilon_0]$, where $\varepsilon_0>0$ is small. We can then study perturbations $Y^\varepsilon\times\lbrace0\rbrace$ of $X_F\times\lbrace0\rbrace$ on $M\times[-\varepsilon_0, \varepsilon_0]$. By \cite[pg 90]{Fenichel1979GeometricSP}, the graph transform argument can still be applied in this case. Hence, $p_\varepsilon$ and $W^{u, loc}_{Y^\varepsilon}(p_\varepsilon)$ will vary in $C^k$ families parametrized by $\varepsilon$. We remark that since all of these manifolds can be seen as $\varepsilon$-parametrized graphs over $p_0$ and a disc in $\mathcal{N}^u$, they are in fact isotopic submanifolds of $M$. We will also use this trick (of augmenting the vector field) in Theorem \ref{thm_fs_global} for the fast-slow systems.

Let $\widetilde{W}^u_{X_F}(p_0)$ be a (compact) $\codim 0$ overflowing submanifold of the global unstable manifold $W^u_{X_F}(p_0)$. In particular, $W^{u, loc}_{X_F}(p_0)\subset\widetilde{W}^u_{X_F}(p_0)$. Hence we can encode every point $y_0\in\widetilde{W}^u_{X_F}(p_0)\setminus W^{u, loc}_{X_F}(p_0)$ as $\phi^T_{X_F}(x_0)$ for some $T>0$ and $x_0\in\partial (W^{u, loc}_{X_F}(p_0))$. Recall that by Remark \ref{rem_smooth_data} the flow depends smoothly on the initial data $(x_0, X_F)$, because $M$ is compact. So the $C^k$ smooth perturbation $Y^\varepsilon$, will perturb $\widetilde{W}^u_{X_F}(p_0)$ also $C^k$ smoothly by a $C^k$ family of $C^k$ maps. Since embeddings are open by Theorem \ref{thm_open_embedd}, we will in fact obtain isotopic overflowing $C^k$ submanifolds of $M$.

To the end, we remark that the smooth family of $W^{u, loc}_{Y^\varepsilon}(p_\varepsilon)$ can also be obtained by Perron's method, see \cite[Apx 1, 4]{Palis1993HyperbolicityAS}.
\end{example}

\begin{defn}\cite{Kuehn2015MultipleTS}A \textbf{Fast-slow system} is a system of ordinary differential equations of the form
\begin{align}\label{eqn_fs_slow}
\begin{split}
\varepsilon \partial_t \theta &=\varepsilon\dot{\theta}=f_1(s, \theta, \varepsilon),\\
\partial_t s &=\,\,\,\dot{s}=f_2(s, \theta, \varepsilon),
\end{split}
\end{align}
where $f_1\in C^k(\R^m\times\R^n\times\R, \R^n)$, $f_2\in C^k(\R^m\times\R^n\times\R, \R^m)$, and $0<\varepsilon\ll 1$. Here, $\theta$ are called the \textbf{fast variables} and $s$ are called the \textbf{slow variables}. We can equivalently put $\tau=t/\varepsilon$ and obtain
\begin{align}\label{eqn_fs_fast}
\begin{split}
\partial_\tau \theta &=\theta^\prime=f_1(s, \theta, \varepsilon),\\
\partial_\tau s &=s^\prime=\varepsilon f_2(s, \theta, \varepsilon).
\end{split}
\end{align}


We will refer to $t$ as the \textbf{slow time} and to $\tau$ as the \textbf{fast time}. The \textbf{fast-slow flow} will be denoted (equivalently) as $\cdot_{f\hbox{-}s\,}t$ or $\fs\tau$.
\end{defn}

\begin{defn}\cite{Kuehn2015MultipleTS}
The \textbf{slow flow} is the solution of the \textbf{slow subsystem}, which is the differential-algebraic equation
\begin{align*}
0&=f_1(s, \theta, 0),\\
\dot{s}&=f_2(s, \theta, 0)
\end{align*}
and will be denoted as $\cdot_{slow}t$.

The \textbf{fast flow} is the solution of the \textbf{fast subsystem}, which is given by
\begin{align*}
\theta^\prime&=f_1(s, \theta, 0),\\
s^\prime&=0
\end{align*}
and will be denoted as $\cdot_{fast}\tau$.
\end{defn}

\begin{defn}\cite{Kuehn2015MultipleTS}
The \textbf{critical set} $C_0$ is defined by
$$C_0:=\lbrace (s, \theta)\in\R^m\times\R^n\,|\,f_1(s, \theta, 0)=0\rbrace.$$
A subset $S_0\subset C_0$ is called \textbf{normally hyperbolic} (for the fast flow) if for each $p\in S_0$ the $(n\times n)$-matrix $D_\theta f_1(p, 0)$ has no eigenvalue with zero real part. We will assume that $S_0$ is a smooth compact $m$-dimensional manifold (possibly with corners).

Points in $C_0$, which do not lie in any normally hyperbolic manifold, are called \textbf{singularities}.
\end{defn}

\begin{rem} Note that the equilibrium points of the fast flow are in bijection with the points of $C_0$.
\end{rem}

\begin{thm}\label{thm_fs_global}\cite{Fenichel1979GeometricSP, Jones1995GeometricSP, Fenichel77} Let $S_0$ be a normally hyperbolic manifold for the fast flow and let $W^{s/u, loc}_0(S_0)$ be its local stable and unstable manifolds. If $\varepsilon_0>0$ is sufficiently small, then there are families $\lbrace S_\varepsilon\rbrace_{\varepsilon\in[0, \varepsilon_0]}$ and $\lbrace W^{s, loc}_\varepsilon(S_\varepsilon)\rbrace_{\varepsilon\in[0, \varepsilon_0]}, \lbrace W^{u, loc}_\varepsilon(S_\varepsilon)\rbrace_{\varepsilon\in[0, \varepsilon_0]}$ such that:
\begin{itemize}
\item[$(i.)$]Families consists of $C^k$ manifolds (with corners) and are $C^{k-1}$ in $\varepsilon$ (the manifolds in each family are isotopic submanifolds of $\R^{m+n}$).
\item[$(ii.)$]$W^{s, loc}_\varepsilon(S_\varepsilon)\cap W^{u, loc}_\varepsilon(S_\varepsilon)=S_\varepsilon$.
\item[$(iii.)$] If $\varepsilon>0$, then $S_\varepsilon, W^{s, loc}_\varepsilon(S_\varepsilon), W^{u, loc}_\varepsilon(S_\varepsilon)$ are locally invariant for the fast-slow flow.
\end{itemize}
\end{thm}

\begin{sproof}
We augment the system $(\ref{eqn_fast_fl})$ by $\varepsilon^\prime=0$ for $|\varepsilon|$ small and would like to apply Theorems \ref{thm_overlfow} and \ref{thm_unstable}. The variable $\varepsilon$ is not too much of a problem due to the behavior of the vector field in the $\varepsilon$ direction. However, slow variables contribute with center-like directions that we need to deal with. For this, we perturb the augmented system by bump functions as discussed in Remark \ref{rem_ext_overflow}.
\end{sproof}



\begin{rem}\label{rem_reg_problem} Due to normal hyperbolicity, we can use the implicit function theorem to write $S_0$ locally as a graph over the slow variable. We will assume that $S_0$ can be globally written as a graph $\lbrace (s, h_0(s))\,|\,s\in K\rbrace$, where $K\subset \R^m$ is a compact, contractible (smooth) manifold with corners.

Then by Theorem \ref{thm_fs_global} the $S_\varepsilon$ can be written as graphs of functions $h_\varepsilon(s)$, which are $C^k$ in $s$ and $C^{k-1}$ in $\varepsilon$.

On $S_\varepsilon$, the fast-slow flow is given as
$$s^\prime=\varepsilon f_2(s, h_\varepsilon(s), \varepsilon).$$
and hence in the slow time as
$$\dot{s}=f_2(s, h_\varepsilon(s), \varepsilon).$$
Observe that as $\varepsilon\rightarrow 0$, the limit of the flow exists and is given by
$$\dot{s}=f_2(s, h_0(s), 0).$$
In particular, this family of flows on $S_\varepsilon$ is $C^{k-1}$ in $\varepsilon$. Hence, the study of the fast-slow flow on $S_\varepsilon$ is a regular perturbation problem in the sense of \cite[Def 5.1.4]{Kuehn2015MultipleTS}.
\end{rem}



\begin{rem} On the first sight, the results from Theorem \ref{thm_fs_global} look a bit weak, since the perturbed manifolds are only locally invariant. But in Remark \ref{rem_reg_problem} we already described the dynamics on $S_\varepsilon$. However, we still do not understand well the dynamics on the perturbations of stable and unstable manifolds. For this, we introduce a certain flow invariant fibration on $W^{s/u, loc}_\varepsilon(S_\varepsilon)$. Which, in particular, will show us that $W^{s/u, loc}_\varepsilon(S_\varepsilon)$ are decaying to $S_\varepsilon$ exponentially and hence, the notation is justified.
\end{rem}


\begin{thm}\label{thm_fen3}\cite{Jones1995GeometricSP, Fenichel1979GeometricSP} Let $S_0$ be a normally hyperbolic manifold for the fast flow. Let $\varepsilon_0>0$ be small. Let $\lbrace S_\varepsilon\rbrace_{\varepsilon\in[0, \varepsilon_0]}$ be parametrized by $v_\varepsilon\in S_\varepsilon$, where the parametrization is $C^k$ in $v$ and $C^{k-1}$ in $\varepsilon$.  Then for each $v_\varepsilon\in S_\varepsilon$ there exists a $C^k$ manifold $W^{s, loc}(v_\varepsilon)$(with corners) such that:
\begin{itemize}
\item[$(i.)$]$W^{s, loc}_\varepsilon(v_\varepsilon)\subseteq W^{s, loc}_\varepsilon(S_\varepsilon)$.
\item[$(ii.)$]$W^{s, loc}_\varepsilon(v_\varepsilon)$ is $C^k$ in $v$ and $C^{k-1}$ in $\varepsilon$.
\item[$(iii.)$]If $\varepsilon=0$, then $W^{s, loc}_0(v_0)$ is a local stable manifold of $v_0\in S_0$ for the fast flow.
\item[$(iv.)$]If $\varepsilon>0$, then $W^{s, loc}_\varepsilon(v_\varepsilon)$ is \textbf{positively invariant}. I.e. if $v_\varepsilon\fs[0, \tau]\in W^{s, loc}_\varepsilon(S_\varepsilon)$, then
$$\big\lbrace w\fs\tau\,|\,w\in W^{s, loc}_\varepsilon(v_\varepsilon), w\fs\tau\subset W^{s, loc}_\varepsilon(S_\varepsilon)\big\rbrace\subset W^{s, loc}_\varepsilon(v_\varepsilon\fs\tau).$$
\item[$(v.)$]For each $\varepsilon>0$ there are constants $C_s>0, L_s<0$ such that if $w\in W^{s, loc}_\varepsilon(v_\varepsilon)$, then
$$d(w\fs\tau-v_\varepsilon\fs\tau)\leq C_s e^{L_s \tau},$$
for all $\tau\geq 0$ for which $w\fs [0, \tau]\subset W^{s, loc}_\varepsilon(S_\varepsilon)$ and $v_\varepsilon\fs [0, \tau]\subset W^{s, loc}_\varepsilon(S_\varepsilon)$.
\end{itemize}

Analogous statement holds for the unstable case, only then in $(iv.)$ and $(v.)$ we consider $\tau\leq0$.
\end{thm}

\begin{sproof} The theorem can be proven similarly to Theorem \ref{thm_fs_global} by a graph transform argument applied to the augmented fast-slow system, which is deformed by bump functions in a neighborhood of $\partial S_0$. However, to obtain a result about the unperturbed system, we need to impose the constraints in $(iv.)$-$(v.)$ to be sure that we do not leave the locally invariant manifolds $W^{s/u, loc}_\varepsilon(S_\varepsilon)$.
\end{sproof}

\begin{rem} Let $D$ be a subset of $S_\varepsilon$. Based on the construction of fibrations from Theorem \ref{thm_fen3}, we introduce
$$W^{s, loc}_\varepsilon(D)=\bigcup_{v\in D}W^{s, loc}_\varepsilon(v)\hbox{ and }W^{s, loc}_\varepsilon(D)=\bigcup_{v\in D}W^{s, loc}_\varepsilon(v).$$
\end{rem}

\begin{defn_thm}\cite{Fenichel1979GeometricSP, Jones2009GeneralizedEL}\label{defn_fen_coord} Let us consider the $C^k$ fast-slow system (\ref{eqn_fs_fast}) with $m$-dimensional normally hyperbolic manifold $S_0\subset \R^{m+n}$. Let $U_0$ be an open (contractible) subset of $S_0$. Then there are \textbf{Fenichel coordinates}. That is a $C^{k-1}$ coordinate change $\Phi^{fen}$ of a neighbourhood of $U_0\times\lbrace0\rbrace\subset \R^{m+n}\times[0, \varepsilon_0)$ such that (\ref{eqn_fs_fast}) takes the \textbf{Fenichel normal form}.
\begin{align}\label{eqn_fenichel_normal_form}
\begin{split}
a^\prime&= A(a, b, x, \varepsilon) a,\\
b^\prime&= B(a, b, x, \varepsilon) b,\\
x^\prime&=\varepsilon\big[h(x, \varepsilon)+ H(a, b, x, \varepsilon)ab\big],
\end{split}
\end{align}
where
\begin{itemize}
\item[(i.)](\ref{eqn_fenichel_normal_form}) is $C^{k-2}$.
\item[(ii.)] $A$ and $B$ are matrices of the form $\ell_u\times\ell_u$ and $\ell_s\times\ell_s$, respectively (recall that $\ell_s+\ell_u=n$). The real parts of the eigenvalues of $A$ are in $(\alpha, \infty)$, and of $B$ are in $(-\infty, -\beta)$, for some $\alpha, \beta>0$. Also $(a, b)\rightarrow H(a, b, x, \varepsilon)ab$ is bilinear.
\item[(iii.)]For each $\varepsilon\in[0, \varepsilon_0)$, it holds in these coordinates that $U_\varepsilon=\lbrace a=0, b=0\rbrace$, $W^{u, loc}_\varepsilon(U_\varepsilon)=\lbrace b=0\rbrace$, and $W^{s, loc}_\varepsilon(U_\varepsilon)=\lbrace a=0\rbrace$. (Here, we were considering $x\in U_{fen}\subset \R^m$, for some $U_{fen}$ diffeomorphic to the original $U_0\subset S_0$.) In particular, along $W^{u, loc}_\varepsilon(U_\varepsilon)$ we have that $(x(\tau), a(\tau), 0)$ is a solution iff $(x(\tau), 0, 0)$ is a solution. Analogously for $W^{s, loc}_\varepsilon(U_\varepsilon)$.
\item[(iv.)]If $W^{u, loc}_\varepsilon(U_\varepsilon)\xrightarrow{\pi_a}U_\varepsilon$ and $W^{s, loc}_\varepsilon(U_\varepsilon)\xrightarrow{\pi_b}U_\varepsilon$ are fibrations induced by the canonical projections $\pi_a, \pi_b$, then the flow of $(\ref{eqn_fenichel_normal_form})$ is a fiber invariant, i.e., the flow maps fibers into fibers. Moreover, let $p\in U_\varepsilon$ and $p_u$ be a lift to the fiber $\pi_a^{-1}(p)$, then the distance of images of $p$ and $p_u$ exponentially decreases under the backward flow. Analogously, for the fibers $\pi_b$, only with the forward flow.
\end{itemize}
This is provided that $\Vert a\Vert \leq\delta, \Vert b\Vert \leq\delta$ and $\varepsilon_0, \delta>0$ are sufficiently small.

By $B_{\delta, U_{fen}}$ we will denote the set $\lbrace(a, b, x)\in\R^{\ell_u}\times\R^{\ell_s}\times\R^m\,|\,\Vert a\Vert \leq\delta, \Vert b\Vert \leq\delta, x\in U_{fen}\rbrace$. We also put $\Phi^{fen}_\varepsilon:=\Phi^{fen}\vert_{\varepsilon\, fixed}$.
\end{defn_thm}

\begin{sproof}First, the coordinate change. Since $U_0$ is contractible, we have a normally hyperbolic splitting of trivial bundles $T\R^{m+n}|_{U_0}=\mathcal{N}^u\oplus\mathcal{N}^s\oplus TU_0$, which gives locally $C^{k-1}$ coordinates $(a, b, x)$ in $\R^{\ell_u}\times\R^{\ell_s}\times\R^m$. By Fenichel's theory (Theorem \ref{thm_fs_global}), we can, for each $\varepsilon\geq0$ sufficiently small, write $W^{s/u, loc}_\varepsilon(U_\varepsilon)$ in these coordinates locally as graphs of some auxiliary functions
$$W^{u, loc}_\varepsilon(U_\varepsilon)=\hbox{Graph}(h^u(a, x, \varepsilon))\hbox{ and }W^{s}(U_\varepsilon)=\hbox{Graph}(h^s(b, x, \varepsilon))$$
such that these functions are $C^{k-1}$ with a domain on $\lbrace x\in U_{fen}\subset R^m, \Vert a\Vert <\delta, \Vert b\Vert <\delta, \Vert \varepsilon\Vert <\varepsilon_0\rbrace$.

Then we change coordinates by straightening $W^{u, loc}_\varepsilon(U_\varepsilon)$. This introduce new variables $(a_1, b_1, x_1)=(a, b-h^u(a, x, \varepsilon), x)$. The Jacobi matrix of the coordinate transformation has a full rank at $a=b=0, \varepsilon=0$, so the transformation is invertible on $\lbrace x\in U_{fen}, \Vert a\Vert <\delta, \Vert b\Vert <\delta, \Vert \varepsilon\Vert <\varepsilon_0\rbrace$ with $\delta, \varepsilon_0>0$ sufficiently small.

Now the argument is similar for the straightening of $W^{s, loc}_\varepsilon(U_\varepsilon)$. Then it remains to straighten the flow invariant fibers of $W^{s, loc}_\varepsilon(U_\varepsilon)$ and $W^{s}(U_\varepsilon)$, which were defined in Theorem \ref{thm_fen3}. This involves changing the $x$-coordinate by the images of the affine maps given by basepoint maps of the fibrations. Again, the rank argument shows that the transformation is invertible on a neighborhood of $a=b=0, \varepsilon=0$.

The Fenichel normal form follows from looking at invariant sets $a=0$ or $b=0$ together with the property that projections along fibers intertwine with the flow.

In particular, the obtained map $H$ is $C^{k-2}$, \cite{Fenichel77}.

The signature of eigenvalues of matrices $A$ and $B$ follows from the normal hyperbolicity of $S_0$ and Taylor expansion estimates; see \cite[Prop 4.1.1]{Kuehn2015MultipleTS}.
\end{sproof}

\begin{defn}\cite{Szmolyan1991TransversalHA}\label{defn_sing_unstabl} Let $\mathfrak{p}_0$ be a hyperbolic equilibrium of the slow flow on a normally hyperbolic manifold for the fast flow. Then we put
$$W^u_{sing}(\mathfrak{p}_0):=\bigcup_{x\in W^u_{slow}(\mathfrak{p_0})}W^u_{fast}(x)$$
and call the set as the \textbf{singular unstable manifold of $\mathfrak{p}_0$}. The restriction of $W^u_{sing}(\mathfrak{p}_0)$ to the local stable manifolds for the fast and slow flows is called the \textbf{local singular unstable manifold of $\mathfrak{p}_0$} and denoted by $W^{u, loc}_{sing}(\mathfrak{p}_0)$. We will also often use the term
$$\widetilde{W}^u_{sing}(\mathfrak{p}_0):=\bigcup_{x\in \widetilde{W}^u_{slow}(\mathfrak{p_0})}W^{u, loc}_{fast}(x),$$
where $\widetilde{W}^u_{slow}(\mathfrak{p_0})$ is a (compact) $\codim 0$ overflowing submanifold of $W^u_{slow}(\mathfrak{p_0})$.

Analogously, we also define the terms in the stable case.
\end{defn}

\begin{cor}\label{cor_fs_local}\cite{Fenichel1979GeometricSP, Szmolyan1991TransversalHA} Let $S_0$ be a normally hyperbolic manifold for the fast flow. Assume also that $p_0\in S_0$ is a normally hyperbolic equilibrium for the slow flow. Then for $\varepsilon>0$ sufficiently small, there are families of normally hyperbolic equilibria $\mathfrak{p}_\varepsilon$ and their local stable and unstable manifolds. Dimensions of the local stable (unstable) manifolds are equal to the sum of stable (unstable) directions in fast and slow flows. The families are $C^{k-1}$ in $\varepsilon$ and describe isotopic $C^k$ manifolds.
\end{cor}

\begin{proof}
Let us consider a Fenichel coordinates around $\mathfrak{p}_0$. By Remark \ref{rem_reg_problem} and the discussion in Example \ref{example_hyperbolicity}, the following holds. The submanifold $W^{u, loc}_{slow}(\mathfrak{p}_0)\subset U_{fen}$ extends to a smooth family of $W^{u, loc}(\mathfrak{p}_\varepsilon)\subset U_{fen}$. Here $W^{u, loc}(\mathfrak{p}_\varepsilon)$ are the local unstable manifold for the (slow time) fast-slow flow restricted to $U_{fen}$.

Hence, we can smoothly extend $W^{u, loc}_{sing}(\mathfrak{p}_0)$ to $W^{u, loc}_{\ffs}(\mathfrak{p}_\varepsilon)$ by taking the unstable fibers over the perturbed $W^{u, loc}_{slow}(\mathfrak{p}_0)$. Note that by the flow invariance of fibers and the exponential decay of solutions in $W^{u, loc}_\varepsilon(U_\varepsilon)$ to $U_\varepsilon$ (Theorem \ref{defn_fen_coord} $(iv.)$) we obtain that in our case $W^{u, loc}_{\ffs}(\mathfrak{p}_\varepsilon)$ are not only locally invariant but actually overflowing.
\end{proof}

\begin{rem}\label{rem_compact_fenichel} By the discussion in Example \ref{example_hyperbolicity} there is also the following extension of Corrollary \ref{cor_fs_local}:  $\widetilde{W}^u_{sing}(\mathfrak{p}_0)\subset S_0$ extends smoothly to isotopic $\widetilde{W}^u_{\ffs}(\mathfrak{p}_\varepsilon)$, where $\widetilde{W}^u_{\ffs}(\mathfrak{p}_\varepsilon)$ which are (compact) $\codim 0$ overflowing submanifolds of $W^u_{\ffs}(\mathfrak{p}_\varepsilon)$.
\end{rem}

\section{Exchange lemma}
\label{s:exch}
\begin{rem}\label{rem_flow_box_coord} Let us consider the set $U_0$ from Definition \ref{defn_fen_coord}. Assume that $U_0$ is rectifiable. Then we can, after a $C^{k-1}$ coordinate change, write the Fenichel normal form of (\ref{eqn_fs_fast}) as
\begin{align}\label{eqn_fenichel_normal_form2}
\begin{split}
a^\prime&= A(a, b, x, \varepsilon) a,\\
b^\prime&= B(a, b, x, \varepsilon) b,\\
x^\prime&=\varepsilon\big[(1, 0, \cdots, 0)^T+ \widetilde{H}(a, b, x, \varepsilon)ab\big].
\end{split}
\end{align}
for $(a, b, x)\in B_{\delta, U_{fen}}\subset \R^{\ell_u}\times\R^{\ell_s}\times\R^m$ and $\varepsilon\in[0, \varepsilon_0]$. Indeed, first we rectify the slow flow along $U_0$. The Flow box theorem, for instance, can accomplish this. And then we make coordinate changes, as in the proof of Theorem \ref{defn_fen_coord}.

By $\phi_{red, \varepsilon}^\tau$ we will denote the \textbf{reduced flow} of $(\ref{eqn_fenichel_normal_form2})$ to $U_{fen}$ (parametrized with the fast time).
\end{rem}

\begin{rem}\label{rem_monot_estim} By Theorem \ref{defn_fen_coord} we know how to track the trajectories $u_\varepsilon(\tau)=(x(\tau), a(\tau), b(\tau))$ of the system $(\ref{eqn_fenichel_normal_form2})$ that start in the stable or the unstable manifold, i.e. in the set $\lbrace a=0\rbrace$ or $\lbrace b=0\rbrace$.

In this section, we will be interested in tracking $u_\varepsilon$ when $a(0)$ and $b(0)$ are nonzero. First observation will be that the norm of $a$ is exponentially decreasing in the forward time and the norm of $b$ is exponentially decreasing in the backward time. This follows from eigenvalue estimates from Theorem \ref{defn_fen_coord} $(ii.)$, for more details see \cite[cor 1]{Jones2009GeneralizedEL}. We also remark that these are only $C^0$ estimates for $u_\varepsilon$.

In order to have more information about the track of $u_\varepsilon$ or even some flow invariant manifolds, we will present so-called Exchange lemmata which are based on the technique of \textbf{Silnikov boundary value problems} \cite{shilnikov1967existence}. That is a quest to find solutions of $(\ref{eqn_fenichel_normal_form2})$ on a time interval $[0, T]$ with boundary value problems
\begin{equation}\label{egn_first_silnikovBVP}
x(0)=x^1, a(0)=a^0, b(T)=b^1
\end{equation}
or
\begin{equation}\label{egn_second_silnikovBVP}
x(T)=x^1, a(0)=a^0, b(T)=b^1.
\end{equation}
We can heuristically see the existence of the solution to the Silnikov BVP as in Figure \ref{figure_Silnikov_Bvp}, see also \cite{Hsu2015AGS}. The Exchange lemmata, derived from the solutions of two Silnikov BVPs $(\ref{egn_first_silnikovBVP})$, $(\ref{egn_second_silnikovBVP})$, give us two different perspectives on how the tracked invariant manifolds can be parametrized. To the Silnikov BVP $(\ref{egn_second_silnikovBVP})$ will correspond Exchange lemmata \ref{lemma_exchangeA} and \ref{lemma_exchangeB} and to BVP $(\ref{egn_first_silnikovBVP})$ corresponds Exchange Lemma \ref{lemma_exchangeC}.
\begin{figure}[!htbp]
\labellist
\pinlabel $0$ at -8 -8
\pinlabel $a^1$ at 200 -8
\pinlabel $b^0$ at -8 200
\endlabellist
\centering
\includegraphics[scale=0.58]{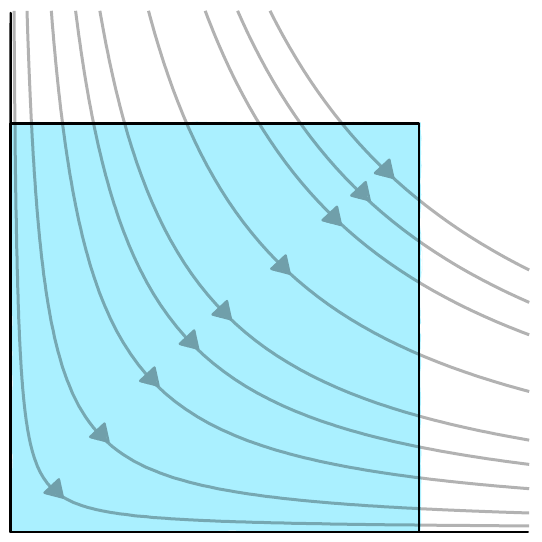}
\vspace{0.3cm}
\caption{Trajectories of the system $a^\prime=a,b^\prime=-b, x^\prime=0$ in the rectangle $\lbrace (a, b)\in[0, a^1]\times[0, b^0]\rbrace$ can be parametrized by $T\geq0$, where $b(0)=b^0$ and $a(T)=a^1$. A solution can stay inside the rectangle for an arbitrarily long time $T$ if it is sufficiently close to the axes.}
\label{figure_Silnikov_Bvp}
\end{figure}
\end{rem}

\begin{rem_not}\label{rem_flow_inv_mfld}By $\bm{M}$ we will denote a $C^{k-2}$ submanifold of $B_{\delta, U_{fen}}\times[0, \varepsilon_0]$, such that $\varepsilon$-level sets are isotopic $C^{k-2}$ submanifolds. We will denote the $\varepsilon$-levels by $M_\varepsilon$. $M^\ast_\varepsilon$ will characterize the evolution of $M_\varepsilon$ under the flow of $(\ref{eqn_fenichel_normal_form2})$.
\end{rem_not}

The following is a simplified version of Schecter's formulation of the Exchange Lemma.

\begin{ExchangeA}\cite{Schecter2008ExchangeL2}\label{lemma_exchangeA} Let us consider the $C^{k-2}$ system $(\ref{eqn_fenichel_normal_form2})$ and the $(\ell_{u}+1)$-dimensional manifold $\bm{M}$ (from Remark \ref{rem_flow_inv_mfld}). Assume that inside $B_{\delta, U_{fen}}$ the manifold $M_0$ intersects $\lbrace a=0\rbrace$ transversally at some point $(b^0, 0, 0)$. Let $\overline{x}>0$ such that
$$u_0(t)=\lbrace (t, 0,\dots, 0)^T\,|\,t\in[0, \overline{x}]\rbrace\subset U_{fen}$$
($u_0$ is parametrized by the slow time) and $D$ be a small neighborhood of $(0, \overline{x})$ in $ax_1$-space. Then for $\varepsilon_0>0$ sufficiently small there exist $C^{k-3}$ functions
\begin{align*}
\widetilde{b}&:D\times [0, \varepsilon_0]\rightarrow\R^{\ell_s},\\
\widetilde{x}&:D\times [0, \varepsilon_0]\rightarrow\R^{m-1},
\end{align*}
such that
\begin{itemize}
\item[$(i.)$]$\widetilde{b}(a, x_1, 0)=0$.
\item[$(ii.)$]$\widetilde{x}(a, x_1, 0)=\widetilde{x}(0, x_1, \varepsilon)=0$.
\item[$(iii.)$] As $\varepsilon\rightarrow 0$, the functions $(\widetilde{b}, \widetilde{x})\rightarrow 0$ exponentially (in uniform norm), along all derivatives up to order $k-3$ with respect to all variables.
\item[$(iv.)$] For $\varepsilon\in(0, \varepsilon_0]$, the set $$Q_\varepsilon=\lbrace (a, b, x)\,|\,(a, x_1)\in D, b=\widetilde{b}(a, x_1, \varepsilon), (x_2,\dots, x_m)=\widetilde{x}(a, x_1, \varepsilon)\rbrace$$
is contained in $M^\ast_\varepsilon$. See also Figure \ref{figure_Exchange}.
\end{itemize}

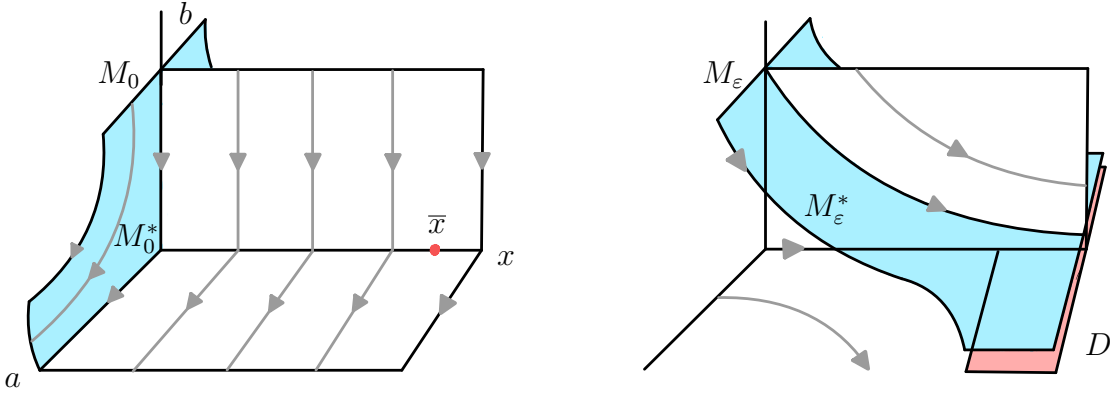
\begin{figure}[!htbp]
\labellist
\pinlabel $\overline{x}$ at 172 749
\pinlabel $x$ at 200 735
\pinlabel $a$ at 0 685
\pinlabel $b$ at 70 835
\pinlabel $M_0$ at 43 810
\pinlabel $M_0^\ast$ at 50 743
\pinlabel $M_\varepsilon$ at 288 810
\pinlabel $M_\varepsilon^\ast$ at 330 755
\pinlabel $D$ at 440 700
\endlabellist
\centering
\includegraphics[scale=0.93]{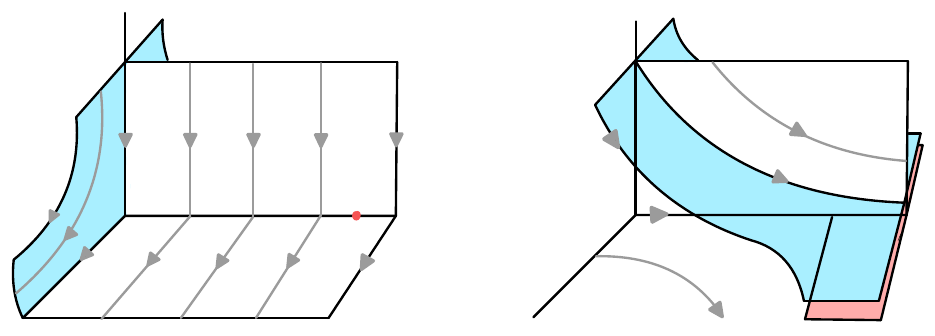}
\vspace{0.3cm}
\caption{$(\ell_u+1)$-Exchange Lemma for $m=\ell_s=\ell_u=1$. \textit{On the left:} evolution of $M_0$ under fast flow. \textit{On the right:} $C^{k-3}$-inclination of $M_\varepsilon$ under the fast-slow flow.}
\label{figure_Exchange}
\end{figure}
\end{ExchangeA}

\begin{rem_not}\label{rem_flow_inv_mfld2}\cite{Schecter2008ExchangeL2}Let us consider $(\ref{eqn_fenichel_normal_form2})$ and the $(\ell_{u}+\sigma+1)$-dimensional manifold $\bm{M}$ (from Remark \ref{rem_flow_inv_mfld}), where $0\leq\sigma\leq m-1$.

Let us assume that
\begin{itemize}
\item[$(i.)$] Inside $B_{\delta, U_{fen}}$ the manifold $M_0$ intersects $\lbrace a=0\rbrace$ transversally in a manifold $N_0$ (of dimension $\sigma$).
\item[$(ii.)$] $N_0$ projects to a $\sigma$-dimensional submanifold $P_0$ of $U_{fen}$.
\item[$(iii.)$] The vector $(1, 0,\dots, 0)^T$ is not tangent to $P_0$.
\end{itemize}
Note that for $\varepsilon_0>0$ small we can (uniquely) define $N_\varepsilon, P_\varepsilon$ that satisfy $(i.)-(iii.)$ with indices $0$ replaced by $\varepsilon\in[0, \varepsilon_0]$.

Hence, similarly to Theorem \ref{defn_fen_coord} and Remark \ref{rem_flow_box_coord} we can introduce ($\varepsilon$ dependent) $C^{k-1}$ coordinates $(a, b, u, v, w)\in\R^{\ell_u}\times\R^{\ell_s}\times\R\times\R^{\sigma}\times\R^{m-1-\sigma}$, which locally take $P_\varepsilon$ to the $v$-space and such that the system (\ref{eqn_fenichel_normal_form2}) takes a form

\begin{align}\label{eqn_fenichel_normal_form3}
\begin{split}
a^\prime&= A(a, b, u, v, w, \varepsilon) a,\\
b^\prime&= B(a, b, u, v, w, \varepsilon) b,\\
u^\prime&=\varepsilon\big[1+ e(a, b, u, v, w, \varepsilon)ab\big],\\
v^\prime&=\varepsilon F(a, b, u, v, w, \varepsilon)ab,\\
w^\prime&=\varepsilon G(a, b, u, v, w, \varepsilon)ab,\\
\end{split}
\end{align}

for $(a, b, u, v, w, \varepsilon)\in B_{\delta, U_{fen}}$ and $\varepsilon\in[0, \varepsilon_0]$.
\end{rem_not}

The following is another simplified version of Schecter's formulation of the Exchange Lemma.

\begin{ExchangeB}\cite{Schecter2008ExchangeL2}\label{lemma_exchangeB} Let us consider the $C^{k-1}$ system $(\ref{eqn_fenichel_normal_form3})$ and the $(\ell_{u}+\sigma+1)$-dimensional manifold $\bm{M}$ (from Remark \ref{rem_flow_inv_mfld2}). Let $\overline{u}>0$ such that
$$\lbrace (t, 0,\dots, 0)\,|\,t\in[0, \overline{u}]\rbrace\subset U_{fen}$$
and $D$ be a small neighborhood of $(0, \overline{u}, 0)$ in $auv$-space. Then for $\varepsilon_0>0$ sufficiently small there exist $C^{k-3}$ functions
\begin{align*}
\widetilde{b}&:D\times [0, \varepsilon_0]\rightarrow\R^{\ell_s},\\
\widetilde{w}&:D\times [0, \varepsilon_0]\rightarrow\R^{m-1-\sigma},
\end{align*}
such that
\begin{itemize}
\item[$(i.)$]$\widetilde{b}(a, u, v, 0)=0$.
\item[$(ii.)$]$\widetilde{w}(a, u, v, 0)=\widetilde{w}(0, u, v, \varepsilon)=0$.
\item[$(iii.)$] As $\varepsilon\rightarrow 0$ the functions $(\widetilde{b}, \widetilde{w})\rightarrow 0$ exponentially, along all derivatives up to order $k-3$ with respect to all variables.
\item[$(iv.)$] For $\varepsilon\in(0, \varepsilon_0]$, the set
$$Q_\varepsilon=\left\lbrace (a, b, u, v, w)\,|\,(a, u, v)\in D, b=\widetilde{b}(a, u, v, \varepsilon), w=\widetilde{w}(a, u, v, \varepsilon)\right\rbrace$$
is contained in $M^\ast_\varepsilon$.
\end{itemize}
\end{ExchangeB}

The following is Brunovsk{\'y}'s formulation of the Exchange Lemma.

\begin{ExchangeB}\cite{Brunovsk1999CrInclinationTF}\label{lemma_exchangeC} Let us consider the $C^{k-1}$ system $(\ref{eqn_fenichel_normal_form3})$ and the $(\ell_{u}+\sigma+1)$-dimensional manifold $\bm{M}$ (from Remark \ref{rem_flow_inv_mfld2}). 

Let $\overline{t}>0$ and let $V_p$ be a small neighborhood of $0$ in the $v$-space such that
$$\phi^t_{slow}(0, V_P, 0)\subset U_{fen}\hbox{ for all }t\in[0, \overline{t}].$$
Let $\underline{t}\in(0, \overline{t})$ and let $\delta_0>0$ be small. Then for $\varepsilon>0$ sufficiently small there is a $C^{k-3}$ function
$$\widetilde{s}^\varepsilon:\lbrace\Vert a\Vert\leq\delta_0\rbrace\times\phi^{[\underline{t}/\varepsilon, \overline{t}/\varepsilon]}_{red, \varepsilon}(0, V_p, 0)\rightarrow\R^{\ell_s+m}$$
such that 
$$(a, b, u, v, w)=\phi^\tau_{\ffs}(a^0, b^0, 0, v^0, 0)$$
for some $(a^0, b^0, 0, v^0, 0)\in M_\varepsilon$ and $\tau\in[\underline{t}/\varepsilon, \overline{t}/\varepsilon]$ if and only if
\begin{equation}\label{eqn_s_funct_ext}
(b, u, v, w)=\widetilde{s}^\varepsilon(a, x),
\end{equation}
where $x\in\phi^\tau_{red, \varepsilon}(0, v^0, 0)$.

Moreover, as $\varepsilon\rightarrow 0$ the functions $\widetilde{s}^{\varepsilon}(a, x)\rightarrow (0, x)$ exponentially, along all derivatives up to order $k-3$ with respect to all variables.

To the end, we recall that by our choice of coordinates we have the simplifications 
$$\phi^{[\underline{t}/\varepsilon, \overline{t}/\varepsilon]}_{red, \varepsilon}(0, V_p, 0)=([\underline{t}/\varepsilon, \overline{t}/\varepsilon], V_p, 0)$$
and 
$$\phi^\tau_{red, \varepsilon}(0, v^0, 0)=(\varepsilon\tau, v^0, 0).$$
\end{ExchangeB}

\begin{rem}\cite[Rem 4.2, Rem 4.5]{Brunovsk1999CrInclinationTF}\label{rem_inclination1} 
By $(\ref{eqn_s_funct_ext})$ the set $\lbrace a\rbrace\times Graph(\widetilde{s}^\varepsilon)$ describe only points of those trajectories through $M_\varepsilon$ over $V_P$ $(=:M_\varepsilon\vert_{V_P})$ that do not leave the set $\lbrace\Vert a\Vert\leq\delta_0\rbrace$. 

Let $\widehat{u}_\varepsilon:[0, \tau]\rightarrow\R^{m+n}$ be such a (fast time) trajectory, then in addition 
\begin{itemize}
\item $\pi_{x}(\widehat{u}_\varepsilon)$ and $\phi^{[0, \tau]}_{red, \varepsilon}\big(\pi_x(\widehat{u}(0))\big)$ are $O(e^{-1/\varepsilon})$ close (recall also that along $U_{fen}$ the $1/\varepsilon$-multiple of the reduced vector field converges to the slow vector field as $\varepsilon\rightarrow 0$).
\item $\Vert a^0\Vert=O(e^{-1/\varepsilon})$.
\end{itemize}
Hence, only the points of $M_\varepsilon\vert_{V_P}$ with exponentially small $\Vert a^0\Vert$ matter. 
\end{rem}

\begin{rem}\label{rem_combination_exchange}From Brunovsk{\'y}'s Exchange Lemma \ref{lemma_exchangeC} we know that the sets $\lbrace a\rbrace\times Graph(\widetilde{s}^\varepsilon)$ exponentially converge to the $auv$-space as $\varepsilon\rightarrow0$. So for $\varepsilon>0$ small the \textit{canonical} orthogonal projection of $\lbrace a\rbrace\times Graph(\widetilde{s}^\varepsilon)$ onto $auv$-space is injective.

Let us pick some $\delta_1>0$ small and $V_p$ as in Brunovsk{\'y}'s Exchange Lemma \ref{lemma_exchangeC}. Hence the set $Graph(\widetilde{b}^\varepsilon, \widetilde{w}^\varepsilon)$ from Schecter's Exchange Lemma \ref{lemma_exchangeB} can be seen as the set of those points of 
$$\phi^{\overline{u}+\delta_1}_{\ffs}(M_\varepsilon\vert_{V_P})\cap\lbrace\Vert a\Vert\leq\delta_0\rbrace$$
that projects orthogonally to $D$ provided that $D$ and $\varepsilon>0$ are small (here $\overline{u}+\delta_1$ is a slow time).
\end{rem}



\begin{rem} We stated the Exchange lemmata in the form of ``inclination theorems.'' In the literature, another popular version of the Exchange lemma is to track manifolds close to ``jumps'' of trajectories combining the fast and slow flows; see, for example, \cite{Schecter2008ExchangeL2, Jones2009GeneralizedEL, Jones1995GeometricSP}.
\end{rem}

\section{Application to Morse flow trees}\label{sec_appl_chord}
\label{s:appl_trees}
In this section, we relate the flows of $-\nabla E_0$ and $-\nabla E_\varepsilon$ with the multiple time scale dynamics and show the desired correspondence between the Morse flow trees. Since the proofs of the correspondence are a bit technical, we kindly recommend that the reader looks at Remark \ref{rem_short_proof_corresp}, where the key steps and ideas of the proofs are outlined. At the end of the section, we briefly discuss some orientation conventions.

\begin{example}\label{ex_fast_slow} The negative gradient vector fields $-\nabla E_\varepsilon$ determine a fast-slow system
\begin{align}\label{eqn_fast_slow_fl}
\begin{split}
\varepsilon\dot{\theta_1}&=\textcolor{teal}{-\langle P, v_1^{\bot}\rangle}+\varepsilon\big[\tau(s_1)d_1^{-1}\langle P+\varepsilon v_2, \dot{\gamma}(s_1)\rangle-\langle v_2, v_1^\bot\rangle\big],\\
\varepsilon\dot{\theta_2}&=\textcolor{teal}{\langle P, v_2^{\bot}\rangle}-\varepsilon\big[\tau(s_2)d_2^{-1}\langle P-\varepsilon v_1, \dot{\gamma}(s_2)\rangle+\langle v_1, v_2^\bot\rangle\big]\\
\dot{s}_1&=-d_1^{-1}\textcolor{blue}{\langle P, \dot{\gamma}(s_1)\rangle}-\varepsilon \big[d_1^{-1}\langle v_2, \dot{\gamma}(s_1)\rangle\big],\\
\dot{s}_2&=d_2^{-1}\textcolor{blue}{\langle P, \dot{\gamma}(s_2)\rangle}+\varepsilon\big[d_2^{-1}\langle v_1, \dot{\gamma}(s_2)\rangle\big],
\end{split}
\end{align}
where $(s_1, \theta_1, s_2, \theta_2)\in (\R/T\mathbb{Z}\times S^1)^2$. Recall that $d_i=1-\varepsilon \cos(\theta_i) \kappa(s_i)$, see also Notation \ref{notat_torus_first}.

Then the slow subsystem is given by
\begin{align}\label{eqn_slow_fl}
\begin{split}
0&=\textcolor{teal}{-\langle P, v_1^{\bot}\rangle},\\
0&=\textcolor{teal}{\langle P, v_2^{\bot}\rangle},\\
\dot{s}_1&=\textcolor{blue}{-\langle P, \dot{\gamma}(s_1)\rangle},\\
\dot{s}_2&=\textcolor{blue}{\langle P, \dot{\gamma}(s_2)\rangle},
\end{split}
\end{align}
and the fast subsystem is given by
\begin{align}\label{eqn_fast_fl}
\begin{split}
\theta_1^\prime&=\textcolor{teal}{-\langle P, v_1^{\bot}\rangle},\\
\theta_2^\prime&=\textcolor{teal}{\langle P, v_2^{\bot}\rangle},\\
s_1^\prime&=0,\\
s_2^\prime&=0.
\end{split}
\end{align}

Note that the critical set $C_0$ is given by the set $\lbrace \overline{\overline{G}}_1=0\wedge \overline{\overline{G}}_2=0\rbrace$ which was described in Remark \ref{rem_sheets_of_G}. Also,
$$\partial_{\theta_i}\langle P, v_i^\bot\rangle=-\langle P, v_i\rangle.$$
Hence, outside of special and diagonal points, $C_0$ has a structure of \textcolor{teal}{normally hyperbolic} $2$-dimenisonal trivial fiber bundle over $(s_1, s_2)$. The typical fiber is given by $4$ points. Over a special point, it holds that $C_0$ is diffeomorphic to the disjoint union of two circles. At the diagonal, $C_0$ is diffeomorphic to $S^1\times S^1$.

By $S^{out-out}_0$ we will denote a subset of $C_0$ consisting of normally hyperbolic points such that $\langle P, v_1\rangle>0$ and $\langle P, v_2\rangle>0$. We moreover impose the condition that $\pi_{s_1, s_2}(S^{out-out}_0)=\overline{S_K}$, where the standard set $S_K$ was introduced in Remark \ref{rem_delta_standard}. In particular, $S^{out-out}_0$ is diffeomorphic to $\overline{S_K}$, and hence it is a compact manifold with corners.

If $x\in S^{out-out}$, then the linearized fast subsystem has at $x$ a positive eigenvalue $\langle P, v_1\rangle$ with eigenvector $\partial_{\theta_1}\vert_x$ and a negative eigenvalue $-\langle P, v_1\rangle$ with eigenvector $\partial_{\theta_2}\vert_x$.
\end{example}

\begin{rem} The special points impose singularities, so-called turning points. For the study of the dynamics around turning points, see, for example, the Exchange lemma in \cite{Schecter2008ExchangeL2}. To the author, it is not known whether anywhere in literature, there exists a study of Multiple-time scale dynamics for singularities of the type of diagonal points, i.e., when all eigenvalues vanish.

For our purposes, we will be, after careful analysis, able to show that the desired Morse flow trees will generically avoid these singularities.
\end{rem}

\begin{lemma}\label{lemma_identif} The projection $\pi_{s_1, s_2}$, when restricted to $S^{out-out}_0$, is a diffeomorphism onto its image. $\pi_{s_1, s_2}$ canonically identifies the slow flow (\ref{eqn_slow_fl}) on $S^{out-out}_0$ with the $-\nabla E_0$ flow outside of weakly diagonal and weakly special points.
\end{lemma}

\begin{proof} As in Remark \ref{rem_reg_problem}, we can, due to the normal hyperbolicity, describe $S^{out-out}_0$ as a graph over the slow variables. More explicitly, the smooth lift to fast variables is uniquely given as in Lemma \ref{lemma_standard_square} by $$v_i=D_i^{-1/2}\langle P, n(s_i)\rangle n(s_i)+\langle P, b(s_i)\rangle b(s_i),$$ where $D_i=\langle P, n(s_i)\rangle^2+\langle P, b(s_i)\rangle^2.$ Then from Example \ref{ex_fast_slow} we know that the differential equations for the slow flow $(\ref{eqn_slow_fl})$ are fortunately only in slow variables. In addition, they coincide with the differential equations for the flow $\phi_{E_0}$.
\end{proof}

\begin{notat}The lift from Lemma \ref{lemma_identif} will be denoted by $\ioo$.
\end{notat}

\begin{defn}We say that the gradient-like vector fields $X_{E_0}$ and $X_{E_\varepsilon}$ are \textbf{$K$-approximations} of $-\nabla E_0$ and $-\nabla E_\varepsilon$ if there is a vector field $Y$ on $(\R/T\mathbb{Z})^2$ such that
$$X_{E_0}=-\nabla E_0+Y\hbox{ and }X_{E_\varepsilon}=-\nabla E_\varepsilon+Y_0,$$
where $Y_0$ is given by equations $(\pi_{s_1, s_2})_\ast Y_0=Y$ and $(\pi_{\theta_1, \theta_2})_\ast Y_0=0$.
\end{defn}

\begin{rem}Note that unlike the approximations of $-\nabla E_\varepsilon$ in Chapter \ref{ch:adiab_conley} now the $K$-approximations are perturbing $-\nabla E_\varepsilon$ only in $s_1$ and $s_2$ variables.
\end{rem}

\begin{rem}\label{rem_K_approx} For $K$-approximations of $-\nabla E_0$ and $-\nabla E_\varepsilon$ it holds that $S^{out-out}_0$ still consists of normally hyperbolic points and the analogy to Lemma \ref{lemma_identif} is still true. Indeed, the algebraic equations from the slow subsystem are unchanged, and the differential equations from the slow subsystem are still dependent only on variables $s_1$ and $s_2$.
\end{rem}

\begin{rem}Let $X_{E_0}$ and $X_{E_\varepsilon}$ be $K$-approximations of $-\nabla E_0$ and $-\nabla E_\varepsilon$. Let $p_0\in Cr(E_0)\setminus \Delta_0$ and $p_\varepsilon$ be the unique corresponding elements from in $Cr(E_\varepsilon)\cap M_{K, \varepsilon}$. Recall that the correspondence is given by Lemma \ref{lemma_morse_corresp}, see also equations $(\ref{eqn_of_crit_pt}).$

By $\mathfrak{p}_0$ and $\mathfrak{p}_\varepsilon$ we denote the unique realizations of $p_0$ and $p_\varepsilon$ as the equilibria of the slow flow on $S^{out-out}_0$ and fast-slow flow, respectively. Then $\mathfrak{p}_0=\mathfrak{p}_\varepsilon$.
\end{rem}

\begin{rem}\label{rem_short_proof_corresp} In danger of future repetitiveness, let us make a \textit{road map of further theorems and the key steps of their proofs.} In Theorem \ref{thm_corresp_traj_gsp}, Lemmata  \ref{lem_tree_2} and \ref{lem_corres_tree3} we are going to show a bijective correspondence $\Psi^\varepsilon$ between the trees on $(\R/T\mathbb{Z})^2$ and $(\R/T\mathbb{Z}\times S^1)^2$ with $m=0, 1, 2$ interior vertices, provided that $\varepsilon>0$ is small. Then we claim that the cases $m=0, 1, 2$ capture all phenomena that will be necessary for the proof of any case with $m>2$. Recall also that so far we defined the trees only for $Crit$ with low indices and outside of the diagonals, see Chapter \ref{ch:file8}. So, let us briefly discuss the proofs.

$m=0:$

By generic a choice of $X_{E_0}$ and index reasons we will be able to lift all trees (orbits) on $(\R/T\mathbb{Z})^2$ to corresponding orbits of the slow flow on the normally hyperbolic critical submanifold $S^{out-out}\subset (\R/T\mathbb{Z}\times S^1)^2$. These orbits of the slow flow can be realized as transverse intersections of singular stable and unstable manifolds. With the help of Fenichel theory and the Stability Lemma \ref{lem_stability}, we perturb compact parts of these singular manifolds. This will construct desired fast-slow orbits and, in particular, give us the definition and injectivity of $\Psi^\varepsilon$. Then, we continue with a slightly laborious discussion of the local uniqueness of these fast-slow orbits. Such a discussion will be beneficial also in the cases $m=1, 2$. Finally, the surjectivity will be similar to the Theorem \ref{thm_almost_bijection}.

$m=1:$

Now, the single interior vertex is present. Since the bifurcations depend also on the embeddings of $K$ and $T_{K, \varepsilon}$ into $\R^3$, it will not be enough to perturb only the singular stable and unstable manifolds. Hence, we will also be perturbing the evaluation maps from Definition \ref{defn_evE}. For the surjectivity, we use Lemma \ref{lem_preim_treeT} to show that the projections of the bifurcations of trees on the torus converge to the bifurcations of the trees on the knot.

$m=2:$

Now, an interior edge appears. Such an edge corresponds to the space of finite-length trajectories. This space is neither a stable, nor an unstable manifold; hence, the Fenichel theory itself will not be sufficient. So we will represent the tree under the backward flow (see Remark \ref{rem_backward_tree}) and overcome the interior edge with the Exchange Lemma.
\end{rem}

\begin{thm}\label{thm_corresp_traj_gsp} Let $p_0\in Cr_1(E_0)\setminus\Delta_0, q_0\in Cr_0(E_0)\setminus\Delta_0$ and $p_\varepsilon, q_\varepsilon$ are their unique corresponding critical points of $E_\varepsilon$ in $M_{K, \varepsilon}$. If $\varepsilon>0$ is sufficiently small and $X_{E_0}$ and $X_{E_\varepsilon}$ are generic $K$-approximations of $-\nabla E_0$ and $-\nabla E_\varepsilon$, then there is a map
$$\Psi^\varepsilon_{p_0, q_0}:\mathcal{M}_{X_{E_0}}(p_0, q_0)\rightarrow \mathcal{M}_{X_{E_\varepsilon}}(p_\varepsilon, q_\varepsilon)$$
which is a bijection. In addition,
$$\mathcal{M}^{out-out}_{X_{E_\varepsilon}}(p_\varepsilon, q_\varepsilon)=\mathcal{M}_{X_{E_\varepsilon}}(p_\varepsilon, q_\varepsilon).$$
\end{thm}

\begin{proof}
Let $X_{E_0}$ be generically approximating $-\nabla E_0$ in the sense of Corollary \ref{cor_perturb_grad_like}, i.e. $X_{E_0}$ is Morse-Bott-Smale and $W^u_{X_{0}}(p_0)$ is transverse to all special points.

\bigskip

\textit{Definition of $\Psi^\varepsilon_{p_0, q_0}$:}

Let $u \in \mathcal{M}_{X_{E_0}}(p_0, q_0)$. By the choice of $X_{E_0}$ it holds that $u\in S_K$, which was defined by $\delta_K>0$ sufficiently small. Hence, by Remark \ref{rem_K_approx}, there is a small contractible open neighborhood $U$ of $\lbrace p_0\cup u\cup q_0\rbrace\subset K\times K$ which can be canonically lifted as a subset $U_0\subset S^{out-out}_0$. In particular, we lifted $u$ to $u_0$ which is an (transverse) intersection of $W^u_{slow}(\mathfrak{p}_0)$ and $W^s_{slow}(\mathfrak{q_0})$ in $U_0$.

Now, we are going to construct certain $\widetilde{W}^u_{slow}(\mathfrak{p}_0), \widetilde{W}^s_{slow}(\mathfrak{q}_0)\subset U_0$, i.e. by Definition \ref{defn_sing_unstabl} compact overflowing/inflowing $\codim 0$ submanifolds of $W^u_{slow}(\mathfrak{p}_0)$ and $W^s_{slow}(\mathfrak{q}_0)$, respectively.

We pick an auxiliary slow-time parametrization of $u_0$. Then $\widetilde{W}^u_{slow}(\mathfrak{p}_0)$ consists of a small local unstable manifold of $\mathfrak{p}_0$ together with $u_0((-\infty, 1])$. Next, $\widetilde{W}^s_{slow}(\mathfrak{q}_0)$ consists of a small local stable manifold of $\mathfrak{q}_0$ together with a small neighborhood of $u_0([-1, \infty))$. Finally, we slightly cut the obtained set (as in Example \ref{example_smooth_index_pair}), such that $\widetilde{W}^s_{slow}(\mathfrak{q}_0)$ is really inflowing invariant. Note that $\widetilde{W}^u_{slow}(\mathfrak{p}_0)$ is overflowing automatically. 

Now we extend $\widetilde{W}^u_{slow}(\mathfrak{p}_0)$ and $\widetilde{W}^s_{slow}(\mathfrak{q}_0)$ to $\widetilde{W}^u_{sing}(\mathfrak{p}_0)$ and $\widetilde{W}^s_{sing}(\mathfrak{q}_0)$ such that $\widetilde{W}^u_{sing}(\mathfrak{p}_0), \widetilde{W}^s_{sing}(\mathfrak{q}_0)\subset (\Phi^{fen}_0)^{-1}(B_{\delta, U_{fen}})$ for some Fenichel coordinates around $U_0$ with small $\delta>0$.

Last, we consider an auxiliary smooth compact $4$-disk $D$ in a neighborhood of $u_0(0)\in (\R/T\mathbb{Z}\times S^1)^2$. We can choose the disk $D$ small enough such that
\begin{itemize}
\item $u_0\pitchfork\partial D$ in exactly two points.
\item $D\cap\big(\partial\widetilde{W}^u_{sing}(\mathfrak{p}_0)\cup\partial\widetilde{W}^s_{sing}(\mathfrak{q}_0)\big)=\emptyset.$
\item $\widetilde{W}^u_{sing}(\mathfrak{p}_0)\pitchfork D$ and $\widetilde{W}^u_{sing}(\mathfrak{p}_0)\pitchfork D$. 
\item $\pi_{s_1, s_2}(D)\cap \widehat{u}=\emptyset$, where $\widehat{u}\in \mathcal{M}_{X_{E_0}}(p_0, q_0)$ such that $u\neq \widehat{u}$ (we recall that $\mathcal{M}_{X_{E_0}}(p_0, q_0)$ is a finite set).
\end{itemize}
See also Figure \ref{figure_fast_slow_intersect}.
\begin{figure}[!htbp]
\labellist
\pinlabel $\widetilde{W}^u_{slow}(\mathfrak{p}_0)$ at 60 90
\pinlabel $\widetilde{W}^s_{slow}(\mathfrak{q}_0)$ at 424 104
\pinlabel $D_0$ at 227 80
\pinlabel $\mathfrak{p}_0$ at 30 40
\pinlabel $\mathfrak{q}_0$ at 355 42
\endlabellist
\centering
\includegraphics[scale=0.70]{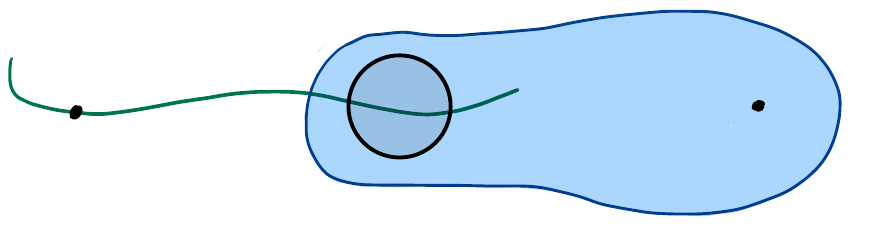}
\vspace{0.3cm}
\caption{The intersection of $\widetilde{W}^u_{slow}(\mathfrak{p}_0)$ and $\widetilde{W}^s_{slow}(\mathfrak{q}_0)$ with $D_0:=D\cap U_0$.}
\label{figure_fast_slow_intersect}
\end{figure}
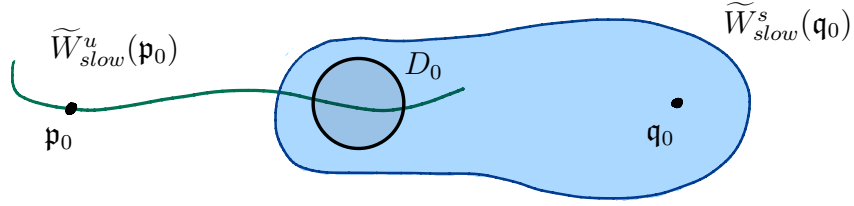

Then, if $\varepsilon_0>0$ is sufficiently small, the Fenichel theory (Corollary \ref{cor_fs_local} and Remark \ref{rem_compact_fenichel}) gives us a $C^k$-isotopy:
$$H:\big(D\times\widetilde{W}^u_{sing}(\mathfrak{p}_0)\times \widetilde{W}^s_{sing}(\mathfrak{q}_0)\big)\times[0, \varepsilon_0]\longrightarrow (\R/T\mathbb{Z}\times S^1)^6$$
from the canonical inclusion to
$$h^{[\varepsilon_0]}\big(D, \widetilde{W}^u_{sing}(\mathfrak{p}_0), \widetilde{W}^s_{sing}(\mathfrak{q}_0)\big)=\big(D, \widetilde{W}^u_{\ffs}(\mathfrak{p}_{\varepsilon_0}), \widetilde{W}^s_{\ffs}(\mathfrak{q}_{\varepsilon_0})\big).$$

On the product of $3$ tori $(\R/T\mathbb{Z}\times S^1)^6$ we consider coordinates $$((r_1, \rho_1, r_2, \rho_2),(s_1, \theta_1, s_2, \theta_2), (y_1, \alpha_1, y_2, \alpha_2)).$$ Then we consider the diagonal $\Delta\subset(\R/T\mathbb{Z}\times S^1)^6$ which is given by $8$ equations
$$\lbrace r_i=s_i=y_i\wedge\rho_i=\theta_i=\alpha_i\rbrace.$$

We stress that $u_0\subset S^{out-out}_0$, which is a normally hyperbolic set. So, by the construction, it holds that $h^{[0]}$ is stratum transverse to the diagonal $\Delta$ and $c_0=\pi_{r_1, \rho_1, r_2, \rho_2}\big((h^{[0]})^{-1}(\Delta)\big)$ defines a compact curve $c_0$ with boundary which is embedded in $(\R/T\mathbb{Z}\times S^1)^2.$ Hence, if $\varepsilon_0>0$ is sufficiently small, Stability Lemma \ref{lem_stability} and Theorem \ref{thm_open_embedd} determine the embedded curves $c_\varepsilon=\pi_{r_1, \rho_1, r_2, \rho_2}\big((h^{[0]})^{-1}(\Delta)\big)\subset(\R/T\mathbb{Z}\times S^1)^2$ for each $\varepsilon\in[0, \varepsilon_0]$.

In particular, the curves are isotopic, and thus each $c_\varepsilon$ is $C^k$ $O(\varepsilon)$-close to $c_0$. Also, for $\varepsilon>0$, $c_\varepsilon$ describes the unique  (transverse) intersection of $\widetilde{W}^u_\ffs(\mathfrak{p}_\varepsilon)$ and $\widetilde{W}^s_\ffs(\mathfrak{q}_\varepsilon)$ inside $D$. Since $\widetilde{W}^u_\ffs(\mathfrak{p}_\varepsilon)$ and $\widetilde{W}^s_\ffs(\mathfrak{q}_\varepsilon)$ are locally invariant, $c_\varepsilon$ is not only a curve, but also a (partial) flow trajectory. By continuing flowing along $W^u_\ffs(\mathfrak{p}_\varepsilon)$ and $W^s_\ffs(\mathfrak{q}_\varepsilon)$, we extend to $c_\varepsilon$ to a heteroclinic orbit $u_\varepsilon$ from $\mathfrak{p}_\varepsilon$ to $\mathfrak{q}_\varepsilon$. In addition, $u_\varepsilon$ is transversely cut out.

Hence, now we put $\Psi^\varepsilon_{p_0, q_0}(u):=u_\varepsilon$ which determines a well defined map.

\bigskip

\textit{Injectivity of $\Psi^\varepsilon_{p_0, q_0}$:} 

If we consider two distinct $u, \widehat{u}\in\mathcal{M}_{X_{E_0}}(p_0, q_0)$, then by Stability Lemma \ref{lem_stability} $\Psi^\varepsilon_{p_0, q_0}(u)$ and $\Psi^\varepsilon_{p_0, q_0}(\widehat{u})$ remain distinct, provided that $\varepsilon>0$ is small.

\bigskip

\textit{Closeness of $u_0$ and $\Psi^\varepsilon_{p_0, q_0}(u)$ (i.e. $u_\varepsilon$):}

We are going to show that for every $\delta_0>0$ there is a $\varepsilon_0>0$ such that for every $\varepsilon\in(0, \varepsilon_0)$ it holds that 
$$d_{haus}(u_0, u_\varepsilon)<\delta_0.$$

It will be enough to describe the bound for the distance of $\phi^{[0, \infty)}_{slow}(c_0)$ from $\phi^{[0, \infty)}_{\ffs}(c_\varepsilon)$, since the bound in the backward time can be done analogously. From the above, we know that $c_0$ and $c_\varepsilon$ are $O(\varepsilon)$ close. Moreover, $c_0\subset U_0$ and $c_\varepsilon\subset U_\varepsilon$.

Hence, by Gronwall estimates in Fenichel coordinates, we see that after a finite time $T>0$ $\phi^{T}_{slow}(c_0)$ and $\phi^{T}_{\ffs}(c_\varepsilon)$ are still $O(\varepsilon)$-close (here recall that Feinchel coordinates depend smoothly on $\varepsilon$).

We also know by Remark \ref{rem_compact_fenichel} (and proof of Corollary \ref{cor_fs_local}) that $\widetilde{W}^s_{slow}(\mathfrak{q}_0)$ and $\widetilde{W}^s(\mathfrak{q}_\varepsilon)$ (defined by the perturbed slow flow) are isotopic. Hence we can take a small $4$-dimensional ball $B_{\delta_0}(\mathfrak{q}_0)$ with the following property. If $\varepsilon>0$ is small, then by Stability Lemma \ref{lem_stability} the intersections $B_{\delta_0}(\mathfrak{q}_0)\pitchfork \widetilde{W}^s_{slow}(\mathfrak{q}_0)$ and $B_{\delta_0}(\mathfrak{q}_0)\pitchfork \widetilde{W}^s(\mathfrak{q}_\varepsilon)$ are both inflowing manifolds (with respect to given flows). The inflowing property of the intersection follows from eigenvalue estimates as in \cite[prop 2]{Jones2009GeneralizedEL}. Alternatively, one can use instead of $B_{\delta_0}(\mathfrak{q}_0)$ some $O(\delta)$ neighborhood of $\mathfrak{q}_0$ which is coming from the graph transform construction of the local stable manifolds.

But $c_0\subset\widetilde{W}^s_{slow}(\mathfrak{q}_0)$ and $c_\varepsilon\subset \widetilde{W}^s(\mathfrak{q}_\varepsilon)$. So once $\phi^{T}_{slow}(c_0)$ and $\phi^{[0, \infty)}_{\ffs}(c_\varepsilon)$ enter $B_{\delta_0}(\mathfrak{q}_0)$, they got trapped. Thus, we showed that
$$d_{haus}(\phi^{[0, \infty)}_{slow}(c_0), \phi^{[0, \infty)}_{\ffs}(c_\varepsilon))<\delta_0.$$

\bigskip

\textit{Local uniqueness of $\Psi^\varepsilon_{p_0, q_0}(u)$ (i.e. $u_\varepsilon$):}

In particular, we are going to show that there is a closed neighborhood $V_{u_0}$ of $u_0$ with the following property: If $\varepsilon>0$ is sufficiently small, then $u_\varepsilon$ is the unique element of $\mathcal{M}_{X_{E_\varepsilon}}(p_\varepsilon, q_\varepsilon)$ that lies (entirely) in $V_{u_0}$. 

Pick $x_0\in Int(c_0)$ and $\delta_1>0$ small such that $B_{\delta_1}(x_0)\subset Int(D)$. Then we put
$$V_{u_0}^{\mathfrak{q}_0}:=B_{\delta_1}(x_0)\cup T^{\mathfrak{q}_0}_{\delta_1/2}\cup B_{\delta_1}(\mathfrak{q}_0),$$
where the balls $B_{\delta_1}(x_0)$ and $B_{\delta_1}(\mathfrak{q}_0)$ are ``glued'' by a thin tube $T^{\mathfrak{q}_0}_{\delta_1/2}$ which is just a tubular neighborhood $\nu_{\delta_1/2}(\phi^{[0, \infty)}_{slow}(x_0))$.

We can always assume that $\delta_0\ll\delta_1$. Now, note that if $\delta_1>0$ is sufficiently small, then $\partial V_{u_0}^{\mathfrak{q}_0}\cap \widetilde{W}^s_{sing}(\mathfrak{q}_0)=\emptyset.$ Indeed, $\widetilde{W}^s_{slow}(\mathfrak{q}_0)$ is overflowing, so $\phi^{[0, \infty)}_{slow}(x_0)\subset Int(\widetilde{W}^s_{sing}(\mathfrak{q}_0))$. Since, $\widetilde{W}^s_{sing}(\mathfrak{q}_0)$ and $\widetilde{W}^s_{\ffs}(\mathfrak{q}_\varepsilon)$ are isotopic and compact,
$$\partial V_{u_0}^{\mathfrak{q}_0}\cap \widetilde{W}^s_{\ffs}(\mathfrak{q}_\varepsilon)=\emptyset$$
provided that $\varepsilon>0$ is sufficiently small.

Note a crucial observation
\begin{itemize}
\item[$(\Box)$] If $y_0\in B_{\delta_1}(x_0)\setminus\widetilde{W}^s_{\ffs}(\mathfrak{q}_\varepsilon)$ and $\lim_{t\rightarrow\infty}\phi^t_{\ffs}(y_0)=\mathfrak{q}_\varepsilon$, then $\phi^{[0, \infty)}_{\ffs}(y_0)$ does not lie entirely in $V_{u_0}^{\mathfrak{q}_0}$.
\end{itemize}
Indeed, the trajectory from $y_0$ will enter $\widetilde{W}^s_{\ffs}(\mathfrak{q}_\varepsilon)\setminus \mathfrak{q}_\varepsilon$ at some point. This can be done only through $\partial\widetilde{W}^s_{\ffs}(\mathfrak{q}_\varepsilon)$. But $\partial\widetilde{W}^s_{\ffs}(\mathfrak{q}_\varepsilon)$ does not lie in $V_{u_0}^{\mathfrak{q}_0}$.

Now, we analogously construct $V_{u_0}^{\mathfrak{p}_0}$ (with the same constant $\delta_1>0$). Finally, we put 
$$V_{u_0}:=V_{u_0}^{\mathfrak{p}_0}\cup V_{u_0}^{\mathfrak{q}_0}.$$
We can assume that $\delta_1>0$ is so small that $V_{u_0}\setminus
B_{\delta_1}(x_0)$ is a disconnected set. In particular, any heteroclinic orbit form $\mathfrak{p}_\varepsilon$ to $\mathfrak{q}_\varepsilon$, that stays in $V_{u_0}$, needs to go through $B_{\delta_1}(x_0)$.

Now we would like to traverse $V_{u_0}$ with a certain set $K_{u_0}$ that will have nice properties with respect to slow and fast-slow flows. By the continuity of the flow, there is a $1$-ball $\widehat{B}_{\delta_1}(x_0)\subset U_0$ that is transverse to the slow flow and 
$$V_{u_0}\setminus\lbrace(s_1, \theta_1, s_2, \theta_2)\,\vert\,(s_1, s_2)\in\pi_{s_1, s_2}(\widehat{B}_{\delta_1}(x_0))\rbrace$$
 is disconnected, provided that $\delta_1>0$ is small. We put
 $$K_{u_0}:=\lbrace(s_1, \theta_1, s_2, \theta_2)\,\vert\,(s_1, s_2)\in\pi_{s_1, s_2}(\widehat{B}_{\delta_1}(x_0))\rbrace, $$
see also Figure \ref{figure_fast_slow_trap}. Moreover, if $\varepsilon>0$ is sufficiently small, then the following holds. By the same computation as in $(\ref{eqn_sign_lift})$, the fast-slow flow is transverse to $K_0$ ``in the same direction'' as the slow is transverse to $\widehat{B}_{\delta_1}(x_0)$.
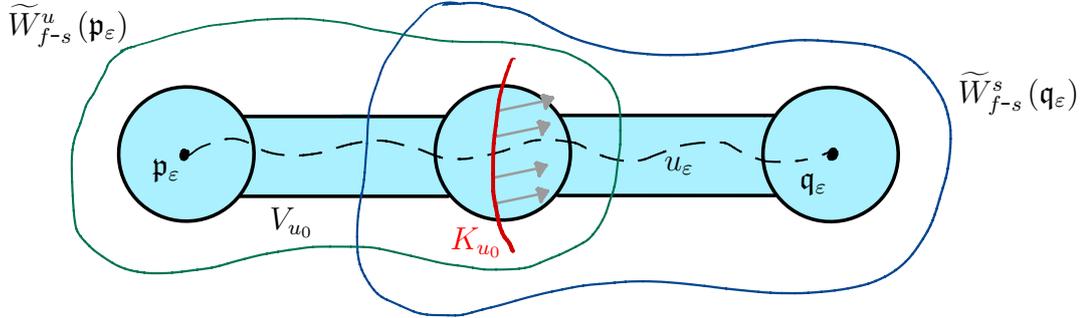
\begin{figure}[!htbp]
\labellist
\pinlabel $\widetilde{W}^u_{\ffs}(\mathfrak{p}_\varepsilon)$ at 7 134
\pinlabel $\widetilde{W}^s_{\ffs}(\mathfrak{q}_\varepsilon)$ at 424 104
\pinlabel $\textcolor{red}{K_{u_0}}$ at 186 38
\pinlabel $\mathfrak{p}_\varepsilon$ at 50 70
\pinlabel $\mathfrak{q}_\varepsilon$ at 334 64
\pinlabel $u_\varepsilon$ at 275 72
\pinlabel $V_{u_0}$ at 105  47
\endlabellist
\centering
\includegraphics[scale=0.85]{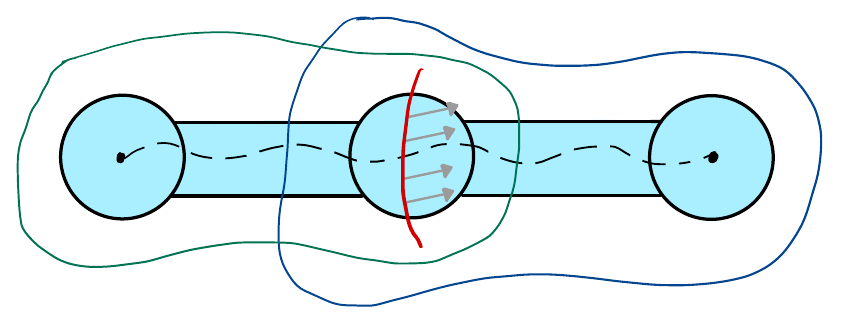}
\vspace{0.3cm}
\caption{A visualization of the splitting of $V_{u_0}$ by $\textcolor{red}{K_{u_0}}$. Note that $\partial\widetilde{W}^u_{\ffs}(\mathfrak{p}_\varepsilon)$ and $\partial\widetilde{W}^u_{\ffs}(\mathfrak{p}_\varepsilon)$ will have now non-empty intersections with $V_{u_0}$. But the intersections will be on the different connected components of the splitting of $V_{u_0}$ by $\textcolor{red}{K_{u_0}}$.}
\label{figure_fast_slow_trap}
\end{figure}

Hence, the following holds
\begin{itemize}
\item If $y_0\in B_{\delta_1}(x_0)\setminus \big(\widetilde{W}^s_{\ffs}(\mathfrak{q}_\varepsilon)\cap \widetilde{W}^u_{\ffs}(\mathfrak{p}_\varepsilon)\big)$ and $y_0$ is contained in some heteroclinic orbit $\widehat{u}_\varepsilon$ from $\mathfrak{p}_\varepsilon$ to $\mathfrak{q}_\varepsilon$, then $\widehat{u}_\varepsilon$ does not lie entirely in $V_{u_0}$.
\end{itemize}

Finally, since $\widetilde{W}^s_{\ffs}(\mathfrak{q}_\varepsilon)$ and $\widetilde{W}^u_{\ffs}(\mathfrak{p}_\varepsilon)$ contribute by the unique heteroclinic orbit - $u_\varepsilon$, the set $V_{u_0}$ has the desired property.

We will also introduce a notation - $N^\varepsilon_{u_0}$ which will be a subset of the $G_{K, \varepsilon}$-strip that projects by $\pi_{s_1, s_2}$ to $\pi_{s_1, s_2}(V_{u_0})$.
\bigskip

\textit{The surjectivity of the map $\Psi^\varepsilon_{p_0, q_0}$:}

This follows from a similar argument as the proof of Lemma \ref{thm_base_ift}; for the sake of completeness, we will still present the proof. 

By contradiction. Let us assume that there is a sequence $\lbrace \widehat{u}_{\varepsilon_n}\rbrace$ of $\varepsilon_n$-solutions with $\varepsilon_n\rightarrow 0$ that are not images of $\Psi^{\varepsilon_n}_{p_0, q_0}$. Which is by the above equivalent to assume that $\widehat{u}_{\varepsilon_n}$ are not in any of the sets $V_{u_0}$.

On the other hand, recall the uniform bounds from Theorem \ref{thm_unif_bound} on $\widehat{u}_{\varepsilon_n}$. Hence, by \cite[Thm 4.7]{frauenfelder2020moduli}, the $\pi_{s_1, s_2}(\widehat{u}_{\varepsilon_n})$ converge in Floer-Gromov $C^1_{loc}$ sense to some broken flow line from $p_0$ to $q_0$. But since $X_{E_0}$ is Morse-Bott-Smale, by the index reason, the limit is in fact a $0$-solution. We will denote it by $u$.  Hence, for $n\gg 0$ it holds that $\pi_{s_1, s_2}(u_{\varepsilon_n})\subset \pi_{s_1, s_2}(V_{u_0})$. Recall also that $u$ is the only $0$-solution inside $\pi_{s_1, s_2}(V_{u_0})$.

Now, since $\pi_{s_1, s_2}(\widehat{u}_{\varepsilon_n})\subset \pi_{s_1, s_2}(V_{u_0})$, in particular $\widehat{u}_{\varepsilon_n}$ avoids special and diagonal points. Hence we can apply Lemma \ref{lemma_solution_strip} and thus each $\widehat{u}_{\varepsilon_n}\subset N^{\varepsilon_n}_{u_0}$. But, for $n\gg 0$, $N^{\varepsilon_n}_{u_0}\subset V_{u_0}$. Also, by the construction $V_{u_0}$ consists the unique $\varepsilon_n$-solution - $u_{\varepsilon_n}$. So $\widehat{u}_{\varepsilon_n}=u_{\varepsilon_n}$. Contradiction.
\end{proof}

\begin{thm}\label{thm_corresp_trees_gsp} Let $m>1$ and $\mathcal{T}\in\clubsuit_m$. Let $p_0\in Cr_1(E_0)\setminus\Delta_0, q_0^1,\dots, q_0^m\in Cr_0(E_0)\setminus\Delta_0$ and $p_\varepsilon, q_\varepsilon^1,\dots, q_\varepsilon^m$ are their unique corresponding critical points of $E_\varepsilon$ in $M_{K, \varepsilon}$. If $\varepsilon>0$ is sufficiently small and $X_{E_0}$ and $X_{E_\varepsilon}$ are generic $K$-approximations of $-\nabla E_0$ and $-\nabla E_\varepsilon$, then there is a map
$$\Psi^\varepsilon_{\mathcal{T}; p_0; q_0^1,\dots, q_0^m}:\clubsuit_{X_{E_0}}(\mathcal{T}; p_0; q_0^1,\dots, q_0^m)\rightarrow \clubsuit_{X_{E_\varepsilon}}(\mathcal{T}; p_\varepsilon; q_\varepsilon^1,\dots, q_\varepsilon^m)$$
which is a bijection. In addition,
$$\clubsuit_{X_{E_\varepsilon}}^{out-out}(\mathcal{T}; p_\varepsilon; q_\varepsilon^1,\dots, q_\varepsilon^m)=\clubsuit_{X_{E_\varepsilon}}(\mathcal{T}; p_\varepsilon; q_\varepsilon^1,\dots, q_\varepsilon^m).$$
\end{thm}

\begin{rem} We are going to show the proof of Theorem \ref{thm_corresp_trees_gsp} for $m= 2, 3$. Then we will claim that we captured all crucial phenomena, and the remaining cases of $m$ are analogous.
\end{rem}

\begin{lemma}\label{lem_tree_2} If $m=2$, then Theorem \ref{thm_corresp_trees_gsp} holds.
\end{lemma}

\begin{proof}
Let us consider a generic $X_{E_0}$ as in the proof of Lemma \ref{lemma_dim_tree_knot}, see also Remark \ref{rem_restr_R}. I.e. if $u\in\clubsuit_{X_{E_0}}(\mathcal{T}; p_0; q_0^1, q_0^2)$, then 
$$u\in S_K\hbox{ and }u\pitchfork\widehat{\widehat{\mathcal{R}}}$$
provided that the constant $\delta_K>0$ (which defines $S_K$) is small.

\bigskip

\textit{Definition of $\Psi^\varepsilon_{\mathcal{T}; p_0; q_0^1, q_0^2}$:}

Let $u\in\clubsuit_{X_{E_0}}(\mathcal{T}; p_0; q_0^1, q_0^2)$.  By the above, we can consider a small neighborhood of $\lbrace p_0\cup u\cup q_0^1\cup q_0^2\rbrace$ in $S_K$ that consists of three (disjoint) contractible components. By Remark \ref{rem_K_approx}, we can lift, with $\ioo$, diffeomorphically this neighborhood to some $U_0\cup U^1_0\cup U^2_0\subset S^{out-out}_0$. Also, $\lbrace p_0\cup u\cup q_0^1\cup q_0^2\rbrace$ is lifted to $\lbrace \mathfrak{p}_0\cup u_0\cup \mathfrak{q}_0^1\cup \mathfrak{q}_0^2\rbrace$, where $u_0$, the lift of $u$, is locally invariant under the slow flow.

Similar to the proof of Theorem \ref{thm_corresp_traj_gsp}, we consider certain manifolds 
\begin{equation}\label{eqn_sing manifolds}
\widetilde{W}^u_{sing}(\mathfrak{p}_0), \widetilde{W}^s_{sing}(\mathfrak{q}_0^1), \widetilde{W}^s_{sing}(\mathfrak{q}_0^2).
\end{equation}
To do this, we take $\widetilde{W}^u_{slow}(\mathfrak{p}_0)$ as an overflowing $\codim 0$ submanifold of $W^u_{slow}(\mathfrak{p}_0)$ that is inside $U_0$ and contains $\lbrace \mathfrak{p}_0\cup (u_0\cap U_0)\rbrace$. And analogously we construct $\widetilde{W}^s_{slow}(\mathfrak{q}_0^1)$ and $\widetilde{W}^s_{slow}(\mathfrak{q}_0^2)$. We require that the manifolds from $(\ref{eqn_sing manifolds})$ lie in the set $(\Phi^{fen}_0)^{-1}(B_{\delta, U_{fen}})$ for some Fenichel coordinates around $U_0\cup U_0^1\cup U_0^2$ with small $\delta>0$. This fully determines $(\ref{eqn_sing manifolds})$.

Now, we would like to construct certain homotopy of the maps that resemble the evaluation maps from Chapter \ref{ch:file8}.

First, we consider a product of three flat tori $(\R/T\mathbb{Z}\times S^1)^6$ with coordinates $((s_1, \theta_1, s_2, \theta_2), (y_1, \alpha_1, y_3, \alpha_3), (z_2, \beta_2, z_3, \beta_3))$. And in these coordinates, we will view $\widetilde{W}^u_{sing}(\mathfrak{p}_0)\times \widetilde{W}^s_{sing}(\mathfrak{q}_0^1)\times \widetilde{W}^s_{sing}(\mathfrak{q}_0^2)$ as an $8$-dimensional submanifold.

In particular, we can locally write
\begin{equation}\label{eqn_w_u}
\widetilde{W}^u_{sing}(\mathfrak{p}_0)=\big(s_1(r), \theta_1(s_1(r), s_2(r)), s_2(r), \theta_2\big),
\end{equation} 
where $(s_1(r), s_2(r))$ is a local parametrization of $1$-dimensional $W^{u}_{X_{E_0}}(p_0)$. Then $\theta_1(s_1(r), s_2(r))$ is a function which is determined by the lift $\ioo$. Analogously, we have the local descriptions
\begin{equation}\label{eqn_w_s1}
\widetilde{W}^s_{sing}(\mathfrak{q}_0^1)=\big(y_1, \alpha_1, y_2, \alpha_3(y_1, y_3)\big),
\end{equation} 
and
\begin{equation}\label{eqn_w_s2}
\widetilde{W}^s_{sing}(\mathfrak{q}_0^2)=\big(z_2, \beta_2(z_2, z_3), z_3, \beta_3\big),
\end{equation} 

Let $x$ be the unique point of $u\pitchfork \mathcal{R}$. Now, by $D_x$ we will denote a silly auxiliary $12$-dimensional ball 
$$B_{\delta_1}\big((\ioo(x), \ioo\circ\varphi^L(x), \ioo\circ\varphi^U(x))\big)\subset(\R/T\mathbb{Z}\times S^1)^6.$$
If $\delta_1>0$ is sufficiently small, then the following holds
\begin{itemize}
\item $D_x\subset Int\big(\widetilde{W}^u_{sing}(\mathfrak{p}_0)\times \widetilde{W}^s_{sing}(\mathfrak{q}_0^1)\times \widetilde{W}^s_{sing}(\mathfrak{q}_0^2)\big)$.
\item $D_x\pitchfork \widetilde{W}^u_{sing}(\mathfrak{p}_0)\times \widetilde{W}^s_{sing}(\mathfrak{q}_0^1)\times \widetilde{W}^s_{sing}(\mathfrak{q}_0^2).$
\item $D_x$ does not contain any other lifted tree from $\clubsuit_{X_{E_0}}(\mathcal{T}; p_0; q_0^1, q_0^2)$ then $u$. In more detail, $$(\emph{Ev}_{\mathcal{T}; p_0; q_0^1, q_0^2})^{-1}(0)\cap[0, 1]\times\pi_{s_1, s_2, y_1, y_3, z_2, z_3}(D_x)$$ consists of a single point (the representative of $u$ from the bijection in Corollary \ref{cor_bij_treesT}). Recall also that $\clubsuit_{X_{E_0}}(\mathcal{T}; p_0; q_0^1, q_0^2)$ is a finite set.
\end{itemize}
See also Figure \ref{figure_bifur_disc}.

\begin{figure}[!htbp]
\vspace{0.2cm}
\labellist
\pinlabel $(\R/T\mathbb{Z})^2$ at 320 235
\endlabellist
\centering
\includegraphics[scale=0.84]{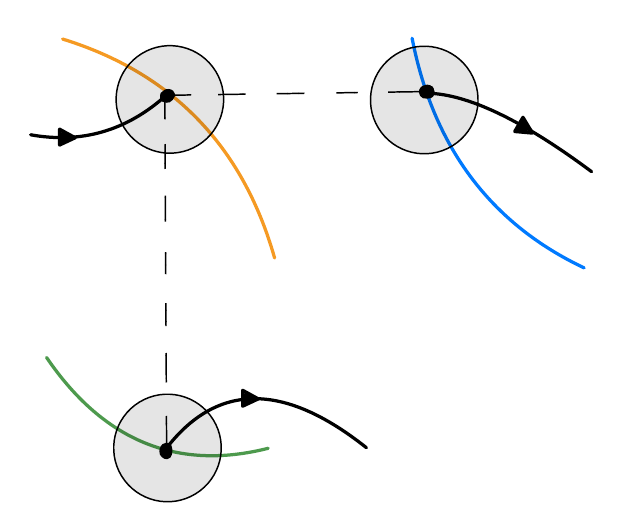}
\vspace{0.4cm}
\caption{A neighborhood of the bifurcation of $u$ depicted in the single configuration space $(\R/T\mathbb{Z})^2$. In more detail, after the canonical identification the curves $\textcolor{orange}{\mathcal{R}}, \textcolor{teal}{\mathcal{R}^L}, \textcolor{blue}{\mathcal{R}^U}$ and the three projections $\pi_{s_1, s_2}(D_x), \pi_{y_1, y_3}(D_x), \pi_{z_3, z_2}(D_x)$ are visualized in the same space.}
\label{figure_bifur_disc}
\end{figure}

We take one more time the product of three flat tori $(\R/T\mathbb{Z}\times S^1)^6$ but now with ``bold'' coordinates $((\bm{s}_1, \bm{\theta}_1, \bm{s}_2, \bm{\theta}_2), (\bm{y}_1, \bm{\alpha}_1, \bm{y}_3, \bm{\alpha}_3), (\bm{z}_2, \bm{\beta}_2, \bm{z}_3, \bm{\beta}_3))$. Then we canonically embed $D_x$ into this product $(\R/T\mathbb{Z}\times S^1)^6$. We will denote the image of the embedding still by $D_x$.

Now, as in Theorem \ref{thm_corresp_traj_gsp}, if $\varepsilon_0>0$ is sufficiently small, the Fenichel theory (Corollary \ref{cor_fs_local} and Remark \ref{rem_compact_fenichel}) gives us a $C^k$-isotopy $H_1$:
\begin{align*}
\big(D_x\times[0, 1]\times\widetilde{W}^u_{sing}(\mathfrak{p}_0)\times \widetilde{W}^s_{sing}(\mathfrak{q}_0^1)\times &\widetilde{W}^s_{sing}(\mathfrak{q}_0^2)\big)\times[0, \varepsilon_0]\longrightarrow\\
 &(\R/T\mathbb{Z}\times S^1)^6\times[0, 1]\times(\R/T\mathbb{Z}\times S^1)^6
\end{align*}
from the canonical inclusion to $h^{[\varepsilon_0]}_1:$
\begin{align*}
\big(D_x, [0, 1], \widetilde{W}^u_{sing}(\mathfrak{p}_0), \widetilde{W}^s_{sing}(\mathfrak{q}_0^1), &\widetilde{W}^s_{sing}(\mathfrak{q}_0^2)\big)\longmapsto\\
&\big(D_x, [0, 1], \widetilde{W}^u_{\ffs}(\mathfrak{p}_{\varepsilon_0}), \widetilde{W}^s_{\ffs}(\mathfrak{q}_{\varepsilon_0}^1), \widetilde{W}^s_{\ffs}(\mathfrak{q}_{\varepsilon_0}^2)\big).
\end{align*}
We remark, that the convention of the coordinates on $(\R/T\mathbb{Z}\times S^1)^6\times[0, 1]\times(\R/T\mathbb{Z}\times S^1)^6$ is the following. On the first factor of $(\R/T\mathbb{Z}\times S^1)^6\times[0, 1]\times(\R/T\mathbb{Z}\times S^1)^6$ we have the ``bold'' coordinates. So by $\pi_{\text{non-bold}}$ we understand the canonical projection onto the remaining coordinates.

Next, we define a smooth homotopy 
$$H_2:(\R/T\mathbb{Z}\times S^1)^6\times[0, 1]\times(\R/T\mathbb{Z}\times S^1)^6\rightarrow(\R/T\mathbb{Z}\times S^1)^6\times (S^1)^3\times(\R/T\mathbb{Z})^3\times \R^3$$
by
$$h^{[\varepsilon]}_2:=\big(\ev_{\text{auxil}}, \ev_\varepsilon\circ\pi_{\text{non-bold}}\big)$$
for every $\varepsilon\in[0, \varepsilon]$. Here, the maps $\ev_\varepsilon$ were introduced in Definition \ref{defn_evE}. The map $\ev_{\text{auxil}}$ is just given as a difference of each ``bold'' coordinate with the corresponding ``non-bold'' coordinate. I.e. for example $\bm{s}_1-s_1$. Note that the coordinate $\ell\in[0, 1]$ does not have a ``bold'' coordinate equivalent, and hence does not appear in $\ev_{\text{auxil}}$.

Finally, we define a $C^k$ homotopy $H:$
\begin{align}
\begin{split}\label{eqn_homotop_tree}
\big(D_x\times[0, 1]\times\widetilde{W}^u_{sing}(\mathfrak{p}_0)\times &\widetilde{W}^s_{sing}(\mathfrak{q}_0^1)\times \widetilde{W}^s_{sing}(\mathfrak{q}_0^2)\big)\times[0, \varepsilon_0]\longrightarrow\\
 &(\R/T\mathbb{Z}\times S^1)^6\times (S^1)^3\times(\R/T\mathbb{Z})^3\times \R^3
\end{split}
\end{align}
as the horizontal composition $H_2\circ H_1$, i.e.
$$h^{[\varepsilon]}:=h^{[\varepsilon]}_2\circ h^{[\varepsilon]}_1.$$

It is crucial to observe that
\begin{itemize}
\item[$(\Delta)$]\label{observation_inters}$h^{[0]}\pitchfork 0$ and $(h^{[0]})^{-1}(0)$ consists of a single point, where in addition  $$\pi_{s_1, \theta_1, s_2, \theta_2}\big((h^{[0]})^{-1}(0)\big)=\ioo(x).$$
\end{itemize}

Now, we would like to enlight $(\Delta)$. Let $\mathbbm{x}_0:=(h^{[0]})^{-1}(0)$. First, note that $h_1^{[0]}$ is just an inclusion that adjusts the domain, where we would like to count the zeros of $h^{[0]}$. By the descriptions $(\ref{eqn_w_u}), (\ref{eqn_w_s1}), (\ref{eqn_w_s2})$ we see that at $\mathbbm{x}_0$ the coordinates $\alpha_1, \theta_2, \beta_3$ are uniquely determined by $\theta_1(s_1(r), s_2(r)), \beta_2(z_2, z_3), \alpha_1(y_1, y_3)$, respectively. Hence, by the definition of $D_x$ and relation (\ref{eqn_ev0}) we obtain that $\mathbbm{x}_0$ is a single point and $\pi_{s_1, \theta_1, s_2, \theta_2}(\mathbbm{x}_0)=\ioo(x).$ In addition, by Lemma \ref{cor_transverse_hitR} it follows that $\ev_0\circ\pi_{\text{non-bold}}\pitchfork0$. Finally, from the definition of $D_x$, we immediately get that $\ev_{\text{auxil}}\pitchfork 0$.

Hence, if $\varepsilon_0>0$ is small, Stability Lemma \ref{lem_stability} determines a $C^k$-family of points 
$$\left\lbrace\mathbbm{x}_\varepsilon:=\pi_{\text{non-bold}}\big((h^{[\varepsilon]})^{-1}(0)\big)\right\rbrace_{\varepsilon\in[0, \varepsilon_0]}$$
which by Corollary \ref{cor_coresp_tree_tor} uniquely determines a family of trees
$$\lbrace u_\varepsilon\in\clubsuit_{X_{E_\varepsilon}}(\mathcal{T}; p_\varepsilon; q_\varepsilon^1, q_\varepsilon^2)\rbrace_{\varepsilon\in[0, \varepsilon_0]}.$$
We will also use a notation 
$$x_\varepsilon:=\pi_{s_1, \theta_1, s_2, \theta_2}(\mathbbm{x}_\varepsilon).$$

Finally, we put $\Psi^\varepsilon_{\mathcal{T}; p_0; q_0^1; q_0^2}(u):=u_\varepsilon$ which determines a well defined map.

\bigskip

\textit{Injectivity of $\Psi^\varepsilon_{\mathcal{T}; p_0; q_0^1, q_0^2}$:} 

It is a straightforward consequence of the Stability Lemma \ref{lem_stability} and Corollary \ref{cor_coresp_tree_tor}.

\bigskip

\textit{Closeness of $u_0$ and $\Psi^\varepsilon_{\mathcal{T}; p_0; q_0^1, q_0^2}(u)$ (i.e. $u_\varepsilon$):}

We would like to show that for every $\delta_0>0$ there is a $\varepsilon_0>0$ such that for every $\varepsilon\in(0, \varepsilon_0)$ it holds that
$$d_{haus}(u_0, u_\varepsilon)<\delta_0.$$
For this, we have to compare the distance of each edge of $u_\varepsilon$ from $u_0$. 

For example, we have to compare the distance of $\phi^{(-\infty, 0]}_{\ffs}(x_\varepsilon)$ from $\phi^{(-\infty, 0]}_{slow}(x_0)$.  By Stability Lemma \ref{lem_stability} we know that $x_0$ and $x_\varepsilon$ are $O(\varepsilon)$ close. But since $x_0\in U_0$ and $x_\varepsilon\in W^u(U_\varepsilon)$, it follows that also $x_0$ and $\pi_a(x_\varepsilon)$ are $O(\varepsilon)$ close (here $\pi_a$ is the base projection from Fenichel fibration - Theorems \ref{thm_fen3} and  \ref{defn_fen_coord}). So by the same argument as in Theorem \ref{thm_corresp_traj_gsp} we can bound the distance of $\phi^{(-\infty, 0]}_\ffs(\pi_a(x_\varepsilon))$ from $\phi^{(-\infty, 0]}_{slow}(x_0)$. Moreover the distance of $x_\varepsilon$ and $\pi_a(x_\varepsilon)$ exponentially decay in backward-time by Theorems \ref{thm_fen3} which imply the bound for the distance of $\phi^{(-\infty, 0]}_\ffs(x_\varepsilon)$ from $\phi^{(-\infty, 0]}_{slow}(x_0)$.

\bigskip

\textit{Local uniqueness of $\Psi^\varepsilon_{\mathcal{T}; p_0; q_0^1, q_0^2}(u)$ (i.e. $u_\varepsilon$):}

In particular, we are going to show that there is a closed neighborhood $V_{u_0}$ of $u_0$ with the following property: If $\varepsilon>0$ is sufficiently small, then $u_\varepsilon$ is the unique element of $\clubsuit_{X_{E_\varepsilon}}(\mathcal{T}; p_\varepsilon; q_\varepsilon^1, q_\varepsilon^2)$ that lies (entirely) in $V_{u_0}$. 

In fact, this is an easier version of the local uniqueness from the proof of Theorem \ref{thm_corresp_traj_gsp}. We just construct $V_{u_0}\subset(\R/T\mathbb{Z}\times S^1)^2$ as the union of small neighborhoods $V^{\mathfrak{p}_0}_{u_0}, V^{\mathfrak{q}_0^1}_{u_0}, V^{\mathfrak{q}_0^2}_{u_0}$ which are around the connected components of $u_0$. Then the local uniqueness will follow from Stability Lemma \ref{lem_stability}, definition of $D_x$ together with the observation $(\Box)$ from the proof of Theorem \ref{thm_corresp_traj_gsp}.

\bigskip

\textit{Trap for trees:}

In particular, using a similar technique as was in Lemma \ref{lemma_solution_strip}, we are going to construct a subset $N^\varepsilon_{u_0}\subset V_{u_0}$ with the following property: if $\varepsilon>0$ is small and $\widehat{u}_\varepsilon\in\clubsuit_{X_{E_\varepsilon}}(\mathcal{T}; p_\varepsilon; q_\varepsilon^1, q_\varepsilon^2)$ such that $\pi_{s_1, s_2}(\widehat{u}_\varepsilon)\subset\pi_{s_1, s_2}(N^\varepsilon_{u_0})$, then $\widehat{u}_\varepsilon\subset N^\varepsilon_{u_0}$.

We define $N^{\varepsilon}_{u_0}$ as the union of the following three sets.
\begin{align*}
N^{\varepsilon, \mathfrak{p}_0}_{u_0}&=\big\lbrace (s_1, \theta_1, s_2, \theta_2)\in M_{K, \varepsilon}\,\vert\,(s_1, s_2)\in\pi_{s_1, s_2}(V^{\mathfrak{p}_0}_{u_0})\wedge\vert\langle P, v_1^\bot\rangle\vert<\varepsilon^{1/2}\wedge\vert\langle P, v_2^\bot\rangle\vert<\varepsilon^{1/3}\big\rbrace,\\
N^{\varepsilon, \mathfrak{q}_0^1}_{u_0}&=\big\lbrace (s_1, \theta_1, s_2, \theta_2)\in M_{K, \varepsilon}\,\vert\,(s_1, s_2)\in\pi_{s_1, s_2}(V^{\mathfrak{q}_0^1}_{u_0})\wedge\vert\langle P, v_1^\bot\rangle\vert<\varepsilon^{1/3}\wedge\vert\langle P, v_2^\bot\rangle\vert<\varepsilon^{1/2}\big\rbrace,\\
N^{\varepsilon, \mathfrak{q}_0^2}_{u_0}&=\big\lbrace (s_1, \theta_1, s_2, \theta_2)\in M_{K, \varepsilon}\,\vert\,(s_1, s_2)\in\pi_{s_1, s_2}(V^{\mathfrak{q}_0^2}_{u_0})\wedge\vert\langle P, v_1^\bot\rangle\vert<\varepsilon^{1/3}\wedge\vert\langle P, v_2^\bot\rangle\vert<\varepsilon^{1/2}\big\rbrace.
\end{align*}
It remains to enlight, why the defined set $N^{\varepsilon}_{u_0}$ has the desired properties.

So let $\widehat{u}_\varepsilon\in\clubsuit_{X_{E_\varepsilon}}(\mathcal{T}; p_\varepsilon; q_\varepsilon^1, q_\varepsilon^2)$ such that $\pi_{s_1, s_2}(\widehat{u}_\varepsilon)\subset\pi_{s_1, s_2}(N^\varepsilon_{u_0})$.  Recall also that $\pi_{s_1, s_2}(N^\varepsilon_{u_0})\subset S_K$ and $\mathfrak{p}_\varepsilon\subset S^{out-out}$. Hence from Figure \ref{figure_laypunov_trap} in the proof of Lemma \ref{lemma_solution_strip} we immediately see that over $\pi_{s_1, s_2}(N^{\varepsilon, \mathfrak{p}_0}_{u_0})$ $\widehat{u}_\varepsilon$ has to satisfy the relations $F_1^{[\varepsilon]}\geq0\wedge\vert\langle P, v_1^\bot\rangle\vert<\varepsilon^{1/2}$. Here we remark that since $\varepsilon^{1/2}>c_G\varepsilon$, the ``inward-pointing behavior'' of $X_{E_\varepsilon}$ along $\langle P, v_1^\bot\rangle=\varepsilon^{1/2}$ and $\langle P, v_1^\bot\rangle=C_G \varepsilon$ is the same, see also Theorem \ref{thm_g} $(iii.)$.

Analogously we obtain that over $\pi_{s_1, s_2}(N^{\varepsilon, \mathfrak{q}_0^1}_{u_0})$ $\widehat{u}_\varepsilon$ has to satisfy the relations $F_2^{[\varepsilon]}\geq0\wedge\vert\langle P, v_2^\bot\rangle\vert<\varepsilon^{1/2}$. In fact, by Lemma \ref{lemma_aux_strip}, $F_2^{[\varepsilon]}\geq\widehat{\delta}$ for some $\widehat{\delta}>0$.

Let $y_\varepsilon\in\widehat{\widehat{\mathcal{R}_\varepsilon}}.$ Hence, by Lemma \ref{lem_aux_intersect_tor_chords}, we see that also the point $y_\varepsilon$ has to satisfy $F_2^{[\varepsilon]}\geq0\wedge\vert\langle P, v_2^\bot\rangle\vert<\varepsilon^{1/3}$ provided that $\varepsilon>0$ is small. Hence the same argument as in Lemma \ref{lemma_solution_strip} implies that over $\pi_{s_1, s_2}(N^{\varepsilon, \mathfrak{p}_0}_{u_0})$ $\widehat{u}_\varepsilon$ has also to satisfy the relations $F_2^{[\varepsilon]}\geq0\wedge\vert\langle P, v_2^\bot\rangle\vert<\varepsilon^{1/3}$. In particular, over $\pi_{s_1, s_2}(N^{\varepsilon, \mathfrak{p}_0}_{u_0})$ $\widehat{u}_\varepsilon$ is contained in $N^{\varepsilon, \mathfrak{p}_0}_{u_0}$. Arguments for the remaining sets are analogous.

\bigskip

\textit{Surjectivity of $\Psi^\varepsilon_{\mathcal{T}; p_0; q_0^1, q_0^2}$:}
The argument is almost the same as in the surjectivity in Theorem \ref{thm_corresp_traj_gsp}. The only difference in the logic is that now $\pi_{s_1, s_2}
(\widehat{u}_{\varepsilon_n})$ does not subsequentially $C^1_{loc}$ converge to some (potentially) broken flow line. However, by Lemma \ref{lem_preim_treeT} we know that
$$\pi_{\ell, s_1, s_2, y_1, y_3, z_2, z_3}\big(\widehat{\widehat{\ev_{\varepsilon_n}}}^{-1}(0)\big)\rightarrow \widehat{\widehat{\evK}}^{-1}(0)$$
as $\varepsilon_n\rightarrow 0$. So by \cite[thm 4.22]{mescher2016perturbedgradientflowtrees} $\pi_{s_1, s_2}
(\widehat{u}_{\varepsilon_n})$ subsequentially $C^1_{loc}$ converge to some broken Morse flow tree $u$ with leaves corresponding to $p_0, q_0^1, q_0^2$ and a single interior vertex. But by generic assumptions on $X_{E_0}$ it follows that $u\in \clubsuit_{X_{E_0}}(\mathcal{T}; p_0; q_0^1, q_0^2)$.
\end{proof}

\begin{lemma}\label{lem_corres_tree3} If $m=3$, then Theorem \ref{thm_corresp_trees_gsp} holds.
\end{lemma}

\begin{proof}
Let us consider a tree $\mathcal{T}\in \clubsuit_3$ as in Figure \ref{figure_tree3}.
\begin{figure}[!htbp]
\labellist
\pinlabel $v_0$ at 50 718
\pinlabel $v_1$ at 378 620
\pinlabel $v_2$ at 378 690
\pinlabel $v_3$ at 378 760
\pinlabel $v_a^{int}$ at 270 635
\pinlabel $v_b^{int}$ at 180 670
\endlabellist
\centering
\includegraphics[scale=0.58]{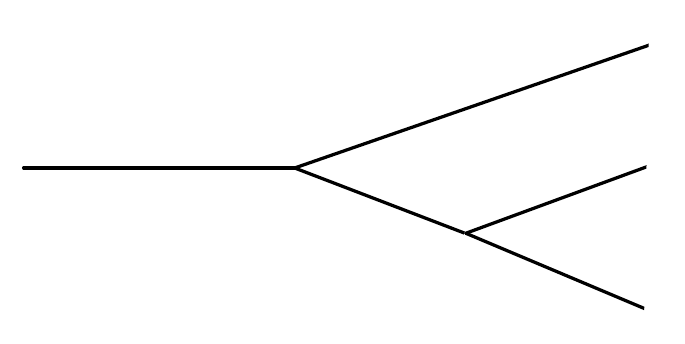}
\vspace{0.3cm}
\caption{The chosen tree $\mathcal{T}\in \clubsuit_3$. $\mathcal{T}$ has one interior edge, the vertex $v_0$ is the root and recall that the exterior vertices $v_1, v_2, v_3$ are ordered.}
\label{figure_tree3}
\end{figure}
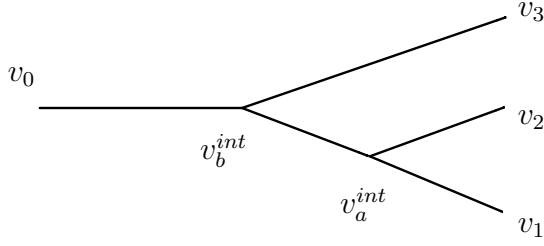

The key difference from Lemma \ref{lem_tree_2} is that now $\mathcal{T}$ contains an interior edge. Interior edges, as we know from Corollaries \ref{cor_bij_treesT} and \ref{cor_coresp_tree_tor}, correspond to the space of positive finite-length trajectories. Such a space is clearly neither a stable, nor an unstable manifold, so the perturbations of such a space can not be fully captured by Fenichel theory. However, we will be able to construct the perturbed trees \textit{iteratively under the backward flow $\phi^{-t}_{X_{E_\varepsilon}}$} as was discussed in Remark \ref{rem_backward_tree}. In particular, in order to track the perturbed trajectories along the interior edge, we will use the Exchange lemma \ref{lemma_exchangeC}.

\bigskip

\textit{Definition of $\Psi^\varepsilon_{\mathcal{T}; p_0; q_0^1, q_0^2, q_0^3}$:}

First, as in Lemma \ref{lem_tree_2}, we take a generic $X_{E_0}$ and $u\in\clubsuit_{X_{E_0}}(\mathcal{T}; p_0; q_0^1, q_0^2, q_0^3)$. We also lift with the map $\ioo$ a small neighborhood of $\lbrace p_0\cup u\cup q^1_0\cup q^2_0\cup q^3_0\rbrace$ to $U_0\cup U_0^{int}\cup U_0^1\cup U_0^2\cup U_0^3\subset S^{out-out}$. The lift of $u$ will be denoted by $u_0$. We construct manifolds
$$\widetilde{W}^u_{sing}(\mathfrak{p}_0), \widetilde{W}^s_{sing}(\mathfrak{q}_0^1), \widetilde{W}^s_{sing}(\mathfrak{q}_0^2), \widetilde{W}^s_{sing}(\mathfrak{q}_0^3)$$
analogously to (\ref{eqn_sing manifolds}). Let $x\in u\cap \mathcal{R}$ be the unique point that corresponds to the bifurcation of $u$ at the vertex $v_a^{int}$. Then we construct a small $12$-dimensional ball $D_x$ as in Lemma \ref{lem_tree_2}. Now, similarly to (\ref{eqn_homotop_tree}) we define a $C^k$ homotopy 
\begin{align*}
H:\big(D_x\times[0, 1]\times(\R/T\mathbb{Z}\times S^1)^2 \times &\widetilde{W}^s_{sing}(\mathfrak{q}_0^1)\times \widetilde{W}^s_{sing}(\mathfrak{q}_0^2)\big)\times[0, \varepsilon_0]\longrightarrow\\
 &(\R/T\mathbb{Z}\times S^1)^6\times (S^1)^3\times(\R/T\mathbb{Z})^3\times \R^3.
\end{align*}
Note that the domain of $H$ now contains a whole copy of the configuration space $(\R/T\mathbb{Z}\times S^1)^2$ instead of $\widetilde{W}^u_{sing}(\mathfrak{p}_0)$ and on this space act $h^{[\varepsilon]}_1$ by the identity. Recall also that we see $(\R/T\mathbb{Z}\times S^1)^2\times \widetilde{W}^s_{sing}(\mathfrak{q}_0^1)\times \widetilde{W}^s_{sing}(\mathfrak{q}_0^2)$ as a submanifold of the manifold $(\R/T\mathbb{Z}\times S^1)^6$ with the coordinates  $((s_1, \theta_1, s_2, \theta_2), (y_1, \alpha_1, y_3, \alpha_3), (z_2, \beta_2, z_3, \beta_3))$. 

Similarly to $(\Delta)$ from Lemma \ref{lem_tree_2} it follows that
\begin{itemize}
\item $h^{[0]}\pitchfork 0$ and $\pi_{s_1, \theta_1, s_2, \theta_2}\big((h^{[0]})^{-1}(0)\big)$ is an embedded $2$-dimensional submanifold of $\pi_{s_1, \theta_1, s_2, \theta_2}(D_x)$ which can be locally written as 
$$\big( s_1(r), \theta_1, s_2(r), \theta_2(s_1(r), s_2(r))\big),$$
where $(s_1(r), s_2(r))$ is a local parametrization of the $1$-dimensional manifold $\mathcal{R}$, $\theta_2(s_1(r), s_2(r))$ is a function determined by $\ioo$ and $\theta_1$ vary in intervals around the function $\theta_1(s_1(r), s_2(r))$ (which is also given by the lift $\ioo$).
\end{itemize}

Hence if $\varepsilon_0>0$ is small, Stability Lemma \ref{lem_stability} and Theorem \ref{thm_open_embedd} determine a $C^k$-family of isotopic submanifolds 
$$\left\lbrace\mathcal{R}_\varepsilon^{\theta_1}:=\pi_{\text{non-bold}}\big((h^{[\varepsilon]})^{-1}(0)\big)\right\rbrace_{\varepsilon\in[0, \varepsilon_0]}$$
of $\pi_{s_1, \theta_1, s_2, \theta_2}(D_x)$.

By the genericity of $X_{E_0}$, we know that at $x$ it holds that $u\pitchfork \mathcal{R}$. Also by Lemma \ref{lem_R_crit} $x\notin Crit(E_0)$. So if $D_x$ is sufficiently small, then $\mathcal{R}_0^{\theta_1}\cap U_0^{int}$ is transverse to the slow flow.

Let $\overline{t}>0$ be the time length of the trajectory of $u$ that corresponds to the interior edge. Now, in a small neighborhood of $\phi^{[-\overline{t}, 0]}_{slow}(\ioo(x))$ we would like to flow from $\mathcal{R}^{\theta_1}_\varepsilon$ with the fast-slow flow for a (slow) time $-\overline{t}-\delta_1$, where $\delta_1>0$ is small. The result will be denoted by $(\mathcal{R}^{\theta_1}_\varepsilon)^\ast$. Since the evolution of $\mathcal{R}^{\theta_1}_\varepsilon$ in time has to be captured by Exchange Lemmata \ref{lemma_exchangeB} or \ref{lemma_exchangeC}, we can understand $(\mathcal{R}^{\theta_1}_\varepsilon)^\ast$ only as a set that comes from some $B_{U_{fen}, \delta}$ given by Fenichel coordinates. In particular, by Remark \ref{rem_inclination1} $(\mathcal{R}^{\theta_1}_\varepsilon)^\ast$ will not contain the points that evolve too far in the $\theta_1$-direction from $U_\varepsilon^{int}$ (or more precisely, $(\mathcal{R}^{\theta_1}_\varepsilon)^\ast$ will not overflow from the set $(\Phi^{fen}_\varepsilon)^{-1}(\Vert b\Vert\leq\delta)$).

Let $\delta_2>0$ be small. Then by Exchange Lemma \ref{lemma_exchangeB} and Remark \ref{rem_combination_exchange} there is a closed $O(\delta_2)$ neighborhood $B_{\delta_2}$ of $\phi^{-\overline{t}}_{slow}(\ioo(x))$ such that $B_{\delta_2}\cap W^{s}_{sing}(U_0^{int})$ together with $\lbrace B_{\delta_2}\cap(\mathcal{R}^{\theta_1}_\varepsilon)^\ast \rbrace_{\varepsilon\in(0, \varepsilon_0]}$ form a $C^k$ family of isotopic $3$-dimensional submanifolds of $(\R/T\mathbb{Z}\times S^1)^2$, provided that $\varepsilon_0>0$ is small. We recall that as $\varepsilon\rightarrow0$ $(\mathcal{R}^{\theta_1}_\varepsilon)^\ast$ are approaching the stable manifold $W^{s}_{sing}(U_0^{int})$, since we are using the backward-flow.

Now, analogously to Lemma \ref{lem_tree_2}, we are in the situation where we want to perturb a certain transverse intersection which corresponds to the interior vertex $v_b^{int}$. In more detail, we want to perturb the ``intersection'' of $\widetilde{W}^{u}_{sing}(\mathfrak{p}_0), \widetilde{W}^{s}_{sing}(\mathfrak{q}_0^3)$ and $B_{\delta_2}\cap W^{s}_{sing}(U_0^{int})$ together with some small auxiliary ($12$-dimensional) ball $\widehat{D}\subset(\R/T\mathbb{Z}\times S^1)^{6}$ (radius of $\widehat{D}$ will be much smaller then $\delta_2>0$). In particular, Stability Lemma \ref{lem_stability} will gives us a $C^k$-family of points 
$$\lbrace\mathbbm{y}_\varepsilon\rbrace_{\varepsilon\in[0, \varepsilon_1]}\subset\widehat{D},$$
provided that $\varepsilon_1>0$ is small.

Finally, note that by Remark \ref{rem_backward_tree} and Corollary \ref{cor_coresp_tree_tor} each of $\mathbbm{y}_\varepsilon$ determines the unique tree $u_\varepsilon\in\clubsuit_{X_{E_\varepsilon}}(\mathcal{T}; p_\varepsilon; q_\varepsilon^1, q_\varepsilon^2, q_\varepsilon^3)$. In particular, the assignment
$u\mapsto u_\varepsilon$ gives us a well defined map $\Psi^\varepsilon_{\mathcal{T}; p_0; q_0^1, q_0^2, q_0^3}$.

\bigskip

\textit{Injectivity of $\Psi^\varepsilon_{\mathcal{T}; p_0; q_0^1, q_0^2, q_0^3}$:} 

By Stability Lemma \ref{lem_stability} and finiteness of $\clubsuit_{X_{E_0}}(\mathcal{T}; p_0; q_0^1, q_0^2, q_0^3)$ we obtain that two different elements of $\clubsuit_{X_{E_0}}(\mathcal{T}; p_0; q_0^1, q_0^2, q_0^3)$ will contribute by two disjoint sets $\mathcal{R}^{\theta_1}_\varepsilon$. Then the injectivity of $\Psi^\varepsilon_{\mathcal{T}; p_0; q_0^1, q_0^2, q_0^3}$ follows.

\bigskip

\textit{Closeness of $u_0$ and $\Psi^\varepsilon_{\mathcal{T}; p_0; q_0^1, q_0^2, q_0^3}(u)$ (i.e. $u_\varepsilon$):}

We would like to show that for every $\delta_0>0$ there is a $\varepsilon_0>0$ such that for every $\varepsilon\in(0, \varepsilon_0]$ it holds that
$$d_{haus}(u_0, u_\varepsilon)<\delta_0.$$
For this, we have to compare the distance of each edge of $u_\varepsilon$ from $u_0$. The cases of edges between $v_b^{int}$ and $v_0$ or $v_3$ are solved as in Lemma \ref{lem_tree_2}. So let us focus on the interior edge.

For $\varepsilon\in[0, \varepsilon_1]$ let us denote by $y_\varepsilon$ the unique projection of $\mathbbm{y}_\varepsilon$ into $(\mathcal{R}^{\theta_1}_\varepsilon)^\ast$. In particular, by the construction of $\mathbbm{y}_\varepsilon$ it follows that $y_\varepsilon$ are $O(\varepsilon)$ close to $y_0$. Also, by Remark \ref{rem_inclination1} the fast-slow trajectory from $y_\varepsilon$ reaches $\mathcal{R}^{\theta_1}_\varepsilon$ after a slow time $\overline{t}_1+O(\varepsilon)$. Let denote such a trajectory by $\widetilde{u}_\varepsilon$.

Recall also that $\mathcal{R}^{\theta_1}_\varepsilon$ and $\mathcal{R}^{\theta_1}_0$ are $O(\varepsilon)$-close in the ``unstable'' direction $\theta_2$. Hence by Remark \ref{rem_inclination1} and the monotonicity estimates from Remark \ref{rem_monot_estim} it follows that $\widetilde{u}_\varepsilon$ and $\widetilde{u}_0$ are $O(\varepsilon)$ close. Hence, also the terminate points of $\widetilde{u}_\varepsilon$ and $\widetilde{u}_0$ are $O(\varepsilon)$ close, and thus also their unique lifts to $(h^{[\varepsilon]})^{-1}(0)$ and $(h^{[0]})^{-1}(0)$ are $O(\varepsilon)$-close. Then the distance of the remaining edges of $u_\varepsilon$ from $u_0$ follows again as in Lemma \ref{lem_tree_2}.

\bigskip

\textit{Local uniqueness of $\Psi^\varepsilon_{\mathcal{T}; p_0; q_0^1, q_0^2, q_0^3}(u)$ (i.e. $u_\varepsilon$):}

In particular, we are going to show that there is a closed neighborhood $V_{u_0}$ of $u_0$ with the following property: If $\varepsilon>0$ is sufficiently small, then $u_\varepsilon$ is the unique element of $\clubsuit_{X_{E_\varepsilon}}(\mathcal{T}; p_\varepsilon; q_\varepsilon^1, q_\varepsilon^2, q_\varepsilon^3)$ that lies (entirely) in $V_{u_0}$. 

Similarly to the proof of Theorem \ref{thm_corresp_traj_gsp} we put
$$V_{u_0}^{int}:=B_{\delta_3}(y_0)\cup T_{\delta_3}\cup B_{\delta_3}(\phi^{\overline{t}}_{slow}(y_0)),$$
where $\delta_3>0$ such that $\delta_3\ll\min\lbrace\delta, \delta_1, \delta_2, \hbox{radii of }D_x\hbox{ and }\widehat{D}\rbrace$. I.e. we ``glued'' the balls $B_{\delta_3}(y_0), B_{\delta_3}(\phi^{\overline{t}}_{slow}(y_0))$ by a thin tube $T_{\delta_3}$ which is just a tubular neighborhood $\nu_{\delta_3}(\phi^{(0, \overline{t})}_{slow}(y_0))$. Recall also that $\phi^{\overline{t}}_{slow}(y_0)=\ioo(x)$.

Now, we claim that
\begin{itemize}
\item[$(\bigcirc)$] If there is a $z\in B_{\delta_3}(y_0)\setminus(\mathcal{R}^{\theta_1}_\varepsilon)^\ast$ and a (slow time) $t_1>0$ such that $\phi^{t_1}_\ffs(z)\in\mathcal{R}^{\theta_1}_\varepsilon\cap B_{\delta_3}(\phi^{\overline{t}}_{slow}(y_0))$, then $\phi^{[0, t_1]}_\ffs(z)$ does not lie entirely in $V^{int}_{u_0}.$
\end{itemize}
In order to show $(\bigcirc)$ we would like to use as in $(\square)$ (Theorem \ref{thm_corresp_traj_gsp}) the local invariance of $(\mathcal{R}^{\theta_1}_\varepsilon)^\ast$. Hence, let us inspect the topological boundary $\partial(\mathcal{R}^{\theta_1}_\varepsilon)^\ast$.

First, it contains $\mathcal{R}^{\theta_1}_\varepsilon$ which has clearly nonempty intersection with $V_{u_0}^{int}$. However, the slow flow is strictly outward-pointing from  $(\mathcal{R}^{\theta_1}_\varepsilon)^\ast$ along $\mathcal{R}^{\theta_1}_\varepsilon$, and thus no trajectory can enter $(\mathcal{R}^{\theta_1}_\varepsilon)^\ast$ along $\mathcal{R}^{\theta_1}_\varepsilon$.

Now we would like to show that the remaining components of the boundary $\partial(\mathcal{R}^{\theta_1}_\varepsilon)^\ast$ lie outside of the $V_{u_0}^{int}$, which will imply $(\bigcirc)$.

By Remarks \ref{rem_inclination1} under the backward-flow $(\mathcal{R}^{\theta_1}_\varepsilon)^\ast$ contains, as a part of the boundary, $(\mathcal{R}^{\theta_1}_\varepsilon)^\ast\cap(\Phi^{fen}_\varepsilon)^{-1}(\Vert b\Vert\leq\delta)$, but since $\delta_3\ll\delta$, we do need to worry about this. Recall also that $(\mathcal{R}^{\theta_1}_\varepsilon)^\ast\cap(\Phi^{fen}_\varepsilon)^{-1}(\Vert a\Vert\leq\delta)=\emptyset$ by Remark \ref{rem_monot_estim}.

The remaining boundary phenomena can be captured by the behavior of the slow variables. In particular, by Remark \ref{rem_inclination1} and the fact that $\delta_3\ll\delta_1$ it follows that $(\mathcal{R}^{\theta_1}_\varepsilon)^\ast\cap\phi^{-\overline{t}-\delta_1}_{\ffs}(\mathcal{R}^{\theta_1}_\varepsilon)$ is outside of $V^{int}_{u_0}$. Since $\delta_3$ is much smaller then the radius of $D_x$ we can treat the case $(\mathcal{R}^{\theta_1}_\varepsilon)^\ast\cap\phi^{(-\overline{t}-\delta_1, 0)}_{\ffs}(\partial\mathcal{R}^{\theta_1}_\varepsilon)$ analogously. And the claim $(\bigcirc)$ follows.

Now, similarly to Lemma \ref{lem_tree_2} the discs $\widehat{D}, D_x$ induce small neighborhoods $V^{\mathfrak{p}_0}_{u_0}, V^{\mathfrak{q}_0^1}_{u_0}, V^{\mathfrak{q}_0^2}_{u_0}, V^{\mathfrak{q}_0^3}_{u_0}$. Finally, their union together with $V_{u_0}^{int}$ induce $V_{u_0}$.

\bigskip

\textit{Trap for trees:}

In particular, we would like to construct a subset $N^\varepsilon_{u_0}\subset V_{u_0}$ with the following property: if $\varepsilon>0$ is small and $\widehat{u}_\varepsilon\in\clubsuit_{X_{E_\varepsilon}}(\mathcal{T}; p_\varepsilon; q_\varepsilon^1, q_\varepsilon^2, q_\varepsilon^3)$ such that $\pi_{s_1, s_2}(\widehat{u}_\varepsilon)\subset\pi_{s_1, s_2}(N^\varepsilon_{u_0})$, then $\widehat{u}_\varepsilon\subset N^\varepsilon_{u_0}$.

This is done by the same trick as in Lemma \ref{lem_tree_2}, one only has ensure that the conditions $\vert\langle P, v^\bot_i\rangle\vert<\varepsilon^j$ are still satisfied when translated over the interior vertices. For this, one can use suitable powers of $\varepsilon$, provided that $\varepsilon>0$ is sufficiently small.

For example, we can define $N^{\varepsilon}_{u_0}$ as the union of the following five sets.
\begin{align*}
N^{\varepsilon, \mathfrak{p}_0}_{u_0}&=\big\lbrace (s_1, \theta_1, s_2, \theta_2)\in M_{K, \varepsilon}\,\vert\,(s_1, s_2)\in\pi_{s_1, s_2}(V^{\mathfrak{p}_0}_{u_0})\wedge\vert\langle P, v_1^\bot\rangle\vert<\varepsilon^{1/2}\wedge\vert\langle P, v_2^\bot\rangle\vert<\varepsilon^{1/3}\big\rbrace,\\
N^{\varepsilon, \mathfrak{q}_0^1}_{u_0}&=\big\lbrace (s_1, \theta_1, s_2, \theta_2)\in M_{K, \varepsilon}\,\vert\,(s_1, s_2)\in\pi_{s_1, s_2}(V^{\mathfrak{q}_0^1}_{u_0})\wedge\vert\langle P, v_1^\bot\rangle\vert<\varepsilon^{1/4}\wedge\vert\langle P, v_2^\bot\rangle\vert<\varepsilon^{1/2}\big\rbrace,\\
N^{\varepsilon, \mathfrak{q}_0^2}_{u_0}&=\big\lbrace (s_1, \theta_1, s_2, \theta_2)\in M_{K, \varepsilon}\,\vert\,(s_1, s_2)\in\pi_{s_1, s_2}(V^{\mathfrak{q}_0^2}_{u_0})\wedge\vert\langle P, v_1^\bot\rangle\vert<\varepsilon^{1/3}\wedge\vert\langle P, v_2^\bot\rangle\vert<\varepsilon^{1/2}\big\rbrace,\\
N^{\varepsilon, \mathfrak{q}_0^3}_{u_0}&=\big\lbrace (s_1, \theta_1, s_2, \theta_2)\in M_{K, \varepsilon}\,\vert\,(s_1, s_2)\in\pi_{s_1, s_2}(V^{\mathfrak{q}_0^3}_{u_0})\wedge\vert\langle P, v_1^\bot\rangle\vert<\varepsilon^{1/4}\wedge\vert\langle P, v_2^\bot\rangle\vert<\varepsilon^{1/2}\big\rbrace,\\
N^{\varepsilon, int}_{u_0}&=\big\lbrace (s_1, \theta_1, s_2, \theta_2)\in M_{K, \varepsilon}\,\vert\,(s_1, s_2)\in\pi_{s_1, s_2}(V^{int}_{u_0})\wedge\vert\langle P, v_1^\bot\rangle\vert<\varepsilon^{1/3}\wedge\vert\langle P, v_2^\bot\rangle\vert<\varepsilon^{1/3}\big\rbrace.
\end{align*}

\bigskip

\textit{Surjectivity of $\Psi^\varepsilon_{\mathcal{T}; p_0; q_0^1, q_0^2, q_0^3}$:}
Now, it follows from the same argument as in Lemma \ref{lem_tree_2}.
\end{proof}

Now, we are going to briefly discuss some \textbf{orientation conventions}.

\begin{rem}\label{rem_sign_conv}Let us consider the canonical orientations on $(\R/T\mathbb{Z})^2$ and $(\R/T\mathbb{Z}\times S^1)^2$ which are given by coordinates $(s_1, s_2)$ and $(s_1, \theta_1, s_2, \theta_2).$ Then on $S^{out-out}$ we induce an orientation by the orientation on $\overline{S_K}\subset(\R/T\mathbb{Z})^2$ such that the map $\ioo$ is orientation preserving.

Now, at each point $p_0\in Crit(E_0)\setminus\Delta_0$ we choose an arbitrary orientation of $W_{X_{E_0}}^u(p_0)$. Then by Lemma \ref{lemma_identif} the canonical isomorphism $T_{p_0}W_{X_{E_0}}^u(p_0)\cong T_{\mathfrak{p}_0} W^u_{slow}(\mathfrak{p}_0)$ induces an orientation on $W^u_{slow}(\mathfrak{p}_0)$ (where $\mathfrak{p}_0\subset S^{out-out}$). Moreover, by Fenichel theory, we have an isomorphism
\begin{equation}\label{eqn_or_unstable}
T_{p_\varepsilon} W^u_{\ffs}(\mathfrak{p}_\varepsilon)\cong T_{\mathfrak{p}_0} W^u_{slow}(\mathfrak{p}_0)\oplus\langle \partial_{\theta_2}\rangle,
\end{equation}
which induces an orientation of $W^u_{\ffs}(\mathfrak{p}_\varepsilon)$. 

Next, the orientation on $W^s_{X_{E_0}}(p_0)$ is given by 
$$T_{p_0}W^s_{X_{E_0}}(p_0)\oplus T_{p_0}W^u_{X_{E_0}}(p_0)=T_{p_0}(\R/T\mathbb{Z})^2,$$
which gives us also an orientation of $W^s_{slow}(\mathfrak{p}_0)$.
Analogously, by Fenichel theory, we have an isomorphism
\begin{equation}\label{eqn_or_stable}
T_{p_\varepsilon} W^s_{\ffs}(\mathfrak{p}_\varepsilon)\cong T_{\mathfrak{p}_0} W^s_{slow}(\mathfrak{p}_0)\oplus\langle \partial_{\theta_1}\rangle,
\end{equation}
which induce an orientation of $W^s_{\ffs}(\mathfrak{p}_\varepsilon)$.

Note that the orientations of stable and unstable manifolds at $\mathfrak{p}_0$ and $\mathfrak{p_\varepsilon}$ match in the following sense: there are orientation-preserving isomorphisms
$$T_{p_0}W^s_{X_{E_0}}(p_0)\oplus\langle\partial_{\theta_1}\rangle\oplus T_{p_0}W^u_{X_{E_0}}(p_0)\oplus\langle\partial_{\theta_2}\rangle\cong T_{\mathfrak{p}_\varepsilon}(\R/T\mathbb{Z}\times S^1)^2$$
and
$$T_{p_\varepsilon} W^s_{\ffs}(\mathfrak{p}_\varepsilon)\oplus T_{p_\varepsilon} W^u_{\ffs}(\mathfrak{p}_\varepsilon)= T_{\mathfrak{p}_\varepsilon}(\R/T\mathbb{Z}\times S^1)^2.$$
\end{rem}

\begin{rem}\label{rem_orient_con_intersect}Let $M$ be an oriented manifold with oriented submanifolds $Q_1, Q_2$. Let moreover $N:=Q_1\pitchfork Q_2$. Then we orient $N$ by the convention that at some $x \in N$ the ordered set
\begin{equation}\label{eqn_orient_con_intersect}
\lbrace e_1,\dots, e_i, e_{i+1},\dots, e_j, e_{j+1},\dots, e_n\rbrace
\end{equation}
is an oriented basis of $T_x M$. Here $\lbrace e_1,\dots, e_j\rbrace$ is an oriented basis of $T_x Q_1$, $\lbrace e_{i+1},\dots, e_n\rbrace$ is an oriented basis of $T_x Q_2$. And finally $\lbrace e_{i+1},\dots, e_j\rbrace$ is a basis of $T_x N$ with the determined orientation. 
\end{rem}

\begin{rem}\label{rem_sign_or} Let $p_0\in Crit_1(E_0)\setminus \Delta_0, q_0\in Crit_0(E_0)\setminus \Delta_0$ and $p_\varepsilon, q_\varepsilon$ will be the corresponding $Crit(E_\varepsilon)\cap M_{K, \varepsilon}$ of indices $2$ and $1$, respectively.

Let $u\in \mathcal{M}_{X_{E_0}}(p_0, q_0)$ and $u_\varepsilon:=\Psi^\varepsilon_{p_0, q_0}(u)\in \mathcal{M}_{X_{E_\varepsilon}}(p_\varepsilon, q_\varepsilon)$, see Theorem \ref{thm_corresp_traj_gsp}. Our goal is to compute and compare the \textbf{orientation signs} $\hbox{sgn}(u)$ and $\hbox{sgn}(u_\varepsilon).$

Let $x\in u$. Then $x_0$ and $x_\varepsilon$ will denote the corresponding points of the slow trajectory $u_0$ and $u_\varepsilon$. Here, the correspondence of $x_0$ and $x_\varepsilon$ is given by Fenichel theory, see Theorem \ref{thm_corresp_traj_gsp}.

By index reasons $\dim(W^s_{X_{E_0}}(q_0))=2$ and $\dim(W^u_{X_{E_0}}(p_0))=1$. We would like to give an orientation to the manifold $u\subset W^s_{X_{E_0}}(q_0)\cap W^u_{X_{E_0}}(p_0)$ by the convention from Remark \ref{rem_orient_con_intersect}. For the oriented basis $\lbrace e_b\rbrace$ of $T_x W^u_{X_{E_0}}(p_0)$ we can find $e_a\in T_x (\R/T\mathbb{Z})^2$ such that
$$\lbrace e_a, e_b\rbrace$$
is an oriented basis of $T_x (\R/T\mathbb{Z})^2$. But then, by the dimension reason, this is also an oriented basis of $T_x W^s_{X_{E_0}}(q_0)$. Hence, $e_b$ induces an orientation of $u$ (as a manifold).

Next, by index reasons $\dim(W^s_{\ffs}(\mathfrak{q}_\varepsilon))=3$ and $\dim(W^u_{\ffs}(\mathfrak{p}_\varepsilon))=2$. Let us give an orientation to the manifold $u_\varepsilon\subset W^s_{\ffs}(\mathfrak{q}_\varepsilon)\cap W^u_{\ffs}(\mathfrak{p}_\varepsilon).$

Now, by $(\ref{eqn_or_unstable})$ and $(\ref{eqn_or_stable})$, there is an orientation \textit{reversing} isomorphism identifying the basis
$$\lbrace\partial_{\theta_1}, e_a, e_b, \partial_{\theta_2}\rbrace$$
with an oriented basis of $T_{x_\varepsilon}(\R/T\mathbb{Z}\times S^1)^2$.  Also $\lbrace\partial_{\theta_1}, e_a, e_b\rbrace$ can be identified by an orientation \textit{preserving} isomorphism with an oriented basis of $T_{x_\varepsilon} W^s_{\ffs}(\mathfrak{q}_\varepsilon)$. And $\lbrace e_b, \partial_{\theta_2}\rbrace$ can be identified by an orientation \textit{preserving} isomorphism with an oriented basis of $T_{x_\varepsilon} W^u_{\ffs}(\mathfrak{p}_\varepsilon)$. Hence, $-e_b$ induces an orientation on $u_\varepsilon$ as a manifold.

Finally, by $C^1$ estimates from Fenichel theory, at $x$ the orientations of $u$ and $X_{E_0}\vert_u$ agree iff at $x_\varepsilon$ the orientations of $u_\varepsilon$ and $X_{E_\varepsilon}\vert_{u_\varepsilon}$ does not agree.

Hence,
\begin{equation}\label{eqn_sqn_corresp}
\hbox{sgn}(u)=-\hbox{sgn}(u_\varepsilon).
\end{equation}
To the end, we refer the reader for the orientation conventions of trees to \cite{mescher2016perturbedgradientflowtrees}. We claim that these orientation conventions can be adapted to our setting and result in the same correspondence as in (\ref{eqn_sqn_corresp}).  
\end{rem}

\chapter{Cord algebras}
\label{ch:cord}
This chapter aim is to introduce Morse models of cord algebras for $T_{K, \varepsilon}$ and relate them to the existing cord algebras for the underlying knot $K$. 
As a consequence of the multiple times scale dynamics, we will show directly the equivalence of cord algebra for $K$ and $T_{K, \varepsilon}$ with $\mathbb{Z}$ coefficients. Moreover, we upgrade the cord algebra for $T_{K, \varepsilon}$ by loop space coefficients $\mathbb{Z}[\lambda^{\pm1}, \mu^{\pm 1}]$ ($\cong\mathbb{Z}\pi_{1}(T_K)\cong H_\bullet(\Omega_\ast T_K)$). We recall that in the knot case, the loop space coefficients were important and the corresponding cord algebra contained much more information about the underlying knot, \cite{Cieliebak2016KnotCH}. In the torus case, the cord algebra with loop space coefficients will be greatly influenced by the geometry of $M_{K, \varepsilon}$ - the space of outward-pointing chords. In particular, on the torus, we will be able to detect an analogy of a bit mysterious skein relation, which was on the knot, replacing short contractible cords by $1-\mu$.

\section{Topological cord algebra for knots}
\begin{assump}From now on, $K$ will be an oriented framed knot in $\R^3$ and the framing $\mathfrak{f}_K$ will induce the choice of the longitude $\mathfrak{L}\subset T_K$. We will choose the meridian $\mathfrak{m}\subset T_K$ such that $\mathfrak{m}$ projects orthogonally onto $\ast\in K$.
\end{assump}

\begin{defn} A \textbf{topological cord} is a continuous map $c:[0, 1]\rightarrow\R^3$ such that $c([0, 1])\cap \mathfrak{L}=\emptyset$ and $c(0), c(1)\in K\setminus \ast$. Two topological cords are \textbf{homotopic} if they are homotopic through topological cords.
\end{defn}

\begin{rem} A linear topological cord is a chord.
\end{rem}

\begin{defn}\label{defn_cord_top}\cite{Cieliebak2016KnotCH} Let $\mathcal{A}^{top}$ be a free unital noncommutative $\mathbb{Z}$-algebra generated by homotopy classes of topological chords and extra generators $\lambda, \mu$ such that $\lambda$ and $\mu$ have inverses and commute together.

Then the \textbf{topological cord algebra of $K$ with loop space coefficients} is the quotient ring
$$Cord^{top}(K; \mathbb{Z}[\lambda^{\pm1}, \mu^{\pm1}])=\mathcal{A}^{top}/\mathcal{I}^{top},$$
where $\mathcal{I}^{top}$ is a two-sided ideal of $\mathcal{A}^{top}$ generated by the skein relation as in Figure \ref{figure_skein} 
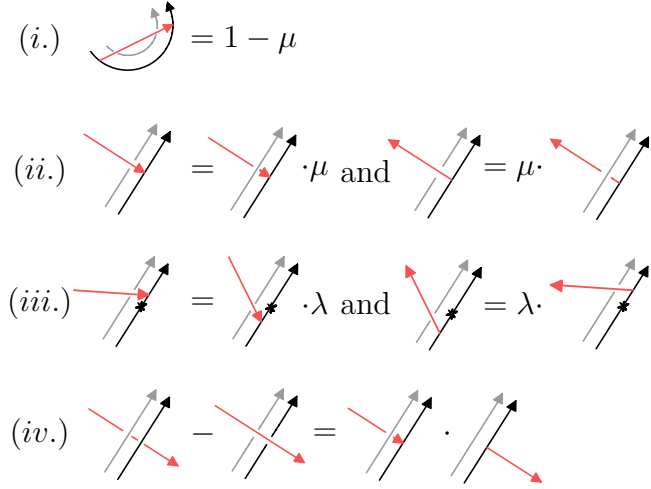
\begin{figure}[!htbp]
\labellist
\pinlabel $=1-\mu$ at 170 800
\pinlabel $=$ at 130 670
\pinlabel $=$ at 130 540
\pinlabel $-$ at 130 410
\pinlabel $\cdot\mu$ at 240 670
\pinlabel $\cdot\lambda$ at 240 540
\pinlabel $=$ at 250 410
\pinlabel $\hbox{and}$ at 290 670
\pinlabel $\hbox{and}$ at 290 540
\pinlabel $=\mu\cdot$ at 435 670
\pinlabel $=\lambda\cdot$ at 435 540
\pinlabel $\cdot$ at 370 410
\pinlabel $(i.)$ at -30 800
\pinlabel $(ii.)$ at -30 670
\pinlabel $(iii.)$ at -30 540
\pinlabel $(iv.)$ at -30 410
\endlabellist
\centering
\includegraphics[scale=0.39]{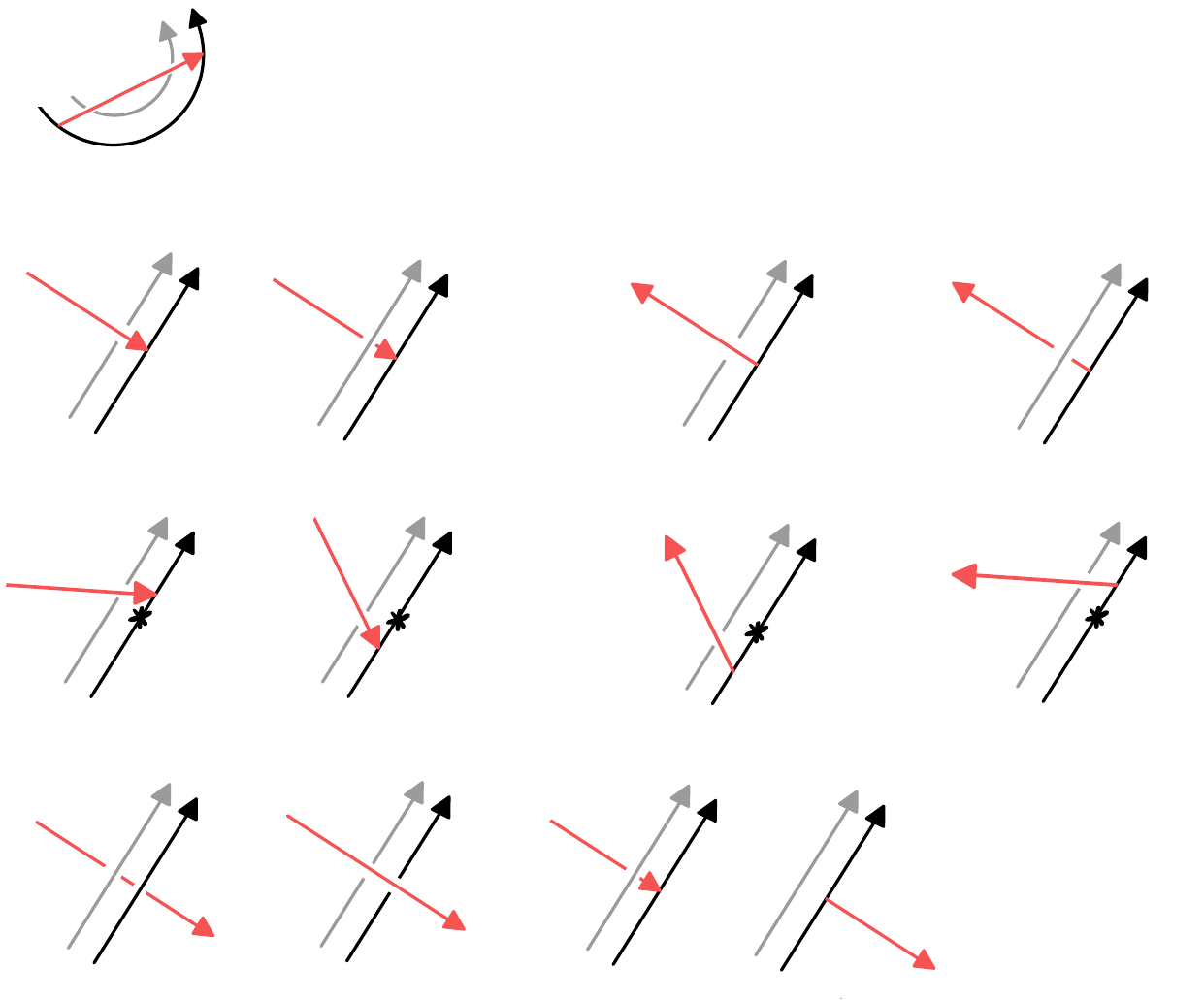}
\vspace{0.3cm}
\caption{The skein relations for $\mathcal{I}^{top}$. The black curves describe the knot $K$, grey the longitude $\mathfrak{L}$, and in red are the cords.
Skein relations describe the following phenomena: $(i.)$ is for contractible cords, $(ii.)$ and $(iii.)$ are for cords crossing the framing (i.e. $\mathfrak{L}$) and the base point, respectively. Finally, $(iv.)$ relates cords as they cross the knot.}
\label{figure_skein}
\end{figure}
\end{defn}

\begin{rem} $Cord^{top}(K; \mathbb{Z}[\lambda^{\pm1}, \mu^{\pm1}])$ does not depend on the choice of $\ast$ and is an invariant under ambient isotopy of $K$.

Using so-called string homology it was shown in the paper \cite{Cieliebak2016KnotCH} that $Cord^{top}(K; \mathbb{Z}[\lambda^{\pm1}, \mu^{\pm1}])$ is isomorphic to the $0$-th degree Legendrian contact homology $LCH_0(\mathcal{L}^\ast K).$

In addition, in the same paper, it was proven that $Cord^{top}(K; \mathbb{Z}[\lambda^{\pm1}, \mu^{\pm1}])$ detects the unknot.
\end{rem}

\begin{rem} If we put $\lambda^{\pm 1}=\mu^{\pm 1}=1$, then $Cord^{top}(K; \mathbb{Z}[\lambda^{\pm1}, \mu^{\pm1}])$ becomes a \textbf{Cord algebra of $K$ over $\mathbb{Z}$}, denoted by $Cord^{top}(K; \mathbb{Z}).$

Note that $Cord^{top}(K; \mathbb{Z})$ is independent of the framing. In more detail, the cords, which are homotopic in $\pi_1(\R^3, K)$, are identified in the cord algebra. Specially, the cords from $[0]\in\pi_1(\R^3, K)$ vanish in the cord algebra. 

Since $\lambda=1$, by the skein relation $(iii.)$ from Figure \ref{figure_skein}, the base point becomes redundant too.

For various versions of $Cord^{top}(K)$ see also \cite{ng2007framedknotcontacthomology, Ng_2005, okamoto2022toward, okamoto2024legendrian}.
\end{rem}

\section{Morse models for Cord algebras over $\mathbb{Z}$}

\begin{defn}\cite{petrak2019definition, Ng_2005, okamoto2024legendrian}\label{def_cordK_Z}Let $K$ be generic and $X_{E_0}$ be a generic approximation of $-\nabla E_0$.

Let $C_0^M(K; \mathbb{Z})$ be a free unital noncommutative $\mathbb{Z}$-algebra generated by $Crit_0(E_0)\setminus \Delta_0$. Let $C_1^M(K; \mathbb{Z})$ be a $\mathbb{Z}$-vector space generated by $Crit_1(E_0)\setminus\Delta_0$. Let $D_K:C_1^M(K; \mathbb{Z})\rightarrow C_0^M(K; \mathbb{Z})$ be a linear map which is defined on generators as
$$D_K(p_0)=\sum_{m>0}\sum_{\substack{\mathcal{T}\in\clubsuit_m\\ q_0^1,\dots,q_0^m\in Crit_0(E_0)\setminus\Delta_0}}\sum_{u\in\clubsuit_{X_{E_0}}(\mathcal{T}; p_0; q_0^1,\dots, q_0^m)}\hbox{sqn}(u)q_0^1\dots q_0^m.$$

Then the \textbf{Morse cord algebra of $K$ over $\mathbb{Z}$} is the quotient ring
$$Cord^M(K; \mathbb{Z})=C_0^M(K; \mathbb{Z})/\mathcal{I}_{K},$$
where $\mathcal{I}_{K}$ is a two-sided ideal of $C_0^M(K; \mathbb{Z})$ generated by the image of the map $D_K$.
\end{defn}

\begin{rem} By Lemma \ref{lemma_dim_tree_knot} it follows that $D_K(p_0)$ is a finite sum.
\end{rem}

\begin{thm}\cite{petrak2019definition} Let $K_0$ and $K_1$ be connected by a generic ambient isotopy. Let $X_{E_0}^0$ and $X_{E_0}^1$ be connected by a generic smooth family of gradient-like vector fields $\lbrace X_{E_0}^r\rbrace_{r\in[0, 1]}$ approximating $\lbrace-\nabla E_0^r\rbrace_{r\in[0, 1]}$. By ``generic'' we understand that only for a finite number of $r\in(0, 1)$ the elements $K_r$ and $X_{E_0}^r$ are not generic in the sense of Definition \ref{def_cordK_Z}. Then
$$Cord^M(K_0; \mathbb{Z})\cong Cord^M(K_1; \mathbb{Z}).$$
\end{thm}

\begin{defn}Let $\varepsilon>0$ be small. Let $K$ be generic and $X_{E_\varepsilon}$ be a $K$-generic approximation of $-\nabla E_\varepsilon$.

Let $C_0^M(T_{K, \varepsilon}; \mathbb{Z})$ be a free unital noncommutative $\mathbb{Z}$-algebra which is generated by $Crit_1(E_\varepsilon\vert_{M_{K, \varepsilon}\setminus \Delta_{\varepsilon}})$. Let $C_1^M(T_{K, \varepsilon}; \mathbb{Z})$ be a $\mathbb{Z}$-vector space generated by $Crit_2(E_\varepsilon\vert_{M_{K, \varepsilon}\setminus\Delta_{\varepsilon}})$. Let $D_{T_{K, \varepsilon}}:C_1^M(T_{K, \varepsilon}; \mathbb{Z})\rightarrow C_0^M(T_{K, \varepsilon}; \mathbb{Z})$ be a linear map which is defined on generators as
$$D_{T_{K, \varepsilon}}(p_\varepsilon)=\sum_{m>0}\sum_{\substack{\mathcal{T}\in\clubsuit_m\\ q_\varepsilon^1,\dots,q_\varepsilon^m\in Crit_0(E_\varepsilon)\setminus\Delta_\varepsilon}}\sum_{u_\varepsilon\in\clubsuit_{X_{E_\varepsilon}}^{out-out}(\mathcal{T}; p_\varepsilon; q_\varepsilon^1,\dots, q_\varepsilon^m)}\hbox{sqn}(u_\varepsilon)q_\varepsilon^1\dots q_\varepsilon^m.$$

Then the \textbf{Morse cord algebra of $T_{K, \varepsilon}$ over $\mathbb{Z}$} is the quotient ring
$$Cord^M(T_{K, \varepsilon}; \mathbb{Z})=C_0^M(T_{K, \varepsilon}; \mathbb{Z})/\mathcal{I}_{T_{K, \varepsilon}},$$
where $\mathcal{I}_{T_{K, \varepsilon}}$ is a two-sided ideal of $C_0^M(T_{K, \varepsilon}; \mathbb{Z})$ generated by the image of the map $D_{T_{K, \varepsilon}}$.
\end{defn}

\begin{thm}Let $\varepsilon>0$ be small and $K$ generic. Let $X_{E_0}$ and $X_{E_\varepsilon}$ be $K$-generic approximations of $-\nabla E_0$ and $-\nabla
E_\varepsilon$, respectively. 

Then for $i=0, 1$ the canonical isomorphisms $\Theta_i^\varepsilon: C_i(K; \mathbb{Z})\rightarrow C_i(T_{K, \varepsilon}; \mathbb{Z})$ satisfy
$$\Theta_0^\varepsilon\circ D_K=D_{T_{K, \varepsilon}}\circ \Theta_1^\varepsilon.$$

In particular, 
$$Cord^M(K; \mathbb{Z})\cong Cord^M(T_{K, \varepsilon}; \mathbb{Z}).$$
\end{thm}

\begin{proof}
By Lemma \ref{lemma_morse_corresp} we obtain the cannonical isomorphisms $\Theta_0^\varepsilon, \Theta_1^\varepsilon$. Then the theorem follows from Theorems \ref{thm_corresp_traj_gsp} and \ref{thm_corresp_trees_gsp}. See also Remarks \ref{rem_sign_conv} and \ref{rem_sign_or} for the sign conventions.
\end{proof}

\begin{cor}Let $\varepsilon>0$ be small. Let $K_0$ and $K_1$ be connected by a generic ambient isotopy. Let $X_{E_\varepsilon}^0$ and $X_{E_\varepsilon}^1$ be connected by a generic smooth family of gradient-like vector fields $\lbrace X_{E_\varepsilon}^r\rbrace_{r\in[0, 1]}$ $K$-approximating $\lbrace-\nabla E_\varepsilon^r\rbrace_{r\in[0, 1]}$. Then
$$Cord^M(T_{K_0, \varepsilon}; \mathbb{Z})\cong Cord^M(T_{K_1, \varepsilon}; \mathbb{Z}).$$
\end{cor}

\section{Morse models for Cord algebras with loop space coefficients}
\subsection{Knot case}
\begin{rem_not}\label{rem_not_extended_trees}There is a canonical extension of Definition \ref{defn_ribbon_tree_K} of  the trees $\clubsuit_{X_{E_0}}(\mathcal{T}; p_0; q_0^1,\dots, q_0^m)$. We would like to allow the leaves be identified not only with the points in $Crit_0(E_0)\setminus \Delta_0$, but also with the Bott-critical submanifold $\Delta_0$. Let us denote such a space of trees by 
$$\clubsuit_{X_{E_0}}(\mathcal{T}; p_0; Q_0^1,\dots, Q_0^m).$$

By the argument from \cite[prop 7.14]{Cieliebak2016KnotCH}, one can still bound the number of potential bifurcations of the trees from above. Moreover, $\Delta_0$ is the unique global minimum of $E_0$, and in particular $\dim(W^s_{X_{E_0}}(\Delta_0))=2$. So, a small adjustment of the argument from Lemma \ref{lemma_dim_tree_knot} implies that for a generic approximation $X_{E_0}$ of $-\nabla E_0$ it holds that $\clubsuit_{X_{E_0}}(\mathcal{T}; p_0; Q_0^1,\dots, Q_0^m)$ is a $0$-dimensional compact manifold. In addition, the set of all of these trees $\clubsuit_{X_{E_0}}$ is finite.
\end{rem_not}

\begin{rem_not}\label{rem_morse_diag_K} Let us consider $\Delta_0$ with the induced metric from the flat metric on $(\R/T\mathbb{Z})^2$. By $$h_{\Delta_0}$$ we denote a Morse function on $\Delta_0$ such that $-\nabla h_{\Delta_0}$ is Morse-Smale.

It follows that $HM_0(h_{\Delta_0}; \mathbb{Z})\cong\mathbb{Z}$ and all elements of $Crit_0(h_{\Delta_0})$ are homologous, see for example \cite[rem 4.5.4]{Audin2013MorseTA}.

Thus for simplicity, we shall assume that $Crit(h_{\Delta_0})$ consists of two critical points - $M_{\Delta_0}, m_{\Delta_0}$ of indices $1, 0$, respectively.

Moreover, since the set $\clubsuit_{X_{E_0}}$ is finite, we impose a generic condition that $\overline{u}\cap Crit(h_{\Delta_0})=\emptyset$ for every $u\in\clubsuit_{X_{E_0}}$. Last, we impose the condition that $\ast\cap Crit(h_{\Delta_0})=\emptyset$.

This will give us a good choice of $h_{\Delta_0}$ for the cord algebra.
\end{rem_not}

\begin{defn} Let $p_0\in Crit_1(E_0)\setminus \Delta_0$ and $q_0^1, \dots, q_0^m\in (Crit_0(E_0)\setminus \Delta_0)\cup \lbrace m_{\Delta_0}\rbrace$, by $\lbrace j_\ell\rbrace_{\ell=1,\dots, k}$ we denote those indices from $\lbrace 1,\dots, m\rbrace$ such that $q_0^{j_\ell}=m_{\Delta_0}$. Let $\mathcal{T}\in\clubsuit_m.$ Let $X_{E_0}$ be a gradient-like vector field adapted to $E_0$. Then $u$ is a \textbf{Cascade Morse flow tree from $p_0$ to $q_0^1, \dots, q_0^m$ modeled on $\mathcal{T}$} if the following holds.

$u$ is a concatenation of a tree $u_{c, 0}$ and partial flow lines $\lbrace u^{j_\ell}\rbrace$, where
\begin{itemize}
\item $u_{c, 0}\in\clubsuit_{X_{E_0}}(\mathcal{T}; p_0; Q_0^1,\dots, Q_0^m)$ with $Q_0^i=q_0^i$ if $i\neq j_\ell$, else $Q_0^i=\Delta_0$.
\item Let $\overline{u_{c, 0}}\cap Q_0^{j_\ell}= x^{j_\ell}$. Then $u^{j_\ell}$ is given by $\phi^{[0, \infty)}_{-\nabla h_{\Delta_0}}(x^{j_\ell})$. In particular, the omega limit of $x^{j_\ell}$ is $m_{\Delta_0}$.
\end{itemize}
The convention is that the pairs of concatenated paths are identified with the single edge of $\mathcal{T}$, which will be called \textbf{broken} with the \textbf{cascade} and the \textbf{diagonal} part. So the leaves are identified with $q_0^1,\dots, q_0^m$.

The \textbf{set of all Cascade Morse flow trees from $p_0$ to $q_0^1, \dots, q_0^m$ modeled on $\mathcal{T}$} is denoted by 
$$\clubsuit_{X_{E_0}}^c(\mathcal{T}; p_0; q_0^1, \dots, q_0^m).$$
The set of all such trees is denoted by $\clubsuit_{X_{E_0}}^c$.
\end{defn}

\begin{defn}\label{defn_intr_fram_pt}\cite{petrak2019definition} We say that the chord $P=\gamma(s_2)-\gamma(s_1)$ \textbf{intersects framing at the starting point}, if the orthogonal projection of $P$ to the normal plane at $\gamma(s_1)$ is a positive multiple of $\mathfrak{f}_K(s_1)$. 
That is, if it holds
$$P-\langle P, \dot{\gamma}(s_1)\rangle \dot{\gamma}(s_1)=\alpha \mathfrak{f}_K(s_1),$$
for some $\alpha>0$. The \textbf{set of all chords intersecting framing at the starting point} is denoted by $\mathcal{F}^S$. Analogously, we define the \textbf{intersection of the framing at the endpoint} and the set $\mathcal{F}^E$.
\end{defn}

\begin{rem}\cite{petrak2019definition} For generic $K$ and generic framing $\mathfrak{f}_K$ the following holds. The closures $\overline{\mathcal{F}^{S}}, \overline{\mathcal{F}^{E}}$ are embedded one dimensional submanifolds of $(\R/T\mathbb{Z})^2$ with boundary of $s_1$- and $s_2$-special points, respectively. Moreover, $\overline{\mathcal{F}^{S}}$ and $\overline{\mathcal{F}^{E}}$ intersect $\Delta_0$ transversely \textit{at the same points}.
\end{rem}

\begin{rem}\label{rem_restr_skein} Definition \ref{defn_intr_fram_pt} allows us to restrict the interpretation of the skein relation $(ii.)$ from Definition \ref{defn_cord_top} to chords that differ by the flow of $X_{E_0}$.

The restriction of the skein relation $(iii.)$ to chords flowing under the flow is clear. 

Recall that $\overline{\mathcal{F}^{S}}, \overline{\mathcal{F}^{E}}$ are well-defined also on $\Delta_0$. So, in fact, we can also interpret the skein relation $(ii.)$ for the trivial chords flowing under the flow of $-\nabla h_{\Delta_0}$. The interpretation of $(iii.)$ for the trivial chords is clear.
\end{rem}

Now we can, in the spirit of \cite{petrak2019definition}, introduce a certain evaluation function $\widehat{D}_u$ on a tree $u\in\clubsuit^c_{X_{E_0}}$. Such a function will just count how many times the tree $u$ intersects the framing and the base point as we flow from the root to the leaves.

\begin{defn}\label{defn_cord_K_dif} Let $u\in\clubsuit_{X_{E_0}}^c(\mathcal{T}; p_0; q_0^1, \dots, q_0^m)$. Recall also the notation for the trees from Remark \ref{rem_graphs}. Then the evaluation function $\widehat{D}_u$ is recursively defined on the vertices of $\mathcal{T}$ as
\begin{equation*}
\widehat{D}_u(v)=
    \begin{cases}
      \mu^{\alpha_1(\emph{e}_r)}\lambda^{\beta_1(\emph{e}_r)}\widehat{D}_u(v^{out}_{\emph{e}_r})\mu^{\alpha_2(\emph{e}_r)}\lambda^{\beta_2(\emph{e}_r)}, & \text{if }v\text{ is the root},\\
      \\
      \mu^{\alpha_1(\emph{e}^L_v)}\lambda^{\beta_1(\emph{e}^L_v)}\widehat{D}_u(v^{out}_{\emph{e}^L_v})\mu^{\alpha_2(\emph{e}^L_v)+\alpha_1(\emph{e}^U_v)} & \text{if }v\text{ is an interior vertex},\\
      \quad\lambda^{\beta_2(\emph{e}^L_v)+\beta_1(\emph{e}^U_v)}\widehat{D}_u(v^{out}_{\emph{e}^U_v})\mu^{\alpha_2(\emph{e}^U_v)}\lambda^{\beta_2(\emph{e}^U_v)}, &\\
      \\
 q_0^j, & \text{if }v\text{ is a leaf represented}\\
        & \text{by the chord }q_0^j,
    \end{cases}
\end{equation*}
where $\alpha_k, \beta_k$ are $\mathbb{Z}$-valued functions on the corresponding edges of the tree. The functions are determined by the skein relations $(ii.), (iii.)$ from Definition \ref{defn_cord_top} and Remark \ref{rem_restr_skein} as we flow along $X_{E_0}$ with the chords in $\R^3$. Here, $k=1$ corresponds to the intersections at the starting points and $k=2$ at the endpoints.
\end{defn}

\begin{defn}\cite{petrak2019definition}\label{defn_cord_morse_loop}Let $K$ be a generic knot with generic framing $\mathfrak{f}_K$ and base point $\ast$. Let $X_{E_0}$ be a generic approximation of $-\nabla E_0$ and let $h_{\Delta_0}$ be generic.

Let $C_0^M(K; \mathbb{Z}[\lambda^{\pm 1}, \mu^{\pm 1}])$ be a free unital noncommutative $\mathbb{Z}$-algebra generated by $(Crit_0(E_0)\setminus \Delta_0)\cup\lbrace m_{\Delta_0}\rbrace$ and extra generators $\lambda, \mu$ such that $\lambda$ and $\mu$ have inverses and commute together. Let $C_1^M(K; \mathbb{Z}[\lambda^{\pm 1}, \mu^{\pm 1}])$ be a $\mathbb{Z}$-vector space generated by $Crit_1(E_0)\setminus\Delta_0$. Let $D_K:C_1^M(K; \mathbb{Z}[\lambda^{\pm 1}, \mu^{\pm 1}])\rightarrow C_0^M(K; \mathbb{Z}[\lambda^{\pm 1}, \mu^{\pm 1}])$ be a linear map which is defined on generators as
\begin{equation}\label{eqn_dif_cordK}
D_K(p_0)=\sum_{u\in \clubsuit_{X_{E_0}}^c}\hbox{sqn}(u)\widehat{D}_u(v^{p_0}).
\end{equation}

Then the \textbf{Morse cord algebra of $K$ with loop space coefficients} is the quotient ring
$$Cord^M(K; \mathbb{Z}[\lambda^{\pm 1}, \mu^{\pm 1}])=C_0^M(K; \mathbb{Z}[\lambda^{\pm 1}, \mu^{\pm 1}])/\mathcal{I}_{K},$$
where $\mathcal{I}_{K}$ is a two-sided ideal of $C_0^M(K; \mathbb{Z}[\lambda^{\pm 1}, \mu^{\pm 1}])$ generated by the image of the map $D_K$ and $\textcolor{blue}{m_{\Delta_0}-1+\mu}$.
\end{defn}

\begin{rem}\label{rem_gen_MK}One can also canonically extend the definition (\ref{eqn_dif_cordK}) of $D_K$ to $M_{\Delta_0}$. I. e. we would like to consider also the flow lines of $-\nabla h_{\Delta_0}$ from $M_{\Delta_0}$ to $m_{\Delta_0}$ and count their intersections with $\mathcal{F}^{S/E}$ and the base point $\ast$. There are two such flow lines with opposite orientations. In addition, we will obtain
$$D_K(M_{\Delta_0})=\mu^\alpha\lambda^\beta m_{\Delta_0}\mu^{-\alpha}\lambda^{-\beta}-\mu^{\widetilde{\alpha}}\lambda^{\widetilde{\beta}} m_{\Delta_0}\mu^{-\widetilde{\alpha}}\lambda^{-\widetilde{\beta}},$$
for some $\alpha, \beta, \widetilde{\alpha}, \widetilde{\beta}\in \mathbb{Z}$, since the intersections are coming in pairs on $\Delta_0$. Moreover, $(m_{\Delta_0}-1+\mu)\sim 0$, so $D_K(M_{\Delta_0})=0$. In particular, the generator $M_{\Delta_0}$ is redundant. 
\end{rem}

\begin{rem} In Definition \ref{defn_cord_morse_loop} we used Cascade Morse flow trees. That is slightly different from \cite{petrak2019definition}, where Andreas P. just perturbed appropriately the Morse-Bott function $E_0$ to the Morse setting. So, the cord algebra $Cord^M(K, \mathbb{Z}[\lambda^{\pm1}, \mu^{\pm1}])$ from \cite{petrak2019definition} was not counting cascade flow trees.

Both approaches are expected to be equivalent as a consequence of Multiple-time scale dynamics. In more detail, one can see $X_{E_0}$ as the fast flow with $\Delta_0$ as the normally hyperbolic critical manifold and $-\nabla h_{\Delta_0}$ as the slow flow. For the correspondence of the fast-slow flow lines and flow lines with cascade, see \cite{banyaga2013cascades}. See also Theorem \ref{thm_corresp_trees_gsp}.

Our approach with cascades has an advantage that it will be in the same spirit as the $Cord^M(T_{K, \varepsilon}; \mathbb{Z}[\lambda^{\pm1}, \mu^{\pm1}])$. A bit surprisingly, in the torus case, the cascades will be a significant simplification.
\end{rem}

\begin{thm}\cite{petrak2019definition} For a generic $(K, \mathfrak{f}_K, \ast, X_{E_0}, h_{\Delta_0})$ the sum in (\ref{eqn_dif_cordK}) is finite and the differential $D_K$ is well-defined. Moreover, $Cord^M(K; \mathbb{Z}[\lambda^{\pm1}, \mu^{\pm1}])$ is the knot invariant under a generic ambient isotopy.
\end{thm}

\subsection{Torus case}
\begin{rem}Now, in the similar flavor to the knot case, we would like to define $Cord^M(T_{K, \varepsilon}; \mathbb{Z}[\lambda^{\pm1}, \mu^{\pm1}])$ as a count of certain Cascade Morse flow trees $\clubsuit_{X_{E_\varepsilon}}^{c, out-out}$. Our definition of $Cord^M(T_{K, \varepsilon}; \mathbb{Z}[\lambda^{\pm1}, \mu^{\pm1}])$ will be then strongly influenced by the fact that the trees from $\clubsuit_{X_{E_\varepsilon}}^{c, out-out}$ lie in $M_{K, \varepsilon}$. Namely, we will see that the cuspidal behavior of $M_{K, \varepsilon}$ near $\Delta_\varepsilon$ shifts the expected dimension of the trees by $-1$ for each leaf in $\Delta_\varepsilon$.

So, first, we define $Cord^M(T_{K, \varepsilon}; \mathbb{Z}[\lambda^{\pm1}, \mu^{\pm1}])$ and then we continue with a \textit{conceptual} discussion that $Cord^M(T_{K, \varepsilon}; \mathbb{Z}[\lambda^{\pm1}, \mu^{\pm1}])$ makes sense and recapture the loop space coefficients in the same way as $Cord^MK; \mathbb{Z}[\lambda^{\pm1}, \mu^{\pm1}]).$
\end{rem}

\begin{rem_not}\label{rem_not_extended_treesT}Similarly to Remark \ref{rem_not_extended_trees} we extend cannonically Definition \ref{defn_ribbon_tree_T} of the trees $\clubsuit_{X_{E_\varepsilon}}^{out-out}(\mathcal{T}; p_\varepsilon; q_\varepsilon^1,\dots, q_\varepsilon^m)$. We would like to allow the leaves to be identified not only with the points in $Crit_1(E_\varepsilon\vert_{M_{K, \varepsilon}\setminus \Delta_\varepsilon})$, but also with $\Delta_\varepsilon$, which is a submanifold of the Bott submanifold $\Delta_{full}\subset(\R/T\mathbb{Z}\times S^1)^2$ (recall Definition \ref{defn_out_out} for $\Delta_\varepsilon$ and $\Delta_{full}$). Let us denote such a space of trees by 
$$\clubsuit_{X_{E_\varepsilon}}^{out-out}(\mathcal{T}; p_\varepsilon; Q_\varepsilon^1,\dots, Q_\varepsilon^m).$$
\end{rem_not}

\begin{rem_not}\label{rem_morse_diag_T} Similarly to Remark \ref{rem_morse_diag_K}, let us consider $\Delta_\varepsilon\subset M_{K, \varepsilon}$ with the induced metric from the flat metric on $(\R/T\mathbb{Z}\times S^1)^2$. Recall that by Remark \ref{rem_descr_diag} it holds that $\Delta_\varepsilon\cong (\R/T\mathbb{Z})\times[0, 1]$. By
 $$h_{\Delta_\varepsilon}$$
we denote a Morse function on $\Delta_\varepsilon$ such that $-\nabla h_{\Delta_\varepsilon}$ is Morse-Smale and strictly inward-pointing on $\partial \Delta_\varepsilon$.

For simplicity, we shall assume that $h_{\Delta_\varepsilon}$ has two critical points - $M_{\Delta_\varepsilon}, m_{\Delta_\varepsilon}$ of indices $1, 0$, respectively.

We will also assume $\overline{u_\varepsilon}\cap \big(Crit(h_{\Delta_\varepsilon})\cup W^s_{-\nabla h_{\Delta_\varepsilon}}(M_{\Delta_\varepsilon})\big)=\emptyset$ for every $u_\varepsilon\in\clubsuit_{X_{E_\varepsilon}}^{out-out}$. To the last, we impose a condition on $h_{\Delta_\varepsilon}$ that $\mathfrak{m}_\varepsilon\cap Crit(h_{\Delta_\varepsilon})=\emptyset$, where $\mathfrak{m}_\varepsilon$ is the meridian curve in $T_{K, \varepsilon}$.
\end{rem_not}

\begin{defn} Let $p_\varepsilon\in Crit_2(E_\varepsilon\vert_{M_{K, \varepsilon}\setminus \Delta_\varepsilon})$ and $q_\varepsilon^1, \dots, q_\varepsilon^m\in Crit_1(E_\varepsilon\vert_{M_{K, \varepsilon}\setminus \Delta_\varepsilon})\cup \lbrace m_{\Delta_\varepsilon}\rbrace$, by $\lbrace j_\ell\rbrace_{\ell=1,\dots, k}$ we denote those indices from $\lbrace 1,\dots, m\rbrace$ such that $q_\varepsilon^{j_\ell}=m_{\Delta_\varepsilon}$. Let $\mathcal{T}\in\clubsuit_m.$ Let $X_{E_\varepsilon}$ be a gradient-like vector field adapted to $E_\varepsilon$. Then $u_\varepsilon$ is a \textbf{Cascade Morse flow tree from $p_\varepsilon$ to $q_\varepsilon^1, \dots, q_\varepsilon^m$ modeled on $\mathcal{T}$} if the following holds.

$u_\varepsilon$ is a concatenation of a tree $u_{c, \varepsilon}$ and partial flow lines $\lbrace u^{j_\ell}\rbrace$, where
\begin{itemize}
\item $u_{c, \varepsilon}\in\clubsuit_{X_{E_\varepsilon}}(\mathcal{T}; p_\varepsilon; Q_\varepsilon^1,\dots, Q_\varepsilon^m)$ with $Q_\varepsilon^i=q_\varepsilon^i$ if $i\neq j_\ell$, else $Q_\varepsilon^i=\Delta_\varepsilon$.
\item Let $\overline{u_{c, \varepsilon}}\cap Q_\varepsilon^{j_\ell}= x^{j_\ell}$. Then $u^{j_\ell}$ is given by $\phi^{[0, \infty)}_{-\nabla h_{\Delta_\varepsilon}}(x^{j_\ell})$. In particular, the omega limit of $x^{j_\ell}$ is $m_{\Delta_\varepsilon}$.
\end{itemize}
The convention is that the pairs of concatenated paths are identified with the single edge of $\mathcal{T}$, which will be called \textbf{broken} with the \textbf{cascade} and \textbf{diagonal} part. So the leaves are identified with $q_\varepsilon^1,\dots, q_\varepsilon^m$.

The \textbf{set of all Cascade Morse flow trees from $p_\varepsilon$ to $q_\varepsilon^1, \dots, q_\varepsilon^m$ modeled on $\mathcal{T}$} is denoted by 
$$\clubsuit_{X_{E_\varepsilon}}^{c, out-out}(\mathcal{T}; p_\varepsilon; q_\varepsilon^1, \dots, q_\varepsilon^m).$$
The set of all such trees is denoted by $\clubsuit_{X_{E_\varepsilon}}^{c, out-out}$.
\end{defn}

Now in analogy to Definition \ref{defn_cord_K_dif}, we would like to introduce an evaluation function $\widehat{D}_{u_\varepsilon}$ on a tree $u_\varepsilon\in\clubsuit^{c, out-out}_{X_{E_\varepsilon}}$. Such a function will just count how many times the realization of $u_\varepsilon$ in $T_{K, \varepsilon}\subset\R^3$ intersects the curves $\mathfrak{L}_\varepsilon$ and $\mathfrak{m}_\varepsilon$.

\begin{defn}\label{defn_cord_K_dif} Let $u_\varepsilon\in\clubsuit_{X_{E_\varepsilon}}^c(\mathcal{T}; p_\varepsilon; q_\varepsilon^1, \dots, q_\varepsilon^m)$. Recall also the notation for the trees from Remark \ref{rem_graphs}. Then the evaluation function $\widehat{D}_{u_\varepsilon}$ is recursively defined on the vertices of $\mathcal{T}$ as
\begin{equation*}
\widehat{D}_{u_\varepsilon}(v)=
    \begin{cases}
      \mu^{\alpha_1(\emph{e}_r)}\lambda^{\beta_1(\emph{e}_r)}\widehat{D}_{u_\varepsilon}(v^{out}_{\emph{e}_r})\mu^{\alpha_2(\emph{e}_r)}\lambda^{\beta_2(\emph{e}_r)}, & \text{if }v\text{ is the root},\\
      \\
      \mu^{\alpha_1(\emph{e}^L_v)}\lambda^{\beta_1(\emph{e}^L_v)}\widehat{D}_{u_\varepsilon}(v^{out}_{\emph{e}^L_v})\mu^{\alpha_2(\emph{e}^L_v)+\alpha_1(\emph{e}^U_v)} & \text{if }v\text{ is an interior vertex},\\
      \quad\lambda^{\beta_2(\emph{e}^L_v)+\beta_1(\emph{e}^U_v)}\widehat{D}_u(v^{out}_{\emph{e}^U_v})\mu^{\alpha_2(\emph{e}^U_v)}\lambda^{\beta_2(\emph{e}^U_v)}, &\\
      \\
 q_\varepsilon^j, & \text{if }v\text{ is a leaf represented}\\
        & \text{by the chord }q_\varepsilon^j,
    \end{cases}
\end{equation*}
where, for $k=1, 2,$  $\alpha_k, \beta_k$ are $\mathbb{Z}$-valued functions on the corresponding edges of the tree $u_\varepsilon$. They are $(-1)^{k+1}$ multiples of the algebraic intersection numbers of $\mathfrak{m}_\varepsilon, \mathfrak{L}_\varepsilon$ with $\Gamma_{\varepsilon}(\pi_{s_k, \theta_k}(u_\varepsilon\vert_{\emph{e}}))$, where $u_\varepsilon\vert_{\emph{e}}$ is the restriction of $u_\varepsilon$ to the given edge $\emph{e}$ of $\mathcal{T}$.
\end{defn}

\begin{defn}\label{defn_cord_morse_loopT}Let $\varepsilon>0$ be small. Let $(K, \mathfrak{f}_K, \ast, h_{\Delta_\varepsilon})$ be generic and $X_{E_\varepsilon}$ be a generic approximation of $-\nabla E_\varepsilon$.

Let $C_0^M(T_{K, \varepsilon}; \mathbb{Z}[\lambda^{\pm 1}, \mu^{\pm 1}])$ be a free unital noncommutative $\mathbb{Z}$-algebra generated by $Crit_1(E_\varepsilon\vert_{M_{K, \varepsilon}\setminus \Delta_\varepsilon})\cup\lbrace m_{\Delta_\varepsilon}\rbrace$ and extra generators $\lambda, \mu$ such that $\lambda$ and $\mu$ have inverses and commute together. Let $C_1^M(T_{K, \varepsilon}; \mathbb{Z}[\lambda^{\pm 1}, \mu^{\pm 1}])$ be a $\mathbb{Z}$-vector space generated by $Crit_2(E_\varepsilon\vert_{M_{K, \varepsilon}\setminus\Delta_\varepsilon})$. 
Let $$D_{T_{K, \varepsilon}}:C_1^M(T_{K, \varepsilon}; \mathbb{Z}[\lambda^{\pm 1}, \mu^{\pm 1}])\rightarrow C_0^M(T_{K, \varepsilon}; \mathbb{Z}[\lambda^{\pm 1}, \mu^{\pm 1}])$$ be a linear map which is defined on generators as
\begin{equation}\label{eqn_dif_cordT}
D_{T_{K, \varepsilon}}(p_\varepsilon)=\sum_{u_\varepsilon\in \clubsuit_{X_{E_\varepsilon}}^{c, out-out}}\hbox{sqn}(u_\varepsilon)\widehat{D}_{u_\varepsilon}(v^{p_\varepsilon}).
\end{equation}

Then the \textbf{Morse cord algebra of $T_{K, \varepsilon}$ with loop space coefficients} is the quotient ring
$$Cord^M(T_{K, \varepsilon}; \mathbb{Z}[\lambda^{\pm 1}, \mu^{\pm 1}])=C_0^M(T_{K, \varepsilon}; \mathbb{Z}[\lambda^{\pm 1}, \mu^{\pm 1}])/\mathcal{I}_{T_{K, \varepsilon}},$$
where $\mathcal{I}_{T_{K, \varepsilon}}$ is a two-sided ideal of $C_0^M(T_{K, \varepsilon}; \mathbb{Z}[\lambda^{\pm 1}, \mu^{\pm 1}])$ generated by the image of the linear map $D_{T_{K, \varepsilon}}$ and $\textcolor{blue}{m_{\Delta_\varepsilon}-1}$.
\end{defn}

\begin{rem}\label{rem_gen_MTE}Due to the analogous argument as in Remark \ref{rem_gen_MK} we dropped from the definition of $Cord^M(T_{K, \varepsilon}; \mathbb{Z}[\lambda^{\pm 1}, \mu^{\pm 1}])$ the redundant generator $M_{\Delta_\varepsilon}$.
\end{rem}

\begin{example}
For simplicity, we restrict the coefficient group ring to $\mathbb{Z}_2[\lambda^{\pm 1}, \mu^{\pm 1}]$.

 We would like to characterize $Cord^M(T_{K, \varepsilon}; \mathbb{Z}_2[\lambda^{\pm 1}, \mu^{\pm 1}])$ in the case when $K$ is an unknot $K_U$ with two binormal chords of Morse index $1$ and two binormal chords of Morse index $2$.

For that, we can choose $K_U$ as an ellipse in the $xy$-plane with (oriented) framing $s\mapsto b(s)$. The base point $\ast$ is chosen to lie outside of the endpoints of the binormal chords of $K$. The meridian curve $\mathfrak{m}_\varepsilon$ has an orientation given by the $\theta$-coordinate.

Next, $(Crit_2(E_\varepsilon)\cap M_{K_U, \varepsilon})\setminus\Delta_\varepsilon$ consists of two critical point - $p_\varepsilon, \widehat{p}_\varepsilon$. And $((Crit_1(E_\varepsilon)\cap M_{K_U, \varepsilon})\setminus \Delta_\varepsilon)\cup\lbrace m_{\Delta_\varepsilon}\rbrace=\lbrace m_{\Delta_\varepsilon}\rbrace$.

Let us pick $p_\varepsilon$. It follows that $D_{T_{K_U, \varepsilon}}(p_\varepsilon)$ counts only cascade flow lines from $\clubsuit_{X_{E_\varepsilon}}^{c, out-out}(p_\varepsilon, m_{\Delta_\varepsilon})$. Let us discuss the expected dimension of $\clubsuit_{X_{E_\varepsilon}}^{c, out-out}(p_\varepsilon, m_{\Delta_\varepsilon})$. We would like to focus on the cascade part of $\clubsuit_{X_{E_\varepsilon}}^{c, out-out}(p_\varepsilon, m_{\Delta_\varepsilon})$ first. Note that on $(\R/T\mathbb{Z}\times S^1)^2$ it holds that 
$$\dim\big(W^u_{X_{E_\varepsilon}}(p_\varepsilon)\big)=2\hbox{ and }\dim\big(W^s_{X_{E_\varepsilon}}(\Delta_{full})\big)=4.$$
So, $\dim W_{X_{E_\varepsilon}}(p_\varepsilon, \Delta_{full})=2$, which is by $1$ dimension higher, then we wanted! However, we are not counting trajectories in $(\R/T\mathbb{Z}\times S^1)^2$ but rather in $M_{K_U, \varepsilon}$.

Now, we would like to recall few facts about $M_{K_U, \varepsilon}$ near $\Delta_\varepsilon$ and the structure of the dynamics on that space.

By Corollary \ref{lem_cusp_fibre}, there is a $O(\varepsilon)$-close subset $\Delta_{\varepsilon}^{cusp}\subset \Delta_\varepsilon$ with the following property. Near $\Delta_{\varepsilon}^{cusp}$, the normal projection $\pi_{\mathcal{N}}$ gives $M_{K_U, \varepsilon}$ a structure of locally trivial stratified fiber bundle $(=:\widehat{M}_{K_U, \varepsilon})$ over $\Delta_{\varepsilon}^{cusp}$, and the typical fiber is a union of two cuspical regions.

Next, in Conjecture \ref{conj_eating_cusp}, we slightly restricted $\widehat{M}_{K_U, \varepsilon}$ to some $\widehat{\widehat{M}}_{K_U, \varepsilon}$, roughly said, by shrinking $\Delta_\varepsilon^{cusp}\rightarrow\Delta_{\varepsilon, \delta}$. Then, by this conjecture, it follows that
$$\textcolor{teal}{\dim\big(W^s_{X_{E_\varepsilon}}(\Delta_{\varepsilon, \delta})\cap M_{K_U, \varepsilon}\big)=3.}$$
We remark that Conjecture \ref{conj_eating_cusp} is a consequence of the cuspical geometry of $M_{K_U, \varepsilon}$, Sternberg's Linearization Theorem \ref{thm_sternberg} and the outward-pointing behavior of $X_{E_\varepsilon}$ on $\partial \widehat{\widehat{M}}_{K_U, \varepsilon}$. For the behavior of $X_{E_\varepsilon}$ on $\partial \widehat{\widehat{M}}_{K_U, \varepsilon}$ we only gave an evidence in Example \ref{exm_diag_grad}, but we expect that a laborious local computation as in Section \ref{sec:closer_diag} will verify the claim.

Hence, for a generic $X_{E_\varepsilon}$ it holds that
$$\dim\big(W_{X_{E_\varepsilon}}(p_\varepsilon, \Delta_{\varepsilon, \delta})\big)=1.$$
Let us put $\delta:=\varepsilon$, so $\Delta_{\varepsilon, \delta}$ is $O(\varepsilon)$-close subset of $\Delta_\varepsilon$. So, for $\varepsilon>0$ small, we expect that 
$$\dim\big(\clubsuit_{X_{E_\varepsilon}}^{out-out}(p_\varepsilon, \Delta_\varepsilon)\big)=0,$$
see also Figure \ref{figure_unstable_diagonal}, where the trajectories were computed numerically with \textit{Mathematica} (in the case that $K$ is a circle in $xy$-plane).
\begin{figure}[!htbp]
\labellist
\pinlabel $c_{p_\varepsilon}$ at 330 120
\pinlabel $\times$ at 178 255
\pinlabel $c_{triv}$ at 190 276
\endlabellist
\centering
\includegraphics[scale=0.57]{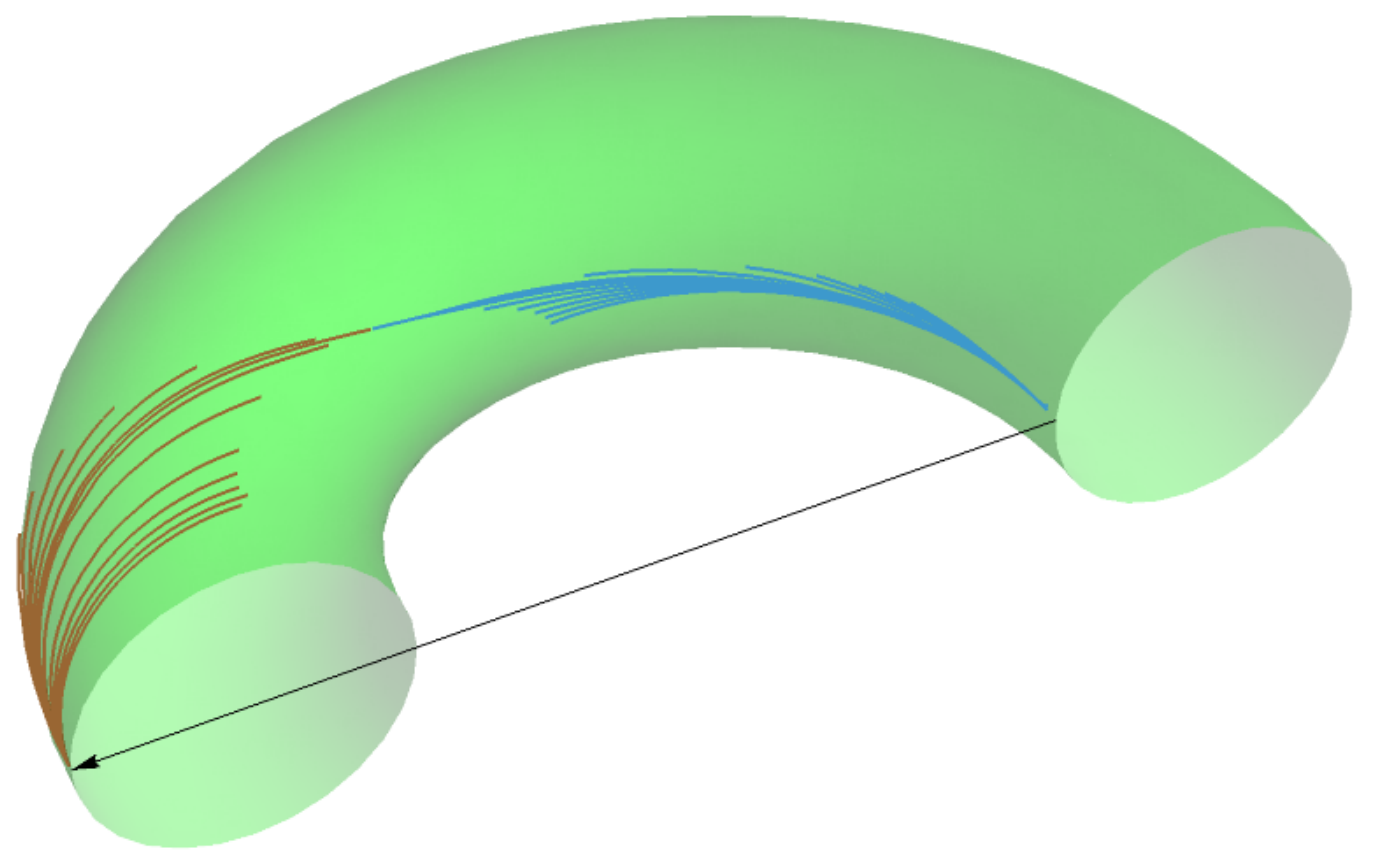}
\vspace{0.3cm}
\caption{A visualization of $W^u_{-\nabla E_\varepsilon}(p_\varepsilon)\cap M_{K_U, \varepsilon}$ on $T_{K_U, \varepsilon}$. The \textcolor{red}{red} curves depict the endpoints of chords emanating from $c_{\varepsilon}$ under the negative gradient flow constraint to $M_{K_U, \varepsilon}$. Similarly, the \textcolor{blue}{blue} curves depict the starting points of chords emanating from $c_{\varepsilon}$ under the negative gradient flow constraint to $M_{K_U, \varepsilon}$. 
Note also a single pair of red and blue curves that are connected by a trivial chord $c_{triv}$; they represent a trajectory $u_{c, \varepsilon}$ from $p_\varepsilon$ to $\Delta_\varepsilon$ that stays the whole time in $M_{K_U, \varepsilon}$.}
\label{figure_unstable_diagonal}
\end{figure}

In Figure \ref{figure_unstable_diagonal} we saw a single $u_{c, \varepsilon}$ which was a cascade part of some trajectory $u_\varepsilon\in\clubsuit^{c, out-out}_{X_{E_\varepsilon}}(p_\varepsilon, m_{\Delta_\varepsilon})$. Or more precisely, we saw a trace of $u_{c, \varepsilon}$ in $T_{K_U, \varepsilon}$. That is the curve, which is given by concatenation of
$$u_{c, \varepsilon}^{trace}:=\textcolor{red}{\Gamma_\varepsilon(\pi_{s_2, \theta_2}(u_{c, \varepsilon}^{backward}))}\ast \textcolor{blue}{\Gamma_\varepsilon(\pi_{s_2, \theta_2}(u_{c, \varepsilon}^{forward}))}.$$
However, then there is also another trajectory $\widetilde{u}_\varepsilon\in\clubsuit^{c, out-out}_{X_{E_\varepsilon}}(p_\varepsilon, m_{\Delta_\varepsilon})$ such that the trace $\widetilde{u}_{c, \varepsilon}^{trace}$ of the cascade part of $\widetilde{u}_{c, \varepsilon}$ is symmetric to $u_{c, \varepsilon}^{trace}$ along the $xy$-plane. See also Figure \ref{figure_two_trajectories}.

\begin{figure}[!htbp]
\labellist
\pinlabel $c_{p_\varepsilon}$ at 330 120
\pinlabel $u_{c, \varepsilon}^{trace}$ at 270 300
\pinlabel $\widetilde{u}_{c, \varepsilon}^{trace}$ at 360 200
\endlabellist
\centering
\includegraphics[scale=0.45]{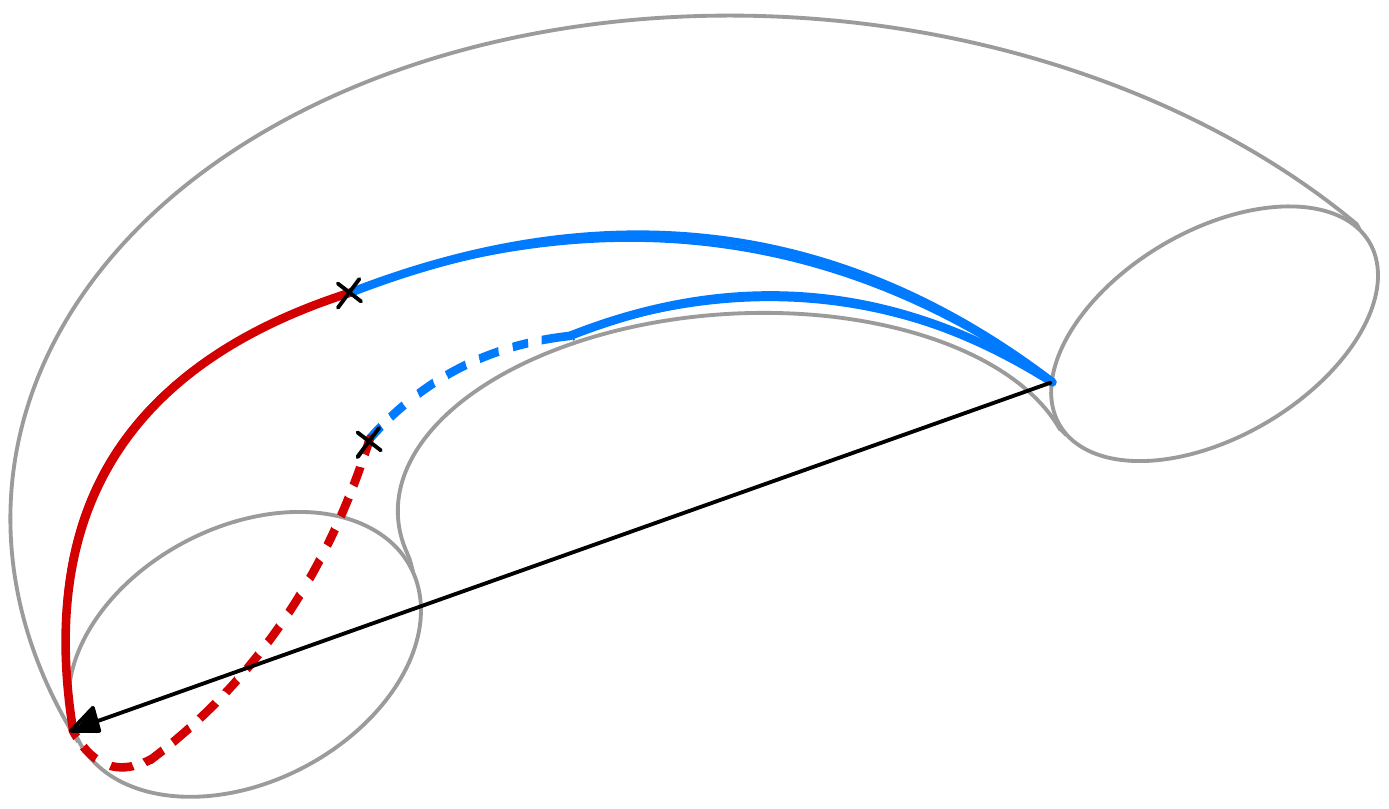}
\vspace{0.3cm}
\caption{Traces of $u_{c, \varepsilon}$ and $\widetilde{u}_{c, \varepsilon}$.}
\label{figure_two_trajectories}
\end{figure}

Let us now compute the contributions of paths $u_\varepsilon$ and $\widetilde{u}_\varepsilon$ to the differential $D_{T_{K_U, \varepsilon}}$. Analogously to Remark \ref{rem_gen_MTE}, the diagonal parts of $u_\varepsilon$ and $\widetilde{u}_\varepsilon$ are redundant in the count. So only the intersections of the traces $u_{c, \varepsilon}^{trace}, \widetilde{u}_{c, \varepsilon}^{trace}$ with $\mathfrak{m}_\varepsilon$ and $\mathfrak{L}_\varepsilon$ matter. Let us consider the orientations as in Figure \ref{figure_two_trajectories_loop}.

\begin{figure}[!htbp]
\labellist
\pinlabel $c_{p_\varepsilon}$ at 330 120
\pinlabel $\textcolor{teal}{\mathfrak{L}_\varepsilon}$ at 240 325
\pinlabel $\textcolor{orange}{\mathfrak{m}_\varepsilon}$ at 465 380
\pinlabel $(\Gamma_\varepsilon)_\ast\partial_s$ at 590 304
\pinlabel $-(\Gamma_\varepsilon)_\ast \partial_\theta$ at 560 220
\endlabellist
\centering
\includegraphics[scale=0.45]{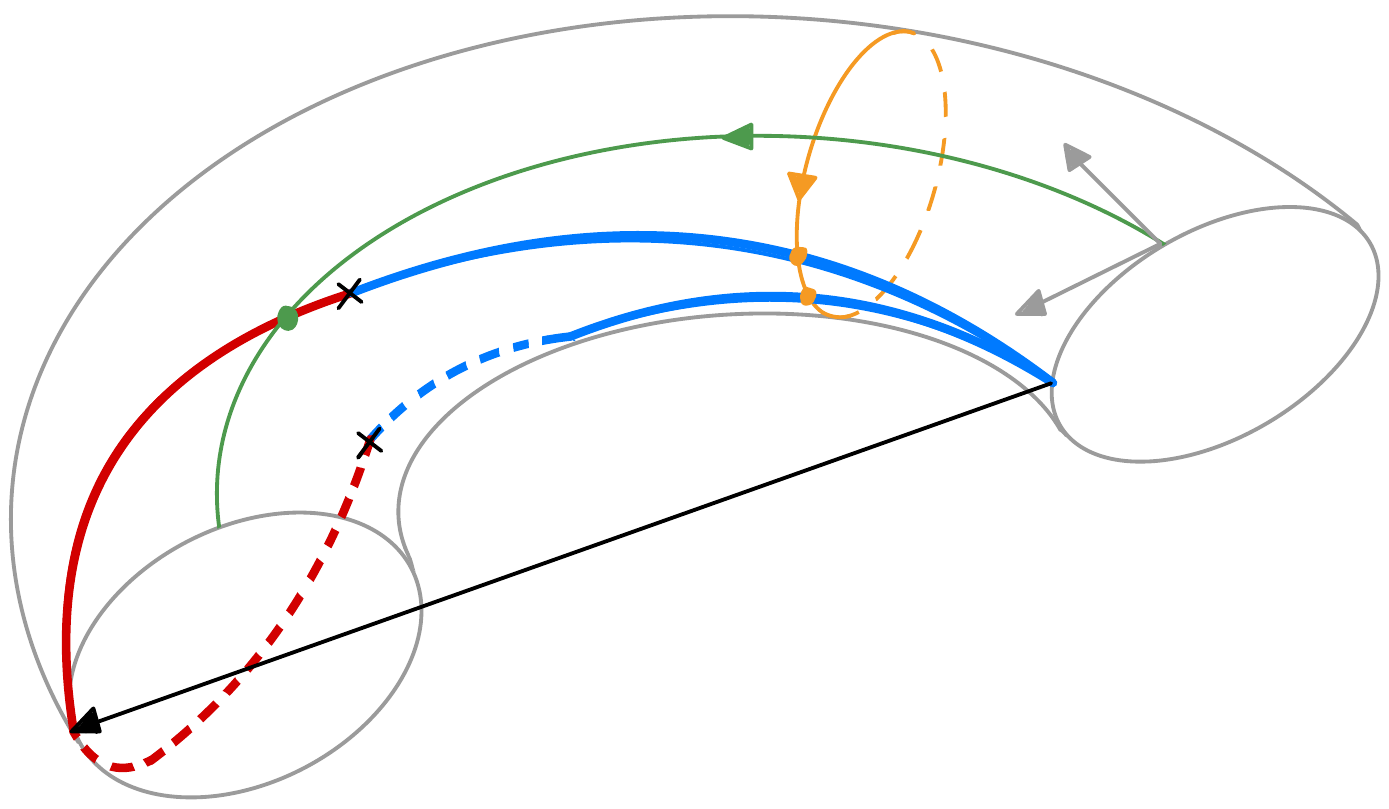}
\vspace{0.3cm}
\caption{Orientations of curves $\textcolor{teal}{\mathfrak{L}_\varepsilon}, \textcolor{orange}{\mathfrak{m}_\varepsilon}$ and their intersections with traces $u_{c, \varepsilon}^{trace}$ and $\widetilde{u}_{c, \varepsilon}^{trace}$. The orientations are induced by the basis $\lbrace \partial_s, -\partial_\theta\rbrace.$}
\label{figure_two_trajectories_loop}
\end{figure}
So the algebraic intersection number of $\mathfrak{m}_\varepsilon$ with the curves $u_{c, \varepsilon}^{trace}, \widetilde{u}_{c, \varepsilon}^{trace}$ is in both cases $+1$. This will contribute in both cases by $\cdot\lambda$.

More intersecting is the algebraic intersection number of $\mathfrak{L}_\varepsilon$ with the curves $u_{c, \varepsilon}^{trace}, \widetilde{u}_{c, \varepsilon}^{trace}$. In the first case, we will get a contribution by $\textcolor{teal}{\cdot\mu}$ and in the second case, there is no intersection, so the contribution is $\textcolor{teal}{\cdot 1}$! Compare this with the skein relation $(i.)$ from Definition \ref{defn_cord_top}.

Since on the other half of the torus there are another two cascade flow lines of $\clubsuit_{X_{E_\varepsilon}}^{c, out-out}(p_\varepsilon, m_{\Delta_\varepsilon})$, it follows that 
\begin{align*}
D_{T_{K_U, \varepsilon}}(p_\varepsilon)&=\lambda(\mu^{\alpha_1}\lambda^{\beta_1}m_{\Delta_\varepsilon}\mu^{-\alpha_1}\lambda^{-\beta_1})\mu\\
&\quad+\lambda(\mu^{\alpha_2}\lambda^{\beta_2}m_{\Delta_\varepsilon}\mu^{-\alpha_2}\lambda^{-\beta_2})\\
&\quad+(\mu^{\alpha_3}\lambda^{\beta_3}m_{\Delta_\varepsilon}\mu^{-\alpha_3}\lambda^{-\beta_3})\mu\\
&\quad+(\mu^{\alpha_4}\lambda^{\beta_4}m_{\Delta_\varepsilon}\mu^{-\alpha_4}\lambda^{-\beta_4}),
\end{align*}
for some numbers $\alpha_i, \beta_i\in\mathbb{Z}$ that come from the diagonal parts of the trajectories.

Since in cord algebra it holds that $m_{\Delta_\varepsilon}=1$, we obtain the relation
$$(\lambda+1)(\mu+1)=0.$$
The same relation will also follow from the cascade flow trajectories from $\widehat{p}_\varepsilon$ to $m_{\Delta_\varepsilon}$. Hence
$$Cord^M(T_{K_U, \varepsilon}; \mathbb{Z}_2[\lambda^{\pm 1}, \mu^{\pm 1}])\cong \mathbb{Z}_2[\lambda^{\pm 1}, \mu^{\pm 1}]/((\lambda+1)(\mu+1)).$$
Which matches the computation of $Cord^M(K_U; \mathbb{Z}_2[\lambda^{\pm 1}, \mu^{\pm 1}])$ from \cite{petrak2019definition}. 
\end{example}
\chapter{Relation to symplectic geometry}
\label{ch:rel_symplectic}
In this chapter we attempt to show that our definitions of cord algebras for $T_{K, \varepsilon}$ did not fall from the sky, but can be expected to arise from a $J$-holomorphic curve invariant for one component of unit conormal bundle of $T_{K, \varepsilon}$ - $\mathcal{L}^\ast_+ T_{K, \varepsilon}$. This invariant will be $0$-th degree Legendrian contact homology $LCH_0(\mathcal{L}^\ast_+ T_{K, \varepsilon})$ (with $\mathbb{Z}$ or $\mathbb{Z}[\pi_1(\mathcal{L}^\ast_+T_K)]$ coefficients).

There are two approaches how to relate $Cord^M(T_{K, \varepsilon})$ with $LCH_0(\mathcal{L}^\ast_+ T_{K, \varepsilon})$. Direct and indirect.

The indirect approach is based on the fact that the identification of $Cord^M(K)$ with $LCH_0(\mathcal{L}^\ast K)$ is known, \cite{Cieliebak2016KnotCH, petrak2019definition}. We recall that this identification uses the so-called $0$-th string homology $H^{string}_0(K)$. In addition, the identification of $LCH_0(\mathcal{L}^\ast K)$ with $LCH_0(\mathcal{L}^\ast_+ T_{K, \varepsilon})$ is dictated by the canonical Legendrian isotopy and the relation $Cord^M(K)$ to $Cord^M(T_{K, \varepsilon})$ was the aim of the last chapters.

So, this chapter will be mainly focused on outlining the direct approach, i.e., the approach that does not involve the theories for the underlying knot. For this we construct propose a model of the $0$-th string homology $H^{string}_0(T_{K, \varepsilon})$ as an intermediate step between $LCH_0(\mathcal{L}^\ast_+ T_{K, \varepsilon})$ and $Cord^M(T_{K, \varepsilon})$. To obtain a chain-level map from $LCH_0(\mathcal{L}^\ast_+ T_{K, \varepsilon})$ to the string homology $H^{string}_0(T_{K, \varepsilon})$ we introduce a moduli space $\mathcal{M}_{T_{K, \varepsilon}}$ of $J$-holomorphic curves with boundaries on the Lagrangian manifold $L^\ast_+ T_{K, \varepsilon}\cup \R^3$ with an arboreal singularity along $T_{K, \varepsilon}$. We will also analyze the geometric information about $J$-holomorphic curves which arises from the local behavior near the arboreal singularity $T_{K, \varepsilon}$. 

\section{Symplectic set-up}

\begin{assump} We shall assume that the reader is familiar with some basic notions of symplectic geometry on the level of Appendices \ref{ch:file1} and \ref{ch:file2}. We will briefly recall some of the definitions and conventions.
\end{assump}

\begin{rem_not}\label{rem_not_symplectic} On the cotangent bundle $T^\ast\R^3$ we consider coordinates $(q, p)$. Then the Liouville one form $\lambda=pdq$ determines the symplectic form $\omega=d\lambda$. So the Liouville vector field $X$ is given by $X=p\partial_{p}$. Moreover, the Hamiltonian $H(x)=(1/2)\Vert x^\sharp\Vert^2_{g_{Euc}}$ determines the Hamiltonian vector field $X_H=p\partial_q$.

The unit cotangent bundle $S^\ast\R^3$ is a contact manifold with a contact form $\lambda_1$ which is given by the restriction of $\lambda$ to $S^\ast\R^3$. $\xi$ denotes the contact structure $\xi=\ker \lambda_1$. Also, the restriction of $X_H$ to $S^\ast\R^3$ is the Reeb vector field $R$.

For a submanifold $Q\subset\R^3$, the conormal bundle $L^\ast Q$ is a Lagrangian submanifold of $T^\ast \R^3$. Then the unit conormal bundle $\mathcal{L}^\ast Q=L^\ast Q\cap S^\ast\R^3$ is a Legendrian submanifold of $S^\ast\R^3$.

Especially, if $Q=T_K$, then the Legendrian submanifold $\mathcal{L}^\ast T_K$ has two connected components: $\mathcal{L}^\ast_+ T_K$ and $\mathcal{L}^\ast_- T_K$. They are both diffeomorphic to $T_K$ by translations in (co)normal directions. $\mathcal{L}^\ast_+ T_K$ is determined by a choice of the co-orientation on $T_K$. Our convention will be that we co-orient $T_K$ by a normal vector that is pointing \textit{inside} of $\nu K$.

Finally, $L^\ast_+ T_K$ and $L^\ast_- T_K$ will be two Lagrangian submanifolds with the boundary $T_K$ such that $\mathcal{L}^\ast_\pm T_K\subset L^\ast_\pm T_K$ and $L^\ast_+ T_K\cup_{T_K} L^\ast_- T_K= L^\ast T_K.$

To any Reeb chord $\bm{a}$ on $\mathcal{L}^\ast_+T_K$ we can assign its degree $\vert \bm{a}\vert$ as in Definition \ref{defn_reeb_pure}. $\vert\bm{a}\vert$ is equal to a well-defined Maslov index of a certain loop of Lagrangian subspaces minus $1$. Roughly speaking, such a loop is constructed by a concatenation of $T\mathcal{L}^\ast_+T_K$ along a capping path from the endpoint of $\bm{a}$ to the starting point of $\bm{a},$ translation of $T\mathcal{L}^\ast_+T_K$ along $\bm{a}$ with the linearized Reeb flow and that all is closed by a positive rotation at the endpoint of $\bm{a}$. The trivialization of the corresponding symplectic bundle is taken over a capping disc.
\end{rem_not}

\begin{rem}\label{rem_symplectomorphism} There is a canonical symplectomorphism
\begin{align*}
\varphi:(\R\times S^\ast \R^3, d(e^s\lambda_1))&\rightarrow(T^\ast \R^3\setminus \R^3, \omega)\\
(s, q, p)&\mapsto (q, e^s p).
\end{align*} 
For the Liouville vector field $\partial_s$ on the symplectization, it holds that $\varphi_\ast \partial_s=X$, hence $\varphi$ is a Liouville diffeomorphism.

Moreover, $\varphi$ induces $\R$-invariant extensions of $\lambda_1, R, \xi$ on $T^\ast N\setminus N$ that will be denoted by the same letters. More precisely,
$$\lambda_1=|p|^{-1}p dq,\quad R=|p|^{-1}p\partial_q\quad\hbox{ and }\quad\xi=\langle R, X\rangle^{\bot \omega}.$$

Also, the $2$-form form $d\lambda_1$ is symplectic on $\xi$ and will be denoted by $\omega_1$. Next, $\varphi^\ast H=(1/2)e^{2s}$, hence $(\varphi^{-1})_\ast X_H=e^s R$. And moreover, $\varphi$ is a canonical diffeomorphism between Lagrangian submanifolds $\R\times\mathcal{L}^\ast Q$ and $L^\ast Q\setminus Q$.
\end{rem}

\begin{rem}\label{rem_proof_eqn_lin_ham_flow_liouv} We would like to deduce Formula (\ref{eqn_lin_ham_flow_liouv}), i.e. $d\phi_H^t(x)X_x=X_y+tR_y$, where $y=\phi^t_{X_H}(x)$.

Since $\varphi$ is a Liouville diffeomorphism, we can make the computation in the symplectization $\R\times S^\ast \R^3$.

Recall the notation $x=(q, p)$. By Remark \ref{rem_symplectomorphism}, in $\R\times S^\ast \R^3$ we have that
$$\phi^t_H:(s, x)\mapsto (s, \phi_R^{e^r t}(x)).$$
Let $x\in S^\ast \R^3$, then we define a curve 
\begin{align*}
\gamma:(-\varepsilon, \varepsilon)&\rightarrow \R\times S^\ast \R^3\\
s &\mapsto (s, x).
\end{align*}
Hence $\gamma$ is an integral curve of $X=\partial_s$. Now we put $u(s):=\phi^t_H(\gamma(s))=(s, \phi_R^{e^r t}(x)).$ Then
$$d\phi^t_H(x, 0) X=\frac{d}{ds}\Bigr|_{s=0}u(s)=t e^0 R_{(y, 0)}+X_{(y, 0)}=tR_{(y, 0)}+X_{(y, 0)},$$
where $(y, 0)=\phi^t_H(x, 0)$. The same computation works for any smooth manifold $N$ instead of $\R^3$.
\end{rem}

\begin{lemma}\label{lemma_isotop_conormal} $\mathcal{L}^\ast K$ and $\mathcal{L}^\ast_+T_{K, \varepsilon}$ are canonically Legendrian isotopic.
\end{lemma}

\begin{proof}
By Remark \ref{rem_proof_eqn_lin_ham_flow_liouv} we have that $\phi^t_R(q, p)=(q+t p, p)$, hence by the construction $\phi^\varepsilon_R(\mathcal{L}^\ast K)=\mathcal{L}^\ast_+T_{K, \varepsilon}$. Next, $\lbrace\phi_R^t\rbrace_{t\in[0, \varepsilon]}$ is an ambient contact isotopy and, in particular, also the desired Legendrian isotopy.
\end{proof}

\begin{assump} From now on, we assume that $K$, and hence also $T_{K, \varepsilon}$, are real analytic. Note that real analyticity can be achieved from a smooth manifold by a $C^1$-perturbation, see \cite[5.8]{Cieliebak2012FromST}. Moreover, recall that for $\varepsilon>0$ small and $K$ generic, the function $E_\varepsilon$ is Morse-Bott by Lemmata \ref{lemma_generic_morse_knot} and \ref{lemma_morse_corresp}. Thus, by Theorem \ref{lemma_index_formula} we will assume that all Reeb chords on $\mathcal{L}^\ast_+T_{K, \varepsilon}\subset S^\ast\R^3$ are nondegenerate.
\end{assump}

\begin{defn}[\cite{Cieliebak2016KnotCH}] \label{defn_adms} An $\omega$-compatible almost complex structure $J$ on $T^\ast\R^3$ is called \textbf{admissible} if
\begin{itemize}
\item[$(i.)$] On $T^\ast\R^3\setminus D^\ast\R^3\cong \R_+\times S^\ast \R^3$ it holds that $J|_\xi$ is compatible with $\omega_1|_\xi$, invariant under translations in the Liouville direction, $J(X)=R$ and $J(R)=-X$.
\item[$(ii.)$] On $T^\ast\R^3\setminus T_K$ it holds that $J|_\xi$ is compatible with $\omega_1|_\xi$ and preserves $\langle X, R\rangle$. Along the zero section $0_{\R^3}$ the almost complex structure $J$ agrees with the standard almost complex structure $J_{st}$ defined by $J_{st}(\partial_{q^i}):=-\partial_{p_i}.$
\item[$(iii.)$] $J$ is integrable near $T_K$ such that $\R^3$ and $T_K$ are real analytic.
\item[$(iv.)$] If $s>1$ and $\gamma_R$ is a Reeb chord in $(\lbrace s\rbrace\times S^\ast\R^3, \lambda_1)$ with endpoints on $\lbrace s\rbrace\times\mathcal{L}^\ast_+T_K$, then $\varphi^\ast J$ is ``adapted'' to $\gamma_R$ and $\lbrace s\rbrace\times\mathcal{L}^\ast_+T_K$. 

Roughly speaking, near $\gamma_R$ there is a Pfaff-Darboux chart such that in the Lagrangian projection $\varphi^\ast J$ is integrable and the sheets of $\mathcal{L}^\ast_+T_{K}$ are real analytic. For details, see \cite{Cieliebak2016KnotCH}.
\item[$(v.)$] There is an exhaustion $\R^3=\bigcup_{k\in\mathbb{N}}V_k$ by compact domains $V_k$ with smooth boundaries such that pullbacks $\pi^{-1}(\partial V_k)$ to $T^\ast\R^3$ are $J$-convex hypersurfaces, for the definition of $J$-convex hypersurfaces see {\normalfont\cite{Cieliebak2012FromST}}.
\end{itemize}
Moreover, a $d(e^s\lambda_1)$-compatible almost complex structure $J_1$ on $\R\times S^\ast\R^3$ is called \textbf{admissible} if $J_1$ satisfies conditions $(i.)$ and $(iv.)$ on the whole symplectization.
\end{defn}

\begin{rem} By \cite[pg 38]{SFT00} a contact manifold $(M, \alpha)$ admits an exhaustion $M=\bigcup_{k\in\mathbb{N}}V_k$ by compact domains $V_k$ with smooth pseudo-convex boundaries, if for each $k$ any trajectory of the Reeb vector field $R|_{V_k}$ does not have an internal tangency with $\partial V_k$. This is satisfied by $(S^\ast\R^3, \lambda_1)$. Indeed, by Example \ref{empl_contc} $(S^\ast\R^3, \lambda_1)$ is contactomorphic to $(T^\ast S^{\hspace{0.07em}2}\times\R, dt\times \lambda_{S^\ast S^{\hspace{0.07em}2}})$, hence $V_k:=D^\ast_kS^{\hspace{0.07em}2}\times[-k, k]$ will be the desired exhaustion, here $D^\ast_k S^{\hspace{0.07em}2}=\lbrace x\in T^\ast S^{\hspace{0.07em}2}\,|\,||x^\sharp||\leq k\rbrace.$
\end{rem}

\begin{example}Let $\rho:[0, \infty)\rightarrow[1, \infty)$ be a smooth function such that $\rho(r)=1$ for $r<1/3$ and $\rho(r)=r$ for $r>2/3$. Let us take $g_{Euc}$ on $\R^3$ together with the induced global coordinates $(q, p)$ on $(T^\ast\R^3, \omega)$ and define an almost complex structure $J_{\rho}$ on $T^\ast\R^3$ by
$$J_{\rho}(\partial_{q^i}):=-\rho(\Vert p\Vert)\partial_{p_i}\hbox{ and }J_{\rho}(\partial_{p_i}):=\rho(\Vert p\Vert)^{-1}\partial_{q^i}.$$
Then $J_\rho$ is $\omega$-compatible and satisfies conditions $(i.)$-$(iii.), (v.)$ from Definition \ref{defn_adms}. After a suitable deformation near infinity through compatible almost complex structures, we will obtain $J_\rho$ that satisfies also $(iv.)$.
\end{example}

\newpage

The following lemma is an analog to \cite[Lemma 8.6]{Cieliebak2016KnotCH}, where suitable charts around $K$ were constructed.

\begin{lemma} \label{lemma_coord} Let $J$ be an almost complex structure on $T^\ast \R^3$ with totally real submanifolds $\R^3$ and $T_K$. Moreover, let $J$ be integrable near $T_K$ such that $\R^3$ and $T_K$ are real analytic. 

Let $\delta>0$ be small, then we put $N:=S^1\times S^1\times (-\delta, \delta)$.  Then there exists a holomorphic embedding from some neighborhood $\mathcal{O}p\big(N, N^{\mathbb{C}}\big)$ to $(T^\ast\R^3, J)$ such that
\begin{itemize}
\item[$(i.)$] $S^1\times S^1\times\lbrace 0\rbrace\times \lbrace 0\rbrace$ is mapped to $T_K$,
\item[$(ii.)$] $S^1\times S^1\times\lbrace 0\rbrace\times (\R\cap(-\delta_1, \delta_1)^2)$ is mapped to $\R^3$,
\item[$(iii.)$] $S^1\times S^1\times\lbrace 0\rbrace\times (i\R_{\geq 0}\cap(-\delta_1, \delta_1)^2)$ is mapped to $L^\ast_+ T_K$,
\end{itemize}
where $\delta_1>0$ is small.

In particular, this gives us local holomorphic coordinates $\mathbb{C}^2\times\mathbb{C}$ around any point $x\in T_K$, we call them the \textbf{standard coordinates around $x$}.
\end{lemma}

\begin{proof}
Let us take a real analytic embedding $\widetilde{\Gamma}:S^1\times S^1\hookrightarrow\R^3$ which determines $T_K$. Next, $w_1$ and $w_2$ will be real analytic unit vector fields tangent to $T_K$ that are induced by the components of $S^1\times S^1$. Then the unit real analytic vector field $n:=w_1\times w_2$ is perpendicular to $T_K$ (with respect to $g_{Euc}$). We can assume that $w_1, w_2$ were chosen such that $\lbrace n, w_1, w_2\rbrace$ is a positively oriented basis frame on $T\R^3$ and $n$ is pointing inside of $\nu K$. 

If $(s, \theta, \rho)$ are coordinates on $N:=S^1\times S^1\times (-\delta, \delta)$, then 
\begin{align*}
\phi: N&\rightarrow\R^3\\
(s, \theta, \sigma)&\mapsto \widetilde{\Gamma}(s, \theta)+\rho n(s, \theta)
\end{align*} 
is a real analytic embedding for $\delta>0$ sufficiently small. Now, by Bruhat-Whitney Theorem, \cite{Cieliebak2012FromST}, $N^\mathbb{C}$ is a well-defined complex manifold such that $N$ is a real analytic totally real submanifold of $N^\mathbb{C}$. Let $\sigma$ be the conjugation on $(T^\ast\R^3, -J)$. Then by \cite[Lemma 5.40]{Cieliebak2012FromST} the composition $\sigma\circ\phi^\mathbb{C}$ gives the desired $(J, i)$-holomorphic embedding from some small neighborhood $\mathcal{O}p\big(N, N^{\mathbb{C}}\big)$ to $(T^\ast\R^3, J)$. 
\end{proof}

\newpage

\section{Moduli spaces}
\begin{defn}\label{defn_mod_switch} Let $m\in\N_0$, $J_1$ be an admissible almost complex structure on $(\R\times S^\ast, d(e^s\lambda_1))$ and $D_{m+1}$ be an unit $(m+1)$-punctured disc with some fixed choice of counterclockwisely placed boundary punctures $p_0,\dots, p_{m+1}$. Then let $\textbf{a}, \textbf{b}_1, \dots, \textbf{b}_m$ be Reeb chords on $\mathcal{L}^\ast_+T_K$ that are parametrized by times $T_0, T_1,\dots, T_m$, respectively. We define a moduli space
$$\mathcal{M}^{sy}_{T_K}(\textbf{a}, \textbf{b}_1, \dots, \textbf{b}_m):=\lbrace u=(a, f)\in C^0(D_{m+1}, \R\times S^\ast\R^3)\,|\,(i.)-(v.)\rbrace/\sim,$$
where $\sim$ is a conformal reparametrization of the domain and
\begin{itemize}
\item[$(i.)$]On the interior the map $u$ is $(J_1, j_\kappa)$-holomorphic, i.e. 
$$du+J_1\circ du\circ j_\kappa=0,$$
where $\kappa$ is some conformal structure on $D_{m+1}$.
\item[$(ii.)$]$Im(u|_{\partial D_{m+1}})\subset\R\times \mathcal{L}^\ast_+ T_K$.
\item[$(iii.)$]For a half-strip neighborhood of $p_0$ there is a constant $a_0$ such that
$$a(\tau, t)-T_0\tau-a_0\rightarrow 0,\quad f(\tau, t)\rightarrow\textbf{a}(T_0t)$$
uniformly in $t$ as $\tau\rightarrow\infty$. We say that $p_0$ is a \textbf{positive puncture asymptotic to} $\textbf{a}$.
\item[$(iv.)$]If $k\in\lbrace 1,\dots, m\rbrace$, then for a half-strip neighborhood of $p_k$ there is a constant $a_k$ such that 
$$a(\tau, t)+T_k\tau-a_k\rightarrow 0,\quad f(\tau, t)\rightarrow\textbf{b}_k(-T_kt)$$
uniformly in $t$ as $\tau\rightarrow\infty$. We say that $p_k$ is a \textbf{negative puncture asymptotic to} $\textbf{b}_k$.
\item[$(v.)$] $u$ has a finite Hofer energy, see \cite{Bourgeois_2003}.
\end{itemize} 
\end{defn}

\begin{thm}\label{thm_dim_symp} \cite{EES05, Cieliebak2016KnotCH, karlsson2019note} 
For a generic admissible almost complex structure $J_1$ the moduli space $\mathcal{M}^{sy}_{T_K}(\textbf{a}, \textbf{b}_1, \dots, \textbf{b}_m)$ is a manifold of the dimension
$$\dim\mathcal{M}^{sy}_{T_K}(\textbf{a}, \textbf{b}_1, \dots, \textbf{b}_m)=|\bm{a}|-\sum_{k=1}^m |\bm{b}_k|.$$
If $\mathcal{L}^\ast T_K$ is spin, there is a choice of coherent orientations for $\mathcal{M}^{sy}_{T_K}(\textbf{a}, \textbf{b}_1, \dots, \textbf{b}_m)$. 
\end{thm}

\begin{rem} Since the almost complex structure $J_1$ is invariant under $\R$-translations, the quotient $\mathcal{M}^{sy}_{T_K}(\textbf{a}, \textbf{b}_1, \dots, \textbf{b}_m)/\R$ becomes an oriented manifold of dimension $|\bm{a}|-\sum_{k=1}^m |\bm{b}_k|-1$.
\end{rem}

\begin{defn}\label{defn_mod_symp}Let $m\in\N_0$, $J$ be an admissible almost complex structure on $(T^\ast\R^3, \omega)$ and $D_{m+1}$ be a unit $(m+1)$-punctured disc with some fixed choice of counterclockwisely placed boundary punctures $p_0,\dots, p_{m+1}$. Then let $\textbf{a}$ be a Reeb chord on $\mathcal{L}^\ast_+T_K$ and $\textbf{n}:=(n_1,\dots, n_m)$ an $m$-tuple of positive half-integers. We define a moduli space
$$\mathcal{M}_{T_K}(\textbf{a}, \textbf{n}):=\lbrace u\in C^0(D_{m+1}, T^\ast\R^3)\,|\,(i.)-(v.)\rbrace/\sim,$$
where $\sim$ is a conformal reparametrization of the domain and
\begin{itemize}
\item[$(i.)$]$u$ is $(J, j_\kappa)$-holomorphic on interior, where $\kappa$ is some conformal structure on $D_{m+1}$.
\item[$(ii.)$]$Im(u|_{\partial D_{m+1}})\subset\R^3\cup \textcolor{blue}{L^\ast_+T_K}$.
\item[$(iii.)$]$p_0$ is a positive puncture asymptotic to $\bm{a}$. Here we used the fact that $T^\ast\R^3\setminus D^\ast\R^3\cong\R_+\times S^\ast\R^3$.
\item[$(iv.)$]If $k\in\lbrace 1,\dots, m\rbrace$, then there is a point $x_k\in T_K$ on a half-strip neighborhood of $p_k$ such that,
$$u(\tau, t)\rightarrow x_k$$
uniformly in $t$ as $\tau\rightarrow\infty$. We say that $p_k$ is a \textbf{switch at} $x_k$ between two (possibly the same) irreducible components of $\R^3\cup L^\ast_+ T_K$.
\item[$(v.)$] Let us take a standard holomorphic coordinates around $x_k$. Let us identify the quadrant $Q:=\lbrace z\in(\mathbb{C}, i)\, |\, Im(z)\geq 0, Re(z)\geq0, 0<|z|\leq1\rbrace$ with a neighborhood of the puncture $p_k$ in $D_{m+1}$. Then $u$ can be locally written in these coordinates as $u(z)=(u_1(z), u_2(z))\in\mathbb{C}^2\times\mathbb{C}$, where
\begin{equation}\label{eqn_local_switch}
u_1(z)=\sum_{r=0}^\infty b_{r; 1}z^r\quad\hbox{and}\quad u_2(z)=\sum_{r=0}^\infty b_{r; 2}z^{r+2n_k}.
\end{equation}
Here all $b_{r; 1}$ belong to $\R^2$, $b_{r; 2}$ are all in $\R$ or all in $i\R$ and $b_{0; 2}\neq 0$. Finally, $n_k\in\frac{1}{2}\mathbb{N}_+$ and is called the \textbf{winding number at} $p_k$. If $n_k\notin\N_+$, we say that $u$ \textbf{switches sheets at} $p_k$ (via the orientation of $\partial D_{m+1}$). $u_1$ will be called the \textbf{ tangent component of} $u$ and $u_2$ the \textbf{normal component of }$u$.
\item[$(vi.)$] $u$ has a finite Hofer energy.
\end{itemize} 
\end{defn}

\begin{thm}\label{thm_dim_symp} \cite{Cieliebak2016KnotCH, Cieliebak2009CompactnessFH} 
For a generic admissible almost complex structure $J$ the moduli space $\mathcal{M}_{T_K}(\textbf{a}, \bm{n})$ is a manifold of the dimension
$$\dim\mathcal{M}_{T_K}(\textbf{a}, \bm{n})=|\bm{a}|+\sum_{k=1}^m(1-n_k).$$
After a choice of the spin structure on $L^\ast_+T_K$ together with the spin structure on $\R^3$, $\mathcal{M}_{T_K}(\textbf{a}, \bm{n})$ becomes naturally oriented.
\end{thm}

\begin{proof}
Even though the $J$-holomorphic curves from $\mathcal{M}_{T_K}(\bm{a}, \bm{n})$ are now asymptotic rather to the arboreal singularity of $\R^3\cup L^\ast_+T_K$ then to the clean Lagrangian intersection of $\R^3\cup L^\ast T_K$, the arguments work in the same manner as in the case of clean intersections.
\end{proof}

\begin{rem}\label{rem_index_torus} Recall that by Remark \ref{rem_corresp_reeb_bin_chord} the cotangent lift gives us a bijective correspondence between unit speed binormal chords on the torus $T$ and Reeb chords on $\mathcal{L}^\ast T$. This restricts to a bijection between strictly inward-pointing chords on $T$ and Reeb chords on $\mathcal{L}^\ast_+ T$ (once we co-orient $T$ by some tubular neighborhood of a knot, see Chapter \ref{ch:file6}).

Let $\bm{a}$ be a (nondegenerate) Reeb chord on $\mathcal{L}^\ast_+T$ and $\gamma_M$ be the corresponding (nondegenerate) binormal chord on $T$. Then by Theorem \ref{lemma_index_formula} it holds that
$$|\bm{a}|=Ind_{\gamma_M}-1.$$
Since $Ind_{\gamma_M}\in\lbrace0, 1, 2, 3, 4\rbrace$, it follows that $|\bm{a}|\in\lbrace-1, 0, 1, 2, 3\rbrace$.

Let us assume that moreover $T=T_K$ for some $K$. Then $Ind_{\gamma_M}\in\lbrace 1, 2, 3\rbrace$. Indeed $\gamma_M$ is strictly inward-pointing, so $\gamma_M^{op}$ is strictly outward-pointing. And since Morse index does not care about the orientation of $\gamma_M$, one can also apply Lemma \ref{lemma_morse_corresp} for the Morse index computations. Finally, we obtain that $|\bm{a}|\in\lbrace0, 1, 2\rbrace$.
\end{rem}

\begin{rem}In the knot case the moduli space $\mathcal{M}_K(\bm{a}, \bm{n})$ was also introduced, see \cite{Cieliebak2016KnotCH}. Such a space counts punctured $J$-holomorphic discs with boundary on the clean Lagrangian intersection $L^\ast K\cup \R^3$. $\mathcal{M}_K(\bm{a}, \bm{n})$ is also a manifold, but now of the dimension 
$$\dim \mathcal{M}_K(\bm{a}, \bm{n})=\vert\bm{a}\vert+\sum_{k=1}^m(1- 2n_k).$$
Moreover, if $\bm{a}$ is a Reeb chord on $\mathcal{L}^\ast K$ and $\gamma_M$ is the corresponding binormal chord on $K$, then by Theorem \ref{lemma_index_formula} it holds that
$$|\bm{a}|=Ind_{\gamma_M},$$
where $Ind_{\gamma_M}\in\lbrace 0, 1, 2\rbrace$.
\end{rem}

\begin{cor} Let $\bm{a}_0$ and $\bm{a}_\varepsilon$ be the corresponding Reeb chords on $\mathcal{L}^\ast K$ and $\mathcal{L}^\ast T_{K, \varepsilon}$, respectively (they are identified via the canonical contact isotopy). Moreover, assume that their degrees are $1$, then
$$\dim\mathcal{M}_K(\bm{a}_0, \tfrac{1}{2}, \tfrac{1}{2})=0\hbox{ and }\dim\mathcal{M}_{T_{K, \varepsilon}}(\bm{a}_\varepsilon, \tfrac{1}{2}, \tfrac{1}{2})=1.$$
\end{cor}

\begin{rem} For the dimension reasons, the moduli spaces $\mathcal{M}_K(\bm{a}_0, \tfrac{1}{2}, \tfrac{1}{2})$ were used in \cite{Cieliebak2016KnotCH} as suitable ``building blocks'' in the identification of Legendrian contact homology $LCH_0(\mathcal{L}^\ast K)$ with so-called string homology $H^{string}_0(K)$.

In our case, $\dim\mathcal{M}_{T_{K, \varepsilon}}(\bm{a}_\varepsilon, \tfrac{1}{2}, \tfrac{1}{2})$ is $1$-dimensional and not compact. In the following section, we are going to compactify the moduli space and interpret $\overline{\mathcal{M}}_{T_{K, \varepsilon}}(\bm{a}_\varepsilon, \tfrac{1}{2}, \tfrac{1}{2})$ geometrically.
\end{rem}

\section{Geometry of switches} 

\begin{rem}\label{rem_model_for_switches}Let us work a bit on the interplay of the conditions $(iv.)$ and $(v.)$ from Definition \ref{defn_mod_symp} of $\mathcal{M}_{T_K}(\bm{a}, \bm{n})$.

\begin{itemize}
\item Let us assume that $p_k$ is switch at $x_k$ from $L^\ast_+ T_K$ to $\R^3$ with winding number $n_k$.

Let us consider the local coordinates of $u=(u_1, u_2):Q\rightarrow\mathbb{C}^2\times \mathbb{C}$. We would like to inspect the normal component $u_2$. Recall that by Lemma \ref{lemma_coord} at $\text{codom}(u_2)$ it holds that $L^\ast_+ T_K$ is viewed as $i\R_{\geq 0}$ and $\R^3$ as $\R$.

Then, after Schwarz Reflexion principle, $u_2$ will become holomorphic also on $Q\cap(\R\cup i\R)$. After inspecting $u_2$ on real and imaginary numbers, it follows that
$$m_k:=2n_k\hbox{ is odd and }b_{2r+1; 2}=0.$$
Moreover, after sending $z\rightarrow 0$ we see that $b_{0; 2}\in(-1)^{m_k}\R_{>0}$. Specially, 
\begin{equation}\label{eqn_switch_LR}
\begin{cases}
\hbox{if }n_k=\frac{1}{2},\hbox{ then }b_{0; 2}\in\R_{>0},\\
\hbox{if }n_k=\frac{3}{2},\hbox{ then }b_{0; 2}\in\R_{<0}.
\end{cases}
\end{equation}
\item Let us assume that $p_k$ is a switch at $x_k$ from $\R^3$ to $L^\ast_+ T_K$ with winding number $n_k$.

After analysis similar as above, we obtain that
$$2n_k\hbox{ is odd and }b_{2r+1; 2}=0.$$
Also, after sending $z\rightarrow 0$ we get that $b_{0; 2}\in i\R_{>0}$. Trivially, 
\begin{equation}\label{eqn_switch_LR}
\begin{cases}
\hbox{if }n_k=\frac{1}{2},\hbox{ then }b_{0; 2}\in i\R_{>0},\\
\hbox{if }n_k=\frac{3}{2},\hbox{ then }b_{0; 2}\in i\R_{>0}.
\end{cases}
\end{equation}

\item Let us assume that $p_k$ is a switch at $x_k$ from $L^\ast_+ T_K$ back to $L^\ast_+ T_K$ with winding number $n_k$.

Again, after inspecting $u_2$ on real and imaginary numbers, we obtain that
$$2n_k\hbox{ is even and }b_{2r+1; 2}=0.$$
However, after sending $z\rightarrow 0$, it turns out that already $n_k$ has to be even. In addition, $b_{0; 2}\in i\R_{>0}$. Trivially, 
\begin{equation}\label{eqn_switch_LL}
\hbox{if }n_k=2,\hbox{ then }b_{0; 2}\in i\R_{>0}.
\end{equation}
\end{itemize} 
\end{rem}

\begin{notat} Let $u\in \mathcal{M}_{T_K}(\textbf{a}, \textbf{n})$. Then for $i\in \lbrace 1,\dots, m-1\rbrace$ we denote by $c^{u}_{p_i, p_{i+1}}:[a_i, a_{i+1}]\rightarrow T^\ast\R^3$ a smooth path induced by $u$ along the boundary arc between the punctures $p_i$ and $p_{i+1}$. The parametrization of $c^{u}_{p_i, p_{i+1}}$ is chosen such that it matches with $3$-jets of $u$ at local coordinates \ref{eqn_local_switch} near punctures $p_i, p_{i+1}$.
\end{notat}

\begin{thm}\label{thm_compact_strip}
Let $\bm{a}$ be a Reeb chord on $\mathcal{L}^\ast_+T_K$ with $\vert \bm{a}\vert=0$ and $J$ be generic. Then the moduli space $\mathcal{M}_{T_K}(\bm{a}, \tfrac{1}{2}, \tfrac{1}{2})$ admits a natural compactification to a compact oriented $1$-dimensional manifold $\overline{\mathcal{M}}_{T_K}(\bm{a}, \tfrac{1}{2}, \tfrac{1}{2})$ with boundary. The boundary consists of
\begin{itemize}
\item[$(i.)$] moduli spaces $\mathcal{M}_{T_K}(\bm{a}, \tfrac{3}{2}, \tfrac{1}{2})$ and $\mathcal{M}_{T_K}(\bm{a}, \tfrac{1}{2}, \tfrac{3}{2})$, where one of $\tfrac{1}{2}$ switches was replaced by a $\tfrac{3}{2}$ switch.
\item[$(ii.)$] moduli spaces $\mathcal{M}_{T_K}(\bm{a}, 2, \tfrac{1}{2}, \tfrac{1}{2})$ and $\mathcal{M}_{T_K}(\bm{a}, \tfrac{1}{2}, \tfrac{1}{2}, 2)$, where a $2$ switch appeared on one of the boundary components that are mapped to $L^\ast_+ T_K$.
\end{itemize}

In addition, the $J$-holomorphic curves of $\overline{\mathcal{M}}_{T_K}(\bm{a}, \tfrac{1}{2}, \tfrac{1}{2})$ have the following geometry along $\R^3$
\begin{itemize}
\item If $u\in \mathcal{M}_{T_K}(\bm{a}, \tfrac{1}{2}, \tfrac{1}{2})$, then $c^u_{p_1, p_2}$ is strictly outward-pointing from $T_{K}$ at the endpoints.
\item If $u\in \mathcal{M}_{T_K}(\bm{a}, \tfrac{3}{2}, \tfrac{1}{2})$, then $c^u_{p_1, p_2}$ is strictly outward-pointing from $T_{K}$ at $a_2$ and tangent to $T_K$ at $a_1$. Moreover, near $a_1$ $c^u_{p_1, p_2}$ cubically deviates into $\nu T_K$.
\item If $u\in \mathcal{M}_{T_K}(\bm{a}, \tfrac{1}{2}, \tfrac{3}{2})$, then $c^u_{p_1, p_2}$ is strictly outward-pointing from $T_{K}$ at $a_1$ and tangent to $T_K$ at $a_2$. Moreover, near $a_2$ the opposite path $-c^u_{p_1, p_2}$ cubically deviates from $\nu T_K$.
\end{itemize}
\end{thm}

\begin{proof}
The points $(i.), (ii.)$ follow from SFT compactness as in \cite{Cieliebak2009CompactnessFH} and gluing analysis similar to \cite{Cieliebak2016KnotCH, EES05}. For the explicit descriptions of possible local degenerations, see Remarks \ref{rem_bad_switch}, \ref{rem_good_switch}, and \ref{rem_geometry_switches}. In these remarks, we will also show the bullet points of Theorem \ref{thm_compact_strip}.
\end{proof}

\begin{rem}\label{rem_bad_switch}We stress out that the boundary $\partial \overline{\mathcal{M}}_{T_K}(\bm{a}, \tfrac{1}{2}, \tfrac{1}{2})$ from Theorem \ref{thm_compact_strip} does not contain $\mathcal{M}_{T_K}(\bm{a}, 1)$. Even though the virtual dimension of $\mathcal{M}_{T_K}(\bm{a}, 1)$ is $0$, such a space is not a good candidate, since it is an empty space. This is due to our restriction of the conormal bundle $L^\ast T_K$ to $L^\ast_+T_K$, see Figure \ref{figure_bad_switch}. Compare with \cite{Cieliebak2016KnotCH}, where such degenerations were possible and corresponded to the so-called $\delta_N$-string operation.
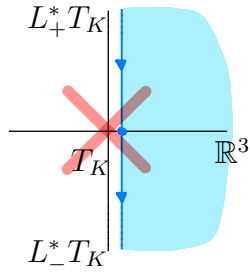
\begin{figure}[!htbp]
\labellist
\pinlabel $\R^3$ at 190 90
\pinlabel ${L^\ast_+T_K}$ at 50 200
\pinlabel ${L^\ast_-T_K}$ at 50 0
\pinlabel $T_K$ at 70 80
\endlabellist
\centering
\includegraphics[scale=0.45]{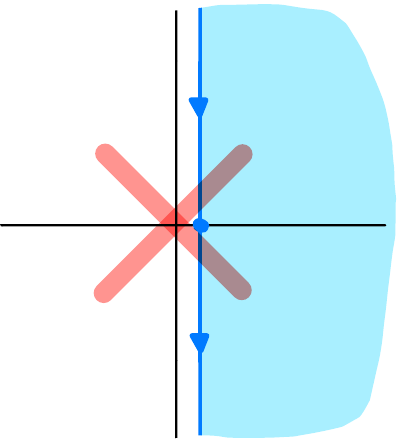}
\vspace{0.3cm}
\caption{A local visualization of the normal component of a potential element $u\in \mathcal{M}_{T_K}(\bm{a}; 1)$ around the puncture $p_1$ with the winding number $1$. But we do not allow $u$ to have a boundary on $L^\ast_-T_K$. So $u$ can not exist.}
\label{figure_bad_switch}
\end{figure}
\end{rem}

\begin{rem}\label{rem_good_switch} We are going to geometrically enlight the boundary phenomenon $\mathcal{M}_{T_K}(\bm{a}, 2, \tfrac{1}{2}, \tfrac{1}{2})\subset\partial \overline{\mathcal{M}}_{T_K}(\bm{a}, \tfrac{1}{2}, \tfrac{1}{2})$. Analogous will hold for $\mathcal{M}_{T_K}(\bm{a}, \tfrac{1}{2}, \tfrac{1}{2}, 2)$. There is a mark point $p\in\partial D_3$ with cyclic order $[p_0, p, p_1]$ and an $1$-parametric family of curves $\widetilde{u}_k\in \mathcal{M}_{T_K}(\bm{a}, \tfrac{1}{2}, \tfrac{1}{2})$ with the following property. With increasing $k$ the points $\widetilde{u}_k(p)$ approach $T_K$ along $L^\ast_+ T_K$. By Remark \ref{rem_model_for_switches}, the resulted curve $u\in \mathcal{M}_{T_K}(\bm{a}, 2, \tfrac{1}{2}, \tfrac{1}{2})$ has the normal component $u_2$ near $p$ of $O(z^4)$. So the family $\widetilde{u}_k$ intersects $T_K$ ``tangentially''. See also Figure \ref{figure_n_string_boundary}.
\begin{figure}[!htbp]
\labellist
\pinlabel $\R^3$ at 217 75
\pinlabel ${L^\ast_+T_K}$ at 60 200
\pinlabel $T_K$ at 142 75
\pinlabel $\textcolor{blue}{\times}$ at 111 125
\endlabellist
\centering
\includegraphics[scale=0.45]{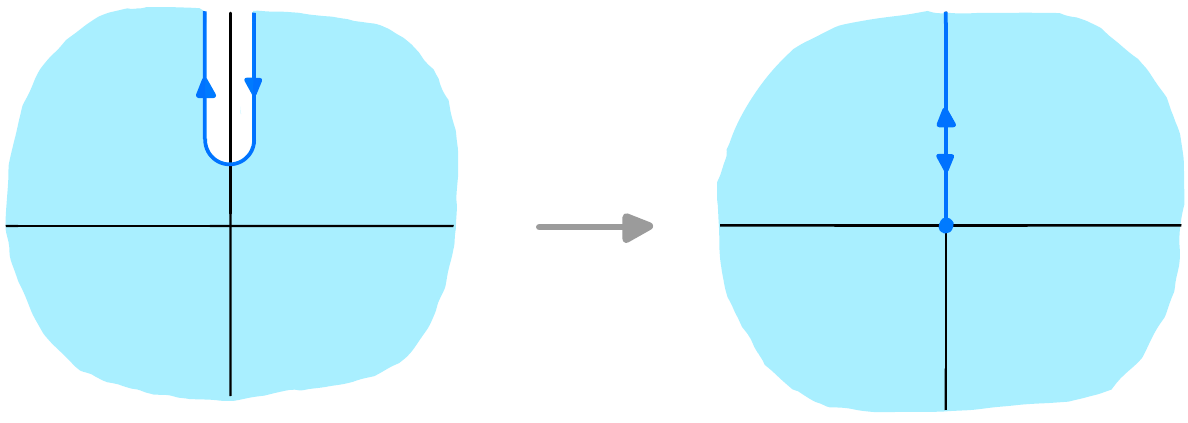}
\vspace{0.3cm}
\caption{\textit{On the left:} $\widetilde{u}_k\in \mathcal{M}_{T_K}(\bm{a}; \tfrac{1}{2}, \tfrac{1}{2})$ around the mark point $p$. \textit{On the right:} $u\in \mathcal{M}_{T_K}(\bm{a}; 2, \tfrac{1}{2}, \tfrac{1}{2})$ around the puncture $p_1$ of the winding number $2$.}
\label{figure_n_string_boundary}
\end{figure}
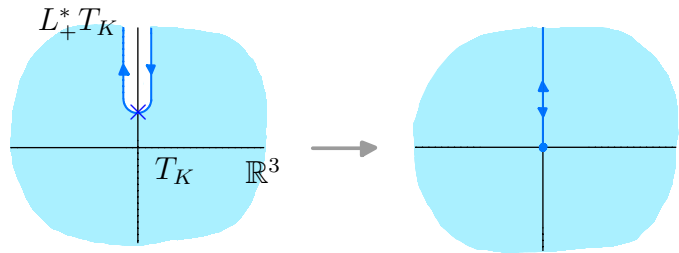
\end{rem}

\begin{rem}\label{rem_geometry_switches}
Now, we are going inspect locally the degenerations of $\overline{\mathcal{M}}_{T_K}(\bm{a}, \tfrac{1}{2}, \tfrac{1}{2})$ as a $\tfrac{1}{2}$ switch changes to a $\tfrac{3}{2}$ switch.

For each switch we write locally $u\in\overline{\mathcal{M}}_{T_K}(\bm{a}, \tfrac{1}{2}, \tfrac{1}{2})$ in the local coordinates $u=(u_1, u_2):Q\rightarrow \mathbb{C}^2\times\mathbb{C}$. Again, we will be interested in the normal components - $u_2$, for which we use Remark \ref{rem_model_for_switches}. 

If $n_1=\tfrac{1}{2}$, then near $p_1$ we can write that $u_2=b_{0; 2}z+O(z^3),$ where $b_{0; 2}\in\R_{>0}$. In particular, $p_1$ switches sheets from $L^\ast_+ T_K$ to $\R^3$. By the order of the leading term of $u_2$, it follows that $c^u_{p_1, p_2}$ is not tangent to $T_K$ at $c^u_{p_1, p_2}(a_1)$. Next, recall that we have chosen the co-orientation of $T_K$ by the normal vector pointing inside of $\nu K$. Hence, since $b_{0; 2}\in\R_{>0}$, $c^u_{p_1, p_2}$ is strictly outward-pointing from $T_K$ at $c^u_{p_1, p_2}(a_1)$, see Figure \ref{figure_normal_wind}.

If $n_2=\tfrac{1}{2}$, then near $p_2$ we can write that $u_2=b_{0; 2}z+O(z^3),$ where $b_{0; 2}\in i\R_{>0}$. By the same argument it follows that $p_2$ switches sheets from $\R^3$ to $L^\ast_+ T_K$ and $c^u_{p_1, p_2}$ is strictly outward-pointing from $T_K$ at $c^u_{p_1, p_2}(a_2)$, see Figure \ref{figure_normal_wind}.
\begin{figure}[!htbp]
\labellist
\pinlabel $\R^3$ at 250 50
\pinlabel ${L^\ast_+T_K}$ at 80 170
\pinlabel $T_K$ at 90 50
\endlabellist
\centering
\includegraphics[scale=0.45]{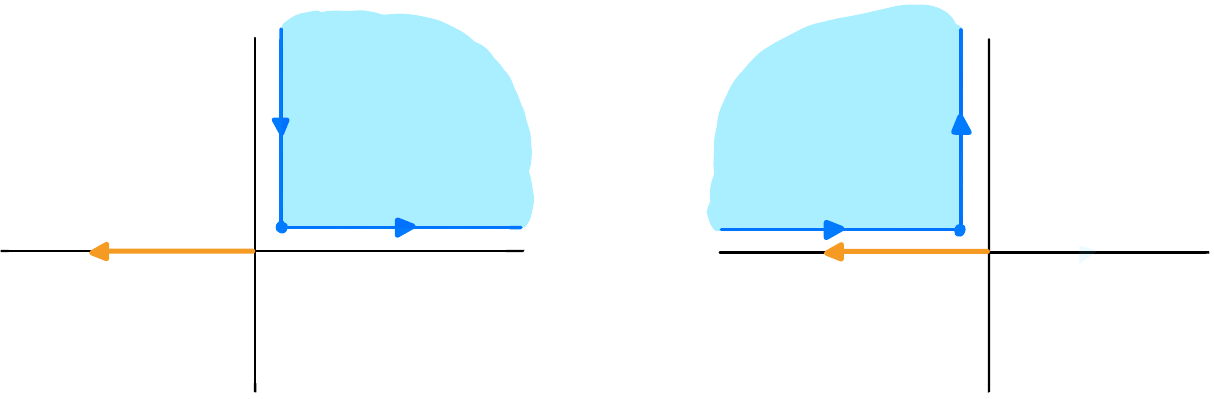}
\vspace{0.3cm}
\caption{\textit{On the left:} A local visualization of the normal component of $u\in \mathcal{M}_{T_K}(\bm{a}; \tfrac{1}{2}, \tfrac{1}{2})$ around the puncture $p_1$ of the winding number $\tfrac{1}{2}$. The blue arrows depict the boundary orientation of $u$, and the orange arrow represents the positive co-orientation of $T_K$. \textit{On the right:} $u\in \mathcal{M}_{T_K}(\bm{a}; \tfrac{1}{2}, \tfrac{1}{2})$ around the puncture $p_2$ of the winding number $\tfrac{1}{2}$.}
\label{figure_normal_wind}
\end{figure}

Now, the winding numbers of the value $\tfrac{3}{2}$. In these cases $u_2=O(z^3)$, so $c^u_{p_1, p_2}$ is tangent to $T_K$ at given endpoints. More precisely, if $n_1=\tfrac{3}{2}$, then near $p_1$ we can write that $u_2=b_{0; 2}z^3+O(z^5),$ where $b_{0; 2}\in \R_{<0}$, And if $n_2=\tfrac{3}{2}$, then near $p_2$ we can write that $u_2=b_{0; 2}z^3+O(z^5),$ where $b_{0; 2}\in i\R_{>0}$. Moreover, from the signs of coefficients of $O(z^3)$ terms we obtain that $c^u_{p_1, p_2}$ deviates from $T_K$ as proposed in Theorem \ref{thm_compact_strip}. See also Figure \ref{figure_tangent_wind}.
\begin{figure}[!htbp]
\labellist
\pinlabel $\R^3$ at 250 90
\pinlabel ${L^\ast_+T_K}$ at 90 200
\pinlabel $T_K$ at 160 65
\endlabellist
\centering
\includegraphics[scale=0.45]{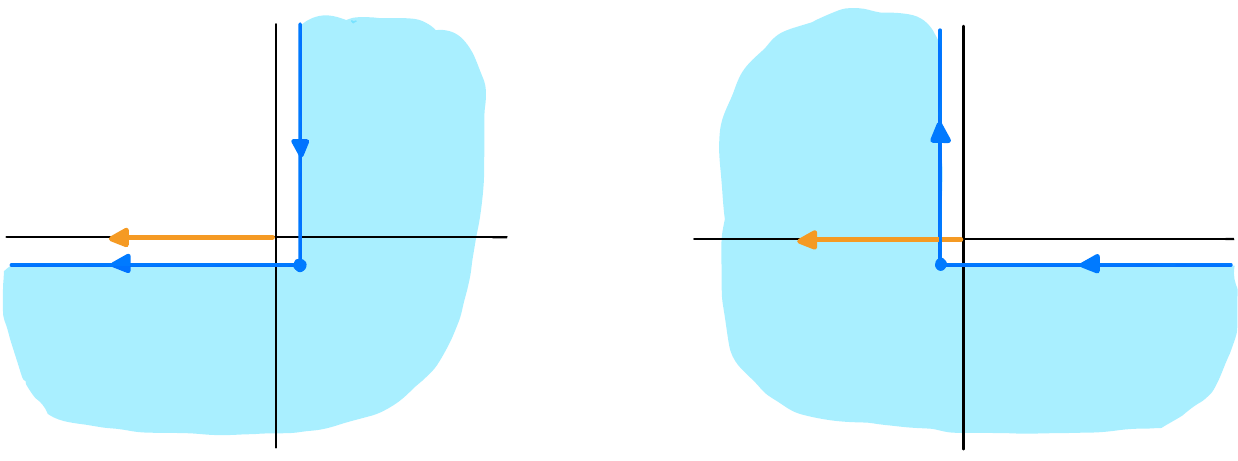}
\vspace{0.3cm}
\caption{\textit{On the left:} $u\in \mathcal{M}_{T_K}(\bm{a}; \tfrac{3}{2}, \tfrac{1}{2})$ around the puncture $p_1$ with winding number $\tfrac{3}{2}$. \textit{On the right:} $u\in \mathcal{M}_{T_K}(\bm{a}; \tfrac{1}{2}, \tfrac{3}{2})$ around the puncture $p_2$ with winding number $\tfrac{3}{2}$.}
\label{figure_tangent_wind}
\end{figure}

Now, let $u\in \mathcal{M}_{T_K}(\bm{a}, \tfrac{3}{2}, \tfrac{1}{2})$ and $\lbrace \widetilde{u}_k\rbrace_k$ be a sequence of elements of $\mathcal{M}_{T_K}(\bm{a}, \tfrac{1}{2}, \tfrac{1}{2})$ that converge to $u$. We would like to geometrically describe how the holomorphic curves are deformed. By the above, we know that the normal components of $\widetilde{u}_k$ are of the form $a_k z+O(z^3)$, where $a_k\rightarrow 0$ as $k\rightarrow \infty$. This means that for $k\gg0$ large each $c^{\widetilde{u_k}}_{p_1, p_2}$ intersects $T_K$ additionally at some time $a_1+\varepsilon_k$. Since $\varepsilon_k\rightarrow 0$ as $k\rightarrow\infty$, we obtain certain vanishing “spikes”. See also Figure \ref{figure_spike}. In addition, these vanishing $a_k$ coefficients will contribute as a gluing parameter which induces a parametrization of $\overline{\mathcal{M}}_{T_K}(\bm{a}, \tfrac{3}{2}, \tfrac{1}{2})$ near $\mathcal{M}_{T_K}(\bm{a}, \tfrac{1}{2}, \tfrac{1}{2})$. Similarly, we can describe $\mathcal{M}_{T_K}(\bm{a}, \tfrac{3}{2}, \tfrac{1}{2})$, see also Figure \ref{figure_spike}.
\begin{figure}[!htbp]
\labellist
\pinlabel $\R^3$ at 255 90
\pinlabel ${L^\ast_+T_K}$ at 90 200
\pinlabel $T_K$ at 105 72
\endlabellist
\centering
\includegraphics[scale=0.45]{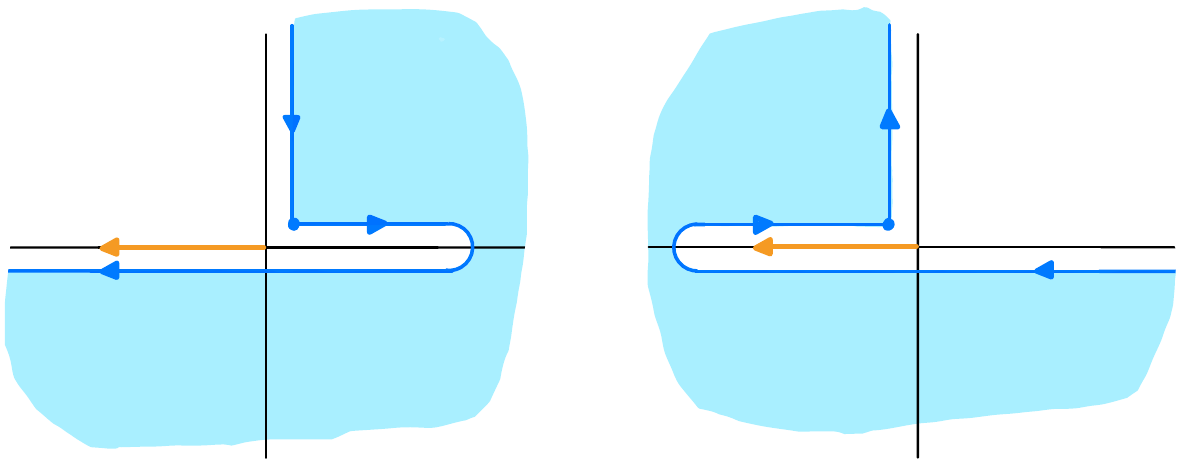}
\vspace{0.3cm}
\caption{\textit{On the left:} A local visualization of the normal component of $\widetilde{u}_k\in \mathcal{M}_{T_K}(\bm{a}; \tfrac{1}{2}, \tfrac{1}{2})$ around the puncture $p_1$ with the winding number $\tfrac{1}{2}$. Note that due to the spike the curve $c^{\widetilde{u_k}}_{p_1, p_2}$ intersects $T_K$ again at some time $\varepsilon_k$. Also as $k\rightarrow \infty$ the spike vanishes and the winding at $p_1$ becomes $\tfrac{3}{2}$. \textit{On the right:} $\widetilde{u}_k\in \mathcal{M}_{T_K}(\bm{a}; \tfrac{1}{2}, \tfrac{1}{2})$ around the puncture $p_2$ with the winding number $\tfrac{1}{2}$ which becomes $\tfrac{3}{2}$ as $k\rightarrow\infty$.}
\label{figure_spike}
\end{figure}
\end{rem}

\begin{conj}\label{cor_moduli_boundary_phenomena}Let $\bm{a}$ be a Reeb chord on $\mathcal{L}^\ast_+ T_K$ and $J$ be generic. Then the moduli space $\mathcal{M}_{T_K}(\bm{a}; \tfrac{1}{2},\dots, \tfrac{1}{2})$ with $2\ell$ switches admits a natural compactification to a compact oriented $(\vert\bm{a}\vert+\ell)$-dimensional manifold $\overline{\mathcal{M}}_{T_K}(\bm{a}, \tfrac{1}{2},\dots, \tfrac{1}{2})$ with corners. The codimension $1$ strata consist of:
\begin{itemize}
\item[$(i.)$] moduli spaces $\mathcal{M}_{T_K}(\bm{a}, \tfrac{1}{2},\dots, \tfrac{1}{2}, \tfrac{3}{2}, \tfrac{1}{2},\dots \tfrac{1}{2})$ with $2\ell$ switches, where one $\tfrac{1}{2}$ switch was replaced by a $\tfrac{3}{2}$ switch.
\item[$(ii.)$] moduli spaces $\mathcal{M}_{T_K}(\bm{a}, \tfrac{1}{2},\dots, \tfrac{1}{2}, 2, \tfrac{1}{2},\dots \tfrac{1}{2})$ with $2\ell+1$ switches, where we added one $2$ switch from $L^\ast_+ T_K$ back to $L^\ast_+ T_K$ (i.e. not from $\R^3$ to $\R^3$).
\item[$(iii.)$] products of moduli spaces
$$\mathcal{M}^{sy}_{T_K}(\bm{a}, \bm{b}_1,\dots, \bm{b}_m)/\R\times\prod_{i=1}^m\mathcal{M}_{T_K}(\bm{b}_i; \tfrac{1}{2},\dots, \tfrac{1}{2}),$$
where the Reeb chords $\bm{b}_i$ have total degree $\vert \bm{a}\vert-1$ and the moduli spaces $\mathcal{M}_{T_K}(\bm{b}_i; \tfrac{1}{2},\dots, \tfrac{1}{2})$ have together $2\ell$ switches.
\item[$(iv.)$] moduli spaces $\mathcal{M}_{T_K}(\bm{a}, \tfrac{1}{2},\dots, \tfrac{1}{2}, 1, \tfrac{1}{2},\dots \tfrac{1}{2})$ with $2\ell-1$ switches, where two consecutive switches $\tfrac{1}{2}$ were replaced by a single $1$ switch, provided that the first replaced $\tfrac{1}{2}$ switch was switching from $\R^3$ to $L^\ast_+ T_K$. See also Figure \ref{figure_string_degeneration}.
\end{itemize}
The codimension $2$ strata consist of corners given by the combinations of phenomena $(i.)-(iv.)$ and
\begin{itemize}
\item[$(v.)$] moduli spaces $\mathcal{M}_{T_K}(\bm{a}, \tfrac{1}{2},\dots, \tfrac{1}{2}, \tfrac{5}{2}, \tfrac{1}{2},\dots \tfrac{1}{2})$ with $2\ell$ switches, where one $\tfrac{1}{2}$ switch was replaced by a $\tfrac{5}{2}$ switch.
\item[$(vi.)$] moduli spaces $\mathcal{M}_{T_K}(\bm{a}, \tfrac{1}{2},\dots, \tfrac{1}{2}, 2, \tfrac{1}{2},\dots \tfrac{1}{2})$ with $2\ell-1$ switches, where two consecutive switches $\tfrac{1}{2}$ were replaced by a single $2$ switch, provided that the first replaced $\tfrac{1}{2}$ switch was switching from $L^\ast_+ T_K$ to $\R^3$. See also Figure \ref{figure_string_degeneration}.
\end{itemize}
\end{conj}

\begin{pproof}
The boundary codimension $1$ configurations $(i)-(ii.)$ were already discussed in Theorem \ref{thm_compact_strip}, and the configuration $(iii.)$ is just a standard breaking phenomenon from Symplectic field theory, see \cite{Cieliebak2016KnotCH, Bourgeois_2003}. The boundary codimension $1$ configuration $(iv.)$ is given by the collapse of two switches into a single one, which is realized by a vanishing “spike” in the conormal direction, see Figure \ref{figure_string_degeneration}. Note that this degeneration is not violating the boundary conditions coming from the arboreal singularity. For the same phenomenon, we refer to \cite{Cieliebak2016KnotCH}.
\begin{figure}[!htbp]
\labellist
\pinlabel $\R^3$ at 235 60
\pinlabel ${L^\ast_+T_K}$ at 85 200
\pinlabel $T_K$ at 95 60
\endlabellist
\centering
\includegraphics[scale=0.47]{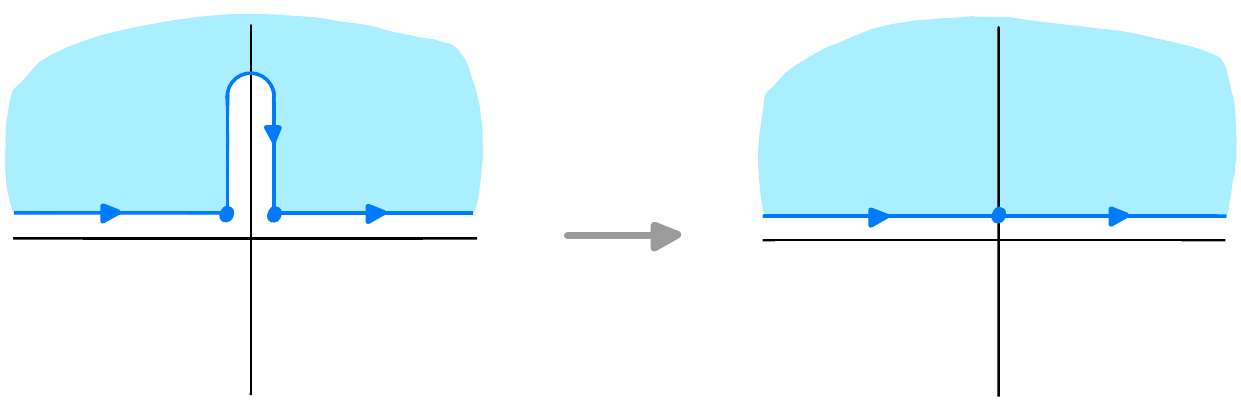}
\caption{\textit{Configuration $(iv.)$:} vanishing conormal “spike”.}
\label{figure_string_degeneration}
\end{figure}

Now, we are going to inspect the codimension $2$ phenomenon $(vi.)$. The configuration $(vi.)$ will appear as an intersection of two configurations $(i.)$. We take an arc from $\partial D_{2\ell+1}$ connecting punctures $p_{i}, p_{i+1}$ of winding numbers $\tfrac{1}{2}, \tfrac{1}{2}$. Then, while collapsing the arc, we obtain a $2$-parametric family of parameters which creates a corner.

Let us take a local model of the normal component of $J$-holomorphic curves near punctures $p_i, p_{i+1}$.  Since we would like to work with two punctures rather then the quadrant $Q$, we take as a domain a punctured half disc $H:=\lbrace \zeta\in(\mathbb{C}, i)\,\vert\, Im(\zeta)\geq 0, 0<\vert \zeta\vert\leq 1, \zeta\neq -\varepsilon\rbrace$. Recall that we can translate the points from $H$ to $Q$ by $\zeta\mapsto\zeta^{1/2}$. Then up to $O(\zeta^{5/2})$, we can locally model the normal components of the $J$-holomorphic curves on the holomorphic functions
$$f_{\varepsilon, \delta}(\zeta):=(\zeta+\varepsilon)^{1/2}(\zeta+\delta)i\zeta^{1/2},$$
where $0\leq\delta\leq\varepsilon$. The function $f_{\varepsilon, \delta}$ has the following properties. The zeros are $-\varepsilon, -\delta, 0$. Along $(-1, -\varepsilon)\cup(0, 1)$, the real part of $f_{\varepsilon, \delta}$ vanishes and the imaginary part matches the conormal boundary conditions. Analogously, along $(-\varepsilon, 0)$ the imaginary part of $f_{\varepsilon, \delta}$ vanishes.

Geometrically, observe that the zero at $-\delta$ corresponds to the intersection of two small spikes with $T_K$ - one spike for each puncture. Also, as $-\delta$ moves to $-\varepsilon$ the winding at $p_i$ changes from $\tfrac{1}{2}$ to $\tfrac{3}{2}$ and similarly for the puncture $p_{i+1}$. See also Figure \ref{figure_sing_holom}. Hence, we can consider as gluing parameters $\varepsilon, \delta$ with $0\leq\delta\leq\varepsilon$.
\begin{figure}[!htbp]
\labellist
\pinlabel $\R^3$ at 250 65
\pinlabel ${L^\ast_+T_K}$ at 80 200
\pinlabel $T_K$ at 165 60
\pinlabel $\textcolor{blue}{\bullet}$ at 478 91
\pinlabel $\textcolor{blue}{\bullet}$ at 122 105
\pinlabel $\textcolor{blue}{\bullet}$ at 148 81
\endlabellist
\centering
\includegraphics[scale=0.45]{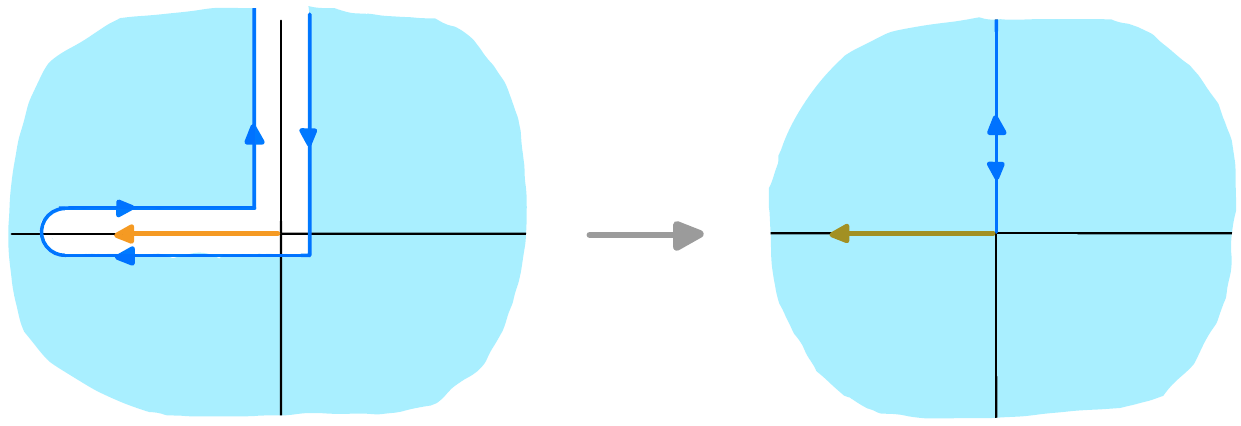}
\caption{\textit{On the left:} $u\in \mathcal{M}_{T_K}(\bm{a}; \tfrac{3}{2}, \tfrac{1}{2})$ around collapsing punctures $p_1$ and $p_2$. \textit{On the right:} $u\in \mathcal{M}_{T_K}(\bm{a}; 2)$ around the puncture $p_1$ with winding number $2$.}
\label{figure_sing_holom}
\end{figure}
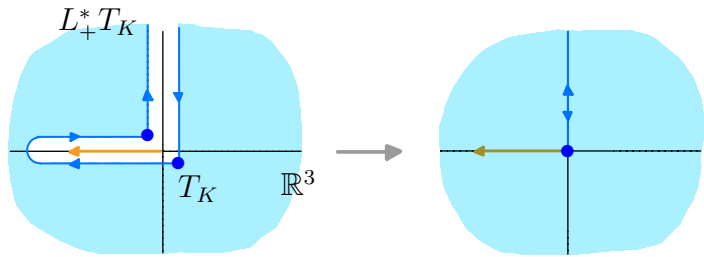

Now, let us consider the codimension $2$ phenomenon $(v.)$. The configuration $(v.)$ can be created from intersection of configurations $(i.)$ and $(ii.)$. I.e. $(v.)$ can appear as a collapse of punctures with winding $2$ and $\tfrac{1}{2}$, or one puncture can change its winding from $\tfrac{3}{2}$ to $\tfrac{5}{2},$ see Figure \ref{figure_5_2_weight}. Since both deformations correspond to different arcs of the punctured disc, we obtain a corner.
\begin{figure}[!htbp]
\labellist
\pinlabel $\R^3$ at 250 65
\pinlabel ${L^\ast_+T_K}$ at 90 200
\pinlabel $T_K$ at 163 60
\pinlabel $\textcolor{blue}{\bullet}$ at 471 99
\pinlabel $\textcolor{blue}{\bullet}$ at 139 79
\endlabellist
\centering
\includegraphics[scale=0.45]{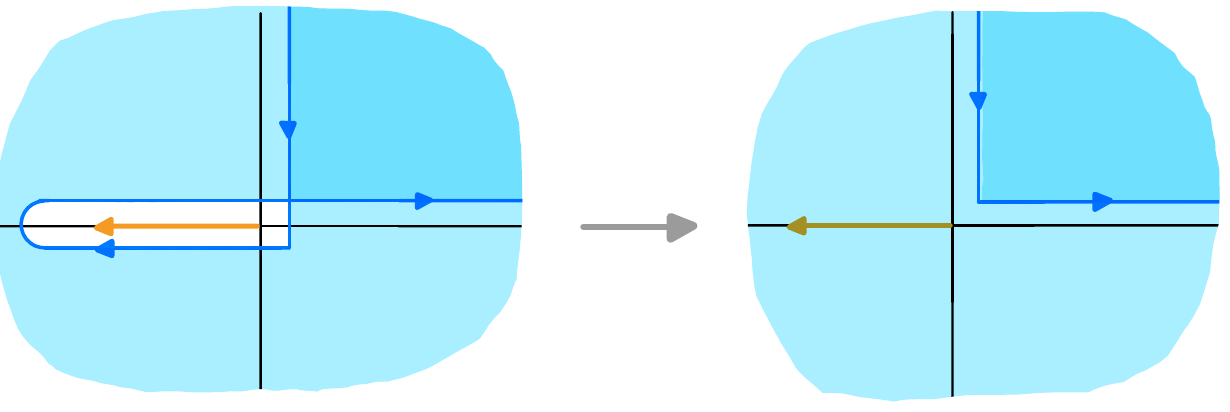}
\vspace{0.3cm}
\caption{\textit{On the left:} $u\in \mathcal{M}_{T_K}(\bm{a}; \tfrac{3}{2}, \tfrac{1}{2})$ around the puncture $p_1$ with winding number $\tfrac{3}{2}$. \textit{On the right:} $u\in \mathcal{M}_{T_K}(\bm{a}; \tfrac{5}{2}, \tfrac{1}{2})$ around the puncture $p_1$ with winding number $\tfrac{5}{2}$.}
\label{figure_5_2_weight}
\end{figure}
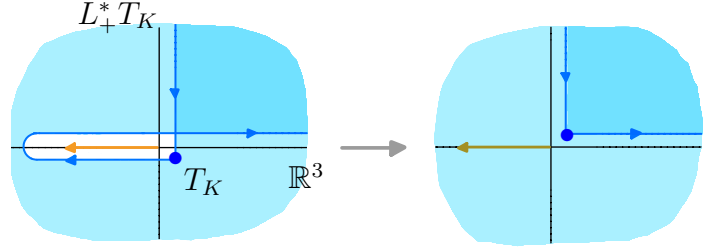
\end{pproof}

\begin{lemma}There is a constant $k_{sw}\in \mathbb{N}$ such that every $\mathcal{M}_{T_K}(\bm{a}, \bm{n})$ is empty, if the number of switches is bigger then $k_{sw}$.
\end{lemma}

\begin{proof}
The lemma follows from similar argument as \cite[thm 1.2]{Cieliebak2009CompactnessFH}.
\end{proof}

\section{Broken strings and chain maps}
In this section, we outline a chain map between the chain complexes for Legendrian contact homology and so-called String homology. Then, we discuss the relation to Morse model of the cord algebra. \textit{The precise proofs of the statements will remain to be done.}

Motivated by broken strings from \cite{Cieliebak2016KnotCH, Ekholm_2017}, we are going to construct a degree zero string homology for $T_{K, \varepsilon}$. See also \cite{okamoto2022toward}.

\begin{notat} In order to agree with the notation from \cite{Cieliebak2016KnotCH, Ekholm_2017}, the ambient space $\R^3$ will be denoted by $Q$.

Next, let us consider an $(\varepsilon/2)$-radius tubular neighborhood $\nu_{\varepsilon/2}T_{K, \varepsilon}$. Recall that the positive co-orientation $T_{K, \varepsilon}$ is given by vectors pointing \textit{inside} of $\nu_\varepsilon K$. Then we put
$$N_+:=\nu_\varepsilon K\cap \nu_{\varepsilon/2}T_{K, \varepsilon}.$$ 
Let us fix $x_0\in(\partial N_+\setminus T_{K, \varepsilon})$ and $v_0\in T_{x_0}N_+$. 
\end{notat}

\begin{figure}[!htbp]
\labellist
\pinlabel $T_{K, \varepsilon}$ at 0 800
\endlabellist
\centering
\includegraphics[scale=0.8]{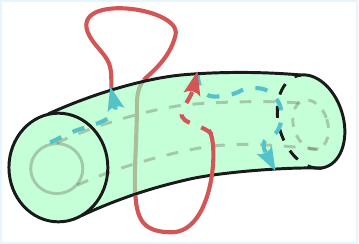}
\vspace{0.3cm}
\caption{A visualization of a broken string on $T_{K, \varepsilon}$ that conists of \textcolor{red}{$Q$-strings} and \textcolor{teal}{$N$-strings}. Note the matching jet condions of $N/Q$-strings at their endpoints.}
\label{figure_strings2}
\end{figure}

\begin{defn}\label{defn_broken_string} A \textbf{broken string with $\ell$ $Q$-strings} on $T_{K, \varepsilon}$ is a tuplet $s=(a_1,\dots, a_{2\ell+1}; s_1, \dots, s_{2\ell+1})$ consisting of real numbers $0=a_0<\dots<a_{2\ell+1}$ and non-constant $C^\infty$-maps 
$$s_{2i+1}:[a_{2i}, a_{2i+1}]\rightarrow \textcolor{blue}{N_+}, \quad s_{2i}:[a_{2i-1}, a_{2i}]\rightarrow Q$$
such that the following holds
\begin{itemize}
\item[$(i.)$]$s_0(0)=s_{2\ell+1}(a_{2\ell+1})=x_0$.
\item[$(ii.)$] for $j=1,\dots, 2\ell$, $s_j(a_j)=s_{j+1}(a_j)\in T_{K, \varepsilon}$.
\item[$(iii.)$] Let $\sigma_j$ denotes the normal component of $s_j$ near its endpoints. Then for every $2k\in\mathbb{N}$ the derivatives $\sigma_j^{(2k)}$ at the endpoints vanish. For $i=1,\dots, \ell$ and every $k\in \mathbb{N}$ it holds that
\begin{align*}
\sigma^{(2k-1)}_{2i}(a_{2i})=(-1^k)\sigma^{(2k-1)}_{2i+1}(a_{2i}),\\
\sigma^{(2k-1)}_{2i-1}(a_{2i-1})=(-1^{k+1})\sigma^{(2k-1)}_{2i}(a_{2i-1}).
\end{align*}
\end{itemize}
Paths $s_{2i}$ are called \textbf{$Q$-strings} and paths $s_{2i+1}$ are called \textbf{$N$-strings}, see Figure \ref{figure_strings2}. We equip $\Sigma_\ell^\circ$; the set of broken strings with $\ell$ $Q$-strings, with $C^\infty$ topology, and by completion extend it by also allowing constant $0$ time $Q/N$-strings; \textbf{trivial $Q/N$-strings}. The resulting space of all strings will be denoted by $\Sigma_\ell$, and we put $\Sigma_\ell^\Delta:=\Sigma_\ell\setminus \Sigma^\circ_\ell$. Moreover, $\Sigma_\ell^{Q_n}$ will consist of broken strings with at least $n>0$ trivial $Q$-strings. The space of all broken strings will be denoted by $\Sigma$.
\end{defn}

\begin{rem}$Q$-strings are outward-pointing from $T_{K, \varepsilon}$. This follows from the condition $(iii.)$ in Definition \ref{defn_broken_string} and the geometry of $N_+$.
\end{rem}

\begin{defn} A broken string from $\Sigma_\ell$ is called a string with a \textbf{$Q$-tangency} if it contains a nontrivial $Q$-string $s_{2i}$ with a vanishing first normal derivative $\sigma^{(1)}_{2i}$ at some endpoint $p$, and $s_{2i}$ can be potentially concatenated at $p$ with some non-trivial local $N$-string. (The second condition just dictates the sign of the first nonzero derivative $\sigma^{(2k-1)}(p)$).

The subset of all of these strings is denoted by $\partial^Q\Sigma_\ell$. And analogously we define $\partial^Q\Sigma^{Q_n}_\ell$.
\end{defn}

\begin{defn} A broken string from $\Sigma_\ell$ is called a string with a \textbf{$N$-tangency} if it contains a nontrivial $N$-string $s_{2i-1}$ touching $T_{K, \varepsilon}$ at some interior point $p$.

The subset of all of these strings is denoted by $\partial^N\Sigma_\ell$. And analogously we define $\partial^N\Sigma^{Q_n}_\ell$.
\end{defn}

\begin{rem_not}Let $\ev_{Q, 2i}:\Sigma_\ell\times(0, 1)\rightarrow Q$ be an evaluation map on interiors of $Q$-strings $s_{2i}$. In more detail,
$$\ev_{Q, 2i}:(s, t)\mapsto s_{2i}((1-t)a_{2i-1}+t a_{2i})).$$
Then, we define a map $F_{Q, 2i}:\Sigma_\ell\times(0, 1)\times T_{K, \varepsilon}\rightarrow Q\times \R$ by
$$(s, t, x)\mapsto\big(\ev_{Q, 2i}(s, t)-x, \langle\partial_t\ev_{Q, 2i}(s, t), v(x)\rangle\big),$$
where $v(x)$ is a \textit{negatively} co-oriented (normalized) normal vector at $x\in T_{K, \varepsilon}$. Compare $F_Q$ also with maps for outward-pointing chords from Definition \ref{defn_funct_out}.

Let $d\geq0$. We also introduce switch evaluation maps $T_{Q, 2i}^d:\Sigma_\ell\rightarrow \R^d$ by
$$s\mapsto\big(\sigma_{2i}^{(2k-1)}(a_{2i})\big)_{1\leq 2k-1\leq d}.$$

Analogously, we introduce also $\ev_{N, 2i-1}$ and $T_{N, 2i-1}^d$.
\end{rem_not}

\begin{defn}\label{defn_gen_chain} Let $d\geq0$ and $\Delta^d=\lbrace(\lambda_1,\dots, \lambda_d)\in\R^d\,\vert\,\lambda_k\geq0, \sum_k \lambda_k\leq1\rbrace.$ A \textbf{weakly generic singular $d$-chain} in $\Sigma_\ell^\circ$ is a smooth map $S:\Delta^d\rightarrow\Sigma^\circ_\ell$ such that
\begin{itemize}
\item[$(i.)$] on $Q$-strings the maps $F_Q\circ S$ and $T^{d+1}_Q\circ S$ are stratum transverse to $0\times[0, \infty)$, respectively jet transverse to $0$. 
\item[$(ii.)$] on $N$-strings the maps $T^{d+1}_N\circ S$ are jet transverse to $0$, $\ev_{N}\circ S(Int(\Delta^d))\cap T_{K, \varepsilon}=\emptyset$ and $\overline{(\ev_{N}\circ S\vert_{\partial\Delta^d})^{-1}(T_{K, \varepsilon})}$ are $(d-1)$-dimensional manifolds with corners. Moreover, the connected components of $\overline{(\ev_{N}\circ S\vert_{\partial\Delta^d})^{-1}(T_{K, \varepsilon})}$ projects to $\partial\Delta^d$ as immersed manifolds.
\end{itemize} 

Finally, on weakly generic singular chains, we introduce the following manifolds with corners
 $$M^d_{\delta_Q(S)} :=\overline{(F_Q\circ S)^{-1}(0\times[0, \infty))}\hbox{ and }M^d_{\delta_N(S)}:=\overline{(\ev_{N}\circ S)^{-1}(T_{K, \varepsilon})},$$ 
In particular, $\dim M^d_{\delta_Q(S)}=\textcolor{blue}{d}$ and $\dim M^d_{\delta_N(S)}=\textcolor{blue}{d-1}$.
\end{defn}

\begin{rem} The condition $(ii.)$ from Definition \ref{defn_gen_chain} looks a bit artificial and will be used in the definition of the String homology, Definition \ref{defn_string_homology}. An obvious way to avoid $(ii.)$ is to forget $N$-strings and work only with $Q$-strings with lengths bounded from below by some small positive number. But, as a consequence, we will then aim for constructing a string homology with coefficients only in $\mathbb{Z}$ and not in $\mathbb{Z}[\pi_1(T_{K, \varepsilon})]$.
\end{rem}

\begin{defn}Using the conditions $(i.), (ii.)$ from Definition \ref{defn_gen_chain} we define also weakly generic singular chains for $\partial^Q\Sigma^\circ_\ell$ and $\partial^N\Sigma^\circ_\ell$. 

Weakly generic singular chains for broken strings $\Sigma_\ell$ with a single trivial $Q/N$-string are introduced as follows. We again use the conditions $(i.), (ii.)$, but not for the trivial string, and also the corner evaluation maps $T^{d+1}$ are not applied for the endpoints surrounding the trivial string.

This induces weakly generic singular $d$-chains on the whole $\Sigma$. We also extend the definitions of $M^d_{\delta_Q(S)}$ and $M^d_{\delta_N(S)}$.
\end{defn}

\begin{rem} We can view the connected components of $M^d_{\delta_Q(S)}$ as weakly generic singular $d$-chains inside $S$, let us denote them by $S_Q$. Analogously, $M^d_{\delta_N(S)}$ will contribute by weakly generic singular $(d-1)$-chains inside $S$, let us denote them by $S_N$.
\end{rem}

\begin{defn}A \textbf{string coproducts $\delta_Q$ and $\delta_N$} are defined on a weakly generic singular $d$-chain $S$ by sums of $S_Q/S_N$ with inserted trivial $N/Q$-strings at $M^d_{\delta_Q(S)}$ and $M^d_{\delta_N(S)}$, respectively.
\end{defn}

\begin{figure}[!htbp]
\labellist
\pinlabel $\longrightarrow$ at 150 767
\pinlabel $\longrightarrow$ at 460 767
\pinlabel $\delta_Q$ at 145 780
\pinlabel $\delta_N$ at 458 780
\endlabellist
\centering
\includegraphics[scale=0.74]{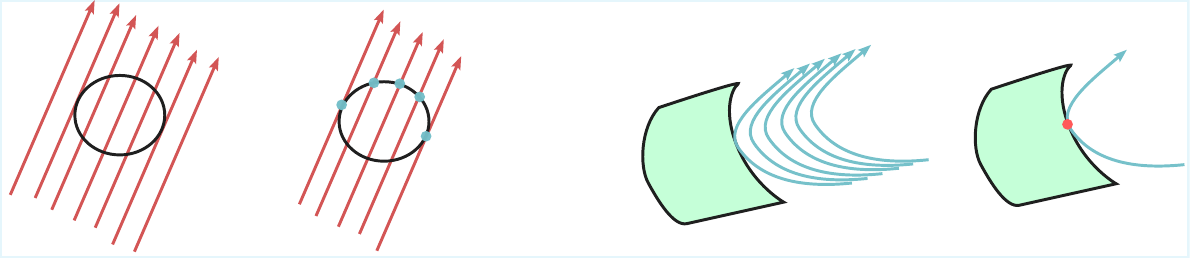}
\vspace{-0.2cm}
\caption{String operations $\delta_Q$ and $\delta_N$.}
\label{figure_coproducts}
\end{figure}
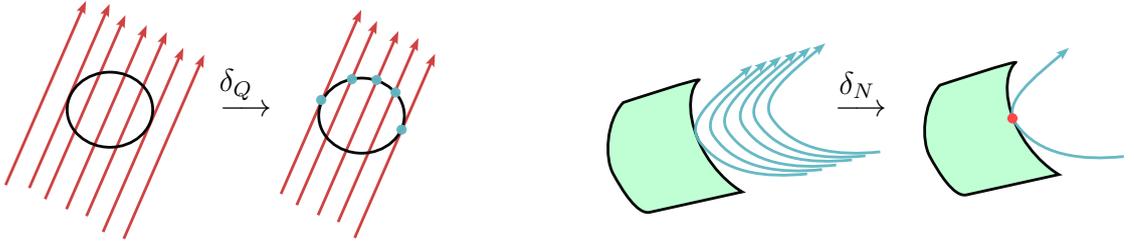

\begin{defn}\label{defn_string_homology}By $C^{sing}_{k+\ell}(\Sigma_\ell, \partial^Q\Sigma_\ell\cup\partial^N\Sigma_\ell)$ we denote a free $\mathbb{Z}$-module generated by weakly generic relative singular $(k+\ell)$-simplices in $\Sigma_\ell$. Similarly, by $C^{sing}_{k+\ell-n}(\Sigma_\ell^{Q_n})$ we understand a free $\mathbb{Z}$-module generated by weakly generic singular $(k+\ell-n)$-simplices in the space of broken strings with $n$ constant and $n-\ell$ non-constant $Q$-strings. We put 
$$C^{string}_k(\Sigma):=\bigoplus_{\ell=0}^{\infty}\big[C^{sing}_{k+\ell}(\Sigma_\ell, \partial^Q\Sigma_\ell\cup\partial^N\Sigma_\ell)\oplus\bigoplus_{n=0}^{\ell}C^{sing}_{k+\ell-n}(\Sigma_\ell^{Q_n}, \partial^Q\Sigma_\ell^{Q_n}\cup\partial^N\Sigma_\ell^{Q_n})\big].$$

Concatenations at the base point induce on $C^{string}_k(\Sigma)$ the structure of a noncommutative (but strictly associative) $\mathbb{Z}$-algebra. 

The singular boundary operator $\partial^{sing}$ and the string coproducts $\delta_Q, \delta_N$ induce a linear map of degree $-1$
$$\partial^{sing}+\delta_Q+\delta_N: C^{string}_1(\Sigma)\rightarrow C^{string}_0(\Sigma),$$
see also Figure \ref{figure_coproducts}.

We define the degree zero \textbf{string homology} $H^{string}_0(T_{K, \varepsilon})$ as
$$H^{string}_0(T_{K, \varepsilon})=C^{string}_0(\Sigma)/\mathcal{I}^{string},$$
where $\mathcal{I}^{string}$ is a two-sided ideal of $C^{string}_0(\Sigma)$ generated by the image of $\partial^{sing}+\delta_Q+\delta_N$.
\end{defn}

\begin{rem} Let us heuristically relate the skein relations arising from $\mathcal{I}^{string}$ to $Cord^M(T_{K, \varepsilon}; \mathbb{Z}[\lambda^{\pm1}, \mu^{\pm1}])$. Hence let us consider a $1$-dimensional family $s^t$ of broken strings in $C_1^{string}(\Sigma)$. Then from $(\partial^{sing}+\delta_Q+\delta_N)(s^t)=0$ we obtain the following skein relations:
\begin{figure}[!htbp]
\labellist
\pinlabel $(i.)\quad 0=$ at 0 750
\pinlabel $(ii.)\quad 0=$ at 0 630
\pinlabel $-$ at 200 750
\pinlabel $-$ at 200 630
\pinlabel $+$ at 370 750
\pinlabel $+$ at 370 630
\endlabellist
\centering
\includegraphics[scale=0.5]{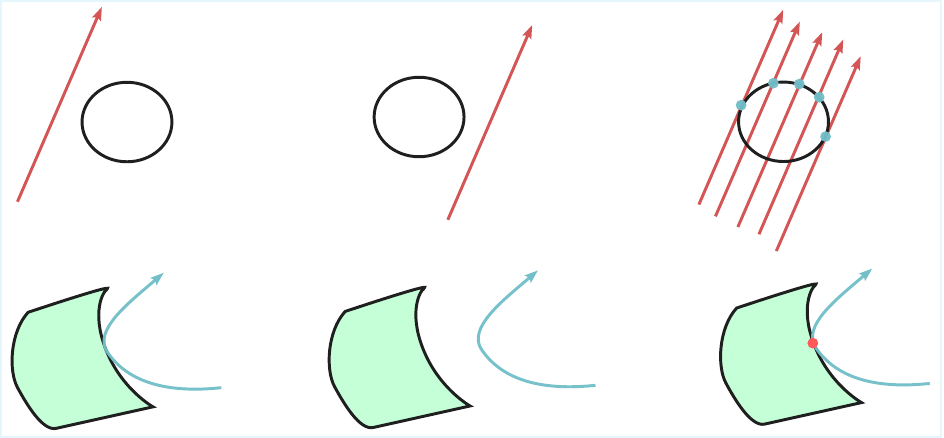}
\label{figure_coproducts_relations}
\end{figure}

Now, we claim that the relation $(i.)$ corresponds in $Cord^M(T_{K, \varepsilon}; \mathbb{Z}[\lambda^{\pm1}, \mu^{\pm1}])$ to non/bifurcations of a gradient-like trajectory of chords on $M_{K, \varepsilon}$. Since the first term of $(ii.)$ vanishes as a relative chain, the relation $(ii.)$ translates to $m_{\Delta_\varepsilon}=1$ in $Cord^M(T_{K, \varepsilon}; \mathbb{Z}[\lambda^{\pm1}, \mu^{\pm1}])$. Compare also with the skein relations for the knot case as in \cite[pg. 9]{Cieliebak2016KnotCH} or \cite[pg. 22]{Ekholm_2017}.
\end{rem}

\begin{assump} We assume that the reader is familiar with the Legendrian contact homology with loop space coefficients on the level of Appendix \ref{ch:file5}. We will briefly recall the concept.
\end{assump}

\begin{rem_not} $C(\mathcal{R})$ is a graded ring, which is as a $\mathbb{Z}$-module generated by the Reeb strings. Reeb strings are alternating words of Reeb chords on $\mathcal{L}^\ast_+ T_K$ and capping paths in $\mathcal{L}^\ast_+ T_K$ concatenating the ends of Reeb chords (and the base point $x_0$). The grading is given by a bidegree consisting of the singular chain degree and the sum of Reeb chord degrees. The differential on $C(\mathcal{R})$ is given by $\partial_{\mathcal{L}}:=\partial^{sing}+\partial^{sy}$, where $\partial^{sy}(\bm{a})$ is counting holomorphic discs $\mathcal{M}^{sy}/\R$ with the Reeb chord $\bm{a}$ at the positive infinity. Then, the homology of $(C(\mathcal{R}), \partial_{\mathcal{L}})$ is called \textbf{Legendrian contact homology}.
\end{rem_not}

\begin{thm}For $i=0, 1$ there are well defined maps $\Phi_i:C_i(\mathcal{R})\rightarrow C_i^{string}(\Sigma
)$ that intertwine with boundary operators, i.e.
\begin{equation}\label{egn_chain_map_holom}
\Phi_0\circ\partial_\mathcal{L}=(\partial^{sing}+\delta_Q+\delta_N)\circ\Phi_1.
\end{equation}
\end{thm}

\begin{sproof}We define maps $\Phi_i$ on Reeb chords by
\begin{equation*}\label{eqn_holom_chain}
\bm{a} \mapsto\sum_{\ell=0}^{k_{sw}}[\overline{\mathcal{M}}_{T_{K, \varepsilon}}(\bm{a}, \underbrace{\tfrac{1}{2},\dots,\tfrac{1}{2}}_{2\ell})],
\end{equation*}
where $[\overline{\mathcal{M}}_{T_{K, \varepsilon}}(\bm{a}, \tfrac{1}{2},\dots,\tfrac{1}{2})]$ denotes the boundaries of the disks from the moduli space, which become chains of broken strings after connecting the ends of $\bm{a}$ by capping paths to the base point $x_0$.

Now, by Conjecture \ref{cor_moduli_boundary_phenomena}, the dimension of $\overline{\mathcal{M}}_{T_{K, \varepsilon}}(\bm{a}, \tfrac{1}{2},\dots,\tfrac{1}{2})$ is $\vert\bm{a}\vert+\ell$. So, $\Phi_i(\bm{a})\in C_{\vert\bm{a}\vert}^{string}(\Sigma)$, and $\Phi_i$ preserves degree.

Let us verify the equation $(\ref{egn_chain_map_holom}).$ We need to inspect how the codimension $1$ strata of $\overline{\mathcal{M}}_{T_{K, \varepsilon}}(\bm{a}, \tfrac{1}{2},\dots,\tfrac{1}{2})$ are contributing to $\partial^{sing}\Phi(\bm{a})$ provided that $\vert\bm{a}\vert=1$. By Conjecture \ref{cor_moduli_boundary_phenomena} (Theorem \ref{thm_compact_strip}), there are four types of configurations:
\begin{itemize}
\item Discs with a $\tfrac{3}{2}$ switch. These discs contribute by $[0]$ as relative chains, since the normal jet conditions of the discs at $\tfrac{3}{2}$ switch are matching with $\partial^Q\Sigma$.
\item Discs with a symplectization building which are counted by $\Phi(\partial_{\mathcal{L}}(\bm{a}))$.
\item Discs with a switch of winding $1$ along $\R^3$ which are counted by $\delta_Q\Phi(\bm{a})$ (this phenomenon appeared already in \cite{Cieliebak2016KnotCH}).
\item Discs with a switch of winding $2$ along $L^\ast_+ T_{K, \varepsilon}$ which are counted by $\delta_N\Phi(\bm{a})$. Note that here we used the shift in grading of $C^{sing}_\ast(\Sigma^{Q_n}_\ell)$ by the number of trivial $Q$-strings.

For example, $\mathcal{M}_{T_{K, \varepsilon}}(\textbf{a}, 2, \tfrac{1}{2}, \tfrac{1}{2})\subset\partial\overline{\mathcal{M}}_{T_{K, \varepsilon}}(\textbf{a}, \tfrac{1}{2}, \tfrac{1}{2})$, so $\dim \mathcal{M}_{T_{K, \varepsilon}}(\textbf{a}, 2, \tfrac{1}{2}, \tfrac{1}{2})=1$. Hence, $[\mathcal{M}_{T_{K, \varepsilon}}(\textbf{a}, 2, \tfrac{1}{2}, \tfrac{1}{2})]\in C^{sing}_1(\Sigma^{Q_1}_2)$. But by the definition we have that $C^{sing}_1(\Sigma^{Q_1}_2)=C^{sing}_{0+2-1}(\Sigma^{Q_1}_2)\subset C^{string}_0(\Sigma)$.
\end{itemize}    
\vspace{-0.8cm}\end{sproof}

Now, let us discuss more the relations of $H_0^{string}(T_{K, \varepsilon})$ to $Cord^M(T_{K, \varepsilon})$ even though the reader is gently warned that the proceeding discussion will be slightly speculative.

\begin{defn}Let $T_0>0$. Then by $\mathcal{P}_{K, \varepsilon}$ we understand the space of paths $\mathcal{H}^3([0, t_0], \R^3; T_{K, \varepsilon})$ that are outward-pointing from $T_{K, \varepsilon}$ at their endpoints.

Next, by $\mathcal{P}_{K, \varepsilon}^{chord}\subset\mathcal{P}_{K, \varepsilon}$ we denote the realization of the space $M_{K, \varepsilon}$ as a space of $T_0$-time chords (with constant speed parametrizations). By $\mathcal{P}^\Delta_{K, \varepsilon}$ we denote the constant paths.
\end{defn}

\begin{rem}$\mathcal{P}_{K, \varepsilon}$ is a Banach manifold with corners. This is followed by similar argument as the proof of Lemma \ref{lemma_morse_hilb}.
\end{rem}

\begin{defn}By $\bm{E}_\varepsilon$ we denote the energy functional on $\mathcal{H}_3([0, T_0], \R^3; T_{K, \varepsilon})$. Next, $\partial^- \mathcal{P}_{K, \varepsilon}$ is defined as the subset of $\partial\mathcal{P}_{K, \varepsilon}$, where $-\nabla \bm{E}_\varepsilon$ is strictly outward-pointing via Definition \ref{def_outward_point}. Then by $\partial_- \mathcal{P}_{K, \varepsilon}$ we denote the closure of $\partial^- \mathcal{P}_{K, \varepsilon}$ in $\partial \mathcal{P}_{K, \varepsilon}$.
\end{defn}

\begin{lemma}\label{lemma_pseudo_def_retract}
The flow of $-\nabla \bm{E}_\varepsilon$ is a ``deformation'' of $\mathcal{P}_{K, \varepsilon}$ to $\mathcal{P}_{K, \varepsilon}^{chord} $ relative to $\partial_- \mathcal{P}_{K, \varepsilon}$.

  More precisely. With respect to the flow of $-\nabla \bm{E}_\varepsilon$ the set $\partial_- \mathcal{P}_{K, \varepsilon}\cup \mathcal{P}_{K, \varepsilon}^{chord}$ is a positively invariant exit set for $\mathcal{P}_{K, \varepsilon}$ (compare with Conley index pairs from Definition \ref{defn_index_pair}). In addition, for every $c\in \mathcal{P}_{K, \varepsilon}$ there is a $t_0\in[0, \infty]$ such that $\phi^{t_0}_{-\nabla\bm{E}_\varepsilon}\cap (\partial_- \mathcal{P}_{K, \varepsilon}\cup\mathcal{P}_{K, \varepsilon}^{chord})\neq\emptyset$.
\end{lemma}

\begin{proof}The lemma follows from the following observations. We know from Remark \ref{rem_crit_path} that 
$$Crit(\bm{E}_\varepsilon\vert_{\mathcal{P}_{K, \varepsilon}})=Crit(\bm{E}_\varepsilon\vert_{\mathcal{P}_{K, \varepsilon}^{chord}\cup \mathcal{P}_{K, \varepsilon}^{\Delta}})$$
and also $\mathcal{P}_{K, \varepsilon}^{\Delta}\subset \partial_- \mathcal{P}_{K, \varepsilon}$. 

Moreover, $\mathcal{P}_{K, \varepsilon}^{chord}$ is locally invariant under the flow $\phi_{-\nabla \bm{E}_\varepsilon}$ and $\phi_{-\nabla \bm{E}_\varepsilon}$ has no periodic orbits.
\end{proof}

\begin{conj}\label{conj_gradient_string}Let $K$ be generic and $\varepsilon>0$ small. Let $c$ be a path from $\partial\mathcal{P}_{K, \varepsilon}$. Then $c\in\partial^-\mathcal{P}_{K, \varepsilon}$ iff $c$ can be reparametrized such that one of the following conditions on normal jets $\sigma$ of $c$ holds:
\begin{itemize}
\item[$(S)$] At the starting point we have $\big(\sigma^{(1)}(0), \sigma^{(2)}(0), \sigma^{(3)}(0)\big)\in0\times0\times(-\infty, 0)$.
\item[$(E)$] At the endpoint we have $\big(\sigma^{(1)}(T_0), \sigma^{(2)}(T_0), \sigma^{(3)}(T_0)\big)\in0\times0\times(0, \infty)$.
\end{itemize}
Note that the curve $c$ satisfying $(S)$ or $(E)$ can be viewed as a $Q$-string in $\partial^Q\Sigma_1$.

Moreover, the analogous holds for paths in $\mathcal{P}_{K, \varepsilon}^{chord}$.
\end{conj}

\begin{hproof} Instead of general paths from $\mathcal{P}_{K, \varepsilon}$ we consider only $\mathcal{P}_{K, \varepsilon}^{chord}$ (i.e. $M_{K, \varepsilon}$). For this, we use the results of Section \ref{s:energy}, where we investigated the behavior of $-\nabla E_\varepsilon$ along $\partial M_{K, \varepsilon}$. We would like to show that the chords $\partial^{-}\mathcal{P}_{K, \varepsilon}^{chord}$ are precisely those boundary chords that appear as limits of interior chords while creating a single vanishing spike as in Figure \ref{figure_spike}.

Let us restrict $M_{K, \varepsilon}$ outside of weakly special and diagonal points. Recall that in Lemmata \ref{lem_mfld_standard} and \ref{lemma_standard_square} we described $M_{K, \varepsilon}\vert_{S_K}$ as a trivial bundle over the standard set $S_K\subset(\R/T\mathbb{Z})^2$  with a typical fiber as the square. Then $-\nabla E_\varepsilon$ is strictly outward-pointing along $\lbrace F_1^{[\varepsilon]}>0\wedge F_2^{[\varepsilon]}=0\rbrace$. Observe that along this set, the chords have an exterior tangent at the endpoints. Moreover, as the chords are approaching the chord with the tangency at the endpoint, they are additionally intersecting \textit{once} $T_{K, \varepsilon}$ close to the end-point. As the chords converge to the tangent one, the additional intersection points converge to the endpoint, see also Figure \ref{figure_contraction_standard} left. In particular, the chords of $\lbrace F_1^{[\varepsilon]}>0\wedge F_2^{[\varepsilon]}=0\rbrace$ will correspond to $(E).$ On the another hand, $-\nabla E_\varepsilon$ is strictly inward-pointing along $\lbrace F_1^{[\varepsilon]}=0\wedge F_2^{[\varepsilon]}>0\rbrace$. However, along this set, the chords have an exterior tangency at the starting points. Hence, these chords do not belong to the groups $(S)$ or $(E)$, see also Figure \ref{figure_contraction_standard} right.
\begin{figure}[!htbp]
\labellist
\endlabellist
\centering
\includegraphics[scale=0.55]{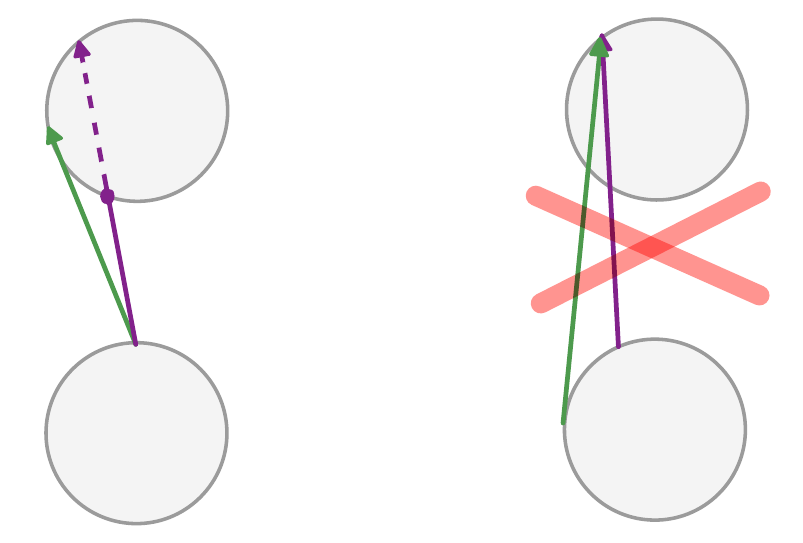}
\vspace{0.3cm}
\caption{A visualization of the chords near to the boundary of $M_{K, \varepsilon}\vert_{S_K}$. \textit{On the left:} the chords approaching the chord with the tangency at the endpoint become shorter, and the tangency at the endpoint is exterior. Note that the purple point will converge to the endpoint - compare with the vanishing spikes from Figure \ref{figure_spike}. \textit{On the right:} the chords approaching the chord with tangency at the starting point become longer, and the tangency at the starting point is exterior.}
\label{figure_contraction_standard}
\end{figure}

Now, let us inspect the behavior of chords from $M_{K, \varepsilon}$ near the diagonal $\Delta_\varepsilon$. From Section \ref{sec:closer_diag} we know that near $\Delta_\varepsilon$ the set $M_{K, \varepsilon}$ looks like a fibration over $\Delta_\varepsilon$ with cuspical fibers. Moreover we expect that $-\nabla E_\varepsilon$ is strictly outward-pointing along both sets $\lbrace F_1^{[\varepsilon]}>0\wedge F_2^{[\varepsilon]}=0\rbrace$ and $\lbrace F_1^{[\varepsilon]}=0\wedge F_2^{[\varepsilon]}>0\rbrace$. We would like to verify that the chords from these sets belong to the groups $(S)$ and $(E)$. For this, see Figure \ref{figure_contraction_diagonal}.
\begin{figure}[!htbp]
\labellist
\pinlabel $x$ at 0 165
\pinlabel $\textcolor{teal}{c_1}$ at 230 105
\pinlabel $\textcolor{teal}{c_2}$ at 190 190
\pinlabel $\textcolor{violet}{\widetilde{c}}$ at 250 165
\endlabellist
\centering
\includegraphics[scale=0.55]{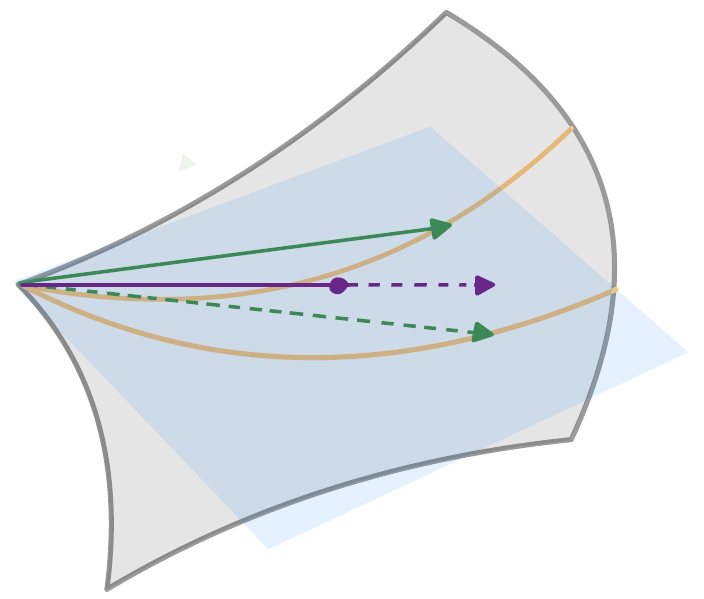}
\vspace{0.3cm}
\caption{A visualization of the outward-pointing chords that starts at some $x\in T_{K, \varepsilon}$ and are short, i.e. close to the diagonal. $T_{K, \varepsilon}$ is drawn as a graph, where the positive co-orientation is pointing below the graph, i.e., below the graph is the interior of $\nu_\varepsilon K$. The blue plane is the plane tangent to $T_{K, \varepsilon}$ at $x$. The orange lines describe the endpoints of the chords that emanate from $x$ and are tangent to $T_{K, \varepsilon}$ at the start/endpoint. The chord $\textcolor{teal}{c_1}$ has an interior tangent at the starting point and lies entirely on the tangent blue plane. $\textcolor{teal}{c_2}$ has an exterior tangent at the endpoint and on its interior lies outside to $\nu_\varepsilon K$. The purple chord $\textcolor{violet}{\widetilde{c}}$ is intersecting $T_{K, \varepsilon}$ also in its interior at the purple point. As $\widetilde{c}$ approaches $c_1$, the purple point converges to the starting point. Also, as $\widetilde{c}$ approaches $c_2$, the purple point converges to the endpoint point. Compare with vanishing spikes from Figure \ref{figure_spike}.}
\label{figure_contraction_diagonal}
\end{figure}
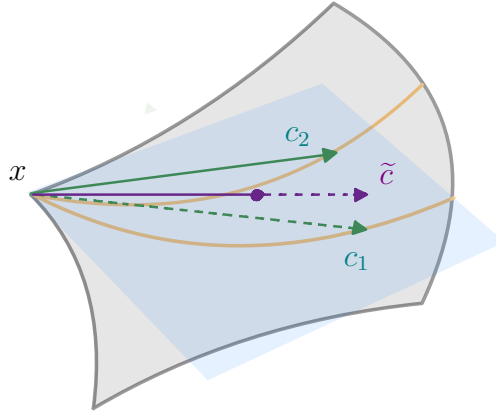
\end{hproof}

\newpage

\begin{rem}\textit{On a chain map from $C_k^{string}(\Sigma)$ to $C_k^M(T_{K, \varepsilon}; \mathbb{Z}[\lambda^{\pm 1}, \mu^{\pm 1}])$}

In the spirit of Lemma \ref{lemma_pseudo_def_retract} and Conjecture \ref{conj_gradient_string} it seems natural to take singular chains in $\Sigma$ and drag the $Q$-strings with the flow of $-\nabla \bm{E}_\varepsilon$ till $Crit(\bm{E}_\varepsilon)$ or $\partial^Q\Sigma$. In particular, $N$-strings will be translated to the words of $\lambda, \mu$. During this process, we will be performing string coproduct operations. However, it is not known how to bound the number of these string operations. 

Therefore, such an obstacle was resolved for the knot case by different deformation of broken strings, see \cite{Cieliebak2016KnotCH}. Their approach was to shrink $Q$-strings to piecewise-linear and then finally to linear. While keeping the (existing) switching points fixed, such a process will have control over the number of performed string operations.

If we try to mimic the approach from \cite{Cieliebak2016KnotCH} in the torus case, we encounter the following problem. Let us try to shrink a single piecewise linear $Q$-string to a linear one while fixing the endpoints. During this process, a $Q$-string $s$ can become tangent to $T_{K, \varepsilon}$ at one of its endpoints. If $s\in\partial^Q\Sigma$, then we can successfully forget about the string, since our theories are relative to this boundary. However, this is clearly not always the case. So the precise linearization process remains unknown.
\end{rem}

\appendix
\pagebreak
\part*{Appendices}
\addcontentsline{toc}{part}{Appendices}
\pagebreak
\chapter{Maslov index in $(\R^{2n}, \omega_0)$}
\label{ch:file1}

 
\section{Linear algebra}

\begin{defn} A \textbf{symplectic vector space} $(V, \omega)$ is a real finite dimensional vector space $V$ together with a \textbf{symplectic form} $\omega$, that is, a skew--symmetric nondegenerate bilinear map $\omega:V\times V\rightarrow \R$.

A map $\varphi:(V_1, \omega_1)\rightarrow(V_2, \omega_2)$ is called a \textbf{linear symplectomorphism} if it is a vector space isomorphism and $\varphi^\ast\omega_2=\omega_1$.

A basis $\lbrace v_i, w_i\rbrace_i$ of $(V, \omega)$ is called a \textbf{symplectic basis} if 
$$\omega(v_j, v_k)=\omega(w_j, w_k)=0\text{ and }\omega(v_j, w_k)=\delta_{jk}.$$
\end{defn}

\begin{rem} \label{rem_basis} Every $(V, \omega)$ has a symplectic basis, in particular $(V, \omega)$ is even-dimensional. Also, the choice of the symplectic basis uniquely determines $\omega$. For details, see \cite[Thm 2.1.3.]{mcduff_salamon_2017}
\end{rem}

\begin{example} $(\R^{2n}, \omega_0)$. Let $(x_i, y_i)$ be standard (Cartesian) coordinates on $\R^{2n}$, then its basis vectors $\lbrace v^{x_i}, w^{y_i}\rbrace$ form a symplectic basis of the \textbf{standard symplectic form} $\omega_0$. Hence, in view of Remark \ref{rem_basis}, every symplectic vector space is symplectomorphic to $(\R^{2n}, \omega_0)$ for some $n$. 

More explicitly, $\omega_0$ is given as
$$\omega_0(v, w):=\langle v, -J_0 w\rangle,$$
$ \text{where } v, w\in R^{2n}, J_0:=\begin{pmatrix}
0 & -1\\
1 & 0
\end{pmatrix}$, and $\langle\, ,\,\rangle$ is the standard inner product. Here $J_0$ is called the \textbf{standard complex structure}, because we can identify $(\R^{2n}, \omega)$ with $(\mathbb{C}^n, i)$ such that the action of $J_0$ is identified with the multiplication by $i$.

The \textbf{symplectic group} $Sp(2n, \R)$ is the group of all linear symplectomorphisms $\varphi:(\R^{2n}, \omega_0)\rightarrow(\R^{2n}, \omega_0)$. $Sp(2n, \R)$ is connected, see \cite[Prop 2.2.4]{mcduff_salamon_2017}. 
\end{example}

\begin{example} \label{example_dual} $(V\times V^\ast, \omega_\times)$, where 
$$\omega_\times ((v_1, v_2^\ast), (w_1, w_2^\ast)):=v_2^\ast(w_1)-w_2^\ast(v_1),$$
for $(v, v^\ast), (w, w^\ast)\in V\times V^\ast.$

Then the \textbf{anti-symplectic involution} on $(V\times V^\ast, \omega_\times)$ is the map $\sigma_{st}:(v_1, v_2^\ast)\mapsto(v_1, -v_2^\ast)$. Note that $\sigma_{st}^{\ast}\omega_\times=-\omega_\times$.
\end{example}

\begin{defn} Let $V$ be a real even dimensional vector space. Then a \textbf{linear complex structure} $J$ is a (real) linear automorphism on $V$ such that 
$$J^{2}=-\mathbbm{1}.$$
We will write $(V, J)$ and call this tuple a \textbf{complex vector space}.
\end{defn}

\begin{rem} For any vector space $(V, J)$ of the real dimension $2n$ there are $n$ vectors $\lbrace v_i\rbrace_i$ such that $\lbrace v_i, Jv_i\rbrace_i$ is a real basis of $V$. In particular, there is a vector space isomorphism $\varphi:V\rightarrow\R^{2n}$ such that $\varphi^\ast J_0=J$. For details, see \cite[Thm 2.5.1.]{mcduff_salamon_2017} 
\end{rem}

\begin{defn}
Let $(V, J)$ be a complex vector space and $\omega$ a symplectic form on $V$. Then $J$ is called \textbf{compatible with} $\omega$ if
$$\omega(Jv, Jw)=\omega(v, w)\text{ for all }v, w\in V$$
and
$$\omega(v, Jv)>0\text{ for all nonzero }v\in V.$$
\textbf{The space of all $J$ compatible with} $\omega$ will be denoted by
$$\mathcal{J}(V, \omega).$$

Also, we define the space 
$$\Omega(V, J):=\lbrace \omega \text{ symplectic form on }V\,|\,J\in\mathcal{J}(V, \omega)\rbrace,$$
or equivalently,
$$\Omega(V, J):=\lbrace g(J\cdot,\cdot)\,|\,g\text{ inner product on }V\text{ such that }J^\ast g=g\rbrace.$$
\end{defn}

\begin{rem} If $J\in\mathcal{J}(V, \omega)$, then $g_J(v, w):=\omega(v, Jw)$ defines an inner product on $V$ and $h:=g_J+i\omega$ defines a Hermitian inner product.
\end{rem}

\begin{lemma}[{\cite[Thm 2.5.5.]{mcduff_salamon_2017}}] \label{Lemma_rotation} $\mathcal{J}(V, \omega)$ is non-empty and contractible.
\end{lemma}

\begin{proof} Due to the existence of the symplectic basis, we only need to consider $\mathcal{J}(R^{2n}, \omega_0)$.

Elements of $\mathcal{J}(\R^{2n}, \omega_0)$ can be characterized by 
\begin{equation}\label{complx} J^2=-\mathbbm{1},\quad J^TJ_0J=J_0, \quad \langle v, -J_0Jv\rangle>0\text{ for all nonzero }v\in \R^{2n}.
\end{equation}
Thus
$$(J_0J)^T=-J^TJ_0=J^TJ_0J^2=J_0J,$$
where in the first equality we used the definition of $J_0$. In the next two equalities, we used consequently the first two identities of (\ref{complx}).
Observe that 
$$P:=-J_0J$$
is a symmetric, positive definite, and symplectic matrix. Hence any $J\in\mathcal{J}(\mathbb{R}^{2n},\omega_{0})$ can be written uniquely in the form $J=J_0P$, where $P$ satisfies the conditions above. Then by \cite[Thm 2.2.3.]{mcduff_salamon_2017} from the eigendecompositions of $P$ and $P^\alpha$ we know that any power $P^\alpha$ is also of this kind for $\alpha\geq 0$. Hence, we can connect $J=J_0P^1$ with $J_0$ by sending the exponent $1\rightarrow 0$. Note that the obtained path lies in $\mathcal{J}(\mathbb{R}^{2n},\omega_{0})$, which concludes the proof.
\end{proof}

\begin{rem} Let $\tilde{g}$ be an inner product on $(V, J)$. Then for $g:=1/2(\tilde{g}+J^\ast \tilde{g})$ we have $J^\ast g=g$. Hence $\Omega(V, J)$ is non-empty. Moreover, $\Omega(V, J)$ is convex and thus contractible.
\end{rem}

\begin{defn} A linear subspace $\Lambda$ of $(V, \omega)$ is called \textbf{Lagrangian} if $\omega|_\Lambda=0$ and $\dim \Lambda=\frac{1}{2}\dim V$.
\end{defn}

\begin{rem} Lagrangian subspaces are the largest linear subspaces on which the symplectic form vanishes.
\end{rem}

\begin{rem}\label{rem_two_lagr} Let $\Lambda_1$ and $\Lambda_2$ be transverse Lagrangians in $(V, \omega)$. If $\lbrace v_1,\dots, v_n\rbrace$ is a basis of $\Lambda_1$, then there is an unique basis $\lbrace w_1,\dots, w_n\rbrace$ of $\Lambda_2$ such that $\lbrace v_1,\dots, v_n, w_1,\dots, w_n\rbrace$ is a symplectic basis of $(V, \omega)$. See \cite{otto_2019}.
\end{rem}

\begin{defn} A linear subspace $R$ of $(V, J)$ is called \textbf{totally real} if $R\cap JR=\lbrace 0\rbrace$ and $\dim R=\frac{1}{2}\dim V$.
\end{defn}

\begin{lemma}For Lagrangian and totally real subspaces it holds:
\begin{itemize}
\item[$(i.)$] Let $\Lambda$ be a Lagrangian subspace in $(V, \omega)$. Then $\Lambda\perp_{g_J}J\Lambda$ for every $J\in\mathcal{J}(V, \omega)$. In particular, $\Lambda$ is totally real.
\item[$(ii.)$] Let $R$ be a totally real subspace in $(V, J)$. Then there is $\omega\in\Omega(V, J)$ such that $R$ is $\omega$-Lagrangian. Moreover, the space of such $\omega$ is contractible.
\end{itemize}
\end{lemma}

\begin{proof}
Ad $(i.)$: Let $v_1, v_2\in\Lambda$. Then $g_J(Jv_1, v_2)=\omega(Jv_1, Jv_2)=\omega(v_1, v_2)=0.$

Ad $(ii.)$: Let $\tilde{g}$ be any inner product on $R$. It induces an inner product on $JR$ by $J^\ast\tilde{g}$. This naturally induces the unique inner product $g$ on $V$ such that $R\perp_g JR$. Hence $J^\ast g=g$ and $g(J\cdot, \cdot)$ is our desired symplectic form. Since the space of inner products on $R$ is convex, the lemma follows.
\end{proof}

\begin{rem} For a pair of totally real subspaces in $(V, J)$ there might not exist
$\omega\in\Omega(V, J)$ such that they become both Lagrangian.

Indeed, let $R_1=\langle (1, 0)^T, (0, 1)^T\rangle$ and $R_2=\langle (i, 0)^T, (1, i)^T\rangle$ in $(\mathbb{C}^2, i)$. Assume that $R_1$ is Lagrangian for some $\omega\in\Omega(\mathbb{C}^2, i)$. Then
\begin{align*}
\omega((i, 0)^T, (1, i)^T)&=\omega((i, 0)^T, (1, 0)^T)+\omega((i, 0)^T, (0, i)^T)\\
&=\omega(i(1, 0)^T, (1, 0)^T)+\omega(i(1, 0)^T, i(0, 1)^T)\\
&=-\omega((1, 0)^T, i(1, 0)^T)<0.
\end{align*} 
\end{rem}

\section{Maslov index}
In this section, we start with the Maslov index for loops in the Lagrangian Grassmannian. Then we generalize this notion to loops in the totally real Grassmannian and paths in the Lagrangian Grassmannian. The section will be finished by introducing the notion of K{\"a}hler angles.

\begin{defn} The \textbf{Lagrangian Grassmanian} $\Lambda(n)$ is the manifold consisting of all Lagrangian subspaces of $(\mathbb{R}^{2n},\omega_0)$.
\end{defn}

\begin{lemma}[\cite{arnold_67}]\label{Lemma_arnold_maslov} $\pi_1(\Lambda(n))\cong\mathbb{Z}$.
\end{lemma}

\begin{proof} First, fibre bundles can be seen as examples of fibrations, see \cite[Prop 4.48]{hatcher_2019}.
Hence, the determinant map $\det$ on $O(n)$ and $U(n)$ induce the following fibration sequences
$$SO(n)\longrightarrow O(n)\xrightarrow{\det} S^0\text{ and }SU(n)\longrightarrow U(n)\xrightarrow{\det} S^1.$$

Also, $U(n)$ acts on $\Lambda(n)$ transitively with stabilizer $O(n)$ at $\R^n$, and in particular $\Lambda$ really has the structure of a manifold. Indeed, for an $\Lambda\in\Lambda(n)$, observe that $\Lambda$ and $J_0\Lambda$ are orthogonal with respect to the standard inner product. Next, if $B, B^\prime$ are orthonormal bases of $\Lambda, \Lambda^\prime\in \Lambda(n)$ respectively, then a linear map that maps $B$ to $B^\prime$ and $J_0B$ to $J_0B^\prime$ is a unitary automorphism. Here we used the fact that $GL(n, \mathbb{C})\cap O(2n)=U(n)$. From the ambiguity of the choice of orthonormal basis, it follows that the stabilizer of this action at $\R^n$ is $O(n)$. Hence, we obtain a fibration
$$O(n)\longrightarrow U(n)\longrightarrow \Lambda(n).$$

Recall that $\det$ of any orthogonal automorphism is equal to $\pm 1$. Hence, we take $\det^2$ to have well defined map on $\Lambda(n)$ that assigns to any $\Lambda\in\Lambda(n)$ the determinant square of the unitary automorphism mapping the $\lbrace x_i\rbrace_i$-plane $\subset\R^{2n}$ to $\Lambda$. This induces a fibration 
$$S\Lambda(n)\longrightarrow \Lambda(n)\xrightarrow{\det^2} S^1,$$
where $S\Lambda(n)$ is the set of $\Lambda\in \Lambda(n)$ such that $\det^2(\Lambda)=1$. Since $SU(n)$ acts transitively on $S\Lambda(n)$ with stabilizer $SO(n)$, $S\Lambda(n)\cong SU(n)/SO(n)$. Hence, we obtain the following commutative diagram of six fibrations
\[
\begin{tikzcd}
SO(n) \arrow[r] \arrow[d] & O(n) \arrow[r, "\det"] \arrow[d]           & S^0 \arrow[d]        \\
SU(n) \arrow[d] \arrow[r] & U(n) \arrow[d] \arrow[r, "\det"] & S^1 \arrow[d, "z^2"] \\
S\Lambda(n) \arrow[r]     & \Lambda(n) \arrow[r, "\det^2"]             & S^1\nospaceperiod              
\end{tikzcd}
\]

From this, we conclude the following LES of homotopy groups, see \cite[Thm 4.41]{hatcher_2019}:
\begin{equation}\label{LES}\dots\longrightarrow\pi_1\left(\frac{SU(n)}{SO(n)}\right)\longrightarrow\pi_1(\Lambda(n))\xrightarrow{(\det^2)_\ast}\pi_1(S^1)\longrightarrow\pi_0\left(\frac{SU(n)}{SO(n)}\right)\longrightarrow\dots
\end{equation}
It remains to study ${SU(n)}/{SO(n)}$. Since $SO(n)$ is path connected and $SU(n)$ is simply connected, the LES of homotopy groups shows us that ${SU(n)}/{SO(n)}$ is simply connected and, in particular, path connected. Plugging this to the sequence (\ref{LES}) we obtain that 
$$\pi_1(\Lambda(n))\xrightarrow[isom]{(\det^2)_\ast} \pi_1(S^1)\xrightarrow[isom]{\deg}\mathbb{Z}.$$
\end{proof}

\begin{defn} \label{defn_masl_1} Let $\gamma$ be a continuous loop in $\Lambda(n)$. Then the \textbf{Maslov index of} $\gamma$ is defined as
$$\mu(\gamma):=\degr\circ\,{\det}^2(\gamma).$$
\end{defn}

\begin{defn_thm}[{\cite{otto_2019}}]\label{defn_maslov} Let $\Lambda_1\in\Lambda(n)$. Then $\Lambda(n)$ can be written as the following disjoint union
$$\Lambda(n)=\bigcup_{k\geq 0}\Sigma_{k}(\Lambda_1),$$
where each $\Sigma_{k}(\Lambda_1)$ is the subset of $\Lambda(n)$ of all Lagrangian subspaces that have $k$-dimensional intersection with $\Lambda_1$. It follows that $\Sigma_{k}(\Lambda_1)$ is a submanifold of $\Lambda(n)$ such that 
$\dim\Lambda(n)=n(n+1)/2$ and $\codim\Sigma_{k}(\Lambda_1)=k(k+1)/2$.
Then \textbf{Maslov cycle (with respect to $\Lambda_1$)} is the following stratified space
$$\Sigma(\Lambda_1):=\overline{\Sigma_1(\Lambda)}=\bigcup_{k\geq 1}\Sigma_{k}(\Lambda_1).$$
$\Sigma(\Lambda_1)$ has the structure of a singular algebraic variety of codimension $1$. Here, lower strata are at least of codimension $2$ in $\Sigma(\Lambda_1)$ (specially codimension $1$ boundary of $\Sigma(\Lambda_1)$ is empty).
\end{defn_thm}

\begin{rem}\label{rem_cross_form} Let $\Lambda_2$ be an element of $\Lambda(n)$. Pick any $\Lambda_1\in\Sigma_0(\Lambda_2)$. Since $\Lambda_2$ and $\Lambda_1$ are transverse, they can be identified respectively with the $\lbrace x_i\rbrace_i\text{--plane}$ and the $\lbrace y_i\rbrace_i\text{--plane}$ in $(\R^{2n}, \omega_0)$. Hence any $\Lambda\in\Sigma_0(\Lambda_1)$ can be uniquely expressed as a graph $\lbrace(a, Sa)\,|\,a\in\lbrace x_i\rbrace_i\text{--plane}\subset\R^{2n}\rbrace$ for some linear map $S:\R^n\rightarrow\R^n$. Since $\Lambda$ is Lagrangian, we deduce that $S$ is represented by a symmetric matrix. 
In fact, this gives us a diffeomorphism
\begin{align*}
\Sigma_0(\Lambda_1)&\longrightarrow S^2(\Lambda_2)\\
\Lambda&\longmapsto\langle \cdot, S\cdot\rangle,
\end{align*}
where $S^2(\Lambda_2)$ is the vector space of all quadratic forms on $\Lambda_2$.

Next, observe that any smooth curve 
$\gamma(t):(-\epsilon, \epsilon)\rightarrow\Lambda(n)$, such that $\gamma(0)=\Lambda_2,$
 can be expressed as $\lbrace(a, S(t)a)\,|\,a\in\R^n\rbrace$. Hence the tangent vector $\overset{\bm .}{\gamma}(0)=\widehat{\Lambda}_2$ can be represented by the quadratic form $\langle \cdot, \overset{\bm .}{S}(0)\cdot\rangle$. This gives us a vector space isomorphism
\begin{align}\label{eqv_tangent}
\begin{split}
T_{\Lambda_2}\Lambda(n)&\longrightarrow S^2(\Lambda_2)\\
\widehat{\Lambda}_2&\longmapsto\langle \cdot, \overset{\bm .}{S}(0)\cdot\rangle=:Q(\Lambda_2, \widehat{\Lambda}_2).
\end{split}
\end{align}
This isomorphism is well defined, i.e. $Q(\Lambda_2, \widehat{\Lambda}_2)$ does not depend on the choice of $\Lambda_1$, see \cite[Thm 1.1]{robbin_salamon_1993} .
\end{rem}


\begin{defn} Let $\gamma: [0, 1]\rightarrow \Lambda(n)$ be a smooth path. Then the number $t\in[0, 1]$ is called a \textbf{crossing} if $\gamma(t)\in\Sigma_k(\Lambda_1)$ for some $k\in\lbrace 1,\dots, n\rbrace$. At each crossing $t$ the \textbf{crossing form} is defined as
$$\Gamma(\gamma, \Lambda_1, t):=\left.Q(\gamma(t), \dot{\gamma}(t))\right|_{\gamma(t)\cap\Lambda_1}.$$
\end{defn}

\begin{rem} For any $k\in\lbrace 0,\dots, n\rbrace$, we could express the tangent space to $\Sigma_k(\Lambda_1)$ at any point $\Lambda_2\in\Sigma_k(\Lambda_1)$ as
\begin{equation}\label{eqv_subtangent}
T_{\Lambda_2}\Sigma_k(\Lambda_1)=\big\lbrace \widehat{\Lambda}_2\in T_{\Lambda_2}\Lambda(n)\,|\left.Q(\Lambda_2, \widehat{\Lambda}_2)\right|_{\Lambda_2\cap\Lambda_1}=0\big\rbrace.
\end{equation}
See \cite[Lemma 10.2.4]{otto_2019}.

Moreover, any smooth path $\gamma$ in $\Lambda(n)$ is tangent to $\Sigma_k(\Lambda_1)$ at the crossing $t$ if and only if $\gamma(t)\in\Sigma_k(\Lambda_1)$ and $\Gamma(\gamma, \Lambda_1, t)=0$. This motivates the following definition, which can be seen as a certain generalization of transversality.
\end{rem}

\begin{defn} A crossing is called \textbf{regular} if the corresponding crossing form is non-singular. 
\end{defn}

\begin{rem} Observe that Formulae (\ref{eqv_tangent}) and (\ref{eqv_subtangent}) allow us to coorient $\Sigma_1(\Lambda_1)$ using the quadratic forms. Hence, if $t$ is a regular crossing for $\gamma(t)\in \Sigma_1(\Lambda_1)$, then the signature of the crossing form is equal to the algebraic intersection number of $\gamma$ and $\Sigma_1(\Lambda_1)$ at $\gamma(t)$. 
\end{rem}

\begin{defn} Let $\gamma$ be a smooth path $\gamma:[0, 1]\rightarrow\Lambda(n)$ with only regular crossings. Then the \textbf{relative Maslov index of $\gamma$ with respect to $\Lambda_1$} is defined as
\begin{equation}\label{rel_maslov}
\mu(\gamma, \Lambda_1):=\frac{1}{2}\sign \Gamma(\gamma, \Lambda_1, 0)+\sum_{0<t<1}\sign \Gamma(\gamma, \Lambda_1, t)+\frac{1}{2}\sign \Gamma(\gamma, \Lambda_1, 1).
\end{equation}
\end{defn}

\begin{rem} Regular crossings are isolated. Consequently, there is only a finite number of nonzero summands in Formula (\ref{rel_maslov}), see \cite[pg.~9]{robbin_salamon_1993}.
\end{rem}

\begin{rem} There is a natural extension of the relative Maslov index to all continuous paths in $\Lambda(n)$.

Indeed, every continuous path in $\Lambda(n)$ is homotopic with fixed endpoints to some smooth path with only regular crossings, see \cite[Lemma 2.1]{robbin_salamon_1993}. If any two smooth paths in $\Lambda(n)$ have only regular crossings and are smoothly homotopic with fixed endpoints, then they have the same Maslov index, see \cite[Lemma 2.2]{robbin_salamon_1993}. The rest follows from the Whitney approximation theorem, see \cite{lee_2013}.
\end{rem}

\begin{lemma} [{\cite[Thm 2.3~and~Rem 2.6]{robbin_salamon_1993}}] The relative Maslov index satisfies the following properties:
\begin{itemize}
\item[(Naturality)] For $\Psi\in Sp(2n, \R)$ we have
$$\mu(\Psi \gamma, \Psi \Lambda_1)=\mu(\gamma, \Lambda_1).$$
\item[(Concatenation)] For $a<b<c$ we have
$$\mu(\gamma|_{[a, c]}, \Lambda_1)=\mu(\gamma|_{[a, b]}, \Lambda_1)+\mu(\gamma|_{[b, c]}, \Lambda_1).$$
\item[(Product)] If $n'+n''=n$, then $\Lambda(n')\times\Lambda(n'')$ can be seen as a submanifold of $\Lambda(n)$ and 
$$\mu(\gamma'\times \gamma'', \Lambda_1'\times\Lambda_1'')=\mu(\gamma', \Lambda_1')+\mu(\gamma'', \Lambda_1'').$$
\item[(Homotopy)] Two paths in $\Lambda(n)$ with the same endpoints are homotopic with fixed endpoints if and only if they have the same Maslov index.
\item[(Zero)] Let $k\in\lbrace 1,\dots, n\rbrace$. Then for any path $\gamma$ in $\Sigma_k(\Lambda_1)$ we have 
$$\mu(\gamma, \Lambda_1)=0.$$
\item[(Loop)] If $\gamma$ is a loop in $\Lambda(n)$, then
$$\mu(\gamma, \Lambda_1)=\mu(\gamma).$$
In particular, the relative Maslov index of a loop does not depend on the choice of the Maslov cycle.
\end{itemize}
\end{lemma}

\begin{rem}\label{rem_opposite_path} If $\gamma(t)$, where $t\in[0, 1]$, is a path in $\Lambda(n)$, then $-\gamma$ denotes the \textbf{opposite path} $\gamma(1-t)$. Thus $\gamma\ast-\gamma$ is a contractible loop and consequently $\mu(\gamma, \Lambda_1)=-\mu(-\gamma, \Lambda_1)$ for any $\Lambda_1\in\Lambda(n)$.
\end{rem}

\noindent{\it For the simplicity, we will use the following (natural) notation. If $(V, J)$ is a complex vector space and $V=V^1\oplus\dots\oplus V^k$, then the expression $e^J V^i$ means $e^{J|_{V^i}}V^i.$}

\begin{defn} [\cite{EES05}] Let $\Lambda_1, \Lambda_2\in\Lambda(n)$ and $J\in\mathcal{J}(\R^{2n}, \omega_0)$. Then the \textbf{K{\"a}hler angle} $\theta(\Lambda_1, \Lambda_2)=(\theta_1, \dots, \theta_n)\in[0, \pi)^n$ and orthogonal decomposition $R^{2n}=W^0\oplus\dots\oplus W^s$ are defined inductively as follows. Take the Hermitian inner product $h=g_J+i\omega_0$.

 If $\dim(\Lambda_1\cap\Lambda_2)=r_0$, then put $$\theta_1=\dots=\theta_{r_0}=0\hbox{ and }W^0=(\Lambda_1\cap\Lambda_2)\oplus J(\Lambda_1\cap\Lambda_2).$$

Then we define $\Lambda^1_i:=(\Lambda_1\cap\Lambda_2)^{\bot_h}\cap\Lambda_i\in\Lambda(n-r^0)$ for $i=1, 2$. Next, if $\alpha_1\in(0, \pi)$ is the smallest angle such that $\dim(e^{\alpha_1 J}\Lambda^1_1\cap\Lambda^1_2)=r_1>0$, then we set $$\theta_{r_0+1}=\dots=\theta_{r_0+r_1}=\alpha_1\hbox{ and }W^1=(e^{\alpha_1 J}\Lambda_1^1\cap\Lambda_2^1)\oplus J(e^{\alpha_1 J}\Lambda_1^1\cap\Lambda_2^1).$$

Then we define $\Lambda^2_1:=(e^{\alpha_1 J}\Lambda_1^1\cap\Lambda_2^1)^{\bot_h}\cap e^{\alpha_1 J}\Lambda_1^1$ and $\Lambda^2_2:=(e^{\alpha_1 J}\Lambda_1^1\cap\Lambda_2^1)^{\bot_h}\cap \Lambda_2^1$. Here $\Lambda_1^2, \Lambda_2^2\in\Lambda(n-r_0-r_1)$. Next, if $\alpha_2\in(0, \pi)$ is the smallest angle such that $\dim(e^{\alpha_2 J}\Lambda^2_1\cap\Lambda^2_2)=r_2>0$, then we set $$\theta_{r_0+r_1+1}=\dots=\theta_{r_0+r_1+r_2}=\alpha_1+\alpha_2\hbox{ and }W^2=(e^{\alpha_2 J}\Lambda_1^2\cap\Lambda_2^2)\oplus J(e^{\alpha_2 J}\Lambda_1^2\cap\Lambda_2^2).$$
Now, we repeat until $\theta_n$ is defined.

Next, using the Hermitian splitting of $\R^{2n}$ we define $\theta(\Lambda_1,, \Lambda_2)J$ as $\theta_{r_0+\dots +r_k}J|_{W^k}$ on $W^k$.

Then by the \textbf{positive rotation of $\Lambda_1$ to $\Lambda_2$ by $J$} we understand the path 
\begin{align*}
rot\>\![\Lambda_1, \Lambda_2; J,+]:[0, 1]&\rightarrow\Lambda(n)\\
t&\mapsto e^{t\theta(\Lambda_1, \Lambda_2) J}\Lambda_1.
\end{align*}

Also by the \textbf{negative rotation of $\Lambda_1$ to $\Lambda_2$ by $J$} we understand the path 
\begin{align*}
rot\>\![\Lambda_1, \Lambda_2; J,-]:[0, 1]&\rightarrow\Lambda(n)\\
t&\mapsto e^{t(\theta(\Lambda_1, \Lambda_2)-\pi) J}\Lambda_1,
\end{align*}
where $\theta(\Lambda_1, \Lambda_2)-\pi:=(\theta_1-\pi,\dots, \theta_n-\pi)$.

Finally, when $J$ is clear from the context, we can emphasize the angle by saying only \textbf{rotation of $\Lambda_1$ to $\Lambda_2$ by $\theta(\Lambda_1, \Lambda_2)$ or by $\theta(\Lambda_1, \Lambda_2)-\pi$}.
\end{defn}

\begin{lemma}\label{lemma_kahler_def} Let $\Lambda_1, \Lambda_2\in\Lambda(n)$ and $J\in\mathcal{J}(\R^{2n}, \omega_0)$. If $\Lambda_1\pitchfork\Lambda_2$, then there is $\alpha\in(0, \pi)$ such that $e^{\alpha J}\Lambda_1\cap\Lambda_2\neq\lbrace 0\rbrace.$ In particular, K{\"a}hler angle is well-defined.
\end{lemma}

\begin{proof}
By the choice of the symplectic basis, we can assume that $\Lambda_1$ is the $\lbrace x_i\rbrace_i$-plane and $J=J_0$. If $J_0\Lambda_1\cap\Lambda_2\neq\lbrace 0\rbrace,$ we can put $\alpha=\pi/2$.

So let us consider the case when $J_0\Lambda_1\cap\Lambda_2=\lbrace 0\rbrace.$ Hence, as in the Remark \ref{rem_cross_form}, $\Lambda_2$ is a graph $\lbrace (a,Sa)\,|\,a\in\lbrace x_i\rbrace_i\hbox{-plane and }S\hbox{ is a symmetrix matrix}\rbrace$.
Moreover 
\begin{align*}
e^{\alpha J_0}&=\sum_{k=0}^\infty\frac{1}{k!}(\alpha J_0)^k\\
&=\sum_{k=0}^\infty\frac{1}{(2k)!}a^{2k}(-\mathbbm{1})^k+\sum_{k=0}^\infty\frac{1}{(2k+1)!}a^{2k+1}(-\mathbbm{1})^k J_0\\
&=\cos (\alpha)+ \sin(\alpha) J_0.
\end{align*}
Hence, $e^{\alpha J_0}\Lambda_1\cap\Lambda_2$ is non-trivial if and only if the system of equations
\begin{align*}
\begin{cases}
 \,\,\,\, a=\cos (\alpha) b,\\
 Sa=\sin(\alpha) b.
 \end{cases}
\end{align*}
has a solution for $a, b\in\R^n$ and $\alpha\in(0, \pi)\setminus\lbrace\pi/2\rbrace$.
This is equivalent to 
\begin{equation}\label{eqn_rotation_positive}
Sa=\tan(\alpha) a.
\end{equation}
But $S$ is symmetric, hence has a basis of eigenvectors with real eigenvalues. Thus $\alpha$ exists and we are done.
\end{proof}

\begin{lemma}\label{rem_positive_negative} Let $\Lambda_1, \Lambda_2\in\Lambda(n)$ and $J\in\mathcal{J}(\R^{2n}, \omega_0)$, then $\rot\>\![\Lambda_1, \Lambda_2; J,+]=-\rot\>\![\Lambda_2, \Lambda_1; J,-]$.
\end{lemma}

\begin{proof}
Let us show the lemma in the case when $\Lambda_1\pitchfork\Lambda_2$ and $J_0\Lambda_1\cap\Lambda_2=\lbrace 0\rbrace$. Then the general case easily follows. Moreover, similarly as in Lemma \ref{lemma_kahler_def} we can assume that
$\Lambda_1$ is the $\lbrace x_i\rbrace_i$-plane and $J=J_0$. And hence analogously $\Lambda_2=\lbrace(a, Sa)\rbrace$. 

Now, $e^{\beta J_0}\Lambda_1\cap\Lambda_2$, for some $\beta\in(0, \pi)\setminus\lbrace\pi/2\rbrace$, is non-trivial if and only if the system of equations
\begin{align*}
\begin{cases}
 \cos(\beta)a-\sin(\beta)Sa=b,\\
 \cos(\beta)Sa+\sin(\beta)a=0.
 \end{cases}
\end{align*}
has a solution for $a, b\in\R^n$ and $\beta\in(0, \pi)\setminus\lbrace\pi/2\rbrace$.

This is equivalent to 
$$Sa=-\tan(\beta)a$$ 
and hence to
\begin{equation}\label{eqn_rotation_opposite}
Sa=\tan(\pi-\beta)a.
\end{equation}
Now, we compare Formulas (\ref{eqn_rotation_positive}) and (\ref{eqn_rotation_opposite}). If we take a $g_{J_0}$-orthogonal eigenbasis of $S$, then eigenvectors generate (as $J_0$-complex vector spaces) same orthogonal decompositions of $\R^{2n}$ for $\theta(\Lambda_1, \Lambda_2)$ and $\theta(\Lambda_2, \Lambda_1)$. And the lemma follows.
\end{proof}

\begin{defn} Let $\Lambda_1, \Lambda_2\in\Lambda(n)$. Then $\Sigma(\Lambda_1, \Lambda_2)$ denote the set of all $\Lambda\in\Lambda(n)$ such that $\Lambda\cap\Lambda_1=\Lambda_2\cap\Lambda_1$.
\end{defn}

\begin{lemma} \label{lemma_contr} $\Sigma(\Lambda_1, \Lambda_2)$ is a contractible manifold.
\end{lemma}

\begin{proof}
Observe that when $\Lambda_1$ and $\Lambda_2$ are transverse, then by the Remark \ref{rem_cross_form} $\Sigma(\Lambda_1, \Lambda_2)$ is identified with the space of symmetric $n\times n$-matrices and hence is contractible.

Suppose that $\Lambda_1$ and $\Lambda_2$ are not transverse. We take a compatible $J$ and the induced Hermitian inner product $h$. Since $\Lambda\cap\Lambda_1$ is fixed for any $\Sigma(\Lambda_1, \Lambda_2)$, we will restrict ourselves only to the Hermitian complement $(\Lambda_2\cap\Lambda_1)^{\bot_h}$. Hence, analogously as before, $\Sigma(\Lambda_1, \Lambda_2)$ is identified with the space of symmetric $(n-k)\times (n-k)$-matrices, where $k=\dim(\Lambda_2\cap\Lambda_1)$, and we are done.
\end{proof}

\begin{lemma}[\cite{Cappell1994OnTM}]\label{lemma_rotation_homotop}  Let $\Lambda_1, \Lambda_2\in\Lambda(n)$. If $J_1, J_2\in\mathcal{J}(\R^{2n}, \omega_0)$, then the paths $\rot\>\![\Lambda_1, \Lambda_2; J_1, \pm]$ and $\rot\>\![\Lambda_1, \Lambda_2; J_2, \pm]$ are homotopic with fixed endpoints.
\end{lemma}

\begin{proof} Due to Remark \ref{rem_positive_negative} we will inspect only the case of positive rotations. First observe that for any $J\in\mathcal{J}(\R^{2n}, \omega_0)$ the positive rotation of $\Lambda_1$ to $\Lambda_2$ by $J$ stays in $\Sigma(\Lambda_1, \Lambda_2)$ for $t\in(0, 1]$.

Since $\mathcal{J}(R^{2n}, \omega_0)$ is contractible, there is a continuous family $\lbrace J^{s}\rbrace_{s\in[1, 2]}$ such that $J^{1}=J_1$ and $J^{2}=J_2$. Since by the Gram-Schmidt process the orthogonal complement depends continuously on the inner product, $\lbrace J^{s}\rbrace_s$ induces a continuous family of $g_{J^s}$-orthogonal splittings $\Lambda_1=(\Lambda_1\cap\Lambda_2)\oplus C^s$. In more detail, let us consider the space $\lbrace P\in\End(\Lambda_1)\,|\,P^2=P, P(\Lambda_1)=\Lambda_1\cap\Lambda_2\rbrace$. Next, if $P_J$ denotes orthogonal projection induced by $J$, then $\lbrace P_{J^s}\rbrace_s$ is the continuous family, where $\ker P_{J^s}=C^s$.

Hence for $\varepsilon>0$ \textit{small} is the map 
\begin{align*}
\gamma: [1, 2]\times[0, 1]&\rightarrow \Lambda(n)\\
(s, t)&\mapsto P_{J^s}(\Lambda_1)+ e^{\varepsilon t J^{s}}(\mathbbm{1}-P_{J^s})(\Lambda_1),
\end{align*}
is continuous, $\gamma|_{[1, 2]\times\lbrace0\rbrace}=\Lambda_1$ and $\gamma|_{[1, 2]\times(0, 1]}\in\Sigma(\Lambda_1, \Lambda_2)$.

Since $\Sigma(\Lambda_1, \Lambda_2)$ is by the Lemma \ref{lemma_contr} contractible, \textit{any} pair of paths from $\gamma(1, 1)$ and $\gamma(2, 1)$ to $\Lambda_2$ makes together with $\gamma(1, t)$ a contractible loop, see Figure \ref{figure_rotation_homotopy}.

\begin{figure}[!htbp]
\labellist
\pinlabel $\Lambda(n)\setminus\Sigma(\Lambda_1,\,\Lambda_2)$ at 51 23
\pinlabel $\Lambda_1$ at 310 -6
\pinlabel $\Lambda_2$ at 360 430
\pinlabel $\gamma(1,\,1)$ at 130 210
\pinlabel $\gamma(2,\,1)$ at 470 230
\pinlabel $\Sigma(\Lambda_1,\,\Lambda_2)$ at 50 350
\endlabellist
\centering
\includegraphics[scale=0.450]{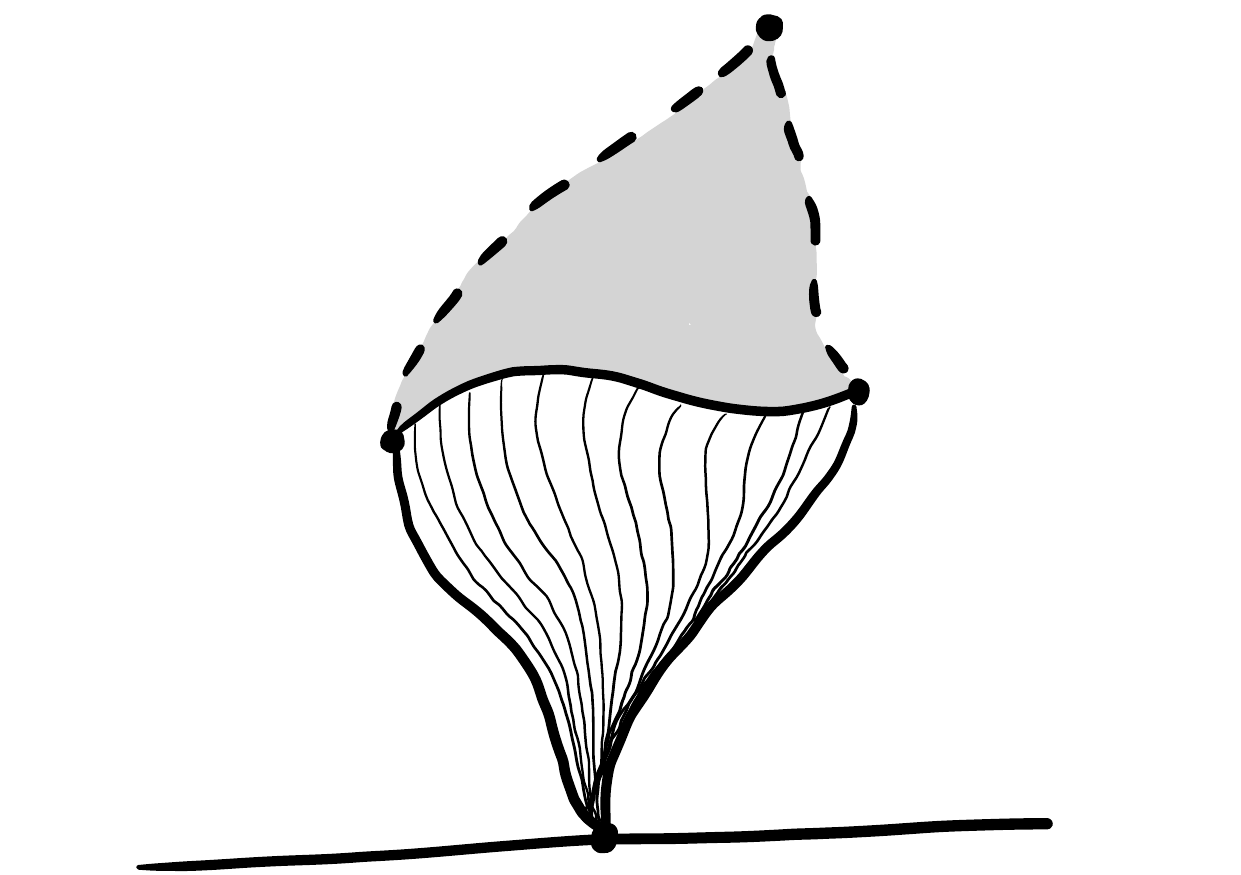}
\caption{The extension of $\gamma$ to the homotopy between $\rot[\Lambda_1, \Lambda_2; J_1, +]$ and $\rot[\Lambda_1, \Lambda_2; J_2, +]$}
\label{figure_rotation_homotopy}
\end{figure}
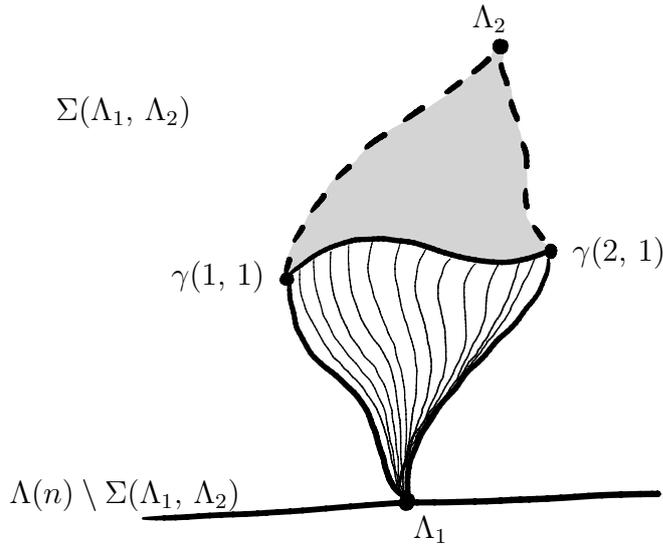

Hence there is a homotopy with fixed endpoints between $\rot\>\![\Lambda_1, \Lambda_2; J_1, +]$ and $\rot\>\![\Lambda_1, \Lambda_2; J_2, +]$.
\end{proof}

\begin{lemma} \label{lemma_maslov_rot} Let $\Lambda_1, \Lambda_2\in\Lambda(n)$ and $J\in\mathcal{J}(R^{2n}, \omega_0)$. Then 
$$\pm\mu(\rot\>\![\Lambda_1, \Lambda_2; J, \pm], \Lambda_1)=\pm\mu(\rot\>\![\Lambda_1, \Lambda_2; J, \pm], \Lambda_2)=\frac{n-k}{2},$$
where $k=\dim(\Lambda_1\cap\Lambda_2)$.
\end{lemma}

\begin{proof}
We will prove only that $\mu(\rot\>\![\Lambda_1, \Lambda_2; J, +], \Lambda_1))=(n-k)/2$. Then $\mu(\rot\>\![\Lambda_1, \Lambda_2; J, -], \Lambda_2))=(k-n)/2$ can be shown by the analogous computation, and the rest follows from the Remarks \ref{rem_opposite_path}, \ref{rem_positive_negative}.

First, we would like to show that $\mu(\rot\>\![\Lambda_1, \Lambda_2; J, +], \Lambda_1))=(n-k)/2$.

Note that that there is a symplectic basis $\lbrace v_1,\dots, v_n, w_1,\dots, w_n\rbrace$ such that $\Lambda_1=\lbrace v_1,\dots, v_n\rbrace$ and $\Lambda_2=\lbrace v_{n-k+1},\dots, v_n, w_1,\dots, w_{n-k}\rbrace$. By the product and naturality properties of the Maslov index we can ignore the part with the constant rotation and consider only $\Lambda_1$ as the $\lbrace x_i\rbrace_{i=1}^{n-k}\text{--plane}$ and $\Lambda_2$ as the $\lbrace y_i\rbrace_{i=1}^{n-k}\text{--plane}$. Moreover, by Lemma \ref{lemma_rotation_homotop} we can assume that $J=J_0$. 

Hence we have
$$\mu(\rot\>\![\Lambda_1, \Lambda_2; J, +], \Lambda_1))=\mu(e^{tJ_0}\lbrace x_i\rbrace_i\text{--plane}, \lbrace x_i\rbrace_i\text{--plane}).$$
Observe that the only crossing appears at $t=0$. Next, let $\varepsilon>0$ small, then for $t\in [0, 1)$ we have
\begin{align*}
\rot\>\![\Lambda_1, \Lambda_2; J, +](t)&=\lbrace (\cos(t)a, \sin(t)\,a)\,|\,a\in\lbrace x_i\rbrace_i\text{--plane}\rbrace\\
&=\lbrace (a, \tan(t)\,a)\,|\,a\in\lbrace x_i\rbrace_i\text{--plane}\rbrace.
\end{align*}
Hence, the corresponding crossing form at $t=0$ is
$$\langle\cdot, \left. \frac{d}{dt} \right\rvert_{t = 0}\tan(t)\cdot\rangle=\langle\cdot,\cdot\rangle$$
and $\mu(\rot\>\![\Lambda_1, \Lambda_2; J, +], \Lambda_1))=(n-k)/2$ by Formula (\ref{rel_maslov}).

The case $\mu(\rot^{\,\Lambda_1,\Lambda_2}_{J, -}, \Lambda_1))$ follows from the analoguos computations. In the remaining cases, we may choose a symplectic basis such that $\Lambda_2$ is the $\lbrace x_i\rbrace_{i=1}^{n-k}\text{--plane}$ (instead of $\Lambda_1$) and the rest is analogous.
\end{proof}

\begin{defn} The \textbf{totally real Grassmanian} $R(n)$ is the manifold consisting of all totally real subspaces of $(\mathbb{C}^{n}, i)$.
\end{defn}

\begin{rem} Analogously as for $\Lambda(n)$ in Lemma \ref{Lemma_arnold_maslov} we see that $R(n)$ has the structure of the homogenuous space. That is
$$R(n)=GL(n, \mathbb{C})/GL(n, \R).$$

Let us naively inspect the stratification of $R(n)$. Let $\Sigma_k^{real}(R)$ denote the subset of $R(n)$ of all totally real subspaces that have $k$-dimensional intersection with the totally real subspace $R$. Take $\Sigma_0^{real}(i\R^2)$, then from looking on graphs over $\R^2$ we have
$$\Sigma_0^{real}(i\R^2)\cong\lbrace A\in\hbox{Mat}_{2\times 2}(\R)\,|\,\pm i\notin\hbox{spectrum}(A)\rbrace.$$
But $\hbox{Mat}_{2\times 2}(\R)\cong\R^4$ and from the characteristic polynomial we see that the space of $A$ such that $\pm i\in \hbox{spectrum}(A)$ has the dimension equal to $2$. Hence $\Sigma_0^{real}(i\R^2)$ is neither contractible nor even only simply connected.

However, by the Gram-Schmidt process, the natural inclusion $\Lambda(n)\hookrightarrow R(n)$ is a deformation retract.

Now we define a map
\begin{align*}
\rho:R(n)&\longrightarrow S^1\\
[A]&\longmapsto\frac{\det(A^2)}{\det(AA^\ast)},
\end{align*}
where $A\in GL(n, \mathbb{C})$ and $\ast$ denotes conjugate transpose. Observe that $\rho$ is well defined.
\end{rem}

\begin{defn}[{\cite[C.3.]{mcduff_salamon_2012}}] \label{defn_masl_2} Let $\gamma$ be a continuous loop in $R(n)$. Then the \textbf{Maslov index of} $\gamma$ is defined as
$$\mu(\gamma):=\degr\circ\,\rho(\gamma).$$
\end{defn}

\begin{rem} \cite{Oh2015Symplectic} Maslov index of a loop in $R(n)$ is also well defined on homotopy classes of loops. In particular, when $\gamma$ is a loop in $\Lambda(n)\subset R(n)$, then Definitions \ref{defn_masl_1} and \ref{defn_masl_2} coincide.
\end{rem}



\chapter{Maslov index of Reeb chords}
\label{ch:file2}
 
\section{Vector bundles}\label{sect_vect}
\begin{defn} An \textbf{orthogonal vector bundle} $(E\rightarrow M, g)$ is a continuous real vector bundle with an \textbf{orthogonal structure} $g$, that is a continuous assignment of an inner product $g_x$ to the fibers $E_x$, for each $x\in M$.

Moreover, if $F$ is a subbundle of $E$, then the subbundle $F^\bot$ is defined fiberwise as
$$F^\bot:=\lbrace v\in E_x\,|\,g_x(v, w)=0\hbox{ for all }w\in E_x\rbrace.$$
\end{defn}

\begin{rem}\label{rem_complement} Observe that $F\oplus F^\bot\cong E.$
\end{rem}

\begin{defn} \label{defn_vect_bundle} A \textbf{symplectic vector bundle} $(E\rightarrow M, \omega)$ is a continuous real vector bundle with a \textbf{symplectic structure $\omega$}, that is a continuous assignment of (linear) symplectic forms $\omega_x$ to the fibers $E_x$ for each $x\in M$.

A \textbf{Lagrangian subbundle} $(L\rightarrow M)$ is a subbundle of $E$ with fibers consisting of Lagrangian subspaces.

Two symplectic vector bundles $(E_1\rightarrow M_1, \omega_1)$ and $(E_2\rightarrow M_2, \omega_2)$ are \textbf{isomorphic}, if there is a vector bundle isomorphism $\psi:E_1\rightarrow E_2$ such that $\psi^\ast\omega_2=\omega_1$.

Next, two Lagrangian subbundles $L_1\subset E_1$ and $L_2\subset E_2$ are \textbf{isomorphic} if there is an isomorphism $\psi$ of symplectic bundles $E_1$ and $E_2$ such that $\psi L_1=L_2$.

A \textbf{complex vector bundle} $(E\rightarrow M, J)$ is a continuous real vector bundle with a \textbf{complex structure} $J$, that is a continuous assignment of linear complex structures $J_x$ to $E_x$ for each $x\in M$.

Two complex vector bundles $(E_1\rightarrow M_1, J_1)$ and $(E_2\rightarrow M_2, J_2)$ are \textbf{isomorphic} if there is a vector bundle isomorphism $\psi:E_1\rightarrow E_2$ such that $\psi^\ast J_2=J_1$.

A complex structure $J$ is called \textbf{compatible with $\omega$}, if $J_x$ is compatible with $\omega_x$ on $E_x$ for each $x\in M$. The \textbf{Space of all $J$ compatible with} $\omega$ will be denoted by
$\mathcal{J}(E, \omega).$ Analogously we define the space 
$$\Omega(E, J):=\lbrace \omega\hbox{ symplectic structure on }E\,|\,J\in\mathcal{J}(E, \omega)\rbrace.$$

A \textbf{Hermitian vector bundle} $(E, \omega, J)$ is a real vector bundle $E$ with a symplectic structure $\omega$ and compatible complex structure $J$. Then a \textbf{Hermitian structure} $h$ is the Hermitian inner product on each fibre induced by $\omega$ and $J$.
\end{defn}

\begin{rem} \cite[ch. 2.6]{mcduff_salamon_2017} Note that symplectic vector bundles have local trivializations with fibers $(\R^{2n}, \omega_0)$. This is done by a parametric version of the choice of a symplectic basis as in the case $(V, \omega)$. Then the transition functions are elements of $Sp(2n, \R)$. Conversely, $Sp(2n, \R)$-vector bundles uniquely determine symplectic vector bundles up to isomorphism. The analogous statement is also true for the complex vector bundles and $GL(n, \mathbb{C})$-vector bundles, for Hermitian vector bundles and $U(n)$-vector bundles, etc.
\end{rem}

\begin{rem} \cite[ch. 2.6]{mcduff_salamon_2017}\label{rem_contr} $\mathcal{J}(E, \omega)$ and $\Omega(E, J)$ are nonempty and contractible. Also, if $R\subset(E, J)$ is a totally real subbundle, then the space of $\omega\in\Omega(E, J)$ such that $J$ is $\omega$-Lagrangian is nonempty and contractible.
\end{rem}

\begin{rem} \cite{hirsch_1997}\label{rem_pullback} Since we would like to work with various trivializations of pullback bundles, it will be useful to recall the definition of the pullback bundle together with its charts.

Let $(E\xrightarrow{\pi} M)$ be a (real) rank $k$ vector bundle together with the maximal atlas $\mathcal{A}$ and let $f:N\rightarrow M$ be a continuous map between topological spaces $N, M$. Then the \textbf{pullback bundle} $f^\ast E$ is given as a vector bundle $(E_0\xrightarrow{\pi_0} M_0)$ together with a maximal atlas $\mathcal{A}_0$, where
\begin{align*}
& E_0:=\lbrace (e, n)\in E\times N\,|\,\pi(e)=f(n)\rbrace,\\
& M_0:=N,\\
& \pi_0:(e, n)\mapsto n.
\end{align*}
The maximal atlas $\mathcal{A}_0$ contains all charts $(\varphi_0, U_0)$ that are constructed as follows. If $(\varphi, U)\in\mathcal{A}$ and $f^{-1}(U)$ is non empty, then $U_0:=f^{-1}(U)$ and
$$\varphi_0:\pi_0^{-1}(U_0)\rightarrow N\times\R^k, (e, n)\mapsto(n, pr_2\circ\varphi(e)).$$
The analogous construction also holds for $G$-vector bundles with any structure group $G$.
\end{rem}

\begin{rem} \label{rem_compos_pullback} Let $E\rightarrow A_3$ be a $G$-vector bundle and consider continuous maps $f_1: A_1\rightarrow A_2$, $f_2:A_2\rightarrow A_3$ and $f_3:A_1\rightarrow A_3$ such that $f_3=f_2\circ f_1$. Then
\begin{equation*}
f^\ast_1 f^\ast_2 E=f_3^\ast E.
\end{equation*}
\end{rem}

\begin{rem}\label{rem_sum_pullback} If $E_1$ and $E_2$ are two $G$-vector bundles over $M$ and $f$ is a continuous map $f:N\rightarrow M$, then
$$f^\ast (E_1\oplus E_2)\cong f^\ast(E_1)\oplus f^\ast(E_2).$$
\end{rem}

\begin{thm}[\cite{bott_tu_2011, Steenrod_23}]\label{thm_triv_g_bundles} Let $M, N$ be topological spaces, where $N$ is moreover a topological manifold, and $f_1, f_2$ are two continuous maps from $N$ to $M$. If $f_1$ and $f_2$ are homotopic and $E\rightarrow M$ is a real vector bundle, then the bundles $f_1^\ast E$ and $f_2^\ast E$ are isomorphic. Moreover, the same statement also holds for any $G$-vector bundles. 
\end{thm}

\begin{rem} \label{rem_triv_homotop} Based on Theorem \ref{thm_triv_g_bundles} we can make the following useful observation.

Let $A$ be a manifold that is a deformation retract onto $B$, where $i$ and $r$ are the corresponding inclusion and retraction, respectively. If $E\rightarrow A$ is a $G$-vector bundle then
$$E\cong (\mathbbm{1}_A)^\ast E\cong (i\circ r)^\ast E=r^\ast i^\ast E=r^\ast E|_B.$$
Hence, any trivialization of $E|_B$ (if it exists) induces a trivialization of $E$ that agrees with the first one when restricted to $B$.

In particular, vector bundles over contractible manifolds are trivial (take for $B$ a point). Moreover, if the structure group is connected, then all trivializations are homotopic.
\end{rem}

\begin{defn} \label{defn_maslov_index_lagr} Let $(E\rightarrow M, \omega)$ be a symplectic vector bundle and $\gamma$ be a continuous positively oriented contractible loop $\gamma: S^1=\lbrace (x_1, x_2)\in\R^2\,|\,x_1^2+x_2^2= 1\rbrace\rightarrow M$, where the orientation of $S^1$ is induced from the canonical orientation of $\R^2$. Let $L^{\gamma}$ be a Lagrangian subbundle of $(\gamma^{\ast}E, \gamma^\ast\omega)$. Take some capping disc $\overline{\gamma}$ of $\gamma$, that is a continuous map $\overline{\gamma}:D=\lbrace (x_1, x_2)\in\R^2\,|\,x_1^2+x_2^2\leq 1\rbrace\rightarrow M$ such that $\left.\overline{\gamma}\right|_{S^1}=\gamma$. 
Next, take a symplectic trivialization $\Phi:(\overline{\gamma}^\ast E, \overline{\gamma}^\ast \omega)\rightarrow(D\times\R^{2n}, \omega_0)$.  We define the \textbf{Maslov index of } $L^\gamma$ \textbf{with respect to} $\Phi$ \textbf{and} $\overline{\gamma}$ as
$$\mu(L^\gamma):=\mu(\Phi L^\gamma).$$

The Maslov index of $L^\gamma$ is \textbf{well-defined} if it does not depend on the choice of $\Phi$ and $\overline{\gamma}$.
\end{defn}

\begin{rem} \label{rem_maslov_totally_real} Definition \ref{defn_maslov_index_lagr} can be analogously rewritten in the language of complex vector bundles and totally real subbundles.
\end{rem}

\begin{rem} \label{rem_well_def} First, we would like to discuss the well--definiteness of $\mu(L^\gamma)$ for fixed $\overline{\gamma}$. But, since $D$ is contractible and $Sp(2n, \R)$ is connected, all trivializations $\Phi$ are homotopic. Hence, all $\mu(\Phi L^\gamma)$ are equal.

Moreover, if $\psi: E\rightarrow E^\prime$ is an isomorphism of symplectic vector bundles over $D$, then $\mu(L)=\mu(\psi L)$. Indeed, $\Phi\psi^{-1}$ is a symplectic trivialization of $E^\prime$ and $\mu(\psi L)$ does not depend on the choice of the trivialization.

However, to investigate the dependence of the index on different choices of capping discs, we first need to introduce the first Chern number.
\end{rem}

\begin{rem} We take a complex vector bundle over the sphere $(E\rightarrow S^2=\lbrace (x_1, x_2, x_3)\in\R^3\,|\,x_1^2+x_2^2+x_3^2=1\rbrace, J)$ of real rank $2n$. Note that $S^2$ is a union of the upper and the lower hemisphere: $S^2=S^+\cup S^-$. Also, the intersection of these hemispheres is a circle $S^1$, which can be parametrized with respect to the canonical orientation of $S^2$ induced from $\R^3$ in $t\in[0,1]$ as $(\cos t, \sin t, 0)$.

Then we choose $\omega\in\mathcal{J}(E, \omega)$ and obtain a Hermitian bundle $(E, \omega, J)$. Next, since the hemispheres $S^\pm$ are contractible, there are unitary trivializations $\Phi^\pm:(E|_{S^\pm}, J)\rightarrow S^\pm\times (\mathbb{R}^{2n}, J_0)$. Restricting to $S^1=\R/\mathbb{Z}$, we have a continuous family of unitary linear maps
$\Phi^\pm_t:E_{[t]}\rightarrow \mathbb{R}^{2n}.$ Hence, for any $t\in [0, 1]$, $\Phi^+_t(\Phi^-_t)^{-1}\in U(n)$.
\end{rem}

\begin{defn} The \textbf{first Chern number} of the complex vector bundle $(E\rightarrow S^2, J)$ is defined as
$$c_1(E):=\deg\circ\det\left(\Phi^+_t(\Phi^-_t)^{-1}\right),$$
for $t\in [0, 1]$ parametrizing the circle $S^1$.
\end{defn}

\begin{rem}From the fibration $SU(n)\rightarrow U(n)\xrightarrow{\det} S^1$ and the LES for homotopy groups, we have that the map $\det_\ast:\pi_1(U(n))\rightarrow\pi_1(S^1)$ is an isomorphism.

By Remark \ref{rem_contr} and contractibility of $S^\pm$ we obtain that all choices of $\omega, \Phi^+$ and $\Phi^-$ are homotopic. Hence, we have a well-defined $c_1$ for isomorphism classes of complex vector bundles over $S^2$.

If $E$ is trivial, then $\Phi^+_t(\Phi^-_t)^{-1}=\mathbbm{1}$ and $c_1$ vanishes. 
\end{rem}

\begin{rem} By Remark \ref{rem_contr} we can analoguosly define $c_1$ for symplectic bundles over $S^2$.
\end{rem}

\begin{rem} \label{rem_ext_c1}Let $(E\rightarrow M, \omega)$ be a symplectic bundle and $f:S^2\rightarrow M$ an orientation preserving homeomorphism, then we put $c_1(E):=c_1(f^\ast E)$.
\end{rem}

\begin{rem} \label{rem_chern} Let $(E\rightarrow S^2, \omega)$ be a symplectic vector bundle and $L$ be a Lagrangian subbundle of $E|_{S^1}$. Then the unitary trivializations $\Phi^\pm$ give rise to the loops $\Phi^\pm_tL_{[t]}\in\Lambda(n)$ for $t$ parametrizing $S^1$. Now we would like to investigate how are the Maslov indices $\mu(\Phi^\pm L)$ related?

By the naturality and homotopy property of the Maslov index, we can assume that $\Phi^+_0=\Phi^-_0$ and $\Phi_1^-L_{[1]}=\lbrace x_i\rbrace_i$--plane in $\R^{2n}$.

Next, observe that for each $t\in[0, 1]$ it holds $$\Phi^+_tL_{[t]}=\Phi^+_{t}\left(\Phi^-_{t}\right)^{-1}\Phi^-_{t}L_{[t]}.$$
Hence, we can compute
\begin{align*} \mu(\Phi^+L)&=\mu(\Phi^+_{t}L_{[t]})\\
&=\mu\left(\Phi^+_{t}\left(\Phi^-_{t}\right)^{-1}\Phi^-_{t}L_{[t]}\right)\\
&=\mu\left(\Phi^+_0\left(\Phi^-_0\right)^{-1}\Phi^-_{t}L_{[t]}\right)+\mu\left(\Phi^+_{t}\left(\Phi^-_{t}\right)^{-1}\Phi^-_1L_{[1]}\right)\\
&=\mu\left(\Phi^-_{t}L_{[t]}\right)+\deg\circ\det\left(\Phi^+_{t}\left(\Phi^-_{t}\right)^{-1}\right)^2\\
&=\mu\left(\Phi^-_{t}L_{[t]}\right)+2\deg\circ\det\left(\Phi^+_{t}\left(\Phi^-_{t}\right)^{-1}\right)\\
&=\mu\left(\Phi^-_{t}L_{[t]}\right)+2c_1(E)\\
&=\mu\left(\Phi^-L\right)+2c_1(E),
\end{align*}
where in the third equality, we used the concatenation and homotopy property of the Maslov index.
\end{rem}

Now, we are ready to continue the discussion from Remark \ref{rem_well_def}.

\begin{rem} As in Definition \ref{defn_maslov_index_lagr} we take a symplectic vector bundle $(E\rightarrow M, \omega)$, a loop $\gamma$ in $M$ and a Lagrangian subbundle $L^\gamma$ of $\gamma^\ast E$. How does the Maslov index $\mu(L^\gamma)$ depend on the choice of the capping disc?

Let $\overline{\gamma}^+$, $\overline{\gamma}^-$ denote two arbitrarily capping discs. They are continuous maps
$\overline{\gamma}^\pm:D\rightarrow M$ such that $\overline{\gamma}^+(\cos t, \sin t)=\overline{\gamma}^-(\cos t, \sin t)=\gamma(\cos t, \sin t)$ for each $t\in[0, 1]$ (again, here we used parametrizations with respect to the canonical orientation of $D$). 

Next, since $\overline{\gamma}^+$ and $\overline{\gamma}^-$ agree on $\partial D=S^1$, by the Pasting lemma, there exists the following gluing along $S^1$. We introduce a continuous function 
$$\pmb{\bm{\gamma}}:=\overline{\gamma}^+\cup_{S^1}(\overline{\gamma}^-)^{op}:D\cup_{S^1}D^{op}\rightarrow M,$$
where $(\overline{\gamma}^-)^{op}$ denotes capping disc $\overline{\gamma}^-$ with opposite orientation. Clearly, there is a canonical orientation-preserving homeomorphism between $D\cup_{S^1}D^{op}$ and $S^2$. Hence, by Remark \ref{rem_ext_c1} we get to the situation as in Remark \ref{rem_chern}. Hence, we obtain the following corollary.
\end{rem}


\begin{cor} \label{cor_chern} Let $(E\rightarrow M, \omega)$ be a symplectic vector bundle, $\gamma$ be a contractible loop in $M$ and $L^\gamma$ be a Lagrangian subbundle of $(\gamma^\ast E, \gamma^\ast\omega)$. If for any two capping discs $\gamma^\pm$ it holds that $c_1(\pmb{\bm{\gamma}}^\ast E)=0$, then the Maslov index $\mu(L^\gamma)$ does not depend on the choice of the capping disc and is well-defined.

\end{cor}

\begin{rem} \label{rem_whitney_sum} If $E_1$ and $E_2$ are two symplectic vector bundles over $S^2$, then a straightforward computation shows that for the Whitney sum $E_1\oplus E_2$ it holds that $c_1(E_1\oplus E_2)=c_1(E_1)+c_1(E_2)$ (see for example \cite{otto_2019}).
\end{rem}

\section{Basics of symplectic and contact manifolds}
\begin{defn} The pair $(M, \omega)$ is called a \textbf{symplectic manifold}, if $\omega$ is a \textbf{symplectic form} on the smooth manifold $M$, that is a closed nondegenerate $2$--form on $M$.
\end{defn}

\begin{defn} The \textbf{Liouville form} $\lambda$ on the symplectic manifold $(M, \omega)$ is a one form $\lambda$ such that $\omega=d\lambda$. The \textbf{Liouville vector field $X$ associated to $\lambda$} is the unique solution of the equation $i_X\omega=\lambda$.
\end{defn}

\begin{defn} Two symplectic manifolds $(M_1, \omega_1)$ and $(M_2, \omega_2)$ are caled \textbf{symplectomorphic}, if there is a diffeomorphism $\psi:M_1\rightarrow M_2$ such that $\psi^\ast\omega_2=\omega_1$. Here $\psi$ is called a \textbf{symplectomophism}.

Moreover, let $\lambda_1$ and $\lambda_2$ are Liouville forms on $M_1$ and $M_2$, respectively. If $\psi^\ast \lambda_2=\lambda_1$, then $\psi$ is called a \textbf{Liouville diffeomorphism}.
\end{defn}

\begin{Darboux}[\cite{mcduff_salamon_2017}]
Let $(M_1, \omega_1)$ and $(M_2, \omega_2)$ be two symplectic manifolds. If $x_1$ and $x_2$ are two points in $M_1$ and $M_2$ respectively, then there exist two neighbourhoods $U_1$ and $U_2$ of these points in $M_1$ and $M_2$ respectively such that $(U_1, \omega_1|_{U_1})$ is symplectomorphic to $(U_2, \omega_2|_{U_2})$. In other words, all symplectic manifolds of the same dimension are locally symplectomorphic.
\end{Darboux}

\begin{defn} Let $H$ be a smooth function on the symplectic manifold $(M, \omega)$. Then the \textbf{Hamiltonian vector field $X_H$ associated to $H$} is the unique solution of $\omega(X_H, \cdot)=-dH$. The induced flow of $X_H$ is called a \textbf{Hamiltonian flow asociated to $H$} and denoted by $\phi^t_H.$
\end{defn}

\begin{lemma} \label{lemma_ham_flow} $\lbrace \phi^t_H\rbrace_t$ is a smooth family of symplectomorphisms.
\end{lemma}

\begin{proof}
We would like to show that $\frac{d}{dt}(\phi^t_H)^\ast\omega=0$. Then it follows that $(\phi^t_H)^\ast\omega=(\phi^0_H)^\ast\omega=\omega$.

We compute
\begin{align*}
\frac{d}{dt}\Bigr|_{t=t_0}(\phi^t_H)^\ast\omega&=(\phi^{t_0}_H)^\ast\left(\frac{d}{dt}\Bigr|_{t=0}(\phi^t_H)^\ast\omega\right)\\
&=(\phi^{t_0}_H)^\ast\left(di_{X_H}\omega+i_{X_H}d\omega\right)\\
&=0,
\end{align*}
where in the second equality we used Cartan´s magic formula.
\end{proof}

\begin{defn} Let $(M, \omega)$ be a $2n$-dimensional symplectic manifold. Then its submanifold $L$ is called \textbf{Lagrangian} if $\dim L=n$ and $\omega|_{TL}=0$.
\end{defn}

\begin{defn} Let $(M, \omega)$ be a symplectic manifold with two Lagrangian submanifolds $L_1, L_2$. We say that \textbf{$L_1$ and $L_2$ intersect cleanly along a submanifold $K$} if 
$L_1\cap L_2=K$ and $T_xK=T_xL_1\cap T_xL_2$ for each $x\in K$.
\end{defn}

\begin{defn} Let $M$ be a smooth $(2n+1)$--dimensional manifold. The \textbf{contact form} $\alpha$ is a $1$--form on $M$ such that
$$\alpha\wedge(d\alpha)^n\neq 0.$$
A hyperlane field $\xi$ on $M$ is called a \textbf{contact structure}, if for every point $x\in M$ there is a neighborhood $U_x\subseteq M$ and a contact form $\alpha$ defined on $U_x$ such that $\xi|_{U_x}=\ker \alpha\subseteq TU_x$. 

The pair $(M, \xi)$ is called a \textbf{contact manifold}.
\end{defn}

\begin{lemma}[\cite{geiges_2008}]\label{lemma_contact}
A contact structure $\xi$ is a kernel of some global contact form if and only if $\xi$ is coorientable.

If the contact form $\xi$ is a kernel of some contact form $\alpha$, then any nonvanishing function $f$ defines a contact form $f\alpha$ with $\xi=\ker f\alpha$. Moreover, $(\xi, d\alpha|_\xi)$ is a symplectic vector bundle.
\end{lemma}

\noindent \textit{In the view of Lemma \ref{lemma_contact}, we will be studying only coorientable contact structures. Hence, we can write $(M, \alpha)$ as a contact manifold.}

\begin{defn}
Two contact manifolds $(M_1, \xi_1=\ker\alpha_1)$ and $(M_2, \xi_2=\ker\alpha_2)$ are called \textbf{contactomorphic}, if there is a diffeomorphism $\psi:N_1\rightarrow N_2$ such that $\psi_\ast(\xi_1)=\xi_2$. Here $\psi$ is called a \textbf{contactomorphism}.

Equivalently, we say that the diffeomorphism $\psi$ is a \textbf{contactomorphism}, if $\psi^\ast\alpha_2=f\alpha_1$ for some nonvanishing function $f$ on $M_1$.

Moreover, if $f=1$, then $\psi$ is called a \textbf{strict contactomorphism}.
\end{defn}

\begin{Pfaff}[\cite{geiges_2008}]\label{thm_pfaff}
All contact manifolds of the same dimension are locally strictly contactomorphic.
\end{Pfaff}

\begin{defn} Let $(M, \alpha)$ be a contact manifold. Then the \textbf{Reeb vector field $R$ asociated to $\alpha$} is the unique solution of:
\begin{itemize}
\item[(i.)] $i_R\alpha=1,$
\item[(ii.)] $i_Rd\alpha=0.$
\end{itemize}
The induced flow of $R$ is called the \textbf{Reeb flow} and denoted by $\phi^t_R.$
\end{defn}

\begin{lemma} $\lbrace \phi^t_R\rbrace_t$ is a smooth family of contactomorphisms.
\end{lemma}

\begin{proof}
The proof is analogous to the proof of the Lemma \ref{lemma_ham_flow}. Only here we have
$$\frac{d}{dt}\Bigr|_{t=0}(\phi^t_R)^\ast\alpha=di_R\alpha+i_Rd\alpha=0.$$
\end{proof}

\begin{defn} Let $(M, \alpha)$ be a $(2n+1)$-dimensional contact manifold. Then its submanifold $\mathcal{L}$ is called \textbf{Legendrian} if $\dim \mathcal{L}=n$ and $\alpha|_{T\mathcal{L}}=0$.

Two Legendrian submanifolds are called \textbf{Legendrian isotopic} if there is a smooth ambient isotopy between them through Legendrian submanifolds.
\end{defn}

\begin{defn}
Let $(M, \alpha)$ be a contact manifold with a Legendrian submanifold $\mathcal{L}.$ A curve $\gamma_R:[0, T]\rightarrow M$ is called a \textbf{Reeb chord on $\mathcal{L}$} if it has endpoints on $\mathcal{L}$ and it is an integral of the Reeb flow $\phi_R^t$.

Next, $\gamma_R$ is called \textbf{pure} if its endpoints lie on the same connected component of the Legendrian submanifold $\mathcal{L}$. Otherwise $\gamma_R$ is called \textbf{mixed (from $\mathcal{L}_0$ to $\mathcal{L}_1$)}, here $\mathcal{L}_0, \mathcal{L}_1$ are two different connected components of $\mathcal{L}$ such that $\gamma_R(0)\in\mathcal{L}_0$, $\gamma_R(T)\in\mathcal{L}_1$.

Moreover, $\gamma_R$ is called \textbf{nondegenerate} if $T_{\gamma_{R}(T)}\mathcal{L}\pitchfork d\phi_R^TT_{\gamma_{R}(0)}\mathcal{L}$ in $\xi_{\gamma_{R}(T)}$.
\end{defn}

\begin{rem} \label{rem_reeb_triv} Let $\gamma_{R}$ be a pure Reeb chord on $\mathcal{L}$. Then there is a continuous (positively oriented) loop $\gamma:S^1\rightarrow M$, parametrized by $t\in\R/\mathbb{Z}=S^1$, that satisfies the following:
\begin{itemize}
\item[$(i.)$] if $t\in[0, 1/3],$ then $\gamma(t)=c(3t)$, where $c$ is a continuous ``capping'' path $c:[0, 1]\rightarrow\mathcal{L}$ with $c(0)=\gamma_R(T)$ and $c(1)=\gamma_R(0)$,
\item[$(ii.)$] if $t\in[1/3, 2/3],$ then $\gamma(t)=\gamma_R(3Tt-T)$,
\item[$(iii.)$] if $t\in[2/3, 1],$ then $\gamma(t)=\gamma_R(T)$.
\end{itemize}
We take the pullback bundle $(\gamma^\ast\xi, \gamma^\ast(d\alpha|_\xi))$. Next, we define its Lagrangian subbundle $L^\gamma$ as follows (see also Figure \ref{figure_reeb_grad}):
\begin{itemize}
\item[$(i.)$] for $t\in[0, 1/3]$ we put  $L^\gamma_{[t]}:=T_{\gamma(t)}\mathcal{L}$
\item[$(ii.)$] for $t\in(1/3, 2/3]$ we put $L^\gamma_{[t]}:=d\phi_R^{(3Tt-T)}T_{\gamma_R(0)}\mathcal{L}$
\item[$(iii.)$] for $t\in(2/3, 1)$ we put $L^\gamma_{[t]}:=\rot\>\![d\phi_R^{T}T_{\gamma_R(0)}\mathcal{L}, T_{\gamma_R(T)}\mathcal{L}; J, +](3t-2)$, where $J$ is a complex structure in $\mathcal{J}(\xi_{\gamma_R(T)}, d\alpha|_{\xi_{\gamma_R(T)}})$.
\end{itemize}

\begin{figure}[!htbp]
\labellist
\pinlabel $\mathcal{L}$ at 500 440
\pinlabel $c$ at 60 120
\pinlabel $\gamma_R$ at 460 190
\pinlabel $\overline{\gamma}$ at 270 200
\endlabellist
\centering
\includegraphics[scale=0.45]{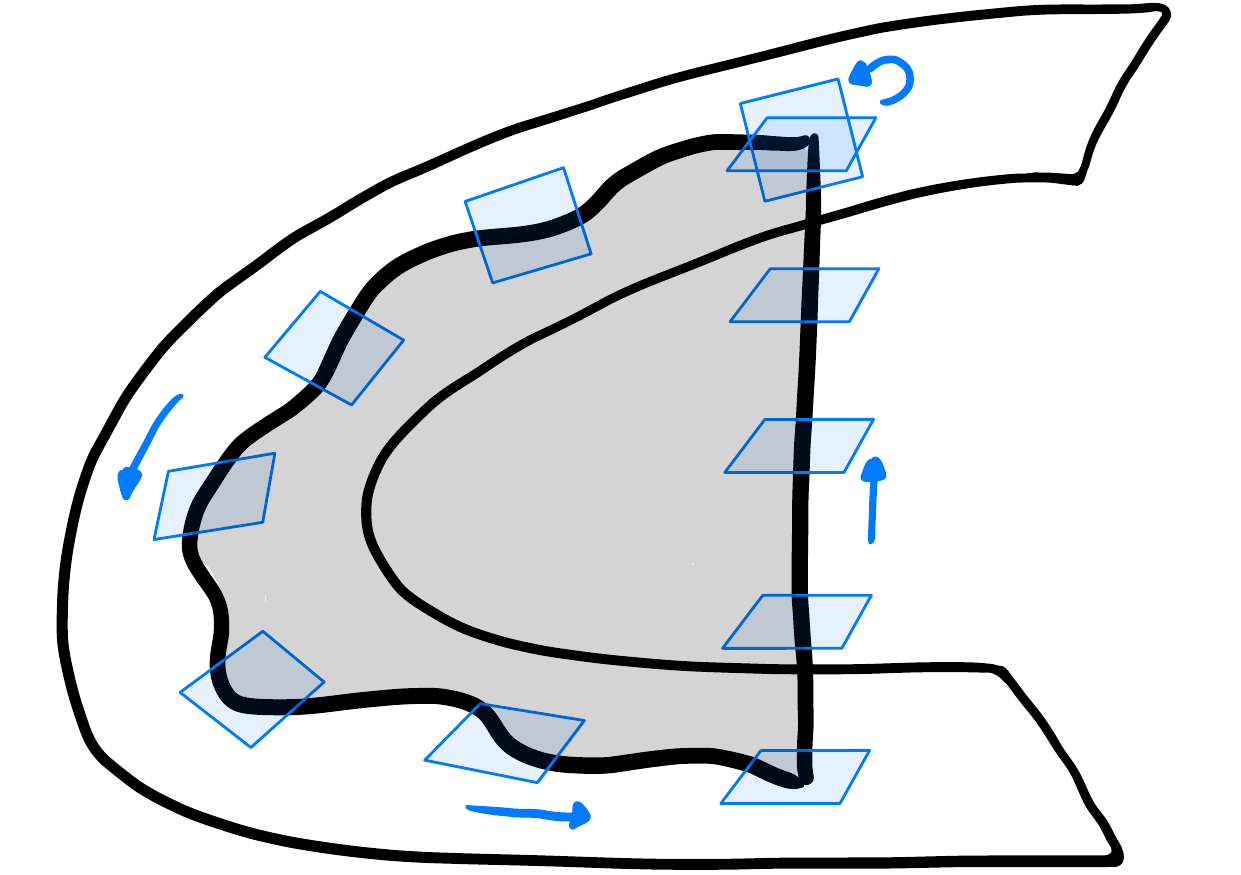}
\caption{The Lagrangian subbundle $L^\gamma$ of $\gamma^\xi$ (in blue) and some capping disc $\overline{\gamma}$ of $\gamma$ (in grey).}
\label{figure_reeb_grad}
\end{figure}
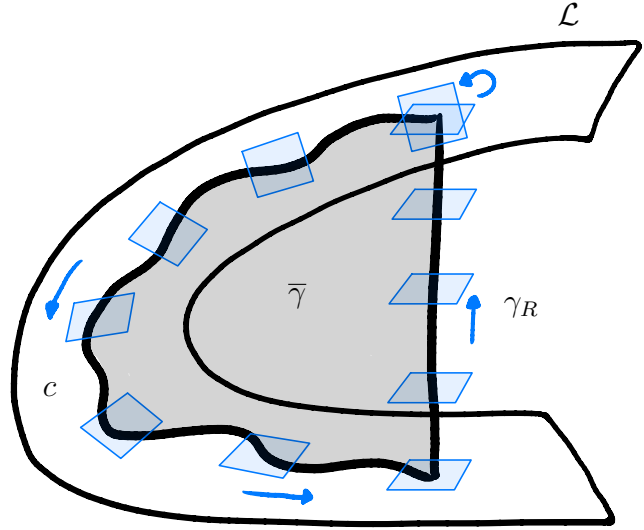

\end{rem}

\begin{defn} \label{defn_reeb_pure} Let $\gamma_R$ and $L^\gamma$ are as in Remark \ref{rem_reeb_triv}. Then the \textbf{Maslov index (or degree) of the pure Reeb chord} $\gamma_{R}$ is defined as 
$$|\gamma_{R}|=\mu(L^\gamma)-1.$$
We say that $|\gamma_R|$ is \textbf{well-defined}, if $\mu(L^\gamma)$ is well-defined in the sense of Definition \ref{defn_maslov_index_lagr} and moreover $\mu(L^\gamma)$ does not depend on the choices of the capping path $c$ and the complex structure $J$.
\end{defn}

\begin{rem} From Lemma \ref{lemma_rotation_homotop} we obtain that the Maslov index of a pure Reeb chord does not depend on $J$.
\end{rem}

\section{On cotangent bundles}
\begin{example} A classical example of the symplectic manifold is a \textbf{cotangent bundle} $(T^{\ast}N\xrightarrow{\pi}N)$. To see this, we introduce first a local coordinates $(q^i\circ\pi, p_i)_i$ on $T^\ast U\subset T^\ast N$. Let $(q^i)_i$ be local coordinates on an open neighborhood $U\subset N$. Then $(p_i)_i$ are defined as follows. If $x\in T^\ast U$, then $\lbrace dq^i|_{q}\rbrace_i$ form a basis of $T_{q}^\ast N$ for $q=\pi(x)$. Therefore there are unique functions $p_i:x\mapsto x(\partial_{q^i}|_q)$. For simplicity, we will write $q^i:= q^i\circ\pi$. Also, we will denote points of $T^\ast N$ as $x=(q, p)$.

Then we will define on $T^\ast N$ a $1$-form $\lambda$. For each $x\in T^\ast N$, the linear map $\lambda_x:T_xT^\ast N\rightarrow\R$ is given by $\lambda_x:=x\circ d\pi$. It will be useful to have a description of $\lambda$ in local coordinates. Since $d\pi(x)\partial_{q^i}|_x=\partial_{q^i}|_q$ and $d\pi(x)\partial_{p_i}|_x=0$, we conclude that $\lambda=\sum_i p_idq^i$.

Now, we put $\omega=d\lambda$, which can be locally written as $\omega=\sum_i dp_i\wedge dq^i$. Since $\omega$ is closed and nondegenerate, it is a symplectic form and $\lambda$ is a Liouville form on $T^\ast N$.
\end{example}

\begin{rem} \label{rem_ham} Since the cotangent bundle $T^\ast N$ has the global Liouville form $\lambda$, we can introduce a global Liouville vector field $X$, which is locally given by $X=\sum_i p_i\partial_{p_i}$.

Let $g$ be a Riemannian metric on $N$. Then we define a function $H$ on $T^\ast N$ as 
\begin{equation}\label{eqn_ham_fce}
H(x):=(1/2)||\,x^\sharp\,||_g^2.
\end{equation}
The corresponding Hamiltonian vector field can be expressed as follows. Let $x\in T^\ast N$, then we take a normal coordinates $(q^i)_i$ centered at $q$ and induced coordinates $(p_i)_i$. Then $(X_H)_x=\sum_i p_i(x)\partial_{q^i}|_x$. 
\end{rem}

\noindent\textit{The only Hamiltonian function on $T^\ast N$ we will consider is the function (\ref{eqn_ham_fce}).}

\begin{defn} Let $Q$ be a submanifold without boundary in $N$. Then the \textbf{conormal bundle of $Q$} is defined as 
$$L^\ast Q:=\lbrace x=(q, p)\in T^\ast N\,|\,q\in Q, x(v)=0 \text{ for all }v\in T_qQ \rbrace.$$
\end{defn}

\begin{rem} \label{rem_conormal} Since $TL^\ast Q$ is subbundle of $TT^\ast N|_{L^\ast Q}$ and $\pi(T^\ast N|_Q)=Q$, we have a projection along fibers $d\pi(TL^\ast Q)=TQ$. Hence, on $TL^\ast Q$ vanish the Liouville form $\lambda$. Since also $\dim L^\ast Q=\dim T^\ast N/2$, $L^\ast Q$ is a Lagrangian submanifold of $T^\ast N$.
\end{rem}

\begin{defn} The \textbf{unit cotangent bundle} on $(N, g)$ is defined as
$$S^\ast N:=\lbrace x\in T^\ast N\,|\,\,||\,x^\sharp\,||_g=1\rbrace.$$
Moreover, let $Q$ be a submanifold without boundary in $N$. Then the \textbf{unit conormal bundle of $Q$} is defined as
$$\mathcal{L}^\ast Q:=S^\ast N\cap L^\ast Q.$$
\end{defn}

\begin{rem} \label{rem_contact_objects} Note that $S^\ast N\subseteq T^\ast N$ is a contact manifold with a contact structure $\xi=\ker \lambda_1$, where $\lambda_1:=\lambda|_{S^\ast N}$. Indeed, by coordinate description of $X$, we see that $X$ is transverse to $S^\ast N$. Hence $\lambda_1\wedge(d\lambda_1)^{n-1}=(1/n)i_X(d\lambda)^n\neq0$ on $S^\ast N$, where $n=\dim N$.

At each point $x\in S^\ast N$, we take a coordinate description of $X_H$ as in Remark \ref{rem_ham}. Observe that $X_H|_x$ satisfies equations for the Reeb vector field on $S^\ast N$. And hence $X_H|_{S^\ast N}=R$.

Thus, we conclude the following. Let $\phi_H^t$ and $\phi_{R}^t$ be Hamiltonian and Reeb flows on $T^\ast N$ and $S^\ast N$ respectively. Then 
\begin{equation}\label{eqn_ham_reeb}
\phi_H^t|_{S^\ast N}= \phi_{R}^{t},
\end{equation}
and hence there is a bijection between integral curves of the Hamiltonian flow $\phi^t_H$ on $S^\ast N\subset T^\ast N$ and integral curves of the Reeb flow $\phi^t_{R}$. These curves are called \textbf{corresponding}.
\end{rem}

\begin{rem} From the coordinate descriptions of $R$ and $X$ we deduce the following. At each point $x$ of $S^{\ast}N$ and any submanifold $L^\ast Q$ of $T^{\ast}N$, we have two splittings:
\begin{align*} T_{x}T^{\ast}N &= \langle X_{x}\rangle\times \langle R_{x}\rangle\times \xi_{x},\\
T_{x}L^\ast Q &= \langle X_{x}\rangle\times T_{x}\mathcal{L}^\ast Q.
\end{align*}
By the Lemma \ref{lemma_contact}, $(\xi, \omega|_\xi)$ is a symplectic vector bundle. By looking on coordinates, we see that $X_x$ and $R_x$ form a symplectic basis of $(\eta_x:=\langle X_{x}\rangle\times \langle R_{x}\rangle, \omega|_{\eta_x})$. This extends to the symplectic vector bundle $(\eta, \omega|_\eta)$.

Also $\langle X_{x}\rangle, \langle R_{x}\rangle$ and $T_{x}\mathcal{L}^\ast Q$ are Lagrangian subspaces when restricted to $\eta_x$ and $\xi_x$, respectively. And analogously these subspaces extend to the Lagrangian subbundles $\langle X\rangle, \langle R\rangle$ and $T\mathcal{L}^\ast Q$ of $\eta$ and $\xi$ respectively.
\end{rem}

\begin{rem} Let $x, y\in S^\ast N$ such that $y=\phi_H^t(x)$ for some $t\in\R$. Then by Formula (\ref{eqn_ham_reeb}) we have $d\phi_H^t(x)R_x=R_y$ and $d\phi_H^t(x)\xi_x=\xi_y$.

Moreover
\begin{equation} \label{eqn_lin_ham_flow_liouv}
d\phi_H^t(x)X_x=X_y+tR_y.
\end{equation}
This is immediate when $N=\R^n$, since then $\phi^t_H(q, p)=(q+t p, p)$. In the general case, it is convenient to do the computation in the symplectization $\R\times S^\ast N$, see also the computation in Remark \ref{rem_proof_eqn_lin_ham_flow_liouv} which works for any $N$. 
\end{rem}

\begin{rem}\cite{geiges_2008} As a contact manifold, $S^\ast N$ can be identified with the space of cooriented contact elements of $N$, which has a natural contact structure. In particular, different choices of $g$ on $N$ give us contactomorphic $S^\ast N$.
\end{rem}

\begin{example}\label{empl_contc} $T^\ast S^{n-1}\times\R$ is a contact manifold with a contact form $dt-\lambda|_{S^\ast S^{n-1}},$ where $t$ denotes $\R$-direction and $\R^n$ is equipped with the standard inner product $\langle\cdot,\cdot\rangle$.

Then $T^\ast S^{n-1}$ can be writen as 
$$T^\ast S^{n-1}=\lbrace(a, b)\in\R^{n}\times\R^{n}\,|\,|a|=1, \langle a, b\rangle=0\rbrace.$$
Now, we can define a map 
\begin{align*}
\psi: S^\ast\R^n&\rightarrow T^\ast S^{n-1}\times \R\\
(q, p)&\mapsto (p, q-\langle q, p\rangle p, \langle q, p\rangle),
\end{align*}
where $(q, p)$ is written with respect to the standard coordinates on $T^\ast \R^n$. Here $q-\langle q, p\rangle p$ is an orthogonal projection of the vector $q_p$ to $T_p^\ast S^{n-1}$. $\langle q, p\rangle$ is then an orthogonal projection of $q_p$ to $N_p^\ast S^{n-1}$.

It follows that $\psi^\ast(dt-\lambda|_{S^\ast S^{n-1}})=\lambda|_{S^\ast \R^n},$ hence it is a (strict) contactomorphism.
\end{example}

\begin{rem}\cite{ferme2017invariants} By the construction, a smooth ambient isotopy between two submanifolds $Q_1, Q_2\in N$ induces a Legendrian isotopy between $\mathcal{L}^\ast Q_1, \mathcal{L}^\ast Q_2 \in S^\ast N.$
\end{rem}

\begin{rem}\cite{ferme2017invariants} Let $N$ be oriented manifold and $Q$ be a $\codim 1$ oriented connected submanifold of $N$. Then $\mathcal{L}^\ast Q$ has two connected components $\mathcal{L}^\ast_+ Q$ and $\mathcal{L}^\ast_- Q$. Both are diffeomorphic to $Q$. Here $\mathcal{L}^\ast_+ Q$ is determined by the positive co-orientation of $Q$.

Naturally, $L^\ast_+ Q$ and $L^\ast_- Q$ will be two Lagrangian submanifolds with boundary $Q$ such that $\mathcal{L}^\ast_\pm Q\subset L^\ast_\pm Q$ and $L^\ast_+ Q\cup_Q L^\ast_- Q= L^\ast Q.$
\end{rem}

\section{Binormal and Hamiltonian chords}
In this section, we would like to describe relations between binormal and Hamiltonian chords. We begin with the discussion of the Morse index of binormal chords. Then we describe the splitting of $TT^\ast N$ into the horizontal and vertical parts. The section will be finished by introducing the Maslov index for Hamiltonian chords.

\begin{lemma}\label{lemma_morse_hilb} Let $Q$ be a submanifold of $(N, g)$ and $T>0$. Then the set
$$\mathcal{H}^1([0, T], N; Q):=\lbrace\gamma\in\mathcal{H}^1([0, T], N)\,|\,(\gamma(0), \gamma(T))\in Q\times Q\rbrace$$
is a smooth Hilbert manifold and at $\gamma\in \mathcal{H}^1([0, T], N; Q)$ the tangent space is
$$T_\gamma\mathcal{H}^1([0, T], N; Q)=\lbrace v\in\mathcal{H}^1([0, T], \gamma^\ast TN)\,|\,(v(0), v(T))\in T_{(\gamma(0), \gamma(T))}(Q\times Q)\rbrace.$$
\end{lemma}

\begin{sproof}For details \cite{Klingenberg1995} or in more general setting \cite{Elasson1967GeometryOM}. Let $\dim N=n$ and $\dim Q=k$.

First, recall the following general fact. If $f:U\subset\R^n\rightarrow U^\prime\subset\R^n$ is a diffeomorphism and $\gamma:[0, T]\rightarrow\R^n$ is of class $\mathcal{H}^1$, then $f\circ\gamma:[0, T]\rightarrow U^\prime$ is also of class $\mathcal{H}^1$.

Let us pick a curve $\gamma\in  \mathcal{H}^1([0, T], N; Q)$. We would like to construct a chart around $\gamma$. First, we take a trivialization $\Phi:\gamma^\ast TN\rightarrow[0, T]\times\R^n$ such that $\Phi(T_{\gamma(0)}Q)=\Phi(T_{\gamma(T)}Q)=\R^k$. 

Now, we pick an auxiliary Riemannian metric $\widetilde{g}$ on $N$. For the simplification, $\widetilde{g}$ will have the property that some \textit{neighbourhoods} of $\gamma(0)$ and $\gamma(T)$ in $Q$ are totally geodesic in $N$. Observe that such a $\widetilde{g}$ always exists. Indeed, first we pick a chart on $N$ around $u(0)$ that maps $Q$ to some vector space. Then the pullback of the Euclidean metric is locally our desired $\widetilde{g}$. Next, we make the same for $u(T)$. Finally, outside some small neighbourhoods of $\gamma(0)$ and $\gamma(T)$ in $N$ we interpolate the pullback metric with any Riemannian metric on $N$, and hence we obtain $\widetilde{g}$ on the whole $N$.

Let $V\subset\gamma^\ast TN$ be a neighbourhood of the zero section such that the map $\exp^{\widetilde{g}}:T_{\gamma(t)}N\subset V\rightarrow N$ is a diffeomorphism onto its image for each $t\in[0, T]$.

Now we put
\begin{align*}
U_\gamma:=\lbrace v\in\mathcal{H}^1([0, T], \R^n)\,|\,(t, v(t))\in\Phi(V)\hbox{ for each }t\in[0, T&],\\
\Phi(v(0))=\Phi(v(T))=\R^k&\rbrace.
\end{align*} 
Observe that $U_\gamma$ is an open subset in 
$$\lbrace v\in\mathcal{H}^1([0, T], \R^n)\,|\,v(0)=v(T)=\R^k\rbrace.$$
The latter is a closed subspace of $\mathcal{H}^1([0, T], \R^n)$ and hence also a Hilbert space. 

Now the charts 
\begin{align*}
\varphi_\gamma:U_\gamma&\longrightarrow\mathcal{H}^1([0, T], N; Q)\\
v&\longmapsto\exp^{\widetilde{g}}(\Phi^{-1}(v))
\end{align*}
define a smooth atlas on $\mathcal{H}^1([0, T], N; Q)$. Moreover, the induced structure does not depend on $\widetilde{g}$.
\end{sproof}

\begin{defn}\label{def_energ_1} We define the \textbf{energy functional} $\bm{E}$ on $\mathcal{H}_1([0, T], N; Q)$ as
$$\bm{E}(\gamma)=\int_0^T||\dot{\gamma}(t)||_g^2dt.$$
\end{defn}

\begin{defn} \label{def_bin_pure} A (non-constant) curve $\gamma_M:[0, T]\rightarrow (N, g)$ is called a \textbf{binormal chord on $Q$} if it is a geodesic with endpoints at $Q$ and at these endpoints $\gamma_M$ is perpendicular to $Q$.
\end{defn}

\begin{rem}\cite{Klingenberg1995} $\bm{E}$ is well defined. That is, $\bm{E}$ does not depend on the choice of charts on $\mathcal{H}_1$. Moreover, $\bm{E}$ is smooth.
\end{rem}

\begin{rem}\cite{Klingenberg1995}\label{rem_crit_path} Critical points of $\bm{E}$ are precisely the binormal chords and constant paths on $Q$.
\end{rem}

\begin{defn}If $\gamma_M\in\mathcal{H}^1([0, T], N; Q)$ is a critical point of $\bm{E}$, then the \textbf{Morse index of} $\gamma_M$ is defined as
\begin{align*}
\hbox{Ind}_{\gamma_M}:=\sup\lbrace\dim V\,|&\,V\hbox{ linear subspace of }T_{\gamma_M}\mathcal{H}^1([0, T], N; Q)\hbox{ on which }\\
&\hbox{the symmetric bilinear form}\
D^2E(\gamma_M)\hbox{ is negative definite}\rbrace.
\end{align*}
Moreover, $\gamma_M$ is called \textbf{nondegenerate} if
\begin{align*}
0=\sup\lbrace\dim V\,|&\,V\hbox{ linear subspace of }T_{\gamma_M}\mathcal{H}^1([0, T], N; Q)\\
&\hbox{ such that }V\hbox{ is a }0\hbox{-eigenspace of }
D^2E(\gamma_M)\rbrace.
\end{align*}
\end{defn}

\begin{Morse}[\cite{Ambrose1961, Klingenberg1995, bolton1977morse}] \label{thm_morse_index}
Let $\gamma_M\in\mathcal{H}^1([0, T], N; Q)$ be a binormal chord. Then
$$\hbox{Ind}_{\gamma_M}=index(H_0-H_T)+\sum_{0<t<T}\lbrace\hbox{multiplicity of the conjugate point }\gamma_M(t)\rbrace,$$
where $H_0, H_T$ are the second fundamental forms of $Q$ at $\gamma_M(0), \gamma_M(T),$ respectively, in directions of the normal vectors $\dot{\gamma}_M(0), \dot{\gamma}_M(T)$, respectively. For definitions of the second fundamental form and the multiplicity of a conjugate point, see \cite{Klingenberg1995}.
\end{Morse}

\begin{defn} \label{def_ham_pure} A curve $\gamma_H:[0, T]\rightarrow T^\ast N$ is called a \textbf{Hamiltonian chord on $L^\ast Q$} if it has endpoints at $L^\ast Q$ and it is an integral curve of the Hamiltonian flow $\phi_H^t$.

Moreover, a Hamiltonian chord $\gamma_H$ is called \textbf{nondegenerate} if $T_{\gamma_{H}(T)}L^\ast Q\pitchfork d\phi_H^TT_{\gamma_{H}(0)}L^\ast Q$ in $T_{\gamma_{H}(T)}T^\ast N$.
\end{defn}
 
\begin{lemma} [{\cite[Thm~2.3.1]{otto_2019}}]\label{lemma_corresp_ham_bin}
If $\gamma_M(t)$ is a binormal chord on $Q$, then the map $t\mapsto(\gamma_M(t), \dot{\gamma}_M(t)^\flat)\in T^\ast N$ defines a Hamiltonian chord $\gamma_H(t)$ on $L^\ast Q$. Moreover, such an assignment defines a bijection between the sets of binormal chords on $Q$ and Hamiltonian chords on $L^\ast Q$. These chords are called \textbf{corresponding}.
\end{lemma}

\begin{rem}\label{rem_split_tangent} For details see \cite[Apx 3]{spivak}. Let $(E\xrightarrow{\pi} N)$ be a vector bundle. 
Applying the tangent functor $T$, we get a vector bundle $(TE\xrightarrow{d\pi} TN)$. 

The \textbf{Vertical subbundle} is the subbundle of $TE$ defined as $T^{v}E:=\ker d\pi$. For any $q\in N$ and $x\in E_q$ we define a map $I_{x}^v$ as
\begin{align*} I_{x}^v:E_q &\longrightarrow T_{x}E\\
\widetilde{x}&\longmapsto\frac{d}{dt}\Big|_{t=0}(x+t\widetilde{x}).
\end{align*}
Observe that $I_{x}^v$ is a vector space isomorphism onto its image, which is $T^v_{x}E$. Hence, we obtain a vector bundle isomorphism $I^v:\pi^\ast(E)\rightarrow T^vE$.

A \textbf{Horizontal subbundle} is a choice of a subbundle $T^hE$ of $TE$ such that $TE=T^vE\oplus T^hE$. Take a Koszul connection $\nabla:\Gamma(E)\rightarrow\Gamma(T^\ast N\otimes E)$. Then for any $q\in N, x\in E_q$ we define a map
\begin{align*} I_{x}^h:T_qN &\longrightarrow T_{x}E\\
v&\longmapsto ds(q)v-(\nabla_vs)(q),
\end{align*}
where $s$ is a section $s:N\rightarrow E$ such that $s(q)=x$ and the map $I^h_x$ does not depend on the choice of $s$. This leads to the vector bundle monomorphism $I^h:\pi^\ast(TN)\rightarrow TE$, where the image of $I^h$ is a horizontal subbundle. Moreover, every horizontal subbundle $T^hE$ is determined by the unique Kozsul connection on $E$.

Hence, we obtain a vector bundle isomorphism $I=I^h\times I^v:\pi^\ast(TN \oplus E)\rightarrow TE$.

By $P=P^h\times P^v:TE\rightarrow\pi^\ast(TN\oplus E)$ we will denote the inverse map to $I$.
\end{rem}

\begin{rem} \label{rem_cotg_split} In our special case we have a vector bundle $(E\rightarrow N)=(T^\ast N\rightarrow N)$ with a symplectic form $\omega$. Let $g$ be a Riemannian metric on $N$. Fix a point $x=(q, p)\in T^\ast N$. Now, we take a normal coordinates $(q^i)_i$ centered at $q$ and induced coordinates $(p_i)_i$. Then 
$\lbrace\partial _{q^i}|_x, -\partial_{p_i}|_x \rbrace_i$ is a symplectic basis of $(T_xT^\ast N, \omega_x)$. Next, $g_q$ induce dual pairing on $T_q N\times T_q^\ast N$. As in Example \ref{example_dual}, this give us a symplectic form $(\omega_\times)_q$ on $T_q N\times T_q^\ast N$ with a symplectic basis $\lbrace \partial_{q^i}|_q,  -dq^i|_q\rbrace_i$.

Moreover, $g_q$ induce on $(T_q N\times T_q^\ast N, (\omega_{\times})_q)$ a Hermitian structure, where the compatible complex structure $(J^g)_q$ and the inner product $(g^D)_q$ are
\begin{equation}\label{eqn_alm_str_metric}
(J^g)_q=\begin{pmatrix}
0 & g_q^\ast\\
-g_q & 0
\end{pmatrix}\hbox{ and }
(g^D)_q=\begin{pmatrix}
g_q & 0\\
0 & g_q^\ast
\end{pmatrix}.
\end{equation}
Here $(\omega_{\times})_q=(g^D)_q((J^g)_q\cdot, \cdot)$.
 
We would like to describe $I_x$ with respect to the bases above. For this, we choose $\nabla$ as the dual connection to the Levi-Civita connection on $(N, g)$. For $s:N\rightarrow T^\ast N$ we take a local section such that $s(\widetilde{q})=\sum_i p_i(x)dq^i|_{\widetilde{q}}$. Note that $s(q)=x$ and $(\nabla s)(q)=0$. Then we obtain linear maps 
\begin{align}\label{split}
\begin{split}
I_{x}^h: &T_q N \longrightarrow T_{x}T^\ast N \quad\quad\quad\quad\quad\quad I_{x}^v:T_q^\ast N \longrightarrow T_{x}T^\ast N\\
&\partial_{q^i}|_q \,\longmapsto\partial_{q^i}|_{x}\,\quad\quad\quad\raisebox{15pt}{ and }  \quad\quad \quad\,\,  dq^i|_q \,\longmapsto\partial_{p_i}|_{x}.
\end{split}
\end{align}
Hence, $I_x=I_{x}^h\times I_{x}^v:T_q N\times T_q^\ast N\rightarrow T_{x}T^\ast N$ is a linear symplectomorphism.

It follows that the map $I:\pi^\ast(T N\oplus T^\ast N)\longrightarrow TT^\ast N$ is an isomorphism of symplectic vector bundles. Moreover, the pullbacks $P^\ast(J^g)$ and $P^\ast(g^D)$ induce a Hermitian structure on $TT^\ast N$.
\end{rem}

\begin{rem} \label{rem_ham_liouv} Observe that $P$ ``splits'' $X$ and $R$ at each point $x=(p, q)\in S^\ast N$. Indeed, by the formulas (\ref{split}) we have
\begin{alignat*}{3}
P_{x}^h(X_{x})&=0&&\text{ and }\quad P_{x}^v(X_{x})&&=x,\\
P_{x}^h(R_{x})&=x^\sharp\quad &&\text{ and }\quad P_{x}^v(R_{x})&&=0.
\end{alignat*}
\end{rem}

\begin{defn} Let $V$ be a linear subspace of $\mathbb{R}^n$. Then the \textbf{conormal space of} $V$ is defined as
$$L^\ast V:=V\times (V^\bot)^\ast=\lbrace(a, b)\in \mathbb{R}^{2n}=\R^n\times(\R^n)^\ast\,|\,a\in V, b(v)=0 \text{ for all }v\in V \rbrace.$$
\end{defn}

\begin{rem} \label{rem_tangent} Observe that if we restrict $T_{x} L^\ast Q$ to horizontal and vertical subspaces of $T_{x}T^\ast N$ then we get
$$P^h_{x}:T_{x}^h L^\ast Q\mapsto T_{q} Q\text{ and }P^v_{x}:T_{x}^v L^\ast Q\mapsto L_{q}^\ast Q.$$
Indeed, since $P^h=d\pi$, the first mapping was already described in Remark \ref{rem_conormal}. Because $P^h_x: T_{x}^h L^\ast Q\rightarrow T_{q} Q$ is an isomorphism, we get that $\dim T_{x}^v L^\ast Q=n-\dim Q=\dim L_{q}^\ast Q$. Since $P_x^v:T_x^vT^\ast N\rightarrow T^\ast_q N$ is an isomorphism with an inverse $I^v_x$, it is enough to show that $I^v_x$ maps every element of $L_{q}^\ast Q$ to $T_{x}^v L^\ast Q$. But, by the definition of $I^v_x$, $I^v_x(L_{q}^\ast Q)$ are elements of $T_x L^\ast Q$ that are in $\ker d\pi$, which is $T_x^v L^\ast Q$.
\end{rem}

\begin{rem}\label{rem_orb} We would like to define the Maslov index of a Hamiltonian chord $\gamma_H$ on $L^\ast Q\subset T^\ast N$. First, we would like to trivialize the symplectic bundle $\gamma_H^\ast TT^\ast N$ appropriately.

Let $$\lbrace\Phi_t:(T_{\gamma_H(t)}T^\ast N, \omega_{\gamma_H(t)})\longrightarrow (\mathbb{R}^{2n}, \omega_0)\rbrace_{t\in[0, T]}$$
be a continuous family of linear symplectomorphisms that satisfies the following:
\begin{itemize}
\item[(a)]$\Phi_t:T_{\gamma_H(t)}^vT^\ast N\longmapsto L^\ast\lbrace 0\rbrace$ for all $t\in[0, T]$,
\item[(b)]$\Phi_t:T_{\gamma_H(t)} L^\ast Q\longmapsto L^\ast V_t$ for some vector subspaces $V_t\subset\R^n $, where $t\in\lbrace 0, T\rbrace$.
\end{itemize}
Note that $\Phi_td\phi^t_H(\gamma_H(0))\Phi_0^{-1}$ is a path in $Sp(2n, \R)$ and hence maps Lagrangian subspaces to Lagrangian subspaces.
\end{rem}

\begin{defn}\label{def_maslov_ind_ham_pure} Let $\Phi$ be a trivialization as in Remark \ref{rem_orb}. Then the \textbf{Maslov index of the Hamiltonian chord} $\gamma_H$ is defined as
\begin{align*}\mu_{\gamma_H}&=\mu\left(\Phi_td\phi^t_H(\gamma_H(0))\Phi_0^{-1}\Phi_0T_{\gamma_H(0)} L^\ast Q, \Phi_TT_{\gamma_H(T)} L^\ast Q \right)\\
&=\mu\left(\Phi_td\phi^t_H(\gamma_H(0))T_{\gamma_H(0)} L^\ast Q, \Phi_TT_{\gamma_H(T)} L^\ast Q \right).
\end{align*}
\end{defn}

\begin{lemma}[{\cite[Prop 3.2]{APS}}] $\mu_{\gamma_H}$ does not depend on the choice of the trivialization $\Phi$.
\end{lemma}

\begin{thm}[{\cite[Prop 4.1]{APS}}]\label{thm_Morse_Maslov} Let $\gamma_H$ and $\gamma_M$ be corresponding chords on $Q$ and $L^\ast Q$, respectively. Then $\gamma_H$ is nondegenerate if and only if $\gamma_M$ is nondegenerate. 

Moreover, if $\gamma_H$ and $\gamma_M$ are nondegenerate, then it holds
\begin{equation}\label{Morse_Maslov}
\hbox{Ind}_{\gamma_M}=\mu_{\gamma_H}+\dim Q-\frac{n}{2}.
\end{equation}
 
\end{thm}

\section{Binormal and Reeb chords}

We will start this section with some discussion of the Maslov index in $T^\ast N$. Then we conclude with the goal of this chapter, that is, a relation between the Morse index of a binormal chord and the Maslov index of the corresponding pure Reeb chord.

\begin{lemma} \label{lemma_lift} Let $\dim N=n>2$ and assume that every loop in $Q$ is contractible in $N$. Then every loop in $\mathcal{L^\ast Q}$ is contractible in $S^\ast N$.
\end{lemma}

\begin{proof}
Let $\gamma:S^1\rightarrow\mathcal{L}^\ast Q$. We put $\gamma_0:=\pi\circ\gamma:S^1\rightarrow N$. Since $\gamma_0$ is contractible, there exists a capping disc $\overline{\gamma}_0$ of $\gamma_0$. 

Then $\gamma$ is a continuous section of $\gamma_0^\ast S^\ast N$. Since $\overline{\gamma}_0^\ast S^\ast N\cong D\times S^{n-1}$ for $S^{n-1}$ simply connected, there is an extension of $\gamma$ to a continuous section $\overline{\gamma}:S^1\rightarrow S^\ast N$.
\end{proof}

\noindent\textit{From now on, $N$ and $Q$ will satisfy the assumptions of Lemma \ref{lemma_lift}.}

\begin{lemma} \label{lemma_chern_split} The first Chern number satisfies the following:
\begin{itemize}
\item[(i.)] Let $\pmb{\bm{\gamma}}:S^2\rightarrow T^\ast N$ be a continuous map, then $c_1(\pmb{\bm{\gamma}}^\ast TT^\ast N)=0$.
\item[(ii.)] Let $\pmb{\bm{\gamma}}:S^2\rightarrow S^\ast N$ be a continuous map, then $c_1(\pmb{\bm{\gamma}}^\ast \xi)=c_1(\pmb{\bm{\gamma}}^\ast \eta)=0$. 
\end{itemize}
\end{lemma}

\begin{proof} Ad $(i.)$: Let $\pmb{\bm{\gamma}}_0:=\pi\circ \pmb{\bm{\gamma}}:S^2\rightarrow N$. We pick a Riemannian metric $g$ on $N$. Then by the fact that $I$ is an isomorphism of symplectic vector bundles and Remark \ref{rem_compos_pullback} about pullback maps, we have
$$c_1\big((\pmb{\bm{\gamma}}^\ast TT^\ast N, \omega)\big)=c_1\big((\pmb{\bm{\gamma}}_0^\ast (TN\oplus T^\ast N), \omega_\times, J^g)\big).$$

Now $g$ induce $O(n)$-structures on $\pmb{\bm{\gamma}}_0^\ast TN$ and $\pmb{\bm{\gamma}}_0^\ast T^\ast N$. Pick $S^+\subset S^2$. Then we take an isometric trivialization 
$$\varphi^h_+:\pmb{\bm{\gamma}}_0^\ast TN|_{S^+}\rightarrow S^+\times \R^n,$$
and find the isometric trivialization
$$\varphi^v_+:\pmb{\bm{\gamma}}_0^\ast T^\ast N|_{S^+}\rightarrow S^+\times (\R^n)^\ast,$$
which is defined pointwise as $\varphi^v_+:=\big((\varphi^h_+)^\ast\big)^{-1}$. We also define an automorphism $\widetilde{\sigma_+}$ on $S^+\times (\R^n\times(\R^n)^\ast)$ as $(Id, \sigma_{st})$, where $\sigma_{st}$ is the anti-symplectic involution defined in Example \ref{example_dual}. Then $\Phi^+:=\widetilde{\sigma_+}\circ(\varphi^h_+, \varphi^v_+)$ is an unitary trivialization of $(\pmb{\bm{\gamma}}_0^\ast (TN\oplus T^\ast N), \omega_\times, J^g)$. Now we construct analogously $\Phi^-$.
 
But then $\det(\Phi^+_t(\Phi^-_t)^{-1})$ is a real number for each $t\in S^1$. Hence, the corresponding loop has a zero degree.

Ad $(ii.)$: Observe that $\eta$ has a global frame spanned by nonvanishing Liouville and Reeb vector fields. Hence, also the pullback bundle $\pmb{\bm{\gamma}}^\ast \eta$ has a global frame. Thus $\pmb{\bm{\gamma}}^\ast \eta$ is trivial symplectic vector bundle and $c_1(\pmb{\bm{\gamma}}^\ast \eta)=0.$

Since $S^\ast N$ is embedded into $T^\ast N$, then, by the above, also $c_1(\pmb{\bm{\gamma}}^\ast TT^\ast N)=0$. Next, $\pmb{\bm{\gamma}}^\ast TT^\ast N$ can be written as a Whitney sum of $\pmb{\bm{\gamma}}^\ast \xi$ and $\pmb{\bm{\gamma}}^\ast \eta$. By the property of $c_1$ on the Whitney sum, the lemma follows.
\end{proof}

\begin{lemma} \label{lemma_maslov_class} Let $\gamma$ be a continuous loop $\gamma:S^1\rightarrow \mathcal{L}^\ast Q\subset S^\ast N$ with a capping disc $\overline{\gamma}$, then $\mu(\gamma^\ast T\mathcal{L}^\ast Q)=0$.
\end{lemma}

\begin{proof}
Recall that $\mathcal{L}^\ast Q$ is embedded in $L^\ast Q$.  First, we would like to show that $\mu(\gamma^\ast TL^\ast Q)=0$.

Let $\overline{\gamma}_0:=\pi\circ \overline{\gamma}:S^1\rightarrow N$. Pick a Riemannian metric $g$ on $N$. Then as in Lemma \ref{lemma_chern_split} there is an isomorphism mapping $(\overline{\gamma}^\ast TT^\ast N, \omega)$ to $(\overline{\gamma}_0^\ast (TN\oplus T^\ast N), \omega_\times)$. For simplicity, we will still denote it by $P$. Then $P$, when restricted to $S^1=\partial D$, also identify Lagrangian subbundles $\gamma^\ast TLQ$ and $\gamma_0^\ast(TQ\oplus L^\ast Q)$, see Remark \ref{rem_tangent}.

Then as in Lemma \ref{lemma_chern_split}, isometric trivializations 
$$\varphi^h:\overline{\gamma}_0^\ast TN\rightarrow S^1\times\R^n\hbox{ and }\varphi^v:\overline{\gamma}_0^\ast T^\ast N\rightarrow S^1\times(\R^n)^\ast$$
define together with anti-symplectic involution $\sigma_{st}$ a symplectic trivialization $\Phi:=\widetilde{\sigma}\circ(\varphi^h, \varphi^v)$ of $\overline{\gamma}_0^\ast (TN\oplus T^\ast N)$. 

Now we pick $\Sigma:=\overline{\Sigma_1(\R^n\times\lbrace 0\rbrace)}$. But, if $k:=\dim Q$, then the loop $\Phi\gamma_0^\ast(TQ\oplus L^\ast Q)$ lies completely in the $k$-th stratum of $\Sigma$. Hence, by the zero property of the Maslov index $\mu(\gamma^\ast TL^\ast Q)=\mu(\Phi P\gamma^\ast TL^\ast Q)=0$.

Now, we would like to show that $\mu(\gamma^\ast \langle R\rangle)=0$. As in the proof of the Lemma \ref{lemma_chern_split}, nonvanishing Liouville and Reeb vector fields imply the existence of some symplectic trivialization $\Psi_1$ of $\gamma^\ast\eta$ that maps $\gamma^\ast \langle R\rangle$ to the constant Lagrangian subspace in $\Lambda(1)$. Hence $\mu(\gamma^\ast \langle R\rangle)=\mu(\Psi_1(\gamma^\ast \langle R\rangle))=0$.

Next, since $\gamma^\ast\xi$ is over a contractible disc, there exists a symplectic trivialization $\Psi_2$ of $\gamma^\ast\xi$. Then $\Psi_1\times\Psi_2$ defines a symplectic trivialization of $\gamma^\ast\eta\oplus\gamma^\ast\xi=\gamma^\ast(\eta\oplus\xi)=\gamma^\ast TT^\ast N$. Now, we compute
\begin{align*}
0&=\mu(\gamma^\ast TL^\ast Q)\\
&=\mu\big((\Psi_1\times\Psi_2)(\gamma^\ast TL^\ast Q)\big)\\
&=\mu\big((\Psi_1\times\Psi_2)(\gamma^\ast\langle R\rangle\oplus\gamma^\ast T\mathcal{L}^\ast Q)\big)\\
&=\mu\big(\Psi_1\gamma^\ast\langle R\rangle\big)+\mu(\Psi_2\gamma^\ast T\mathcal{L}^\ast Q)\\
&=\mu(\Psi_2\gamma^\ast T\mathcal{L}^\ast Q)\\
&=\mu(\gamma^\ast T\mathcal{L}^\ast Q).
\end{align*}
Which finishes the proof.
\end{proof}

\begin{cor} \label{lemma_degree_reeb} If $\gamma_{R}$ is a pure Reeb chord on $\mathcal{L}^\ast Q\subset S^\ast N$, then $|\gamma_{R}|$ exists and is well defined. 
\end{cor}
\begin{proof} 
First, the existence of $|\gamma_{R}|$. The existence of $c, J$ and $\overline{\gamma}$ immediately follows from our assumptions and Remark \ref{rem_contr}. Hence, we can apply the Definition \ref{defn_maslov_index_lagr} of the Maslov index $\mu(L^\gamma)$.

Next, the ``well--definiteness'' of $|\gamma_{R}|$. By Corollary \ref{cor_chern} and Lemma \ref{lemma_chern_split} we need to inspect only the dependence of $|\gamma_{R}|$ on $J$ and $c$. 

For the various complex structures, it is enough to note the following. If we fix any trivialization of $\gamma^\ast\xi|_{[2/3, 1]}$, then by Lemma \ref{lemma_rotation_homotop} any two complex structures from Remark \ref{rem_reeb_triv} induce two paths in $\Lambda(n)$ that are homotopic with fixed end points.

Now, we would like to treat the case of two different capping paths. Let us take two discs $\overline{\gamma_1}, \overline{\gamma_2}$ of the Reeb chord $\gamma_R$ that are bounding two different capping paths. If we parametrize their boundary in the standard way, we see that they agree on $[1/3, 1]$. Next, the gluing $\pmb{\bm{\gamma}}:=\overline{\gamma_1}\cup_{[1/3, 1]}\overline{\gamma_2}$ is a map with a contractible domain. Hence $\pmb{\bm{\gamma}}^\ast\xi$ is trivial and the trivialization of $\pmb{\bm{\gamma}}^\ast\xi$ restricts to trivializations of $\overline{\gamma_1}^\ast\xi$ and $\overline{\gamma_2}^\ast\xi$. With respect to these trivializations it follows that $\mu(L^{\gamma_1})-\mu(L^{\gamma_2})$ is equal to the Maslov index of a loop in $\mathcal{L}^\ast Q$, which vanishes by Lemma \ref{lemma_maslov_class}.
\end{proof}

\begin{rem} \label{rem_corresp_reeb_bin_chord} By Remark \ref{rem_contact_objects} and Lemma \ref{lemma_corresp_ham_bin} we know that for any Reeb chord $\gamma_R:[0, T]\rightarrow S^\ast N$ with endpoints on $\mathcal{L}^\ast Q$ there is the unique binormal chord $\gamma_M:[0, T]\rightarrow N$ with endpoints on $Q$ that is parametrized by the arc length, and vice versa.
\end{rem}

\begin{thm} \label{lemma_index_formula} 
Let $\gamma_R$ and $\gamma_M$ be corresponding chords on $\mathcal{L}^\ast Q$ and $Q$, respectively. Then $\gamma_R$ is nondegenerate if and only if $\gamma_M$ is nondegenerate. 

Moreover, if $\gamma_R$ and $\gamma_M$ are nondegenerate, then it holds
\begin{equation}\label{Morse_Reeb}
|\gamma_R|=\text{Ind}_{\gamma_M}-\dim Q+n-2.
\end{equation}
\end{thm}

\begin{proof}
First, the correspondence of nondegenerate chords follows from Theorem \ref{thm_Morse_Maslov} and Formula (\ref{eqn_lin_ham_flow_liouv}).

Next, by Theorem \ref{thm_Morse_Maslov} it remains to show that
\begin{equation}\label{e_ham_reeb}
\mu_{\gamma_H}=|\gamma_{R}|-\frac{n}{2}+2.
\end{equation}

Let $\overline{\gamma}:D\rightarrow S^\ast N$ and $L^{\gamma, 2}$ be a capping disc and a Lagrangian subbundle from the definition of $|\gamma_R|$. Then we put $\overline{\gamma}_0:=\pi\circ \overline{\gamma}$ and $L^{\gamma, 1}$ is a Lagrangian subbundle of $\gamma^\ast\eta$ defined as follows. Since $\gamma$ is a positively oriented loop, we can take a parametrization by $S^1=\R/\mathbb{Z}$ and for $t\in \R/\mathbb{Z}$ we put 
\begin{itemize}
\item[$(i.)$] for $t\in[0, 1/3]$ we put  $L^{\gamma, 1}_{[t]}:=\langle R_{\gamma(t)}\rangle$
\item[$(ii.)$] for $t\in(1/3, 2/3]$ we put $L^{\gamma, 1}_{[t]}:=d\phi_H^{(3Tt-T)}\langle R_{\gamma_R(0)}\rangle$
\item[$(iii.)$] for $t\in(2/3, 1)$ we put $L^{\gamma, 1}_{[t]}:=\rot\>\![d\phi_H^{T}\langle R_{\gamma_R(0)}\rangle, \langle R_{\gamma_R(T)}\rangle; J, +](3t-2)$, where $J$ is a complex structure in $\mathcal{J}(\eta_{\gamma_R(T)}, \eta|_{\omega_{\gamma_R(T)}})$.
\end{itemize}
Observe that $\mu(L^{\gamma, 1})$ is well-defined.

Now we would like to compute $\mu(L^\gamma)$, where $L^\gamma:=L^{\gamma, 1}\oplus L^{\gamma, 2}\subset \gamma^\ast TT^\ast N$. The procedure will be similar to the proofs of Lemmata \ref{lemma_maslov_class}, \ref{lemma_chern_split}. So $\mu(L^\gamma)=\mu(PL^\gamma)$ for $P$ mapping $(\overline{\gamma}^\ast TT^\ast N, \omega)$ to $(\overline{\gamma}_0^\ast (TN\oplus T^\ast N), \omega_\times)$ as before, but $\Phi=\widetilde{\sigma}\circ(\varphi^h, \varphi^v)$ will be constructed as follows.

We would like to use Remarks \ref{rem_complement}, \ref{rem_sum_pullback} and \ref{rem_triv_homotop}. First, note that $[-1/3, 1/3]$ is contractible in $S^1=\R/\mathbb{Z}$ and 
$$(\gamma_0|_{[-1/3, 1/3]})^\ast TQ\oplus (\gamma_0|_{[-1/3, 1/3]})^\ast TQ^\bot\cong (\gamma_0|_{[-1/3, 1/3]})^\ast TN.$$
Hence there exists 
$$\varphi^h:(\gamma_0|_{[-1/3, 1/3]})^\ast TN\rightarrow [-1/3, 1/3]\times \R^n$$
 and some linear subspace $V\subset\R^n$ such that $\varphi^h:T_{\gamma_0(t)}Q\mapsto V\subset\R^n$ for each $t\in[-1/3, 1/3]$. Moreover, there exists a trivialization of $\overline{\gamma}_0^\ast TN$ that extends $\varphi^h$. We will denote this trivialization by the same symbol $\varphi^h$. This determines $\varphi^v: \overline{\gamma}_0^\ast T^\ast N\rightarrow D\times(\R^n)^\ast$ such that $\varphi^v:L^\ast_{\gamma_0(t)} Q\mapsto (V^\bot)^\ast$ for each $t\in[-1/3, 1/3]$. Finally, $\widetilde{\sigma}=(Id, \sigma_{st})$ is the automorphism of $D\times\R^n\times(\R^n)^\ast$, where $\sigma_{st}$ is the anti-symplectic on $\R^n\times(\R^n)^\ast$.   

Let us denote $\widetilde{\Phi}:=\Phi\circ P$, then we compute
\begin{align*}
\begin{split}
\mu(L^\gamma)&=\mu(\widetilde{\Phi}L^\gamma)\\
&=\mu(\widetilde{\Phi}L^\gamma, \widetilde{\Phi}L^\gamma_{[0]})\\
&=\mu(\widetilde{\Phi}L^\gamma_{[0, 1/3]}, \widetilde{\Phi}L^\gamma_{[0]})+\mu(\widetilde{\Phi}L^\gamma_{[1/3, 2/3]}, \widetilde{\Phi}L^\gamma_{[0]})+\mu(\widetilde{\Phi}L^\gamma_{[2/3, 1]}, \widetilde{\Phi}L^\gamma_{[0]}).
\end{split}
\end{align*}
Here because $\widetilde{\Phi}L^\gamma_{[0, 1/3]}$ is a constant path in $\Lambda(n),$ we have 
\begin{equation*}\label{e1}
\mu(\widetilde{\Phi}L^\gamma_{[0, 1/3]}, \widetilde{\Phi}L^\gamma_{[0]})=0.
\end{equation*}
Next, 
\begin{equation*}
\mu(\widetilde{\Phi}L^\gamma_{[1/3, 2/3]}, \widetilde{\Phi}L^\gamma_{[0]})=\mu_{\gamma_H}
\end{equation*}
and by Lemma \ref{lemma_maslov_rot}
\begin{equation*}
\mu(\widetilde{\Phi}L^\gamma_{[2/3, 1]}, \widetilde{\Phi}L^\gamma_{[0]})=\frac{n}{2}.
\end{equation*}

Hence
\begin{equation}\label{eqn1}
\mu(L^\gamma)=\mu_{\gamma_H}+\frac{n}{2}.
\end{equation}

On the other hand, we have
\begin{align*}
\mu(L^\gamma)=\mu(L^{\gamma, 1}\oplus L^{\gamma, 2})\\
=\mu(L^{\gamma, 1})+\mu(L^{\gamma, 2}).
\end{align*}
But 
$$\mu(L^{\gamma, 2})=|\gamma_R|+1.$$
In order to compute $\mu(L^{\gamma, 1})$, we take a symplectic trivialization $\Psi$ of $\overline{\gamma}^\ast\eta$ that sends a nonvanishing frame spanned by $R$ and $X$ to fixed basis vectors. It follows that
\begin{align*}
\mu(L^{\gamma, 1})&=\mu(\Psi L^{\gamma, 1})\\
&=\mu(\Psi L^{\gamma, 1}, \Psi L^{\gamma, 1}_{[0]})\\
&=\mu(\Psi L^{\gamma, 1}_{[0, 1/3]}, \Psi L^{\gamma, 1}_{[0]})+\mu(\Psi L^{\gamma, 1}_{[1/3, 2/3]}, \Psi L^{\gamma, 1}_{[0]})+\mu(\Psi L^{\gamma, 1}_{[2/3, 1]}, \Psi L^{\gamma, 1}_{[0]}).\\
\end{align*}

Because $\Psi L^{\gamma, 1}_{[0, 1/3]}$ is constant path in $\Lambda(1)$, we have
$$\mu(\Psi L^{\gamma, 1}_{[0, 1/3]}, \Psi L^{\gamma, 1}_{[0]})=0.$$
Next, by Formula \ref{eqn_lin_ham_flow_liouv} the only intersection of $\Psi L^{\gamma, 1}_{[1/3, 2/3]}$ with $\Psi L^{\gamma, 1}_{[0]}$ appears at $t=1/3$. Moreover, this intersection contributes with $+1/2$. Hence
$$\mu(\Psi L^{\gamma, 1}_{[1/3, 2/3]}, \Psi L^{\gamma, 1}_{[0]})=\frac{1}{2}.$$
And by Lemma \ref{lemma_maslov_rot}
$$\mu(\Psi L^{\gamma, 1}_{[2/3, 1]}, \Psi L^{\gamma, 1}_{[0]})=\frac{1}{2}.$$

Thus
\begin{equation}\label{eqn2}
\mu(L^{\gamma})=|\gamma_R|+2.
\end{equation}
Hence from Formulas (\ref{eqn1}) and (\ref{eqn2}) we obtain Formula (\ref{e_ham_reeb}), which finishes the proof.
\end{proof}
\chapter{Legendrian contact homology}
\label{ch:file5}

\begin{rem}\cite{felix2001rational}Let $X$ be a smooth manifold. We define a \textbf{normalized singular chain complex} as graded free abelian groups $C_\ast(X)$ together with the differential $\partial^{sing}:C_\ast(X)\rightarrow C_{\ast-1}(X)$, where
$$C_i(X)=\frac{\mathbb{Z}\langle C^\infty(\Delta^i, X)\rangle}{\mathbb{Z}\langle \hbox{degenerate maps}\rangle},\hspace{1cm}i\in\N_0.$$
Here, $\Delta^i$ is the convex hull of the standard basis $e_0,\dots, e_i$ of $\R^{i+1}$ and degenerate maps are those maps that are independent of at least one coordinate of $\R^{i+1}$. And for generators $\alpha\in C_i(X)$ we put $\partial^c\alpha=\sum_k(-1)^k\alpha\circ\delta_{k},$ where $\delta_{k}$ is the inclusion of the $k$-the face $\Delta^{i-1}\rightarrow\Delta^i$.

Next, if $X_1, X_2$ are smooth manifolds, then there is an Eilenberg-Zilber map \begin{align*}
\times:C_i(X_1)\otimes C_j(X_2)&\longrightarrow C_{i+j}(X_1\times X_2),
\end{align*}
which assigns to $\Delta^i\otimes \Delta^j$ the canonical subdivision of $\Delta^{i+j}$. For details, see \cite{felix2001rational}.
Moreover, $\times$ satisfies associativity. 
\end{rem}

\begin{rem}\cite{felix2001rational} Let $X$ be a smooth manifold and $x_0, x_1\in X$. Then $\Omega_{x_0, x_1}X$ denote a \textbf{Moore path space}. That is a set
\begin{align*}
\Omega_{x_0, x_1}X:=\lbrace \gamma\in C^\infty([0, T], X)\,|\,&T\in\R_{\geq0}, \gamma(0)=x_0, \gamma(R)=x_1\hbox{ and}\\
&\hbox{and all derivatives at endpoints vanish}\rbrace
\end{align*}
equipped with the $C^\infty$-compact-open topology and the multiplication $\ast:\Omega_{x_0, x_1}X\times\Omega_{x_1, x_2}X\rightarrow\Omega_{x_0, x_2}X$ given by concatenation as
$$(\gamma_1, \gamma_2)\longmapsto\gamma_1\ast\gamma_2(t)\equiv
\begin{cases}
\gamma_1(t)\hspace{2cm} \hbox{if }t\in[0, T_1],\\
\gamma_2(t-T_1)\hspace{1.1cm} \hbox{if }t\in[T_1, T_1+T_2].
\end{cases}$$
Here the maps $\gamma_i$ have domains $[0, T_i]$. If moreover $x_0=x_1$, then $\Omega_{x_0, x_0}X$ is called a \textbf{Moore based loop space} and denoted $\Omega_{x_0}X$. Last, if $x_1$ is allowed to vary on $X$, then the resulting space is called a \textbf{free-end Moore path space} and denoted $\overline{\Omega}_{x_0}X$. 

Observe that $\ast$ is associative.

If $X$ is path-connected, then the map evaluating paths at endpoints gives us a fibration
\begin{equation}\label{sequence_loop}
\Omega_{x_0}X\longrightarrow \overline{\Omega}_{x_0}X\longrightarrow X
\end{equation}
and $\overline{\Omega}_{x_0}X$ is contractible.
\end{rem}

\begin{rem}\cite{felix2001rational} The composition of $\times$ and $\ast$ gives $\big(C(\Omega_{x_0}X), \partial^c\big)$ a structure of a (strictly) associative differential graded ring.
\end{rem}

\begin{rem}\label{rem_groupring_loop} Now we would like to describe the structure of $H(\Omega_{x_0}\mathcal{L}^\ast_+T_K)$. 

First recall that $\mathcal{L}^\ast_+T_K\cong T_K$. Since $T_K$ has the universal cover $\R^2$, by the lifting property $\mathcal{L}^\ast_+T_K$ is an Eilenberg-MacLane space $K(G, 1)$, \cite{hatcher_2019}. Here $G=\mathbb{Z}^2$.

From fibration (\ref{sequence_loop}) and LES of homotopy groups we obtain that $\pi_{0}(\Omega_{x_0}\mathcal{L}^\ast_+T_K)\cong\pi_1(\mathcal{L}^\ast_+T_K)$ is the only nonzero homotopy group of $\mathcal{L}^\ast_+T_K$.

By Hurewicz theorem, \cite{felix2001rational}, the only nonzero homology group of $\Omega_{x_0}\mathcal{L}^\ast_+T_K$ is concentrated at the degree $0$. Here
$$H_0(\Omega_{x_0}\mathcal{L}^\ast_+T_K)\cong\bigoplus_{\pi_{1}(\mathcal{L}^\ast_+T_K)}\mathbb{Z}$$
and the multiplicative structure at $H_0(\Omega_{x_0}\mathcal{L}^\ast_+T_K)$ is induced by the group structure of $\pi_{1}(\mathcal{L}^\ast_+T_K)$.

Hence
$$H_0(\Omega_{x_0}\mathcal{L}^\ast_+T_K)\cong\mathbb{Z}[\pi_{1}(\mathcal{L}^\ast_+T_K)].$$
\end{rem}

\begin{defn} Let $x_0$ be a fixed point in $\mathcal{L}^\ast_+T_K$ and $\ell\in\mathbb{N}_0$. A \textbf{Reeb string with $\ell$ chords} is a word 
$$\alpha_1\textbf{a}_1\alpha_2\textbf{a}_2 \dots \textbf{a}_\ell\alpha_{\ell+1},$$
where
\begin{itemize}
\item[$(i.)$] For $k\in\lbrace 1,\dots, \ell\rbrace$, $\textbf{a}_i$ are Reeb chords on $\mathcal{L}^\ast_+T_K$ parametrized by times $T_i$.
\item[$(ii.)$] $\alpha_1\in\Omega_{x_0, \textbf{a}_1(T_1)}\mathcal{L}^\ast_+T_K$.
\item[$(iii.)$] For $k\in\lbrace 2,\dots, \ell\rbrace$, $\alpha_k\in\Omega_{\textbf{a}_{k-1}(0), \textbf{a}_k(T_k)}\mathcal{L}^\ast_+T_K$.
\item[$(iv.)$] $\alpha_{\ell+1}\in\Omega_{\textbf{a}_\ell(T_0), x_0}\mathcal{L}^\ast_+T_K$.
\end{itemize}
See also Figure \ref{figure_reeb_string}.

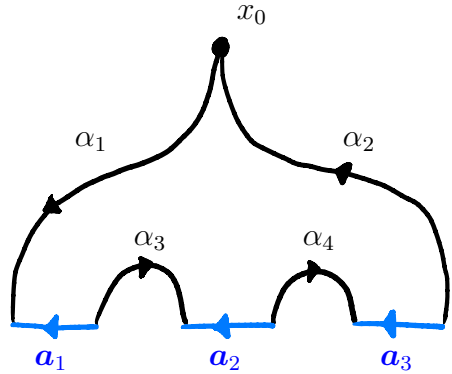
\begin{figure}[!htbp]
\labellist
\pinlabel $x_0$ at 300 420
\pinlabel $\textcolor{blue}{\bm{a}_1}$ at 225 290
\pinlabel $\textcolor{blue}{\bm{a}_2}$ at 290 290
\pinlabel $\textcolor{blue}{\bm{a}_3}$ at 354 290
\pinlabel $\alpha_1$ at 240 372
\pinlabel $\alpha_2$ at 340 372
\pinlabel $\alpha_3$ at 262 335
\pinlabel $\alpha_4$ at 325 335
\endlabellist
\centering
\includegraphics[scale=1]{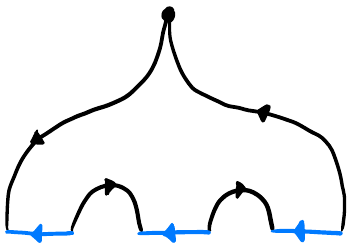}
\vspace{0.1cm}
\caption{The Reeb string 
$\alpha_1{\bm{a}_1}\alpha_2{\bm{a}_2}\alpha_3{\bm{a}_3}\alpha_4$. Observe that Reeb chords (in blue) are traversed in the opposite direction.}
\label{figure_reeb_string}
\end{figure}

Next, $\mathcal{R}^\ell$ denotes the set of Reeb strings with $\ell$ chords equipped with $C^m$-compact-open topology. We put $\mathcal{R}:=\coprod_{\ell\geq 0}\mathcal{R}^\ell$. Then the concatenation $\ast$ at $x_0$ gives $\mathcal{R}$ a structure of a $H$-space. 
\end{defn}

\begin{rem} The sub-$H$-space $\mathcal{R}^0$ agrees with $\Omega_{x_0}\mathcal{L}^\ast_+T_K$.
\end{rem}

\begin{rem} The composition of $\times$ and $\ast$ gives $C(\mathcal{R})$ the structure of a noncommutative but associative ring, which is graded by the degree of the singular chains. Next, we can extend the grading of $C(\mathcal{R})$ by adding the sum of degrees of the Reeb chords that correspond to given Reeb strings. Since the degree of a Reeb chord is independent of a capping path, $C(\mathcal{R})$ still has the structure of a strictly associative graded ring.
\end{rem}

\begin{defn} Let $\textbf{a}, \textbf{b}_1, \dots, \textbf{b}_m$ be Reeb chords on $\mathcal{L}^\ast_+T_K$ that are parametrized by times $T_0, T_1,\dots, T_m$, respectively. If $u\in\mathcal{M}^{sy}(\textbf{a}, \textbf{b}_1, \dots, \textbf{b}_m)/\R$, then $\partial(u)$ is a word
$$\partial(u):=\beta_1\textbf{b}_1\beta_2\textbf{b}_2\dots\textbf{b}_m\beta_{m+1},$$
where $\beta_k$ are boundary arcs of $u$ projected to $\mathcal{L}^\ast_+T_K$ and are in counterclockwise order and orientation of $\partial D_{m+1}$. Hence
\begin{itemize}
\item[$(i.)$] $\beta_1\in\Omega_{\textbf{a}(T_0), \textbf{b}_1(T_1)}\mathcal{L}^\ast_+T_K$.
\item[$(ii.)$] For $k\in\lbrace 2,\dots, m\rbrace$, $\beta\in\Omega_{\textbf{b}_{k-1}(0), \textbf{b}_k(T_k)}\mathcal{L}^\ast_+T_K$.
\item[$(iii.)$] $\beta_{m+1}\in\Omega_{\textbf{b}_m(0), \textbf{a}(0)}\mathcal{L}^\ast_+T_K$.
\end{itemize}
\end{defn}

\begin{defn} The differential $\partial^{sy}$ on $C(\mathcal{R})$ is defined as follows. Let $W\in C(\mathcal{R})$ be a chain of a type $W=\alpha_1\textbf{a}_1\dots \textbf{a}_\ell\alpha_{\ell+1}$, then
$$\partial^{sy}(W):=\sum_{k=1}^\ell\sum_{\substack{m\geq 0\\ \textbf{b}_1,\dots, \textbf{b}_m}}\sum_{\substack{\dim\mathcal{M}^{sy}(\textbf{a}, \textbf{b}_1, \dots, \textbf{b}_m)/\R=0 \\ u\in\mathcal{M}^{sy}(\textbf{a}, \textbf{b}_1, \dots, \textbf{b}_m)/\R}}\varepsilon(-1)^{d+|\textbf{a}_1|+\cdots+|\textbf{a}_{k-1}|}\partial(u)\cdot_k W.$$
Here $\varepsilon$ is a sign from orientation of $0$-dimensional manifold $\mathcal{M}^{sy}(\textbf{a}, \textbf{b}_1, \dots, \textbf{b}_m)/\R$ (see also Theorem \ref{thm_dim_symp}), $d$ is a singular chain degree and $\partial(u)\cdot_k W$ is a word
$$\partial(u)\cdot_k W:=\alpha_1\textbf{a}_1\dots\alpha_k\partial(u)\alpha_{k+1}\dots\textbf{a}_\ell\alpha_{\ell+1},$$
where paths $\alpha_k$ and $\alpha_{k+1}$ are concatenated with paths $\beta_1$ and $\beta_m$, respectively. See Figure \ref{figure_reeb_string_enhanced}.

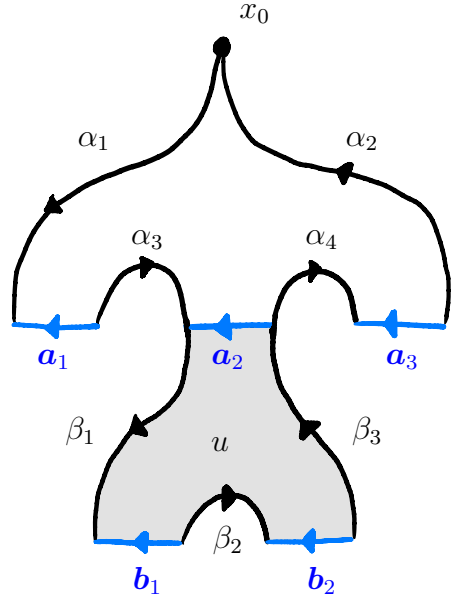
\begin{figure}[!htbp]
\labellist
\pinlabel $x_0$ at 275 425
\pinlabel $\textcolor{blue}{\bm{a}_1}$ at 200 295
\pinlabel $\textcolor{blue}{\bm{a}_2}$ at 265 295
\pinlabel $\textcolor{blue}{\bm{a}_3}$ at 330 295
\pinlabel $\alpha_1$ at 215 377
\pinlabel $\alpha_2$ at 315 377
\pinlabel $\alpha_3$ at 235 340
\pinlabel $\alpha_4$ at 300 340
\pinlabel $\beta_1$ at 210 268
\pinlabel $\beta_2$ at 265 227
\pinlabel $\beta_3$ at 317 268
\pinlabel $\textcolor{blue}{\bm{b}_1}$ at 235 212
\pinlabel $\textcolor{blue}{\bm{b}_2}$ at 300 212
\pinlabel $u$ at 262 263
\endlabellist
\centering
\includegraphics[scale=1]{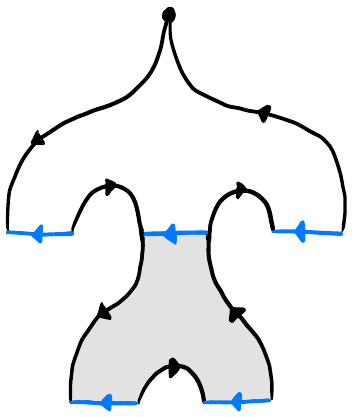}
\vspace{0.3cm}
\caption{The Reeb string $\partial(u)\cdot_2\alpha_1{\bm{a}_1}\alpha_2{\bm{a}_2}\alpha_3{\bm{a}_3}\alpha_4$, where $\partial(u)=\beta_1\bm{b}_1\beta_2\bm{b}_2\beta_3$.}
\label{figure_reeb_string_enhanced}
\end{figure}
\end{defn}

\begin{rem}By Stokes' Theorem $L(\bm{a})-\sum_{k=1}^m L(\bm{b}_k)\geq0$. Hence, the sum in the definition of $\partial^{sy}$ is finite.
\end{rem}

\begin{defn} $\partial_{\mathcal{L}}:=\partial^{sing}+\partial^{sy}:C(\mathcal{R})\rightarrow C(\mathcal{R})$.
\end{defn}

\begin{rem}$\partial_{\mathcal{L}}$ has degree $-1$.
\end{rem}

\begin{defn_thm}[\cite{Cieliebak2016KnotCH, EES05}] $\partial_{\mathcal{L}}^2=0$, the homology of $(C(\mathcal{R}), \partial_{\mathcal{L}})$ is called \textbf{Legendrian contact homology $LCH(\mathcal{L}^\ast_+T_K)$} and is an invariant under all choices. In particular, $LCH(\mathcal{L}^\ast_+T_K)$ is an invariant under Legendrian isotopy of $\mathcal{L}^\ast_+T_K$ in a fixed chosen spin structure.
\end{defn_thm}

\begin{defn} Let $(R, \partial_R)$ be a differential graded (dg) ring. An \textbf{$(R, \partial_R)$-NC-algebra} is a dg ring $(S, \partial_S)$ together with a dg ring homomorphism $(R, \partial_R)\rightarrow(S, \partial_S)$. 

Two $(R, \partial_R)$-NC-algebras $(S_1, \partial_{S_1}), (S_2, \partial_{S_2})$ are \textbf{isomorphic} if there is a dg ring isomorphism $(S_1, \partial_{S_1})\rightarrow(S_2, \partial_{S_2})$ that intertwines with maps $(R, \partial_R)\rightarrow(S_1, \partial_{S_1})$ and $(R, \partial_R)\rightarrow(S_2, \partial_{S_2})$. 

Here ``NC'' means ``noncommutative''.
\end{defn}

\begin{example}By Remark \ref{rem_groupring_loop}, $LCH(\mathcal{L}^\ast_+T_K)$ is an NC-algebra over the group ring $\mathbb{Z}[\pi_1(\mathcal{L}^\ast_+T_K)].$
\end{example}
\chapter{Proof of Lemma \ref{lemma_aux_diag}}
\label{apx:proof_lemma}

First, we describe $K$ in suitable coordinates around $\gamma(\overline{s})$. For this, we can assume that $\kappa(\overline{s})>0$. We take the following isothermal coordinates $(x, y, z)\in(\R^3, (1/2)\kappa(\overline{s})g_{Euc})$. We would like to find $\delta_{\overline{s}}$ and $\varepsilon_{0, \overline{s}}$ in these coordinates. Note that the actual $\delta_{\overline{s}}$ and $\varepsilon_{0, \overline{s}}$ need to be scaled by $\kappa(\overline{s})/2$ and $2/\kappa(\overline{s})$, respectively.

The coordinates are characterized by:
\begin{itemize}
\item $\gamma(\overline{s})=(0, 0, 0)$.
\item Around $\gamma(\overline{s})$ the curve $\gamma$ is reparametrized by $t\in(-\overline{\delta}, \overline{\delta})$ and given as
\begin{align*}
x&=t,\\
y&=t^2+Y_1 t^3+Y_2 t^4+O(t^5),\\
z&=\quad\quad Z_1 t^3+Z_2 t^4+O(t^5),
\end{align*}
for some real constants $Y_1, Y_2, Z_1, Z_2$.
\end{itemize}
Here $\overline{\delta}>0$ was chosen sufficiently small, such that the above chart exists. Note that in these coordinates, $\gamma$ is not necessarily parametrized by arc length, but for purposes of this lemma, it is not an issue.

Let $\varepsilon>0$. Then $\Gamma_\epsilon:(-\overline{\delta}, \overline{\delta})\times \R/2\pi\mathbb{Z}\rightarrow\R^3$ will still describe a parametrization of the torus $T_{K, \varepsilon}$, but now in our local coordinates. By symmetry, we can restrict ourself to $[0, \overline{\delta})$.

First, we treat the cases $(i.)$ and $(ii.)$.

We introduce a smooth family of functions $\lbrace g_{t, \varepsilon}\rbrace_{t\in[0,\overline{\delta})}$, which are defined as
\begin{align*}
g_{t, \varepsilon}:\R/2\pi\mathbb{Z}&\rightarrow\hspace{2cm}\R\\
\theta\hspace{0.45cm}&\mapsto \pi_{y}(\Gamma_{\varepsilon}(t, \theta))^2+\pi_{z}(\Gamma_{\varepsilon}(t, \theta))^2
\end{align*}
Here $\pi_{y}, \pi_{z}$ are the canonical projections on $\R^2_{y, z}$. The minimum of $g_{t, \varepsilon}$ describes the square of the distance of the circle $\Gamma_{\varepsilon}|_t$ from the $x$-axis. Recall that if the minimum is $\varepsilon^2$, then $\Gamma_\varepsilon|t$ is tangent to $\mathcal{C}_{(0, 0, 0), \varepsilon}$

Now, we inspect the uniqueness of the minimum of $g_{t, \varepsilon}$. Clearly, for $t=0$, there is a whole circle of minima. We can translate the problem into finding the minima of the distance between an ellipse $E$ and a point $p$ in the plane. Observe that if $p$ is lying outside of the $E$ then there is the unique point on $E$ achieving the minima. It corresponds to the touching point of $E$ and some circle with the center $p$. And such a point is unique by the convexity of the circle and the ellipse. If $p$ lies inside $E$, there could potentially be two minima, but we are not interested in this particular case. To achieve the translation of the problem we just project the circle $\Gamma_{\varepsilon}|_t$ to $\pi_{y, z}(\Gamma_{\varepsilon}|_s)\subset\R^2_{y, z}$ and $x$-axis to $(0, 0)$. By shrinking $\overline{\delta}$, if necessary, we can assume that $\pi_{y, z}$ is an embedding for any $t\in(0, \overline{\delta})$.

First, let us find when $\pi_{y, z}(\Gamma_{\varepsilon}|_s)$ lies outside of the point $(0, 0)$. It will be sufficient to solve the inequality
\begin{equation}\label{ineq_dist1}
\pi_{y}(\gamma(t))>\varepsilon.
\end{equation}
For this we split inequality (\ref{ineq_dist1}) into two inequalities, where one will contain $\varepsilon$-terms and the second $O(t^3)$-terms. That is
\begin{alignat}{1}
 (3/4)t^2>\varepsilon \label{ineq_dist2},\\
 (1/4)t^2+O(t^2)>0\label{ineq_dist3}.
\end{alignat}

First, inequality (\ref{ineq_dist2}). There exists $\varepsilon_1>0$ such that for every $\varepsilon\in(0,\varepsilon_1)$ it holds that for every $t\in(\sqrt{4/3}\,\varepsilon^{1/2}, \overline{\delta})$ inequality (\ref{ineq_dist2}) is satisfied.

Next, inequality (\ref{ineq_dist3}). Here we obtain immediately that there exists $\delta_1\in(0, \overline{\delta})$ such that inequality (\ref{ineq_dist3}) holds for $t\in(0, \delta_1)$.

Altogether, inequality (\ref{ineq_dist1}) holds for $t\in(\sqrt{4/3}\,\varepsilon^{1/2}, \delta_1)$. Note that for $\delta_1$ we can find $\varepsilon_1>0$ such that for every $\varepsilon\in(0, \varepsilon_1]$ the intervals $(\sqrt{4/3}\,\varepsilon^{1/2}, \delta_1)$ are nonempty. And specially, the corresponding $g_{t, \varepsilon}$ has the unique minimum.

Now, we are going to explicitly find the minimum of $g_{t, \varepsilon}$.

\begin{figure}[!htbp]
\labellist
\pinlabel ${\hbox{distance}^2}$ at -9 140
\pinlabel $\varepsilon^2$ at 11 82
\pinlabel $t$ at 250 0
\pinlabel $g_{t, \varepsilon}(\theta(t))$ at 250 120
\pinlabel $0$ at 13 0
\endlabellist
\centering
\includegraphics[scale=0.98]{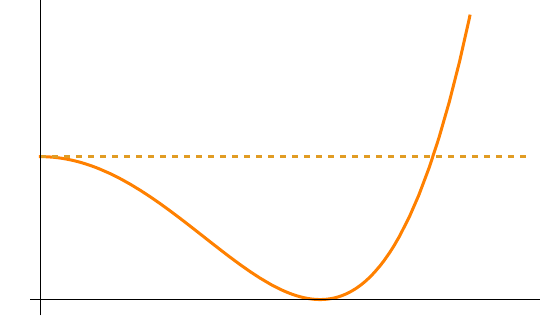}
\vspace{0.3cm}
\caption{Graph of the function $g_{t, \varepsilon}(\theta(t))$.}
\label{figure_distance}
\end{figure}
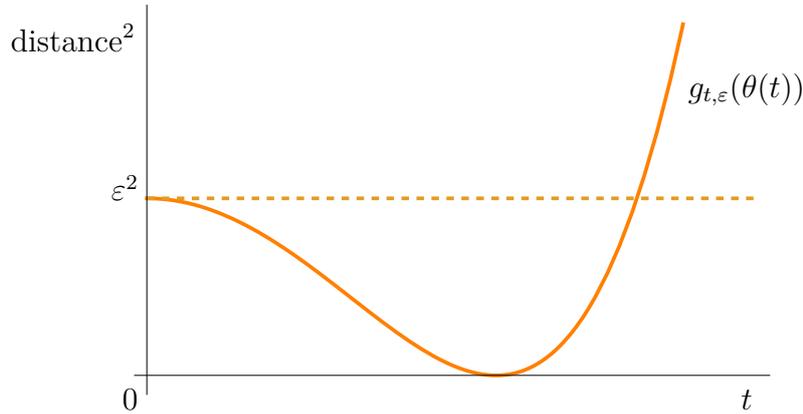

Let us assume that around $(0, 0, 0)$ $\gamma$ contains no torsion (i.e. $z=0$). Then we can quickly find the minimum of $g_{t, \varepsilon}$. First note that $\gamma$ lies in the $xy$-plane, and the $xy$-plane is orthogonal to the circle $\Gamma_{\varepsilon}|_t$ and to the $x$-axis. Hence by Lagrangian multipliers, the minimum is attended at $\theta(t)=\pi$ (or in another words, at $\gamma(t)-\varepsilon n(t)$). In particular, if $\gamma$ does not contain any higher order terms of $t$ than of the second order, then $g_{t, \varepsilon}(\theta(t))$ is given as on Figure \ref{figure_distance} and explicitly as
\begin{align*}
g_{t, \varepsilon}(\theta(t))&=\left(t^2-\frac{\varepsilon}{\sqrt{1+4t^2}}\right)^2\\
&=\varepsilon^2-2(\varepsilon+2\varepsilon^2)t^2+(1+4\varepsilon+16\varepsilon^2)t^4+O(t^5).
\end{align*}

When $\gamma$ also includes torsion, the computation of the desired $\theta(t)$ is more complicated. One needs to solve $\partial_\theta g_{t, \varepsilon}(\theta)=0$, which is a quartic equation in $\ell$ after the substitution $\cot(\theta/2)=\ell$. Our aim is to compute the Taylor expansion of $g_{t, \varepsilon}(\theta(t))$ up to order $4$. This can be done by cutting $g_{t, \varepsilon}(\theta(t))$ into smaller pieces and computing the sufficiently high Taylor expansions of these pieces first. Using \textit{Mathematica} we obtain that
\begin{align*}
&g_{t, \varepsilon}(\theta(t))=\\
&=\varepsilon^2-2(\varepsilon+2\varepsilon^2)t^2-2Y_1\varepsilon(1+6\varepsilon)t^3\\
&+\frac{1+(8-Z_1^2+2Y_2)\varepsilon+(32-9Y^2_1+24Y_2-12Z_1^2)\varepsilon^2+(64-36Y^2_1+64Y_2-36Z_1^2)\varepsilon^4}{1+4\varepsilon}t^4\\
&+O(t^5).
\end{align*}
Which can be rewritten as
$$g_{t, \varepsilon}(\theta(t))=\varepsilon^2-2(\varepsilon+2\varepsilon^2)t^2-\varepsilon P_1t^3+(1+\varepsilon P_2)t^4+O(t^5),$$
where $P_1$ and $P_2$ are some expressions that depend on $\lbrace\varepsilon, Y_1, Y_2, Z_1\rbrace$ and have the property that $\varepsilon^{1/2}P_1$ and $\varepsilon^{1/2}P_2$ are $O(\varepsilon^{1/2})$.

Now we are going to inspect analytical properties of $g_{t, \varepsilon}(\theta(t))$ relative to $\varepsilon^2$. In more detail, will describe the monotonicity of $g_{t, \varepsilon}(\theta(t))$ and upper/lower bounds to $g_{t, \varepsilon}(\theta(t))$, which will help us to localize the solutions of $g_{t, \varepsilon}(\theta(t))=\varepsilon^2$.

First, we inspect when $g_{t, \varepsilon}(\theta(t))$ is growing. Hence, we study the inequality
\begin{equation}\label{ineq_dist4}
\partial_t\big(g_{t, \varepsilon}(\theta(t))\big)>0.
\end{equation}
Writing $1\cdot t^3=(1/4)t^3+(3/4)t^3$, we, analogously to inequality (\ref{ineq_dist1}), split inequality (\ref{ineq_dist4}) into two inequalities, where the first one is
\begin{equation*}
4(1/4)t^3+O(t^4)>0.
\end{equation*}
Which holds for $t\in(0, \delta_2)$, where $\delta_2>0$ is sufficiently small.

Now we take the second inequality of the splinting of inequality (\ref{ineq_dist4}). Here we substitute $t:=\varepsilon^{1/2}t_1$ and obtain
$$-4\varepsilon^{3/2}t_1-8\varepsilon^{7/2}t_1-\varepsilon^{2}3P_1t_1^2+4\varepsilon^{3/2}t_1^3+\varepsilon^{5/2}4P_2t_1^3>0.$$
Then we divide the inequality by $\varepsilon^{3/2}$, put $\varepsilon:=0$ and obtain
$$-4t_1+4(3/4)t_1^3>0.$$

Hence, altogether, there exists $\delta_2, \varepsilon_2>0$ such that for every $\varepsilon\in(0, \varepsilon_2]$ the function $g_{t, \varepsilon}(\theta(t)$ is growing along $(\sqrt{4/3}\,\varepsilon^{1/2}, \delta_2)$.

Now, we inspect when $g_{t, \varepsilon}(\theta(t))$ is decreasing, that is the inequality
\begin{equation}\label{ineq_dist8}
\partial_t\big(g_{t, \varepsilon}(\theta(t))\big)<0.
\end{equation}
Again we put $t:=\varepsilon^{1/2}t_1$, so $O(t^4)$ becomes $\varepsilon^2 O(t_1^4)$. Then divide the inequality by $\varepsilon^{3/2}>0$ and put $\varepsilon:=0.$ So we are left with
$$-4t_1+4t_1^3<0.$$
Hence there is $\varepsilon_3>0$ such that for every $\varepsilon\in(0, \varepsilon_3]$ it holds that for every $t\in(0, \varepsilon^{1/2})$ inequality (\ref{ineq_dist8}) is satisfied.

Next, we are going to inspect when $g_{t, \varepsilon}(\theta(t))$ is smaller/bigger than $\varepsilon^2$.

First, the inequality
\begin{equation}\label{ineq_dist5}
g_{t, \varepsilon}(\theta(t))<(\varepsilon/2)^2.
\end{equation}
We substitute $t:=\varepsilon^{1/2}t_1$. So $O(t^{5})$ becomes $\varepsilon^{5/2}O(t_1^5)$. Then we consequently divide inequality (\ref{ineq_dist5}) by $\varepsilon^2$, put $\varepsilon:=0$ and obtain
$$1-2t^2+t^4<1/4.$$
Hence, there exists $\varepsilon_4>0$ such that for every $\varepsilon\in[\varepsilon_4, 0)$ it holds that for every $t\in(\sqrt{1/2}\,\varepsilon^{1/2}, \sqrt{3/2}\,\varepsilon^{1/2})$ inequality (\ref{ineq_dist5}) is satisfied.

Now the inequality
\begin{equation}\label{ineq_dist6}
g_{t, \varepsilon}(\theta(t))>(2\varepsilon)^2.
\end{equation}
We write $1\cdot t^4=(1/2)t^4+(1/2)t^4$. As before, we split inequality (\ref{ineq_dist6}) into two inequalities. The first one is
\begin{equation*}
(1/2)t^4+O(t^5)>0.
\end{equation*}
This holds for $t\in(0, \delta_3)$, where $\delta_3>0$ is sufficiently small.

Now we take the second inequality of the splinting of inequality (\ref{ineq_dist6}). Here we consequently substitute $t:=\varepsilon^{1/2}t_1$, divide the inequality by $\varepsilon^{3/2}$, put $\varepsilon:=0$ and obtain
$$1-2t^2+(1/2)t^4>4.$$
Hence, altogether, there exists $\delta_3, \varepsilon_5>0$ such that for every $\varepsilon\in(0, \varepsilon_5]$ it holds that for every $t\in((23/10)\varepsilon^{1/2}, \delta_3)$ inequality (\ref{ineq_dist6}) is satisfied.

Next, we put $\delta_{\overline{s}}:=\hbox{min}\lbrace \delta_1, \delta_2, \delta_3\rbrace$. Then for every $\delta_0\in(0, \delta_{\overline{s}}]$ we put $\varepsilon_{0, \overline{s}}:=\hbox{min}(\varepsilon_1, \varepsilon_2, \varepsilon_3, \varepsilon_4, \varepsilon_5, \delta_0, \varepsilon_{good})$. This finishes the poofs of $(i.)$.

We will solve the cases $(ii.)$ analogously.

Since $\gamma$ is transverse to the $yz$-plane, the intersection $E_{t, \varepsilon}:=C_{t, \varepsilon}\cap\lbrace x=0\rbrace$ is an ellipse. The parametrization of $\Gamma_{\varepsilon}|_t$ by $\theta$ naturally induces the parametrization of $E_{\varepsilon, t}$, which will also be denoted by $\theta$.

Let $\lbrace h_{t, \varepsilon}\rbrace_{t\in[0, \delta_0)}$ be a smooth family of functions, where
\begin{align*}
h_{t, \varepsilon}:\R/2\pi\mathbb{Z}&\rightarrow\hspace{2cm}\R\\
\theta\hspace{0.45cm}&\mapsto \pi_{y}(E_{t, \varepsilon}(\theta))^2+\pi_{z}(E_{t, \varepsilon}(\theta))^2.
\end{align*}

Similarly as before, we are interested in finding the minimum of $h_{t, \varepsilon}$.

First we inspect the uniqueness. Let $c:[0, \delta_0)\rightarrow\R^3$ be a curve assigning to each $t$ the intersection of the line $\lbrace\gamma(t)+\langle\dot{\gamma}(t)\rangle\rbrace$ with the plane $\lbrace x=0\rbrace$. Then we are interested in solving the inequality
\begin{equation}\label{ineq_dist7}
\pi_y(c(t))<-\varepsilon.
\end{equation}
Which can be rewritten as
$$t^2+O(t^3)>\varepsilon.$$
Hence, similarly to inequality (\ref{ineq_dist1}), there exist $\delta_4>0$ and $\varepsilon_6>0$ such that for every $\varepsilon\in(0, \varepsilon_6)$ and every $t\in(\varepsilon^{1/2}, \delta_4)$ inequality (\ref{ineq_dist7}) is satisfied. And specially such a $h_{t, \varepsilon}$ has the unique minimum.

Also, if $\gamma$ does not contain any higher order terms of $t$ than of the second order, then by Lagrangian multipliers we compute $\theta(t)=0$ and
\begin{align*}
h_{t, \varepsilon}(\theta(t))&=\left(t^2-\frac{\varepsilon(1+4t^2)}{\sqrt{1+4t^2}}\right)^2\\
&=\varepsilon^2-2(\varepsilon-2\varepsilon^2)t^2+(1-4\varepsilon)t^4+O(t^5).
\end{align*}

In the general case, by \textit{Mathematica}, we obtain that
$$h_{t, \varepsilon}(\theta(t))=\varepsilon^2-2(\varepsilon-2\varepsilon^2)s^2-\varepsilon P_3s^3+(1+\varepsilon P_4)s^4+O(s^5),$$
where $P_3$ and $P_4$ are some expressions that depend on $\lbrace\varepsilon, Y_1, Y_2, Z_1\rbrace$ and have the property that $\varepsilon^{1/2}P_3$ and $\varepsilon^{1/2}P_4$ are $O(\varepsilon^{1/2})$.

Analogously, as above, we obtain the following analytical properties. First, there exist $\delta_5>0$ and $\varepsilon_7>0$ such that for every $\varepsilon\in[\varepsilon_7, 0)$ the function $h_{t, \varepsilon}(\theta(t))$ is growing on $(\sqrt{4/3}\,\varepsilon^{1/2}, \delta_5)$ and decreasing on $(0, \varepsilon^{1/2})$. Also there exists $\varepsilon_8>0$ such that for every $\varepsilon\in[\varepsilon_8, 0)$ it holds that for every $t\in(\sqrt{1/2}\,\varepsilon^{1/2}, \sqrt{3/2}\,\varepsilon^{1/2})$ the inequality $h_{t, \varepsilon}(\theta(t))<(\varepsilon/2)^2$ is satisfied. And there exist $\delta_6>0$ and $\varepsilon_9>0$ such that for every $\varepsilon\in(0, \varepsilon_9]$ it holds that for every $t\in((23/10)\,\varepsilon^{1/2}, \delta_6)$ the inequality $h_{t, \varepsilon}(\theta(t))>(2\varepsilon)^2$ is satisfied.

Now, we put $\delta_{\overline{s}}:=\hbox{min}\lbrace \delta_1,\dots, \delta_6\rbrace$. Then for $\delta_0\in(0, \delta_{\overline{s}}]$ we put $\varepsilon_{0, \overline{s}}:=\hbox{min}(\varepsilon_1,\dots, \varepsilon_9, \delta_0, \varepsilon_{good})$. This finishes the proof of $(i.)-(ii.)$.

\chapter*{List of Symbols}
\label{ch:symbols}
\addcontentsline{toc}{chapter}{List of Symbols}

\begin{tabular}{cp{0.74\textwidth}}
$T^\ast N, S^\ast N$ & cotangent and unit cotangent bundles over $N$, respect.\\[0.5ex]
$L^\ast Q, \mathcal{L}^\ast Q$ & conormal and unit conormal bundles over $Q$, respect.\\[0.8ex]
$L^\ast_+ Q, \mathcal{L}^\ast_+ Q$ & one half of (unit) conormal bundles over the $\codim 1$ submanifold $Q$, see Remark \ref{rem_not_symplectic}\\[0.8ex]
$T_{K, \varepsilon}$ & $\varepsilon$-radius torus around the knot $K$\\[0.8ex]
$M_{K, \varepsilon}$ & space of outward-pointing chords on $T_{K, \varepsilon}$\\[0.8ex]
$\Delta_0, \Delta_{full}, \Delta_\varepsilon$ & various diagonals of $K\times K, T_{K, \varepsilon}\times T_{K, \varepsilon}, M_{K, \varepsilon}$\\[0.8ex]
$X_f$ & gradient like vector field adapted to $f$\\[0.8ex]
$E_0, E_\varepsilon$ & restrictions of the standard energy function $E:\R^3\times \R^3\rightarrow\R$ to the (parametrizations) of $K\times K$ and $T_{K, \varepsilon}\times T_{K, \varepsilon}$\\[0.8ex]
$\clubsuit_{X_{E_\varepsilon}}^{c, out-out}$ & cascade Morse flow tress of $X_{E_\varepsilon}$ that are restricted to $M_{K, \varepsilon}$\\[0.8ex]
$\mathfrak{m}, \mathfrak{L}$ & meridian and longitude\\[0.8ex]
$\mu, \lambda$ & generators of the cord algebra corresponding to $\mathfrak{m}$ and  $\mathfrak{L}$, respect.\\[0.8ex]
$\nu Q$ & tubular neighborhood of $Q$\\[0.8ex]
$\gamma(s), \Gamma_\varepsilon(s, \theta)$ & parametrizations of $K$ and $T_{K, \varepsilon}$\\[0.8ex]
$v(s, \theta), v^{\bot}(s, \theta)$ & normal vector to $T_{K, \varepsilon}$ at $\Gamma_\varepsilon(s, \theta)$ and its orthogonal complement $\partial_\theta v= v^\bot$\\[0.8ex]
$P$ & chord $\gamma(s_2)-\gamma(s_1)$\\[0.8ex]
$d_i$ & $1-\cos(\theta_i)\kappa(s_i)$\\[0.8ex]
$F^{[\varepsilon]}_1, F^{[\varepsilon]}_2$ & functions of outward-pointing constraints for chord\\[0.8ex]
$S_K$ & standard set, i.e., $(\R/T\mathbb{Z})^2$ without weakly special point $W_{(\overline{s}_1, \overline{s}_2)}$ and weakly diagonal points\\[0.8ex]
$\kappa_M, \kappa_m, w_M, w_m$ & principal curvatures and principal directions\\[0.8ex]
$\langle\cdot, \cdot\rangle, \widetilde{\langle\cdot, \cdot\rangle}$ & pullback and flat metrics on $(\R/T\mathbb{Z}\times S^1)^2$\\[0.8ex]
$G_{K, \varepsilon}$-strip & Definition \ref{defn_g_strip}\\[0.8ex]
$\mathcal{R}$ & space of chords intersecting $K$ in its interior\\[0.8ex]
$\widehat{\widehat{\mathcal{R_\varepsilon}}}$ & space of chords intersecting $T_{K, \varepsilon}$ in its interior, but the intersection is not too close to the endpoints, Definition \ref{defn_evE}\\[0.8ex]
$W^u_{sing}$ & singular unstable manifold, Definition \ref{defn_sing_unstabl}\\[0.8ex]
$\cdot_{f\hbox{-}s\,}t, \fs\tau$ & fast-slow flow in slow/fast time scales\\[0.8ex]
$B_{\delta, U_{fen}}$ & box given by the codomain of Fenichel chart\\[0.8ex]
$X_H, R$ & Hamiltonian and Reeb vector fields\\[0.8ex]
$\vert \bm{a}\vert$ & Reeb chord degree\\[0.8ex]
$\mathcal{M}^{sy}_{T_K}, \mathcal{M}_{T_K}$ & Moduli spaces of $J$-holomorphic curves in the symplectization and with switches, respect.\\

\end{tabular}\\

\renewcommand*{\bibfont}{\scriptsize}
\setlength{\bibsep}{1pt plus 0.4ex}
\bibliography{priklady_literatury}{}

\openright

\end{document}